\documentclass[11pt]{amsart}

\usepackage[alphabetic]{amsrefs} 
\usepackage[letterpaper,top=1in,bottom=1in,left=1in,right=1in]{geometry}
\usepackage{amsthm,amsfonts,amssymb,amsmath,amsxtra,url}

\usepackage{mathrsfs}

\usepackage{xr-hyper}
\usepackage[pdfencoding=auto, psdextra,colorlinks=true,linkcolor=blue, citecolor=teal!50!green,urlcolor=cyan]{hyperref}
\usepackage{verbatim}

\usepackage{tikz-cd}
\usepackage{array,booktabs}
\usepackage{xspace}
\usepackage[all,cmtip]{xy}
\usepackage{dynkin-diagrams}
\usepackage{enumerate}
\usepackage{bm}
\usepackage{enumitem}
\usepackage{xcolor}
\usepackage{multicol}

\newtheorem{lem}{Lemma}[section]
\newtheorem{definition}[lem]{Definition}
\newtheorem{thm}[lem]{Theorem}
\newtheorem{rk}[lem]{Remark}
\newtheorem{prop}[lem]{Proposition}
\newtheorem{cor}[lem]{Corollary}
\newtheorem{convention}[lem]{Convention}
\newtheorem{construction}[lem]{Construction}

\newtheorem{Definition and Proposition}[lem]{Definition/Proposition}
\numberwithin{equation}{section}

\newtheorem*{lemS}{Lemma}
\newtheorem*{definitionS}{Definition}

\newtheorem*{rkS}{Remark}
\newtheorem*{propS}{Proposition}

\newtheorem*{thmA}{Theorem A}
\newtheorem*{thmB}{Theorem B}

\newenvironment{proofof}[1][]{{\noindent\it Proof of  {#1}.}\hspace{.5em}}{\hfill $\square$\par}

\newcommand{\bb}[1]{\mathbb{#1}}

\newcommand{\ca}[1]{\mathcal{#1}}

\newcommand{\und}[1]{\und{#1}}
\newcommand{\sbst}{\subset}
\newcommand{\spst}{\supset}
\newcommand{\lra}{\longrightarrow}

\newcommand{\mbf}[1]{\mathbf{#1}}

\newcommand{\iso}{\cong}
\newcommand{\mrm}[1]{\mathrm{#1}}
\newcommand{\wat}{\widehat}

\newcommand{\et}{\text{{\'e}t}}

\newcommand{\dr}{\mathrm{dR}}

\newcommand{\der}{\mrm{der}}
\newcommand{\ad}{\mrm{ad}}
\newcommand{\zpp}{ZP_\Phi}
\newcommand{\x}{\ca{X}}
\newcommand{\sh}{\mrm{Sh}}
\newcommand{\wdtd}{\widetilde}

\newcommand{\A}{\bb{A}_f}
\newcommand{\bss}{\backslash}
\newcommand{\disju}{\coprod}
\newcommand{\stb}{\mrm{Stab}}
\newcommand{\phip}{\Phi^\prime}
\newcommand{\zbkp}{{\bb{Z}_{(p)}}}

\DeclareMathOperator{\lie}{\mrm{Lie}}
\newcommand{\ull}{\underline{\mrm{Lie}}}
\newcommand{\mmin}{\mrm{min}}
\newcommand{\cusp}{\mrm{Cusp}}
\newcommand{\bd}[4]{\xrightarrow[#4]{({#1},{#2})_{#3}}}

\newcommand{\dpr}{{\prime\prime}}
\newcommand{\Ap}{\bb{A}_f^p}
\newcommand{\tw}[1]{{#1}^{\ca{P}}}
\newcommand{\coker}{\operatorname{coker}}
\newcommand{\tg}{\tilde{\gamma}}
\DeclareMathOperator{\isom}{\underline{Isom}}
\newcommand{\zbkpt}{\bb{Z}_{(p)}^\times}
\newcommand{\zhp}{\wat{\bb{Z}}^p}
\newcommand{\zbkppt}{\bb{Z}_{(\square)}^\times}
\newcommand{\V}{\mbf{V}}
\newcommand{\ion}{\iota^{\tg}_N}
\newcommand{\ionx}{\iota^{\tg}_{N,x}}
\newcommand{\ff}{\mbf{f}^{\tg}_N}
\newcommand{\ffab}{\mbf{f}^{\tg,ab}_N}
\newcommand{\ffsab}{\mbf{f}^{\tg,sab}_N}
\newcommand{\ffet}{\mbf{f}^{\tg,et}_N}
\newcommand{\fftr}{\mbf{f}^{\tg,tr}_N}

\DeclareMathOperator{\spec}{\mrm{Spec}}
\newcommand{\gal}{\mrm{Gal}}
\newcommand{\aut}{\underline{\mrm{Aut}}}
\newcommand{\Aut}{\mrm{Aut}}
\newcommand{\Endo}{\mrm{End}}
\newcommand{\uend}{\ul{\mrm{End}}}
\newcommand{\Hom}{\mrm{Hom}}
\newcommand{\uhom}{\ul{\mrm{Hom}}}
\newcommand{\G}{\ca{G}}
\newcommand{\hol}{\mrm{hol}}
\DeclareMathOperator{\gr}{\mrm{Gr}}
\newcommand{\lcj}[2]{{}^{#2}{#1}}
\newcommand{\rcj}[2]{{#1}^{#2}}
\DeclareMathOperator{\Frac}{\mrm{Frac}}
\newcommand{\W}{\ca{W}}
\newcommand{\Q}{\ca{Q}}
\newcommand{\T}{\ca{T}}
\newcommand{\bml}{\bm{\lambda}}
\newcommand{\pol}{\mrm{pol}}

\DeclareMathOperator{\spf}{\mrm{Spf}}

\DeclareMathOperator{\im}{\mrm{Im}}
\DeclareMathOperator{\Ext}{\mrm{Ext}}
\newcommand{\id}{\mrm{Id}}
\DeclareMathOperator{\Pic}{\mrm{Pic}}
\DeclareMathOperator{\biext}{\mrm{Biext}}
\newcommand{\ul}[1]{\underline{#1}}
\newcommand{\e}{\mrm{e}}
\newcommand{\PEL}{\mrm{PEL}}

\usepackage[OT2,T1]{fontenc}
\DeclareSymbolFont{cyrletters}{OT2}{wncyr}{m}{n}
\DeclareMathSymbol{\Sha}{\mathalpha}{cyrletters}{"58}
\newcommand{\mxh}{\mbf{M}_{\mrm{MH}}}
\newcommand{\imxh}{\bb{M}_{\mrm{MH}}}

\newcommand{\can}{\mrm{can}}
\newcommand{\p}{\square}
\newcommand{\zbkpp}{\bb{Z}_{(\square)}}
\newcommand{\zhpp}{\wat{\bb{Z}}^\square}
\DeclareMathOperator{\chara}{\mrm{char}}
\newcommand{\zhppt}{\wat{\bb{Z}}^{\square,\times}}
\newcommand{\App}{\A^\square}
\newcommand{\z}{\mbf{z}}
\newcommand{\oo}{{\ca{O}_F}}
\newcommand{\op}{{\ca{O}_{F,(p)}}}
\DeclareMathOperator{\ob}{\mrm{Ob}}
\DeclareMathOperator{\mor}{\mrm{Mor}}

\newcommand{\pp}{\mathscr{P}}
\newcommand{\psii}{\bm{\psi}}
\newcommand{\fo}{\mrm{for}}
\newcommand{\red}{\mrm{red}}
\newcommand{\zn}{\lcj{Z}{(n)}}
\newcommand{\zf}{\lcj{Z}{f}}

\newcommand{\wn}{\lcj{W}{(n)}}
\newcommand{\wf}{\lcj{W}{f}}
\newcommand{\wj}{\lcj{W}{j}}
\newcommand{\nt}{\natural}
\newcommand{\cc}{{\Phi,\sigma}}
\newcommand{\cpl}[2]{{(#1)}{}^{\wedge{\ }}_{#2}}

\newcommand{\ag}{\mathscr{A}}
\newcommand{\agsb}{\mathscr{A}^\circ}
\newcommand{\atd}{\wdtd{\mathscr{A}}}
\newcommand{\atdsb}{\wdtd{\mathscr{A}}^\circ}

\newcommand{\eg}{\mathscr{E}}

\newcommand{\K}{\wdtd{K}}
\newcommand{\qbar}{\overline{\bb{Q}}}

\newcommand{\jg}{\mathscr{J}}

\newcommand{\fj}{\underline{\mrm{FJ}}}

\newcommand{\geom}[1]{\bar{#1}}

\usepackage[symbols,nogroupskip,toc=false]{glossaries-extra}
\usepackage{glossary-mcols}
\setglossarystyle{mcolindex}
\glsxtrRevertTocMarks
\makenoidxglossaries
\glsxtrnewsymbol[description={\hfill a cusp label representative \hfill}]{Phi}{\ensuremath{\Phi}}

\glsxtrnewsymbol[description={\hfill Pink's $P_1$ and Madapusi's $Q_\Phi$\hfill}]{PPhi}{\ensuremath{P_\Phi}}

\glsxtrnewsymbol[description={\hfill Pink's $Q$ and Madapusi's $P_\Phi$\hfill}]{QPhi}{\ensuremath{Q_\Phi}}

\glsxtrnewsymbol[description={\hfill Pink's $W_1$ and Madapusi's $U_\Phi$\hfill}]{WPhi}{\ensuremath{W_\Phi}}

\glsxtrnewsymbol[description={\hfill Pink's $U_1$ and Madapusi's $W_\Phi$\hfill}]{UPhi}{\ensuremath{U_\Phi}}

\glsxtrnewsymbol[description={\hfill}]{H0}{\ensuremath{H_0}}

\glsxtrnewsymbol[description={\hfill}]{uQx}{\ensuremath{u^{Q}_x}}

\glsxtrnewsymbol[description={\hfill}]{DPhi}{\ensuremath{D_\Phi}}

\glsxtrnewsymbol[description={\hfill}]{hbar}{\ensuremath{\hbar}}

\glsxtrnewsymbol[description={\hfill rational boundary component \hfill}]{rbc}{\ensuremath{(P_\Phi,D_\Phi)}}

\glsxtrnewsymbol[description={\hfill}]{shKPhiF}{\ensuremath{\sh_{K_\Phi,F}}}

\glsxtrnewsymbol[description={\hfill}]{solshKPhiF}{\ensuremath{\overline{\sh}_{K_\Phi,F}}}

\glsxtrnewsymbol[description={\hfill}]{shKPhihF}{\ensuremath{\sh_{K_\Phi,h,F}}}

\glsxtrnewsymbol[description={\hfill}]{CLRGX}{\ensuremath{\ca{CLR}(G,X)}}

\glsxtrnewsymbol[description={\hfill}]{CuspKGX}{\ensuremath{\mrm{Cusp}_K(G,X)}}

\glsxtrnewsymbol[description={\hfill}]{CuspKGXSigma}{\ensuremath{\mrm{Cusp}_K(G,X,\Sigma)}}

\glsxtrnewsymbol[description={\hfill}]{IQ}{\ensuremath{I(Q)}}

\glsxtrnewsymbol[description={\hfill}]{Ep}{\ensuremath{E^p}}

\glsxtrnewsymbol[description={\hfill}]{tau}{\ensuremath{\tau}}

\glsxtrnewsymbol[description={\hfill}]{zmrmZPhiKC}{\ensuremath{\mrm{Z}_{[\Phi_i],K}(\bb{C})}}

\glsxtrnewsymbol[description={\hfill}]{DeltaPhiK}{\ensuremath{\Delta_{\Phi_i,K}}}

\glsxtrnewsymbol[description={\hfill}]{zmrmZPhiK}{\ensuremath{\mrm{Z}_{[\Phi],K}}}

\glsxtrnewsymbol[description={\hfill}]{DeltaPhiKcirc}{\ensuremath{\Delta_{\Phi,K}^\circ}}

\glsxtrnewsymbol[description={\hfill}]{zmrmZUpsilon}{\ensuremath{\mrm{Z}_{\Upsilon,K}}}

\glsxtrnewsymbol[description={\hfill}]{shKPhisubsigma}{\ensuremath{\sh_{K_\Phi,\sigma}}}

\glsxtrnewsymbol[description={\hfill}]{shKPhibksigma}{\ensuremath{\sh_{K_\Phi}(\sigma)}}

\glsxtrnewsymbol[description={\hfill}]{shKPhibkSigma}{\ensuremath{\sh_{K_\Phi}(\Sigma)}}

\glsxtrnewsymbol[description={\hfill}]{PPhi+}{\ensuremath{\mbf{P}_{\Phi}^+}}

\glsxtrnewsymbol[description={\hfill}]{PP}{\ensuremath{\mbf{P}_\Phi}}

\glsxtrnewsymbol[description={\hfill cocharacter group\hfill}]{SKPhivee}{\ensuremath{\mbf{S}_{K_\Phi}^\vee}}

\glsxtrnewsymbol[description={\hfill}]{mPhi}{\ensuremath{m(\Phi)}}

\glsxtrnewsymbol[description={\hfill}]{EKPhi}{\ensuremath{\mbf{E}_{K_\Phi}}}

\glsxtrnewsymbol[description={\hfill}]{LambdaKPhi}{\ensuremath{\mbf{\Lambda}_{K_\Phi}}}

\glsxtrnewsymbol[description={\hfill}]{zpp}{\ensuremath{ZP_\Phi}}

\glsxtrnewsymbol[description={\hfill}]{simZP}{\ensuremath{\sim_{ZP}}}

\glsxtrnewsymbol[description={\hfill}]{CuspZPKGXSigma}{\ensuremath{\mrm{Cusp}^{ZP}_K(G,X,\Sigma)}}

\glsxtrnewsymbol[description={\hfill}]{CuspZPKGX}{\ensuremath{\mrm{Cusp}^{ZP}_K(G,X)}}

\glsxtrnewsymbol[description={\hfill}]{DeltaPhiZP}{\ensuremath{\Delta^{ZP}_{\Phi,K}}}

\glsxtrnewsymbol[description={\hfill}]{DeltaZPcircPhiK}{\ensuremath{\Delta^{ZP,\circ}_{\Phi,K}}}

\glsxtrnewsymbol[description={\hfill}]{ZPPhi}{\ensuremath{[{ZP}(\Phi)]}}

\glsxtrnewsymbol[description={\hfill}]{ZPPhisigma}{\ensuremath{[{ZP}(\Phi,\sigma)]}}

\glsxtrnewsymbol[description={\hfill}]{fjKPhi2l}{\ensuremath{\fj_{K_{\Phi_2}}(l)}}

\glsxtrnewsymbol[description={\hfill}]{Hc}{\ensuremath{H^c}}

\glsxtrnewsymbol[description={\hfill 1-motives \hfill}]{Qca}{\ensuremath{\ca{Q}}}

\glsxtrnewsymbol[description={\hfill}]{IGG0}{\ensuremath{I_{G/G_0}}}

\glsxtrnewsymbol[description={\hfill}]{Ec}{\ensuremath{E^c}}

\glsxtrnewsymbol[description={\hfill}]{Ewdtd}{\ensuremath{\wdtd{E}}}

\glsxtrnewsymbol[description={\hfill}]{FKZ}{\ensuremath{F_{K_Z}}}

\glsxtrnewsymbol[description={\hfill}]{EKZwdtd}{\ensuremath{\wdtd{E}_{K_Z}}}

\glsxtrnewsymbol[description={\hfill}]{Eprime}{\ensuremath{E'}}

\glsxtrnewsymbol[description={\hfill}]{E0}{\ensuremath{E_0}}

\glsxtrnewsymbol[description={\hfill}]{E2}{\ensuremath{E_2}}

\glsxtrnewsymbol[description={\hfill}]{nuGddag}{\ensuremath{\nu_{G^\ddag}}}

\glsxtrnewsymbol[description={\hfill Hodge-type datum\hfill}]{G0X0}{\ensuremath{(G_0,X_0)}}

\glsxtrnewsymbol[description={\hfill Siegel-type datum\hfill}]{GddagXddag}{\ensuremath{(G^\ddag,X^\ddag)}}

\glsxtrnewsymbol[description={\hfill}]{VbbZ}{\ensuremath{V_{\bb{Z}}}}

\glsxtrnewsymbol[description={\hfill}]{SKddag}{\ensuremath{\ca{S}_{K^\ddag}}}

\glsxtrnewsymbol[description={\hfill}]{shKddag}{\ensuremath{\sh_{K^\ddag}}}

\glsxtrnewsymbol[description={\hfill}]{rgddag}{\ensuremath{r(g^\ddag)}}

\glsxtrnewsymbol[description={\hfill}]{hgddag}{\ensuremath{h(g^\ddag)}}

\glsxtrnewsymbol[description={\hfill}]{SKPhiddag}{\ensuremath{\ca{S}_{K_{\Phi^\ddag}}}}

\glsxtrnewsymbol[description={\hfill}]{SKPhi0}{\ensuremath{\ca{S}_{K_{\Phi_0}}}}

\glsxtrnewsymbol[description={\hfill}]{ZcaUpsilon}{\ensuremath{\ca{Z}_{\Upsilon,K}}}

\glsxtrnewsymbol[description={\hfill}]{SK0}{\ensuremath{\ca{S}_{K_0}}}

\glsxtrnewsymbol[description={\hfill}]{SK0Sigma0}{\ensuremath{\ca{S}_{K_0}^{\Sigma_0}}}

\glsxtrnewsymbol[description={\hfill}]{SKddagSigmaddag}{\ensuremath{\ca{S}_{K^\ddag}^{\Sigma^\ddag}}}

\glsxtrnewsymbol[description={\hfill}]{EtdKPhi}{\ensuremath{\mbf{E}_{\K_\Phi}}}

\title
[Compactifications of Integral Models of Abelian Type]  
{Arithmetic compactifications of integral models of Shimura varieties of abelian type}
\author{Peihang Wu}
\address{BICMR, Peking University, Beijing 100871, China} 
\email{wuph@pku.edu.cn} 
\date{}

\begin{document}

\begin{abstract}
    In this paper, we construct good toroidal and minimal compactifications in the sense of Lan-Stroh for integral models of abelian-type Shimura varieties. 
    We start with finding suitable types of cusp labels and cone decompositions which are compatible with those of the associated Hodge-type Shimura varieties.
    We then study the action of $\bb{Q}$-points of the adjoint group on boundary charts and toroidal compactifications of Hodge-type integral models. In particular, we extend the twisting construction of Kisin and Pappas to boundary charts.
    Finally, up to taking refinements of cone decompositions, we construct an abelian-type toroidal compactification as an open and closed algebraic subspace of a quotient from a disjoint union of Hodge-type toroidal compactifications and construct minimal compactifications with a similar method. Furthermore, we show results on nearby cycles of these compactifications and verify Pink's formula when the level at $p$ is an intersection of $n$ quasi-parahoric subgroups.
\end{abstract}

\maketitle     
\tableofcontents 

\section*{Introduction}
In this paper, we will study the Shimura varieties associated with abelian-type Shimura data and construct compactifications of certain integral models associated with them.\par
To recall the definitions, a pair $(G,X)$ of a connected reductive $\bb{Q}$-group $G$ together with a $G(\bb{R})$-conjugacy class of homomorphisms $\bb{S}:=R_{\bb{C}/\bb{R}}\bb{G}_m\to G_\bb{R}$ is called a \emph{Shimura datum} if it satisfies \cite[2.1.1.1-2.1.1.3]{Del79} (see also \cite[Def. 5.5]{Mil04}). 
Choosing a neat open compact subgroup $K\sbst G(\A)$, the complex double coset space $\sh_K(G,X)(\bb{C}):=G(\bb{Q})\bss X\times G(\A)/K$ uniquely algebraizes to a complex algebraic variety $\sh_K(G,X)_\bb{C}$ by \cite{BB66} and \cite{Bor72}.
Furthermore, by the works of Shimura, Deligne, Milne, Borovoi and others, the inverse system of algebraic varieties $\{\sh_K(G,X)_\bb{C}\}_{K}$ with varying neat open compact subgroups $K$ canonically descends to an inverse system of algebraic varieties $\{\sh_K(G,X)\}_K$ over a number field (called the \emph{reflex field}) $E(G,X)$ contained in $\bb{C}$, such that the Galois action on special points is given by a certain reciprocity law of the class field theory. \par
Abelian-type Shimura varieties and their integral models are important in number theory, and there are interesting problems and applications related to them. For an abelian-type Shimura datum $(G_2,X_2)$, the weight cocharacter and the center of the reductive group are more complicated than those of a Hodge-type Shimura datum; more precisely, the Axioms ($\mrm{SV4}$), ($\mrm{SV}5$) and ($\mrm{SV}6$) in \cite[Sec. 5]{Mil04} are not satisfied in general. To study an abelian-type Shimura variety, one can study a Hodge-type Shimura variety associated with it first; the two varieties are related by a quotient over geometrically connected components and by an \emph{induction construction} of Deligne \cite{Del79}.\par
Let $p$ be a prime number. Let $(G_2,X_2)$ be an abelian-type Shimura datum. By definition (see, e.g., \cite[Intro.]{Kis10}), it is associated with a Hodge-type Shimura datum $(G_0,X_0)$. Assume that $K_2=K_{2,p}K^p_2$, where $K_{2,p}$ is an open compact subgroup of $G_2(\bb{Q}_p)$ and $K_2^p$ is a neat open compact subgroup of $G_2(\Ap)$. When $G_2$ is quasi-split and unramified at $p$ and $K_{2,p}$ is hyperspecial, it was shown by Kisin and Kim-Madapusi that $\sh_{K_2}(G_2,X_2)$ extends to a \emph{smooth integral model} $\ca{S}_{K_2}$, and $\ca{S}_{K_{2,p}}:=\varprojlim \ca{S}_{K_{2,p}K^p_2}$ is \emph{canonical} in the sense that it satisfies \emph{the extension property} (see \cite[2.3.7]{Kis10} and \cite{KM15}).
When $p>2$ and $K_{2,p}$ is parahoric, Kisin, Pappas, and Zhou showed that $\sh_{K_{2,p}}(G_2,X_2)$ extends to an integral model $\ca{S}_{K_{2,p}}(G_2,X_2)$, which satisfies a weaker extension property and has a \emph{local model diagram} (see \cite{KP15} and \cite{KPZ24}).
Moreover, these integral models $\ca{S}_{K_{2,p}}$ are canonical in the sense that they satisfy the \emph{conjecture of Pappas and Rapoport} (see \cite[Conj. 4.2.2]{PR24}, \cite{DvHKZ24} and \cite[Thm. 4.10]{DY25}).
\subsubsection*{Background on compactifications}
Shimura varieties are not proper in general, and they admit smooth and projective toroidal compactifications and projective minimal compactifications over reflex fields (see \cite{BB66}, \cite{AMRT10} and \cite{Pin89}). In many occasions, studying the integral models of these compactifications is very beneficial.\par 
It is conjectured that smooth integral models of Shimura varieties of abelian type have toroidal and minimal compactifications with the same properties mentioned above extending those over the reflex fields (see \cite[Conj. 2.18]{Mil92}). Such compactifications of smooth integral models were constructed by Chai and Faltings \cite{FC90} in the Siegel-type case, by Lan \cite{Lan13} in the PEL-type case and by Madapusi \cite{Mad19} in the Hodge-type case.\par 
In fact, one can expect a similar story for more general integral models. Let us fix a prime number $p$ and continue with the conventions above. As before, $(G_0,X_0)$ by definition admits an embedding (called Hodge embedding in this paper) $\iota:(G_0,X_0)\hookrightarrow (G^\ddag,X^\ddag)$, i.e., an injective group homomorphism $G_0\hookrightarrow G^\ddag:=\mrm{GSp}(V,\psi)$, which induces an embedding from $X_0$ to the union of Siegel upper and lower half-spaces $X^\ddag:=\bb{S}^\pm$. Suppose that we choose a self-dual lattice $V_{\bb{Z}_p}$ of $V_{\bb{Q}_p}$ and an open compact subgroup $K_{0,p}\sbst G_0(\bb{Q}_p)$ such that $K_{0,p}\sbst K_p^\ddag:=\stb_{G^\ddag(\bb{Q}_p)}(V_{\bb{Z}_p}).$ 
Let $K_0^p$ be a neat open compact subgroup of $G_0(\Ap)$ contained in a neat open compact subgroup $K^{\ddag,p}$ of $G^\ddag(\Ap)$. Set $K^\ddag:=K_p^\ddag K^{\ddag,p}$ and $K_0:=K_{0,p}K^{p}_0$. Let $E_0$ be the reflex field of $(G_0,X_0)$. Let $\ca{O}:=\ca{O}_{E_0,(v)}$ for a place $v|p$. Let $\ca{S}_{K_0}$ be the normalization of $\ca{S}_{K^\ddag,\ca{O}}$ in $\sh_{K_0}:=\sh_{K_0}(G_0,X_0)$, where $\ca{S}_{K^\ddag}$ is the integral model of $\sh_{K^\ddag}(G^\ddag,X^\ddag)$ over $\zbkp$ constructed as a moduli space of polarized abelian schemes with level structures. Then $\ca{S}_{K_0}$ is a quasi-projective normal scheme that is flat over $\ca{O}$, and its generic fiber is $\sh_{K_0}$. 
For a choice of Hodge embedding 
$$(G_0,X_0,K_{0,p})\hookrightarrow(G^\ddag,X^\ddag,K_p^\ddag)$$ such that $K_p^\ddag\cap G_0(\bb{Q}_p)\neq K_{0,p}$, or a choice of open compact subgroup $K_{0,p}$ beyond quasi-parahoric subgroups, the scheme $\ca{S}_{K_{0,p}}:=\varprojlim_{K_0^p}\ca{S}_{K_0}$ is \emph{not} canonical in the sense of satisfying the extension property or the conjecture of Pappas and Rapoport.\par 
Despite this, since \cite{Lan16b} and \cite{Mad19}, it has been understood that there are good toroidal and minimal compactifications associated with reasonably defined integral models of Shimura varieties of Hodge type, even when the integral models are not constructed with (quasi-)parahoric level structures or with any definition of canonicity mentioned above.\par 
In fact, Madapusi proved that \emph{all} Hodge-type integral models with hyperspecial or (quasi-)parahoric levels studied in \cite{Kis10}, \cite{KP15}, \cite{KPZ24}, \cite{PR24} and \cite{DvHKZ24} have good toroidal compactifications, and his results can also be applied in a broader context beyond quasi-parahoric levels with minimal restrictions on $K_{0,p}\sbst K_p^\ddag$ (see \cite[Sec. 3.1]{Mad19}). Let $\Sigma^\ddag$ be a smooth admissible rational polyhedral cone decomposition for $(G^\ddag,X^\ddag,K^\ddag)$. By \cite{FC90} and \cite{Lan13}, there is a smooth toroidal compactification $\ca{S}_{K^\ddag}^{\Sigma^\ddag}$ for $\ca{S}_{K^\ddag}$. Let $\Sigma_0$ be the cone decomposition for $(G_0,X_0)$ induced by $\Sigma^\ddag$. Define $\ca{S}_{K_0}^{\Sigma_0}$ to be the normalization of $\ca{S}_{K^\ddag}^{\Sigma^\ddag}$ in $\sh_{K_0}$. 
Madapusi showed that $\ca{S}_{K_0}^{\Sigma_0}$ has the right properties as a toroidal compactification, that is, $\ca{S}_{K_0}^{\Sigma_0}$ has a good stratification and the complete local rings of $\ca{S}^{\Sigma_0}_{K_0}$ at its boundary points are complete localizations of toric torsors, which can be explicitly described by integral models of mixed Shimura varieties (see \cite[Thm. 3.4.3 and Thm. 4.1.5]{Mad19}). The generic fiber of $\ca{S}_{K_0}^{\Sigma_0}$ is $\sh_{K_0}^{\Sigma_0}$, the toroidal compactification of $\sh_{K_0}$ associated with the cone decomposition $\Sigma_0$ by Pink's characteristic $0$ theory \cite{Pin89}. Similar results were also proved by Lan in \cite{Lan16b} when $(G_0,X_0)$ is of PEL type.\par
Hence, one expects that, for all abelian-type Shimura data, a good compactification theory should exist for integral models defined in an appropriate way with arbitrary levels at $p$. 
\subsubsection*{Main results}
Let $\Sigma_2$ be an admissible rational polyhedral cone decomposition for $\sh_{K_2}:=\sh_{K_2}(G_2,X_2)$. Denote by $\sh_{K_2}^{\Sigma_2}$ the toroidal compactification of $\sh_{K_2}$ associated with $\Sigma_2$ by Pink's theory. Denote $E_2:=E(G_2,X_2)$. Fix any place $v_2|p$ of $E_2$.
The main goal of this paper is to construct good toroidal compactifications $\ca{S}_{K_2}^{\Sigma_2'}$ over $\ca{O}_{E_2,v_2}$ associated with suitable refinements $\Sigma_2'$ of $\Sigma_2$ based on Madapusi's results in the Hodge-type case.\par 
The main results can be stated in the following compressed form. The readers will find precise statements in Theorem \ref{maintheorem}.
\begin{thmA}[{Theorem \ref{maintheorem}}]
Assume that $K_{2,p}$ is an open compact subgroup of $G_2(\bb{Q}_p)$ and $K_2^p$ is neat open compact in $G_2(\Ap)$. Let $K_2=K_{2,p}K_2^p$. Let $v_2|p$ be a place of $E_2$ over $p$.\par
For any admissible cone decomposition $\Sigma_2$ for $(G_2,X_2)$ and $K_2$, there is a refinement $\Sigma_2'$ of $\Sigma_2$, which can be made smooth and/or projective, such that there is a proper normal flat model $\ca{S}_{K_2}^{\Sigma_2^\prime}$ over $\ca{O}_{E_2,v_2}$ extending the toroidal compactification $\sh_{K_2}^{\Sigma'_2}$ constructed by Pink in \cite{Pin89}. 
The compactification $\ca{S}^{\Sigma_2'}_{K_2}$ has a good stratification, which extends the stratification obtained from Pink's theory from the generic fiber to the integral model.\par
There is an open dense subscheme $\ca{S}_{K_2}$ in $\ca{S}_{K_2}^{\Sigma_2'}$ extending $\sh_{K_2}$, and the reduced complement $D$ of $\ca{S}_{K_2}$ in $\ca{S}^{\Sigma_2'}_{K_2}$ is a relative effective Cartier divisor. The complete local ring of $\ca{S}^{\Sigma_2'}_{K_2}$ at any point $x$ of $D$ is isomorphic to a complete local ring of a toric scheme. These toric schemes are constructed from finite quotients of the toric schemes associated with certain integral models of boundary mixed Shimura varieties of abelian type.\par
When $G_2$ is quasi-split and unramified at $p$, and $K_{2,p}$ is hyperspecial, the towers of integral models of boundary mixed Shimura varieties mentioned in the last paragraph can be constructed to satisfy the extension property. In this case, $\ca{S}^{\Sigma_2'}_{K_2}$ and its strata are defined over $\ca{O}_{E_2,(v_2)}$.
\end{thmA}
We can also construct a minimal (i.e., Satake-Baily-Borel) compactification $\ca{S}_{K_2}^\mmin$ for $\ca{S}_{K_2}$. More precisely, 
\begin{thmB}[{Theorem \ref{thm-main-theorem-min}}]
There is a normal projective model $\ca{S}_{K_2}^\mmin$ over $\ca{O}_{E_2,v_2}$ (or over $\ca{O}_{E_2,(v_2)}$ if $G_2$ is quasi-split and unramified at $p$, and $K_{2,p}$ is hyperspecial) extending the minimal compactification $\sh_{K_2}^\mmin$ of $\sh_{K_2}$ over $E_2$. There is a proper surjective morphism $\oint_{K_2}^{\Sigma_2'}: \ca{S}^{\Sigma_2'}_{K_2}\to \ca{S}_{K_2}^\mmin$ with geometrically connected fibers, under which the stratifications on the source and the target compatibly extend the stratifications on the generic fibers. 
\end{thmB}
As the aforementioned works \cite{Lan16b} and \cite{Mad19} by Lan and Madapusi, our results also allow the level at $p$ to be arbitrarily high (and low). More precisely, to apply our results, one can pick an open compact subgroup $K_{2,p}$ and choose $n$ Bruhat-Tits stabilizer subgroups $K_{2,p}^i$ containing $K_{2,p}$ from possibly distinct apartments (which means that $K_{2,p}$ can be any proper open compact subgroup of $\cap_{i=1}^n K_{2,p}^i$), then a good compactification theory can be established by the constructions in this paper. Note that the compactifications and their stratifications obtained will depend on the choice of the set of Bruhat-Tits stabilizers, since in our article the models of the compactifications and stratifications at level $K_{2,p}K_2^p$ are constructed by taking relative normalizations from the models at level $(\cap_{i=1}^n K_{2,p}^i) \cdot K_2^p$, and the choice of the latter level structure can be flexible. \par
\begin{rkS}
Since the integral models of compactifications in this paper have good stratifications, one can take the open dense strata (i.e., the interiors) of compactifications to get integral models of Shimura varieties.\par 
In particular, this paper constructed integral models of abelian-type Shimura varieties with quasi-parahoric level structures at $p>0$. In a forthcoming work joint with Shengkai Mao (see Remark \ref{rk-future-work}), we show these models are canonical integral models in the sense of Pappas and Rapoport. In fact, we show that $p$-adic shtukas admit ``canonical extensions'' to integral models of toroidal compactifications of abelian-type Shimura varieties.
\end{rkS}
Another advantage of our results is that we have proved that the compactifications of abelian-type integral models also satisfy (a slightly generalized version of) the qualitative descriptions of good compactifications in the sense of Lan-Stroh (see, e.g., \cite[Prop. 2.2]{LS18i}). The only difference is that the ``$\Gamma_\sigma$'' in the notation of \emph{loc. cit.} is in general not trivial in our case. In other words, in the Hodge-type or PEL-type case, the formal completions at the boundary can be described by (integral models of) mixed Shimura varieties defined by some irreducible mixed Shimura data, while for the abelian-type case, those completions are (also explicitly) described by \emph{finite quotients} of (integral models of) those mixed Shimura varieties (see Remark \ref{rk-field-of-def}). Hence, one can immediately extend many results from the Hodge-type case to the abelian-type case, once the arguments to prove them only involve a good description of the boundary.\par
As an application, we generalize the nearby cycle results in \cite{LS18} to the abelian-type case at all levels.
\begin{propS}[{Proposition \ref{prop-nearby} and Corollary \ref{cor-nearby-intersection}}]
    Fix a prime $l\neq p$. 
Let $\ca{V}$ be a lisse $\overline{\bb{Q}}_l$-sheaf associated with an algebraic representation $\xi$ of $G^c_2$ on a finite-dimensional $\overline{\bb{Q}}_l$-vector space $V_\xi$ or a finite $\bb{F}_l$-sheaf equipped with an action of an open compact subgroup $G_2^c(\bb{Z}_l)$.\par
Assume that the projection of $K_2$ into $G^c_2(\bb{Q}_l)$ factors through $G^c_2(\bb{Z}_l)$. Then there are natural isomorphisms induced by adjunctions that are equivariant under the actions of the absolute Galois group $\gal(\overline{\eta}/\eta)$:
\begin{enumerate}
    \item $R\Psi_{\ca{S}_{K_2}^{\Sigma_2}}RJ_{\eta,*}^{\Sigma_2}\ca{V}\xrightarrow{\sim} RJ_{\overline{s},*}^{\Sigma_2}R\Psi_{\ca{S}_{K_2}}\ca{V}$ and $J_{\overline{s},!}^{\Sigma_2}R\Psi_{\ca{S}_{K_2}}\ca{V}\xrightarrow{\sim}R\Psi_{\ca{S}_{K_2}^{\Sigma_2}}J_{\eta,!}^{\Sigma_2}\ca{V}$; similar results also hold for minimal compactifications.
    \item $R\Psi_{\ca{S}_{K_2}^\mmin}J_{,\eta,!*}^{\mmin}\ca{V}[d]\xrightarrow{\sim}J_{\overline{s},!*}^\mmin R\Psi_{\ca{S}_{K_2}}\ca{V}[d]$.
\end{enumerate}
\end{propS}
As another application, we prove Pink's formula when the level at $p$ is an intersection of $n$ quasi-parahoric subgroups.
\begin{propS}[{See Proposition \ref{prop-pink-formula} for details}]
Fix a positive integer $n$.
Fix $n$ Bruhat-Tits stabilizer subgroups $K^i_{2,p}$ of $G_2(\bb{Q}_p)$. Let $K_{2,p}$ be an open compact subgroup of $G_2(\bb{Q}_p)$ such that $\cap_{i=1}^n K^{i,\circ}_{2,p}\sbst K_{2,p}\sbst \cap_{i=1}^n K^i_{2,p}$. Then Pink's formula holds for the minimal compactification $\ca{S}^\mmin_{K_2}$ constructed as the relative normalization from the integral model of the minimal compactification at level $\cap_{i=1}^n K^i_{2,p}$, where the latter one is constructed as in Case ($\mrm{STB}_n$) for $\{K_{2,p}^i\}_{i=1}^n$.
\end{propS}
\subsubsection*{Method.}
From the definition of abelian-type Shimura data, we make the following construction.
\begin{lemS}[{See \S\ref{sss-conclusion}}]
 There is a connected reductive group $G$ and two Shimura data $(G,X_a)$ and $(G,X_b)$, such that $G^\der=G_2^\der$, such that $X_a$ and $X_b$ differ only by a homomorphism from $\bb{S}$ to $Z_{G,\bb{R}}$ (which we denote by $(G,X_a)\sim (G,X_b)$ below to refer to this similarity), and such that there is a diagram 
 \begin{equation*}
    \begin{tikzcd}
&& (G_0,X_0)\arrow[d,"\pi^b"]\\
(G_2,X_2)\arrow[r,hook,"\pi^a"]& (G,X_a)\arrow[r,phantom,description,"\sim"]& (G,X_b)          \end{tikzcd}
\end{equation*} 
 where $\pi^a$ is an embedding such that $\pi^a(G_2^\der)=G^\der$ and $\pi^b$ is a map whose kernel is the kernel of the central isogeny $G_0^\der\to G^\der=G_2^\der$.   
\end{lemS}
This is a variant of the construction in \cite{Del79}. For the topic of this paper, we found some advantages of doing the construction above: First of all, since the difference between the descent data of the objects associated with $(G,X_a)$ and $(G,X_b)$ is essentially induced by a central homomorphism $c:\bb{S}\to Z_{G,\bb{R}}$ (which is explicitly determined by the reciprocity law mentioned above) and since there is a direct map $(G_0,X_0)\to (G,X_b)$, we do not need to break up the integral model of $\sh_{K_0}(G_0,X_0)$ into geometrically connected components to construct strata and compactifications for $(G_2,X_2)$. \par
Secondly, to construct integral models with the strategy above, we do not rely on certain parahoric $\bb{Z}_p$-group schemes associated with $G$, $G_0$ and $G_2$ to construct the groups used in the induction construction. This gives us the flexibility we need in order to find suitable integral models when the level structure is arbitrary.\par
Another advantage is that the morphism from $(G_2,X_2)$ to $(G,X_a)$ is an embedding. By choosing a suitable level $K\sbst G(\A)$, one can show that the morphisms induced by $(G_2,X_2)\hookrightarrow(G,X_a)$ between compactifications and strata are also embeddings. In particular, results for $(G_2,X_2)$ can be immediately deduced from those for $(G,X_a)$.\par
Let us now sketch some crucial steps in this article.
\subsubsection*{Cusp labels and cone decompositions.}
Let $(G,X)$ be any Shimura datum. For each cusp label representative $\Phi$ (see \cite{Pin89} and \cite{Mad19}), Pink associated with $\Phi$ a $\bb{Q}$-subgroup $P_\Phi$ (which is denoted by $P_1$ in \cite{Pin89}), which is a normal subgroup of the admissible $\bb{Q}$-parabolic $Q_\Phi$ associated with $\Phi$ (see \cite[Def. 4.5 and Prop. 4.6]{Pin89}).\par
Note that there is in general no direct map between $P_\Phi$'s for $(G_0,X_0)$ and $(G_2,X_2)$ since there is no morphism from $(G_0,X_0)$ to $(G_2,X_2)$.
To address the issue, let $ZP_\Phi$ be the identity component of the group generated by the center $Z_G$ of $G$ and the group $P_\Phi$.
We use the following weaker equivalence relation to create coarser classes of cusp labels. 
\begin{definitionS}[{Definition \ref{equi-z}}] Let $\Phi_1$ and $\Phi_2$ be two cusp label representatives. Then we say 
$\Phi_1\preceq\Phi_2$ if there are $q^\prime\in ZP_{\Phi_2}(\bb{A}_f)$ and $\gamma\in G(\bb{Q})$ such that
$\gamma P_{\Phi_1}\gamma^{-1}\sbst P_{\Phi_2}$, $P_{\Phi_1}(\bb{Q})\gamma X^+_{\Phi_1}=P_{\Phi_2}(\bb{Q})X^+_{\Phi_2}$ and $\gamma g_{\Phi_1}\equiv q^\prime g_{\Phi_2}$ modulo $K$.
\end{definitionS}
A stratum defined by the equivalence relation above is actually a disjoint union of some strata in Pink's definition; by considering this kind of disjoint union, we can uniformize the construction of strata and mixed Shimura varieties for $(G,X_a)$ and $(G,X_b)$.
\begin{rkS}
In fact, this is an example of an alternative definition of boundary components, where the mixed Shimura data produced are not irreducible (see \cite[Rmk. 4.11]{Pin89} and \S\ref{subsubsec-mixedshimura}). 
\end{rkS}
More specifically, we can find cone decompositions as follows:
\begin{propS}[{Proposition \ref{zp-cones}}]
We can choose a common cone decomposition $\Sigma$ for both $(G,X_a)$ and $(G,X_b)$, such that $\Sigma$ induces a cone decomposition $\Sigma_2'$, which is a refinement of $\Sigma_2$ that can be made smooth, projective, or both smooth and projective, and such that $\Sigma$ also induces cone decompositions $\Sigma_0^\alpha$ for $(G_0,X_0) $ under the pullback of $\pi^b(\alpha)$ for each $\alpha$ in an index set $I_{G/G_0}
:=\stb_{G(\bb{Q})}(X_0)\pi^b (G_0(\A))\bss G(\A)/K$.
\end{propS}

\subsubsection*{Construction.}
We now have obtained the following diagram:
\begin{equation*}
    \begin{tikzcd}
&& (G_0,X_0);\Sigma_0^\alpha\arrow[d]\\
(G_2,X_2);\Sigma_2'\arrow[r,hook]& (G,X_a);\Sigma\arrow[r,phantom,description,"\sim"]& (G,X_b);\Sigma.        
    \end{tikzcd}
\end{equation*}
Note that one cannot get a good theory if one applies Deligne's induction construction directly to toroidal compactifications; one reason is that one cannot choose smooth cone decompositions for general level $K_2$ if cone decompositions over all geometrically connected components of the Shimura variety $\sh(G_2,X_2)$ are chosen the same. Hence, we first construct the toroidal compactification $\ca{S}^{\Sigma}_K(G,X_b)$ of $(G,X_b)$ as a disjoint union of certain quotients of Hodge-type toroidal compactifications associated with possibly different cone decompositions $\Sigma_0^\alpha$. We then construct $\ca{S}_{K_2}^{\Sigma_2'}$ as an open and closed algebraic subspace of $\ca{S}_{K}^{\Sigma}(G,X_a)$. The latter one is isomorphic to $\ca{S}_K^{\Sigma}(G,X_b)$ after a possibly ramified normalized base change to a ring $\ca{O}_{K_Z}$, but the ramification of this normalized base change is controlled by $c$ and the intersection of the center $Z_G(\bb{Q}_p)$ and the level group $K_p\sbst G(\bb{Q}_p)$ at $p$. For example, this (normalized) base change is unramified if $Z_G(\bb{Q}_p)\cap K_p$ contains the parahoric subgroup of $Z_G(\bb{Q}_p)$.\par
We can show the following statement: 
\begin{propS}[{Corollary \ref{cor-maintheorem}}]\label{prop-kummer-etale}
After taking the normalized base change to $\ca{O}_{K_Z}$, $\ca{S}_{K_2,\ca{O}_{K_Z}}^{\Sigma_2'}$ is an open and closed subspace of $\ca{S}_{K,\ca{O}_{K_Z}}^{\Sigma}(G,X_b)$, and $\ca{S}_{K,\ca{O}_{K_Z}}^{\Sigma}(G,X_b)$, equipped with the log structure defined by the complement of $\ca{S}_{K,\ca{O}_{K_Z}}(G,X_b)$ in $\ca{S}_{K,\ca{O}_{K_Z}}^\Sigma(G,X_b)$, admits a finite Kummer {\'e}tale cover of disjoint union of Hodge-type toroidal compactifications $\ca{S}^{\mrm{TOR}}_{K_0,\ca{O}_{K_Z}}:=\disju_{\alpha\in I_{G/G_0}}\ca{S}^{\Sigma_0^\alpha}_{K_0^\alpha,\ca{O}_{K_Z}}$ (with $K_{0,p}^\alpha$ chosen depending on $K_p$). The induced map between interiors $\disju_{\alpha\in I_{G/G_0}}\ca{S}_{K_0^\alpha,\ca{O}_{K_{Z}}}\to\ca{S}_{K,\ca{O}_{K_Z}}(G,X_b)$ is finite {\'e}tale.
\end{propS}
Note that we show the statement above for any choice of $(G_0,X_0)$ attached to $(G_2,X_2)$ that satisfies a condition on reflex fields (see \cite[Lem. 4.6.22 (3)]{KP15}). \par
There are two main difficulties in the quotient construction.\par
The first difficulty is to find certain $\ag$- and $\agsb$-type groups for mixed Shimura varieties, with which one can construct integral models of toroidal compactifications and integral models of boundary mixed Shimura varieties as (disjoint unions of) quotients. 
Here we have to use a method that is different from \cite{Kis10} and \cite{KP15} since we should construct the integral models and describe them at every finite level $K_pK^p$ instead of describing at level $K_p$ with $K_p$ a parahoric subgroup.
An important input to achieve this is constructing a collection of integral models $\{\ca{S}_{\lcj{K_{\Phi_0^\alpha}}{\gamma}}\}$ such that the whole $G^\ad(\bb{Q})$ acts on it. 
In {Section} \ref{sec-tor-abelian}, our proof of this crucially relies on the canonicity results in \cite{Kis10}, \cite{KM15}, \cite{PR24}, \cite{DY25} and \cite{DvHKZ24} and a detailed computation of the level groups associated with the boundary components in \cite{Mao25}. 
Even though there might be other ways to obtain this, the advantage of this approach is that we can establish a stronger statement along the way on how the construction of the toroidal compactification $\ca{S}_{K_0^\alpha}^{\Sigma_0^\alpha}$ is independent of the choice of Hodge embedding with prescribed formation. \par
The second difficulty is to show that the integral models of boundary mixed Shimura varieties are correct torus torsors. 
To do this, we use twisting of $1$-motives to give an explicit moduli description of certain actions following the idea of Kisin in \cite{Kis10}. 
The main result on twisting can be described as follows:
\begin{propS}[{Proposition \ref{prop-twt-1}}]
Let $G_{0,\zbkp}$ (resp. $Z_\zbkp$) be the closure of $G_0$ (resp. $Z_{G_0}$) in $G^\ddag_\zbkp$. Let $G^\ad_\zbkp:=G_{0,\zbkp}/Z_\zbkp$. For any $\gamma \in G^\ad_{0,\zbkp}(\zbkp)_1$ and any cusp label representative $\Phi=(Q,X_0^+,g^\ddag)$ as in \cite{Mad19}, there is a twisted cusp label representative $\Psi$ defined as $\Psi:=(\gamma Q\gamma^{-1},\gamma(X_0^+),\gamma g^\ddag\gamma^{-1})$, such that the action $$\gamma^{-1}:\sh_{K_\Phi}(\bb{C})\xrightarrow{\sim}\sh_{K_{\Psi}}(\bb{C}),$$ 
sending $\wp=[(x,\mbf{p}g^\ddag)]$ to $\wp^\gamma=[(\gamma(x),\gamma \mbf{p}g^\ddag\gamma^{-1})]$, extends to a twisting construction between $1$-motives $$(\Q,\bml,[u]_{(K_\Phi)^{g^\ddag}})\mapsto (\Q^\gamma,\bml^\gamma,[u^{\gamma}]_{(K_\Phi)^{g^\ddag\gamma^{-1}}}).$$
This action is compatible with twisting of abelian schemes on $\ca{S}_K$ constructed by Kisin \cite{Kis10} and Kisin-Pappas \cite{KP15}, under the isomorphism of formal completions in \cite[Thm. 4.1.5 (5)]{Mad19} for Hodge-type Shimura varieties.
\end{propS}
With this, we show that the integral model of the boundary mixed Shimura variety associated with $ZP_\Phi$ for any cusp label representative $\Phi$ of $(G,X_a)$ is a torus torsor when $K_p\cap Z_G(\bb{Q}_p)$ contains the parahoric subgroup of $Z_G(\bb{Q}_p)$ (see Lemma \ref{lem-twt-free-action}, Lemma \ref{lem-twist-admissible} and Proposition \ref{prop-mixsh-canonicity}). Then we extend this result to arbitrary level structures essentially by using an explicit description of the torus part of the Siegel-type boundary mixed Shimura varieties (see \S\ref{subsubsec-main-3}).
\subsection*{Organization of the paper}
In {Section} \ref{sec-supplement}, we review Pink's characteristic $0$ theory for canonical models of toroidal compactifications of mixed Shimura varieties. Subsequently, we introduce a definition of coarser cusp labels, which helps us with finding compatible choices of cone decompositions for $(G_0, X_0)$ and $(G_2, X_2)$. This selection of cusp labels and cones also enables us to construct $\ca{S}_{K_2}^{\Sigma_2}$ as the pullback of some quotient in later steps. Some group-theoretic lemmas are established along the way.\par
In {Section} \ref{sec-1-motives-and-deg}, the main goal is to review the theory of degeneration of abelian schemes developed in \cite{Mum72}, \cite{FC90} and \cite{Lan13}. The main difference is that we review it in the language of $1$-motives (cf. \cite{Str10} and \cite{Mad19}); some propositions are recorded for later use.\par
In {Section} \ref{sec-twt-construction}, we extend Kisin-Pappas' construction of twisting abelian schemes up to prime-to-$p$ isogenies to the twisting of $1$-motives and semiabelian schemes with additional structures. We then explain that the twisting constructions on the interior $\ca{S}_{K_0}$ and on the boundaries are compatible and induce an action on toroidal compactifications.\par
In {Section} \ref{sec-tor-abelian}, we first generalize Deligne's induction to construct boundary mixed Shimura varieties associated with $\sh^{\Sigma_2'}_{K_2}$. We then construct $\ca{S}_{K_2}^{\Sigma_2'}$ as the ``subquotients'' described above and present the main theorem (see Theorem \ref{maintheorem}). \par
To complete the construction of a good compactification theory, we also construct the minimal compactification $\ca{S}_{K_2}^\mmin$ of $\ca{S}_{K_2}$. We finish the paper by proving Pink's formula and by showing that the adjunction morphisms of nearby cycles are isomorphisms (see {Section} \ref{sec-min-cpt}).\par \S\ref{subsec-zp-cusp}, \S\ref{subsec-twt-1-mot}, \S\ref{subsec-ext}, \S\ref{subsec-con-stra-comp-quo} and \S\ref{subsec-main-thm} are the main technical {subsections} of this paper.\par
A list of symbols appearing in the first and the second {sections} is attached as an appendix for the convenience of the readers.\par 
Here is a diagram showing the rough relationships among all sections and many key references in the literature. We apologize to the readers that we are not able to exhibit all the references we have used or to describe all the details within one diagram for now.
\begin{equation*}
    \begin{tikzcd}
    \fbox{\cite{Pin89}}\arrow[d]&\fbox{\cite{Del79}}\arrow[d]&\fbox{\cite{Kis10}\&\cite{KP15}}\arrow[d]&\fbox{\cite{FC90}\&\cite{Lan13}}\arrow[d]& \fbox{\cite{Mad19}}\arrow[dl]\arrow[d]\\
    \fbox{\S 1}\arrow[rrr,leftrightarrow,bend left=15, dashed,"{+App. \ref{cpr-pel-cl}\text{ compatible}}" description]\arrow[r]\arrow[ddr,"\S 1.4"]& \fbox{\S 4.1}\arrow[d]&\fbox{\S 3.1}\arrow[d]&\fbox{\S 2.1-2.4}\arrow[r]\arrow[dl]& \fbox{\S 2.5}\arrow[ddlll,bend left]\\
    & \fbox{\S 4.2}\arrow[d]& \fbox{\S 3.2-3.3}\arrow[l]\arrow[d]&&\\
    &\fbox{\S 4.3}\arrow[dd,bend right]\arrow[d]& \fbox{\S 3.4}\arrow[l]&&\\
    &\fbox{\S 5.1-5.3}\arrow[d]&&&\\
    \fbox{\cite{Pin92},\cite{LS18b},\cite{LS18}}\arrow[r]&\fbox{\S 5.4-5.5}&&&
    \end{tikzcd}
\end{equation*}
\subsection*{Acknowledgments}
This is a major revision of my dissertation submitted to the University of Minnesota in August 2024. 
I am very grateful to my academic advisor Kai-Wen Lan for his guidance. I would like to thank him for pointing out many errors in the preliminary drafts of this paper and providing invaluable advice on doing and writing math. He also encouraged me to generalize the results to Case ($\mathrm{STB}_n$) $(n>1)$ and its deeper levels in this revision. I would like to thank Shengkai Mao for many fruitful discussions throughout the years.\par 
Part of this work was presented at the Morningside Center of Mathematics in October 2024, and I am very grateful to Xu Shen for his invitation.  
Moreover, I want to thank Pol van Hoften, Kentaro Inoue, Kai-Wen Lan, Shengkai Mao, Georgios Pappas, Xu Shen, Sug Woo Shin, Liang Xiao, Alex Youcis, Zhiyu Zhang, and Yihang Zhu for their interest and encouragement during the preparation of the paper and many helpful communications.
\subsection*{Notation and Conventions}\label{n-a-c}
\subsubsection*{General notation}
\begin{itemize}
\item In this article, all rings have identities. We fix a prime number $p$. 
\item Let $\p$ be either $\{p\}$ or $\emptyset$. For example, if $\p=\emptyset$, $\zbkpp=\bb{Q}$, $\zhpp=\wat{\bb{Z}}=\varprojlim_{n\in \bb{Z}_{>0}}\bb{Z}/n\bb{Z}$, and $\App=\A$, the finite adeles over $\bb{Q}$; if $\p=\{p\}$, $\zhpp:=\varprojlim_{n\in \bb{Z}_{>0} \text{ and }p\nmid n}\bb{Z}/n\bb{Z}$ and $\App=\Ap$ the finite adeles away from $p$.
\item In the notation ``$\otimes_R$'' for tensor products over a commutative ring $R$, we will omit the subscript ``$R$'' when $R=\bb{Z}$, or when $R$ is clear in the context.  
\item\textbf{(Conjugation)} Let $G$ be a group acting on a set $S$. Let $S_1$ be a subset of $S$. We denote the image in $S$ of the action of $g\in G$ on $S_1$ by $\lcj{S_1}{g}$ (resp. $S_1^g$) if the action of $G$ is a left action (resp. right action). In particular, let $K$ be a subgroup of a group $G$ and let $g\in G$ be any element in $G$, then we write $\lcj{K}{g}:=gKg^{-1}$ and write $\rcj{K}{g}:=g^{-1}Kg$.
\item Let $E$ be a number field and let $v$ be a place of $E$. We denote by $\ca{O}_{E,(v)}$ the localization of $\ca{O}_E$ at the prime corresponding to $v$ and by $\ca{O}_{E,v}$ its completion. We write the subscript of $E_v$ out to avoid multiple subscripts in the notation.
  \item Let $\ca{F}$ be an fppf sheaf of $B$-modules on a scheme $S$. Let $A$ be a $B$-algebra. Denote $\aut_{A}\ca{F}:=(\uend_S \ca{F}\otimes_B A)^\times$. Let $C$ be an $A$-algebra, and let $T$ be an object in $S_{fppf}$, we can form $\aut_A(\ca{F}(T))(C):=((\uend_S \ca{F})(T)\otimes_B C)^\times$, the group of $C$-valued automorphisms of $\ca{F}|_T$. To save the notation, we denote $\aut_A(\ca{F})(A)$ by $\aut_A\ca{F}$, which is a sheaf of $A$-modules. We will omit $A$ in the above notation if $A=\bb{Z}$.
  \end{itemize}
\subsubsection*{Notation related to Shimura varieties and compactifications}
  \begin{itemize}
  \item For a connected reductive group $G$,  we denote by $G^\der$ the derived group (resp. $G^\ad$ the adjoint group) of $G$. For an algebraic group $G$ over $\bb{R}$ or $\bb{C}$, denote by $G^\circ$ the identity component of $G$. For a Lie group $\mathscr{H}$ over $\bb{R}$ or $\bb{C}$, denote by $\mathscr{H}^+$ its identity component. So it makes sense to write $G(\bb{R})^+$ or $G(\bb{C})^+$. We will also use the symbol $G(\bb{R})_+$ to denote the pullback to $G(\bb{R})$ of the identity component $G^{\ad}(\bb{R})^+$ under the natural homomorphism $G(\bb{R})\to G^{\ad}(\bb{R})$. Let $G(\bb{Q})?:=G(\bb{R})?\cap G(\bb{Q})$, where $?={}^+,{}_+$. Let $G^\ad(\bb{R})_1$ be the subgroup of $G^\ad(\bb{R})$ that stabilizes $X\hookrightarrow X^\ad$, that is, the image of $G(\bb{R})$ in $G^\ad(\bb{R})$. Define $G^\ad(\bb{Q})_1:=G^\ad(\bb{R})_1\cap G^\ad(\bb{Q})$, etc.
  \item\textbf{(Sign conventions)} We follow Pink's sign conventions for Hodge and weight cocharacters (see \cite[Ch. 1]{Pin89}), and the class field theory isomorphism (see \cite[11.3]{Pin89}). 
  The Hodge cocharacter $\mu_\bb{C}$ is defined such that $z=(z,1)\in \bb{G}_m(\bb{C})\sbst \bb{S}(\bb{C})$ acts on $V^{-1,0}_\bb{C}$ by $z$. The split torus $\bb{G}_{m,\bb{R}}$ is embedded into $\bb{S}$ by $i_w:\bb{G}_{m,\bb{R}}\hookrightarrow \bb{S}$, $r\in \bb{R}^\times\hookrightarrow \bb{C}^\times$. Let $V_\bb{R}$ be an $\bb{R}$-vector space. Suppose there is a map $h: \bb{S}\to \mrm{GL}(V_\bb{R})$. Define the weight cocharacter $\omega$ to be $\omega: \bb{G}_{m,\bb{R}}\xrightarrow{i_w}\bb{S}\to \mrm{GL}(V_\bb{R})$. We say that $h$ has weight $-n$ if $\omega(r)$ sends $v\in V_\bb{R}$ to $r^n v$. 
Moreover, we choose the class field theory isomorphism $\mrm{rec}_E:\gal(E^{ab}/E)\iso \pi_0(\bb{G}_m(\bb{A}_E)/\bb{G}_m(E))$ such that, for any non-archimedean place $v$ of $E_v$ with a uniformizer $\pi_v$, and any finite abelian field extension $L/E$ such that $L/E$ is unramified at $v$, $(\cdots,\pi_v^{-1},\cdots)$ corresponds to (arithmetic) Frobenius of $\gal(L/E)$.
 \item\textbf{(Subgroups)} Let $G$ be a linear algebraic group over $\bb{Q}$. A subgroup $\Gamma$ in $G(\bb{Q})$ is called a congruence subgroup if it is of the form $\Gamma= G(\bb{Q})\cap K$ for some open compact subgroup $K\sbst G(\A)$. A subgroup $\Gamma'$ of $G(\bb{Q})$ is called an arithmetic subgroup if it is commensurable with (any congruence subgroup) $\Gamma$.\par
 Assume that $G$ is connected reductive. We call an open compact subgroup $K_p$ of $G(\bb{Q}_p)$ a Bruhat-Tits stabilizer subgroup if $K_p$ is the group of $\bb{Z}_p$-points of a Bruhat-Tits stabilizer group scheme $\mathscr{G}$. Denote by $K_p^\circ$ the parahoric subgroup associated with $K_p$, i.e., the $\bb{Z}_p$-points of the parahoric group scheme $\mathscr{G}^\circ$ associated with $\mathscr{G}$. An open compact subgroup $K_p''\sbst G(\bb{Q}_p)$ is called a quasi-parahoric subgroup (associated with $K_p$) if $K_p^\circ\sbst K_p'' \sbst K_p$.
\item\textbf{(Neatness)} Let $G$ be a linear algebraic group over $\bb{Q}$. Let $\square=\{p\},$ or $\emptyset$. Following Pink's convention, we say an open compact subgroup $K^\square\sbst G(\A^\square)$ is \textbf{neat} if, for some (and for any) faithful representation $i:G\hookrightarrow \mrm{GL}(V)$ and any $g=(g_\ell)_{\ell\notin \square}\in K^\p$, the intersection $$\bigcap\limits_{\ell\notin \square} (\overline{\bb{Q}}^\times\cap \Gamma_\ell)_{\mrm{torsion}}$$ is trivial. In the intersection above, $\Gamma_\ell$ is defined to be the group generated by the eigenvalues of $g_\ell$ in $\overline{\bb{Q}}_\ell^\times$.\par
  If $K^\square$ is neat, then for any $g\in G(\A)$ and any open compact subgroup $K\sbst G(\A)$ such that the projection of $K$ to $G(\A^\square)$ is $K^\square$, the intersection $g K g^{-1}\cap G(\bb{Q})$ is a neat (congruence) subgroup. See \cite[pp.12-13]{Pin89}. \emph{Although neatness will be assumed throughout the paper, we will re-emphasized it when this assumption is crucial in the argument.}
\item (\textbf{Extension property}) Let $\ca{O}$ be a Dedekind domain with fraction field $K$ and let $X$ be an $\ca{O}$-scheme. Following \cite{Kis10}, we say that $X$ has the extension property if, for any regular and formally smooth $\ca{O}$-scheme $S$, any morphism $S\otimes_{\ca{O}}K\to X$ extends to a morphism $S\to X$.
\item (\textbf{Cones}) Let $\mbf{E}\iso \bb{G}_m^n$ be a split algebraic torus over $\bb{Z}$.  The rational polyhedral cones in this paper are \emph{closed}, i.e., a rational polyhedral cone $\sigma$ in the cocharacter group $X_*(\mbf{E})_\bb{R}$ is a subset determined by $\sigma=\bb{R}_{\geq 0}v_1+ \cdots+\bb{R}_{\geq 0}v_m$ for $v_i\in X_*(\mbf{E})_\bb{Q}$ and $1\leq i\leq m\leq n$. This convention coincides with those in \cite{Pin89} and \cite{Mad19}. The $\sigma$ in \cite[Def. 6.1.1.5]{Lan13} is written as $\sigma^\circ$ in this paper. Here, $\sigma^\circ$ denotes the interior of the closed cone $\sigma$. The closure $\overline{\sigma}$ in \emph{loc. cit.} is written as $\sigma$ in this paper.\par
Let $\ca{T}$ be an $\mbf{E}$-torsor over a scheme $\ca{S}$. Denote by $\ca{T}(\sigma)$ the twisted (relatively) affine torus embedding and by $\ca{T}_{\sigma}$ the $\sigma$-stratum (see \cite[5.2]{Pin89}, \cite[Sec. 6.1.2]{Lan13} and \cite[2.1.17]{Mad19}). These two terms defined in \cite[5.2]{Pin89} and \cite[Sec. 6.1.2]{Lan13} are the same because the definitions of $\sigma^\vee$ and $\sigma^\perp$ remain the same for the two different definitions of $\sigma$. 
\item When writing an embedding of (integral models of) Shimura varieties, we often omit the reflex fields. That is, $\sh_{K_1}(G_1,X_1)\hookrightarrow \sh_{K_2}(G_2,X_2)$ means $\sh_{K_1}(G_1,X_1)\hookrightarrow \sh_{K_2}(G_2,X_2)\otimes_{E(G_2,X_2)} E(G_1,X_1)$. We will also omit the reflex fields when saying relative normalizations. That is, the relative normalization in $\sh_{K_1}(G_1,X_1)$ of $\ca{S}_{K_2}(G_2,X_2)$ in fact means the relative normalization in $\sh_{K_1}(G_1,X_1)$ of $\ca{S}_{K_2}(G_2,X_2)_{\ca{O}_{E_1(G_1,X_1),(v_1)}}$ or $\ca{S}_{K_2}(G_2,X_2)_{\ca{O}_{E_1(G_1,X_1),v_1}}$.
\end{itemize}

\newpage
%
%
%
%
\section{Characteristic zero theory}\label{sec-supplement}
This {section} has two main goals. The first one (see \S\ref{subsec-bound-chart-char-0} and \S\ref{subsec-bound-chart-conn-comp}) is to review and summarize some key definitions and theorems in Pink's thesis \cite{Pin89}. See also a summary of it in \cite[Sec. 2.1]{Mad19}. The second one (see \S\ref{subsec-zp-cusp}) is a preparation for our main construction in {Section} \ref{sec-tor-abelian}. Namely, we shall introduce our strategy of passing from Hodge-type Shimura data to abelian-type ones, and we shall explain how to make a good choice of cone decompositions. 
\subsection{Boundary charts in characteristic zero theory}\label{subsec-bound-chart-char-0}
In \S\ref{subsec-bound-chart-char-0}, we collect some definitions and results in \cite{Pin89} and \cite[Sec. 2.1]{Mad19} due to Pink and Madapusi.
\subsubsection{}\label{subsubsec-mixedshimura} 
Let $P$ be a connected linear algebraic group over $\bb{Q}$. Let $W$ be the unipotent radical of $P$ and let $U$ be a normal subgroup of $W$. Let $\ca{X}$ be a $P(\bb{R})U(\bb{C})$-homogeneous space together with a $P(\bb{R})U(\bb{C})$-equivariant morphism $\gls{hbar}:\ca{X}\to \mrm{Hom}(\bb{S}_{\bb{C}},P_{\bb{C}})$. The triple $(P,\ca{X},\hbar)$ is called a \textbf{mixed Shimura datum} if it satisfies the conditions listed in \cite[Def. 2.1]{Pin89}. We will omit the morphism $\hbar$ in the notation of a mixed Shimura datum when $\hbar$ is clear in the context. We denote $\hbar(P,\x):=(P,\hbar(\x))$.\par 
When $W=1$, the triple is called a \textbf{pure Shimura datum}. We do not omit the word ``pure'' to distinguish it from the usual definition of Shimura data as in \cite{Del71}, \cite{Del79} and so on. 
In other words, if $(P,\ca{X},\hbar)$ is a pure Shimura datum, $\hbar(P,\ca{X})$ is a Shimura datum in the sense of \emph{loc. cit}., but $(P,\ca{X},\hbar)$ itself is a Shimura datum if and only if $\hbar$ is a homeomorphism.\par
Let $P^\prime$ be a normal $\bb{Q}$-subgroup of $P$. Define $T:=P/P^\prime$ with a canonical quotient map $\phi:P\to T$. As in \cite[Prop. 2.9]{Pin89}, there is a canonical \textbf{quotient mixed Shimura datum} $(T,\ca{Y})$ of $(P,\ca{X})$ by $P^\prime$, and a canonical \textbf{quotient morphism}, which we abusively denote by $\phi:(P,\ca{X})\to (T,\ca{Y})$. 
In fact, there is an isomorphism $\ca{Y}\iso T(\bb{R})\phi(U(\bb{C}))/\phi(\mrm{Stab}_{P(\bb{R})U(\bb{C})}(x))$ for any $x\in\ca{X}$. Sometimes we will abusively write $\ca{Y}$ as $\ca{X}/P'$, although $\ca{Y}$ is not a quotient of $\ca{X}$ in general. In particular, if $T$ is a torus, $(T,\ca{Y})$ being a mixed Shimura datum implies that $\ca{Y}$ is finite and discrete. Hence, $\phi(\stb_{P(\bb{R})U(\bb{C})}(x))$ is a subgroup of $T(\bb{R})\phi(U(\bb{C}))$ containing $(T(\bb{R})\phi(U(\bb{C})))^+$.\par
A mixed Shimura datum $(P,\ca{X},\hbar)$ is \textbf{irreducible} if $P$ is the only normal $\bb{Q}$-subgroup $Q$ of $P$ such that all homomorphisms $h_x:\bb{S}_\bb{C}\to P_\bb{C}$ defined by $x\in\hbar(\ca{X})$ factor through $Q_\bb{C}$. 
\subsubsection{}\label{subsubsec-can-model}
Let $(P,\ca{X})$ be a mixed Shimura datum. For any neat open compact subgroup $K$ in $P(\A)$, define
$$\sh_K(P,\ca{X})(\bb{C}):= P(\bb{Q})\bss \ca{X}\times P(\A)/K.$$
The group $P(\bb{Q})$ acts on $\ca{X}\times P(\A)$ diagonally on the left, while the group $K$ acts only on $P(\A)$ on the right.
By \cite[Prop. 3.3 (b)]{Pin89}, the double coset $\sh_K(P,\ca{X})(\bb{C})$ is a complex manifold; by \cite[Prop. 6.26 and Prop. 9.24]{Pin89}, it admits a canonical structure of a smooth quasi-projective variety $\sh_K(P,\ca{X})_\bb{C}$ over $\bb{C}$.\par
Let $p_f\in P(\A)$ and choose $K'$ such that $p_f^{-1}K'p_f\sbst K$, we can define a morphism
$$\pi_{K',K}(p_f): \sh_{K'}(P,\ca{X})(\bb{C})\lra \sh_K(P,\ca{X})(\bb{C})$$
sending $[(x,p)]$ to $[(x,p\cdot p_f)]$. Then the map $\pi_{K',K}(p_f)$ also algebraizes to a morphism between mixed Shimura varieties. Any morphism between mixed Shimura data $f:(P_1,\ca{X}_1)\to (P_2,\ca{X}_2)$, with a compatible choice of neat open compact subgroups $K_1\sbst P_1(\A)$ and $K_2\sbst P_2(\A)$ such that $K_1\sbst K_2$, also induces a morphism between mixed Shimura varieties $f: \sh_{K_1}(P_1,\ca{X}_1)_\bb{C}\to \sh_{K_2}(P_2,\ca{X}_2)_\bb{C}$.\par
Let us briefly recall the definition of canonical models of mixed Shimura varieties (see \cite[Ch. 11]{Pin89} and also \cite[1.5.3]{KSZ21}). 
If the mixed Shimura datum is of the form $(T,\ca{Y})$ where $T$ is a torus and $\ca{Y}$ is a finite set, there is a morphism 
$$R_{E_T/\bb{Q}}\bb{G}_{m,E_T}\xrightarrow{R_{E_T/\bb{Q}}\mu_y}R_{E_T/\bb{Q}}T_{E_T}\xrightarrow{\mrm{Norm}_{E_T/\bb{Q}}}T,$$
where $E_T=E(T,\ca{Y})$ and $\mu_y$ is the Hodge cocharacter associated with the homomorphism $h_y:\bb{S}_\bb{C}\to T_\bb{C}$ determined by any $y\in \ca{Y}$. This is independent of the choice of $y$ since $\hbar(\ca{Y})$ is a single point. 
The morphism above induces a continuous homomorphism $\pi_0(\bb{G}_m(\bb{A}_{E_T})/\bb{G}_m(E_T))\to \pi_0(T(\bb{A})/T(\bb{Q}))$. Pre-composing this homomorphism with the reciprocity law $\mrm{rec}_{E_T}: \mrm{Gal}(E^{ab}_T/E_T)\xrightarrow{\sim}\pi_0(\bb{G}_m(\bb{A}_{E_T})/\bb{G}_m(E_T))$, we obtain $$r_{E_T}(T,\ca{Y}):\mrm{Gal}(\overline{E}_T/E_T)\lra \mrm{Gal}(E^{ab}_T/E_T)\lra \pi_0(T(\bb{A})/T(\bb{Q})).$$ 
Let $K_T$ be any neat open compact subgroup in $T(\A)$, since $\sh_{K_T}(T,\ca{Y})(\overline{E}_T)=T(\bb{Q})\bss \ca{Y}\times T(\A)/K_T$ is a finite set, the left action of $T(\bb{A})/T(\bb{Q})$ on it factors through $\pi_0(T(\bb{A})/T(\bb{Q}))$. Hence, there is a well-defined left action of $\mrm{Gal}(\overline{E}_T/E_T)$ on $\sh_{K_T}(T,\ca{Y})(\overline{E}_T)$ through $r_{E_T}(T,\ca{Y})$. This action forms a descent datum, which determines a canonical model $\sh_{K_T}(T,\ca{Y})$ for $\sh_{K_T}(T,\ca{Y})_\bb{C}$ over $E_T$.\par
By \cite[Thm. 11.18]{Pin89}, for every mixed Shimura datum $(P,\ca{X})$ and $K$ neat open compact subgroup of $P(\A)$, $\sh_K(P,\ca{X})_\bb{C}$ admits a \textbf{canonical model} $\sh_K(P,\ca{X})$ over its reflex field $E:=E(P,\ca{X})$, and this canonical model is uniquely characterized by the following two properties:
\begin{itemize}
    \item For any $p_f$ and $K'$ as above, the morphism $\pi_{K',K}(p_f)$ descends to $E$;
    \item For any embedding $\iota:(T,\ca{Y})\hookrightarrow (P,\ca{X})$ and any $\iota(K_T)\sbst K$, the induced morphism $\iota: \sh_{K_T}(T,\ca{Y})_\bb{C}\to \sh_K(P,\ca{X})_\bb{C}$ descends to $E_T$.
\end{itemize}
\subsubsection{}\label{subsubsec-boundary} 
Let $G$ be a connected reductive group over $\bb{Q}$ and $X$ be a $G(\bb{R})$-conjugacy class of homomorphisms of the form $h:\bb{S}\to G_\bb{R}$. Let $(G,X)$ be a Shimura datum, which means that $(G,X)$ satisfies \cite[2.1.1.1-2.1.1.3]{Del79}. We fix a triple $(H_0,h_0,h_\infty)$ as in \cite[Sec. 4.3]{Pin89}. More precisely, let $\gls{H0}$ be the reference group defined as in \emph{loc. cit.} over $\bb{R}$. In fact, $H_0$ is characterized by 
$$H_0(\bb{R}):=\{(z,\alpha)\in\bb{S}(\bb{R})\times \mrm{GL}_{2}(\bb{R})||z|^2=\det \alpha\},\text{ and}$$
$$H_0(\bb{C}):=\{((z_1,z_2);\alpha)\in \bb{S}(\bb{C})\times \mrm{GL}_2(\bb{C})|z_1z_2=\det \alpha\}.$$
Let $h_0:\bb{S}\to H_0$ and $h_\infty:\bb{S}\to H_0$ be homomorphisms from $\bb{S}$ to $H_0$ chosen as in \emph{loc. cit.}, with associated weight cocharacters $\omega_0$ and $\omega_\infty$ obtained by pre-composing with the natural embedding $\bb{G}_{m,\bb{R}}\hookrightarrow \bb{S}$.\par
Let us recall the definition of \textbf{admissible $\bb{Q}$-parabolic subgroups} (see \cite[Def. 4.5]{Pin89}). The adjoint group $G^\ad$ of $G$ can be decomposed as a product $G^\ad=\prod_{i=1}^k G_i$ of $\bb{Q}$-simple reductive groups $G_i$. An admissible $\bb{Q}$-parabolic subgroup $Q$ is the preimage in $G$ of a subgroup $Q^\ad=\prod_{i=1}^k Q_i$, where $Q_i\sbst G_i$ is either $G_i$ or a maximal proper $\bb{Q}$-parabolic subgroup of $G_i$. \par
Let $Q$ be any admissible $\bb{Q}$-parabolic subgroup of $G$. Then the homomorphisms $h_0$ and $h_\infty$ are chosen such that $\omega_0\cdot\omega_\infty^{-1}$ factors through $H_0^{\der}$. The difference $\lambda_{H_0}:=\omega_\infty\cdot \omega_0^{-1}$ is a cocharacter $\lambda_{H_0}:\bb{G}_{m,\bb{R}}\to H_0$. Let $B_0$ (resp. $U_0$) be the Borel subgroup (resp. unipotent radical of a Borel subgroup) of $H_0$ corresponding to the subspace of nonnegative (resp. positive) weights under the adjoint representation of $\lambda_{H_0}$ on $\lie H_0$.
For any admissible $\bb{Q}$-parabolic subgroup $Q$ of $G$ and any $x\in X$, there is a unique homomorphism $\gls{uQx}: H_{0,\bb{C}}\to G_{\bb{C}}$ characterized by the conditions listed in \cite[Prop. 4.6]{Pin89} such that $u_x^Q\circ h_0=h_x$, and such that $u_x^Q\circ h_{\infty}$ and $\omega_{x,\infty}^Q:=u_x^Q\circ\omega_\infty$ factor through $Q_{\bb{C}}$. Let $\omega_{x,0}^Q:=u_x^Q\circ \omega_0=\omega_x$. Then $\omega_{x,0}^Q\cdot\omega_{x,\infty}^{Q,-1}$ factors through $G_{\bb{C}}^{\der}$.\par
A \textbf{cusp label representative} $\gls{Phi}$ of $(G,X)$ is defined by a triple $\Phi=(Q_\Phi,X^+_\Phi,g_\Phi)$ as in \cite[2.1.7]{Mad19}, where $Q_\Phi$ is an admissible $\bb{Q}$-parabolic subgroup of $G$, $X^+_\Phi$ is a connected component of $X$ and $g_\Phi$ is an element in $G(\bb{A}_f)$. For any $x\in X$, we can denote by $[x]$ its image under $X\to \pi_0(X)$, which represents a connected component of $X$. 
Let $\gls{PPhi}$ be the smallest (connected) normal $\bb{Q}$-subgroup of $\gls{QPhi}$, such that $u^{Q_\Phi}_x\circ h_\infty$ factors through $P_{\Phi,\bb{C}}$. We can also denote $P_\Phi$ by $P_{Q_\Phi}$ because it only depends on $Q_\Phi$. \par
Let $W_\Phi$ be
the unipotent radical of $Q_\Phi$ and let $\gls{UPhi}$ be the center of $\gls{WPhi}$ (see \cite[III. Cor. 4.4, p.147]{AMRT10}).
Consider the $Q_\Phi(\bb{R})U_\Phi(\bb{C})$-equivariant map
\begin{equation}\label{qruc}
    \gls{tau}:X\lra \pi_0(X)\times \mrm{Hom}(\bb{S}_\bb{C}, P_{\Phi,\bb{C}})
\end{equation}
defined by mapping $x\in X$ to $([x],u_x^{Q_\Phi}\circ h_\infty)$. By \cite[Sec. 4.11]{Pin89}, the $P_\Phi(\bb{R})U_\Phi(\bb{C})$-orbit of the image of any $x\in X$ under the map above only depends on the image of $x$ in $\pi_0(X)$.\par
Let $\gls{DPhi}$ be the $P_{\Phi}(\bb{R})U_{\Phi}(\bb{C})$-orbit of the image of $x_\Phi\in X^+_\Phi$ under (\ref{qruc}). The morphism over connected components $\pi_0(\tau^{-1}(D_\Phi))\xrightarrow{\sim}\pi_0(D_\Phi)$ induced by $\tau$ is a bijection. We can also denote it by $D_{Q_\Phi,X^+_\Phi}$, because it only depends on $Q_\Phi$ and $X^+_\Phi$.\par 
By \cite[Sec. 4.11]{Pin89}, the pair $\gls{rbc}$ is a mixed Shimura datum, which is called the \textbf{rational boundary component of $(G,X)$ associated with $\Phi$}. Define $\overline{P}_{\Phi}:= P_{\Phi}/U_\Phi$ and $P_{\Phi,h}:=P_{\Phi}/W_\Phi$. Then we have mixed Shimura data defined by $(\overline{P}_{\Phi}, \overline{D}_\Phi):=(P_{\Phi},D_\Phi)/U_\Phi$ and $(P_{\Phi,h},D_{\Phi,h}):=(P_{\Phi},D_\Phi)/W_\Phi$ as in \cite[Prop. 2.9]{Pin89}. Note that $(P_{\Phi,h},D_{\Phi,h})$ is a pure Shimura datum.\par
For any neat open compact subgroup $K$ of $G(\bb{A}_f)$, define $K_\Phi:=g_\Phi K {g_\Phi}^{-1}\cap P_{\Phi}(\bb{A}_f).$ Define $\overline{K}_\Phi$ (resp. $K_{\Phi,h}$) to be the image of $K_\Phi$ in $\overline{P}_{\Phi}$ (resp. $P_{\Phi,h}$) via the obvious quotient maps.
Define double coset spaces  $\mrm{Sh}_{K_\Phi}(P_{\Phi},D_\Phi)(\bb{C})$ $:=P_{\Phi}(\bb{Q})\backslash D_\Phi\times P_{\Phi}(\bb{A}_f)/ K_\Phi$, $\mrm{Sh}_{\overline{K}_\Phi}(\overline{P}_{\Phi},\overline{D}_\Phi)(\bb{C}):=\overline{P}_{\Phi}(\bb{Q})\backslash \overline{D}_\Phi\times\overline{P}_{\Phi}(\bb{A}_f)/\overline{K}_\Phi$ and $\mrm{Sh}_{K_{\Phi,h}}({P}_{\Phi,h},{D}_{\Phi,h})(\bb{C}):={P}_{\Phi,h}(\bb{Q})\backslash D_{\Phi,h}\times{P}_{\Phi,h}(\bb{A}_f)/{K}_{\Phi,h}$, as in \cite[2.1.7]{Mad19}. They admit canonical structures of smooth quasi-projective algebraic varieties over $\bb{C}$ and they admit canonical models over the reflex field $E:=E(G,X)$ associated with $(G,X)$ (see \S\ref{subsubsec-can-model} above). Denote their canonical models by $\mrm{Sh}_{K_\Phi}(P_\Phi,D_\Phi)$, $\mrm{Sh}_{\overline{K}_\Phi}(\overline{P}_{\Phi},\overline{D}_{\Phi})$ and $\mrm{Sh}_{K_{\Phi,h}}(P_{\Phi,h},D_{\Phi,h})$, respectively. There is a tower $\sh_{K_\Phi}\to \overline{\sh}_{K_\Phi}\to\sh_{K_\Phi,h}$ over $\spec E$.\par
We call $(P_\Phi,D_\Phi)$ (resp. $\sh_{K_\Phi}(P_\Phi,D_\Phi)$) the \textbf{boundary mixed Shimura datum} (resp. \textbf{boundary mixed Shimura variety}) associated with $(G,X;\Phi)$ (resp. $(G,X;\Phi;K)$). 
\begin{rk}
Note that our system of notation is slightly different from those in \cite{Pin89} and \cite{Mad19}, and the systems of notation in the just-mentioned references are also slightly different. See a brief comparison in the List of Symbols.
\end{rk}
\begin{convention}\label{conv-mixed-sh}To avoid overloaded notation, we will frequently adopt the following abbreviations: $\sh_{K_\Phi}:=\sh_{K_\Phi}(P_\Phi,D_\Phi)$, $\overline{\sh}_{K_\Phi}:=\sh_{\overline{K}_{\Phi}}(\overline{P}_\Phi,\overline{D}_\Phi)$ and $\sh_{K_\Phi,h}:=\sh_{K_{\Phi,h}}(P_{\Phi,h},D_{\Phi,h})$; similarly, we write $\gls{shKPhiF}$, $\gls{solshKPhiF}$ and $\gls{shKPhihF}$ for the base change to a field extension $F$ of $E=E(G,X)$; we write $\sh_{K_\Phi}(F)$, $\overline{\sh}_{K_\Phi}(F)$ and $\sh_{K_\Phi,h}(F)$ for $F$-valued points. In other words, we will frequently omit the data in the brackets, which is usually harmless, since the information is hidden in the given data $\Phi$ and the sub/super-scripts of the level groups.
\end{convention}
Since $D_{Q_\Phi,X_\Phi^+}$ is a $P_\Phi(\bb{R})U_\Phi(\bb{C})$-orbit, by \cite[4.14]{Pin89}, there is a natural continuous map $\im: D_{Q_\Phi,X_\Phi^+}\to U_{\Phi}(\bb{R})(-1):=(2\pi\sqrt{-1})^{-1}U_\Phi(\bb{R})$, mapping $x\in D_{Q_\Phi,X_\Phi^+}$ to its imaginary part $u_x\in U_\Phi(\bb{R})(-1)$. As explained in \cite[4.15]{Pin89}, one associates with $(Q_\Phi,X_\Phi^+)$ an open homogeneous self-adjoint non-degenerate cone $\gls{PPhi+}$ in the sense of \cite[Ch. II]{AMRT10}.
\subsubsection{}\label{cusp-label}
Let $\gls{CLRGX}$ be the set of all cusp label representatives for $(G,X)$ defined as in \S\ref{subsubsec-boundary}. 
Note that $G(\bb{Q})$ acts on $\ca{CLR}(G,X)$ by sending $\Phi$ to $$\gamma\Phi:=(\gamma Q_\Phi\gamma^{-1},\gamma X^+_\Phi,\gamma g_\Phi)$$ for any $\gamma\in G(\bb{Q})$.\par 
More importantly, we can also view the action of $\gamma$ as the one induced by the conjugation 
$$\gamma:G\to G;\ g\mapsto \gamma g\gamma^{-1}.$$
This is a homomorphism between algebraic groups. This homomorphism induces a morphism between Shimura data and Shimura varieties
$$\gamma:(G,X)\to(G,X);\ \sh_{K}(G,X)\to \sh_{\lcj{K}{\gamma}}(G,X).$$
By \cite[4.16]{Pin89}, this induces an automorphism of $\ca{CLR}(G,X)$, sending $\Phi$ to a cusp label representative on the \emph{target} $\gamma\Phi:=(\gamma Q_\Phi \gamma^{-1},\gamma X_\Phi^+,\gamma g_\Phi \gamma^{-1})$.\par
For both viewpoints, they induce the same morphism between mixed Shimura varieties. The second viewpoint is more useful since it can be generalized to any conjugation of $\gamma\in G^\ad(\bb{Q})$. \par
Indeed, any $\gamma\in G^\ad(\bb{Q})$ induces a conjugation of $G$, and $\gamma$ acts on $X^\ad$, so it sends $X$ to $\gamma\cdot X$ viewing $X$ as a subspace of $X^\ad$. Then $\gamma$ induces a morphism between Shimura varieties, 
$$\gamma: \sh_K(G,X)\to \sh_{\lcj{K}{\gamma}}(G,\gamma\cdot X),$$
sending $[(x,g)]$ to $[\gamma\cdot x,\gamma g\gamma^{-1}]$.
Then we can define the action of $G^\ad(\bb{Q})$ on $\ca{CLR}(G,X)$, i.e., defining $\gamma\Phi:=(\gamma Q_\Phi \gamma^{-1},\gamma X_\Phi^+,\gamma g_\Phi \gamma^{-1})$.\par
Moreover, $G(\bb{A}_f)$, and in particular, $P_\Phi(\A)$ and $K$, can act on $\ca{CLR}(G,X)$ on both left and right: if we write them as left actions, then $(g_1,g_2)\cdot\Phi:=(Q_\Phi,X^+_\Phi,g_1g_\Phi g_2^{-1})$ for any $g_1,g_2\in G(\bb{A}_f)$.\par 
Let $\Phi_1=(Q_{\Phi_1},X^+_{\Phi_1},g_{\Phi_1})$ and $\Phi_2=(Q_{\Phi_2},X^+_{\Phi_2},g_{\Phi_2})$ be two cusp label representatives. 
As \cite[2.1.14]{Mad19}, we say $\Phi_1$ and $\Phi_2$ are \textbf{equivalent} if there is a $\gamma\in G(\bb{Q})$ such that $\gamma Q_{\Phi_1}\gamma^{-1}=Q_{\Phi_2}$, $P_{\Phi_1}(\bb{Q})\gamma X^+_{\Phi_1}=P_{\Phi_2}(\bb{Q})X^+_{\Phi_2}$, and there is a $q\in P_{\Phi_2}(\bb{A}_f)$ such that $\gamma g_{\Phi_1}\equiv q g_{\Phi_2}$ modulo $K$; we denote this equivalence by $$\Phi_1\xrightarrow[\sim]{(\gamma,q)_K} \Phi_2.$$ 
An equivalent interpretation of the last condition is that $\gamma g_{\Phi_1}$ and $g_{\Phi_2}$ have the same image in $P_{\Phi_2}(\bb{A}_f)\backslash G(\bb{A}_f)/K$.  
The notion of equivalence depends on the open compact group $K$. Denote by $\gls{CuspKGX}$ the set of equivalence classes of cusp label representatives, and the elements $[\Phi]$ in it are called \textbf{cusp labels}. 
We also write $\Phi_1\xrightarrow{(\gamma,q)_K}\Phi_2$ if $\gamma P_{\Phi_1}\gamma^{-1}\sbst P_{\Phi_2}$, $P_{\Phi_1}(\bb{Q})\gamma X_{\Phi_1}^+=P_{\Phi_2}(\bb{Q})X^+_{\Phi_2}$ and $\gamma g_{\Phi_1}\equiv qg_{\Phi_2}$ modulo $K$. Define a partial order $\preceq$ on $\cusp_K(G,X)$ by writing $[\Phi_1]\preceq [\Phi_2]$ if and only if $\Phi_1\xrightarrow{(\gamma,q)_K}\Phi_2$ for some $\gamma$ and $q$.
\begin{lem}\label{dep-g}
Let $\Phi_1$ and $\Phi_2$ be two cusp label representatives. Suppose that $Q_{\Phi_1}=Q_{\Phi_2}=Q$. Then there is an element $q\in Q(\bb{Q})$ such that $qX^+_{\Phi_1}=X^+_{\Phi_2}$.
\end{lem}
\begin{proof}
Since $Q(\bb{Q})G(\bb{Q})_+\sbst G(\bb{Q})$ and $G(\bb{Q})/G(\bb{Q})_+\iso G(\bb{R})/G(\bb{R})_+$ by real approximation (see \cite[2.1.2]{Del79}), it suffices to show $Q(\bb{R})G(\bb{R})_+=G(\bb{R})$. Note $Q$ is a parabolic subgroup and $G(\bb{R})_+$ contains $K_\infty:=\mrm{Stab}_{G(\bb{R})}(x)$ for some $x\in X^+$. Then the statement follows from the Iwasawa decomposition (see \cite[Thm. 3.9, p.131]{PR94}).
\end{proof}
\begin{prop}[{cf. \cite[6.3]{Pin89} and \cite[p.221]{Pin92}}]\label{prop-iphi}
Fix an admissible $\bb{Q}$-parabolic subgroup $Q$ and a connected component $X^+$ of $X$. The set $\gls{IQ}$ defined as $$\{[\Phi]\in\mrm{Cusp}_K(G,X)| \Phi=(Q,X^+_\Phi,g_\Phi),\ X^+_\Phi\sbst X\text{ connected component, }g_\Phi\in G(\A)\}$$ consisting of cusp labels in $\mrm{Cusp}_K(G,X)$ that have representatives of the form $\Phi=(Q,X^+_\Phi,g_\Phi)$ is in bijection with
$$\mrm{Stab}_{Q(\bb{Q})}(D_{Q,X^+})P_Q(\bb{A}_f)\backslash G(\bb{A}_f)/K.$$
\end{prop}
\begin{proof}
By the definition of cusp labels, there is a well-defined map from \begin{equation}\mrm{Stab}_{Q(\bb{Q})}(D_{Q,X^+})P_Q(\bb{A}_f)\backslash G(\bb{A}_f)/K\end{equation}
to $\mrm{Cusp}_K(G,X)$ mapping $[g]$ to $[(Q,X^+,g)]$. In fact, we only have to show that $(Q,X^+,g)\sim (Q,X^+,q\cdot p\cdot g\cdot k)$ for any $q\in\mrm{Stab}_{Q(\bb{Q})}(D_{Q,X^+})$, $p\in P_Q(\A)$ and $k\in K$. It suffices to check that $(Q,X^+,q\cdot p\cdot g\cdot k)\bd{q^{-1}}{p}{K}{\sim}(Q,X^+,g)$. Since $q^{-1}(Q,X^+,q\cdot p\cdot g\cdot k)=(Q,q^{-1}X^+,p\cdot g\cdot k)$ and $p\cdot g\cdot k\equiv p\cdot g$ mod $K$, it suffices to show that $P_Q(\bb{Q})q^{-1}X^+=P_Q(\bb{Q})X^+$. Since $\tau$ induces a bijection over connected components, it suffices to show that $P_Q(\bb{Q})q^{-1}\tau(X^+)=P_Q(\bb{Q})\tau(X^+)=D_{Q,X^+}$, and this follows from the condition that $q\in\stb_{Q(\bb{Q})}(D_{Q,X^+})$.\par 
Then we check that this map is injective. If $[g]$ and $[g^\prime]$ map to the same element in $\mrm{Cusp}_K(G,X)$, then $\gamma\cdot g=p\cdot g^\prime \cdot k$ for some $\gamma\in Q(\bb{Q})=\mrm{Stab}_{G(\bb{Q})}(Q(\bb{Q}))$, $p\in P_Q(\A)$ and $k\in K$. Moreover, $\gamma$ stabilizes $X^+$, so $\gamma\in \stb_{Q(\bb{Q})}(D_{Q,X^+})$.\par 
Finally, we check that this map is surjective. For any cusp label $[\Phi]=[(Q,X^+_\Phi,g_\Phi)]$, there is a $q\in Q(\bb{Q})$ such that $qX^+_\Phi=X^+$ by Lemma \ref{dep-g}, so such an equivalence class is in the image of the map.
\end{proof}
\subsection{Boundary mixed Shimura varieties of Hodge type}\label{subsec-bound-mix-sh-hodge}
\subsubsection{}In this subsection, we show that $(P_{\Phi,h},D_{\Phi,h})$ is a Shimura datum in the usual sense when $(G,X)$ is a Shimura datum of Hodge type or in slightly more general cases. Our proof does not rely on any detailed computation of the groups such as $P_\Phi$ and $Q_\Phi$.\par 
Denote by $N:\bb{S}\to \bb{G}_{m,\bb{R}}$ the morphism induced by taking the norm over $\bb{R}$-points. Then $N(\bb{S}(\bb{R}))=\bb{R}^\times_+=\bb{G}_{m}(\bb{R})^+.$ As in \cite[Example 2.8]{Pin89}, we can define the $0$-dimensional Siegel Shimura datum to be the pair $(\bb{G}_m,\bb{H}^\pm_0)$, where $\bb{H}^\pm_0$ is the $2$-point set defined by the set of isomorphisms between $\bb{Z}$ and $\bb{Z}(1)$. Let $(G_0,X_0)$ be a Shimura datum of Hodge type with an embedding into a Siegel Shimura datum $\iota:(G_0,X_0)\hookrightarrow (G^\ddag,X^\ddag)$.\par
There is a canonical \emph{similitude character} $\nu:=$\gls{nuGddag}$:(G^\ddag,X^\ddag):=(\mrm{GSp}(V,\psi),\bb{H}_g^\pm)\to (\bb{G}_m,\bb{H}_0^\pm)$ from any Siegel Shimura datum associated with a symplectic space $(V,\psi)$ of dimension $2g>0$ to $(\bb{G}_m,\bb{H}_0^\pm)$. In fact, for any $x\in \bb{H}_g^\pm$, the point $\nu(x)$ is the isomorphism $\lambda$ such that $\lambda\circ \psi$ is a polarization with respect to the Hodge structure associated with $h_x$; the map, abusively denoted by $\nu:\mrm{GSp}(V,\psi)\to \bb{G}_m$, is defined to be the quotient map of $\mrm{GSp}(V,\psi)$ by its derived group. 
From the explanation in \S\ref{subsubsec-mixedshimura}, the quotient $\nu:\mrm{GSp}(V,\psi)\to\bb{G}_m$ induces a quotient of $(\mrm{GSp}(V,\psi),\bb{H}_{g}^\pm)$ by $\mrm{Sp}(V,\psi)$. \par
\begin{convention}\label{conv-pm-pol}
We fix a choice of $\sqrt{-1}$, and therefore fix a choice of $\bb{Z}(1)\iso \bb{Z}$. We denote by $X^{\ddag,+}$ the connected component of $X^\ddag$ such that $\psi$ is a (\emph{positive}) polarization with respect to $h_x$ for one of (and therefore all of) the points $x\in X^{\ddag,+}$.\par
Fix any cusp label representative $\Phi^\ddag$. Recall that there is a $Q_{\Phi^\ddag}(\bb{R})$-equivariant map $\tau:X^\ddag\to \pi_0(X^\ddag)\times \Hom(\bb{S}_\bb{C},P_{\Phi^\ddag,\bb{C}})$, where $D_{\Phi^\ddag}$ is a $P_{\Phi^\ddag}(\bb{R})U_{\Phi^\ddag}(\bb{C})$-orbit of $X^{\ddag,+}$. Since $P_{\Phi^\ddag}(\bb{R})$ acts transitively on $\pi_0(X^\ddag)\iso \{\pm 1\}$, $D_{\Phi^\ddag}$ does not depend on the choice of $X^{\ddag,+}$, and $\pi_0(D_{\Phi^\ddag})=\pi_0(X^\ddag)$. Unless otherwise noted, we always let $D_{\Phi^\ddag}^+:=\tau (X^{\ddag,+})$. We see in \S\ref{subsec-mxhg-cont} that points on $D^+_{\Phi^\ddag}$ correspond to mixed Hodge structures of $V$ such that $\psi$ is a (\emph{positive}) polarization.
\end{convention}
Let $\nu_G:(G_0,X_0)\to (T^\prime,Y)$ be the morphism between pure Shimura data induced by the quotient of $(G_0,X_0)$ by $G^\der_0$, where $T^\prime:=G_0/G^\der_0$. Note that there is a commutative diagram
\begin{equation}
    \begin{tikzcd}
    (G_0,X_0)\arrow[rr]\arrow[d,"\nu_G"]&& (G^\ddag,X^\ddag)\arrow[d,"\nu"]\\
    (T^\prime,Y)\arrow[rr,"\nu_G^\prime"]&& (\bb{G}_m,\bb{H}_0^\pm),
    \end{tikzcd}
\end{equation}
where $\nu_G^\prime$ is defined by the quotient of $(T^\prime,Y)$ by $\nu_G(G_0\cap G^{\ddag,\der})$.
\begin{lem}\label{injective}The map $\nu^\prime_G: Y\to \bb{H}_0^\pm$ is injective. The stabilizer $T^\prime(\bb{R})_+$ of any point of $Y$ in $T^\prime(\bb{R})$ is connected.
\end{lem}
\begin{proof}
Let $\omega:\bb{G}_m\to G_0$ be the weight cocharacter. Then $\nu_G^\prime\circ\nu_G\circ\omega$ is isomorphic to $\bb{G}_m\xrightarrow{x\mapsto x^2}\bb{G}_m$.
Denote by $C\iso \mu_2$ the kernel of this composition. \par
Let $G_0^\prime:=G^{\ddag,\der}\cap G_0$. There is a commutative diagram of algebraic groups
\begin{equation}
    \begin{tikzcd}
        &\bb{G}_m\arrow[drr,two heads,"\nu\circ\omega"]\arrow[d,hook,"\omega"]&&\\       G_0^\prime\arrow[r,hook]&G_0\arrow[rr,two heads,"\nu_G^\prime\circ\nu_G"]&& \bb{G}_m,
    \end{tikzcd}
\end{equation}
where the horizontal row is exact.
Hence, there is an almost-direct product decomposition $G_0\iso G_0^\prime\times \omega(\bb{G}_m)/C$ for $G_0$. In what follows, we split the proof into two cases, according to whether or not the two-point finite group $C$ is contained in $G_0^{\der}$.\par
If $C$ lies in $G^\der_0$, then $T^\prime\iso \bb{K}\times \bb{G}_m$, where $\bb{K}$ is the kernel of $\nu_G^\prime$. In fact, $T^\prime\iso \bb{K}\times (\nu_G\circ\omega(\bb{G}_m))/C^\prime$, where $C^\prime$ is the image of $C$ under $\nu_G$, which is trivial in this case. So $\nu_G^\prime$ is the projection of $T^\prime\iso \bb{K}\times \bb{G}_m$ to the second factor. Note that $\bb{K}$ is a \emph{connected and $\bb{R}$-anisotropic} torus since $T^\prime$ is connected and $h_x(i)$ induces a Cartan involution of $(G^\ddag/\bb{G}_m)_{\bb{R}}$, so $\bb{K}(\bb{R})$ is also a \emph{connected compact} Lie group by a theorem of Chevalley
(see, e.g., \cite[Sec. 3.1, Thm. 3.1, and Sec. 3.2, Cor. 1]{PR94}). Hence, $\nu^\prime_G: Y\to \bb{H}_0^\pm$ is isomorphic to the morphism $\pi_0(T^\prime(\bb{R}))\to\pi_0(\bb{G}_m(\bb{R}))$ induced by the natural projection $T^\prime_\bb{R}\to \bb{G}_{m,\bb{R}}$, and therefore is injective since $\bb{K}(\bb{R})$ is connected. The second statement follows from the injectivity proved above, and the connectedness of $\bb{K}(\bb{R})$.\par
If $C$ does not lie in $G^\der_0$, it follows that $C^\prime\iso \mu_2$. Let $U^1$ be the kernel of $N$. Choose any $x\in X_0$. Then $U^\prime$, the image of $U^1$ under $U^1\hookrightarrow \bb{S}\xrightarrow{\nu_G\circ h_x} T^\prime_{\bb{R}}$, is not trivial since $C^\prime$ lies in it; so $U^\prime$ is isomorphic to $U^1$. Since $T^\prime$ is connected, it follows that $\bb{K}$ has at most $2$ connected components. But we claim $\bb{K}$ is connected. In fact, $\bb{K}/C^\prime$ is connected and $C^\prime$ is contained in a compact connected subgroup $U^1$ in the compact algebraic group $\bb{K}$, so $C^\prime$ lies in the identity component of $\bb{K}$, and then $\bb{K}$ is connected. Hence, both $\bb{K}/C^\prime(\bb{R})$ and $\bb{K}(\bb{R})$ are connected by Chevalley's theorem (see, e.g., \cite[Sec. 3.1, Thm. 3.1, and Sec. 3.2, Cor. 1]{PR94}). \par
Furthermore, $T^\prime(\bb{R})$ is connected. To see this, we note that there is an extension of algebraic groups 
$$0\to \bb{G}_m\to T^\prime\to \bb{K}/C^\prime\to 0.$$
Since $\bb{K}/C^\prime(\bb{R})$ is connected and $T^\prime\to \bb{K}/C^\prime$ is surjective, the induced homomorphism between $\bb{R}$-points $T^\prime(\bb{R})\to (\bb{K}/C^\prime)(\bb{R})$ is surjective. Then we find a fibration $\bb{G}_m(\bb{R})\to T^\prime(\bb{R})\to (\bb{K}/C^\prime)(\bb{R})$.  
We can consider the long exact sequence of $\pi_0$- and $\pi_1$- groups
$$\cdots\to \pi_1((\bb{K}/C^\prime)(\bb{R}))\xrightarrow{\ \partial\ }\pi_0(\bb{G}_m(\bb{R}))\to \pi_0(T^\prime(\bb{R}))\to \pi_0((\bb{K}/C^\prime)(\bb{R}))\to 0.$$
Since $U^\prime(\bb{R})$ is a loop in $\bb{K}(\bb{R})\sbst T^\prime(\bb{R})$ containing $C^\prime(\bb{R})\iso\mu_2(\bb{R})=\{\pm1\}\sbst \nu_G\circ\omega(\bb{G}_m(\bb{R}))\iso\bb{G}_m(\bb{R})$, where the last isomorphism is the natural isomorphism induced by the quotient of $\omega(\bb{G}_m)$ by $G_0^\der\cap\omega(\bb{G}_m)=1$. Then there is a half-loop $U^+\iso[-1,1]$ in $U^\prime(\bb{R})$ connecting $1$ and $-1$ in $C^\prime(\bb{R})$ whose image in $(\bb{K}/C^\prime)(\bb{R})$ is a loop, denoted by $\overline{U}^+$. Hence, the monodromy action of the homotopy class of $\overline{U}^+$ interchanges the two elements of $\pi_0(\bb{G}_m(\bb{R}))$; that is, $\partial$ is surjective. Then we find that $T^\prime(\bb{R})$ is connected in the second case.\par
As a result, $Y$ has only one point in this case, and therefore $\nu_G^\prime:Y\to \bb{H}_0^\pm$ is injective. The second statement follows from the connectedness of $T^\prime(\bb{R})$.
\end{proof}
\begin{lem}\label{conn}
Let $(G,X)$ be a Shimura datum of abelian type. Let $G(\bb{R})_+$ be the inverse image of $G^{\ad}(\bb{R})^+$ under $G(\bb{R})\to G^\ad(\bb{R})$. Then $G(\bb{R})_+$ is connected if $(G,X)$ satisfies the following condition:
\begin{itemize}
    \item There is a Hodge-type Shimura datum $(G_0,X_0)$ and a central isogeny  $\pi: G_0\to G$ such that the kernel of $\pi$ is in $G^{\der}_0$ and $\pi$ induces a morphism between Shimura data $(G_0,X_0)\to(G,X)$.
\end{itemize}
\end{lem}
\begin{proof}
Suppose that $(G,X),$ $(G_0,X_0)$ and $\pi$ are as in the statement. Then they fit into a commutative diagram:
\begin{equation}
    \begin{tikzcd}
    (G_0,X_0)\arrow[rrrr]\arrow[d,"\pi"]&&&&(G^\ddag,X^\ddag)\arrow[d,"\nu"]\\
    (G,X)\arrow[rr]&&(T^\prime,Y)\arrow[rr,"\nu^\prime_G"]&& (\bb{G}_m,\bb{H}_0^\pm).
    \end{tikzcd}
\end{equation}
By Lemma \ref{injective}, it suffices to show that $G^\der(\bb{R})_+:=G(\bb{R})_+\cap G^\der(\bb{R})$ is connected.
Fix any $x\in X$. Let $X^+$ be the connected component of $X$ containing $x$. Since $G^\der(\bb{R})^+$ maps surjectively to $G^\ad(\bb{R})^+$ and is contained in $G(\bb{R})^+\sbst G(\bb{R})_+$, $G^\der(\bb{R})_+$ acts on $X^+$ transitively. \par 
Therefore, it suffices to show that the centralizer $\mrm{Cent}_{G^\der(\bb{R})}(h_x)$ is connected. This fact can be extracted from the proof of \cite[Prop. 2.11]{Pin89}. Since $G^\der_\bb{R}$ is connected, the algebraic group $\mrm{Cent}_{G^\der_\bb{R}}(h_x)$ is connected. The Lie algebra of $\mrm{Cent}_{G^\der_\bb{R}}(h_x)$ is fixed by the Cartan involution induced by $h_x(i)$, which implies that $\mrm{Cent}_{G^\der_\bb{R}}(h_x)$ is compact. Then the desired assertion follows by Chevalley's theorem (see, e.g., \cite[Sec. 3.1, Thm. 3.1, and Sec. 3.2, Cor. 1]{PR94}).
\end{proof}
\begin{prop}\label{pure shimura}
Let $(G,\ca{X},\hbar)$ be a pure Shimura datum with $\hbar:\ca{X}\to X=\mrm{Hom}(\bb{S},G_\bb{R})$. Suppose that $(G,X)$ satisfies the condition in Lemma \ref{conn}. Then for any neat open compact subgroup $K$, the pure Shimura variety $\mrm{Sh}_K(G,\ca{X})$ is a Shimura variety.
\end{prop}
\begin{proof}
By definition, it suffices to show that the $G(\bb{R})$-equivariant projection $\ca{X}\to \mrm{Hom}(\bb{S},G_\bb{R})$ is injective. Let $\ca{X}^+$ be a connected component of $\ca{X}$, which maps to a connected component $X^+$ of $X$. By \cite[Prop. 2.11 and Cor. 2.12]{Pin89}, $\hbar$ is a finite covering and $\ca{X}$ is a union of hermitian symmetric domains with $G(\bb{R})^+$ acting on each connected component transitively. Then the stabilizer of $\ca{X}^+$ in $G(\bb{R})$ is a subgroup of $G(\bb{R})_+$ containing $G(\bb{R})^+$. By Lemma \ref{conn}, $G(\bb{R})^+=G(\bb{R})_+$, which implies the desired statement.
\end{proof}
\begin{cor}\label{cor-pure-shimura}
Let $(G,X)$ be a Shimura datum satisfying the condition in Lemma \ref{conn}. Let $\Phi$ be any cusp label representative. Then $(P_{\Phi,h},D_{\Phi,h})$ defined in \S\ref{subsubsec-boundary} is a Shimura datum in the usual sense. If $(G,X)=(G_0,X_0)$ is of Hodge type, $(P_{\Phi,h},D_{\Phi,h})$ is of Hodge type and, in particular, $\sh_{K_\Phi,h}$ is a Hodge-type Shimura variety for any neat open compact subgroup $K\sbst G(\A)$.
\end{cor}
\begin{proof}
We first assume that $(G,X)$ is a Hodge-type Shimura datum. Let $\iota: (G,X)\hookrightarrow (G^\ddag,X^\ddag)$ be the embedding into the Siegel Shimura datum. Let $\Phi=(Q_\Phi,[x_\Phi],g_\Phi)$ be a cusp label representative of $(G,X)$. 
There is a unique minimal admissible $\bb{Q}$-parabolic subgroup $Q^\ddag_\Phi$ of $G^\ddag$ containing $Q_\Phi$. Let $[\iota(x_\Phi)]$ be the unique connected component of $X^\ddag$ containing $\iota(x_\Phi)$. Then $(Q^\ddag_\Phi,[\iota(x_\Phi)],\iota(g_\Phi))$ is a cusp label representative of $(G^\ddag,X^\ddag)$, denoted by $\Phi^\ddag$. Levi decomposition induces an injective homomorphism $P_{\Phi,h}\hookrightarrow P_{\Phi^\ddag,h}$. By \cite[4.25]{Pin89}, the group $P_{\Phi^\ddag,h}$ is a symplectic similitude group and $\hbar(D_{\Phi^\ddag,h})$ is the Siegel upper and lower half-spaces associated with $P_{\Phi^\ddag,h}$. By Proposition \ref{pure shimura}, $(P_{\Phi,h},D_{\Phi,h})=\hbar(P_{\Phi,h},D_{\Phi,h})\hookrightarrow \hbar(P_{\Phi^\ddag,h},D_{\Phi^\ddag,h})=(P_{\Phi^\ddag,h},D_{\Phi^\ddag,h})$ is a Hodge embedding since $\Hom(\bb{S}, P_{\Phi,h,\bb{R}})\hookrightarrow\Hom(\bb{S}, P_{\Phi^\ddag,h,\bb{R}})$.\par
Now we only assume that $(G,X)$ satisfies the condition in Lemma \ref{conn}. Since $\pi:G_0\to G$ is a central isogeny, $Q$ is an admissible $\bb{Q}$-parabolic subgroup of $G$ if and only if $\pi^{-1}(Q)$ is an admissible $\bb{Q}$-parabolic subgroup of $G_0$. Since this question is not related to $g_\Phi$, we can choose a cusp label representative $\Phi=(Q_\Phi,[x_\Phi],g_\Phi)$ of $(G,X)$ such that there is a connected component $[x_{\Phi_0}]$ of $X_0$ mapping to $[x_\Phi]$ and a lifting $g_{\Phi_0}$ of $g_\Phi$ to $G_0(\bb{A}_f)$. Then $\Phi_0=(\pi^{-1}(Q),[x_{\Phi_0}],g_{\Phi_0})$ is a cusp label representative of $(G_0,X_0)$. By the last paragraph and Proposition \ref{pure shimura}, $({P}_{\Phi_0,h},D_{\Phi_0,h})$ is a Hodge-type Shimura datum, and therefore $(P_{\Phi,h},D_{\Phi,h})$ satisfies the condition in Lemma \ref{conn}. Hence, the statement follows from Proposition \ref{pure shimura}. 
\end{proof}
\begin{cor}\label{cor-pure-embedding}
With the conventions above, if $(G,X)=(G_0,X_0)$ is of Hodge type and if $\Phi$ maps to $\Phi^\ddag\in\ca{CLR}(G^\ddag,X^\ddag)$ (see the conventions in Corollary \ref{cor-pure-shimura}), the morphisms $(P_\Phi,D_\Phi)\to (P_{\Phi^\ddag},D_{\Phi^\ddag})$, $(\overline{P}_\Phi,\overline{D}_\Phi)\to (\overline{P}_{\Phi^\ddag},\overline{D}_{\Phi^\ddag})$ and $(P_{\Phi,h},D_{\Phi,h})\to (P_{\Phi^\ddag,h},D_{\Phi^\ddag,h})$ induced by the Hodge embedding $(G,X)\hookrightarrow (G^\ddag,X^\ddag)$ are embeddings. The last morphism $(P_{\Phi,h},D_{\Phi,h})\hookrightarrow(P_{\Phi^\ddag,h},D_{\Phi^\ddag,h})$ is also a Hodge embedding.
\end{cor}
\begin{proof}
From the proof above, we have that $(P_{\Phi,h},D_{\Phi,h})=\hbar(P_{\Phi,h},D_{\Phi,h})\hookrightarrow \hbar(P_{\Phi^\ddag,h},D_{\Phi^\ddag,h})=(P_{\Phi^\ddag,h},D_{\Phi^\ddag,h})$ is a Hodge embedding. By \cite[Cor. 2.12]{Pin89}, every connected component of $D_\Phi$ (resp. $D_{\Phi^\ddag}$) maps isomorphically to its image in $\hbar(D_\Phi)$ (resp. $\hbar(D_{\Phi^\ddag})$). Since the fibers of projection $D_\Phi\to D_{\Phi,h}$ (resp. $D_{\Phi^\ddag}\to D_{\Phi^\ddag,h}$) are connected, this implies that $D_\Phi=\hbar(D_\Phi)$ (resp. $D_{\Phi^\ddag}=\hbar(D_{\Phi^\ddag})$). The same argument also applies to $\overline{D}_\Phi$ and $\overline{D}_{\Phi^\ddag}$.
\end{proof}
\subsubsection{}\label{subsubsec-Harris-group}
We compare the group $P_\Phi$ with the group $P''$ in \cite[2.7]{Har89}; from now on, denote the latter group by $P_\Phi''$.
More precisely, for any $\Phi$ in $\ca{CLR}(G,X)$, let $P''_{\Phi}$ be the maximal $\bb{Q}$-subgroup of $Q_\Phi$ such that the homomorphism $P''_\Phi\to \mrm{GL}(\lie U_\Phi)$ defined by the adjoint representation of $Q_\Phi$ on $\mrm{GL}(\lie U_\Phi)$ factors through $P''_\Phi\to \bb{G}_m\hookrightarrow \mrm{GL}(\lie U_\Phi)$, where the second homomorphism $\bb{G}_m\hookrightarrow \mrm{GL}(\lie U_\Phi)$ is given by the $\bb{G}_m$-action on $\lie U_\Phi$ by homotheties. The group $P''_\Phi$ is a normal subgroup of $Q_\Phi$ containing the center $Z_G$ of $G$. By \cite[Prop. 2.14(a)]{Pin89}, $P_\Phi\sbst P''_\Phi$.
Following \cite{Har89}, let $P_\Phi'$ be the kernel of the adjoint representation $Q_\Phi\to \mrm{GL}(\lie U_\Phi)$. Let $P_\Phi^*$ be the smallest normal $\bb{Q}$-subgroup of $Q_\Phi$ such that $u_x^{Q_\Phi}\circ h_\infty(i)$ factors through it and the unipotent radical of $P_\Phi^*$ is $W_\Phi$.
\begin{lem}\label{lem-Harris-group-compare}
Assume that $(G,X)$ is a Shimura datum such that there is a Hodge-type Shimura datum $(G_0,X_0)$ and a morphism $\pi:(G_0,X_0)\to (G,X)$ such that the kernel of $\pi:G_0\to G$ is finite and is in $G_0^\der$. Then the quotient group $P''_\Phi/Z_G\cdot P_\Phi$ is of compact type. 
\end{lem}
\begin{proof}In fact, we will imitate some arguments in \cite[Lem. 4.9]{Pin89}.
Let $x$ be any point on $X$.
As the boundary weight cocharacter $\omega^{Q_\Phi}_{x,\infty}$ is defined over $\bb{Q}$ by assumption, we have $P''_\Phi=\omega^{Q_\Phi}_{x,\infty}(\bb{G}_m)\cdot P'_\Phi$ and $P_\Phi\supset\omega^{Q_\Phi}_{x,\infty}(\bb{G}_m)\cdot P_\Phi^*$. Note that, from the construction of $h_\infty$ and $u_x^{Q_\Phi}$ (see \cite[Prop. 4.6]{Pin89}), $P_\Phi^*\sbst P'_\Phi$ since the adjoint representation of $u_x^{Q_\Phi}\circ h_\infty(i)$ on $\lie U_\Phi$ is trivial. This implies that $\omega^{Q_\Phi}_{x,\infty}(\bb{G}_m)\cdot P_\Phi^*$ is an almost direct product.\par
Therefore, we are reduced to showing that $P'_\Phi/Z_G\cdot P^*_\Phi$ is of compact type, and it suffices to show that $h_{x,\infty}(i):=u_x^{Q_\Phi}\circ h_\infty(i)$ induces a Cartan involution on $(P'_\Phi/Z_G\cdot W_\Phi)_\bb{R}$ since $h_{x,\infty}(i)$ factors through $P_{\Phi,\bb{C}}^*$. 
Let $\wdtd{G}$ be the centralizer of $\omega^{Q_\Phi}_{x,\infty}$ in $Q_\Phi$, which is a connected reductive subgroup that is naturally isomorphic to the Levi quotient $Q_\Phi/W_\Phi$. The quotient $P'_\Phi/Z_G\cdot W_\Phi$ is a normal subgroup of $Q_\Phi/Z_G\cdot W_\Phi$. 
By the definition of $\wdtd{G}$, $h_{x,\infty}(i)$ factors through $\wdtd{G}_\bb{C}$.\par 
Finally, we claim that $h_x(i)$ and $h_{x,\infty}(i)$ induce the same Cartan involution on $(\wdtd{G}/Z_G)_\bb{R}$: To see this, note that the difference $h_x(i)\cdot h_{x,\infty}(i)^{-1}$ lies in the image of $u_x^{Q_\Phi}(H_{0,\bb{C}}^\der)$. 
Moreover, note that the adjoint representation of $u_x^{Q_\Phi}(H_{0,\bb{C}}^\der)$ on $G$ induces a trivial restriction on $\wdtd{G}$, since the kernel of this adjoint representation is a normal subgroup of a semisimple group, since the conjugation of $\omega_{x,\infty}^{Q_\Phi}\cdot \omega_{x,0}^{Q_\Phi,-1}$ on $\wdtd{G}$ is trivial, and since the image of $u_x^{Q_\Phi}(U_0)$ is in $W_{\Phi,\bb{R}}$. Hence, $h_{x,\infty}(i)$ induces a Cartan involution on $(Q_\Phi/Z_G\cdot W_\Phi)_\bb{R}$ and on $(P'_\Phi/Z_G\cdot W_\Phi)_\bb{R}$.
\end{proof}
\subsection{Compactifications and connected components}\label{subsec-bound-chart-conn-comp}
Let $K\sbst G(\bb{A}_f)$ be a neat open compact subgroup. Let $X^+$ be a connected component of $X$. Define $G(\bb{Q})_+:=G(\bb{Q})\cap G(\bb{R})_+$. We have $$\sh_K(G,X)(\bb{C})=G(\bb{Q})_+\backslash X^+\times G(\bb{A}_f)/K\iso\disju_{g_i\in I} X^+/ \Gamma(g_i).$$ 
Here,  $\Gamma(g_i):=G(\bb{Q})_+\cap g_iKg_i^{-1}$ and $I$ is a set of representatives of the double coset $\pi_0(\mrm{Sh}_K(G,X)_\bb{C})\iso G(\bb{Q})_+\backslash G(\bb{A}_f)/K$. This disjoint union admits an algebraic structure $\mrm{Sh}_K(G,X)_{\bb{C}}\iso\disju_{g_i\in I}X^+/\Gamma(g_i)^{\mrm{alg}}$ over $\bb{C}$.\par 
In this {subsection}, we recall some key definitions and statements in Pink's thesis on minimal and toroidal compactifications of $\sh_K:=\sh_K(G,X)$. Along the way, we also study the geometrically connected components of the strata of the toroidal and minimal compactifications of a single connected component $X^+/\Gamma(g_i)^{\mrm{alg}}$ in the strata of the corresponding compactifications of $\mrm{Sh}_K(G,X)_{\bb{C}}$.
\subsubsection{}\label{dbcoset} We summarize a situation that will repeatedly occur.\par
Let $(X,X_1)$ be a pair of topological spaces with an open and closed embedding $h:X\hookrightarrow X_1$. Let $(G,G_1)$ be a pair of groups equipped with an injective homomorphism $i:G\hookrightarrow G_1$. Let $K$ (resp. $K_1$) be a subgroup of $G$ (resp. $G_1$). Let $H$ (resp. $H_1$) be a subgroup of $G$ (resp. $G_1$). Suppose that $H$ (resp. $H_1$) acts on $X$ (resp. $X_1$) by $c: H\times X\to X$ (resp. $c_1: H_1\times X_1\to X_1$), and $c$, $c_1$, $h$ and $i$ satisfy the following commutative diagram 
\begin{equation}
\begin{tikzcd}
H\times X\arrow{d}{(i,h)}\arrow[rr,"c"]&& X\arrow[d,"h"]\\
H_1\times X_1\arrow[rr,"c_1"]&& X_1,
\end{tikzcd}
\end{equation}
In addition, we assume that $H_1 X=X_1$, and that $H_1$ stabilizes $G$ by conjugation.\par 
We have left group actions of $H_1$ on $X_1$ and $G_1$, and a right group action of $K_1$ on $G_1$. 
Then we have the following decomposition of double cosets
$$H_1\bss X_1\times G_1/K_1= H_1\bss \disju_{h\in I_1} h\cdot X\times G_1/K_1\iso \stb_{H_1}(X)\bss X\times G_1/K_1,$$
where $I_1$ is a set of representatives of the coset $H_1/ \stb_{H_1}(X)$.
Moreover, $$\stb_{H_1}(X)\bss X\times G_1/K_1=\stb_{H_1}(X)\bss X\times \disju_{g\in I_2}\stb_{H_1}(X)GgK_1/K_1,$$
where $I_2$ is a set of representatives of the double coset $\stb_{H_1}(X)G\bss G_1/K_1.$ If $i(K)\sbst gK_1g^{-1}$, we also note that there is a  canonical map from $H\bss X\times G/K$ to $\stb_{H_1}(X)\bss X\times\stb_{H_1}(X)GgK_1/K_1.$ \par
Since $H_1$ stabilizes $G$ by assumption, we can further write $$\stb_{H_1}(X)\bss X\times\stb_{H_1}(X)GgK_1/K_1\iso \stb_{H_1}(X)\cap GgK_1g^{-1}\bss X\times GgK_1/K_1.$$
\subsubsection{}
First, we consider the minimal compactification.
For any cusp label representative $(Q_\Phi,X^+_\Phi,g_\Phi)$, from \S\ref{subsubsec-boundary}, we associate a $P_\Phi(\bb{R})U_\Phi(\bb{C})$-homogeneous space $D_{\Phi}$ with $\Phi$, and $D_{\Phi,h}$ is called a \emph{rational boundary component} of $X$ (see \cite[3.5]{BB66}). Let $X^*=\bigcup D_{\Phi,h}$ be the union of $D_{\Phi,h}$ running over all rational boundary components of $X$ (also including $X$ itself), endowed with the so-called \emph{Satake topology} (see \cite[4.8]{BB66} and \cite[6.2]{Pin89}). Note that the assignment $\Phi\mapsto D_{\Phi,h}$ is not injective: for any two cusp label representatives $\Phi=(Q_\Phi,X^+_\Phi,g_\Phi)$ and $\Phi^\prime=(Q_{\Phi^\prime},X^+_{\Phi^\prime},g_{\Phi^\prime})$, the associated rational boundary components $D_{\Phi,h}=D_{\Phi^\prime,h}$ if and only if
$D_{\Phi}=D_{\Phi^\prime}$ if and only if 
$Q_{\Phi}=Q_{\Phi^\prime}$ and $P_{\Phi}(\bb{Q})X^+_\Phi=P_{\Phi^\prime}(\bb{Q})X^+_{\Phi^\prime}$.\par
The \emph{minimal compactification} of $\mrm{Sh}_K(G,X)(\bb{C})$ is defined as $\mrm{Sh}^{\mrm{min}}_K(G,X)(\bb{C}):=G(\bb{Q})\backslash X^*\times G(\bb{A}_f)/K$, and it possesses a canonical structure of a normal variety  $\mrm{Sh}^{\mrm{min}}_K(G,X)$ over the reflex field $E:=E(G,X)$. See \cite[12.3]{Pin89}.\par 
Consider the $G(\bb{Q})$-orbit $\Psi=G(\bb{Q})\Phi$ of $\Phi$ and let $D_\Psi$ be the union of $D_{g\Phi}$ for $g\in G(\bb{Q})$. Similarly, let $D_{\Psi,h}$ be the union of subspaces $D_{g\Phi,h}$ of $X^*$ for $g\in G(\bb{Q})$.
Denote by $\mrm{Z}_{\Psi,K}(\bb{C}):=G(\bb{Q})\backslash D_{\Psi,h}\times G(\bb{A}_f)/K$ the locally closed subspace of $\mrm{Sh}_K^{\mrm{min}}(G,X)(\bb{C})$ induced by the inclusion $D_{\Psi,h}\hookrightarrow X^*$. 
By \cite[6.3]{Pin89} and Proposition \ref{prop-iphi}, we have $\mrm{Z}_{\Psi,K}(\bb{C})=\mrm{Stab}_{Q(\bb{Q})}(D_{\Phi,h})\backslash D_{\Phi,h}\times G(\bb{A}_f)/K=\disju_{[\Phi_i]\in I(Q_\Phi)} \mrm{Z}_{[\Phi_i],K}(\bb{C})$, where each $$\gls{zmrmZPhiKC}=\Delta_{\Phi_i,K}\backslash\mrm{Sh}_{K_{\Phi_i,h}}(G_{\Phi_i,h},D_{\Phi_i,h})(\bb{C}),$$ and 
$$\gls{DeltaPhiK}:=\mrm{Stab}_{Q_{\Phi_i}(\bb{Q})}(D_{\Phi_i})\cap P_{\Phi_i}(\bb{A}_f)g_{\Phi_i}Kg_{\Phi_i}^{-1}/P_{\Phi_i}(\bb{Q}).$$
Indeed, by Proposition \ref{prop-iphi}, there is an isomorphism $\mrm{Stab}_{Q(\bb{Q})}(D_\Phi)P_\Phi(\bb{A}_f)\backslash G(\bb{A}_f)/K\iso I(Q_\Phi)$ mapping $[g_i]$ in the double coset to $[\Phi_i]=[(Q_\Phi,X^+_\Phi,g_i)]$, and $\Psi$ only depends on the $G(\bb{Q})$-conjugacy class $[Q_\Phi]$ of $Q_\Phi$. Since the strong approximation theorem holds for unipotent groups, we have $\mrm{Z}_{[\Phi_i],K}(\bb{C})=\Delta_{\Phi_i,K}\bss (P_\Phi(\bb{Q})\bss D_{\Phi,h}\times P_\Phi(\A)/K_{\Phi_i})$. By \S\ref{dbcoset}, we can write the disjoint union as above. If we fix a set of representatives for $I(Q_\Phi)$ as above such that $\Phi_i=(Q_\Phi,X_\Phi^+,g_i)$, then
$$\Delta_{\Phi_i,K}\iso \mrm{Stab}_{Q_\Phi(\bb{Q})}(D_\Phi)\cap P_\Phi(\bb{A}_f)g_iKg_i^{-1}/P_\Phi(\bb{Q}).$$

In other words, for any admissible $\bb{Q}$-parabolic subgroup $Q$, we can choose any $\Phi=(Q,X^+_\Phi,g_\Phi)$ and define $\Psi$ as above. Then $\mrm{Z}_{\Psi,K}(\bb{C})$ is the union of all $\mrm{Z}_{[\Phi_i],K}(\bb{C})$ associated with equivalence classes of cusp label representatives formed by the $G(\bb{Q})$-conjugacy class $[Q]$ and other data. Hence, it makes sense to denote $\mrm{Z}_{\Psi,K}(\bb{C})$ as $\mrm{Z}_{[Q],K}(\bb{C})$.\par
Then \begin{equation}\label{[Q]}
\mrm{Sh}_K^\mrm{min}(G,X)(\bb{C})\iso\disju_{[Q]} \mrm{Z}_{[Q],K}(\bb{C}),\end{equation} 
where the disjoint union runs over the conjugacy classes of admissible $\bb{Q}$-parabolic subgroups of $G$.\par
Let $\Phi=(Q,X^+,g_\Phi)$. Let $D^+_\Phi$ be the connected component of $D_\Phi$ containing $X^+$. Then the image $D^+_{\Phi,h}$ of $D_\Phi^+$ in $D_{\Phi,h}$ is a connected component of $D_{\Phi,h}$. We have $\mrm{Z}_{[\Phi_i],K}(\bb{C})\iso\mrm{Stab}_{Q(\bb{Q})}(D_\Phi)\backslash D_{\Phi,h}\times \stb_{Q(\bb{Q})}(D_\Phi)P_{\Phi}(\bb{A}_f)/K_{\Phi_i}\iso\disju_{c_j\in C(\Phi_i)}D_{\Phi,h}^+/\Gamma_{\Phi_i}(c_j)$, where $C(\Phi_i)$ is a set of representatives of the double coset $\pi_0(\mrm{Z}_{[\Phi_i],K}(\bb{C}))\iso\mrm{Stab}_{Q(\bb{Q})}(D^+_\Phi)\bss \stb_{Q(\bb{Q})}(D_\Phi)P_\Phi(\bb{A}_f)/K_{\Phi_i}$ and $\Gamma_{\Phi_i}(c_j):=\mrm{Stab}_{Q(\bb{Q})}(D_\Phi^+)\cap c_jK_{\Phi_i}c_j^{-1}.$\par
By \cite[Prop. 9.24 and Thm. 11.8]{Pin89} again, $\mrm{Z}_{[\Phi],K}(\bb{C}),$ $\mrm{Z}_{\Psi,K}(\bb{C})$ and $\mrm{Z}_{[Q],K}(\bb{C})$ admit algebraic structures over the reflex field $E$ of $(G,X)$, and we denote them by $\gls{zmrmZPhiK}$, $\mrm{Z}_{\Psi,K}$ and $\mrm{Z}_{[Q],K}$, respectively. We also adopt similar notation for their base changes (cf. Convention \ref{conv-mixed-sh}).\par
We have the following observation from the discussions above:
\begin{lem}\label{comp-min}Denote by $\overline{E}$ the algebraic closure of $E:=E(G,X)$. 
Fix a connected component $X^+$ of $X$ so that we can define an isomorphism $\pi_0(\mrm{Sh}_K(G,X)_{\overline{E}})\iso G(\bb{Q})_+\backslash G(\bb{A}_f)/K$ mapping $G(\bb{Q})_+\backslash X^+\times G(\bb{Q})_+K/K$ to $[1]$.
Denote the connected component of $\mrm{Sh}_K^{\mrm{min}}(G,X)_{\overline{E}}$ corresponding to $[g]$ by  $\mrm{Sh}_{K,\overline{E}}^{\mrm{min},[g]}$ for $[g]\in G(\bb{Q})_+\backslash G(\bb{A}_f)/K$.
Let $\mrm{Z}_{[Q],K,\overline{E}}^{[g]}$ be the pullback of $\mrm{Z}_{[Q],K,\overline{E}}\hookrightarrow \mrm{Sh}_K^{\mmin}(G,X)_{\overline{E}}$ to $\mrm{Sh}_{K,\overline{E}}^{\mrm{min},[g]}$.
For any $Q$, there is an isomorphism $\pi_0(\mrm{Z}_{[Q],K,\overline{E}})\iso \mrm{Stab}_{Q(\bb{Q})}(D^+_{Q,X^+})\backslash G(\bb{A}_f)/K$ such that $\mrm{Z}_{[Q],K,\overline{E}}^{[g]}\iso\disju_{[g^\prime]\mapsto [g]} \mrm{Z}_{[Q],K,\overline{E}}^{[g^\prime]}$, where the disjoint union runs over the preimages of $[g]$ under the natural quotient  $$\pi_{Q,X^+}:\mrm{Stab}_{Q(\bb{Q})}(D_{Q,X^+}^+)\backslash G(\bb{A}_f)/K\lra G(\bb{Q})_+\backslash G(\bb{A}_f)/K.$$

Additionally, in the quotient above, the contribution of $\pi_0(\mrm{Z}_{[\Phi],K,\overline{E}})$ is a subquotient
$$\pi_{\Phi}:\stb_{Q(\bb{Q})}(D^+_{Q,X^+})\cap P_\Phi(\A)g_\Phi Kg^{-1}_\Phi\bss P_\Phi(\A)g_\Phi K/K\lra G(\bb{Q})_+\bss G(\A)/K.$$
\end{lem}
\begin{proof}It suffices to check the statement over $\bb{C}$. 
By the definition of Satake topology, the closure $X^{+,*}$ of $X^+$ in $X^*$ has a stratification consisting of connected components $D^+_{\Phi,h}$ of $D_{\Phi,h}$ such that $D^+_{\Phi}$ contains $\tau(X^+)$. Hence, $\mrm{Sh}_K^{\mrm{min}}(G,X)(\bb{C})=G(\bb{Q})_+\backslash X^{+,*}\times G(\bb{A}_f)/K$, and the strata in (\ref{[Q]}) can be written as $G(\bb{Q})_+\backslash G(\bb{Q})_+D^+_{\Phi,h}\times G(\bb{A}_f)/K=\mrm{Z}_{[Q],K}(\bb{C})$, for some connected component $D^+_{\Phi}$ containing $\tau(X^+)$. Then the map between the $\pi_0$ of $\mrm{Z}_{[Q],K,\bb{C}}$ and $\mrm{Sh}^\mrm{min}_K(G,X)(\bb{C})$ can be written as $\pi_0(\mrm{Z}_{[Q],K,\bb{C}})\iso G(\bb{Q})_+\backslash (G(\bb{Q})_+/\mrm{Stab}_{Q(\bb{Q})}(D^+_{Q,X^+}))\times G(\bb{A}_f)/K\to G(\bb{Q})_+\backslash (G(\bb{Q})_+\backslash G(\bb{Q})_+)\times G(\bb{A}_f)/K\iso \pi_0(\mrm{Sh}^{\mrm{min}}_K(G,X)_{\bb{C}}).$ We get the desired statement.
\end{proof}
For any cusp label $[\Phi]$ with a representative $\Phi$ and any geometrically connected component $\sh^+_K$ of $\sh_K(G,X)$ defined over $\overline{E}$, we denote by $\sh_K^{+,\mrm{min}}$ the schematic closure of $\sh_K^+$ in $\sh_{K,\overline{E}}^{\mrm{min}}$ and denote by $\mrm{Z}_{[\Phi],K}^+$ the pullback $\mrm{Z}_{[\Phi],K,\overline{E}}\times_{\sh^\mrm{min}_{K,\overline{E}}}\sh^{+,\mrm{min}}_K$.
\begin{rk}
Although in the theory of compactifications of connected Shimura varieties over $\bb{C}$ (see \cite[Ch. III]{AMRT10}), the strata defined for minimal compactifications are all connected, the strata $\mrm{Z}_{[\Phi],K}^+$ induced from Pink's theory for every connected component $\sh_K^{+,\mrm{min}}$ are not connected in general. 
\end{rk}
\subsubsection{}\label{tor-comp}
Now we recall Pink's work on toroidal compactifications of $\sh_K(G,X)$. Let $Z:=Z_G$ be the center of $G$. Let $\Phi$ be any cusp label representative. Denote by $Z_\Phi$ the center of $P_\Phi$. Let $Z_\Phi(\bb{Q})^\circ$ be the subgroup of $Z_\Phi(\bb{Q})$ that acts trivially on $D_\Phi$. Let $K_{\Phi,W}$ be the projection of $(Z_\Phi(\bb{Q})^\circ\times W_\Phi(\bb{A}_f))\cap K_\Phi$ to the second factor $W_\Phi(\bb{A}_f)$, and let $K_{\Phi,U}:=K_{\Phi,W}\cap U(\bb{A}_f)$. So $K_{\Phi,U}$ is the projection of $(Z_\Phi(\bb{Q})^\circ\times U_\Phi(\A))\cap K_\Phi$ to the second factor.\par 
Moreover, define $K'_{\Phi,U}$ (resp. $K'_{\Phi,W}$) to be the projection of $(Z(\bb{Q})\times U_\Phi(\A))\cap g_\Phi K g_\Phi^{-1}$ (resp. $(Z(\bb{Q})\times W_\Phi(\A))\cap g_\Phi K g_\Phi^{-1}$) to the second factor $U_\Phi(\A)$.
\begin{lem}\label{lem-large-lattices}
If $K$ is neat, $K_{\Phi,U}$ (resp. $K_{\Phi,W}$) is contained in $K_{\Phi,U}'$ (resp. $K_{\Phi,W}'$) as a finite index subgroup. Moreover, if we assume that $Z$ is isogenous to a product of $\bb{Q}$-split tori and $\bb{R}$-anisotropic tori, then $K_{\Phi,U}=K_{\Phi,U}'$ and $K_{\Phi,W}=K_{\Phi,W}'$.
\end{lem}
\begin{proof}We only prove the result for $K_{\Phi,U}$ and $K_{\Phi,U}'$. 
By the proof of \cite[Cor. 4.10]{Pin89}, the center of $ZP_\Phi/ZW_\Phi$ is isogenous to a product of $\bb{Q}$-split tori and compact-type tori. Hence, if $K_\Phi$ is neat, then $K_{\Phi,U}$ is contained in the projection of $(Z(\bb{Q})\cdot Z_\Phi(\bb{Q})^\circ\times U_\Phi(\bb{A}_f))\cap g_\Phi Kg_\Phi^{-1}$ to the second factor $U_\Phi(\bb{A}_f)$, and the latter is isomorphic to $K_{\Phi,U}'$. By neatness assumption of $K$, the only neat arithmetic subgroup of $Z(\bb{Q})$ is the trivial group, so $K_{\Phi,U}=K'_{\Phi,U}$.
\end{proof}
Define the lattice $\gls{SKPhivee}:=(U_\Phi(\bb{Q})\cap K_{\Phi,U})(-1)\sbst U_\Phi(\bb{R})(-1)$, and let $\gls{EKPhi}$ be the algebraic torus whose cocharacter is the lattice $\mbf{S}_{K_\Phi}^\vee$ in $U_\Phi(\bb{Q})(-1)$. Let $\gls{PP}$ be the union of all $\mrm{int}(\gamma^{-1})\mbf{P}_{\Phi'}^+\sbst U_\Phi(\bb{R})(-1)$ such that $\Phi\xrightarrow[\sim]{(\gamma,q)_K}\Phi'$. (Note that $\mbf{P}_\Phi$ is not a rational polyhedral cone.) Define $\gls{LambdaKPhi}:=(U_\Phi(\bb{Q})\cap K_{\Phi,U}')(-1)\sbst U_\Phi(\bb{R})(-1)$. By \cite[3.13, Prop. 11.10 and Thm. 11.8]{Pin89}, $\sh_{K_\Phi}$ is a torus torsor under $\mbf{E}_{K_\Phi}$ over $\overline{\sh}_{K_\Phi}$.\par
Let us recall the definitions related to cone decompositions.
\begin{definition}[{\cite[6.4 and 7.12]{Pin89} and \cite[2.1.23]{Mad19}}]
An \textbf{admissible} rational polyhedral cone decomposition $\Sigma$ for $(G,X,K)$ is an association $\Phi\in \ca{CLR}(G,X)\mapsto \Sigma(\Phi)$, where $\Sigma(\Phi)$ is a rational polyhedral cone decomposition of $\mbf{P}_\Phi$, such that: if $\Phi_1\xrightarrow{(\gamma,p)_K}\Phi_2$, $\mrm{int}(\gamma^{-1})\Sigma(\Phi_2)=\Sigma(\Phi_1)|_{\mrm{int}(\gamma^{-1})\mbf{P}_{\Phi_2}}$ under the natural embedding $\mrm{int}(\gamma^{-1}):U_{\Phi_2}(\bb{R})(-1)\hookrightarrow U_{\Phi_1}(\bb{R})(-1)$.
\end{definition}
Let $\Sigma$ be any admissible rational polyhedral cone decomposition. Then $\Sigma(\Phi)=\Sigma((q,k)\Phi)$ for any $q\in P_\Phi(\bb{A}_f)$ and $k\in K$.
If $\Phi_1\xrightarrow[\sim]{(\gamma,q)_K}\Phi_2$, then $\gamma U_{\Phi_1}(\bb{R})\gamma^{-1}=U_{\Phi_2}(\bb{R})$. Note that any $\gamma\in G(\bb{Q})$ has a left action on $\Sigma$ such that for any cone $\sigma\in\Sigma(\Phi_1)$ where $\sigma^\circ\sbst \mbf{P}^+_{\Phi_1}$, the conjugation $\gamma\sigma\gamma^{-1}$ induces a cone in $\Sigma(\Phi_2)$ whose interior lies in $\mbf{P}^+_{\Phi_2}.$ 
Hence, $G(\bb{Q})$ acts on the set of pairs $(\Phi,\sigma)$ where $\sigma\in\Sigma(\Phi)$ such that $\sigma^\circ\sbst \mbf{P}^+_\Phi$. As in \cite[2.1.26]{Mad19}, we say any two such pairs $(\Phi_1,\sigma_1)$ and $(\Phi_2,\sigma_2)$ are \textbf{equivalent} if $\Phi_1\xrightarrow[\sim]{(\gamma,q)_K}\Phi_2$ and $ \gamma \sigma_1\gamma^{-1}=\sigma_2$; we denote this equivalence by $(\Phi_1,\sigma_1)\xrightarrow[\sim]{(\gamma,q)_K}(\Phi_2,\sigma_2)$. We denote the set of equivalence classes in the form $[(\Phi,\sigma)]$ by \gls{CuspKGXSigma}, the elements in it are called \textbf{cusp labels with cones}. There is a partial order on $\cusp_K(G,X,\Sigma)$ defined as follows: $[(\Phi_1,\sigma_1)]\preceq[(\Phi_2,\sigma_2)]$ if and only if  $\Phi_1\xrightarrow{(\gamma,q)_K}\Phi_2$ and $\gamma^{-1}\sigma_2\gamma$ is a face of $\sigma_1$.\par
Define $\Sigma^+(\Phi)$ to be the cones $\sigma\in \Sigma(\Phi)$ satisfying $\sigma^\circ\sbst \mbf{P}^+_\Phi$.\par
From Proposition \ref{prop-iphi}, we can see that when we fix a $\Phi$, the stabilizer of $[\Phi]$ in $Q_\Phi(\bb{Q})$ is $\wdtd{\Delta}_{\Phi,K}:=\stb_{Q_\Phi(\bb{Q})}(D_{\Phi})\cap P_\Phi(\A)g_\Phi Kg^{-1}_\Phi$. Let $\Delta_{\Phi,K}:=\wdtd{\Delta}_{\Phi,K}/P_\Phi(\bb{Q})$. This means that for any $\Upsilon=[(\Phi,\sigma)]\in\mrm{Cusp}_K(G,X,\Sigma)$, $\Upsilon$ is determined by a pair $([\Phi],[\sigma])$, where $[\Phi]\in \mrm{Cusp}_K(G,X)$ and $[\sigma]$ is a $\Delta_{\Phi,K}$-orbit of $\sigma$ in $U_\Phi(\bb{R})$. Let $\gls{DeltaPhiKcirc}$ be the normal subgroup of $\Delta_{\Phi,K}$ that stabilizes $\sigma$. By \cite[2.1.19]{Mad19}, when $K$ is neat, $\Delta^\circ_{\Phi,K}$ is independent of the choice of $\sigma$ such that $\sigma^\circ\sbst \mbf{P}^+_\Phi$, and it fixes $\mbf{P}_\Phi$.\par
\begin{definition}
Let $\Sigma$ be any admissible rational polyhedral cone decomposition for $(G,X,K)$, we can further define the following terminologies:
\begin{itemize}
\item The decomposition $\Sigma$ is called \textbf{complete} if, for any $\Phi$, $\mbf{P}_\Phi=\bigcup_{\sigma\in \Sigma(\Phi)}\sigma^\circ$. \par
\item The decomposition $\Sigma$ is called \textbf{finite} if $G(\bb{Q})\bss\Sigma/K$ is finite.\par
\item The decomposition $\Sigma$ is a decomposition \textbf{without self-intersections} if, for any $\tau,\sigma\in \Sigma(\Phi)$ such that $\tau$ is a face of $\sigma$ and $\sigma^\circ\sbst \mbf{P}_{\Phi}^+$, and any $\gamma\in G(\bb{Q})$, $q\in P_\Phi(\A)$ and $k\in K$ such that $ (\gamma,q,k)\tau$ is also a face of $\sigma$, we have that $(\gamma,q,k)\tau=\tau$. (This condition is equivalent to \cite[Condition 6.2.5.25]{Lan13} if $K$ is neat; see \cite[Rmk. 6.2.5.26]{Lan13} and the discussion above.)\par
\item The decomposition $\Sigma$ is \textbf{smooth} (resp. \textbf{projective}) if $\Sigma(\Phi)$, as a cone decomposition of $\mbf{P}_\Phi\sbst U_\Phi(\bb{R})(-1)$, is smooth (resp. projective) with respect to $\mbf{\Lambda}_{K_\Phi}$ for any $\Phi$; note that, if $(G,X)$ is a Hodge-type Shimura datum, we can replace the $\mbf{\Lambda}_{K_\Phi}$ in the last sentence with $\mbf{S}_{K_\Phi}^\vee$ by Lemma \ref{lem-large-lattices}.\par
In the last paragraph, being ``smooth with respect to $\mbf{\Lambda}_{K_\Phi}$'' has the following precise meaning: We require that, for any $\sigma\in \Sigma(\Phi)$, there is a $\bb{Z}$-basis $\{v_1,v_2,...,v_n\}$ of $\mbf{\Lambda}_{K_\Phi}$ such that $\sigma=\sum_{i=1}^k\bb{R}_{\geq 0}v_i$ for some $k\leq n$.\par
Moreover, being ``projective with respect to $\mbf{\Lambda}_{K_\Phi}$'' has the following precise meaning: We require that there is a polarization function $\mrm{pol}_{K_\Phi}:\mbf{P}_{\Phi}\to \bb{R}_{\geq 0}$ on $\mbf{P}_\Phi$ such that
\begin{itemize}
    \item $\mrm{pol}_{K_\Phi}$ is continuous and $\Delta_{\Phi,K}$-invariant, and is positive on $\mbf{P}_\Phi\bss \{0\}$, 
    \item the restriction $\mrm{pol}_{K_\Phi}|_{\mbf{\Lambda}_{K_{\Phi}}\cap \mbf{P}_\Phi}$ takes values in $\bb{Z}_{\geq 0}$,
    \item $\mrm{pol}_{K_\Phi}$ is piecewise linear, and is linear on a rational polyhedral cone $\sigma\sbst \mbf{P}_\Phi$ if and only if $\sigma\in \Sigma(\Phi)$, and
    \item we have $\pol_{K_\Phi}(x+y)\geq \pol_{K_\Phi}(x)+\pol_{K_\Phi}(y)$ for any $x$ and $y$ in $\mbf{P}_\Phi$.
\end{itemize}
\end{itemize}
\end{definition}
\begin{convention}\label{conv-cone}
    We always choose $\Sigma$ to be admissible, complete, finite, and without self-intersections. From now on, whenever we say $\Sigma$ is an admissible cone decomposition, it will mean that $\Sigma$ is an admissible rational polyhedral cone decomposition that is complete, finite, and without self-intersections. The refinements of cone decompositions we take in this paper are also admissible, complete, finite, and without self-intersections.
\end{convention}
\begin{definition}[{cf. \cite{Lan19}}]\label{def-func-cones}
Let $f:(G_1,X_1)\to (G_2,X_2)$ be a morphism of Shimura data. Let $g\in G_2(\A)$. Let $K_1\sbst G_1(\A)$ and $K_2\sbst G_2(\A)$ be neat open compact subgroups such that $f(K_1)\sbst g K_2g^{-1}$. By \cite[3.5 and Thm. 12.4]{Pin89}, there is a morphism over the reflex field of $(G_1,X_1)$: 
$$f(g)_{K_1,K_2}:\sh_{K_1}(G_1,X_1)\to \sh_{K_2}(G_2,X_2).$$
\begin{enumerate}
\item Suppose that $(G_1,X_1)=(G_2,X_2)=(G,X)$ and $f=\mrm{id}$. Then, in this case, $f(g)_{K_1,K_2}=[g]_{K_1,K_2}$ is the right action of $g$. Suppose that there is an admissible cone decomposition $\Sigma$ for $(G,X,K_2)$, by \cite[6.7]{Pin89}, the pullback of $\Sigma$, denoted by $[g]_{K_1,K_2}^*(\Sigma)$ is an admissible cone decomposition for $(G,X,K_1)$. We say $[g]_{K_1,K_2}^*(\Sigma)$ is \textbf{induced} by $\Sigma$ under $[g]_{K_1,K_2}$.
\item Suppose that there are admissible cone decompositions $\Sigma_1$ for $(G_1,X_1,K_1)$ and $\Sigma_2$ for $(G_2,X_2,K_2)$. Suppose that for any $\Phi_1\in \ca{CLR}(G_1,X_1)$ mapping to some $\Phi_2\in \ca{CLR}(G_2,X_2)$ under $f(g)_{K_1,K_2}$, the induced map $U_{\Phi_1}(\bb{R})(-1)\to U_{\Phi_2}(\bb{R})(-1)$ is injective. Then we say that $\Sigma_1$ is \textbf{induced} by $\Sigma_2$ under $f(g)_{K_1,K_2}$ if, for any $\sigma_1\in \Sigma_1(\Phi_1)$, $\sigma_1$ is exactly the preimage of some $\sigma_2\in\Sigma_2(\Phi_2)$. Denote $\Sigma_1:=f(g)_{K_1,K_2}^*(\Sigma_2)$.
If the image of $\sigma_1$ is exactly $\sigma_2\in \Sigma_2(\Phi_2)$ itself for any $\sigma_1\in \Sigma_1(\Phi_1)$ and any pair $(\Phi_1,\Phi_2)$ described above, then we say that $\Sigma_1$ and $\Sigma_2$ are \textbf{strictly compatible} with each other.
\item Suppose that there are admissible cone decompositions $\Sigma_1$ for $(G_1,X_1,K_1)$ and $\Sigma_2$ for $(G_2,X_2,K_2)$. We say that $\Sigma_1$ and $\Sigma_2$ are \textbf{compatible} with each other if, for any $\Phi_1\in\ca{CLR}(G_1,X_1)$ mapping to $\Phi_2\in\ca{CLR}(G_2,X_2)$ under $f(g)_{K_1,K_2}$, the image of any $\sigma_1\in \Sigma_1(\Phi_1)$ under the map $\mbf{P}_{\Phi_1}\to \mbf{P}_{\Phi_2}$ is contained in some $\Sigma_2\in \Sigma_2(\Phi_2)$. 
\end{enumerate}
\end{definition}
Let $\Sigma$ be an admissible cone decomposition for $(G,X,K)$. Assume that $K$ is a neat open compact subgroup of $G(\A)$. By the main theorems of \cite{Pin89}, there is a toroidal compactification $\sh_K^\Sigma:=\sh_K^\Sigma(G,X)$ of $\sh_K$ over the reflex field $E:=E(G,X)$, such that $\sh_K^\Sigma$ is a normal algebraic space that is proper over $\spec E$. If $\Sigma$ is smooth (resp. projective), then $\sh_K^\Sigma$ is smooth with a normal crossings boundary divisor $D:=(\sh_K^\Sigma\bss \sh_K)_{\mrm{red}}$ (resp. is projective and representable by a scheme). There is a good stratification labeled by $\cusp_K(G,X,\Sigma)$: 
$$\sh_K^\Sigma\iso \disju\limits_{\Upsilon\in \cusp_K(G,X,\Sigma)}\mrm{Z}_{\Upsilon,K}.$$ 
For any $\Upsilon=[(\Phi,\sigma)]\in \mrm{Cusp}_K(G,X,\Sigma)$, the stratum $\gls{zmrmZUpsilon}$ in $\sh_K^\Sigma$ is a normal variety over $E$ and is isomorphic to $\Delta^\circ_{\Phi,K}\backslash\sh_{K_\Phi,\sigma}$, where $\gls{shKPhisubsigma}$ is the $\sigma$-stratum of the twisted affine torus embedding $\sh_{K_\Phi}\hookrightarrow \gls{shKPhibkSigma}$ over $\overline{\sh}_{K_\Phi}$ defined by $\Sigma^+(\Phi)$. Note that $\sh_{K_\Phi,\sigma}$ is closed in $\sh_{K_\Phi}(\sigma)$ and is locally closed in $\sh_{K_\Phi}(\Sigma)$.\par
Moreover, there is a canonical isomorphism:
$$(\sh_K^\Sigma)^{\wat{\ }}_{\mrm{Z}_{\Upsilon,K}}\iso \Delta^\circ_{\Phi,K}\bss\sh_{K_\Phi}(\sigma)^{\wat{\ }}_{\sh_{K_\Phi,\sigma}}.$$\par
Let $\widetilde{\Delta}^\circ_{\Phi,K}$ be the preimage of $\Delta^\circ_{\Phi,K}$ in $\mrm{Stab}_{Q(\bb{Q})}(D_\Phi)\cap P_\Phi(\bb{A}_f)g_{\Phi}Kg_{\Phi}^{-1}$. Since $\sh_{K_\Phi,\sigma}$ is the closed stratum of the twisted affine torus embedding $\sh_{K_\Phi}\hookrightarrow\gls{shKPhibksigma}$ over an abelian scheme torsor $\overline{\sh}_{K_\Phi}$ over $\sh_{K_\Phi,h}$ (see \cite[3.21]{Pin89} and \cite[2.1.10]{Mad19}), $\pi_0(\sh_{K_\Phi,\sigma,\bb{C}})\iso \pi_0(\sh_{K_\Phi,h,\bb{C}})\iso\stb_{P_\Phi(\bb{Q})}(D_\Phi^+)\backslash P(\bb{A}_f)/K_\Phi.$ Moreover, $\pi_0(\sh_{K,\bb{C}}^\Sigma)\iso\pi_0(\mrm{Sh}_{K,\bb{C}})$ since $\mrm{Sh}_{K}$ is dense in $\sh_K^\Sigma$. \par 
We have the following observation:
\begin{lem}\label{comp-tor}
Let $X^+$ and $Q$ be as in Lemma \ref{comp-min}. 
For any $[g]\in G(\bb{Q})_+\backslash G(\bb{A}_f)/K$ and any $\Upsilon=[(\Phi,\sigma)]\in\mrm{Cusp}_K(G,X,\Sigma)$ with $\Phi=(Q,X^+,g_\Phi)$ so that $[\Phi]\in I(Q)$, denote by $\mrm{Z}_{\Upsilon,\overline{E}}^{[g]}$ (resp. $\sh_{K_\Phi,\sigma,\overline{E}}^{[g]}$) the pullback of $\mrm{Z}_{\Upsilon,\overline{E}}\hookrightarrow \sh_{K,\overline{E}}^\Sigma$ (resp. $\sh_{K_\Phi,\sigma,\overline{E}}\to \sh^\Sigma_{K,\overline{E}}$) to the connected component $\sh_{K,\overline{E
}}^{\Sigma,[g]}$ corresponding to $[g]$. Then $ \mrm{Z}_{\Upsilon,\overline{E}}^{[g]}=\disju_{[g^\prime]\mapsto [g]}\mrm{Z}_{\Upsilon,\overline{E}}^{[g^\prime]}$ and $\sh^{[g]}_{K_\Phi,\sigma,\overline{E}}=\disju_{[g^{\prime\prime}]\mapsto [g]}\sh_{K_\Phi,\sigma,\overline{E}}^{[g^{\prime\prime}]}$, 
where the first disjoint union runs over the preimage of $[g]$ under 
\begin{equation}\label{pi-upsilon}\pi_0(\mrm{Z}_{\Upsilon,\overline{E}})\iso\mrm{Stab}_{\widetilde{\Delta}^\circ_{\Phi,K}}(D^+_{Q,X^+})\backslash P_\Phi(\bb{A}_f)g_\Phi K/K\lra \pi_0(\mrm{Sh}_{K,\overline{E}}),\end{equation} and 
where the second disjoint union runs over the preimage of $[g]$ under
\begin{equation}\pi_0(\sh_{K_\Phi,\sigma,\overline{E}})=\mrm{Stab}_{P_\Phi(\bb{Q})}(D^+_{Q,X^+})\backslash P_\Phi(\bb{A}_f)g_\Phi K/K\lra \pi_0(\mrm{Sh}_{K,\overline{E}}).\end{equation}
\end{lem}
\begin{proof}It suffices to check this over $\bb{C}$. Since $g_\Phi$ is a representative of $\mrm{Stab}_{Q(\bb{Q})}(D_\Phi)P_\Phi(\bb{A}_f)\backslash G(\bb{A}_f)/K$, and 
$\mrm{Stab}_{P_\Phi(\bb{Q})}(D_\Phi^+)\backslash P_\Phi(\bb{A}_f)/K_\Phi$ is isomorphic to $ \mrm{Stab}_{P_\Phi(\bb{Q})}(D_\Phi^+)\backslash P_\Phi(\bb{A}_f)g_\Phi K/K$, 
the map between connected components $\pi_0(\sh_{K_\Phi,\sigma,\bb{C}})\to\pi_0(\mrm{Z}_{[Q],K,\bb{C}})$ induced by the natural morphism $\sh_{K_\Phi,\sigma,\bb{C}}\to \mrm{Z}_{[Q],K,\bb{C}}$ is the map   $\mrm{Stab}_{P_\Phi(\bb{Q})}(D_\Phi^+)\backslash P_\Phi(\bb{A}_f)g_\Phi K/K\to \mrm{Stab}_{Q(\bb{Q})}(D_\Phi^+)\backslash G(\bb{A}_f)/K$ induced by $P_\Phi(\bb{A}_f)g_\Phi K\hookrightarrow G(\bb{A}_f)$ and $\mrm{Stab}_{P_\Phi(\bb{Q})}(D_\Phi^+)\hookrightarrow \mrm{Stab}_{Q(\bb{Q})}(D_\Phi^+)$. Then the statement follows from Lemma \ref{comp-min} and the paragraph above. 
\end{proof}
\begin{rk}
Let $\sh_K^+$ be a geometrically connected component of $\sh_K(G,X)$. Denote by $\sh_K^{+,\Sigma}$ and $\sh_K^{+,\mrm{min}}$ its closure in $\sh_K^\Sigma(G,X)$ and $\sh_K^{\mrm{min}}(G,X)$ respectively. The components $\sh_K^{+,\Sigma}$ and $\sh_K^{+,\mrm{min}}$ might not intersect with all strata of the toroidal and minimal compactifications nontrivially (see the appendix of \cite{LW15} for an affirmative answer about this for the PEL type A and C case); moreover, from (\ref{pi-upsilon}), we can see that for toroidal compactifications, each stratum appearing in some $\sh_K^{+,\Sigma}$, if not trivial, might not be geometrically connected, since $\Delta^\circ_{\Phi,K}$ might be trivial and the action of $\Delta_{\Phi,K}$ might not be trivial, and since $\mrm{Z}_{[\Phi],K}^+$ might not be connected.
\end{rk}
\begin{convention}\label{conv-K-sub-in-delta}
If $K$ is clear in the context, we will omit the subscript $K$ in $\Delta_{\Phi,K}$, $\Delta_{\Phi,K}^\circ$, $\mrm{Z}_{[\Phi],K}$, $\mrm{Z}_{[Q],K}$ and $\mrm{Z}_{[\Phi],K}^+$.
\end{convention}
\subsubsection{} We have the following generalization of \cite[Lem. 2.1.2]{Kis10} for mixed Shimura varieties.
\begin{lem}\label{lem-emb-mix}
Let $f:(P^1,\x_1)\hookrightarrow(P^2,\x_2)$ be an embedding of two mixed Shimura data. For $i=1$ or $2$, assume that $K^i:=K_p^iK^{i,p}$, where $K_p^i$ is an open compact subgroup of $P^i(\bb{Q}_p)$ and $K^{i,p}\sbst P^i(\bb{A}_f^p)$ is neat open compact. Assume that $K^2_p\cap P^1(\bb{Q}_p)=K^1_p$. Then, for any $K^{1,p}$, there is a $K^{2,p}$ containing $K^{1,p}$ such that $f$ induces a closed embedding 
\begin{equation}\label{eq-emb-mix-sh}
    \mrm{Sh}_{K^1}(P^1,X_1)\lra\mrm{Sh}_{K^2}(P^2,X_2).
\end{equation}
\end{lem}
\begin{proof}
This is essentially \cite[Lem. 2.1.2]{Kis10}, with some inputs in \cite{Del71} replaced by those in \cite{Pin89}. 
In fact, since (\ref{eq-emb-mix-sh}) is finite, it suffices to show that 
$$\varprojlim_{K^{1,p}}\mrm{Sh}_{K^{1,p}K_{p}^1}(P^1,\x_1)\lra \varprojlim_{K^{2,p}}\mrm{Sh}_{K^{2,p}K_{p}^2}(P^2,\x_2)$$
is injective, where the inverse limits run over all neat open compact $K^{i,p}$ for $i= 1, 2$. 
Let $Z^i(\bb{Q})^\circ$ be the subgroup of the center $Z^i(\bb{Q})$ of $P^i(\bb{Q})$ whose action on $\x_i$ is trivial, and let $\Gamma_Z^i$ be any arithmetic subgroup of $Z^i(\bb{Q})^\circ$. Denote by $(\Gamma^i_Z)^-$ the closure of $\Gamma^i_Z$ in $P^i(\bb{A}_f)$.
By \cite[Lem. 3.7(b)]{Pin89}, $$ \varprojlim_{K^{1,p}}\mrm{Sh}_{K^{1,p}K_{p}^1}(P^1,\x_1)(\bb{C})\iso P^i(\bb{Q})\backslash \x_i\times P^i(\bb{A}_f)/(\Gamma^i_Z)^- K_p^i.$$
Note that the identity component of the center $Z^i$ of $P^i$ is a torus, since the Lie algebra $\mrm{Lie}\ W_i$ of the unipotent radical $W^i$ of $P^i$ has nontrivial weight by \cite[Def. 2.1]{Pin89} and $G^i:=P^i/W^i$ is reductive. By a theorem of Chevalley (see \cite[Thm. 1]{Che51} and \cite[2.0.10]{Del79}), any finite index subgroup of the group of units $U^i$ of $Z^i(\bb{Q})$ is open in $U^i$ under the topology induced by $P^i(\bb{A}_f^p)$ and $P^i(\bb{A}_f)$, so $(\Gamma_Z^i)^-$ is the completion $(\Gamma_Z^i)^{\widehat{}}$ with respect to finite index subgroups in $U^i$, and $(\Gamma_Z^i)^-$ is also the closure of $\Gamma_Z^i$ in $P^i(\bb{A}_f^p)$. 
Let $\widehat{P}^i(\bb{Q}):=\varprojlim P^i(\bb{Q})/U^{i,n}$ be the completion with respect to $U^{i,n}$. Note that the image of $Z^1(\bb{Q})$ in $P^2(\bb{Q})$ is contained in the $\bb{R}$-points of the centralizer of $\bb{S}\xrightarrow{h_{x,\infty}^1}P^1_\bb{R}\to P^2_\bb{R}$, which is compact modulo $Z^2(\bb{R})W^2(\bb{R})$. Then there is a well-defined morphism $\widehat{P}^1(\bb{Q}) \backslash P^1(\bb{A}_f^p)\to \widehat{P}^2(\bb{Q})\backslash P^2(\bb{A}_f^p)$. From the argument in \cite[Lem. 2.1.2]{Kis10} verbatim, the morphism above and the morphism $\x_1\times P^1(\bb{Q}_p)/K_p^1\to \x_2\times P^2(\bb{Q}_p)/K_p^2$ are injective, and the desired morphism is also injective. Note that the argument of \cite[1.15.3]{Del71} used in the argument of \cite[Lem. 2.1.2]{Kis10} also works, if we replace $G^1$ and $G^2$ in the proof of \cite[1.15.3]{Del71} with $P^1$ and $P^2$, respectively.
\end{proof}
Let $(G,\ca{X},\hbar)$ be a pure Shimura datum. Recall that $G$ is a connected reductive group over $\bb{Q}$. Assume that $K_p$ is an open compact subgroup of $G(\bb{Q}_p)$ that contains a parahoric group $K_p^\circ:=\mathscr{G}^\circ(\bb{Z}_p)$, where $\mathscr{G}^\circ$ is a parahoric group scheme of $G_{\bb{Q}_p}$ over $\bb{Z}_p$. Let $K^p$ be a neat open compact subgroup of $G(\Ap)$ as before. Let \gls{Ep} be the maximal field extension of $E:=E(G,\ca{X})$ that is unramified at all primes dividing $p$. We record the following lemma:
\begin{prop}\label{prop-par-geo-conn}
With the assumptions above, the geometrically connected components of $\sh_K:=\sh_{K_pK^p}(G,\ca{X})$ are defined over $E^p$.    
\end{prop}
\begin{proof}
By \cite[Prop. 2.11]{Pin89}, $(G,\ca{X})$ embeds into a pure Shimura datum $(G_1, \ca{X}_1)$, which has the same reflex field $E$ as $(G,\ca{X})$. In fact, $(G_1,\ca{X}_1):=(T,\ca{Y})\times \hbar(G,\ca{X})$, where $(T,\ca{Y}):=(G,\ca{X})/G^\der$ is a pure Shimura datum associated with the cocenter $T:=G/G^\der$. 
Let $K_1\sbst G_1(\bb{A}_f)$ be a neat open compact subgroup such that $K_1:=K_{1,p}K^p_1$ with $K_{1,p}\sbst G_1(\bb{Q}_p)$ of the form $K_{1,p}=K_{T,p}\times K_p$, $K_{1,p}\cap G(\bb{Q}_p)=K_{p}$, and $K_1^p$ neat open compact in $G_1(\Ap)$. 
We can further choose $K_{T,p}\sbst T(\bb{Q}_p)$ to be the unique parahoric subgroup of $T(\bb{Q}_p)$. 
By Lemma \ref{lem-emb-mix}, there is a neat open compact subgroup $K^{p}_1\sbst G_1(\Ap)$ such that the induced morphism $\sh_{K}\hookrightarrow \sh_{K_1}$ is a closed embedding between smooth varieties of equal dimensions, so it is open and closed. Moreover, this embedding is defined over $E$ by \cite[Prop. 11.10]{Pin89}.
Then we are reduced to show the statement for $\hbar(G,\ca{X})$ and for $(T,\ca{Y})$. 
For the former one, the proof of \cite[Lem. 4.3.2 and Cor. 4.3.9]{KP15} works verbatim. For the latter one, note that $(T,\ca{Y})$ is a pure Shimura datum with $T$ a torus and $\ca{Y}$ finite. In this case, the proof similar to \emph{loc. cit}. also works. 
Indeed, by the discussion in \S\ref{subsubsec-can-model}, the Galois group $\mrm{Gal}(\overline{E}_T/E_T)$ acts on $\sh_{K_T}(\overline{E}_T)$ via $r_{E_T}(T,\ca{Y})/K_{T,p}: \mrm{Gal}(\overline{E}_T/E_T)\to T(\bb{Q})^+\bss T(\A)/K_{T,p}$. 
For any place $v|p$ of $E_T\sbst E$, the image of $\ca{O}^\times_{E_{T,v}}$ in $T(\bb{Q})^+\bss T(\A)/K_{T,p}$ is trivial under $r_{E_T}(T,\ca{Y})\circ\mrm{rec}_{E_T}^{-1}$, since $K_{T,p}$ is the unique Iwahori subgroup of $T(\bb{Q}_p)$, since $\ca{O}^\times_{E_T,v}$ is the unique Iwahori subgroup of $R_{E_{T,v}/\bb{Q}_p}\bb{G}_m$ by \cite[Lem. 2.5.18]{KP23}, and since $r_{E_T}(T,\ca{Y})\circ \mrm{rec}_{E_T}^{-1}|_{E_{T,v}}$ is induced by a homomorphism between algebraic tori $$R_{E_{T,v}/\bb{Q}_p}\bb{G}_m\xrightarrow{R_{E_{T,v}/\bb{Q}_p}\mu}R_{E_{T,v}/\bb{Q}_p}T_{E_{T,v}}\xrightarrow{\mrm{Norm}_{E_{T,v}/\bb{Q}_p}}T_{\bb{Q}_p},$$
which sends the Iwahori subgroup to the Iwahori subgroup by \cite[Lem. 2.5.19]{KP23}. This means that $\mrm{Gal}(\overline{E}_T/E_T^p)$ acts trivially on $\sh_{K_T}(T,\ca{Y})_{\overline{E}_T}$. Hence, we have the desired statement for both $\sh_{K_T}$ and $\sh_{K}(G,\hbar(\ca{X}))$ and finally for $\sh_{K}(G,\ca{X})$.
\end{proof}
\subsubsection{}We come back to the setup that $(G,X)$ is a Shimura datum. Let $K=K_pK^p$ be an open compact subgroup of $G(\A)$ such that $K_p\sbst G(\bb{Q}_p)$ and $K^p$ is neat open compact in $G(\Ap)$.
In this subsection, we assume one of the following assumptions: 
\begin{enumerate}
\item $G$ is quasi-split and unramified at $p$, and $K_p=G_{\zbkp}(\bb{Z}_p)$ is a hyperspecial subgroup associated with some smooth reductive model $G_{\zbkp}$ of $G$ over $\zbkp$. 
\item $K_p=\mathscr{G}^\circ(\bb{Z}_p)$ is a parahoric subgroup of $G(\bb{Q}_p)$.
\end{enumerate}
\begin{prop}\label{comp-phi}With one of the assumptions above,
for any cusp label representative $\Phi$, the geometrically connected components of $\sh_{K_\Phi,h}$, $\overline{\sh}_{K_\Phi}$ and $\sh_{K_\Phi}$ are defined over $E^p$.
\end{prop}
\begin{proof}
It suffices to show that the geometrically connected components of $\sh_{K_\Phi,h}$ are defined over $E^p$. In the first case, $K_{\Phi,h}$ is also hyperspecial; in the second case, $K_{\Phi,h}$ is quasi-parahoric by \cite[Prop. 2.67]{Mao25}. Hence, the statement follows from Proposition \ref{prop-par-geo-conn} above. In the first case, we can also replace the reference \cite[Cor. 4.3.9]{KP15} in the proof of the proposition above with \cite[Prop. 2.2.4]{Kis10}.
\end{proof}
\subsection{$ZP$-cusps}\label{subsec-zp-cusp}
We explain how to construct some cusp labels that are \emph{not} irreducible (see \S\ref{subsubsec-mixedshimura} for the definition), which will be helpful to our purpose. We do not assume that the Shimura data being considered are of abelian type until \S\ref{sss-conclusion}.
\subsubsection{}
Let $Z:=Z_G$ be the center of $G$ as before. We denote by $\gls{zpp}$ the identity component of the group generated by $Z$ and $P_\Phi$, for any cusp label representative $\Phi$. Note that $ZP_\Phi$ is also a normal subgroup of $Q_\Phi$ over $\bb{Q}$. The group $\zpp(\bb{A}_f)$ acts on cusp label representatives by sending $\Phi=(Q_\Phi,X^+_\Phi,g_\Phi)$ to $q^\prime\Phi:=(Q_\Phi,X^+_\Phi,q^\prime g_\Phi)$ for any $q^\prime\in \zpp(\bb{A}_f)$. 
\begin{definition}
We say an admissible cone decomposition (see Convention \ref{conv-cone}) $\Sigma$ is \textbf{$ZP$-invariant} if $\Sigma(q^\prime\Phi)=\Sigma(\Phi)$ for any cusp label representative $\Phi$ and for any $q^\prime\in \zpp(\bb{A}_f)$.
\end{definition}
The following lemma is purely combinatorial and is similar to \cite[Prop. 6.3.3.5]{Lan13}:
\begin{lem}\label{lem-zinv-ref}
Let $\Sigma$ be an admissible cone decomposition with the property $\ca{P}$, where $\ca{P}=$ $\emptyset$, ``smooth'', ``projective'' or ``smooth and projective''. Then there is a refinement $\Sigma^\prime$ of $\Sigma$ which is $ZP$-invariant and has the property $\ca{P}$.
\end{lem}
\begin{definition}For any cusp label representative $\Phi$, we define the \textbf{magnitude} of $\Phi$ to be the dimension of the torus $\mbf{E}_{K_\Phi}$, denoted by $\gls{mPhi}$.
\end{definition}
\begin{definition}\label{equi-z}
There is a partial order $\preceq$ associated with the set of cusp label representatives forming a poset $(\{\Phi\},\preceq)$ defined by the following: $\Phi_1\preceq\Phi_2$ if there are $q^\prime\in ZP_{\Phi_2}(\bb{A}_f)$ and $\gamma\in G(\bb{Q})$ such that
$\gamma P_{\Phi_1}\gamma^{-1}\sbst P_{\Phi_2}$, $P_{\Phi_1}(\bb{Q})\gamma X^+_{\Phi_1}=P_{\Phi_2}(\bb{Q})X^+_{\Phi_2}$ and $\gamma g_{\Phi_1}\equiv q^\prime g_{\Phi_2}$ modulo $K$, denoted also by $\Phi_1\xrightarrow{(\gamma,q^\prime)_K}\Phi_2$. We write $\Phi_1\sim_{ZP}\Phi_2$ if $\Phi_1\preceq \Phi_2$ and $\Phi_2\preceq\Phi_1$. Then ``$\gls{simZP}$'' defines an equivalence relation. We denote by $\gls{CuspZPKGX}$ the set of equivalence classes $\mrm{Cusp}_K(G,X)/\sim_{ZP},$ called the set of \textbf{$ZP$-cusp labels}.\par
Similarly, we can define the set of \textbf{$ZP$-cusp labels with cones} $\gls{CuspZPKGXSigma}$ to be the set of $ZP$-equivalence classes of $\mrm{Cusp}_K(G,X,\Sigma)$. Indeed, we fix any pairs $(\Phi_1,\sigma_1)$ and $(\Phi_2,\sigma_2)$ such that $\sigma_i\in \Sigma(\Phi_i)$ and $\sigma_i^\circ\sbst \mbf{P}_{\Phi_i}^+$, for $i=1$ and $2$. We define $(\Phi_1,\sigma_1)\sim_{ZP}(\Phi_2,\sigma_2)$, if there are $\gamma\in G(\bb{Q})$ and $q'\in ZP_{\Phi_2}(\A)$, such that $\Phi_1\xrightarrow[\sim]{(\gamma,q')_K}\Phi_2$ and $\mrm{int}(\gamma)\sigma_1=\sigma_2$. Then we can define $\cusp_K^{ZP}(G,X,\Sigma):=\cusp_K(G,X,\Sigma)/\sim_{ZP}$. There is also a partial order on $\cusp_K^{ZP}(G,X,\Sigma)$ defined similarly: $[(\Phi_1,\sigma_1)]\preceq_{ZP}[(\Phi_2,\sigma_2)]$ if and only if  $\Phi_1\xrightarrow{(\gamma,q')_K}\Phi_2$ and $\gamma^{-1}\sigma_2\gamma$ is a face of $\sigma_1$.
\end{definition}
\begin{lem}\label{izpq}
Let $Q$ be an admissible $\bb{Q}$-parabolic subgroup of $G$, and let $X^+$ be any connected component of $X$. Let $I^{G,X^+,ZP}(Q):=\{[\Phi]\in\mrm{Cusp}_K^{ZP}(G,X)|\Phi\sim_{ZP} (Q,X^+,g_\Phi),g_\Phi\in G(\bb{A}_f)\}$. Then $I^{G,X^+,ZP}(Q)\iso\mrm{Stab}_{Q(\bb{Q})}(D_{Q,X^+})ZP_Q(\bb{A}_f)\backslash G(\bb{A}_f)/K$.\end{lem}
We will omit the superscript $G$ and $X^+$, and write $I^{G,X^+,ZP}(Q)$ as $I^{ZP}(Q)$, if $G$ and $X^+$ are clear in the context.

\begin{proofof}[Lemma \ref{izpq}]
This follows verbatim from the proof of Proposition \ref{prop-iphi}, after replacing $P_Q(\A)$ there with $ZP_Q(\A)$ here and replacing $\mrm{Cusp}_K(G,X)$ there with $\mrm{Cusp}_K^{ZP}(G,X)$ here.
\end{proofof}
Note that \cite[Prop. 2.14(a)]{Pin89} can be strengthened naively. 
\begin{lem}[{\cite[Prop. 2.14(a)]{Pin89}}]\label{lem-pin214-str}
The conjugation of $ZP_\Phi$ on $U_\Phi$ is a scalar multiplication: $c_\Phi: ZP_\Phi\to \bb{G}_m$.
\end{lem}
\begin{proof}
Since $Z$ acts on $U_\Phi$ trivially, since the quotient of $ZP_\Phi$ by $Z\cap ZP_\Phi$ is $P_\Phi$, and since $(P_\Phi,D_\Phi)$ is irreducible, we can reduce our statement to \emph{loc. cit.} by taking the quotient of $Z\cap ZP_\Phi$.
\end{proof}
\begin{proofof}[Lemma \ref{lem-zinv-ref}]
Note that we only need to consider one cusp label representative $\Phi$ in each equivalence class of the relation $\sim_{ZP}$: In fact, $K'_{\Phi,U}$ is a scalar multiplication of $q^\prime K'_{\Phi,U}(q^\prime)^{-1}$ by Lemma \ref{lem-pin214-str} or the proof of \cite[Prop. 2.14 (a)]{Pin89} for $q^\prime\in \zpp(\bb{A}_f)$. Suppose that $q'K'_{\Phi,U}(q')^{-1}=c_\Phi(q')K'_{\Phi,U}$ for some $c_\Phi(q')\in \A^\times$. Since we have the decomposition $\A^\times=\bb{Q}^\times\cdot \wat{\bb{Z}}^\times$ and since $U_\Phi$ is connected and unipotent, $\mbf{\Lambda}_{K_\Phi}$ is a $\bb{Q}^\times$-scalar multiplication of $\mbf{\Lambda}_{K_{q'\Phi}}$. Therefore, $\Sigma(\Phi)$ has property $\ca{P}$ with respect to $\mbf{\Lambda}_{K_\Phi}$ if and only if $\Sigma(\Phi^\prime)$ has the property $\ca{P}$ with respect to $\mbf{\Lambda}_{K_{\Phi'}}$ for any cusp label repressentatives $\Phi$ and $\Phi'$ such that $\Phi\sim_{ZP}\Phi^\prime$. \par 
Next, we fix a complete set of representatives $\ca{R}$ in $\ca{CLR}(G,X)$ of the equivalence classes $\cusp_K^{ZP}(G,X)$ of $ZP$-cusp labels. We prove the statement by induction on the magnitude of cusp label representatives $\Phi\in\ca{R}$. If $m(\Phi)=0$, then $\Phi=(G,X^+,g_\Phi)$ and there is nothing to show since $\mbf{P}_\Phi$ is trivial.
Assume that we have refined $\Sigma$ by some cone decompositon $\Sigma_n$ for cusp label representatives $\Phi$ of magnitude $m(\Phi)\leq n$ with the property $\ca{P}$ such that the collection $\{\Sigma_n(\Phi)\}_{m(\Phi)\leq n}$ is invariant under the left action of $\zpp(\bb{A}_f)$. 
For any two cusp label representatives $\Phi_1 \xrightarrow{(\gamma,q)_K}\Phi_2$ and $q^\prime\in ZP_{\Phi_2}(\bb{A}_f)$, $q^\prime\Phi_1\xrightarrow{(\gamma,q^\prime q(q^\prime)^{-1})_K}q^\prime\Phi_2$. So for any cusp label representative $\Phi\in\ca{R}$ of magnitude $m(\Phi)=n+1$, the cone decomposition $\Sigma_n$ defines equal cone decompositions for $\mbf{P}_{q^\prime\Phi}-\mbf{P}_{q^\prime\Phi}^+$ for any $q^\prime\in \zpp(\bb{A}_f)$. 
There is an admissible cone decomposition $\Sigma_{n+1}(\Phi)$ for $\mbf{P}_{\Phi}$ with the property $\ca{P}$ \emph{extending} the cone decomposition for $\mbf{P}_\Phi-\mbf{P}^+_\Phi$ defined by $\Sigma_{n}$, and we can define cone decompositions $\Sigma(q^\prime\Phi):=\Sigma(\Phi)$ for all $q^\prime\in \zpp(\bb{A}_f)$. 
Therefore, we can extend the definition of $\Sigma_{n+1}$ to all cusp labels of magnitudes $\leq n+1$, then we can get the desired refinement $\Sigma^\prime$ of $\Sigma$ by induction. For the results of extendability of admissible cone decompositions for $\mbf{P}_\Phi-\mbf{P}^+_{\Phi}$ to $\mbf{P}_\Phi$ preserving the property $\ca{P}$, see, e.g., \cite[pp. 32-35]{KKMS73} and \cite[5.20-5.25 and 9.20]{Pin89}.
\end{proofof}
\begin{rk}\label{zp-orb} We can see from Lemma \ref{izpq} that a $ZP$-cusp label with cone is determined by a $ZP$-cusp label $[\Phi]$ and a $\gls{DeltaPhiZP}$-orbit $[\sigma]$ in $U_\Phi(\bb{R})$, where $\Delta_{\Phi,K}^{ZP}:=\stb_{Q(\bb{Q})}(D_\Phi)\cap ZP(\A)g_\Phi Kg_\Phi^{-1}/ZP_\Phi(\bb{Q})$. We will omit the subscript $K$ in $\Delta^{ZP}_{\Phi,K}$ if $K$ is clear in the context.
\end{rk}
\begin{definition}
Fix an admissible cone decomposition $\Sigma$. Let $Q$ be an admissible $\bb{Q}$-parabolic subgroup of $G$, and let $X^+$ be any connected component of $X$. Define $I^{G,X^+,ZP}(Q,\Sigma):=\{[(\Phi,\sigma)]\in\cusp_K^{ZP}(G,X,\Sigma)|(\Phi,\sigma)\sim_{ZP}(Q,X^+,g_\Phi;\sigma'),g_\Phi\in G(\A),\sigma'\in\Sigma(Q,X^+,g_\Phi)\}$, and similarly define $I^{G,X^+}(Q,\Sigma):=\{[(\Phi,\sigma)]\in\cusp_K(G,X,\Sigma)|(\Phi,\sigma)\sim(Q,X^+,g_\Phi;\sigma'),g_\Phi\in G(\A),\sigma'\in\Sigma(Q,X^+,g_\Phi)\}$.\par
We shall omit the superscript ``$G,X^+$'' if it is clear in the context.
\end{definition}
We will use the following propositions later:
\begin{prop}\label{bij}
Let $\Phi$ be a cusp label representative, and let $\sigma\in\Sigma(\Phi)$. Define 
$$\gls{ZPPhi}:=\{[\Phi^\prime]\in\mrm{Cusp}_K(G,X)|\Phi^\prime\sim_{ZP}\Phi\},$$
and define 
$$\gls{ZPPhisigma}:=\{[(\Phi^\prime,\sigma^\prime)]\in\mrm{Cusp}_K(G,X,\Sigma)|(\Phi^\prime,\sigma^\prime)\sim_{ZP}(\Phi,\sigma)\}.$$
Then $[ZP(\Phi)]$ is bijective to $\wdtd{\Delta}^{ZP}_{\Phi,K}\bss( P_\Phi(\A)\bss ZP_\Phi(\A)/\wdtd{K}_\Phi)$ and $[ZP(\Phi,\sigma)]$ is bijective to $$\wdtd{\Delta}^{ZP,\circ}_{\Phi,K}\bss (P_\Phi(\A)\bss ZP_\Phi(\A)/\wdtd{K}_\Phi),$$ where $\wdtd{K}_\Phi:=ZP(\A)\cap g_\Phi K g_\Phi^{-1}$, $\wdtd{\Delta}_{\Phi,K}^{ZP}:=\stb_{Q_\Phi(\bb{Q})}(D_\Phi)\cap ZP_\Phi(\A)g_\Phi Kg_\Phi^{-1}$ and $\wdtd{\Delta}_{\Phi,K}^{ZP,\circ}:=\stb_{Q_\Phi(\bb{Q})}(D_\Phi,\sigma)\cap ZP_\Phi(\A)g_\Phi Kg_\Phi^{-1}.$
\end{prop}
\begin{proof}
By definition, $[ZP(\Phi)]$ is the fiber of $I(Q_\Phi)\to I^{ZP}(Q_\Phi)$ at $[\Phi]$. This fiber is bijective to the fiber of 
$$\stb_{Q(\bb{Q})}(D_\Phi) P_\Phi(\A)\bss G(\A)/K\to \stb_{Q(\bb{Q})}(D_\Phi) ZP_\Phi(\A)\bss G(\A)/K$$ at $[g_\Phi]$, which is $\wdtd{\Delta}^{ZP}_{\Phi,K}\bss( P_\Phi(\A)\bss ZP_\Phi(\A)g_\Phi K/K)\iso \wdtd{\Delta}_{\Phi,K}^{ZP}\bss(P_\Phi(\A)\bss ZP_\Phi(\A)/K_\Phi)$. The other statement can also be proved in a similar way. The set $[ZP(\Phi,\sigma)]$ is the fiber of 
$$\stb_{Q(\bb{Q})}(D_\Phi) P_\Phi(\A)\bss \Sigma(\Phi)\times G(\A)/K\to \stb_{Q(\bb{Q})}(D_\Phi) ZP_\Phi(\A)\bss \Sigma(\Phi)\times G(\A)/K$$
at $(\sigma,g_\Phi)$, such that $(Q_\Phi,X^+_\Phi,g_\Phi;\sigma)$ represents $[(\Phi,\sigma)]$. Since the stabilizer in $\stb_{Q(\bb{Q})}(D_\Phi)$ of $P_\Phi(\A)\bss \sigma\times ZP_\Phi(\A)g_\Phi K/K$ is $\wdtd{\Delta}^{ZP,\circ}_{\Phi,K}$, $[ZP(\Phi,\sigma)]$ is bijective to $\wdtd{\Delta}^{ZP,\circ}_{\Phi,K}\bss( P_\Phi(\A)\bss ZP_\Phi(\A)/K_\Phi)$. 
\end{proof}
Let $\gls{DeltaZPcircPhiK}:=\wdtd{\Delta}^{ZP,\circ}_{\Phi,K}/ZP_\Phi(\bb{Q})$. 
Again, we shall omit the subscript $K$ in $\wdtd{\Delta}^{ZP}_{\Phi,K}$, $\wdtd{\Delta}^{ZP,\circ}_{\Phi,K}$ and $\Delta^{ZP,\circ}_{\Phi,K}$ if it is clear in the context.
\begin{prop}\label{fiber}For any $[\Phi^\prime]\in[ZP(\Phi)]$, there is a cone $\sigma^\prime\in\Sigma(\Phi^\prime)$ such that $[(\Phi^\prime,\sigma^\prime)]\in[ZP(\Phi,\sigma)]$.  
The fiber of the natural projection $\ell: [ZP(\Phi,\sigma)]\to [ZP(\Phi)]$ at $[\Phi^\prime]\in[ZP(\Phi)]$ is bijective to
$\Delta_{\Phi^\prime}$-orbits of the $\Delta^{ZP}_{\phip}$-orbit $[\sigma^\prime]_{ZP}$. Moreover, we can choose a representative $\Phi^\dpr\sim\Phi^\prime$ such that $\Phi^\dpr=(Q_\Phi,X^+_\Phi,g_{\Phi^\dpr}g_\Phi)$ for some $g_{\Phi^\dpr}$ lifting an element in $\wdtd{\Delta}^{ZP}_\Phi \bss(P_\Phi(\A)\bss ZP_\Phi(\A)/\wdtd{K}_\Phi)$, then this fiber is bijective to $\Delta_\Phi$-orbits of the $\Delta^{ZP}_\Phi$-orbit $[\sigma]_{ZP}$.
\end{prop}
\begin{proof}
If $\Phi\bd{\gamma}{q}{K}{\sim}\Phi^\prime$ for $\gamma\in G(\bb{Q})$ and for $q\in ZP_\Phi(\A)$, then we let $\sigma^\prime:=\gamma\sigma$ and $(\Phi,\sigma)\bd{\gamma}{q}{K}{\sim}(\Phi^\prime,\sigma^\prime)$. (Here we are using an abbreviation: $\gamma\sigma=\mrm{int}(\gamma)\sigma$, and the same for the paragraph below.)\par 
For the second statement, we consider the commutative diagram
\begin{equation}
\begin{tikzcd}
I(Q_\Phi,\Sigma)\arrow{rr}\arrow{d}&& I^{ZP}(Q_\Phi,\Sigma)\arrow{d}\\
I(Q_\Phi)\arrow{rr}&& I^{ZP}(Q_\Phi)
\end{tikzcd}
\end{equation}
Suppose that $[\Phi^\prime]\in I(Q_\Phi)$ maps to $[\Phi]\in I^{ZP}(Q_\Phi)$. The fiber of $[\Phi']$ along $I(Q_\Phi,\Sigma)\to I(Q_\Phi)$ is $\{([\Phi'],[\sigma^\dpr])|[\sigma^\dpr]\in \Delta_{\Phi'}\bss\Sigma(\Phi')\}$. We still assume that $\Phi\bd{\gamma}{q}{K}{\sim}\Phi^\prime$. Then $([\Phi^\prime],[\sigma^\dpr])$ maps to $([\Phi],[\gamma^{-1}\sigma^\dpr]_{ZP})$ via $I(Q_\Phi,\Sigma)\to I^{ZP}(Q_\Phi,\Sigma)$, where $[\gamma^{-1}\sigma^\dpr]_{ZP}$ is the $\Delta^{ZP}_\Phi$-orbit of $\gamma^{-1}\sigma^\dpr$. Then the fiber of $\ell$ at $[\Phi^\prime]$ is isomorphic to the elements of $\Delta_{\Phi^\prime}\bss\Sigma(\Phi^\prime)$ that map to $[\sigma]_{ZP}$ via \begin{equation}\begin{split}\Delta_{\Phi^\prime}\bss\Sigma(\Phi^\prime)\lra\Delta_{\Phi}^{ZP}\bss\Sigma(\Phi^\prime)\\
[\sigma^\dpr]\mapsto[\gamma^{-1}\sigma^\dpr]_{ZP},\end{split}\end{equation}
which is the pullback of $[\sigma]_{ZP}$ along the map above. So it is a $\gamma\Delta^{ZP}_{\Phi}\gamma^{-1}$-orbit of $\sigma^\prime=\gamma\sigma$ modulo the action of $\Delta_{\Phi^\prime}$. Since $\gamma\Delta^{ZP}_\Phi\gamma^{-1}=\Delta^{ZP}_{\Phi^\prime}$, the fiber of $\ell$ at $[\Phi^\prime]$ is bijective to  $\Delta_{\Phi^\prime}$-orbits of $[\sigma^\prime]_{ZP}$.\par
The last statement follows from Proposition \ref{bij}.
\end{proof}
\subsubsection{}A $ZP$-invariant cone decomposition $\Sigma$ can extend over certain central extensions. We have the following variant of Lemma \ref{lem-emb-mix}:
\begin{prop}\label{ext-imp}
Let $G_1$ be a connected reductive group over $\bb{Q}$ with $Z_1$ its center. Let $Z_2$ be a $\bb{Q}$-torus containing $Z_1$ and let $G_2:=G_1\times^{Z_1}Z_2$ be a connected reductive group with center $Z_2$. Denote by $\iota: G_1\to G_2$ the embedding induced by the inclusion $Z_1\hookrightarrow Z_2$.
We abusively denote by $\iota:(G_1,X_1)\hookrightarrow (G_2,X_2)$ the embedding of Shimura data induced by $\iota: G_1\to G_2$. 
Let $K_1$ be an open compact subgroup of $G_1(\A)$ of the form $K_1=K_{1,p}K_1^p$ where 
$K_{1,p}$ is open compact and $K_1^p\sbst G_1(\Ap)$ is neat open compact. Then: 
\begin{enumerate}
\item There is a $\bb{Q}$-subgroup $Z_3$ of $Z_2$ of multiplicative type such that $Z_2=Z_1\cdot Z_3$ and $Z_1\cap Z_3$ is finite.
\item Let $K_2:=K_{2,p}K^p_2$. There is a neat open compact subgroup $K_2^{p,\prime}\sbst G_2(\Ap)$ such that, for any neat open compact $K_2^p\sbst K_2^{p,\prime}$ and $K_2\cap G_1(\A)=K_1$, the following statements are true: $\iota$ induces an open and closed embedding of Shimura varieties $$\iota:\sh_{K_1}(G_1,X_1)\hookrightarrow\sh_{K_2}(G_2,X_2);$$
moreover, let $\Phi_1$ be any cusp label representative of $(G_1,X_1)$ which maps to a cusp label representative $\Phi_2$ of $(G_2,X_2)$ via $\iota$, we have isomorphisms 
$U_{\Phi_1}\iso U_{\Phi_2}$, $\mbf{\Lambda}_{K_{\Phi_1}}\iso \mbf{\Lambda}_{K_{\Phi_2}}$, and $$\iota_{\Phi_1}:\sh_{K_{\Phi_1}}(P_{\Phi_1},D_{\Phi_1})\to\sh_{K_{\Phi_2}}(P_{\Phi_2},D_{\Phi_2}).$$ 
\item If $Z_3^\circ$, the identity component of $Z_3$, is isogenous to a product of $\bb{Q}$-split tori and compact-type tori, then the results in the last part hold for any neat open compact $K_2$. 
\end{enumerate}
\end{prop}
\begin{proof}
In fact, there is a $\bb{Q}$-subtori $T_3$ of $\overline{Z}_2:=Z_2/Z_{G_1^\der}$ such that the intersection of $\overline{Z}_1:=Z_1/Z_{G_1^\der}$ and $T_3$ is finite, and such that $\overline{Z}_2=\overline{Z_1}\cdot T_3$. We then find such a $Z_3$ by taking $Z_3:=Z_2\times_{\overline{Z}_2}T_3$.\par
The first statement in the second part follows directly from \cite[Lem. 2.1.2]{Kis10}. 
Since $G_1$ is a normal subgroup of $G_2$ and $G_1^\der=G_2^\der$, $P_{\Phi_1}\iso P_{\Phi_2}$, $U_{\Phi_1}\iso U_{\Phi_2}$ and $D_{\Phi_1}\iso D_{\Phi_2}$ by definition. Since $K_{\Phi_1}= P_{\Phi_1}(\A)\cap g_{\Phi_1}K_1 g_{\Phi_1}^{-1}=P_{\Phi_2}(\A)\cap g_{\Phi_1}K_2 g_{\Phi_1}^{-1}\cap G_1(\A)=P_{\Phi_2}(\A)\cap g_{\Phi_1}K_2 g_{\Phi_1}^{-1}=K_{\Phi_2}$, we have isomorphisms between boundary mixed Shimura varieties. Moreover, $\mbf{\Lambda}_{K_{\Phi_1}}\sbst \mbf{\Lambda}_{K_{\Phi_2}}$ is of finite index, we can choose $K_2$ such that $K_2/K_1$ sufficiently small and $\mbf{\Lambda}_{K_{\Phi_1}}= \mbf{\Lambda}_{K_{\Phi_2}}$. Note that this is possible because the number of cusp labels of $\sh_{K_1}$ is finite and remains unchanged if we shrink $K_2$. Since we only need to choose $K_2$ to be sufficiently small but containing $K_1$, we can choose a $K_2$ of the form $K_2=K_{2,p}K_2^p$.\par
When $Z_3^\circ$ is isogenous to a product of $\bb{Q}$-split tori and compact-type tori, and $K_2$ is neat, $\Delta_{G_2/G_1,K_2}$ itself is trivial.
\end{proof}
\begin{prop}\label{ext-zp}With the conventions in Lemma \ref{lem-zinv-ref} and Proposition \ref{ext-imp}, we suppose that $K_2$ is chosen as in Proposition \ref{ext-imp}. Let $\Sigma_1$ be an admissible cone decomposition for $\sh_{K_1}(G_1,X_1)$. Then there is an admissible $ZP$-invariant cone decomposition $\Sigma_2$ for $\sh_{K_2}(G_2,X_2)$ with property $\ca{P}$ such that $\iota^*(\Sigma_2)$ is a refinement of $\Sigma_1$. Moreover, $\iota^*(\Sigma_2)$ also has the property $\ca{P}$.
\end{prop}
\begin{proof}It suffices to find a $\Sigma_2$ for $\sh_{K_2}(G_2,X_2)$ such that $\iota^*(\Sigma_2)$ refines $\Sigma_1$, because we can refine $\Sigma_2$ further to obtain a cone decomposition $\Sigma_2'$ with property $\ca{P}$ by Lemma \ref{lem-zinv-ref}. The induced cone decomposition $\iota^*(\Sigma_2)$ will also have the property $\ca{P}$ since $U_{\Phi_1}\iso U_{\Phi_2}$ and $\mbf{\Lambda}_{K_{\Phi_1}}\iso \mbf{\Lambda}_{K_{\Phi_2}}$ from the choice of $K_2$ above.
Let $Q_1$ be an admissible parabolic subgroup of $G_1$, then there is a unique admissible parabolic subgroup $Q_2$ of $G_2$ containing $Q_1$.
Since $Z_2\cap G_1=Z_1$, $Q_1=G_1\cap Q_2$ and $P_{Q_1}=P_{Q_2}\cap G_1$, and since our choice of $K_2$ as in Proposition \ref{ext-imp}, we see that $I(Q_1)\hookrightarrow I(Q_2)$. 
We then do induction on the magnitude of cusp label representatives $\Phi_2$ of $(G_2,X_2)$. 
Fix any complete set of representatives $\ca{R}_1$ of $\mrm{Cusp}_{K_1}(G_1,X_1)$, we can extend $\ca{R}_1$ to a complete set of representatives $\ca{R}_2$ of $\mrm{Cusp}_{K_2}(G_2,X_2)$ since $I(Q_1)\hookrightarrow I(Q_2)$. 
Assume that we have an admissible cone decomposition $\Sigma_2^n$ with property $\ca{P}$ for cusp labels of $(G_2,X_2)$ with magnitude $\leq n$. 
Let $\Phi_2$ be any cusp label representative in $\ca{R}_2$ of magnitude $n+1$. Then $\Phi_2$ represents an element $[g_{\Phi_2}]\in I(Q_{\Phi_2})$. For any $\Phi_2^\prime\succeq \Phi_2$, there are $\gamma\in G(\bb{Q})$ and $q^\prime\in P_{\Phi^\prime_2}$ such that $\Phi_2\xrightarrow{(\gamma,q^\prime)_{K_2}}\Phi_2^\prime$. Let $Q_{\Phi_2^\prime}^\gamma:=\gamma^{-1}Q_{\Phi_2^\prime}\gamma$. Then $\Phi_2^\prime$ represents the image of $[g_{\Phi_2}]$ under $I(Q_{\Phi_2})\to I(Q^\gamma_{\Phi_2^\prime})$ induced by $Q_{\Phi_2}\hookrightarrow Q_{\Phi_2^\prime}^\gamma$. 
If $\Phi_2\in\ca{R}_1$, then the image of $[g_{\Phi_2}]$ in $I(Q^\gamma_{\Phi_2^\prime})$ lies in $I(Q^\gamma_{\Phi_2^\prime}\cap G_1)$, the cone decomposition for $\mbf{P}_{\Phi_2}-\mbf{P}_{\Phi_2}^+$ refines $\Sigma_1$ by induction hypothesis, and we can choose an admissible cone decomposition for $\Phi_2$ that refines $\Sigma_1(\Phi_2)$ and extends $\Sigma_2^n|_{\mbf{P}_{\Phi_2}-\mbf{P}_{\Phi_2}^+}$; if $\Phi_2$ is not in $\ca{R}_1$, we can find an admissible cone decomposition for $\mbf{P}_{\Phi_2}$ extending that of $\mbf{P}_{\Phi_2}-\mbf{P}_{\Phi_2}^+$. Then we get desired result by induction process.
\end{proof}
Let $\iota:(G_1,X_1)\to(G_2,X_2)$ be any morphism between Shimura data. For any $g\in G_2(\A)$ and $K_1\sbst K_2$, there is a morphism
\begin{equation}
\iota(g): \sh_{gK_1g^{-1}}(G_1,X_1)\to \sh_{K_2}(G_2,X_2)
\end{equation}
defined over the reflex field of $(G_1,X_1)$. Over complex points, it is described by $[(x,g_1)]\mapsto [(x,g_1g)]$.
For any admissible cone decomposition $\Sigma_2$ for $\sh_{K_2}(G_2,X_2)$, denote by $g^*\Sigma_2$ the induced cone decomposition for $\sh_{gK_1g^{-1}}(G_1,X_1)$ by pulling back via $\iota(g)$.
\begin{cor}\label{ext-cpt-imp} Let $\iota:(G_1,X_1)\hookrightarrow (G_2,X_2)$ be the embedding as in Proposition \ref{ext-imp}. Let $K_1$ and $K_2$ be the neat open compact subgroups chosen there. 
Let $\Sigma_2$ be the admissible $ZP$-invariant cone decomposition chosen as in Proposition \ref{ext-zp}. Choose a set of representatives $\{g_i\}_{i\in \ca{T}}$ of $\ca{T}:=G_2(\bb{Q})G_1(\A)\bss G_2(\A)/K_2$. As the last paragraph, there are morphisms between Shimura varieties 
\begin{equation}
\iota(g_i): \sh_{g_iK_1g_i^{-1}}(G_1,X_1)\to \sh_{K_2}(G_2,X_2),    
\end{equation}
defined by $[(x,g)]\mapsto [(x,gg_i)]$ over complex points,
which induce morphisms between minimal and 
toroidal compactifications
\begin{equation}
   \disju_{g_i\in\ca{T}}\iota(g_i)^{\mmin}: \disju_{g_i\in\ca{T}}\sh_{g_iK_1g_i^{-1}}^\mrm{min}\to \sh_{K_2}^\mrm{min}
\end{equation}
and
\begin{equation}
    \disju_{g_i\in\ca{T}}\iota(g_i)^{\mrm{tor}}:\disju_{g_i\in\ca{T}}\sh_{g_iK_1g_i^{-1}}^{g_i^*\Sigma_2}\lra \sh_{K_2}^{\Sigma_2}.
\end{equation}
In the morphism above, $g_i^*:=\iota(g_i)^{*}$. Moreover, the morphisms restricted to $\iota^\mmin=\iota(1)^\mmin$ and $\iota^{\mrm{tor}}=\iota(1)^{\mrm{tor}}$, $\iota^\mmin:\sh^{\mrm{min}}_{K_1}\to \sh_{K_2}^\mrm{min}$ and 
$\iota^{\mrm{tor}}:\sh^{\iota^*\Sigma_2}_{K_1}\to \sh^{\Sigma_2}_{K_2}$, are open and closed embeddings over the reflex field of $(G_1,X_1)$.
\end{cor}
\begin{proof}We only need to show the last claim on open and closed embeddings. Define $\Delta_{G_2/G_1,K_2}:=G_2(\bb{Q})\cap G_1(\A)K_2/G_1(\bb{Q})$. 
From \cite[Prop. 7.10]{Pin89} for trivial cone decompositions and the choice of $K_2$ in Proposition \ref{ext-imp}, the group $\Delta_{G_2/G_1, K_2}$ acts on $\sh_{K_1}$
trivially. Since $\sh_{K_1}$ is dense in its minimal and toroidal compactifications, $\Delta_{G_2/G_1,K_2}$ acts on the minimal and toroidal compactifications of $\sh_{K_1}$ trivially. From \cite[Prop. 7.10]{Pin89}, we see the statement is true for toroidal compactifications. (This can also imply that $\iota=\iota(1)$ is an open and closed embedding.) We can check over $\bb{C}$ to see that the schematic closure $(\sh_{K_1})^{\overline{\ }}$ of $\sh_{K_1}$ in $\sh_{K_2}^\mmin$ under $\iota$ is normal and proper. Since $\iota^\mmin$ is proper, quasi-finite, and factors through $(\sh_{K_1})^{\overline{\ }}$, we deduce the claim for minimal compactifications by Zariski's main theorem.
\end{proof}
\subsubsection{}\label{a,b}
Let $G$ be a connected reductive group over $\bb{Q}$ with connected center $Z$. Fix $(G^\ad,X^\ad)$, a Shimura datum associated with its adjoint group. Throughout \S\ref{a,b}, we work under the following assumptions:
\begin{itemize}
\item Suppose that $(G^\ad,X^\ad)$ lifts to two different Shimura data $(G,X_a)$ and $(G,X_b)$ such that the natural quotient from $G$ to $G^\ad$ induces isomorphisms $(G^\ad,X_a^\ad)\xrightarrow{\sim}(G^\ad,X^\ad)\xleftarrow{\sim}(G^\ad,X_b^\ad).$ Note that the maps $X_i\to X_i^\ad$ are injective for $i=a$ and $b$. 
\item Suppose that for one (and therefore any) $x\in X^\ad$, $h_x:\bb{S}\to G^\ad_\bb{R}$ lifts to unique cocharacters $h^i_x:\bb{S}\to G_\bb{R}$ in $X_i$. 
\end{itemize}
Then $c:=h_x^a\cdot (h_x^{b})^{-1}$ is a homomorphism that factors through the center $Z_\bb{R}$ of $G_\bb{R}$, and therefore it is independent of the choice of $x\in X^\ad$.\par
We fix homomorphisms $h_0$ and $h_\infty$ as in \S\ref{subsubsec-boundary}. Since $(H_0/H_0^\der)_\bb{C}\iso \bb{S}_\bb{C}$, we can define $u^c: H_{0,\bb{C}}\to G_\bb{C}$ by $H_{0,\bb{C}}\twoheadrightarrow(H_0/H_0^\der)_\bb{C}
\xrightarrow{c_\bb{C}} Z_\bb{C}\hookrightarrow G_\bb{C}$. For $i=a$ and $b$, and for any admissible $\bb{Q}$-parabolic subgroup $Q$, let $u^{Q,i}_x:H_{0,\bb{C}}\to G_\bb{C}$ be the unique homomorphisms such that $u^{Q,i}_x\circ h_0=h_x^i$, which satisfy \cite[Prop. 4.6 (b)]{Pin89}. 
\begin{lem}\label{uq-mul}
$u^{Q,a}_x=u^{Q,b}_x\cdot u^c$.
\end{lem}
Note that the multiplication in the equation above makes sense because $u^c$ factors through $Z_\bb{C}$.\par
\begin{proofof}[Lemma \ref{uq-mul}]
The statement is proved by verifying conditions in \cite[Prop. 4.6 (b)]{Pin89}. Let $v^{Q,a}_x:=u^{Q,b}_x\cdot u^c$. The homomorphism $\pi^\prime\circ v^{Q,a}_x: H_{0,\bb{C}}\to Q_\bb{C}\twoheadrightarrow(Q/U_Q)_\bb{C}$ is defined over $\bb{R}$ since $u^c$ is defined over $\bb{R}$. Since $\omega^{Q,b}_{x,\infty}$ is of the form $\omega_x^b\cdot \lambda$ for some cocharacter $\lambda$ depending only on $Q$ as \cite[Prop. 4.6 (b)iii]{Pin89} and $u^c$ factors through $Z_\bb{C}$, we have $\omega^{Q,b}_{x,\infty}\cdot (u^c\circ \omega_0)=\omega^b_x\cdot (u^c\circ\omega_0)\cdot \lambda$, whose adjoint action on $(\lie{G})_\bb{C}$ is identical to that of $\omega_{x,\infty}^{Q,b}$. So $v^{Q,a}_x$ satisfies conditions in \cite[Prop. 4.6 (b)]{Pin89}. Hence, the statement follows by the uniqueness part of \textit{loc. cit}.
\end{proofof}
For $i=a,b$, the map $X_i\to X^\ad$ is injective and sends connected components to connected components. So for any cusp label representative $\Phi$, the connected component $X_{a,\Phi}^+$ of $X_a$ associated with $\Phi$ corresponds to a connected component $X_{b,\Phi}^+$ of $X_b$; we denote them by $X_\Phi^+$ abusively. So there is a natural one-to-one correspondence between cusp label representatives of $(G,X_a)$ and those of $(G,X_b)$, sending each of three data identically. Hence, any cone decomposition $\Sigma$ for $(G,X_a)$ can also be viewed as a cone decomposition for $(G,X_b)$. We can abusively denote by $\Phi$ a cusp label representative of both $(G,X_a)$ and $(G,X_b)$.
\begin{prop}\label{ppt-p}Fix a neat open compact subgroup $K\sbst G(\A)$. With the conventions in Lemma \ref{lem-zinv-ref}, let $\Sigma$ be any admissible $ZP$-invariant cone decomposition for $(G,X_a,K)$ with property $\ca{P}$. Then $\Sigma$ also defines an admissible $ZP$-invariant cone decomposition for $(G,X_b,K)$ with property $\ca{P}$.
\end{prop}
\begin{proof}
By Lemma \ref{lem-large-lattices}, the definition of smoothness or projectivity does not depend on the center of $P_\Phi$, so $\Sigma$ also defines a cone decomposition for $(G,X_b)$ with property $\ca{P}$. Fix an admissible $\bb{Q}$-parabolic subgroup and a connected component $X^+$. By Lemma \ref{izpq}, 
we have $I_{ZP}(Q)\iso\mrm{Stab}_{Q(\bb{Q})}(D_{Q,X^+})ZP_Q(\bb{A}_f)\backslash G(\bb{A}_f)/K$. Then the rest of the statement follows from Lemma \ref{zp} below.
\end{proof}
\begin{lem}\label{zp} Let $P_\Phi^i$ be the smallest normal subgroup of $Q$ that $u_x^{Q,i}\circ h_\infty$ factors through, for $i=a,b$. Then  $\zpp^a=\zpp^b$. In particular, $\zpp^a(\bb{A}_f)=\zpp^b(\bb{A}_f)$.
\end{lem}
\begin{proof}
Let $m:Z\times Q\to Q$ be the multiplication homomorphism of $Z\times Q$. Since $m\circ (u^c\times u^{Q,b}_x)\circ h_\infty=u^{Q,a}_x\circ h_\infty$ by Lemma \ref{uq-mul}, and $\zpp^b$ is a normal subgroup of $Q$, we know $\zpp^b$ contains $P^a_\Phi$. Hence, $\zpp^a\sbst\zpp^b$. Symmetrically, we can get the opposite direction of inclusion, so the desired statement is proved.
\end{proof}

Let \gls{Ec} be the field of definition of the Hodge cocharacter $\mu^c$ associated with the homomorphism $c$. Note that $(Z,\{c\})$ can be viewed as a Shimura datum whose associated Shimura varieties are $0$-dimensional. Let $E_i:=E(G,X_i)$ be the reflex field of $(G,X_i)$. Let $E^*$ be a field extension of $\bb{Q}$ such that all geometrically connected components of $\sh_{K_{\Phi_a}}$ and $\sh_{K_{\Phi_b}}$ are defined over $E^*$ for all cusp label representatives $\Phi_a$ and $\Phi_b$ of $(G,X_a)$ and $(G,X_b)$. We have $E^c\sbst \wdtd{E}:=E_a\cdot E_b.$\par 
Let $K_Z:=K_{Z,p}K^p_Z$ be a neat open compact subgroup in $Z(\A)$, where $K_{Z,p}\sbst Z(\bb{Q}_p)\cap K_p$ and $K_Z^p\sbst Z(\A^p)\cap K^p$ is neat open compact. We can consider the zero-dimensional Shimura variety $\sh_{K_Z}$ associated with $(Z,\{c\})$ over $E^c$. Over complex points, $\sh_{K_Z}(\bb{C})=Z(\bb{Q})\bss \{c\}\times Z(\A)/K_Z$.\par
Denote by \gls{FKZ} a finite field extension of $E^c$ such that all geometrically connected components of $\sh_{K_Z}$ are defined over $F_{K_Z}$. The base change $(\sh_{K_Z})_{F_{K_Z}}$ is a disjoint union of finite points over $F_{K_Z}$; in particular, $\sh_{K_Z}$ is proper.\par
Let $S$ be the smallest $\bb{Q}$-subgroup of $Z$ that $c$ factors through. Then the set of cusp labels is represented by $\mrm{Cusp}_{K_Z}(Z,\{c\})\iso Z(\bb{Q})S(\A)\bss Z(\A)/K_Z$. For any $z\in Z(\A)$, we can therefore associate a cusp label $[z]$ in the double coset above.\par
Let $$m^\prime:Z\times G\to Z\times G$$ be the isomorphism induced by $(z,g)\mapsto (z,zg)$ on $R$-value points, where $R$ is any $\bb{Q}$-algebra. 
\begin{prop}\label{zp-ab-min}
Let $i=a$ or $b$. Let $\sh_K^{i,\mmin}:=\sh_K^\mmin(G,X_i)$ be the minimal compactification of $\sh_K^i:=\sh_K(G,X_i)$. With the conventions above, we have the following statements:
\begin{enumerate}
\item The isomorphism $m^\prime$ induces an isomorphism $[m^\prime]: \sh_{K_Z}\times\sh^b_K\to\sh_{K_Z}\times\sh^a_K$ over $\wdtd{E}$. Moreover, $\sh_{K_Z}\times\sh_K^{i,\mmin}\iso(\sh_{K_Z}\times\sh_K^{i})^\mmin$, and $[m^\prime]$ extends to an isomorphism between minimal compactifications 
$[m^\prime]^\mmin:\sh_{K_Z}\times\sh^{b,\mmin}_K\to\sh_{K_Z}\times\sh^{a,\mmin}_K$ over \gls{Ewdtd}.
\item Define $$[ZP^i(\Phi)]:=\{[\Phi^\prime]\in\cusp_K(G,X_i)|\Phi^\prime\sim_{ZP}\Phi\}.$$ 
Let $\mrm{Z}_{[ZP^i(\Phi)],K}:=\disju_{[\Phi']\in [ZP^i(\Phi)]}\mrm{Z}_{[\Phi'],K}$ be the disjoint union of strata in the minimal compactification corresponding to any cusp label $[\Phi^\prime]$ in $[ZP^i(\Phi)]$. Then 
\begin{equation}\label{min-decomp} \mrm{Z}_{[ZP^i(\Phi)],K}\iso \wdtd{\Delta}^{ZP}_{\Phi,K}\bss \sh_{\wdtd{K}_{\Phi,h}^i}(ZP_{\Phi,h}^i,ZP_{\Phi,h}^i(\bb{Q})D_{\Phi,h}),\end{equation}
where $\wdtd{\Delta}^{ZP}_{\Phi,K}:=\stb_{Q_\Phi(\bb{Q})}(D_\Phi)\cap ZP^i_\Phi(\A)g_\Phi Kg_\Phi^{-1},$ where $$\wdtd{K}^i_{\Phi,h}:=(ZP^i(\A)\cap g_\Phi Kg_\Phi^{-1})W_\Phi(\A)/W_\Phi(\A)$$ and where $ZP^i_{\Phi,h}:=ZP^i_\Phi/W_\Phi$.\par
\item $m^\prime:Z\times ZP^b_\Phi\to Z\times ZP^a_\Phi$ induces an isomorphism $[m^\prime]^\Phi:\sh_{K_Z}\times\mrm{Z}_{[ZP^b(\Phi)],K}\to \sh_{K_Z}\times\mrm{Z}_{[ZP^a(\Phi)],K}$ over $\wdtd{E}$.
\end{enumerate}
\end{prop}
\begin{proof}In the proof, we omit the subscripts of the level groups in the symbols of strata and delta groups to simplify the notation. All fiber products are over $\wdtd{E}$ if not specified. 
By \cite[12.3]{Pin89}, the projections from $\wdtd{\sh}^i:=\sh_{K_Z}\times\sh^i_K$ to the first and the second factors induce corresponding morphisms between minimal compactifications $\pi_1^\mmin:\wdtd{\sh}^{i,\mmin}:=(\sh_{K_Z}\times\sh^i_K)^\mmin\to\sh_{K_Z}$ and $\pi_2^\mmin:\wdtd{\sh}^{i,\mmin}\to\sh^{i,\mmin}_K$. 
So there is a canonical morphism $\mrm{id}^\mmin: \wdtd{\sh}^{i,\mmin}\to \sh_{K_Z}\times \sh_K^{i,\mmin}$ induced by the universal property of the fiber product. The morphism $\mrm{id}^\mmin|_{\sh_{K_Z}\times \sh_K^i}$ is an isomorphism, and $\sh_{K_Z}\times \sh_K^{i,\mmin}$ is normal since $\sh_{K_Z}$ is smooth over $\wdtd{E}$.\par
Any cusp label $[\Phi^\times]$ of $\wdtd{\sh}^{i,\mmin}$ maps to a cusp label $[z]$ of $\sh_{K_Z}$ and also to a cusp label of $\sh_K^{i,\mmin}$, so there is a canonical morphism $\pi_{[\Phi^\times]}:\mrm{Z}_{[\Phi^\times]}\to\mrm{Z}_{[z]}\times\mrm{Z}_{[\Phi]}$ for some $[z]\in \cusp_{K_Z}(Z,\{c\})$ and some $[\Phi]\in\cusp_K(G,X_i)$ fitting into the following diagram
\begin{equation}\label{phi-cross-min}
    \begin{tikzcd}
    \sh_{K_{\Phi^\times,h}}(P_{\Phi^\times,h},D_{\Phi^\times,h})\arrow[rr,"\wdtd{\pi}_{[\Phi^\times]}"]\arrow{d}&&\sh_{[z]}\times\sh_{K_{\Phi,h}}(P_{\Phi,h},D_{\Phi,h})\arrow{d}\\
    \mrm{Z}_{[\Phi^\times]}\arrow[rr,"\pi_{[\Phi^\times]}"]&&\mrm{Z}_{[z]}\times\mrm{Z}_{[\Phi]}.
    \end{tikzcd}
\end{equation}
In (\ref{phi-cross-min}) above, the morphism $\wdtd{\pi}_{[\Phi^\times]}$ is locally quasi-finite: In fact, $P_{\Phi^\times}\hookrightarrow S\times P_\Phi$ is a closed embedding, so is the morphism between their Levi parts $P_{\Phi^\times,h}\hookrightarrow S\times P_{\Phi,h}$, then the locally quasi-finiteness follows from \cite[Prop. 3.8]{Pin89}. Moreover, the two vertical morphisms in (\ref{phi-cross-min}) are finite, because $\Delta_{\Phi^\times}$ (resp. $\Delta_{\Phi}$ and resp. $\Delta_{[z]}$) acts on $\sh_{K_{\Phi^\times},h}(P_{\Phi^\times,h},D_{\Phi^\times,h})$ (resp. $\sh_{K_\Phi,h}(P_{\Phi,h},D_{\Phi,h})$ and resp. $\sh_{[z]}$) through a finite group (resp. through a finite group and resp. trivially). Then $\pi_{[\Phi^\times]}$ is locally quasi-finite (see \cite[\href{https://stacks.math.columbia.edu/tag/0GWS}{Tag 0GWS}]{stacks-project}).   
So the first part is deduced from Zariski's main theorem, because $\mrm{id}^\mmin:\wdtd{\sh}^{i,\mmin}\to\sh_{K_Z}\times\sh_K^{i,\mmin}$ is a proper birational quasi-finite morphism between normal schemes, and from \cite[12.3]{Pin89} again.\par
By \S\ref{dbcoset}, 
\begin{equation}\label{dbcoset-min}
    \sh_{\wdtd{K}^i_\Phi}(ZP^i_{\Phi,h},ZP^i_{\Phi,h}(\bb{Q})D_{\Phi,h})\iso \disju_{\Phi^\prime} \wdtd{\Delta}_{ZP,\Phi^\prime,h}\bss\sh_{K_{\Phi'},h}^i.
\end{equation}
In the equation (\ref{dbcoset-min}) above, $\Phi'$ are cusp label representatives of $(G,X_i)$, the varieties $\sh_{K_{\Phi'},h}^i$ denote the boundary mixed Shimura varieties associated with $(P_{\Phi',h}^i,D_{\Phi',h})$ and $K_{\Phi',h}$, and the disjoint union runs over a set of representatives for $$\stb_{ZP^i_{\Phi,h}(\bb{Q})}(D_{\Phi,h})P^i_{\Phi,h}(\A)\bss ZP^i_{\Phi,h}(\A)/\wdtd{K}^i_{\Phi,h};$$
by strong approximation theorem for unipotent groups (see \cite[Lem. 2.7 and Prop. 6.6]{PR94}), the double coset above is isomorphic to 
\begin{equation}\label{label-zp}\stb_{ZP^i_\Phi(\bb{Q})}(D_\Phi)P^i_\Phi(\A)\bss ZP^i_\Phi(\A)/\wdtd{K}^i_\Phi,\end{equation}
which represents a set of cusp label representatives $\Phi^\prime=(Q_\Phi,X^+_\Phi,q_{\Phi^\prime}g_\Phi)$ for $q_{\Phi^\prime}$ consisting of a set of representatives for (\ref{label-zp}). 
Also, in the equation (\ref{dbcoset-min}) above, $\wdtd{\Delta}_{ZP,\Phi^\prime,h}:=\stb_{ZP^i_{\Phi,h}(\bb{Q})}(D_{\Phi,h})\cap P_{\Phi^\prime,h}^i(\A)(q_{\phip}\wdtd{K}^i_\Phi q_{\phip}^{-1}W_{\phip}(\A)/W_{\phip}(\A))$, and also for the same reason as the last sentence, the action of $\wdtd{\Delta}_{ZP,\phip}:=\stb_{ZP^i_\Phi(\bb{Q})}(D_\Phi)\cap P_{\phip}^i(\A)q_{\phip}\wdtd{K}^i_\Phi q_{\phip}^{-1}$ on $\sh^i_{K,h}(\phip)$ factors surjectively through $\wdtd{\Delta}_{ZP,\phip,h}$.\par
The quotient of (\ref{label-zp}) by $\wdtd{\Delta}^{ZP}_\Phi$ is 
\begin{equation}
    \wdtd{\Delta}^{ZP}_\Phi P^i_\Phi(\A)\bss ZP^i_\Phi(\A)/\wdtd{K}_\Phi^i,
\end{equation}
which is in bijection with $[ZP^i(\Phi)]$ by Proposition \ref{bij}. The stabilizer in $\stb_{Q(\bb{Q})}(D_{\Phi})$ of any $\Phi^\prime$ represented in the double coset (\ref{label-zp}) is $\wdtd{\Delta}_{\Phi^\prime}$. 
Hence, the quotient by $\wdtd{\Delta}_\Phi^{ZP}$ of the right-hand side of (\ref{dbcoset-min}) can be written as 
\begin{equation}
\disju\limits_{[\Phi^\prime]\in[ZP^i(\Phi)]} \Delta_{\phip}\wdtd{\Delta}_{ZP,\phip}\bss\sh^i_{K_{\Phi'},h}.\end{equation}
Since $\Delta_{ZP,\phip}:=\wdtd{\Delta}_{ZP,\phip}/P^i_{\phip}(\bb{Q})\sbst\Delta_{\phip}$, the second part follows.\par
For the last part, it suffices to check over complex points that if two points on $\sh_{\wdtd{K}^i_{\Phi,h}}(ZP^i_{\Phi,h},ZP^i_{\Phi,h}(\bb{Q})D_{\Phi,h})$ are equivalent under the relation defined by the action of $\wdtd{\Delta}^{ZP}_\Phi$, then they remain so after the action of any $[z]\in \sh_{K_Z}(\bb{C})$. This can be verified since $Z$ is the center of $G$.
\end{proof}
\begin{lem}\label{zp-iso}With the conventions above, for any cusp label representative $\Phi$ and any $\sigma\in\Sigma(\Phi)$ such that $\sigma^\circ\sbst \mbf{P}^+_\Phi$, we have the following statements:
\begin{enumerate}
\item Let $\stb_{Q_\Phi(\bb{Q})}(D_\Phi,\sigma)$ be the intersection $\stb_{Q_\Phi(\bb{Q})}(D_\Phi)\cap \stb_{Q_\Phi(\bb{Q})}(\sigma)$. Define $\wdtd{\Delta}^{ZP,\circ}_{\Phi,K}:=\stb_{Q_\Phi(\bb{Q})}(D_\Phi,\sigma)\cap ZP_\Phi^i(\A)g_\Phi Kg_\Phi^{-1}$, which is independent of the value of $i$. Define
$$[ZP^i(\Phi,\sigma)]:=\{(\Phi^\prime,\sigma^\prime)\in\cusp_K(G,X_i,\Sigma)|(\Phi^\prime,\sigma^\prime)\sim_{ZP}(\Phi,\sigma)\}.$$
Then 
\begin{equation}\label{tor-decomp-1}\disju\limits_{[(\Phi^\prime,\sigma^\prime)]\in[ZP^i(\Phi,\sigma)]}\Delta^{i,\circ}_{\Phi',K}\bss\sh_{K^i_{\Phi'}}(\sigma^\prime)\iso \wdtd{\Delta}^{ZP,\circ}_{\Phi,K}\bss \sh_{\wdtd{K}_\Phi^i}(ZP_\Phi^i,ZP^i_\Phi(\bb{Q})D_\Phi,\sigma)\end{equation}
over $\wdtd{E}$, where $K_{\Phi'}^i:=P^i(\A)\cap g_{\Phi'}Kg_{\Phi'}^{-1}$ and $\wdtd{K}_\Phi^i:=ZP^i(\A)\cap g_\Phi Kg_\Phi^{-1}$. 
\item The isomorphism (\ref{tor-decomp-1}) can be written as
\begin{equation}\label{tor-decomp-2}
    \disju\limits_{[\phip]\in[ZP^i(\Phi)]}\disju\limits_{\delta\in \Delta^i_{\phip}\bss\Delta^{ZP}_{\phip,K}}\Delta^i_{\phip,K}\bss\sh_{K_{\Phi'}}^i(\Delta^i_{\phip,K}\delta\sigma^\prime)\iso
    \Delta^{ZP}_{\Phi,K}\bss\sh_{\wdtd{K}^i_\Phi}(ZP^i_\Phi,ZP^i_\Phi(\bb{Q})D_\Phi,\Delta^{ZP}_{\Phi,K}\sigma).
\end{equation}
\end{enumerate}
\end{lem}
\begin{proof}Again, we omit the subscripts of the level groups in the symbols of strata and delta groups to simplify the notation.
For the first part, its proof is similar to that of the third part of Proposition \ref{zp-ab-min}. Firstly, $\wdtd{\Delta}^{ZP,\circ}_\Phi\bss \sh_{\wdtd{K}_\Phi^i}(ZP_\Phi^i,ZP_\Phi^i(\bb{Q})D_\Phi,\sigma)$ is isomorphic to
\begin{equation}\label{compl-quot}
\wdtd{\Delta}^{ZP,\circ}_\Phi\bss (\disju_{\Phi^\prime} \Delta_{ZP,\Phi^\prime}\bss\sh^i_{K_{\Phi'}}(\sigma)).
\end{equation}
In (\ref{compl-quot}) above, the disjoint union runs over a set of cusp label representatives $\phip=(Q_\Phi,X_\Phi^+,q_{\phip}g_\Phi)$, where $q_{\phip}$ consisting of a set of representatives for
\begin{equation}\label{good-rep}
\stb_{ZP^i_\Phi(\bb{Q})}(D_\Phi)P^i_\Phi(\A)\bss ZP^i_\Phi(\A)/\wdtd{K}^i_\Phi.\end{equation}
In fact, by \S\ref{dbcoset}, we have the following commutative diagram
\begin{equation}\label{tor-torsor}
    \begin{tikzcd}
    \sh_{\wdtd{K}^i_\Phi}(ZP^i_\Phi,ZP^i_\Phi(\bb{Q})D_\Phi)\arrow["\sim",rr]\arrow{d}&&\disju_{\phip}\Delta_{ZP,\phip}\bss\sh^i_{K_{\Phi'}}\arrow{d}\\
    \sh_{\overline{\wdtd{K}}^i_\Phi}(ZP^i_\Phi/U_\Phi,ZP^i_\Phi(\bb{Q})\overline{D}_\Phi)\arrow["\sim",rr]&&\disju_{\phip}\Delta_{ZP,\phip}\bss\overline{\sh}^i_{K_{\Phi'}}.
    \end{tikzcd}
\end{equation}
Then the corresponding twisted affine torus embedding with respect to $\sigma$ of the two vertical morphisms of (\ref{tor-torsor}) are isomorphic.\par
Since $ZP^i_\Phi(\bb{Q})$ acts on $\sigma$ via a scalar multiplication by \cite[Prop. 2.14 (a)]{Pin89}, we know $\stb_{ZP^i_\Phi(\bb{Q})}(D_\Phi)\sbst \stb_{Q(\bb{Q})}(D_\Phi,\sigma)$. So the quotient of (\ref{good-rep}) by $\wdtd{\Delta}_\Phi^{ZP,\circ}$ is
\begin{equation}
    \wdtd{\Delta}^{ZP,\circ}_\Phi P^i_\Phi(\A)\bss ZP^i_\Phi(\A)/\wdtd{K}^i_\Phi,
\end{equation}
which is in bijection with $[ZP^i(\Phi,\sigma)]$ by Proposition \ref{bij}. 
The stabilizer in $\stb_{Q(\bb{Q})}(D_{\Phi},\sigma)$ of any $\Phi^\prime$ represented in the double coset (\ref{good-rep}) is $\wdtd{\Delta}^\circ_{\Phi^\prime}$. 
Hence, the quotient (\ref{compl-quot}) can be written as \begin{equation}\disju\limits_{[(\Phi^\prime,\sigma^\prime)]\in[ZP^i(\Phi,\sigma)]} \Delta^\circ_{\phip}\Delta_{ZP,\phip}\bss\sh^i_{K_{\Phi'}}(\sigma).\end{equation}
Since $\Delta_{ZP,\phip}\sbst\Delta^\circ_{\phip}$ by \cite[Prop. 2.14 (a)]{Pin89} again, the desired result for representatives as in (\ref{good-rep}) is proved. For another choice of representatives $(\Phi^\dpr,\sigma^\dpr)$ for $[(\Phi^\prime,\sigma^\prime)]$, there is an equivalence 
$(\phip,\sigma^\prime)\bd{\gamma}{p}{K}{\sim}(\Phi^\dpr,\sigma^\dpr)$ for some $\gamma\in G(\bb{Q})$ and $p\in P^i_{\phip}(\A)$, which induces a canonical isomorphism $\Delta^{i,\circ}_{\phip}\bss \sh^i_{K_{\Phi'}}(\sigma')\iso \Delta^{i,\circ}_{\Phi^\dpr}\bss\sh_{K_{\Phi''}}^i(\sigma'')$. So the first part is proved.\par
The second part follows from the first part and Proposition \ref{fiber}.
\end{proof}
\begin{prop}\label{sumup} Let $i=a$ or $b$. Let $\sh_K^{i,\Sigma}:=\sh_{K}^\Sigma(G,X_i)$ be the toroidal compactification of $\sh_K^i:=\sh_{K}(G,X_i)$ associated with a $ZP$-invariant admissible cone decomposition $\Sigma$ with the property $\ca{P}$ as in Proposition \ref{ppt-p}.
Then: 
\begin{enumerate}
\item Let $\pi_2^i$ be the projection $\pi_2^i:\sh_{K_Z}\times\sh_K^i\to\sh_K^i$ to the second factor. Then the pullback of $\Sigma$ along $\pi_2^i$ induces a cone decomposition $\pi^{i,*}_2(\Sigma)$, which we abusively denote by $\Sigma$. Then there is a canonical isomorphism $\mrm{id}^\Sigma:(\sh_{K_Z}\times\sh_K^i)^\Sigma\xrightarrow{\sim}\sh_{K_Z}\times\sh_K^{i,\Sigma}$. So $[m^\prime]$ extends to an isomorphism $[m^\prime]^\Sigma:\sh_{K_Z}\times\sh^{b,\Sigma}_K\to\sh_{K_Z}\times\sh^{a,\Sigma}_K$ over $\wdtd{E}$.
\item 
Define $$\mrm{Z}_{[ZP^i(\Phi,\sigma)],K}:=\disju_{[(\Phi^\prime,\sigma^\prime)]\in[ZP^i(\Phi,\sigma)]} \mrm{Z}^i_{[(\Phi^\prime,\sigma^\prime)],K}.
$$
Then the restriction of $m^\prime$ to $ZP_\Phi^b$, $m^\prime:Z\times ZP^b_\Phi\to Z\times ZP^a_\Phi$, induces an isomorphism denoted by $[m^\prime]^{\Phi,\sigma}$, $$\sh_{K_Z}\times \disju_{[(\Phi^\prime,\sigma^\prime)]\in[ZP^b(\Phi,\sigma)]}\Delta^{b,\circ}_{\Phi',K}\bss\sh_{K_{\Phi'}}^b(\sigma')\xrightarrow{[m']^{\Phi,\sigma}}\sh_{K_Z}\times\disju_{[(\Phi^\prime,\sigma^\prime)]\in[ZP^a(\Phi,\sigma)]}\Delta^{a,\circ}_{\Phi',K}\bss\sh_{K_{\Phi'}}^a(\sigma')$$ over $\wdtd{E}$. We also have a similar isomorphism for the unions of $\sigma'$-closed strata.\par
\item The following diagram commutes:
\begin{equation}\label{diag-zp-a-b}
   \begin{tikzcd}
  {(\sh_{K_Z})^{\wat{}}_{[z]}\times (\sh^{b,\Sigma}_K)^{\wat{}}_{\mrm{Z}_{[ZP^b(\Phi,\sigma)],K}}}\arrow[d,"{[m^\prime]}^\Sigma"]\arrow[r,"\sim"]& 
    {(\sh_{K_Z})^{\wat{}}_{[z]}\times \disju\limits_{[(\Phi^\prime,\sigma^\prime)]\in[ZP^b(\Phi,\sigma)]}\Delta_{\Phi',K}^{b,\circ}\backslash(\sh_{K_{\Phi'}}^b(\sigma'))^{\wat{}}_{\sh^b_{K_{\Phi'},\sigma'}}}\arrow[d,"{[m^\prime]}^{\Phi,\sigma}"]\\
    {(\sh_{K_Z})^{\wat{}}_{[z]}\times (\sh^{a,\Sigma}_K)^{\wat{}}_{\mrm{Z}_{[ZP^a(\Phi,\sigma)],K}}}\arrow[r,"\sim"]& {(\sh_{K_Z})^{\wat{}}_{[z]}\times \disju\limits_{[(\Phi^\prime,\sigma^\prime)]\in[ZP^a(\Phi,\sigma)]}\Delta_{\Phi',K}^{a,\circ}\bss(\sh_{K_{\Phi'}}^a(\sigma'))^{\wat{}}_{\sh^a_{K_{\Phi'},\sigma'}}}.\end{tikzcd}
\end{equation}
\end{enumerate}
\end{prop}
\begin{proof}Again, we omit the subscripts of the level groups in the symbols of strata, cones and delta groups to simplify the notation.
We prove the first part with the same strategy as in Proposition \ref{zp-ab-min}. For any cusp label $[\Phi^\times]$ mapping to $[z]$ and $[\Phi]$ under $\pi_1^\mmin$ and $\pi_2^\mmin$ as in Proposition \ref{zp-ab-min}, we have $\mbf{P}_{\Phi^\times}\xrightarrow{\sim}\mbf{P}_{[z]}\times\mbf{P}_\Phi=\mbf{P}_\Phi$. For any cone $\sigma\in\Sigma(\Phi)$, we abusively denote by $\sigma$ the cone induced by $\mbf{P}_{\Phi^\times}\xrightarrow{\sim}\mbf{P}_{\Phi}$. 
We have the following commutative diagram
\begin{equation}
    \begin{tikzcd}
    \sh_{K_{\Phi^\times}}(P_{\Phi^\times},D_{\Phi^\times})\arrow{rr}\arrow{d}&&\mrm{Z}_{[z]}\times\sh_{K_\Phi}(P_\Phi,D_\Phi)\arrow{d}\\
    \sh_{\overline{K}_{\Phi^\times}}(\overline{P}_{\Phi^\times},\overline{D}_{\Phi^\times})\arrow{rr}&&\mrm{Z}_{[z]}\times \sh_{\overline{K}_\Phi}(\overline{P}_\Phi,\overline{D}_\Phi),
    \end{tikzcd}
\end{equation}
where the horizontal morphisms are locally quasi-finite. 
Then morphism $\mbf{E}_{K_{\Phi^\times}}\to \mbf{E}_{K_\Phi}$ is also finite. So the morphism between twisted affine torus embeddings with respect to $\sigma$, i.e., $\sh_{K_{\Phi^\times}}(P_{\Phi^\times},D_{\Phi^\times},\sigma)\to\mrm{Z}_{[z]}\times\sh_{K_\Phi}(P_\Phi,D_\Phi,\sigma)$, is locally quasi-finite.
So $\mrm{Z}_{[(\Phi^\times,\sigma)]}\to\mrm{Z}_{[z]}\times\mrm{Z}_{[(\Phi,\sigma)]}$ is locally quasi-finite. Then $\mrm{id}^\Sigma$ is an isomorphism by Zariski's main theorem, as it is a proper birational and quasi-finite morphism between normal schemes. Other statements in part one follow from \cite[6.7(b) and 12.4(b)]{Pin89}.\par
The morphism $m^\prime$ induces an isomorphism $$ m^\prime:\sh_{K_Z}\times\sh_{\wdtd{K}_\Phi^b}(ZP_\Phi^b,ZP^b_\Phi(\bb{Q})D_\Phi)\to\sh_{K_Z}\times\sh_{\wdtd{K}_\Phi^a}(ZP^a_\Phi,ZP^a_\Phi(\bb{Q})D_\Phi)$$ by \cite[3.5 (c)]{Pin89}, whose quotient by $\wdtd{\Delta}^{ZP,\circ}_\Phi$ on the second factor is $[m^\prime]^{\Phi,\sigma}$. So the second part is proved by Lemma \ref{zp-iso}.
The third part follows from the first two parts, and \cite[Thm. 2.1.27 and 2.1.28]{Mad19}. Note that the second paragraph also implies a similar isomorphism between closed strata.
\end{proof}
Define \gls{EKZwdtd} to be $\wdtd{E}_{K_Z}:= \wdtd{E}\cdot F_{K_Z}$. Then $\sh_{K_Z}$ splits as a finite union of $\spec \wdtd{E}_{K_Z}$ over $\wdtd{E}_{K_Z}$.
\begin{cor}\label{zp-ab}
With the conventions above, let $s: \spec \wdtd{E}_{K_Z}\to \sh_{K_Z}$ be any section of $(\sh_{K_Z})_{\wdtd{E}_{K_Z}}$ over $\wdtd{E}_{K_Z}$ corresponding to an element $z\in Z(\A)$. Denote by $\pi_2$ the projection of a fiber product to the second factor. Denote by $[g]$ the Hecke action of any element $g\in G(\A)$. We have the following statements:
\begin{enumerate}
\item The compositions $$(m')^\Phi_s:=[z^{-1}]\circ\pi_2\circ (m')\circ(s\times \mrm{id}):\sh_{\K^b_\Phi}(ZP^b_\Phi,ZP^b_\Phi(\bb{Q})D_\Phi)\to\sh_{\K^a_\Phi}(ZP^a_\Phi,ZP^a_\Phi(\bb{Q})D_\Phi)$$ and $$[m^\prime]^\Phi_s:=[z^{-1}]\circ \pi_2\circ[m^\prime]^\Phi\circ (s\times \mrm{id}):\mrm{Z}_{[ZP^b(\Phi)],K}\to\mrm{Z}_{[ZP^a(\Phi)],K}$$ are isomorphisms over $\wdtd{E}_{K_Z}$. Also, $[m^\prime]^\mmin_s:=[z^{-1}]\circ\pi_2^a\circ[m^\prime]^\mmin\circ(s\times\mrm{id}):\sh^{b,\mmin}_K\to\sh^{a,\mmin}_K$ is an isomorphism over $\wdtd{E}_{K_Z}$ extending $[m']_s:=[z^{-1}]\circ\pi_2^a\circ [m']\circ (s\times \mrm{id})$.
\item  Then the composition $[m^\prime]^\Sigma_s:=[z^{-1}]\circ\pi_2\circ[m^\prime]^\Sigma\circ (s\times \mrm{id}):\sh^{b,\Sigma}_K\to\sh^{a,\Sigma}_K$ is an isomorphism over $\wdtd{E}_{K_Z}$. Similarly, for any cusp label representative $\Phi$ and any $\sigma\in\Sigma(\Phi)$, there is an isomorphism $[m^\prime]^{\Phi,\sigma}_s:=[z^{-1}]\circ\pi_2\circ[m^\prime]^{\Phi,\sigma}\circ \disju(s\times \mrm{id}):\disju\limits_{[(\Phi^\prime,\sigma^\prime)]\in[ZP^b(\Phi,\sigma)]}\Delta^{b,\circ}_{\phip}\bss\sh^b_{K_{\Phi'}}(\sigma')\to\disju\limits_{[(\Phi^\prime,\sigma^\prime)]\in[ZP^a(\Phi,\sigma)]}\Delta^{a,\circ}_{\phip}\bss\sh^a_{K_{\Phi'}}(\sigma')$ over $\wdtd{E}_{K_Z}$, and it induces an isomorphism 
$$\disju\limits_{[(\Phi^\prime,\sigma^\prime)]\in[ZP^b(\Phi,\sigma)]}\Delta^{b,\circ}_{\phip,K}\bss\sh^b_{K_{\Phi'},\sigma'}\to\disju\limits_{[(\Phi^\prime,\sigma^\prime)]\in[ZP^a(\Phi,\sigma)]}\Delta^{a,\circ}_{\phip,K}\bss\sh^a_{K_{\Phi'},\sigma'}$$ over $\wdtd{E}_{K_Z}$.\par
Moreover, the following diagram over $\wdtd{E}_{K_Z}$ commutes:
\begin{equation}
   \begin{tikzcd}
  {(\sh^{b,\Sigma}_K)^{\wat{}}_{\mrm{Z}_{[ZP^b(\Phi,\sigma)],K}}}\arrow[d,"{[m^\prime]}^\Sigma_s"]\arrow[r,"\sim"]& 
    {\disju\limits_{[(\Phi^\prime,\sigma^\prime)]\in[ZP^b(\Phi,\sigma)]}\Delta_{\Phi^\prime,K}^{b,\circ}\backslash(\sh_K^b(\Phi^\prime,\sigma^\prime))^{\wat{}}_{\sh^b_{K_{\Phi'},\sigma'}}}\arrow[d,"{[m^\prime]}^{\Phi,\sigma}_s"]\\
    {(\sh^{a,\Sigma}_K)^{\wat{}}_{\mrm{Z}_{[ZP^a(\Phi,\sigma)],K}}}\arrow[r,"\sim"]& {\disju\limits_{[(\Phi^\prime,\sigma^\prime)]\in[ZP^a(\Phi,\sigma)]}\Delta_{\Phi^\prime,K}^{a,\circ}\bss(\sh_K^a(\Phi^\prime,\sigma^\prime))^{\wat{}}_{\sh^a_{K_{\Phi'},\sigma'}}}.\end{tikzcd}
\end{equation}
\item All the isomorphisms $[m']_s$, $[m']^\Phi_s$, $[m']^\mmin_s$, $[m']^\Sigma_s$ and $[m']^{\Phi,\sigma}_s$ above are independent of the choice of $s$.
\end{enumerate}
\end{cor}
\begin{proof}By \cite[12.4(b)]{Pin89} and our construction, all of the morphisms above are defined over $\wdtd{E}_{K_Z}$. Then we can check all the statements above over field extensions.
Since all of $[m']^?$ for $?= \empty,\ \Phi,\ (\Phi,\sigma),\ \Sigma,\ \mmin$ are isomorphisms according to Proposition \ref{zp-ab-min} and Proposition \ref{sumup}, and since $\pi_2$, $s\times \mrm{id}$ and $[z]$ all induce isomorphisms on each geometrically connected component, we only need to count the number of geometrically connected components. 
Over $\overline{\wdtd{E}}$, the morphism $[m^\prime]_s^\Sigma$ (resp. $[m^\prime]_s^{\Phi,\sigma}$) induces isomorphism of connected components via $\pi_0(\sh^{b,\Sigma}_K)\to [z]\times\mrm{Stab}_{G(\bb{Q})}(X^+)\bss G(\A)/K
\to[z]\times\mrm{Stab}_{G(\bb{Q})}(X^+)\bss G(\A)/K\to\pi_0(\sh^{a,\Sigma}_K)$ 

\noindent(resp. $\pi_0(\disju_{[(\Phi^\prime,\sigma^\prime)]\in[ZP^b(\Phi,\sigma)]}\sh^{b}_{K_{\Phi'}}(\sigma'))\to $ $[z]\times\mrm{Stab}_{Q(\bb{Q})}(D^+_{\Phi})\bss ZP^b_{\Phi}(\A)g_\Phi K/K\to$ 

$[z]\times\mrm{Stab}_{Q(\bb{Q})}(D^+_{\Phi})\bss ZP^a_{\Phi}(\A)g_\Phi zK/K
\to$ $\pi_0(\disju_{[(\Phi^\prime,\sigma^\prime)]\in[ZP^a(\Phi,\sigma)]}\sh^{a}_{K_{\Phi'}}(\sigma'))$) (see Lemma \ref{comp-tor}). The morphism in the bracket is an isomorphism because of  Lemma \ref{zp} and Lemma \ref{zp-iso}. Other isomorphisms are checked in a similar way. So the first part is proved. The second part is also true since we can verify it over $\overline{\wdtd{E}}$, since we have (\ref{diag-zp-a-b}) and since we can restrict it to geometrically connected components of $\sh^{i,\Sigma}_K$, and open and closed subschemes of $\sh^b_{K_{\Phi'}}(\sigma')$. Finally, we can check over complex points that $[m']_s$ and $[m']_s^{\Phi}$ are isomorphisms that are independent of the choice of $s$ or $z$. Then the definitions of all the other morphisms are independent of the choice of $s$ or $z$ because of density.
\end{proof}
\begin{convention}\label{conv-m}
    We denote $[m']_s$, $[m']^\Phi_s$, $[m']^\mmin_s$, $[m']^\Sigma_s$ and $[m']^{\Phi,\sigma}_s$ by $m$, $m^\Phi$, $m^\mmin$, $m^\Sigma$ and $m^{\Phi,\sigma}$, respectively, since they do not depend on the choice of sections.
\end{convention}
\subsubsection{}\label{sss-conclusion}
We finally come to the main setup. Let $(G_0,X_0)$ be a Hodge-type Shimura datum. Denote by $\iota: (G_0,X_0)\hookrightarrow (G^\ddag,X^\ddag)$ an embedding into a Siegel Shimura datum. Let $(G_2,X_2)$ be an abelian-type Shimura datum with a central isogeny $\pi^\der:G_0^\der\to G_2^\der$ between the derived groups of $G_0$ and $G_2$ such that $\pi^\der$ induces an isomorphism $(G_0^\ad,X_0^\ad)\iso(G_2^\ad,X_2^\ad)$. 
If $G_{2,\bb{Q}_p}$ is unramified, by \cite[Lem. 3.4.13]{Kis10}, there is a $G_{0}$ such that $G_{0,\bb{Q}_p}$ is also unramified.\par
\begin{construction}\label{cons-group-setting}\upshape
Let $C^\der$ be the kernel of $\pi^\der$. Let $G'_0:=G_0/C^\der$ be the quotient of $G_0$ by $C^\der$. Let $G^*$ be the contracted product $G^*:=G_2\times^{G^\der_2}G'_0$. The center $Z^*$ of $G^*$ can be embedded in some quasi-split $\bb{Q}$-torus $Z'$ by \cite[Prop. 2.1, p. 55]{PR94}. We define $G:=G^*\times^{Z^*}Z'$.\par
By the construction in the last paragraph, there are natural homomorphisms $\pi^a:G_2\to G^*\to G$ and $\pi^b:G_0\to G^*\to G$, where $\pi^a$ is an embedding and $\pi^b$ has finite kernel. Let $X_a$ (resp. $X_b$) be the $G(\bb{R})$-orbit of one of the elements in $X_2$ (resp. $X_0$). 
Since $H^1(\bb{R},Z')$ is trivial by \cite[Lem. 2.4, p.73]{PR94}, $X_a\iso X_b\iso X^\ad$. In particular, the Shimura data $(G,X_a)$ and $(G,X_b)$ satisfy the assumptions in \S\ref{a,b}.
\end{construction}
Fix an open compact subgroup $K_{2,p}$ of $G_2(\bb{Q}_p)$, and a neat open compact subgroup of $K_2^p$ of $G_2(\A)$. Choose an open compact $K_p\sbst G(\bb{Q}_p)$ such that $K_p\cap G_2(\bb{Q}_p)=K_{2,p}$. Let $K_{0,p}$ be an open compact subgroup of $G_0(\bb{Q}_p)$ contained in the preimage of $K_p$ under $\pi^b$. Let $K_{0}:=K_{0,p}K^p_0$ (resp. $K:=K_pK^{p}$), where $K^p_0$ (resp. $K^{p}$) is a neat open compact subgroup of $G_0(\A^p)$ (resp. $G(\A^p)$).\par
By Proposition \ref{ext-imp}, we can choose $K^{p}$ sufficiently small such that $\pi^a$ induces an open and closed embedding $\pi^a:\sh_{K_2}(G_2,X_2)\hookrightarrow \sh_{K}(G,X_a)$ and an isomorphism $\pi^a:\sh_{K_{\Phi_2}}(P_{\Phi_2},D_{\Phi_2})\iso \sh_{K_\Phi}(P_\Phi,D_\Phi),$ for any cusp label representative $\Phi_2$ mapping to $\Phi$ under $\pi^a$.\par
Let $\im G_0(\A)$ be the image of $G_0(\A)$ under $\pi^b$. Define an index set \gls{IGG0} as follows:
$$I_{G/G_0}:=I_{G/G_0,K}:=\stb_{G(\bb{Q})}(X_0)\im G_0(\A)\bss G(\A)/K.$$ 
We omit the subscript $K$ in the notation if it is clear in the context. 
From now on, we fix a set of representatives $\{g_\alpha\}_{\alpha\in I_{G/G_0}}$ of $I_{G/G_0}$; moreover, if $G_0$ and $G$ are all chosen to be unramified at $p$, by \cite[Lem. 2.2.6]{Kis10}, we can choose $g_\alpha$ in $G(\A^p)$. \par 
For any $\alpha$, choose a neat open compact $K_0^{p,\alpha}$ of $G_0(\A^p)$ contained in $\pi^{b,-1}(g_\alpha K^{p}g_\alpha^{-1})$, and choose an open compact subgroup $K_{0,p}^\alpha\sbst G_0(\bb{Q}_p)$ contained in $\pi^{b,-1}(g_\alpha K_pg_\alpha^{-1})$. Denote by $\pi^b(g_\alpha): \sh_{K_0^\alpha}(G_0,X_0)\to\sh_{K}(G,X_b)$ the morphism defined as in \S\ref{subsubsec-can-model}. For any $\alpha$, we also choose $K^{\ddag,\alpha}:=K_p^{\ddag,\alpha}K^{\ddag,\alpha,p}$ such that: 
\begin{itemize}
\item $K_p^{\ddag,\alpha}$ are open compact subgroups of $G^{\ddag}(\bb{Q}_p)$ and $K_p^{\ddag,\alpha}\cap G_0(\bb{Q}_p)\spst K_{0,p}^\alpha$.
\item $K^{\ddag,\alpha,p}$ are neat open compact subgroups of $G^\ddag(\Ap)$ containing $K^{\alpha,p}_0$.
\item $\iota^\alpha: \sh_{K_0^\alpha}(G_0,X_0)\to \sh_{K^{\ddag,\alpha}}(G^\ddag,X^\ddag)$ are finite morphisms.
\end{itemize}
Since $E_a\sbst$\gls{E2}$:=E(G_2,X_2)$ and $E_b\sbst$\gls{E0}$:=E(G_0,X_0)$, we have $\wdtd{E}\sbst$\gls{Eprime}$:=E_0\cdot E_2$.
\begin{prop}\label{zp-cones}
Let $\Sigma_2$ be an admissible cone decomposition for $\sh_{K_2}(G_2,X_2)$. Then by Proposition \ref{ext-zp}, there is an admissible $ZP$-invariant cone decomposition $\Sigma$ for $\sh_{K}(G,X_a)$ and $\sh_{K}(G,X_b)$ such that $\pi^{a,*}(\Sigma)$ refines $\Sigma_2$. \par
For each $g_\alpha$ defined above, let $\Sigma_0^\alpha:=\pi^b(g_\alpha)^*(\Sigma)$. 
Then $\Sigma$ can be chosen such that, for each $\alpha\in I_{G/G_0}$, $\Sigma_0^\alpha$ is induced by an admissible cone decomposition $\Sigma^{\ddag,\alpha}$ for $\sh_{K^{\ddag,\alpha}}(G^\ddag,X^\ddag)$, and such that $\Sigma^{\ddag,\alpha}$ is a refinement of a smooth or smooth and projective admissible cone decomposition for $\sh_{K^{\ddag,\alpha}}(G^\ddag,X^\ddag)$. Moreover, we can further choose $\Sigma$ to be smooth, projective, or to be both smooth and projective.
\end{prop}
\begin{proof}
First, we apply the argument of \cite[Prop. 4.10]{Lan19}. For any $\alpha\in I_{G/G_0}$, we have an embedding $\iota^\alpha: \sh_{K_0^\alpha}\hookrightarrow \sh_{K^{\ddag,\alpha}}$ by construction. For any cusp label representative $\Phi_0^\alpha\in \ca{CLR}(G_0,X_0)$, the map $\iota^\alpha$ sends $\Phi_0^\alpha$ to a cusp label representative $\Phi^{\ddag,\alpha}\in \ca{CLR}(G^\ddag,X^\ddag)$. This assignment induces embeddings of cones $\iota^\alpha_{\mbf{P}}:\mbf{P}_{\Phi_0^\alpha}\hookrightarrow \mbf{P}_{\Phi^{\ddag,\alpha}}$ and $\bb{R}$-spaces $\iota^\alpha_U: U_{\Phi_0^\alpha}(\bb{R})(-1)\hookrightarrow U_{\Phi^{\ddag,\alpha}}(\bb{R})(-1)$. The maps $\iota^\alpha_{\mbf{P}}$ and $\iota^\alpha_U$ are cut out by some hyperplane $H_{\Phi_0^\alpha}\sbst U_{\Phi^{\ddag,\alpha}}(\bb{R})(-1)$. 
By \emph{loc. cit.}, there is a pair of admissible cone decompositions (which are not necessarily smooth or projective) $(\Sigma^{\alpha}_{0,\mrm{scp}},\Sigma^{\ddag,\alpha}_\mrm{scp})$, such that the two cone decompositions are \emph{strictly compatible with each other} (see Definition \ref{def-func-cones}). 
That is, for any $\Phi_0^\alpha$ mapping to $\Phi^{\ddag,\alpha}$, and any $\sigma^{\ddag,\alpha}\in \Sigma^{\ddag,\alpha}_{\mrm{scp}}(\Phi^{\ddag,\alpha})$ such that $\sigma^{\ddag,\alpha,\circ}\cap H_{\Phi_0^\alpha}\neq \emptyset$, we have  $\sigma^{\ddag,\alpha}\sbst H_{\Phi_0^\alpha}$ and $\sigma^{\ddag,\alpha}\in \Sigma^\alpha_{0,\mrm{scp}}(\Phi_0^\alpha)$. In other words, all hyperplanes $H_{\Phi_0^\alpha}$ either intersect with the cone $\sigma^{\ddag,\alpha}$ on a face of its boundaries, or contain the cone, and each $\sigma_0^\alpha\in \Sigma_{0,\mrm{scp}}^{\alpha,+}(\Phi_0^\alpha)$ also maps exactly to a cone in $\Sigma_{\mrm{scp}}^{\ddag,\alpha,+}(\Phi^{\ddag,\alpha})$. For any refinement $\Sigma_0^\alpha$ of $\Sigma_{0,\mrm{scp}}^\alpha$ with property $\ca{P}$, since all the refinements are on the boundaries of each $\sigma^{\ddag,\alpha}$, or are on the whole $\sigma^{\ddag,\alpha}$. By \cite[Prop. 5.23 and Lem. 5.25]{Pin89} and \cite[pp. 32-35]{KKMS73}, there is a refinement $\Sigma^{\ddag,\alpha}$ of $\Sigma^{\ddag,\alpha}_{\mrm{scp}}$ such that $\Sigma_0^\alpha$ and $\Sigma^{\ddag,\alpha}$ are strictly compatible with each other.\par
Then it suffices to choose an admissible cone decomposition $\Sigma$ for $\sh_{K}(G,X_b)$, such that the induced cone decompositions $\pi^{b}(g_\alpha)^*(\Sigma)$ are refinements of $\Sigma_{0,\mrm{scp}}^\alpha$. Indeed, if we can do so, we can refine $\Sigma$ further so that it satisfies the property $\ca{P}$ and is $ZP$-invariant. Denote this refinement by $\Sigma'$. Then the pullbacks $\pi^b(g_\alpha)^*(\Sigma')$ are refinements of $\pi^{b}(g_\alpha)^*(\Sigma)$ and in particular, of $\Sigma_{0,\mrm{scp}}^\alpha$. Then, by the last paragraph, there are refinements $\Sigma^{\ddag,\alpha,\prime}$ of $\Sigma^{\ddag,\alpha}_{\mrm{scp}}$, such that $\pi^b(g_\alpha)^*(\Sigma')$ and $\Sigma^{\ddag,\alpha,\prime}$ are strictly compatible with each other for all $\alpha$. In particular, $\pi^b(g_\alpha)^*(\Sigma')$ are induced by $\Sigma^{\ddag,\alpha,\prime}$ for all $\alpha$.\par
In fact, we can choose such a $\Sigma$. Since the maps $\cusp_{K_0^\alpha}(G_0,X_0)\to\cusp_{K}(G,X_b)$ induced by $\pi^b(g_\alpha)$ are maps between finite sets and with finite fibers and since $|I_{G/G_0}|$ is finite, we only have finitely many prescribed cone decompositions for each $\mbf{P}_{\Phi}$ with $\Phi\in\ca{CLR}(G,X_b)$. For any $\Phi_0^\alpha\in\ca{CLR}(G_0,X_0)$ mapping to $\Phi\in \ca{CLR}(G,X_b)$, the map $\pi^b(g_\alpha)_{\mbf{P}}:\mbf{P}_{\Phi_0^\alpha}\to \mbf{P}_{\Phi}$ induced by $\pi^b(g_\alpha)$ is an isomorphism. We then choose a common refinement of the cone decompositions of all $[\Phi_0^\alpha]$ in $\disju_{\alpha\in I_{G/G_0}}\cusp_{K_0^\alpha}(G_0,X_0)$ mapping to $[\Phi]$; we can do this since there are only finitely many such cone decompositions for $[\Phi]$.
\end{proof}
\begin{definition}\label{def-zp-mix-shimura}
As \S\ref{subsubsec-boundary}, define $\K_\Phi:=\K_\Phi^i:=ZP^i_\Phi(\A)\cap g_\Phi Kg_\Phi^{-1}$, $\overline{\K}_\Phi:=\overline{\K}_\Phi^i:=ZP^i_\Phi/U_\Phi(\A)\cap (g_\Phi K g_\Phi^{-1}/g_\Phi Kg_\Phi^{-1}\cap U_\Phi(\A))$ and $\K_{\Phi,h}:=\K_{\Phi,h}^i:=ZP^i_\Phi/W_\Phi(\A)\cap(g_\Phi Kg_\Phi^{-1}/W_\Phi(\A)\cap g_\Phi Kg_\Phi^{-1})$. Define $\sh_{\K^i_\Phi}:=\sh_{\K_\Phi}(ZP^i_\Phi,ZP^i_\Phi(\bb{Q})D_\Phi)$, 
$\overline{\sh}_{\K^i_\Phi}:=\sh_{\overline{\K}_\Phi}(ZP^i_\Phi/U_\Phi,ZP^i_\Phi(\bb{Q})\overline{D}_\Phi)$, and 
$\sh_{\K^i_\Phi,h}:=\sh_{\K_{\Phi,h}}(ZP^i_\Phi/W_\Phi,ZP_\Phi^i(\bb{Q})D_{\Phi,h})$. 
\end{definition}
Let $\gls{EtdKPhi}$ be the algebraic torus whose group of cocharacters corresponds to the lattice $\bm{\Lambda}_{K_\Phi}\sbst U_\Phi(\bb{R})(-1)$. 
By \cite[3.13]{Pin89} and Lemma \ref{lem-large-lattices}, the variety $\sh_{\K^i_\Phi}$ is a $\mbf{E}_{\K_\Phi}$-torsor over $\overline{\sh}_{\K_\Phi^i}$, and $\overline{\sh}_{\K_\Phi^i}$ is an abelian scheme torsor over $\sh_{\K^i_{\Phi},h}$.

\newpage
\section{One-motives and the theory of degeneration}\label{sec-1-motives-and-deg}
The goal of this {section} is to review and explain Faltings-Chai's theory of degenerating abelian schemes (see \cite{FC90} and \cite{Lan13}). As a preparation for the next {section}, we present the theory in the language of $1$-motives and aim to relate the notions in different references (see, e.g., \cite{Mil90}, \cite{Pin89}, \cite{Str10}, \cite[\S2.2]{Mad19}, \cite{FC90}, \cite{Lan13} and \cite{Lan12b}).\par 
We start with Siegel moduli and $1$-motives (see \S\ref{subsec-siegel} and \S\ref{subsec-1-mot}).
We then review the theory of degeneration (see \S\ref{subsec-deg-data}), and explain the moduli interpretations of boundary mixed Shimura varieties of integral models of Siegel Shimura varieties (see \S\ref{subsec-mxhg-cont}). In particular, we check that for Siegel Shimura varieties, definitions of cusp labels, boundary charts, and level structures of boundary charts given in various references are the same. In \S\ref{subsec-tor-hdg-review}, we review the main theorems of toroidal compactifications of integral models of Siegel/PEL/Hodge-type Shimura varieties.
The results that will be used in subsequent {sections} will be proved throughout the course of the exposition.\par
Since we will only study Hodge-type and Siegel-type Shimura varieties in this {section}, the index set $I_{G/G_0}$ in the previous {section} is not important here; all superscripts $\alpha$ in the conventions of \S\ref{sss-conclusion} will be omitted.
\subsection{Siegel moduli}\label{subsec-siegel}
Recall that $G_0$ is a connected reductive group over $\bb{Q}$ and \gls{G0X0} is a Hodge-type Shimura datum with an embedding $\iota:(G_0,X_0)\hookrightarrow$\gls{GddagXddag}$:=(\mrm{GSp}(V,\psi),\bb{S}^\pm)$. We denote $G^\ddag:=\mrm{GSp}(V,\psi)$ and $X^\ddag:=\bb{S}^\pm$. We denote by $\nu:=\nu_{G^\ddag}: G^\ddag\to \bb{G}_m$ the similitude character of $G^\ddag$. Let $K_0:=K_{0,p}K^p_0$ (resp. $K^\ddag:=K^\ddag_pK^{\ddag,p}$) be a neat open compact subgroup of $G_0(\A)$ (resp. $G^\ddag(\A)$) such that $K_0\sbst K^\ddag\sbst G^\ddag(\A)$, such that $K^p_0$ (resp. $K^{\ddag,p}$) is a neat open compact subgroup of $G_0(\Ap)$ (resp. $G^\ddag(\Ap)$) and such that $K_{0,p}$ (resp. $K^\ddag_p$) is an open compact subgroup of $G_0(\bb{Q}_p)$ (resp. $G^\ddag(\bb{Q}_p)$). 
By Zarhin's trick, after replacing $V$ with $V^{\oplus 8}$ if necessary, we choose a $\bb{Z}$-lattice \gls{VbbZ}$\sbst V$ such that $V_{\bb{Z}}$ is self-dual with respect to $\psi$, such that $K_p^\ddag$ is exactly the stabilizer of $V_{\bb{Z}_p}$ in $G^\ddag(\bb{Q}_p)$, and such that $K^\ddag$ stabilizes $V_{\wat{\bb{Z}}}$.\par 
To avoid repeating similar definitions, let $\p$ be $\emptyset$ or $\{p\}$. If $\p=\emptyset$, then $K^{\ddag,\p}=K^\ddag$ and $K^\ddag_{\p}$ is empty. Denote by $\mbf{M}_{K^\ddag,\zbkpp}^{\mrm{isog}}$ the moduli problem $\mbf{M}_{K^\ddag,\zbkpp}^{\mrm{isog}}:(Sch/\zbkpp)^{\mrm{op}}\to (Grpds)$ defined as follows:\par
For any locally Noetherian and connected $\zbkpp$-scheme $S$, $\mbf{M}_{K^{\ddag,\p}}^{\mrm{isog}}(S)$ is the groupoid of tuples $(\ca{A},\lambda,[\varepsilon^\ddag]_{K^{\ddag,\p}})$, where 
\begin{itemize}
    \item $\ca{A}$ is an abelian scheme over $S$, 
    \item $\lambda$ is a $\zbkppt$-polarization (i.e., $[N]\circ \lambda$ is a (positive) polarization for some positive integer $N$ such that $(N,\p)=1$), and 
    \item $[\varepsilon^\ddag]_{K^{\ddag,\p}}$ is a $\pi_1(S,\overline{s})$-invariant $K^{\ddag,\p}$-orbit of $\App$-equivariant isomorphisms  $\varepsilon^\ddag:V_{\App}\xrightarrow{\sim}V^p\ca{A}_{\overline{s}}$ at some (and therefore any) geometric point $\overline{s}$ of $S$ (this is equivalent to saying that $[\varepsilon^\ddag]_{K^{\ddag,\p}}\in\Gamma(S,\isom(V_{\App},V^\p\ca{A})/K^{\ddag,\p})$), such that, $\psi_{\App}$ is sent to an $\A^{\p,\times}$-multiple of the Weil pairing $\mrm{e}^{\lambda_{\overline{s}}}$ of $\lambda_{\overline{s}}$ via each such $\varepsilon^\ddag$ in the orbit; we will omit the subscript $K^{\ddag,\p}$ when it is clear in the context.
    \end{itemize}
Any pair of such tuples $(\ca{A},\lambda,[\varepsilon^\ddag])$ and $(\ca{A}^\prime,\lambda^\prime,[\varepsilon^{\ddag,\prime}])$ are equivalent in $\mbf{M}^{\mrm{isog}}_{K^\ddag}(S)$, denoted as $(\ca{A},\lambda,[\varepsilon^\ddag])\sim_{\mrm{isog}}(\ca{A}^\prime,\lambda^\prime,[\varepsilon^{\ddag,\prime}])$, if there is a $\zbkppt$-isogeny $f:\ca{A}\to\ca{A}^\prime$ such that $rf^\vee\circ\lambda^\prime\circ f=\lambda$ for some $\bb{Z}_{(\square),>0}^\times$ and such that $V^\p f\circ \varepsilon^\ddag=\varepsilon^{\ddag,\prime}$ modulo the action of $K^{\ddag,\p}$. 
Finally, for general $\zbkpp$-scheme $S$, as the pairs $(\ca{A},\lambda)$ described above are finitely presented, there is a locally Noetherian $\zbkpp$-scheme $S_1$ with a morphism $S\to S_1$ over $\zbkpp$, such that $(\ca{A},\lambda)$ is the pullback of a disjoint union of tuples described above over each connected component of $S_1$. The action of $\pi_1(S,\overline{s})$ for any geometric point $\overline{s}\in S$ on $V^\p\ca{A}$ factors through the action of the {\'e}tale fundamental group $\pi_1(S_1,\overline{s}_1)$ at a geometric point $\overline{s}_1$ of $S_1$ by \cite[Arcata V. Cor. 3.3]{SGA4.5}. Hence, the tuple $(\ca{A},\lambda,[\varepsilon^\ddag])$ is the pullback of a disjoint union of tuples described above over each connected component of $S_1$.\par
\begin{rk}
We will have equivalent interpretations of $\mbf{M}^{\mrm{isog}}_{K^{\ddag,\p}}$ if doing the following changes: (1) We further require that the objects in $\mbf{M}^{\mrm{isog}}_{K^{\ddag,\p}}(S)$ satisfy the (Kottwitz) determinantal condition (see \cite[\S 5]{Kot92} and also \cite[1.3.4]{Lan13}). Since the $\ca{O}$ in \emph{loc. cit.} is $\bb{Z}$, the condition is empty for the Siegel-type case. (2) We can change the term ``$\zbkppt$-polarizations'' to other terms such as ``weak polarizations with respect to $\zbkppt$'' and also allow the $r$ in the equivalence condition above to cover the whole $\zbkppt$. 
\end{rk}
The moduli problem above, defined by $\zbkppt$-isogeny classes of abelian schemes with additional structures, is equivalent to the one defined below using isomorphism classes of abelian schemes with additional structures.\par
Let $(V_{\bb{Z}},\psi_{\bb{Z}})$ be a pair of $\bb{Z}$-lattice $V_\bb{Z}$ with $\psi_\bb{Z}:V_\bb{Z}\times V_\bb{Z}\to \bb{Z}(1)$, the restriction of $\psi$ to $V_\bb{Z}$ which takes its values in $\bb{Z}(1)$. Recall that we require that $V_\bb{Z}$ is self-dual with respect to $\psi$. Define the moduli problem $\mbf{M}^{\mrm{iso}}_{{(V_{\bb{Z}},\psi_\bb{Z})},K^{\ddag,\p}}:
(Sch/\zbkpp)^{\mrm{op}}\to(Grpds)$ of isomorphism classes of abelian schemes with additional structures as follows:\par
For any locally Noetherian and connected $\zbkpp$-scheme $S$, $\mbf{M}_{(V_{\bb{Z}},\psi_\bb{Z}),K^{\ddag,\p}}^{\mrm{iso}}(S)$ is a groupoid whose objects consist of the tuples $(\ca{A},\lambda,[\varepsilon^\ddag_{\zhpp}]_{K^{\ddag,\p}})$, where 
\begin{itemize}
    \item $\ca{A}$ is an abelian scheme over $S$, 
    \item $\lambda$ is a principal polarization, and 
    \item $[\varepsilon^\ddag_{\zhpp}]_{K^{\ddag,\p}}$ is a $\pi_1(S,\overline{s})$-invariant $K^{\ddag,\p}$-orbit of $\wat{\bb{Z}}^\p$-equivariant isomorphisms  $\varepsilon^\ddag_{\zhpp}:V_{\wat{\bb{Z}}^\p}\xrightarrow{\sim}T^\p\ca{A}_{\overline{s}}$ for some (and hence any) geometric point $\overline{s}$ of $S$ (this is equivalent to saying that $[\varepsilon^\ddag_{\zhpp}]_{K^{\ddag,\p}}\in\Gamma(S,\isom(V_{\wat{\bb{Z}}^\p},T^\p\ca{A})/K^{\ddag,\p})$), such that, $\psi_{\wat{\bb{Z}}^\p}$ is sent to a $\wat{\bb{Z}}^{\p,\times}$-multiple of the Weil pairing $\mrm{e}^{\lambda_{\overline{s}}}$ associated with the polarization $\lambda_{\overline{s}}$ via each such $\varepsilon^\ddag_{\zhpp}$; we will omit the subscript $K^{\ddag,\p}$ when it is clear in the context.
    \end{itemize}
Any pair of such tuples $(\ca{A},\lambda,[\varepsilon^\ddag_{\zhpp}])$ and $(\ca{A}^\prime,\lambda^\prime,[\varepsilon^{\ddag,\prime}_{\zhpp}])$ are equivalent in $\mbf{M}^{\mrm{iso}}_{(V_\bb{Z},\psi_\bb{Z}),K^{\ddag,\p}}(S)$, denoted as $(\ca{A},\lambda,[\varepsilon^\ddag_{\zhpp}])\sim_{\mrm{iso}}(\ca{A}^\prime,\lambda^\prime,[\varepsilon^{\ddag,\prime}_{\zhpp}])$, if there is an isomorphism $f:\ca{A}\to\ca{A}^\prime$ such that $f^\vee\circ\lambda^\prime\circ f=\lambda$ and such that $T^\p f\circ \varepsilon^\ddag_{\zhpp}=\varepsilon^{\ddag,\prime}_{\zhpp}$ modulo the action of $K^{\ddag,\p}$. For general $\zbkpp$-scheme $S$, as explained above the objects described above are finitely presented. For any $(\ca{A},\lambda,[\varepsilon^\ddag_{\zhpp}])$ over $S$, there is a locally Noetherian $\zbkpp$-scheme $S_1$ with a morphism $S\to S_1$ over $\zbkpp$, such that $(\ca{A},\lambda,[\varepsilon^\ddag_{\zhpp}])$ is the pullback of a disjoint union of tuples described above over each connected component of $S_1$.\par
\begin{rk}
    Again, it is harmless to do the following modifications to the definition above: (1) Requiring (Kottwitz) determinantal conditions, which are in fact empty for our case. (2) Requiring that one of $\pm\lambda$ to be a polarization, and that $f^\vee\circ\pm\lambda^\prime\circ f=\pm\lambda$ in the definition of the equivalence relation $\sim_{\mrm{iso}}$.
\end{rk}
\begin{rk}\label{rk-change-lattice} The definition of $\mbf{M}^{\mrm{iso}}_{(V_{\bb{Z}},\psi_\bb{Z}),K^{\ddag,\p}}$ depends only on the base change of $(V_\bb{Z},\psi_\bb{Z})$ to $\wat{\bb{Z}}^\p$; that is, if $(V^{(1)}_\bb{Z},\psi^{(1)}_\bb{Z})$ and $(V^{(2)}_\bb{Z},\psi^{(2)}_\bb{Z})$ are symplectic $\bb{Z}$-lattices that are isomorphic over $\wat{\bb{Z}}^\p$ after base change, then the two moduli problems defined by the two symplectic lattices are identical, since only the base change $(V_{\zhpp},\psi_{\zhpp})$ is involved in the definition above. We shall omit the subscript $(V_{\bb{Z}},\psi_\bb{Z})$ if it is clear in the context.\par
In the definition above, the condition that $\lambda$ is a principal polarization follows automatically from the self-duality of $V_{\bb{Z}}$ with respect to $\psi$ and the definition of $\varepsilon^\ddag_{\zhpp}$. In general, if we choose a $V_{\bb{Z}}$ such that $V_\bb{Z}\sbst V_{\bb{Z}}^\vee$, we may require that $\lambda$ is a polarization with degree $d$, where $d:=|V_{\bb{Z}}^\vee/V_{\bb{Z}}|$.
\end{rk}
For any object $(\ca{A},\lambda,[\varepsilon^\ddag_{\zhpp}])$ in $\mbf{M}^{\mrm{iso}}_{K^{\ddag,\p}}(S)$, we can associate an object
$(\ca{A},\lambda,[\varepsilon^\ddag_{\zhpp}\otimes \zbkpp])$ in $\mbf{M}^{\mrm{isog}}_{K^\ddag}(S)$ by viewing $\lambda$ as a $\zbkppt$-polarization and by tensoring $\varepsilon^\ddag_{\zhpp}$ with $\zbkpp$. Then the natural inclusion of equivalence relations from $\mbf{M}^{\mrm{iso}}_{K^{\ddag,\p}}(S)$ to $\mbf{M}^{\mrm{isog}}_{K^{\ddag,\p}}(S)$ induces a natural transformation $t: \mbf{M}^{\mrm{iso}}_{K^{\ddag,\p}}\to\mbf{M}^{\mrm{isog}}_{K^{\ddag,\p}}$, which is in fact an equivalence; see, e.g., \cite[Prop. 1.4.3.4]{Lan13}.\par
Let us recall the construction of the inverse  $t^{-1}:\mbf{M}^{\mrm{isog}}_{K^{\ddag,\p}}(S)\to\mbf{M}^{\mrm{iso}}_{K^{\ddag,\p}}(S)$ following \emph{loc. cit.}: \par 
\begin{construction}\upshape\label{cons-ab}
First, we can assume that $S$ is locally Noetherian and connected over $\zbkpp$. 
Next, for any object $(\ca{A},\lambda,[\varepsilon^\ddag])$ in $\mbf{M}^{\mrm{isog}}_{K^{\ddag,\p}}(S)$ and any representative $\varepsilon^\ddag:V_{\App}\to V^\p\ca{A}_{\overline{s}}$, the image $\varepsilon^\ddag(V_{\zhpp})$ of $V_{\zhpp}$ in $V^\p\ca{A}_{\overline{s}}$ is $\pi_1(S,\overline{s})$-invariant (since for any $\pi\in\pi_1(S,\overline{s})$,  $\pi(\varepsilon^\ddag)=\varepsilon^\ddag\circ k$ for some $k\in K^{\ddag,\p}$ and $K^{\ddag,\p}$ stabilizes $V_{\zhpp}$). Moreover, $\varepsilon^\ddag(V_{\zhpp})$ contains a finite index sub-$\zhpp$-lattice ${N}T^\p\ca{A}_{\overline{s}}$ for some prime-to-$\p$ integer $N>0$. \par
Then there is a homomorphism between abelian schemes $f_1:\ca{A}_{\overline{s}}\to \ca{A}^\prime_{\overline{s}}$ over $\overline{s}$ defined by taking the quotient of $\ca{A}_{\overline{s}}$ by $K_{\overline{s}}:=\frac{1}{N}\varepsilon^\ddag(V_{\zhpp})/T^{\p}\ca{A}_{\overline{s}}\sbst V^{\p}\ca{A}_{\overline{s}}/T^{\p}\ca{A}_{\overline{s}}=\ca{A}_{\overline{s},\mrm{torsion}}$. \par
Since $K_{\overline{s}}$ is $\pi_1(S,\overline{s})$-invariant, it extends to a finite {\'e}tale commutative group scheme $K_S$ over $S$ contained in $\ca{A}[k]$. So $f_1$ is the pullback to $\overline{s}$ of some homomorphism $f:\ca{A}\to \ca{A}'$ defined by taking the quotient of $\ca{A}$ by $K_S$. Hence, there is a $\zbkppt$-isogeny defined by $f^\prime: \ca{A}\xleftarrow{[N]}\ca{A}\xrightarrow{f}\ca{A}^\prime$ such that $V^\p f^\prime\circ\varepsilon^\ddag(V_{\zhpp})=T^\p\ca{A}^\prime$.
Now we have found another object $(\ca{A}^\prime,\lambda^\prime,[\varepsilon^{\ddag,\prime}])$ in $\mbf{M}^{\mrm{isog}}_{K^{\ddag,\p}}(S)$ equivalent to $(\ca{A},\lambda,[\varepsilon^\ddag])$ such that $f^{\prime,\vee}\circ\lambda^\prime\circ f'=\lambda$ and $V^\p f^\prime\circ\varepsilon^\ddag=\varepsilon^{\ddag,\prime}$. 
Finally, since $\A^{\p,\times}=\zbkppt\cdot\wat{\bb{Z}}^{\p,\times}$, and $K^{\ddag,\p}$ stabilizes $V_{\zhpp}$, we can choose a prime-to-$\p$ polarization $\lambda'_0$, which is a $\bb{Z}_{(\p),>0}^\times$-multiple of $\lambda^\prime$, such that $\psi_{\zhpp}$ is sent to a $\wat{\bb{Z}}^{\p,\times}$-multiple of $\mrm{e}^{\lambda_{\overline{s}}}$ via any $\varepsilon^{\ddag,\prime}$ in $[\varepsilon^{\ddag,\prime}]$. Hence, we have found an object $(\ca{A}^\prime,\lambda^\prime_0,[\varepsilon^{\ddag,\prime}])$ in $\mbf{M}^{\mrm{iso}}_{K^{\ddag,\p}}(S)$.
\end{construction}
When $K^{\ddag,\p}$ is neat, the moduli problems $\mbf{M}^{\mrm{isog}}_{K^{\ddag,\p}}$ and $\mbf{M}^{\mrm{iso}}_{(V_\bb{Z},\psi_\bb{Z}),K^{\ddag,\p}}$ are representable by a smooth quasi-projective scheme $\ca{M}_{K^{\ddag,\p}}$ over $\zbkpp$. When $\p$ is empty, $\ca{M}_{K^{\ddag}}$ is (the canonical model of) Siegel Shimura variety \gls{shKddag}$:=\sh_{K^\ddag}(G^\ddag,X^\ddag)$. When $\p=\{p\}$, we denote $\ca{M}_{K^{\ddag,\p}}$ by \gls{SKddag}$:=\ca{S}_{K^\ddag_pK^{\ddag,p}}(G^\ddag,\bb{S}^\pm)$, the integral model of the Siegel Shimura variety $\sh_{K^\ddag}$ associated with the Shimura datum $(G^\ddag,\bb{S}^\pm)$ and the level subgroup $K^\ddag:=K_p^\ddag K^{\ddag,p}$, where $K^\ddag_p$ is hyperspecial. The inverse limit $\ca{S}_{K^\ddag_p}:=\varprojlim_{K^{\ddag,p}} \ca{S}_{K^\ddag_pK^{\ddag,p}}$ is an integral canonical model of $\sh_{K^\ddag_p}:=\varprojlim_{K^{\ddag,p}}\sh_{K^\ddag_pK^{\ddag,p}}$ that satisfies the extension property.
\subsection{1-motives}\label{subsec-1-mot}
\subsubsection{}\label{def-1-motive} We recall the definition of $1$-motives. 
\begin{definition}[{\cite[(10.1.2) and (10.1.10)]{Del74}}; see also {\cite[Sec. 2]{Ray94}} and {\cite[1.1]{Mad19}}]\label{def-1-mtf}
A \textbf{1-motive} \gls{Qca} over a base scheme $S$ is a tuple $\Q=(\underline{Y},\G^{\natural},T,A,\iota,c^\vee)$, where:
\begin{enumerate} 
\item $\underline{Y}$ is an {\'e}tale locally constant finite free $\bb{Z}$-module over $S$, 
\item $\G^{\natural}$ is a semi-abelian scheme over $S$,
\item \label{c} $T$ is a torus over $S$ and $A$ is an abelian scheme over $S$, such that $\G^{\natural} $ is an extension of $A$ by $T$, i.e., there is an exact sequence $0\to T\xrightarrow{\ i\ } \G^{\natural}\xrightarrow{\ \pi\ } A\to 0$,
\item \label{iota}$\iota: \underline{Y}\to \G^{\natural}$ is a homomorphism from $\underline{Y}$ to $\G^{\natural}$, whose composition with the natural quotient $\G^{\natural} \to A$ is $c^\vee:\underline{Y}\to A$.
\end{enumerate}
\end{definition}
Let $\underline{X}:=\mbf{X}^*(T)$ be the sheaf of groups of characters of $T$. Then \ref{c} is equivalent to
\begin{enumerate}
\item[$\mathit{\ref{c}}'$] \textit{There is a homomorphism $c:\underline{X}\to A^\vee$.}
\end{enumerate}
In fact, for any $\chi\in \underline{X}$, the $\bb{G}_m$-torsor produced by the pushout of $\G^{\natural}$ under $-\chi: T\to \bb{G}_m$ corresponds to an invertible sheaf $\ca{O}_\chi$ in $A^\vee:=\underline{\Pic}^0_e(A/S)$; so $\G^\natural$ is associated with a homomorphism $\underline{X}\to \underline{\Pic}_e^0(A/S)$, and this association is bijective (see \cite[Prop. 3.1.5.1]{Lan13}). Note that there is an anti-equivalence $$ \underline{\Pic}_e^0(A/S)\xrightarrow{\ \sim\ }\underline{\Ext}^1_S(A,\bb{G}_m).$$ \par
\begin{rk}
Let $X$ be an $S$-scheme. We will use the same notation to denote a $\bb{G}_m$-extension of $X$ and its corresponding invertible sheaf; we will make the precise meaning clear in the context.
\end{rk}
Let $\ca{P}_A$ be the Poincar{\'e} biextension of $A\times_S A^\vee$ by $\bb{G}_m$.\par
Given $c$ and $c^\vee$, the statement \ref{iota} is equivalent to 
\begin{enumerate}[label=$\mathit{\ref{iota}}'$]
\item\label{eq-4}\textit{There is a trivialization of biextensions $\tau^{-1}: \mbf{1}_{\underline{Y}\times\underline{X}}\to (c^\vee\times c)^*\ca{P}_A$ over $\underline{Y}\times\underline{X}$; or equivalently, a trivialization of biextensions $\tau: \mbf{1}_{\underline{Y}\times\underline{X}}\to (c^\vee\times c)^*\ca{P}_A^{\otimes -1}$, i.e., the inverse of $\tau^{-1}$.}
\end{enumerate}
In fact, after replacing $S$ with an {\'e}tale cover, for any $y\in Y$ and $\chi\in X$, $(c^\vee(y)\times c(\chi))^*\ca{P}_A\iso \ca{O}_\chi|_{c^\vee(y)}$. Then the sections $\{\tau^{-1}(y,\chi):\ca{O}_\chi|_{c^\vee(y)}\to\ca{O}_S\}$ determine and are determined by $\bigoplus_\chi \ca{O}_\chi|_{c^\vee(y)}\xrightarrow{\sum_\chi\tau^{-1}(y,\chi)}\ca{O}_S$, which is a section in $\G^\natural(S)$ that projects to $c^\vee(y)\in A(S)$, since we can write $\G^\natural$ as $\underline{\spec}_A(\bigoplus_{\chi\in X}\ca{O}_\chi)$.\par
Symmetrically, given $c$ and $c^\vee$, $\textit{\ref{iota}}^\prime$ is equivalent to 
\begin{enumerate}[label=$\mathit{\ref{iota}}''$]
\item\label{eq-44} \textit{Let $\G^{\natural,\vee}$ be the semi-abelian scheme over $S$ determined by $c^\vee$ under \cite[Prop. 3.1.5.1]{Lan13}. Then there is a homomorphism $\iota^\vee:\underline{X}\to \G^{\natural,\vee}$. }
\end{enumerate}
Note that the semi-abelian scheme $\G^{\natural,\vee}$ satisfies an exact sequence $0\to T^\vee\xrightarrow{\ i^\vee\ }\G^{\natural,\vee}\xrightarrow{\ \pi^\vee\ }A^\vee\to 0$, which makes it an extension of $A^\vee$ by $T^\vee:=\underline{\Hom} (\underline{Y},\bb{G}_m)$.\par
Hence, a 1-motive $\Q$ determines and is determined by its \textbf{Cartier dual} 
$$\Q^\vee:=(\underline{X},\G^{\natural,\vee}, T^\vee, A^\vee,\iota^\vee,c).$$ 
If there is more than one $1$-motive we are working with, we will use $?_{\Q}^\sharp$, for $?=\underline{Y},\underline{X},\G^{\natural},T,A,c,\iota,\tau$ and for $\sharp=\emptyset, \vee$, to refer to the corresponding objects determined by $\Q$; if the $1$-motive we are referring to is clear in the context, we omit the subscript $\Q$.\par
\begin{rk}
Note that if $S$ is the spectrum of a Noetherian complete local ring, or if $S$ is locally Noetherian and normal, then every torus over $S$ is isotrivial. See \cite[EX. X, 3 and 5]{SGA3}; see also \cite[Rmk. 3.2.5.6]{Lan13}. In these cases, we can require $\underline{Y}$ to be trivialized by a finite {\'e}tale cover in the definition of a $1$-motive $\Q$, and $\underline{X}=\underline{Y}_{\Q^\vee}$ is also trivialized by a finite {\'e}tale cover.
\end{rk}
A \textbf{homomorphism} between 1-motives $f:\Q_1\to \Q_2$ is a pair $f=(f^{et},f^{sab})\in\underline{\Hom}(\underline{Y}_{\Q_1},\underline{Y}_{\Q_2})\times\underline{\Hom}(\G^{\natural}_{\Q_1},\G^{\natural}_{\Q_2})$, such that $\iota_{\Q_2}\circ f^{et}=f^{sab}\circ\iota_{\Q_1}$. The morphism $f^{sab}$ determines and is determined by a pair $(f^{tr},f^{ab})\in \underline{\Hom}(T_{\Q_1},T_{\Q_2})\times\underline{\Hom}(A_{\Q_1},A_{\Q_2})$. Moreover, define the \textbf{Cartier dual $f^\vee:\Q_2^\vee\to \Q_1^\vee$ of $f$} to be $f^\vee:=(f^{tr,\vee},(f^{et,\vee},f^{ab,\vee}))$.\par
\begin{definition}[{\cite[Sec. 1.3]{AB05}; see also \cite[10.1.5]{Del74}}]\label{def-qn}
For any positive integer $n$, define \textbf{$n$-torsion subgroup of $\Q$} to be the quotient sheaf $$\Q[n]:=\frac{\ker ((-\iota+[n]):\underline{Y}\times_S\G^{\natural}\to\G^{\natural})}{\im(([n],\iota):\underline{Y}\to\underline{Y}\times_S\G^{\natural})}.$$ 
\end{definition}
Note that this sheaf $\Q[n]$ is indeed representable by a finite flat commutative group over $S$: To see this, note that there is a canonical exact sequence of fppf sheaves:
\begin{equation}\label{eq-can-exact}
0\to \G^\natural[n]\to \Q[n]\to \ul{Y}/n\ul{Y}\to 0.
\end{equation}
Also, there is a canonical exact sequence of commutative groups over $S$:
\begin{equation}\label{eq-gp-ex}
0\to T[n]\to \G^\natural[n]\to A[n]\to 0.    
\end{equation}
Since finiteness and flatness are fppf-local on the base (see \cite[\href{https://stacks.math.columbia.edu/tag/02LA}{Lem. 02LA}]{stacks-project} and \cite[\href{https://stacks.math.columbia.edu/tag/02L2}{Lem. 02L2}]{stacks-project}), $\G^\natural[n]$ is finite flat over $A[n]$. So $\G^\natural[n]$ is finite flat over $S$. 
In (\ref{eq-can-exact}), the first and the last (non-zero) terms are representable. Since $\G^\natural[n]$ is finite over $S$, $\Q[n]$ is representable by fppf descent of affine morphisms; see 
\cite[Prop. 17.4]{Oor66}, and also \cite[\href{https://stacks.math.columbia.edu/tag/0247}{Lem. 0247}]{stacks-project}; this also implies that $\Q[n]$ is finite flat over $S$ by \cite[\href{https://stacks.math.columbia.edu/tag/02LA}{Lem. 02LA}]{stacks-project} and \cite[\href{https://stacks.math.columbia.edu/tag/02L2}{Lem. 02L2}]{stacks-project} again.
\begin{rk}Note that 
one might also replace $(-\iota+[n])$ with $(\iota+[n])$ and replace $([n],\iota)$ with $([n],-\iota)$ to have a different sign convention for $\Q[n]$. Here we are following the conventions of \cite[10.1.5]{Del74} so that we can define Weil pairings of $1$-motives the same as \cite[10.2.5]{Del74}; moreover, it makes it easy to interpret splittings of the weight filtration $W$ on $\Q[n]$ as liftings $(c_n,c^\vee_n,\iota_n)$ (see Proposition \ref{prop-deg-lvl}, part \ref{bd-lv-3}).
\end{rk}
\begin{rk}\label{rk-qn}
Continuing the previous remark, as in \cite[10.1.5]{Del74}, $\Q$ is represented by a complex $[\underline{Y}\to \G^\natural]\in D^b_{fppf}(S)$, where $\underline{Y}$ is in degree $-1$ and $\G^\natural$ is in degree $0$ and where $D^b_{fppf}(S)$ is the derived category of bounded complexes of fppf-sheaves of abelian groups over $S$. Moreover, by \cite[Prop. 2.3.1]{Ray94}, $\ul{\Hom}_{1-\text{mot}}(\Q_1,\Q_2)=\ul{\Hom}_{D^b_{fppf}(S)}(\Q_1,\Q_2)$. We denote by $\{k\}$ the $k$-shifted complex (since square brackets are occupied by the notation of $n$-torsions). So the multiplication of $\Q\{-1\}$ by $n$ is represented by the following morphism between complexes
\begin{equation}
\begin{tikzcd}
{[\underline{Y}}\arrow[rr,"-\iota"]\arrow[d,"n"]&& {\G^\natural]}\arrow[d,"n"]\\
{[\underline{Y}}\arrow[rr,"-\iota"]&&{\G^\natural]};
\end{tikzcd}
\end{equation}
 the mapping cone of the above morphism is represented by 
$$\ca{C}_n(\Q):=[\underline{Y}\xrightarrow{(n,\iota)} \underline{Y}\times_S\G^\natural\xrightarrow{-\iota+n} \G^\natural],$$
whose degrees are concentrated in $-1,0$ and $1$. Then $\Q[n]=H^0(\ca{C}_n(\Q))\iso \ca{C}_n(\Q).$
If we do not shift $\Q$ by $-1$, we can also formulate $\Q[n]$ as $H^{-1}$ of $[\Q\xrightarrow{n}\Q]$ (see, e.g., \cite[3.1]{Ray94} and \cite[\S1.1.2]{Mad19}).
\end{rk}
\begin{lem}\label{lem-mor-tor}
Let $f:\Q_1\to\Q_2$ be a homomorphism between $1$-motives over $S$. Then $f$ naturally induces a homomorphism, which we abusively denote by $f:\Q_1[n]\to\Q_2[n]$, between $n$-torsion subgroups. 
\end{lem}
\begin{proof}
This is just the consequence of the following commutative diagram:
 \begin{equation}
    \begin{tikzcd}
\ca{C}_n(\Q_1):&\ul{Y}_{\Q_1}\arrow[rrr,"{(n,\iota_{\Q_1})}"]
\arrow[d,"f^{et}"]&&&
{\ul{Y}_{\Q_1}\times_S\G^\natural_{\Q_1}}
\arrow[rr,"{-\iota_{\Q_1}+n}"]
\arrow[d,"{(f^{et},f^{sab})}"]&&
\G^\natural_{\Q_1}\arrow[d,"f^{sab}"] \\
\ca{C}_n(\Q_2):&\ul{Y}_{\Q_2}\arrow[rrr,"({n,\iota_{\Q_2}})"]&&&\ul{Y}_{\Q_2}\times_S\G^\natural_{\Q_2}\arrow[rr,"-\iota_{\Q_2}+n"]&& \G^\natural_{\Q_2}.
    \end{tikzcd}
\end{equation}
\end{proof}
In general,
\begin{definition}\label{def-isog-1-mot}
    A homomorphism $f:\Q_1\to \Q_2$ between $1$-motives is called an \textbf{isogeny} if $f^{et}$ is injective and has finite cokernel, and if $f^{sab}$ is an isogeny, that is, $f^{sab}$ is surjective and finite.\par
    Define the \textbf{kernel} of $f$ to be the quotient sheaf
    $$\ker f:=\frac{\ker ((-\iota_{\Q_2}+f^{sab}):\underline{Y}_{\Q_2}\times_S\G_{\Q_1}^{\natural}\to\G_{\Q_2}^{\natural})}{\im((f^{et},\iota_{\Q_1}):\underline{Y}_{\Q_1}\to\underline{Y}_{\Q_2}\times_S\G^{\natural}_{\Q_1})}.$$
    Equivalently, $\ker f:=H^0(\ca{C}_f(\Q_1,\Q_2))\iso\ca{C}_f(\Q_1,\Q_2)$, where $$\ca{C}_f(\Q_1,\Q_2):=[\underline{Y}_{\Q_1}\xrightarrow{(f^{et},\iota_{\Q_1})} \underline{Y}_{\Q_2}\times_S\G_{\Q_1}^\natural\xrightarrow{-\iota_{\Q_2}+f^{sab}} \G_{\Q_2}^\natural],$$
    the mapping cone of $f\{-1\}$,
whose degrees are concentrated in $-1,0$ and $1$.
\end{definition}
See more definitions about isogenies in Appendix \ref{comm-lem}.\par
As (\ref{eq-can-exact}), there is a canonical exact sequence:
\begin{equation}
    0\to \ker f^{sab}\to \ker f\to \coker f^{et}\to 0.
\end{equation}
By \cite[Prop. 17.4]{Oor66}, \cite[\href{https://stacks.math.columbia.edu/tag/02LA}{Lem. 02LA}]{stacks-project}, \cite[\href{https://stacks.math.columbia.edu/tag/02L2}{Lem. 02L2}]{stacks-project} and \cite[\href{https://stacks.math.columbia.edu/tag/0247}{Lem. 0247}]{stacks-project} again, $\ker f$ is representable by a finite flat commutative group over $S$.
\begin{lem}\label{lem-isog-1-mot}
    Let $f:\Q_1\to \Q_2$ be an isogeny between $1$-motives over $S$. Then there is an integer $N$ and an isogeny $h:\Q_2\to \Q_1$, such that, $h\circ f=[N]_{\Q_1}$ and $f\circ h=[N]_{\Q_2}$. In particular, $\ker f\sbst\Q_1[N]$ and $\ker h\sbst\Q_2[N]$.
\end{lem}
\begin{proof}
    There is an integer $N_1$ and an injective homomorphism $g^{et}:\ul{Y}_{\Q_1}\to \ul{Y}_{\Q_2}$ with finite cokernel such that $g^{et}\circ f^{et}=[N_1]_{\ul{Y}_{\Q_1}}$ and $f^{et}\circ g^{et}=[N_1]_{\ul{Y}_{\Q_2}}$; moreover, there is an integer $N_2$ and an isogeny $g^{sab}:\G^\natural_{\Q_1}\to \G^\natural_{\Q_2}$ such that $g^{sab}\circ f^{sab}=[N_2]_{\G^\natural_{\Q_1}}$ and $f^{sab}\circ g^{sab}=[N_2]_{\G^\natural_{\Q_2}}$. Take $N:=N_1N_2$. Then $h:=([N_2]\circ g^{et},[N_1]\circ g^{sab})$ satisfies the first statement: To see this, we draw the following commutative diagram
    \begin{equation}
        \begin{tikzcd}
&\ul{Y}_{\Q_1}\arrow[r,"\iota_{\Q_1}"]\arrow[dd,"{[N_1]}"]\arrow[dl,"f"]&
\G^\natural_{\Q_1}\arrow[dd,"{[N_1]}"]\arrow[dr,"f"]& \\
\ul{Y}_{\Q_2}\arrow[rrr,dashed,"\iota_{\Q_2}"]\arrow[dd,"{[N_1]}"]\arrow[dr,"g^{et}"]\arrow[dddr,dashed]&&&\G^\natural_{\Q_2}\arrow[dd,"{[N_1]}"]\arrow[dddl,dashed]\\
&\ul{Y}_{\Q_1}\arrow[r,"\iota_{\Q_1}"]\arrow[dd,"{[N_2]}"]\arrow[dl,"f"]&\G^\natural_{\Q_1}\arrow[dr,"f"]\arrow[dd,"{[N_2]}"]&\\
\ul{Y}_{\Q_2}\arrow[rrr,"\iota_{\Q_2}"]\arrow[dd,"{[N_2]}"]&&&\G^\natural_{\Q_2}\arrow[dl,"g^{sab}"]\arrow[dd,"{[N_2]}"]\\
&\ul{Y}_{\Q_1}\arrow[r,"\iota_{\Q_1}"]\arrow[dl,"f"]&\G^\natural_{\Q_1}\arrow[dr,"f"]&\\
\ul{Y}_{\Q_2}&&&\G^\natural_{\Q_2}.
        \end{tikzcd}
    \end{equation}
    In the diagram above, all triangles, parallelograms and trapezoids containing at most one dashed arrow in their edges are commutative. 
    By diagram-chasing, we find $[N_1]\circ g^{sab}\circ\iota_{\Q_2}=\iota_{\Q_1}\circ [N_2]\circ g^{et}$.\par
For the second statement, we only need to show $\ker f\sbst \Q_1[N]$, since the other one can be proved symmetrically. We have the following commutative diagram
\begin{equation}
    \begin{tikzcd}
\ca{C}_f(\Q_1,\Q_2):&\ul{Y}_{\Q_1}\arrow[rrr,"{(f^{et},\iota_{\Q_1})}"]
\arrow[d,equal]&&&
{\ul{Y}_{\Q_2}\times_S\G^\natural_{\Q_1}}
\arrow[rr,"{-\iota_{\Q_2}+f^{sab}}"]
\arrow[d,"{(h^{et},\id)}"]&&
\G^\natural_{\Q_2}\arrow[d,"h^{sab}"] \\
\ca{C}_N(\Q_1):&\ul{Y}_{\Q_1}\arrow[rrr,"({[N],\iota_{\Q_1}})"]&&&\ul{Y}_{\Q_1}\times_S\G^\natural_{\Q_1}\arrow[rr,"-\iota_{\Q_1}+{[N]}"]&& \G^\natural_{\Q_1}.
    \end{tikzcd}
\end{equation}
Since $\ker (-\iota_{\Q_2}+f^{sab})\sbst \ker (-\iota_{\Q_1}+[N])$ and $\im (h^{et}\circ f^{et},\iota_{\Q_1})=\im ([N],\iota_{\Q_1})$, we have the second statement.
\end{proof}
Conversely,
\begin{lem}\label{lem-isog-con}
Let $\Q_1$ and $N$ be defined as in Lemma \ref{lem-isog-1-mot}. Let $K$ be a finite, flat, closed $S$-subgroup of $\Q_1[N]$. \par
Then there is a $1$-motive $\Q_2$ and an isogeny $f:\Q_1\to\Q_2$, such that, $\ker f=K$.
\end{lem}
\begin{proof}
There is a canonical quotient map from $\Q_1[N]$ to $\ul{Y}_{\Q_1}/N\ul{Y}_{\Q_1}$, whose kernel is $\G^\natural_{\Q_1}[N]$. Then $K$ fits into a canonical exact sequence of fppf sheaves
$$0\to K_G\lra K\lra K_Y\lra 0.$$
In the sequence above,  $K_G$ is the kernel of  $K\sbst\Q_1[N]\to \ul{Y}_{\Q_1}/N\ul{Y}_{\Q_1}$, and $K_Y:=K/K_G$ as an fppf sheaf; in fact, $K_Y=K/K_G$ is representable by a scheme since $K_G$ is finite. Since both $K$ and $K_G$ are flat, $K_Y$ is finite and flat in $\underline{Y}_{\Q_1}/N\ul{Y}_{\Q_1}$, and therefore is locally constant.\par
Let $\G^\natural_{\Q_2}:=\G^\natural_{\Q_1}/K_G$ and $\ul{Y}_{\Q_2}$ be the inverse image of $K_Y\sbst \ul{Y}_{\Q_1}/N\ul{Y}_{\Q_1}\iso \frac{1}{N}\ul{Y}_{\Q_1}/\ul{Y}_{\Q_1}$ in $\frac{1}{N}\ul{Y}_{\Q_1}$. 
Let $K_T$ be the kernel of $K_G\sbst \G^\natural_{\Q_1}[N]\to A_{\Q_1}[N]$ and let $K_A:=K_G/K_T$, which is representable by a finite flat group scheme since $K_T$ is finite and flat. There is an exact sequence $0\to K_T\to K_G\to K_A\to 0$. \par
Since $\G^\natural_{\Q_1}$ is an extension of an abelian scheme $A_{\Q_1}$ by a torus $T_{\Q_1}$, and by a theorem of Deligne and Raynaud \cite[Thm. 1.9 (b)]{FC90}, each orbit of $K_G$ is locally on the base contained in an affine open subset of $\G^\natural_{\Q_1}$. Then $\G^\natural_{\Q_2}$ is representable by a scheme. Then $\G^\natural_{\Q_2}$ is an extension of an abelian scheme $A_{\Q_2}:=A_{\Q_1}/K_A$ by a torus $T_{\Q_2}:=T_{\Q_1}/K_T$, and $\ul{Y}_{\Q_2}$ is flat and locally constant.\par
Let $f^{et}:\ul{Y}_{\Q_1}\to \ul{Y}_{\Q_2}$ be the natural inclusion; let $f^{sab}:\G^\natural_{\Q_1}\to \G^\natural_{\Q_2}$ be the natural quotient.\par
Finally, let us define $\iota_{\Q_2}$ such that $f^{sab}\circ\iota_{\Q_1}=\iota_{\Q_2}\circ f^{et}$. Fix any $S^\prime$ over $S$ and any $y\in \ul{Y}_{\Q_2}(S^\prime)$. 
Since $\ul{Y}_{\Q_2}$ is defined to be the fiber product of $K_Y$ and $\frac{1}{N}\ul{Y}_{\Q_1}$ over $\frac{1}{N}\ul{Y}_{\Q_1}/\ul{Y}_{\Q_1}$, we see $y$ projects to some $k_y\in K_Y(S')$. We then lift $k_y$ to $k\in K(S')$ and then to $\wdtd{k}=(\wdtd{y},\wdtd{g})\in \ul{Y}_{\Q_1}\times_S\G_{\Q_1}^\natural(S')$ up to replacing $S'$ with an fppf cover. Since the projection of $\wdtd{y}$ from $\ul{Y}_{\Q_1}(S')$ to $\ul{Y}_{\Q_1}/N\ul{Y}_{\Q_1}(S')\xrightarrow{\frac{1}{N}}\frac{1}{N}\ul{Y}_{\Q_1}/\ul{Y}_{\Q_1}(S')$ coincides with the projection of $y$ to $\frac{1}{N}\ul{Y}_{\Q_1}/\ul{Y}_{\Q_1}(S')$. We have $N y-\wdtd{y}=N y_0$ for $y_0\in \ul{Y}_{\Q_1}(S')$. By adjusting $(N y_0,\iota_{\Q_1}(y_0))$, we have a lifting $(Ny,\wdtd{g}')\in \ul{Y}_{\Q_1}\times_S\G_1^\natural(S')$. We define $\iota_{\Q_2}(y)=\overline{\wdtd{g}}'$, the projection of $\wdtd{g}'$ to $\G_{\Q_2}^\natural(S').$
This definition is independent of the choice of the lifting from $K_Y$ to $K$, since any two liftings differ by an element in $K_G$. Moreover, this definition is independent of the second step of lifting from $K$ to $\ul{Y}_{\Q_1}\times_S\G_{\Q_1}^\natural$, since we will always adjust it by some $(N y_0,\iota_{\Q_1}(y_0))$ so that the first component is $Ny$. Hence, this map is well defined. It is compatible with $\iota_{\Q_1}$ since if $y\in \ul{Y}_{\Q_1}(S')$, the lifting we choose will be $(N y,\iota_{\Q_1}(y))$.
\end{proof}
Moreover, we can define the $p$-adic Tate module by the $p$-adic sheaf $T_p\Q:=\varprojlim_{n\in \bb{Z}_{>0}}\Q[p^n]$ and denote the prime-to-$p$ Tate module by the $\zhp$-sheaf $T^p\Q:=\varprojlim_{n\in\bb{Z}_{>0}\text{ and }p\nmid n}\Q[n]$. Define $\wat{T}\Q:=\varprojlim_{n\in \bb{Z}}\Q[n]$. Define the corresponding rational Tate modules by $V_p\Q:=T_p\Q\otimes\bb{Q}$, $V^p\Q:=T^p\Q\otimes \bb{Q}$ and $\wat{V}\Q:=\wat{T}\Q\otimes \bb{Q}$, respectively.\par
The \textbf{weight filtration $W_\bullet$} on $\Q$ is defined by $W_{-3}\Q=0$, $W_{-2}\Q=T$, $W_{-1}\Q=\G^{\natural}$ and $W_0\Q=\Q$; so $\mrm{Gr}_{-2}^W=T$, $\mrm{Gr}_{-1}^W=A$ and $\gr_0^W=\underline{Y}$. This weight filtration canonically induces weight filtrations on $?\Q$ by setting $W_{i}?\Q:=?W_i\Q$, for $?=T^p$, $T_p$, $V^p$, $V_p$, $\wat{T}$ and $\wat{V}$, and for any $i\in\bb{Z}$; we denote their corresponding graded pieces by $\gr_i^{W_{?\Q}}$. For $n$-torsion subgroups, denote the weight filtration and its graded pieces by $W_i\Q[n]:=(W_i\Q)[n]$ and $\gr_i^{W_{\Q[n]}}$. 
\begin{lem}\label{lem-tate-inj}
Let $\kappa$ be a field. Let $l\neq \chara\kappa$ be a prime number. Let $\Q_1$ and $\Q_2$ be two $1$-motives over $S=\spec \kappa$. Then there is an injective map
\begin{equation}\label{eq-tate-inj}
T_l:\Hom(\Q_1,\Q_2)\hookrightarrow\Hom_{\bb{Z}_l}(T_l\Q_1,T_l\Q_2). 
\end{equation}
\end{lem}
\begin{proof}
Let $f:\Q_1\to \Q_2$ be any homomorphism such that $T_lf=0$. We want to show that $f=0$. Assume that $\kappa$ is algebraically closed. Then, for any integer $n>0$, the morphism $f^{sab}:\G^\natural_{\Q_1}[l^n]\to \G^\natural_{\Q_2}[l^n]$ is $0$. Then the homomorphism $f^{ab}:A_{\Q_1}[l^n]\to A_{\Q_2}[l^n]$ is $0$. So $f^{ab}:A_{\Q_1}\to A_{\Q_2}$ is $0$ since (\ref{eq-tate-inj}) holds for abelian varieties by \cite[Thm. 3, p.176]{Mum74}. Also, (\ref{eq-tate-inj}) holds for tori since tori over an algebraically closed field are copies of $\bb{G}_m$. We then have $f^{tr}=0$ and $f^{sab}=0$.\par
Now we can assume that $f^{sab}=0$. Then $f^{et}(\ul{Y}_{\Q_1})\sbst \ker \iota_{\Q_2}$. So for any $y_1\in \ul{Y}_{\Q_1}(S')$, the projection of $(f^{et}(y_1),0)\in \ul{Y}_{\Q_2}\times_S \G_{\Q_2}^\natural(S')$ to $\ul{Y}_{\Q_2}/l^n\ul{Y}_{\Q_1}(S')$ lies in $\Q_2[l^n]$. But since this projection is trivial for all integers $n>0$, $f^{et}(y_1)$ is trivial. Hence, $f^{et}$ and $f^{sab}$ are both trivial, so is $f$.
\end{proof}
\begin{definition}\label{def-pol} Let $\Q$ be a $1$-motive defined as in Definition \ref{def-1-mtf}. A \textbf{polarization} of a 1-motive $\Q$ is a homomorphism  $\bml:\Q\to\Q^\vee$ such that $\bml^{et}=\phi$ for some injective homomorphism $\phi: \underline{Y}\to\underline{X}$ with finite cokernel, such that $\bml^{tr}=\phi^\vee$, and such that $\bml^{ab}=\lambda_{A}$ for some polarization $\lambda_{A}$ of $A$. \end{definition}
A pair $(\Q,\bml)$ of a 1-motive $\Q$ with a polarization $\bml$ is called a polarized 1-motive. Given a $\bml$, by definition, $c\circ\phi=\bml^{ab}\circ c^\vee$. 
Then we have a trivialization $\varphi:=(\mrm{Id}\times\phi)^*\tau:\mbf{1}_{\underline{Y}\times\underline{Y}}\to (\mrm{Id}\times\phi)^* (c^\vee\times c)^*\ca{P}_A^{\otimes{-1}}\iso (c^\vee\times c^\vee)^*(\mrm{Id}\times\lambda_A)^*\ca{P}_A^{\otimes{-1}}$ of $\bb{G}_m$-biextensions over $\underline{Y}\times\underline{Y}$. Note that by the rigidifications of $\ca{P}_A$ along $e_A\times A^\vee$ and $A\times e_{A^\vee}$, and by partial multiplication laws of biextensions (see \cite[VII, D{\'e}f. 2.1]{SGA7I}), if we replace $S$ with some {\'e}tale cover such that $\underline{X}=X$ and $\underline{Y}=Y$ are constant, then $\tau(Y,0)=\tau(0,X)=1$, and $\tau(-,-)$ (resp. $\varphi(-,-)$) is bilinear over $Y\times X$ (resp. $Y\times Y$).
\subsubsection{}
Next, we recall the definition of biextension of $1$-motives by commutative groups. Using this notion, we recall the definition of Weil pairings for $1$-motives following Deligne \cite[10.2.5]{Del74}.\par 
\begin{definition}[{\cite[10.2.1]{Del74}}]
Let $\Q_1$ and $\Q_2$ be two $1$-motives over $S$. Let $H$ be a commutative group scheme over $S$. A \textbf{biextension $\ca{P}$ of $(\Q_1,\Q_2)$ by $H$} is a biextension $\ca{P}$ of $(\G^\natural_1,\G^\natural_2)$ by $H$ such that 
\begin{enumerate}
\item There is a trivialization $\varrho_1$ (resp. $\varrho_2$) of biextensions $\ca{P}|_{\G^\natural_1\times \underline{Y}_2}$ (resp. $\ca{P}|_{\underline{Y}_1\times\G^\natural_2}$) $\varrho_1:\mbf{1}_{\G^\natural_1\times \underline{Y}_2}\to\ca{P}|_{\G^\natural_1\times \underline{Y}_2}^{\otimes{-1}}$ (resp. $\varrho_2:\mbf{1}_{\underline{Y}_1\times\G^\natural_2}\to\ca{P}|_{\underline{Y}_1\times\G^\natural_2}^{\otimes{-1}}$).
\item $\varrho_1 $ and $\varrho_2$ coincide when they are pulled back to $\underline{Y}_1\times \underline{Y}_2$ under $(\iota_1\times \mrm{Id}_{\underline{Y}_2})$ and $(\mrm{Id}_{\underline{Y}_1}\times\iota_2)$, respectively.
\end{enumerate}
\end{definition}
\begin{rk}Recall that a $1$-motive $\Q$ can also be interpreted as a complex of group schemes $\Q=[\underline{Y}\xrightarrow{\ \iota\ }\G^\natural]$; here, we let $\underline{Y}$ be in degree $-1$ and let $\G^\natural$ be in degree $0$. As remarked in \cite[10.2.1]{Del74}, $\mrm{Biext}^1(\Q_1,\Q_2;H)$, the isomorphism classes of biextensions of $(\Q_1,\Q_2)$ by $H$, is isomorphic to $\underline{\Ext}^1(\Q_1\otimes^L\Q_2,H)$; by hom-tensor adjunction \cite[\href{https://stacks.math.columbia.edu/tag/0A65}{Lem. 0A65}]{stacks-project}, this group is isomorphic to $\underline{\Ext}^1(\Q_1,R\underline{\Hom}(\Q_2,H)).$ 
\end{rk}
\begin{lem}\label{lem-tri}
Let $\ca{P}_\Q$ be the pullback of $\ca{P}_A$ to $\G^\natural\times\G^{\natural,\vee}$. Then $\ca{P}_\Q$ can be canonically viewed as a biextension of $(\Q,\Q^\vee)$ by $\bb{G}_m$.
\end{lem}
\begin{proof}
Let $S^\prime$ be any $S$-scheme such that $\underline{X}(S^\prime)=X$ and $\underline{Y}(S^\prime)=Y$. Giving a section $x\in \G^\natural(S^\prime)$ is equivalent to giving a section $x_A\in A(S^\prime)$ together with a section $x_f$ of the fiber $\G^\natural|_{x_A}$ over $S^\prime$; a section $x_f$ of the trivial $T$-torsor $\G^\natural|_{x_A}$ is determined by a map $\chi\in X(S^\prime)\mapsto x_f(\chi)\in \ca{O}_\chi|_{x_A}(S')$ (here $\ca{O}_\chi$ means the corresponding $\bb{G}_m$-bundle), a section of the pushout of $\G^\natural|_{x_A}$ along $-\chi$. We show that $\varrho_1(x,\chi):=x_f(\chi)^{-1}$ is the desired trivialization of the biextension $\ca{P}_\Q^{\otimes -1}|_{\G^\natural\times \underline{X}}$.\par 
To see this, we have to show that $\varrho_1(-,-)$ can be canonically viewed as a trivialization of $\ca{P}_\Q|_{(x,\chi)}$ for any $x\in\G^\natural(S^\prime)$ and any $\chi\in X(S^\prime)$, and that $\varrho_1(-,-)$ is bilinear. We let $S=S^\prime$. 
Note that any $x\in \G^\natural(S)$ induces a section $x\in \G^\natural\times_A\G^\natural(S)$ under the diagonal embedding $\Delta:\G^\natural\to\G^\natural\times_A\G^\natural$. Then $\varrho_1(x,\chi)$ is defined by (the inverse of) the pushout of $\G^\natural\times_A\G^\natural|_{(x,x)}$ from the second factor along $-\chi$. The case for $\varrho_2(-,-)$ is constructed symmetrically.\par
Next, let us explain the bilinearity. The linearity with respect to the first factor follows from the multiplicative structure of $\G^\natural$ and the diagonal embedding; the linearity of the second factor follows from the multiplicative structure of pushouts with respect to the tensor product: $\ca{O}_\chi\otimes\ca{O}_{\chi^\prime}\xrightarrow{\sim}\ca{O}_{\chi+\chi^\prime}$.\par
Finally, $\varrho_1$ and $\varrho_2$ coincide when pulled back to $Y\times X$: this can be seen from the equivalence of $\tau$, $\iota$ and $\iota^\vee$; see \ref{def-1-motive}, \ref{iota}, \ref{eq-4} and \ref{eq-44}. 
\end{proof}
We will keep on using the notation $\varrho_1(g,h)$ (resp. $\varrho_2(g,h)$) defined in Lemma \ref{lem-tri}, for $g\in \G^\natural(S^\prime)$ and $h\in\iota^\vee(\ul{X})(S^\prime)$ (resp. for $g\in \iota(\ul{Y})(S^\prime)$ and $h\in \G^{\natural,\vee}$), where $S^\prime$ is an $S$-scheme; if there is more than one $1$-motive involved, we will use $\varrho_1^\Q$ and $\varrho_2^\Q$ to specify that the trivializations $\varrho_1$ and $\varrho_2$ are defined for $\ca{P}_\Q$. If $g$ (resp. $h$) is trivial, we denote $\mrm{rig}_h:=\varrho_2(0,h)$ (resp. $\mrm{rig}_g:=\varrho_1(g,0)$).\par
Now we recall Deligne's definition of Weil pairings.\par
For any two valued points $x\in\G^\natural(S^\prime)$ and $y\in\G^{\natural,\vee}(S^\prime)$ for some $S$-scheme $S^\prime$, we denote $\ca{P}_\Q|_{(x,y)}$ by $\ca{P}_{x,y}$. Moreover, we denote $(x\times\id_{\G^{\natural,\vee}})^*\ca{P}_\Q$ (resp. $(\id_{\G^\natural}\times y)^*\ca{P}_\Q$) by $\ca{O}_x$ (resp. $\ca{O}_y$). 
\begin{definition}[{\cite[10.2.5]{Del74}}]\label{def-wp}
The \textbf{Weil pairing for a $1$-motive $\Q$} over $S$ is a perfect pairing ${\e}_{\Q[n]}: \Q[n]\times \Q^\vee[n]\to \bb{G}_m$, which is defined as the difference of two trivializations described as follows:\par
Let $q_1\in \Q[n](S')$ and $q_2\in\Q^\vee[n](S')$. Suppose that $q_1$ lifts to $(y,g)\in \underline{Y}(S^\prime)\times \G^\natural(S^\prime)$ and $q_2$ lifts to $(x,h)\in \underline{X}(S^\prime)\times \G^{\natural,\vee}(S^\prime)$. Then $\mrm{e}_{\Q[n]}(q_1,q_2)$ fits into the commutative diagram
\begin{equation}
\begin{tikzcd}
&\ca{P}_{\iota(y),h}\arrow[rr,"{\varrho_2(y,h)}"]&& \bb{G}_m\\
\ca{P}^{\otimes n}_{g,h}\arrow[ur,"\mrm{can.}"]\arrow[dr,"\mrm{can.}"]&&  &\\
&\ca{P}_{g,\iota^\vee(x)}\arrow[rr,"{\varrho_1(g,x)}"]&&\bb{G}_m\arrow[uu,"{\mrm{e}_{\Q[n]}(q_1,q_2)}"].
\end{tikzcd}
\end{equation}
Symbolically, we write ${\e}_{\Q[n]}(q_1,q_2):=\frac{\varrho_2(y,h)}{\varrho_1(g,x)}$.\par
Moreover, for any polarization $\bml$ of $\Q$, $\bml(\Q[n])\sbst \Q^\vee[n]$ (see Lemma \ref{lem-mor-tor}). Define ${\e}^{\bml}: \Q[n]\times \Q[n]\to \bb{G}_m$ to be ${\e}^{\bml}(q_1,q_2):={\e}_{\Q[n]}(q_1,\bml(q_2))$ for $q_1,q_2\in \Q[n]$.
\end{definition}
The following lemma is immediate:
\begin{lem}\label{lem-wp-ele}
    Let $(\Q,\bml)$ be a polarized $1$-motive and let $f:\Q_1\to\Q$ be an isogeny between $1$-motives. Then $\bml_1:=f^\vee\circ\bml\circ f$ is a polarization of $\Q_1$. \par 
    Let $q_1\in \Q_1[n]$ and $q_2\in \Q^\vee[n]$. Set $q_1=(y,g)\in \ul{Y}_{\Q_1}(S^\prime)\times\G^\natural_{\Q_1} (S^\prime)$ and $q_2=(x,h)\in \ul{X}_\Q(S^\prime)\times \G^{\natural,\vee}_\Q(S^\prime)$. Then $\varrho_1^\Q(f^{sab}(g),x)=\varrho_1^{\Q_1}(g,f^{\vee,et}(x))$, and 
    $\varrho_2^{\Q}(f^{et}(y),h)=\varrho_2^{\Q_1}(y,f^{\vee,sab}(h))$.\par 
    Consequently, we have the following equalities: $\e_{\Q[n]}(f(q_1),q_2)=\e_{\Q_1[n]}(q_1,f^\vee(q_2))$; if $q_1,q_2\in\Q_1[n](S^\prime)$, $\e^{\bml}(f(q_1),f(q_2))=\e^{\bml_1}(q_1,q_2)$; for any $y_1\in \ul{Y}_{\Q_1}(S^\prime)$ and any $x\in \ul{X}_{\Q}(S^\prime)$, $\tau_\Q(f^{\et}(y),x)=\tau_{\Q_1}(y,f^{\vee,et}(x))$.
\end{lem}
\begin{proof}
The first paragraph follows directly from the definition. The second paragraph follows from the construction of $\varrho_1$ and $\varrho_2$ in Lemma \ref{lem-tri}, as the pushout of $g\in\G^\natural_{\Q_1}(S^\prime)$ under $-f^{\vee,et}(x)$ factors through the pushout of $f^{sab}(g)\in\G^\natural_\Q(S^\prime)$ under $-x$, and as the pushout of $f^{\vee,sab}(h)\in\G^{\natural,\vee}_{\Q_1}(S^\prime)$ under $-y$ factors through the pushout of $h\in\G^{\natural,\vee}_\Q(S^\prime)$ under $-f^{et}(y)$. The third paragraph follows from the second paragraph.
\end{proof}
\subsubsection{}\label{subsubsec-biext-rel}
The statements in \S\ref{subsubsec-biext-rel} and \S\ref{subsubsec-n-tor-alt} are not used in the remaining part of this paper. These two subsections are retained to ensure consistency with the submitted version of the thesis.\par
We make some explanations on how the definitions above are compatible with the general theory of Ext- and Biext- groups studied in \cite[VII and VIII]{SGA7I}.\par
First, let us consider the complex $[\underline{Y}\xrightarrow{\ c^\vee\ }A]$ concentrated in degrees $-1$ and $0$. Then, as in \cite[VII, 3.1 and Thm. 3.2.5]{SGA7I}, the group $\underline{\Ext}^1([\ul{Y}\to A],\bb{G}_m)$ can be described as follows:
\begin{thm}[{\cite[VII, Thm. 3.2.5]{SGA7I}}] Let $L=[\cdots\rightarrow L_2\xrightarrow{d_2} L_1\xrightarrow{d_1} L_0\rightarrow 0]$ be a complex of fppf sheaves of commutative groups over $S$ such that $L_0$ is in degree $0$. Let $H$ be a commutative group over $S$. Then
$$\ul{\Ext}^1(L,H)\iso\{\begin{array}{c} (E,\alpha)| E \text{ is an extension of }L_0 \text{ by }H,\ \alpha \text{ is a trivialization}\\ \text{of the pullback extension }d_1^*(E)\text{ such that }\alpha\circ d_2=0\end{array}\}/\iso,
$$
where an isomorphism $f:(E,\alpha)\to (E^\prime,\alpha^\prime)$ is an isomorphism of extensions of $L_0$ by $H$, $f:E \to E^\prime$, such that $f\circ\alpha=\alpha^\prime$. Moreover, $\ul{\Ext}^0(L,H)\iso\ul{\Hom}_{S-gp}(H^0(L),H)$.
\end{thm}
As there is an exact sequence 
$$0\to [0\rightarrow A]\to [\ul{Y}\rightarrow A]\to [\ul{Y}\rightarrow 0]\to 0,$$
applying $R\ul{\Hom}(-,\bb{G}_m)$, we see that there is an exact sequence
$$0\to T^\vee\to \ul{\Ext}^1([\ul{Y}\to A],\bb{G}_m)\to A^\vee\to 0,$$
and $R^0\ul{\Hom}([\ul{Y}\to A],\bb{G}_m)=0$.\par
On the other hand, since the arguments in \cite[VIII, 1.1.4 and 1.4.2]{SGA7I} involve only general homological algebra and also work in the case of $1$-motives, and since $R^0\underline{\Hom}([\ul{Y}\to A],\bb{G}_m)$ is trivial, we have an isomorphism $\biext^1(\ca{T},[\underline{Y}\xrightarrow{c^\vee}A];\bb{G}_m)\xrightarrow{\sim}\underline{\Hom}(\ca{T},\underline{\Ext}^1([\underline{Y}\xrightarrow{c^\vee}A],\bb{G}_m))$ for any fppf sheaf of commutative groups $\ca{T}$. Then the functor $\biext^1(-,[\ul{Y}\to A];\bb{G}_m)$ from $Sh(S_{fppf})^{\mrm{op}}$ to $Sh(S_{fppf})$ is represented by $\ul{\Ext}^1([\ul{Y}\to A],\bb{G}_m)$; the fppf sheaf $\ul{\Ext}^1([\ul{Y}\to A],\bb{G}_m)$ is representable by a commutative group scheme by fppf descent of affine morphisms (see \cite[Prop. 17.4]{Oor66}; see also \cite[\href{https://stacks.math.columbia.edu/tag/0247}{Lem. 0247}]{stacks-project}).\par
More precisely,
\begin{prop}\label{prop-dual-g} $\ul{\Ext}^1([\ul{Y}\to A],\bb{G}_m)$ is representable by $\G^{\natural,\vee}$. In particular, if we assume that $\ul{Y}$ splits, for simplicity, then $\ul{\Ext}^1([\ul{Y}\to A],\bb{G}_m)\iso \bigoplus_{y\in Y}\ca{O}_{c^\vee(y)}$, where $\ca{O}_{c^\vee(y)}:=(c^\vee(y)\times\id_{A^\vee})^*\ca{P}_A$. Similar statements are true if we replace $c^\vee, \ul{X}$ and $A$ with $c$, $\ul{Y}$ and $A^\vee$, respectively.
\end{prop}
\begin{proof}
Let $\ca{P}$ be the $\bb{G}_m$-extension over $A\times\G^{\natural,\vee}$ defined by the pullback of $\ca{P}_A$ under $\id_A\times \pi^\vee$. Then as Lemma \ref{lem-tri}, the diagonal embedding $\Delta:\G^{\natural,\vee}\to \G^{\natural,\vee}\times_{A^\vee}\G^{\natural,\vee}$ defines a canonical trivialization of $\mbf{t}^\Delta:\mbf{1}_{\ul{Y}\times \G^{\natural,\vee}}\to(c^\vee\times\id_{\G^{\natural,\vee}})^*\ca{P}$ over $\ul{Y}\times \G^{\natural,\vee}$ by pushing out sections $-y\in \ul{Y}(S^\prime)$ for any $S$-scheme $S^\prime$; we claim that $(\ca{P},\mbf{t}^\Delta)$ is the tautological biextension (with the sign of trivializations opposite to the definition appearing previously) of the functor $\biext^1(-,[\ul{Y}\to A];\bb{G}_m)$.\par
On the one hand, for any commutative group scheme $H$ over $S$, any $x\in \G^{\natural,\vee}(H)$, the pair $((\id_A\times x)^*\ca{P},(\id_A\times x)^*\mbf{t}^\Delta)$ defines a desired biextension in $\biext^1(H,[\ul{Y}\xrightarrow{c^\vee}A];\bb{G}_m)$. On the other hand, for any biextension $(E,\alpha)$ of $([\ul{Y}\to A],H)$ by $\bb{G}_m$, the image of $(E,\alpha)$ from the top-left corner to the bottom-right corner of the commutative diagram
\begin{equation}
\begin{tikzcd}
\biext^1(H,[\underline{Y}\xrightarrow{c^\vee}A];\bb{G}_m)\arrow[rr,"\sim"]\arrow[d,"\mrm{can.}"]&&\underline{\Hom}(H,\underline{\Ext}^1([\underline{Y}\xrightarrow{c^\vee}A],\bb{G}_m))\arrow[d,"\mrm{can.}"]\\
\biext^1(H,A;\bb{G}_m)\arrow[rr,"\sim"] &&\ul{\Hom}(H,A^\vee)
\end{tikzcd}
\end{equation}
determines a homomorphism $f_E: H\to A^\vee$, and this induces an isomorphism $E\iso (\id_A\times f_E)^*\ca{P}_A$. Then $\alpha$ is a trivialization that assigns each $y\in \ul{Y}(S^\prime)$ a trivialization $\alpha(y)$ of $f_E^*\ca{O}_{c^\vee(y)}(S')$ for any $S$-scheme $S^\prime$, and this is equivalent to a homomorphism $a: H\to H\times_{A^\vee}\G^{\natural,\vee}$. This proves the claim.
\end{proof}
With the definition above, we can then construct $\iota^\vee$ from $\iota$ by pushing out $\G^\natural$ along $-\chi$. By doing this, there is a $\bb{G}_m$-extension $\ca{O}_{\chi}$ with a section $\ul{Y}\to \G^\natural\xrightarrow{-\chi}\ca{O}_{\chi}$ for any $\chi$, which induces a map $\iota^\vee: \ul{X}\to \G^{\natural,\vee}$ by the proposition above.
\subsubsection{}\label{subsubsec-n-tor-alt}The following lemma is self-explanatory due to the definition of $n$-torsion points.
\begin{lem}[{See also \cite[III., Cor. 7.3]{FC90}}]\label{lem-n-tor-ext}
Recall that there are two extensions induced by the weight filtration of $\Q[n]$
\begin{equation}
\begin{split}
0\lra \G^\natural[n]\lra \Q[n]\lra \ul{Y}/n\ul{Y}\lra 0,\\
0\lra T[n]\lra \G^\natural[n]\lra A[n]\lra 0. 
\end{split}
\end{equation}
The first extension is determined by the image of $\iota$ under the canonical connecting homomorphism $\partial:\uhom(\ul{Y},\G^\natural)\to \ul{\Ext}^1(\ul{Y},\G^\natural[n])=\ul{\Ext}^1(\ul{Y}/n\ul{Y},\G^\natural[n])$ induced by applying $R\uhom(\ul{Y},-)$ to $0\to \G^\natural[n]\to \G^\natural\to \G^\natural\to 0$; the second extension is the dual of the  extension determined by the image of $c$ under the canonical connecting homomorphism $\partial:\uhom(\ul{X},A^\vee)\to \ul{\Ext}^1(\ul{X},A^\vee[n])=\ul{\Ext}^1(\ul{X}/n\ul{X},A^\vee[n])$ induced by applying $R\uhom(\ul{X},-)$ to $0\to A^\vee[n]\to A^\vee\to A^\vee\to 0$.
\end{lem}
The construction of the Cartier dual $\Q^\vee$ (see \cite[10.2.11]{Del74}) can be viewed as the following statement:
\begin{prop}\label{prop-cd}
With the notation and conventions introduced in Remark \ref{rk-qn},
the Cartier dual $\Q^\vee$ is isomorphic to the $1$-motive represented by $(\tau^{\leq 1}R\underline{\Hom}(\Q,\bb{G}_m))\{1\}$, where $\tau^{\leq 1}$ denotes the canonical truncation and where $\{-\}$ denotes the shifting of complexes.
\end{prop}
\begin{proof}
There is an exact sequence 
$$0\lra T\lra \Q\lra \Q^\circ\lra 0,$$
 where $\Q^\circ:=[\underline{Y}\to A]$. 
Applying $R\underline{\Hom}(-,\bb{G}_m)$, there is a long exact sequence
\begin{align*}
0\lra R^{0}\underline{\Hom}(\Q^\circ,\bb{G}_m)\lra R^{0}\underline{\Hom}(\Q,\bb{G}_m)\lra \underline{X}\xrightarrow{\ \delta\ } R^1\underline{\Hom}(\Q^\circ,\bb{G}_m)\lra \\
R^1\underline{\Hom}(\Q,\bb{G}_m)\lra 0.\end{align*}
The first term $R^{0}\underline{\Hom}(\Q^\circ,\bb{G}_m)$ is trivial by the second paragraph before Proposition \ref{prop-dual-g}; moreover, the fourth term $R^1\underline{\Hom}(\Q^\circ,\bb{G}_m)=\underline{\Ext}^1([\underline{Y}\xrightarrow{c^\vee}A],\bb{G}_m)\iso \G^{\natural,\vee}$ by Proposition \ref{prop-dual-g}. Moreover, $\delta=-\iota^\vee$ since $\delta$ is constructed by pushing out.\par
There is a distinguished triangle
$$R\ul{\Hom}(\Q^\circ,\bb{G}_m)\lra R\ul{\Hom}(\Q,\bb{G}_m)\lra R\ul{\Hom}(T,\bb{G}_m)\xrightarrow[+1]{\delta};$$
this triangle shifts to a distinguished triangle
$$R\ul{\Hom}(T,\bb{G}_m)\{-1\}\xrightarrow{\ -\delta\{-1\}\ } R\ul{\Hom}(\Q^\circ,\bb{G}_m)\lra R\ul{\Hom}(\Q,\bb{G}_m)\xrightarrow{+1}.$$
From Lemma \ref{lem-cd-hom-alg}, since $R^1\ul{\Hom}(T,\bb{G}_m)=\ul{\Ext}^1(T,\bb{G}_m)$ is trivial, there is a distinguished triangle
$$\tau^{\leq 1}(R\ul{\Hom}(T,\bb{G}_m))\xrightarrow{\ -\delta\ } (\tau^{\leq 1}R\ul{\Hom}(\Q^\circ,\bb{G}_m))\{1\}\lra (\tau^{\leq 1}R\ul{\Hom}(\Q,\bb{G}_m))\{1\}\xrightarrow{+1}.$$
Since $\tau^{\leq 1}R\ul{\Hom}(T,\bb{G}_m)\iso \ul{X}$ and $(\tau^{\leq 1}R\ul{\Hom}(\Q^\circ,\bb{G}_m))\{1\}\iso \G^{\natural,\vee}$, $(\tau^{\leq 1}R\ul{\Hom}(\Q,\bb{G}_m))\{1\}\iso [\ul{X}\xrightarrow{-\delta\ }\G^{\natural,\vee}]$, which is concentrated in degrees $-1$ and $0$. 
\end{proof}
\begin{lem}\label{lem-cd-hom-alg}
Let $f:A^\bullet\to B^\bullet$ be a morphism between complexes of fppf sheaves of commutative groups over $S$. Assume that $A^i=B^i=0$ for any $i<0$, and assume that $H^2(A^\bullet)=0$. Let $\ca{C}_f$ be the mapping cone of $f$. Then $\tau^{\leq 1}\ca{C}_f$ is quasi-isomorphic to the mapping cone of ${\tau^{\leq 1}}f:\tau^{\leq 1}A^\bullet\to \tau^{\leq 1}B^\bullet$.
\end{lem}
\begin{proof}
Write $A^\bullet=[0\to A^0\to A^1\to A^2\to\cdots]$ and $B^\bullet=[0\to B^0\to B^1\to B^2\to \cdots]$. Let $K^1:=\im(d_A^1:A^1\to A^2)$ and $R^1:=\im(d_B^1:B^1\to B^2)$. By assumption, $K^1=\ker (d^2_A:A^2\to A^3)$. Then the mapping cone $\ca{C}_f$ of $f$ can be written as $(\ca{C}_f^i,d_C^i)$, where $\ca{C}_f^i=A^{i+1}\oplus B^i$ and $d_C^i=(-d_A^{i+1},f^{i+1}+d_B^i)$; the first four terms can be written as $[0\to A^0\to A^1\oplus B^0\to A^2\oplus B^1\to A^3\oplus B^2\to \cdots]$, where $A^0$ is concentrated in degree $-1$.\par
Then $\tau^{\leq 1}\ca{C}_f=[0\to A^0\to A^1\oplus B^0\to A^2\oplus B^1\to \im(d_C^1:\ca{C}_f^1\to \ca{C}_f^2)\to 0]$. On the other hand, the mapping cone of $[0\to A^0\to A^1\to K^1\to 0]\to [0\to B^0\to B^1\to R^1\to 0]$ is $\ca{C}_{\tau^{\leq 1} f}:=[0\to A^0\to A^1\oplus B^0\to K^1\oplus B^1\to 0\oplus R^1\to 0]$. There is an obvious morphism $\omega:\ca{C}_{\tau^{\leq 1}f}\to \tau^{\leq 1}\ca{C}_f$ induced by inclusions $K^1\to A^2$ and $R^1\to \im (d_C^1:\ca{C}_f^1\to \ca{C}_f^2)$. We can check that $\omega$ is a quasi-isomorphism: We need to check the isomorphism of $H^1$. If $d_C^1(a_2,b_1)=0$ for some $a_2\in A^2$ and $b_1\in B^1$, then $d_A^2(a_2)=0$, so $a_2\in K^1$ by assumption. Since the image of $d_C^0$ factors through $K^1\oplus B^1$, we deduce that $H^1(\ca{C}_{\tau^{\leq 1}f})=H^1(\tau^{\leq1}\ca{C}_f)$; for the same reason, we can also see that $H^0(\ca{C}_{\tau^{\leq 1}f})=H^0(\tau^{\leq 1}\ca{C}_f)$.
\end{proof}
Finally, the duality of $n$-torsion subgroups of $1$-motives can also be described as follows:
\begin{prop}\label{prop-dual-tor}
$\Q^{\vee}[n]\iso \ul{\Hom}_{1-mot}(\Q[n],\bb{G}_m)$. 
\end{prop}
\begin{proof}
We have the following two diagrams:
\begin{equation}
    \begin{tikzcd}
        \ca{Q}^\vee\{-1\}\arrow[r,"n"]\arrow[d,equal]&
        \ca{Q}^\vee \{-1\}\arrow[r]\arrow[d]&
        \ca{C}_n(\ca{Q}^\vee)\arrow[r,"+1"]\arrow[d,dashed]&{}\\
    \ca{Q}^\vee\{-1\}\arrow[r]& R\ul{\Hom}(\ca{Q},\bb{G}_m)\arrow[r]\arrow[d]& \ca{Y}\arrow[r,"+1"]\arrow[dl,dashed]&{}\\
    & \tau^{>1}R\ul{\Hom}(\ca{Q},\bb{G}_m)\arrow[d,"+1"]\arrow[uur,dashed,"+1"']& &,\\
   {}& {}&&
    \end{tikzcd}
\end{equation}
and
\begin{equation}
    \begin{tikzcd}
        \ca{Q}^\vee\{-1\}\arrow[r]\arrow[d,equal]&
        R\ul{\Hom}(\ca{Q},\bb{G}_m)\arrow[r]\arrow[d,"n"]&
        \tau^{>1}R\ul{\Hom}(\ca{Q},\bb{G}_m)\arrow[r,"+1"]\arrow[d,dashed]&{}\\
    \ca{Q}^\vee\{-1\}\arrow[r]& R\ul{\Hom}(\ca{Q},\bb{G}_m)\arrow[r]\arrow[d]& \ca{Y}\arrow[r,"+1"]\arrow[dl,dashed]&{}\\
    & \wdtd{\ca{C}}_n\arrow[d,"+1"]\arrow[uur,dashed,"+1"']& &.\\
    {}& {}&&
    \end{tikzcd}
\end{equation}
From the Octahedral Axiom (TR4; see \cite[\href{https://stacks.math.columbia.edu/tag/0145}{Def. 0145}]{stacks-project}), we find that the two dashed triangles are distinguished. Then $\Q^\vee[n]\iso H^0(\ca{C}_n(\Q^\vee))\iso H^0(\ca{Y})\iso H^0(\wdtd{\ca{C}}_n)$.
On the other hand, by \cite[Prop. 2.3.1]{Ray94}, $\ul{\Hom}_{1-mot}(\Q[n],\bb{G}_m)\iso R^0\ul{\Hom}(\Q[n],\bb{G}_m)$, and $R^0\ul{\Hom}(\Q[n],\bb{G}_m)\iso H^0(R\ul{\Hom}(\ca{C}_n(\Q),\bb{G}_m))\iso H^0(Cone(R\ul{\Hom}(\Q,\bb{G}_m)\xrightarrow{n}R\ul{\Hom}(\Q,\bb{G}_m)))\iso H^0(\wdtd{\ca{C}}_n)$. So we have the desired result.
\end{proof}
\subsection{Theory of degeneration}\label{subsec-deg-data}
We review the notions of degenerating families, degeneration data, and their level structures. Along the way, we will explain that the definitions of level structures associated with the degeneration data in \cite[Sec. 6.2.3]{Lan13} and in \cite[Sec. 2.2]{Mad19} are equivalent (see Proposition \ref{prop-deg-lvl}). The aim of \S\ref{subsubsec-deg-data-1}-\S\ref{subsubsec-deg-data-4} is to provide a very brief summary; see \cite{FC90} and \cite{Lan13} for all the details that we omit.
\subsubsection{}\label{subsubsec-deg-data-1}
Let $R$ be a Noetherian normal domain that is complete with respect to a radical ideal $I\sbst R$. Let $K:=\Frac{R}$ be the fraction field. Denote $R_i:=R/I^{i+1}$, $\W:=\spec R$, $\W_i:=\spec R_i$ and $\eta:=\spec K$. Denote $\W_\fo:=\spf (R,I)$.
\begin{definition}[{\cite[Def. 4.4.2]{Lan13}}] Define the category $\mbf{DEG}_{\pol}(R,I)$ as follows:
\begin{enumerate}
\item An object of $\mbf{DEG}_{\pol}(R,I)$ is a pair $(\G,\lambda)$ of a semi-abelian scheme $\G$ over $\W$ with a homomorphism $\lambda:\G\to \G^\vee$. 
\item The generic fiber $\G_\eta$ is an abelian scheme, and $\lambda_\eta:\G_\eta\to\G_\eta^\vee$ is a polarization of the abelian scheme $\G_\eta$. The pullback $\G_0:=\G\times_\W \W_0$ is an extension of an abelian scheme $A_0$ by an isotrivial torus $T_0$ over $\W_0$.
\item A morphism $f:(\G_1,\lambda_1)\to (\G_2,\lambda_2)$ in $\mbf{DEG}_{\pol}(R,I)$ is a homomorphism $f$ from $\G_1$ to $\G_2$ respecting the given polarizations over $\eta$, i.e., $f^\vee_\eta\circ\lambda_{2,\eta}\circ f_\eta=\lambda_{1,\eta}$.
\end{enumerate}
\end{definition}
Let $(V_\bb{Z},\psi_\bb{Z})$ be as in \S\ref{subsec-siegel}. Let $S$ be a normal locally Noetherian scheme, and let $U$ be a dense open subscheme of $S$.
\begin{definition}[{\cite[Def. 5.3.2.1]{Lan13}}]
Define the groupoid $\mbf{DEG}_{(V_\bb{Z},\psi_\bb{Z}),K^{\ddag,\p}}(S,U)$ of \textbf{degenerating families of type $\mbf{M}^{\mrm{iso}}_{(V_\bb{Z},\psi_\bb{Z}),K^{\ddag,\p}}$} as follows:
\begin{enumerate}
\item An object in $\mbf{DEG}_{(V_\bb{Z},\psi_\bb{Z}),K^{\ddag,\p}}(S,U)$ is a tuple $(\G,\lambda,[\varepsilon^{\ddag}_{\zhpp}]_{K^{\ddag,\p}})$, where $\G$ is a semi-abelian scheme over $S$ and $\lambda:\G\to\G^\vee$ is a homomorphism. The pullback $(\G_U,\lambda_U)$ is a principally polarized abelian scheme over $U$ such that $[\varepsilon^{\ddag}_{\zhpp}]_{K^{\ddag,\p}}$ is an integral level structure for $(\G_U,\lambda_U)$ of the moduli problem $\mbf{M}^{\mrm{iso}}_{(V_\bb{Z},\psi_\bb{Z}),K^{\ddag,\p}}$. That is, the pullback $(\G_U,\lambda_U,[\varepsilon^{\ddag}_{\zhpp}]_{K^{\ddag,\p}})$ is an object in $\mbf{M}^{\mrm{iso}}_{(V_\bb{Z},\psi_\bb{Z}),K^{\ddag,\p}}(U)$. 
\item Two tuples $(\G,\lambda,[\varepsilon^{\ddag}_{\zhpp}]_{K^{\ddag,\p}})$ and $(\G^\prime,\lambda^\prime,[\varepsilon^{\ddag,\prime}_{\zhpp}]_{K^{\ddag,\p}})$ are equivalent if and only if $(\G_U,\lambda_U,[\varepsilon^{\ddag}_{\zhpp}]_{K^{\ddag,\p}})\iso (\G^\prime_U,\lambda^\prime_U,[\varepsilon^{\ddag,\prime}_{\zhpp}]_{K^{\ddag,\p}})$ in $\mbf{M}^{\mrm{iso}}_{(V_\bb{Z},\psi_{\bb{Z}}),K^{\ddag,\p}}(U)$.
\end{enumerate}
\end{definition}
\subsubsection{}
\begin{definition}[{\cite[p.212]{Lan13}}]Let $\Q$ be a $1$-motive. Define $\bml:=(\phi,\lambda_A)$ to be a pair consisting of an injective homomorphism $\phi:\ul{Y}\to \ul{X}$ with finite cokernel, and a polarization $\lambda_A:A\to A^\vee$ of $A=\Q^{ab}$, such that $\lambda_A\circ c^\vee=c\circ\phi$. With conventions introduced in Definition \ref{def-pol} and the paragraph below it, 
we say that some trivialization $\tau$ of biextension $(c^\vee\times c)^*\ca{P}_A^{\otimes -1}$ is symmetric (with respect to $\bml$) if, for any {\'e}tale $S$-scheme $S^\prime$ and any $y,y^\prime\in \underline{Y}(S^\prime)$, we have $\varphi(y,y^\prime)=\varphi(y^\prime,y)$, where $\varphi:= (\id\times \phi)^*\tau$. 
\end{definition}
\begin{lem}\label{lem-sym}
Let $\tau=\tau_\Q$ be the trivialization of $\ca{P}_A^{\otimes -1}$ canonically constructed from $\iota_\Q$. Let $\bml$ be a polarization of $\Q$. Then $\tau$ is symmetric with respect to $(\bml^{et},\bml^{ab})$.
\end{lem}
\begin{proof}
This follows from Definition \ref{def-pol} and the last statement of Lemma \ref{lem-wp-ele}. 
\end{proof}
\begin{definition}[{\cite[Def. 4.4.6]{Lan13}}]Define the category $\mbf{DD}_{\pol}(R,I)$ of degeneration data over $(R,I)$ (without level structures) as follows:
\begin{enumerate}
\item An object in $\mbf{DD}_{\pol}(R,I)$ is a polarized 1-motive $(\Q_\eta,\bml_\eta)$ over $\eta$. Moreover, we denote $?^\sharp_{\eta}:=?^\sharp_{\Q_\eta}$ for $?=\underline{Y},\underline{X},\G^{\natural},T,A,c,\iota,\tau$ and for $\sharp=\emptyset, \vee$. Then $\underline{X}_{\eta}$ (resp. $\underline{Y}_{\eta}$) extends to an {\'e}tale locally finite free $\bb{Z}$-module $\underline{X}$ (resp. $\underline{Y}$) over $\W$; $\G^{\natural}_{\eta}$ (resp. $G^{\natural,\vee}_{\eta}$) extends to a semi-abelian scheme $\G^{\natural}$ (resp. $G^{\natural,\vee}$) over $\W$, which is an extension of abelian scheme $A$ (resp. $A^\vee$) extending $A_\eta$ (resp. $A^\vee_\eta$) by a torus $T$ (resp. $T^\vee$) extending $T_\eta$ (resp. $T^\vee_\eta$); $c_{\eta}$ (resp. $c^\vee_{\eta}$) extends to a homomorphism $c$ (resp. $c^\vee$). Consequently, $\bml^{et}_\eta=\phi_\eta$ and $\bml_\eta^{ab}=\lambda_{A,\eta}$ canonically extend to homomorphisms $\bml^{et}=\phi:\underline{Y}\to\underline{X}$ and $\bml^{ab}=\lambda_A:A\to A^\vee$.
\item The trivialization $\tau_\eta$ satisfies the \textbf{positivity conditions}, i.e., after replacing $\W$ with an {\'e}tale cover of $\W$ such that $\underline{Y}$ splits, for any $y\in\underline{Y}$, $\tau_\eta(y,\phi(y))$ extends to a section of $(\mrm{Id}\times \phi)^*(c^\vee\times c)^*\ca{P}_A^{\otimes {-1}}$ over $\W$, and this section factors through $I$ if $y$ is nonzero. Moreover, $\tau_\eta$ is \emph{necessarily} symmetric with respect to $(\phi_\eta,\lambda_{A,\eta})$ (see Lemma \ref{lem-sym}).
\item A morphism $f:(\Q_{\eta,1},\bml_{\eta,1})\to(\Q_{\eta,2},\bml_{\eta,2})$ in $\mbf{DD}_{\pol}(R,I)$ is a homomorphism $f:\Q_{\eta, 1}\to\Q_{\eta,2}$ respecting the given polarizations, i.e., $f^\vee\circ\bml_{\eta,2}\circ f=\bml_{\eta,1}$.
\end{enumerate}
\end{definition}
\begin{construction}[{\cite[Ch. II]{FC90} and \cite[Sec. 4.1-4.3]{Lan13}}]\upshape
Let us recall the construction of the functor $\mbf{F}_{\pol}(R,I)$ that associates with any object $(\G,\lambda)$ in $\mbf{DEG}_{\pol}(R,I)$ an object $(\Q_\eta,\bml_\eta)$ in $\mbf{DD}_{\pol}(R,I)$:\par
\begin{enumerate}[label={(\arabic*)}]
\item First, we can write $\G_0:=\G_{\W_0}$ as an extension of an abelian scheme $A_0$ by an isotrivial torus $T_0$ over $\W_0$: $1\to T_0\to\G_0\to A_0\to 1$. Then $T_0$ uniquely extends to a torus $T_i$ over $\W_i$ and $T_i$ induces an extension: $0\to T_i\to\G_i\to A_i\to 0$, where $\G_i:=\G_{\W_i}$ and $A_i$ is an abelian scheme over $\W_i$. Then this induces an extension: $0\to T_{\fo}\to \G_{\fo}\to A_{\fo}\to 0$ over $\W_{\fo}$. Since there is a cubical ample invertible sheaf over $\G$, the extension of $\G_{\fo}$ over $\W_{\fo}$ algebraizes uniquely to a so-called \textbf{Raynaud extension}: $0\to T\xrightarrow{i} \G^\natural\xrightarrow{\pi} A\to 0$ over $\W$. Similarly, there is a canonical Raynaud extension for $\G^\vee$: $0\to T^\vee\xrightarrow{i^\vee} \G^{\vee,\natural}\xrightarrow{\pi^\vee} A^\vee\to 0$. See \cite[Sec. 3.3.3 and 3.4.4]{Lan13}. Thus, $\lambda$ canonically defines two morphisms $\phi^\vee:T\to T^\vee$ and $\lambda_A:A\to A^\vee$, where the former one determines a dual morphism $\phi: \underline{Y}\to \underline{X}$. Then it suffices to construct $\tau_\eta$.
\item Fix an auxiliary ample cubical invertible sheaf $\ca{L}$ over $\G$ which induces some polarization $\lambda'$ (possibly different from $\lambda$).
Then $\ca{L}$ induces an ample cubical invertible sheaf $\ca{L}_{\fo}$ over the formal completion $\G_{\fo}$ and algebraizes to an ample cubical invertible sheaf $\ca{L}^\natural$ over $\G^\natural$. Replacing $\ca{W}$ with a finite {\'e}tale cover, we assume that $T$ and $T^\vee$ split with character groups $\mbf{X}^*(T)=X$ and $\mbf{X}^*(T^\vee)=Y$. 
Then we can fix a cubical trivialization $\mbf{tri}:i^*\ca{L}^\natural\iso \ca{O}_T$ and find an ample invertible sheaf $\ca{M}$ over $A$ such that $\ca{L}^\natural\iso \pi^*\ca{M}$ under this trivialization; moreover, we can write $\pi_*\ca{L}^\natural$ as a direct sum $\pi_*\ca{L}^\natural\iso\bigoplus_{\chi\in X}\ca{M}_\chi$, where $\ca{M}_\chi:=\ca{M}\otimes_{\ca{O}_A} \ca{O}_\chi$, and $\ca{O}_\chi$ is the invertible sheaf corresponding to the $\bb{G}_m$-torsor defined by the pushout of $\G^\natural$ along $-\chi$. By \cite[Prop. 3.1.5.1]{Lan13}, we see that $\ca{O}_\chi\iso (\mrm{Id}_A\times c(\chi))^*\ca{P}_A$. Denote by $p_\chi:\pi_*\ca{L}^\natural\to \ca{M}_\chi$ the projection to the direct factor $\ca{M}_\chi$ of $\pi_*\ca{L}^\natural$ canonically determined by the fixed trivialization $\mbf{tri}$; this projection determines a morphism $\sigma_\chi: \Gamma(\G_\eta,\ca{L}_\eta)\hookrightarrow \Gamma(\G_{\fo,\eta},\ca{L}_{\fo,\eta})\iso \wat{\bigoplus}_{\chi\in X}\Gamma(A_\eta,\ca{M}_{\chi,\eta})\to \Gamma(A_\eta,\ca{M}_{\chi,\eta})$. We then need to show a technical fact that $\sigma_\chi\neq 0$ for all $\chi\in X$ using addition formulas for theta representations (see \cite[Sec. 4.3.3]{Lan13}).
\item There is an isomorphism $T^*_{c^\vee(y)}\ca{M}_\chi\iso \ca{M}_{\chi+\phi(y)}\otimes \ca{M}_\chi(c^\vee(y)).$ Let $K(\ca{M}_{\chi,\eta})$ be the kernel of $\lambda_{\ca{M}_{\chi,\eta}}$, the polarization determined by $\ca{M}_{\chi,\eta}$. Then there is a group $\G(\ca{M}_{\chi,\eta})$, which is an extension of $K(\ca{M}_{\chi,\eta})$ by $\bb{G}_{m,\eta}$, defined by the Weil pairing over $K(\ca{M}_{\chi,\eta})\times K(\ca{M}_{\chi,\eta})$ (see \cite[Sec. 3.2.4]{Lan13}). For any $y\in Y$, there is a canonical isomorphism $\G(\ca{M}_{\chi,\eta})\iso \G(\ca{M}_{\chi+\phi(y),\eta})$ depending only on the trivialization $\mbf{tri}$ (see \cite[Lem. 4.3.2.4]{Lan13}). 
For any $\overline{\chi}\in X/\phi(Y)$, there is a weight-$\overline{\chi}$ space $\Gamma(\G_\eta,\ca{L}_\eta)_{\overline{\chi}}$, whose choice depends only on $\mbf{tri}$, such that $\Gamma(\G_\eta,\ca{L}_\eta)\iso \bigoplus_{\overline{\chi}\in X/\phi(Y)}\Gamma(\G_\eta,\ca{L}_\eta)_{\overline{\chi}}$ and such that $\sigma_\chi$ factors as a sequence $\sigma_\chi: \Gamma(\G_\eta,\ca{L}_\eta)\to \Gamma(\G_\eta,\ca{L}_\eta)_{\overline{\chi}}\xrightarrow{\overline{\sigma}_\chi} \Gamma(A_\eta,\ca{M}_{\chi,\eta})$, where the first map is the natural projection and the second map $\overline{\sigma}_\chi$ is $\G(\ca{M}_{\chi,\eta})$-equivariant. See \cite[Lem. 4.3.2.3 and Lem. 4.3.2.5]{Lan13}. 
In fact, it can be shown that $\overline{\sigma}_\chi$ is an isomorphism of irreducible representations of $\G(\ca{M}_{\chi,\eta})$ (see \cite[Lem. 4.3.2.9]{Lan13}). 
Thus, by comparing $T^*_{c^\vee(y)}\circ \overline{\sigma}_\chi:\Gamma(\G_\eta,\ca{L}_\eta)_{\overline{\chi}}\to \Gamma(A_\eta,\ca{M}_{\chi+\phi(y),\eta})\otimes \ca{M}_\chi(c^\vee(y))_\eta$ and $\overline{\sigma}_{\chi+\phi(y)}:\Gamma(\G_\eta,\ca{L}_\eta)_{\overline{\chi}}\to\Gamma(A_\eta,\ca{M}_{\chi+\phi(y),\eta})$ as isomorphisms of \emph{nonzero} and  \emph{irreducible} representations of $\G(\ca{M}_{\chi,\eta})$, there is a well-defined function $\psii(-,-)\in \ca{M}_\chi(c^\vee(y))_\eta^{\otimes -1}$ such that $\psii(\chi,y)(T^*_{c^\vee(y)}\circ \overline{\sigma}_\chi)=\overline{\sigma}_{\chi+\phi(y)}$ (see \cite[Prop. 4.3.2.10]{Lan13}).
Let $\psii(y):=\psii(0,y)$ and let $\tau_\eta(\chi,y)=\psii(y)^{-1}\psii(\chi,y)$. Then it can be shown that $\tau_\eta$ is symmetric with repsect to $\phi$ and satisfies the positivity conditions (see \cite[Prop. 4.3.1.9]{Lan13}).
\item One can show that the definition of $\tau_\eta$ does not depend on auxiliary choices, that is, it does not depend on the choice of $\ca{L}$ and the trivialization $\mbf{tri}$. Eventually, it can be shown that $\tau_\eta$ does not depend on $\lambda'$. See \cite[Cor. 4.3.4.4 and Prop. 4.5.5.1]{Lan13}.
\end{enumerate}
\end{construction}
\begin{thm}[{\cite[III. Cor. 7.2]{FC90}; see also \cite[Thm. 4.4.16]{Lan13}}]\label{thm-fc}
The functor $$\mbf{F}_{\pol}(R,I):\mbf{DEG}_{\pol}(R,I)\lra\mbf{DD}_{\pol}(R,I)$$ sending $(\G,\lambda)$ to the pair $(\Q_\eta,\bml_\eta)$ constructed as above is an equivalence of categories.
\end{thm}
\begin{rk}The quasi-inverse of the functor $\mbf{F}_\pol(R,I)$ is the so-called ``Mumford's construction'' $\mbf{M}_\pol(R,I)$. Note that the independence of $\tau_\eta$ on the choice of $\lambda$ is actually proved after the construction of $\mbf{M}_\pol(R,I)$; see \cite[Sec. 4.5.5]{Lan13}.
\end{rk}
\subsubsection{}
Fix a pair $(R,I)$ as before. 
Let us briefly recall the construction of the quasi-inverse of the functor $\mbf{F}_\pol(R,I)$, the so-called \textbf{Mumford's construction} $\mbf{M}_\pol(R,I)$; see \cite{Mum72}, \cite[Ch. III]{FC90} and \cite[Sec. 4.5]{Lan13}. \par
Let $\mbf{DEG}(R,I)$ (resp. $\mbf{DD}(R,I)$) be the category whose objects are $\G$ (resp. $\Q_\eta$) that can be extended to an object $(\G,\lambda)\in \ob\mbf{DEG}_\pol(R,I)$ (resp. $(\Q_\eta,\bml_\eta)\in\ob\mbf{DD}_\pol(R,I)$). The morphisms in $\mbf{DEG}(R,I)$ (resp. $\mbf{DD}(R,I)$) are the morphisms between semi-abelian schemes (resp. $1$-motives) that can be extended to $\mor (\mbf{DEG}(R,I))$ (resp. $\mor (\mbf{DD}(R,I))$).
If we forget polarizations, it turns out that there is also a functor $\mbf{F}(R,I):\mbf{DEG}(R,I)\to\mbf{DD}(R,I)$ that is compatible with $\mbf{F}_\pol$, see below.
\begin{construction}[{\cite{Mum72}, \cite[Ch. III]{FC90} and \cite[Sec. 4.5]{Lan13}}]\label{const-mumford}\upshape
Fix some $(\Q_\eta,\bml_\eta)\in \ob \mbf{DD}_\pol(R,I)$. The so-called ``Mumford's construction'' is established in the following steps:
\begin{enumerate}[label={(\arabic*)},start=1]
\item Suppose that we are given an ample invertible sheaf $\ca{M}$ over $A$. Let $\ca{L}^\natural:=\pi^{*}\ca{M}$, the pullback of $\ca{M}$ via $\pi:\G^\natural\to A$. Suppose that there is a cubical trivialization $\psii:\ul{Y}\to \ca{L}^{\natural,\otimes -1}$ and an injective $\phi:\ul{Y}\to \ul{X}$ with finite cokernel. Furthermore, we can assume that $\ul{Y}=Y$ and $\ul{X}=X$, i.e., $\ul{Y}$ and $\ul{X}$ split. The idea is to construct some ``compactification'' $(\pp^\natural,\ca{L}^\natural)$.
\item The first fact is that, \emph{under some technical condition} on $(\ca{L}^\natural,\psii,\phi)$ (see \cite[Ch. III, p. 62]{FC90}, there is a \emph{relatively complete model} $\pp^\natural$. To be precise, $\pp^\natural$ is an integral and locally of finite type scheme over $A$, which contains $\G^\natural$ as an open dense subscheme and satisfies the following properties (see \cite[Ch. III, Def. 3.1]{FC90}):
\begin{enumerate}
\item The invertible sheaf $\ca{L}^\natural$ extends to an invertible sheaf over $\pp^\natural$, which we also denote by $\ca{L}^\natural$. Moreover, $\ca{L}^\natural|_V$ is ample over any quasi-compact open subscheme $V$ of $\pp^\natural$. 
\item The translation action $T:\G^\natural\times\G^\natural\to \G^\natural$ extends to an action of $\G^\natural$ on $\pp^\natural$, i.e., a $\G^\natural$-action $T:\G^\natural\times \pp^\natural\to \pp^\natural$. Denote by $T_g$ the translation action of $g\in\G^\natural(S^\prime)$ for any $S^\prime$ over $\W$. Moreover, the translation action on $\ca{O}_{\G^\natural}$ extends to an action $T^*_g\ca{N}\to \ca{N}$ on $\ca{N}:=\ca{L}^\natural\otimes_{\ca{O}_{\pp^\natural}}\pi^*\ca{M}^{\otimes -1}$.
\item We assume that $\ul{Y}=Y$ splits. Recall that there is a $Y$-action $(\iota,\psii)$ on $(\G^\natural_\eta,\ca{L}^{\natural}_\eta)$. Then there is a $Y$-action on $(\pp^\natural,\ca{L}^\natural)$ extending $(\iota,\psii)$; we denote by $(S,\wdtd{S})$ this extended action, and by $(S_y,\wdtd{S}_y)$ the evaluation of this action at any point $y\in Y$.
\item There is a $\G^\natural$-invariant open subscheme $U$ of $\pp^\natural$, such that $\pp^\natural=\bigcup\limits_{y\in Y}S_y(U)$.
\item $\pp^\natural$ satisfies the \emph{completeness condition} as in \cite[Ch. III, 3, p.61]{FC90}.
    \end{enumerate}
Note that not all $(\ca{L}^\natural,\psii,\phi)$ satisfy the \emph{technical condition} on \cite[Ch. III, p. 62]{FC90}, but for any $\Q_\eta$, one can find some $(\ca{L}^\natural,\psii,\phi)$ satisfying this condition.
\item Denote $\pp^\natural_i:=\pp^\natural\times_\W\W_i$. Every irreducible component of $\pp^\natural_i$ is proper. Moreover, $S_y$ acts on each $\pp_i^\natural$, and the quotient fpqc sheaf $\pp_i:=\pp_i^\natural/Y$ is representable by a projective scheme. The morphism $\pp_i^\natural\to \pp_i$ is {\'e}tale and surjective. The invertible sheaf $\ca{L}^\natural_i:=\ca{L}^\natural\otimes_R R_i$ descends to an ample sheaf $\ca{L}_i$. For different $i$, such pairs $(\pp_i,\ca{L}_i)$ are compatible, and therefore there is a pair $(\pp_{\fo},\ca{L}_\fo)$ consisting of the formal completion $\pp_\fo$ of $\{\pp_i\}$ and an ample sheaf $\ca{L}_\fo$ on it. This pair algebraizes to a pair $(\pp,\ca{L})$ of projective scheme $\pp$ with an ample invertible sheaf $\ca{L}$ on it.
\item Let $\G^{\natural,*}:=\bigcup\limits_{y\in Y}S_y(\G^\natural)\sbst \pp^\natural$. Then $\G^{\natural,*}$ is an open subscheme of $\pp^\natural$. Denote $\G_i:=\G^{\natural,*}\times_\W\W_i$. Consider the complement $C^\natural:=(\pp^\natural- \G^{\natural,*})_{\red}$. Note that $C_i$, the quotient by $Y$ of the pullback of $C^\natural$ to $\W_i$, are reduced closed subschemes of $\pp_i$. We can form a $C_\fo$ over $\W$ whose pullback to $\W_i$ is $C_i$, and this $C_\fo$ algebraizes to a closed subscheme $C$ of $\pp$. Let $\G:=\pp-C$.
Then it can be shown that $\G$ is a semi-abelian scheme, with $\G_\eta$ an abelian scheme, and with $\ca{L}:=\ca{L}|_\G$ an ample invertible sheaf. See \cite[Cor. 4.5.3.9 and Cor. 4.5.3.13]{Lan13}.
\item It can be shown that $\G$ does not depend on the auxiliary choice of $(\ca{L}^\natural,\psii,\phi)$ and $\pp^\natural$. Thus, one can form a functor $\mbf{M}(R,I):\mbf{DD}(R,I)\to \mbf{DEG}(R,I)$, which is an equivalence of category. The quasi-inverse of it is denoted by $\mbf{F}(R,I)$. Moreover, we have $\mbf{M}(R,I)(\Q_\eta^\vee)=\G^\vee$. Thus, there is a $\lambda:\G\to \G^\vee\in \mor(\mbf{DEG}(R,I))$ that is defined by applying $\mbf{M}(R,I)$ to $\bml_\eta$. In fact, it can be shown that $\lambda$ is a polarization, so we can form the functor $\mbf{M}_\pol:\mbf{DD}_\pol(R,I)\to \mbf{F}_\pol(R,I)$; $\mbf{M}_\pol(R,I)$ is the quasi-inverse of $\mbf{F}_\pol$.
\end{enumerate}
\end{construction}
Let us also recall the construction of relative complete models for the graph of a morphism $f:\Q_{\eta,1}\to \Q_{\eta,2}$ in $\mbf{DD}(R,I)$. See \cite[Thm. 4.6]{Mum72}, \cite[Ch. III, Thm. 5.5]{FC90} and \cite[Thm. 4.5.3.6]{Lan13}.\par
Choose  $(\ca{M}_i,\ca{L}^\natural_i,\psii_i,\phi_i)$ satisfying the condition on \cite[Ch. III, p.63]{FC90} as in Construction \ref{const-mumford} for $\Q_{\eta,1}$ and $\Q_{\eta,2}$, respectively. Then there are two relative complete models $\pp^\natural_1$ and $\pp^\natural_2$ for $\Q_{\eta,1}$ and $\Q_{\eta,2}$, respectively.
$f^{sab}$ extends to a morphism $f^{sab}:\G^\natural_1\to \G^\natural_2$ over $\W$. If we let $W_1^\nt$ be the schematic closure of the graph $H^\nt$ of $f^{sab}$ in $\pp_1^\nt\times_\W\pp_2^\nt$. 
Then $W_\fo^\nt:=\bigcup\limits_{(y_1,y_2)\in Y_1\times Y_2} S_{(y_1,y_2)}W_{1,\fo}^\nt$ is a locally finite union. 
Then it can be shown that $W_\fo:=W_\fo^\nt/Y_1\times Y_2$ is a formal subscheme of $P_{1,\fo}\times_{\W_\fo} P_{2,\fo}$; moreover, $W_\fo$ algebraizes to a unique $W\sbst P_1\times_\W P_2$. Then it can be shown that $H=W\cap(\G_1\times_\W \G_2)$ is the graph of a uniquely determined $g:G_1\to G_2$, and we let $g=\mbf{M}(f)$.
\subsubsection{}\label{subsubsec-deg-data-4} The $n$-torsion points of $\G_\eta$ and $\Q_\eta$ are isomorphic. More precisely, they are canonically isomorphic under $\mbf{M}$.
\begin{lem}\label{lem-mum-tor}
Suppose that $n$ is prime to $p$ and that $\mbf{M}(\Q_\eta)=\G$. 
Then there is a canonical isomorphism $\Q_\eta[n]\iso \G_\eta[n]$ induced by $\mbf{M}$.
\end{lem}
\begin{proof}
We can see this in the proof of \cite[Thm. 4.10]{Mum72}, \cite[Thm. 5.9]{FC90} and \cite[Thm. 4.5.3.10]{Lan13}. Let us repeat this for the convenience of the reader.\par
By {\'e}tale descent, we assume that the locally constant group part of all $1$-motives we are working on are constant, i.e., $\ul{X}=X$ and $\ul{Y}=Y$. 
Let $\pp^\natural$ be a relative complete model for $\Q_\eta$. 
Let $\zn^\natural_\eta$ be $Y\times_{\iota,\G^\natural_\eta,[n]}\G^\natural_\eta$. Let $\zn^\natural_{y,\eta}$ be $y\times _{\iota,\G^\natural_\eta,[n]}\G^\natural_\eta$. We see that $\zn_{y,\eta}^\natural\iso\zn_{y+nz,\eta}^\natural$ for any $z\in Y$. Then, if we choose any set of representatives of $Y/nY$ in $Y$, $\disju_{y\in Y/nY}\zn^\natural_{y,\eta}$ is isomorphic to $\Q_\eta[n]$. Let $\sigma_y:\W\to\G^{\natural,*}$ be the unique section that $\sigma_y(\eta)=\iota(y)$. 
Let $\zn_y^\natural:=\W\times_{\sigma_y,\G^{\natural,*},[n]}\G^{\natural,*}$. Since $\G^{\natural,*}_\eta=\bigcup_y\iota(y)(\G^\natural_\eta)=\G^\natural_\eta$, $\zn^\natural_{y,\eta}$ is dense in $\zn^\natural_y$.\par
Denote by $\wn^\nt_1$ the schematic closure of the graph of $[n]$ in $\pp^\nt\times \pp^\nt$. Recall that $\iota$ extends to an action $S$ over $\pp^\nt$. Denote $S_{(a,b)}:=S_a\times S_b$. Denote $\wn^\nt:=\bigcup\limits_{y\in Y}S_{(0,y)}\wn_1^\nt$.\par
Then 
$$\bigcup\limits_{y\in Y}(\overline{\zn^\nt_y})_\fo=\wn^\nt_\fo\cap(\pp_\fo^\nt\times_{\W_\fo}\sigma_0\W_\fo),$$
and both the left-hand side and the right-hand side are locally finite unions.\par
Let $(\overline{\zn})_\fo:=\bigcup\limits_{y\in Y}(\overline{\zn^\nt_y})_\fo/Y=\wn_\fo\cap(\pp_\fo\times_{\W_\fo}\sigma_0\W_\fo),$ where $\wn_\fo:=\wn_\fo^\nt/Y\times Y$. Then
$(\overline{\zn})_\fo$ algebraizes to a closed subscheme $\overline{\zn}$ of $P$, and 
$$\overline{\zn}=\wn\cap(P\times\sigma_0\W).$$
On the other hand, for any finite subset $Y_0\sbst Y$, there is a morphism
$$\overline{q}_\fo:\bigcup\limits_{y\in Y_0}(\overline{\zn^\nt_y})_\fo\lra (\overline{\zn})_\fo.$$
Since $\overline{q}_\fo$ is a morphism between formal completions of proper schemes over $\W$, $\overline{q}_\fo$ algebraizes to a morphism $\overline{q}$.
Since $\overline{q}_\fo^{-1}(\overline{\zn}_\fo\cap\pp_\fo\times C_\fo)\sbst \pp^\nt_\fo\times C^\nt_\fo$, 
$\overline{q}^{-1}(\overline{\zn}\cap\pp\times C)\sbst \pp^\nt\times C^\nt$.
Then $\overline{q}$ restricts to 
$$q:\bigcup\limits_{y\in Y_0}(\zn^\nt_y)\lra (\zn).$$
Choose $Y_0$ to be any set of representatives of $Y/nY$ in $Y$. Then it turns out that $q$ is {\'e}tale, surjective and of degree one; see the proof of \emph{loc. cit}. In particular, $q$ is an isomorphism. Hence, $q_\eta$ is the desired isomorphism induced by Mumford's construction.
\end{proof}
\subsubsection{}\label{subsec-deg-data-lvl}
Next, we recall the definition of level structures of degeneration data, which can be used to describe the boundary mixed Shimura varieties of $\ca{S}_{K^\ddag}$. Let $(V_\bb{Z},\psi_\bb{Z})$ be a pair as in \S\ref{subsec-siegel}, where $V_\bb{Z}$ is a $\bb{Z}$-lattice and $\psi_\bb{Z}:V_\bb{Z}\times V_\bb{Z}\to \bb{Z}(1)$ is the restriction of $\psi$ to $V_\bb{Z}$, which takes its values in $\bb{Z}(1)$. Moreover, we can and we will assume that $V_\bb{Z}$ is self-dual with respect to $\psi$. Set $\p=\emptyset$ or $\{p\}$.\par 
Fix any cusp label representative $\Phi_0=
(Q,X^+_0,g_0)$ of $(G_0,X_0)$. Let $g_0^\p\in G_0(\App)\sbst G^\ddag(\App)$ be the projection of $g_0$ to $\App$-points of $G_0$. 
Suppose that $\Phi_0$ maps to a cusp label representative $\Phi^\ddag=(Q^\ddag,X^{\ddag,+},g^\ddag)$, where $Q^\ddag$ is the unique minimal admissible $\bb{Q}$-parabolic of $G^\ddag$ containing $Q$, $X^{\ddag,+}$ is the unique connected component of $X^\ddag$ containing $X^+_0$, and $g^\ddag=g_0$. 
For any $K^{\ddag}$ containing $K$, recall that $K_{\Phi^\ddag}:=P_{\Phi^\ddag}(\A)\cap g^\ddag K^\ddag (g^{\ddag})^{-1}$. We assume $K^\ddag$ is as in \S\ref{subsec-siegel}. \par
Recall some conventions in Appendix \ref{cpr-pel-cl}.  The parabolic subgroup $Q^\ddag$ defines a filtration $W^\ddag_\bullet=\{W^\ddag_i\}_{i=-2}^0$ on $V_\bb{Q}$, and the filtration $W^\ddag_\bullet$ determines an admissible and fully symplectic filtration $Z^{(g^{\ddag})}:=V_{\zhpp}\cap g^{\ddag,-1}(W^\ddag_i)\otimes_\bb{Q}\App$ on 
$V_{\zhpp}$.
Denote $V^{(g^\ddag)}_\bb{Z}:=g^\ddag V_{\wat{\bb{Z}}}\cap V_\bb{Z}$. Let $F_\bullet^{(g^\ddag)}:=V_\bb{Z}^{(g^\ddag)}\cap W_\bullet^\ddag$ be the filtration on $V^{(g^\ddag)}_\bb{Z}$. Recall that $\nu$ is the similitude character for $G^\ddag$. For any $g^\ddag\in G^\ddag(\A)$, there is a unique decomposition $\nu(g^\ddag)=r(g^\ddag)\cdot h(g^\ddag)$, where \gls{rgddag}$\in \bb{Q}_{>0}^\times$ and \gls{hgddag}$\in \wat{\bb{Z}}^\times$.
\begin{definition}[{see \cite[2.2.3]{Mad19}}]\label{def-deg-lvl}
Let $S$ be a $\zbkpp$-scheme. We fix a cusp label representative $\Phi^\ddag$. Suppose that there is a polarized $1$-motive $(\Q,\bml)$ over $S$. 
A \textbf{level structure of type $(V_\bb{Z},\psi_\bb{Z},\Phi^\ddag, K^{\ddag,\p})$} for $(\Q,\bml)$ is 
a section $[u]_{(K_{\Phi^\ddag}^\p)^{g^\ddag}}\in \Gamma(S,\isom(V_{\wat{\bb{Z}}^\p},T^\p\Q)/\rcj{(K^\p_{\Phi^\ddag})}{g^\ddag})$, which, over any geometric point $\overline{s}\in S$, is a $\pi_1(S,\overline{s})$-invariant $\rcj{(K^\p_{\Phi^\ddag})}{g^\ddag}$-orbit of an isomorphism of $\zhpp$-modules $u_{\overline{s}}:V_{\wat{\bb{Z}}^\p}\xrightarrow{\sim}T^\p\Q_{\overline{s}} $ with an isomorphism $v_{\overline{s}}: \zhpp(1)\xrightarrow{\sim}T^\p\bb{G}_m$ such that the following conditions are satisfied:
\begin{enumerate} 
\item $u_{\overline{s}}$ maps $v_{\overline{s}}\circ\psi_{\zhpp}$ to $\e^{\bml_{\overline{s}}}$, the Weil pairing $\mrm{e}^{\bml_{\overline{s}}}$ over $T^\p\Q_{\overline{s}}\times T^\p\Q_{\overline{s}}$ induced by $\bml_{\overline{s}}$.
\item $u_{\overline{s}}$ maps $Z^{(g^\ddag)}$ to the weight filtration $W_\Q$ on $T^\p\Q_{\overline{s}}$. 
\item There are isomorphisms over $\bb{Z}$, $\alpha_{-2}^{(g^\ddag)}:\gr^{F^{(g^\ddag)}}_{-2}\to \Hom (\ul{X},\bb{Z}(1))$ and $\alpha_0^{(g^\ddag)}:\gr^{F^{(g^\ddag)}}_0\to \ul{Y}$, such that, $v_{\overline{s}}\circ(\alpha_{-2}^{(g^\ddag)}\otimes 1)\circ h(g^\ddag)^{-1} \circ\gr(g^\ddag)= \gr_{-2}u_{\overline{s}}$ and $(\alpha_{0}^{(g^\ddag)}\otimes 1)\circ \gr(g^\ddag)= \gr_0u_{\overline{s}}$; the isomorphisms $\alpha_{-2}^{(g^\ddag)}\otimes 1$ and $\alpha_{0}^{(g^\ddag)}\otimes 1$ are isomorphisms defined by tensoring $\zhpp$.
\end{enumerate}
\end{definition}

We have the following statement relating Definition \ref{def-deg-lvl} with \cite[Ch. 5]{Lan13}:
\begin{prop}[{cf. \cite[Ch. 5]{Lan13}}]\label{prop-deg-lvl} Let $(\Q,\bml)$ be a polarized $1$-motive over a $\zbkpp$-scheme $S$. Fix a cusp label representative $\Phi^\ddag$ as above. We denote $\Q=(\underline{Y},\G^{\natural},T,A,\iota,c^\vee)$, denote $\Q^\vee=(\underline{X},\G^{\natural,\vee},T^\vee,A^\vee,\iota^\vee,c)$ and denote $\tau_\Q$ by $\tau$, for simplicity. Then the following objects associated with $(\Q,\bml)$ are equivalent:
\begin{enumerate}
\item\label{bd-lv-1} We associate a level structure $[u]_{(K^\p_{\Phi^\ddag})^{g^\ddag}}$ of type $(V_\bb{Z},\psi_\bb{Z},\Phi^\ddag,K^{\ddag,\p})$.
\item\label{bd-lv-2} Under the map $\mrm{CL}^\p$ defined as in Appendix \ref{cpr-pel-cl}, $\Phi^\ddag$ defines a PEL cusp label which admits a representative $(Z^{(g^{\ddag})}_{K^{\ddag,\p}},\Phi_{K^{\ddag,\p}}^{(g^{\ddag})},\delta^{(g^{\ddag})}_{K^{\ddag,\p}})$. In particular, $g^\ddag$ determines a $K^{\ddag,\p}$-orbit of $(\varphi_{-2}^{(g^\ddag)},\varphi_0^{(g^\ddag)})$.
We associate a tuple $(\alpha_0^\prime, \alpha_{-2}^\prime, v_{K^{\ddag,\p}}, \beta_{-1,K^{\ddag,\p}}^{(g^{\ddag})},\hat{\varsigma})$ defined as follows:
\begin{itemize}
\item $\alpha_{-2}^\prime: \Hom_{\wat{\bb{Z}}^\p}(\ul{X}\otimes\wat{\bb{Z}}^\p,T^\p \bb{G}_m)\xrightarrow{\sim}
\Hom_{\wat{\bb{Z}}^\p}(\ul{X}\otimes\wat{\bb{Z}}^\p,T^\p\bb{G}_m)$ is an automorphism of $\Hom_{\wat{\bb{Z}}^\p}(\ul{X}\otimes \wat{\bb{Z}}^\p,T^\p\bb{G}_m)$ that is induced by an automorphism of $\ul{X}$.
\item $\alpha_0^\prime: \gr^{W_{T^\p\Q}}_0\xrightarrow{\sim}\gr^{W_{T^\p\Q}}_0$ is an automorphism of $\gr^{W_{T^\p\Q}}_0\iso \ul{Y}\otimes \zhpp$ that is a $\wat{\bb{Z}}^{\p,\times}$-multiple of a base change of an automorphism of $\ul{Y}$ to $\zhpp$.
\item $v_{K^{\ddag,\p}}$ is a $K^{\ddag,\p}$-orbit of isomorphisms $v:\zhpp(1)\xrightarrow{\sim} T^\p\bb{G}_m$.
\item Set $\beta_0^{(g^\ddag)}:=\alpha^\prime_0\circ\varphi_0^{(g^\ddag)}$ and $\beta_{-2}^{(g^\ddag)}:=\alpha^\prime_{-2}\circ v \circ\varphi_{-2}^{(g^\ddag)}$. Moreover, $\beta_{-1,K^{\ddag,\p}}^{(g^{\ddag})}$ is a $K^{\ddag,\p}$-orbit of $\beta_{-1}^{(g^{\ddag})}:\mrm{Gr}^{Z^{(g^{\ddag})}}_{-1}\xrightarrow{\sim} A$. 
\item $\hat{\varsigma}: \bigoplus\limits_{i=0}^2\mrm{Gr}^{W_{T^\p\Q}}_{-i}\xrightarrow{\sim}T^\p\Q$ is a splitting of $W_{T^\p\Q}$, such that, 
\begin{enumerate}[label=(\ref{bd-lv-2}.\arabic*)]
\item\label{2-*}$(\bigoplus\limits_{i=0}^{2}\beta_{-i}^{(g^\ddag)})_* (\delta^{(g^{\ddag})})^{*} (v\circ\psi_{\zhpp})=\hat{\varsigma}^* (\mrm{e}^{\bml})$.
\end{enumerate}
\end{itemize}
\item\label{bd-lv-3} Under the conventions of the second part, we associate a tuple $(\alpha_0',\alpha_{-2}',\beta_{-1,K^{\ddag,p}}^{(g^{\ddag})},c_{K^{\ddag,\p}},c^\vee_{K^{\ddag,\p}},\tau_{K^{\ddag,\p}})$, where:
\begin{enumerate}[label=(\ref{bd-lv-3}.\arabic*)]
\item$\beta_{-1,K^{\ddag,\p}}^{(g^{\ddag})}$ is defined as in the second part, such that, $\beta_{-1}^{(g^{\ddag})}$ maps $\psi_{11}$, the restriction of $\psi_{\zhpp}$ to $\gr_{-1}^{Z^{(g^{\ddag})}}$ is a $\zhppt$-multiple of $\mrm{e}^{\lambda_{A}}$; that is, there is an isomorphism $v:\zhpp(1)\xrightarrow{\sim}T^\p\bb{G}_m$ such that $\beta^{(g^{\ddag}),*}_{-1}(\mrm{e}^{\lambda_{A}})=v\circ \psi_{11}$.
\item $\alpha_0^\prime$ and $\alpha^\prime_{-2}$ are defined exactly as in the second part, so $\beta_0^{(g^\ddag)}$ and $\beta_{-2}^{(g^\ddag)}$ are also defined as above.
\item $c_{K^{\ddag,\p}}$ (resp. $c^\vee_{K^{\ddag,\p}}$) is a homomorphism $\wdtd{c}:\underline{X}\otimes \zbkpp\to A^\vee$ (resp. $\wdtd{c}^\vee: \underline{Y}\otimes \zbkpp\to A$) that lifts $c$ (resp. $c^\vee$) and remains unchanged under the action of $K^{\ddag,\p}$. Note that we can write $\underline{Y}\otimes \zbkpp$ as $\varinjlim_{n,\p\nmid n}\frac{1}{n}\underline{Y}$; the same is true for $\underline{X}$. We can therefore write $\wdtd{c}$ (resp. $\wdtd{c}^\vee$) as an inverse limit of $c_n:\frac{1}{n}\ul{X}\to A^\vee$ (resp. $c_n^\vee:\frac{1}{n}\ul{Y}\to A$) over all $n$, $\p\nmid n$.
\item $\tau_{K^{\ddag,\p}}$ is a trivialization of biextensions $\wdtd{\tau}:\mbf{1}_{(\underline{Y}\otimes\zbkpp)\times \underline{X}}\to (\wdtd{c}^\vee\times c)^*\ca{P}_A^{\otimes {-1}}$ that lifts $\tau$ and remains unchanged under the action of $K^{\ddag,\p}$. We can write $\wdtd{\tau}$ as an inverse limit of $\tau_n:\frac{1}{n}\ul{Y}\times\ul{X}\to (c_n^\vee\times c)^*\ca{P}_A^{\otimes -1}$ over all $n$, $\p\nmid n$.
\item\label{3-**} Denote by $\phi_n: \frac{1}{n}\underline{Y}\to \frac{1}{n}\underline{X}$ the base change of $\phi$ to $\frac{1}{n}\bb{Z}$, and by $\hat{\phi}:\gr^{W_{T^\p\Q}}_{0}\to (\gr^{W_{T^\p\Q}}_{-2})^\vee=\ul{\Hom}(\gr^{W_{T^\p\Q}}_{-2},\bb{G}_m)$ the base change of $\phi$ to $\zhpp$. Moreover, $\hat{\phi}$ naturally induces a pairing $\hat{\phi}^*:\gr_{-2}^{W_{T^\p\Q}}\times\gr_0^{W_{T^\p\Q}}\to \bb{G}_m$.\par
Denote by $\hat{\mbf{d}}_{10}: \gr_{-1}^{W_{T^\p\Q}}\times \gr_0^{W_{T^\p\Q}}\to T^\p\bb{G}_m$ the pairing defined as follows: For any $a=(a_n)\in \gr_{-1}^{W_{T^\p\Q}}=T^\p A$ and any $\hat{y}=(y_n)\in \gr_0^{W_{T^\p\Q}}=\varprojlim\ul{Y}/n\ul{Y}$, define $\hat{\mbf{d}}_{10}:(a,\hat{y})\mapsto \varprojlim_{n,\p\nmid n}\mbf{d}_{10,n}(a_n,y_n)$, where $\mbf{d}_{10,n}(a_n,y_n):=\e_{A[n]}(a_n,(\lambda_{A}\circ c_n^\vee-c_n\circ\phi_n)(\frac{1}{n}\wdtd{y}_n))$, which is independent of the choice of $\wdtd{y}_n\in\ul{Y}$ lifting $y_n$.\par 
Denote by $\hat{\mbf{d}}_{00}:\gr_0^{W_{T^\p\Q}}\times\gr_0^{W_{T^\p\Q}}\to T^\p\bb{G}_m$ the pairing defined as follows: for any $\hat{y}_1=(y_{1,n})_n,\hat{y}_2=(y_{2,n})_n\in\gr_0^{W_{T^\p\Q}}=\varprojlim \ul{Y}/n\ul{Y}$, define $\hat{\mbf{d}}_{00}:(\hat{y}_1,\hat{y}_2)\mapsto \varprojlim_{n,\p\nmid n}\mbf{d}_{00,n}(y_{1,n},y_{2,n})$, where $\mbf{d}_{00,n}(y_{1,n},y_{2,n}):=\tau_n(\frac{1}{n}\wdtd{y}_{1,n},\phi(\wdtd{y}_{2,n}))\tau_n(\frac{1}{n}\wdtd{y}_{2,n},\phi(\wdtd{y}_{1,n}))^{-1}$, which is independent of liftings $\wdtd{y}_{1,n}$ and $\wdtd{y}_{2,n}$ chosen in $\ul{Y}$.\par
Let $\hat{\mbf{d}}_{22}:\gr_{-2}^{W_{T^\p\Q}}\times\gr_{-2}^{W_{T^\p\Q}}\to T^\p\bb{G}_m$ and $\hat{\mbf{d}}_{21}:\gr_{-2}^{W_{T^\p\Q}}\times\gr_{-1}^{W_{T^\p\Q}}\to T^\p\bb{G}_m$ be trivial pairings. Let $\hat{\mbf{d}}_{20}=\hat{\phi}^*$, the pairing of $\gr^{W_{T^\p\Q}}_{0}\times\gr^{W_{T^\p\Q}}_{-2}$ induced by $\hat{\phi}$. Let $\hat{\mbf{d}}_{11}=\e^{\lambda_A}$. Let $\hat{\mbf{d}}_{ij}=-\hat{\mbf{d}}_{ji}^{T}$, the negative of the transpose of $\hat{\mbf{d}}_{ji},$ for $i,j=0,1,2$.\par
Denote by $\psi_{ij}$ the pullback of $\psi_{\zhpp}$ under $(\delta_{-i}^{(g^{\ddag})}\times \delta^{(g^{\ddag})}_{-j})$. 
Then $\beta^{(g^{\ddag}),*}(\hat{\mbf{d}}_{ij})=v\circ \psi_{ij}$, for $i,j=0,1,2$. 
\end{enumerate}
\end{enumerate}
\end{prop}
\begin{rk}Before presenting the proof, we make the following remarks:
\begin{enumerate}
\item One might not deduce this equivalence \emph{verbatim} from the proof of \cite[Thm. 5.2.3.14 and Cor. 5.2.3.15]{Lan13} since the pair $(\Q,\bml)$ might not come from $\mbf{F}_{\pol}(R,I)$, in particular, we are not assuming positive conditions for $\tau$ imposed in $\mbf{DD}_\pol(R,I)$. Nevertheless, we will closely follow the idea and the arguments of \emph{loc. cit.} in the proof below.  
\item Note that Part \ref{bd-lv-3} is helpful in the construction and the \emph{detailed} description of boundary charts (see \cite[Ch. 6]{Lan13}). For our purpose, it is more convenient for us to directly use Part \ref{bd-lv-1} in order to study the twisting action on the family of $1$-motives with additional structures in the next section.
\item One can change the sign conventions in the definition of Weil pairings (see Definition \ref{def-wp}) to get a possibly different sign in the pairings of the proposition above. We choose a convention that makes our results compatible with \cite{Lan13}.
\end{enumerate}
\end{rk}

We explain the proof of Proposition \ref{prop-deg-lvl} in steps. \par
\begin{proofof}[Proposition \ref{prop-deg-lvl}, part \ref{bd-lv-1}$\Longleftrightarrow$part \ref{bd-lv-2}]\par
Part \ref{bd-lv-1}$\implies$part \ref{bd-lv-2} can be directly deduced from \cite[Lem. 5.2.2.14]{Lan13}. Since $G^\ddag(\A)$ has a \emph{right} action on $(Z^{(g^{\ddag})},\Phi^{(g^{\ddag})})$ whose stabilizer is $g^{\ddag,-1}P_{\Phi^\ddag}(\A)g^{\ddag}$, part \ref{bd-lv-2}$\implies$part \ref{bd-lv-1} also follows from \emph{loc. cit.}; while applying \emph{loc. cit.}, one can replace $\overline{\eta}$ there with a general base $S$ and replace $(G,\lambda)$ there with $(\Q,\bml)$. 
\end{proofof}
Next, we show part \ref{bd-lv-2}$\Longleftrightarrow$part \ref{bd-lv-3}.\par
\begin{proofof}[Proposition \ref{prop-deg-lvl}, except that \ref{2-*}$\Longleftrightarrow$\ref{3-**}]\par
First, a splitting $\hat{\varsigma}_{12}:\gr^{W_{T^\p \Q}}_{-1}\oplus\gr^{W_{T^\p \Q}}_{-2}\xrightarrow{\sim}T^\p \G^{\natural}$ determines and is determined by an inverse system of splittings $\varsigma_{12,n}:\gr^{W_{T^\p \Q}}_{-1}[n]\oplus\gr^{W_{T^\p \Q}}_{-2}[n]\xrightarrow{\sim}\G^{\natural}[n]$, for $p\nmid n$ if $\p=\{p\}$. \par 
Write $\G^{\natural}$ as the extension $1\rightarrow T\xrightarrow{i}\G^{\natural}\xrightarrow{\pi}A\rightarrow 1$. For any positive integer $n$, the splitting $\varsigma_{12,n}:\gr^{W_{T^\p \Q}}_{-1}[n]\oplus\gr^{W_{T^\p \Q}}_{-2}[n]\xrightarrow{\sim}\G^{\natural}[n]$ determines and is determined by a projection $p_T:\G^\natural[n]\to T[n]$ such that $p_T\circ i=\id$; this is equivalent to an embedding $\overline{c}_n:\underline{X}/n\underline{X}\to\underline{\Hom}_{S-gp}(\G^\natural[n],\bb{G}_m)\iso \Q^{\vee,\circ}[n]$ such that $\overline{p}_X\circ\overline{c}_n=\id$, after taking $\underline{\Hom}_{S-gp} (-,\bb{G}_m)$. In the last sentence, $\Q^{\vee,\circ}:=\Q^\vee/T^\vee=[\ul{X}\xrightarrow{\ c\ }A^\vee]$, and $\overline{p}_X$ is a quotient of the natural projection to the first factor $\Q^{\vee,\circ}[n]\xrightarrow{p_1}\ul{X}\xrightarrow{\ }\ul{X}/n\ul{X}$.\par 
Then we can see that such a splitting $\overline{c}_n$ determines and is determined by a lifting $c_n:\frac{1}{n}\underline{X}\to A^\vee$ of $c$ such that the following diagram commutes:
\begin{equation}
\begin{tikzcd}
\underline{X}\arrow[r,"\sim"]\arrow[rrr,bend left,"c"]&\frac{1}{n}
\underline{X}\arrow[drr, "c_n"]&& A^\vee \\ 
\underline{X}\arrow[u,"{[n]}"]\arrow[ur,hook]\arrow[rrr,"c"]&&& A^\vee\arrow[u,"{[n]}"],
\end{tikzcd}
\end{equation}
where the natural isomorphism $\underline{X}\iso\frac{1}{n}\underline{X}$ is induced by ``dividing by $n$'' of coefficients $\bb{Z}\iso \bb{Z}\cdot\frac{1}{n}$. \par
In fact, this follows from the definition of $n$-torsion subgroup of $1$-motives. We can assume that $\ul{X}=X$. Fix any $S$-scheme $S^\prime$. If there is a splitting $\overline{c}_n:X/nX\to \Q^{\vee,\circ}[n]$, for any $\overline{\chi}\in X/nX(S^\prime)$, we can therefore find a pair $(\chi,a^\vee)$, lifting $\overline{c}_n(\overline{\chi})$, such that $\chi\in X(S^\prime)$ lifting $\overline{\chi}$, $a^\vee\in A^\vee(S^\prime)$, and $na^\vee=c(\chi)$. Define $c_n(\frac{1}{n}\chi)=a^\vee$. 
For any other choice of the pair $(\chi^\prime,a^{\vee,\prime})$ lifting $\overline{c}_n(\overline{\chi})$, $c(\frac{1}{n}(\chi^\prime-\chi))=a^{\vee,\prime}-a^\vee$. Then there is a well-defined homomorphism $c_n:\frac{1}{n}X\to A^\vee$ lifting $c$. 
Conversely, assume that there is a lifting $c_n:\frac{1}{n}X\to A^\vee$. Choose any set of representatives $\{\frac{1}{n}\chi\}$ of $\frac{1}{n}X/X$ in $\frac{1}{n}X$. 
Then define $\overline{c}_n:X/nX\to \Q^{\vee,\circ}[n],$ by $\overline{\chi}\mapsto[(\chi,c_n(\frac{1}{n}\chi))]$, which maps any $\overline{\chi}\in X/nX$ represented by $\chi$ to the class represented by $(\chi,c_n(\frac{1}{n}\chi))$ in $\Q^{\vee,\circ}[n]$.
\par
A splitting $\hat{\varsigma}: \bigoplus\limits_{i=0}^2\mrm{Gr}^{W_{T^\p \Q}}_{-i}\xrightarrow{\sim}T^\p \Q$ determines and is determined by $\hat{\varsigma}_{12}$ and the restriction $\hat{\varsigma}_0:=\hat{\varsigma}|_{\mrm{Gr}^{W_{T^\p \Q}}_0}$. And $\hat{\varsigma}_0$ is an inverse limit, over all $n$ such that $\p\nmid n$, of splittings $\varsigma_{0,n}:\ul{Y}/n\ul{Y}\to \Q[n]$ of the exact sequence $0\to \G^\natural[n]\to \Q[n]\to \ul{Y}/n\ul{Y}\to 0$. Then, by the same argument as in the third and the fourth paragraph of this proof, a splitting $\varsigma_{0,n}$ determines and is determined by liftings $\iota_n:\frac{1}{n}\ul{Y}\to \G^\natural$ of $\iota$ such that the following diagram commutes:
\begin{equation}
\begin{tikzcd}
\underline{Y}\arrow[r,"\sim"]\arrow[rrr,bend left,"\iota"]&\frac{1}{n}
\underline{Y}\arrow[drr, "\iota_n"]&& \G^\natural \\ 
\underline{Y}\arrow[u,"{[n]}"]\arrow[ur,hook]\arrow[rrr,"\iota"]&&& \G^\natural\arrow[u,"{[n]}"].
\end{tikzcd}
\end{equation}
Then $c_n^\vee:=\pi\circ \iota_n$ is a lifting of $c^\vee$, which corresponds to the splitting $\pi\circ\varsigma_{0,n}:\ul{Y}/n\ul{Y}\to \Q[n]/T[n]=\Q^\circ[n]$ of the exact sequence $0\to A[n]\to \Q^\circ[n]\to \ul{Y}/n\ul{Y}\to 0$.\par
Let $\wdtd{c}:=\varprojlim_{n,\p\nmid n}c_n$, $\wdtd{c}^\vee:=\varprojlim_{n,\p\nmid n} c_n^\vee$ and $\wdtd{\iota}:=\varprojlim_{n,\p\nmid n}\iota_n$; $\wdtd{\iota}$ corresponds to a trivialization $\wdtd{\tau}:\mbf{1}_{(\ul{Y}\otimes \zbkpp)\times \ul{X}}\to (\wdtd{c}^\vee\times c)^*\ca{P}_A^{\otimes -1}$ of biextensions.
Then we have shown that a splitting $\hat{\varsigma}$ is equivalent to a tuple of liftings $(\wdtd{c},\wdtd{c}^\vee,\wdtd{\iota})$ which lifts the tuple $(c,c^\vee,\iota)$; this tuple of liftings is further equivalent to a tuple $(\wdtd{c},\wdtd{c}^\vee,\wdtd{\tau})$.
\end{proofof}
Note that we will complete the proof of Proposition \ref{prop-deg-lvl} if we show \ref{2-*}$\Longleftrightarrow$\ref{3-**}.

\begin{proofof}[Proposition \ref{prop-deg-lvl}, \ref{2-*}$\Longleftrightarrow$\ref{3-**}: first reductions]\par
From now on, to the end of the proof of Proposition \ref{prop-deg-lvl}, we fix an $S$-scheme $S^\prime$ such that $\ul{X}=X$ and $\ul{Y}=Y$ over $S^\prime$.\par
We can write the splitting $\hat{\varsigma}: \bigoplus\limits_{i=0}^{2}\mrm{Gr}^{W_{T^\p \Q}}_{-i}\xrightarrow{\sim}T^\p \Q$ as an inverse limit, over all $n$ such that $\p\nmid n$, of $\varsigma_n: \bigoplus\limits_{i=0}^2\mrm{Gr}^{W_{\Q[n]}}_{-i}\xrightarrow{\sim}\Q[n]$. Let $\e_{ij,n}:=\e^{\bml}|_{\varsigma_n(\gr_{-i}^{W_{\Q[n]}}),\varsigma_n(\gr_{-j}^{W_{\Q[n]}})}$.\par 
To show \ref{2-*}$\Longleftrightarrow$\ref{3-**}, it suffices to show that $\hat{\e}_{ij}:=\varprojlim_{n,\p\nmid n}\e_{ij,n}=\hat{\mbf{d}}_{ij}$.\par
Firstly, we can show the claim above when $(i,j)=(2,1),(2,2)$ and $(1,1)$. In fact, since $\mrm{rig}_g:\ca{P}_{g,0}\to \ca{O}_{S^\prime}$ (resp. $\mrm{rig}_h:\ca{P}_{0,h}\to \ca{O}_{S^\prime}$) factors through $\mrm{rig}_{\pi(g)}:\ca{P}_{\pi(g),0}\to \ca{O}_{S^\prime}$ (resp. $\mrm{rig}_{\pi^\vee(h)}:\ca{P}_{0,\pi^\vee(h)}\to\ca{O}_{S^\prime}$) for any $g\in \G^\natural[n](S^\prime)$ (resp. $h\in\G^{\natural,\vee}[n](S^\prime)$), from Definition \ref{def-wp}, we find that 
\begin{equation}\label{eq-wp-g}\e_{\G^\natural[n]}(q_1,q_2)=\e_{A[n]}(\pi(q_1),\pi^\vee(q_2))\end{equation}
for any $q_1\in\G^\natural[n](S^\prime)$ and $q_2\in\G^{\natural,\vee}[n](S^\prime)$.\par
Then $\e_{21,n}=\e_{22,n}=0$ and $\e_{11,n}=\e^{\lambda_A}$.
\end{proofof}
Before completing the final steps, let us make some preparation on the Weil pairings and trivializations induced by the tuple $(c_n,c^\vee_n,\iota_n)$ and introduce some notation.\par

The Cartier dual of the splitting $\varsigma_n$ induces a splitting $$\varsigma^\vee_n:=\ul{\Hom}(\varsigma_n,\bb{G}_m)^{-1}:\bigoplus\limits_{i=0}^2\gr_{-i}^{W_{\Q^\vee[n]}}\xrightarrow{\sim}\Q^\vee[n].$$ 
As the proof above, the restriction $\varsigma^\vee_{0,n}:=\varsigma^\vee_n|_{\gr_0^{W_{\Q^\vee[n]}}}:\ul{X}/n\ul{X}\to \Q^\vee[n]$ is equivalent to a lifting $\iota_n^\vee:\frac{1}{n}\ul{X}\to \G^{\natural,\vee}$ of $\iota^\vee$ satisfying a similar commutative diagram:
\begin{equation}
\begin{tikzcd}
\underline{X}\arrow[r,"\sim"]\arrow[rrr,bend left,"\iota^\vee"]&\frac{1}{n}
\underline{X}\arrow[drr, "\iota_n^\vee"]&& \G^{\natural,\vee} \\ 
\underline{X}\arrow[u,"{[n]}"]\arrow[ur,hook]\arrow[rrr,"\iota^\vee"]&&& \G^{\natural,\vee}\arrow[u,"{[n]}"].
\end{tikzcd}
\end{equation}
In particular, $\pi^\vee\circ\iota^\vee_n=c_n$.\par
Denote by $\ca{R}^\vee$ the $1$-motive represented by $[\frac{1}{n}\ul{X}\xrightarrow{\iota^\vee_n}\G^{\natural,\vee}]$. Let $\ca{R}$ be the Cartier dual of $\ca{R}^\vee$, which is represented by $[\ul{Y}\xrightarrow{\wdtd{\iota}}\G^{\natural,\prime}]$. More precisely, $\G^{\natural,\prime}=\ul{\Ext}^1([\frac{1}{n}\ul{X}\xrightarrow{c_n}A],\bb{G}_m)$. Then $\ca{R}$ and $\ca{R}^\vee$ fit into the following commutative diagram:
\begin{equation}\label{dia-new-isog}
    \begin{tikzcd}
\ul{Y}\arrow[rr,"\wdtd{\iota}"]\arrow[drr,"\iota"]\arrow[d,"\phi"]&&\G^{\natural,\prime}\arrow[d,"{[n]_T}"]\\
\ul{X}\arrow[drr,"\iota^\vee"]\arrow[d,hook]&&\G^\natural\arrow[d,"{\bml}^{sab}"]\\
\frac{1}{n}\ul{X}\arrow[rr,"\iota_n^\vee"]&&\G^{\natural,\vee}.
    \end{tikzcd}
\end{equation}
In the diagram above, the homomorphism $[n]_T: \G^{\natural,\prime}\to \G^\natural$ is the homomorphism defined by the Cartier dual of 
\begin{equation*}
    \begin{tikzcd}
    \ul{X}\arrow[rr,"c"]\arrow[d,hook]&&A^{\vee}\\
    \frac{1}{n}\ul{X}\arrow[urr,"c_n"],&&    
    \end{tikzcd}
\end{equation*}
which is $T\xrightarrow{x\mapsto x^n}T$ over torus parts of $\G^{\natural,\prime}$ and $\G^{\natural}$.
\begin{lem}\label{lem-lft}
The image of $\varsigma_n(A[n])\sbst\G^\natural[n]$ lifts to $n$-torsion points in $\G^{\natural,\prime}$.
\end{lem}
\begin{proof}
The homomorphism $\overline{c}_n$ determines a homomorphism $\overline{\pi}^\vee: \Q^{\vee,\circ}[n]\to A^\vee[n]$: Indeed, for any $[(\chi,a)]\in \Q^{\vee,\circ}[n]$ such that $\chi\in X$ and $a\in A^\vee(S^\prime)$, $a-\overline{c}_n(\overline{\chi})\in A^\vee(S^\prime)$. Then $\overline{\pi}^\vee:[(\chi,a)]\mapsto a-\overline{c}_n(\overline{\chi})$ is a well-defined homomorphism. Since $\overline{c}_n$ is a splitting of $0\to A^\vee[n]\to\Q^{\vee,\circ}[n]\to \ul{X}/n\ul{X}\to 0$, $\ker\overline{\pi}^\vee=\im \overline{c}_n$ and $\overline{\pi}^\vee|_{A^\vee[n]}=\id_{A^\vee[n]}$. 
On the other hand, from the second step of the proof of Proposition \ref{prop-deg-lvl}, the image of $\overline{c}_n$ in $\Q^{\vee,\circ}[n]$ is $\im (n,c_n)/\im(n,c)$. So the image of $\overline{\pi}^\vee:\Q^{\vee,\circ}[n]\to A^\vee[n]$ is isomorphic to the image of the natural homomorphism $\Q^{\vee,\circ}[n]\to\ca{R}^{\vee,\circ}[n]$, where $\ca{R}^{\vee,\circ}:=[\frac{1}{n}\ul{X}\xrightarrow{c_n}A^\vee]$. Then, by Cartier duality, we find that $\varsigma_n(A[n])\sbst [n]_T(\G^{\natural,\prime}[n])$.
\end{proof}
In all diagrams below, all ``$\ca{P}$'' and ``$\ca{O}$'' denote biextensions or bundles, but not invertible sheaves.\par
\begin{proofof}[Proposition \ref{prop-deg-lvl}, \ref{2-*}$\Longleftrightarrow$\ref{3-**}: $\mbf{d}_{10,n}=\e_{10,n}$]\par
Let $\overline{y}\in Y/nY$ and $a\in A[n](S^\prime)$. Then $\varsigma_n(\overline{y})=(y,g)\in \Q[n]$ for some $y\in Y$ lifting $\overline{y}$ and some $g\in\G^\natural(S^\prime)$ such that $\iota(y)=ng$. And $\varsigma_n(a)=(0,h)$ such that $h\in\G^\natural[n](S^\prime)$ and $\pi(h)=a$. To compute $\e^{\bml}(\varsigma_n(a),\varsigma_n(\overline{y}))$, let us analyze the following diagram:
\begin{equation}\label{dia-d10}
    \begin{tikzcd}
        \ca{P}^{\otimes n}_{h,{\bml}^{sab}(g)}\arrow[rr,"\mrm{can}."]\arrow[d,"\mrm{can}."]&&\ca{P}_{h,{\bml}^{sab}\circ \iota(y)}\arrow[d,"\iso"]\\
        \ca{P}_{0,{\bml}^{sab}(g)}\arrow[d,"\iso"]&&\ca{O}_{\iota^\vee\circ\phi(y)}|_h\arrow[dd,"{\varrho_1(h,\phi(y))}"]\\
        \ca{O}_{{\bml}^{sab}\circ\iota_n(\frac{1}{n}y)}|_0\arrow[d,"{\varrho_2(0,{\bml}^{sab}\circ\iota_n(\frac{1}{n}y))}"']&&\\
        \ca{O}_{S^\prime}&&\ca{O}_{S^\prime}\arrow[ll,
        "{\e^{\bml}(\varsigma_n(a),\varsigma_n(\overline{y}))}"
        ].
    \end{tikzcd}
\end{equation}
By Lemma \ref{lem-lft}, there is an $h^\prime\in\G^{\natural,\prime}[n](S^\prime)$ such that $h=[n]_T(h^\prime)$. 
Let $[n]_A:\G^{\natural}\to\G^{\natural,\prime}$ be the pushout of $\G^\natural$ by $[n]_A:A\to A$. Applying Lemma \ref{lem-wp-ele} to the diagram
\begin{equation*}
    \begin{tikzcd}
\ul{Y}\arrow[rr,"\wdtd{\iota}"]\arrow[drr,"\iota"]\arrow[d,"{[n]}"]&&\G^{\natural,\prime}\arrow[d,"{[n]_T}"]\\
\ul{Y}\arrow[drr]&&\G^\natural\arrow[d,"{[n]_A}"]\\
&&\G^{\natural,\prime},
    \end{tikzcd}
\end{equation*}
we find that $\varrho_1^{\Q}(h,\phi(y))=\varrho_1^{\ca{R}}(h^\prime,\phi_n(y))=\varrho_1^{\ca{R}}(nh^\prime,\phi_n(\frac{1}{n}y))=\varrho_1^{\ca{R}}(0,\phi_n(\frac{1}{n}y))$.\par
The equality $\varrho_1^{\Q}(h,\phi(y))=\varrho_1^{\ca{R}}(h^\prime,\phi_n(y))$ amounts to the commutativity of 
\begin{equation*}
    \begin{tikzcd}
\ca{P}_{h^\prime,\bml^{sab}\circ\iota(y)}\arrow[d,"\mrm{can}."]\arrow[rrr,"{\varrho_1^{\ca{R}}(h^\prime,\phi_n(y))}"]&&&\ca{O}_{S^\prime}\\
\ca{P}_{h,\bml^{sab}\circ\iota(y)};\arrow[urrr,"{\varrho_1^{\Q}(h,\phi(y))}" description]
    \end{tikzcd}
\end{equation*}
also, $\ca{P}^{\otimes n}_{h^\prime,\bml^{sab}(g)}\xrightarrow{\mrm{can}.}\ca{P}_{0,\bml^{sab}(g)}$ factors through $\ca{P}^{\otimes n}_{h^\prime,\bml^{sab}(g)}\xrightarrow{\mrm{can}.}\ca{P}^{\otimes n}_{h,\bml^{sab}(g)}$. 
As a result, (\ref{dia-d10}) can be written as 
\begin{equation}\tag{$\ref{dia-d10}^\prime$}\label{dia-d10-prime}
     \begin{tikzcd}
        \ca{P}^{\otimes n}_{h^\prime,{\bml}^{sab}(g)}\arrow[rr,"\mrm{can}."]\arrow[d,"\mrm{can}."]&&\ca{P}_{h^\prime,{\bml}^{sab}\circ \iota(y)}\arrow[d,"\iso"]\\
        \ca{P}_{0,{\bml}^{sab}(g)}\arrow[d,"\iso"]&&\ca{O}_{\iota^\vee_n\circ\phi_n(y)}|_{h^\prime}\arrow[dd,"{\varrho_1^{\ca{R}}(h^\prime,\phi_n(y))}"]\\
        \ca{O}_{{\bml}^{sab}\circ\iota_n(\frac{1}{n}y)}|_0\arrow[d,"{\varrho_2(0,{\bml}^{sab}\circ\iota_n(\frac{1}{n}y))}"']&&\\
        \ca{O}_{S^\prime}&&\ca{O}_{S^\prime}\arrow[ll,
        "{\e^{\bml}(\varsigma_n(a),\varsigma_n(\overline{y}))}"
        ].
    \end{tikzcd}
\end{equation}
On the other hand, $\e_{\G^{\natural,\prime}[n]}(h^\prime,(\bml^{sab}\circ\iota_n-\iota_n^\vee\circ\phi_n)(\frac{1}{n}y))$ is computed by the following commutative diagram:
\begin{equation}\label{dia-e10}
    \begin{tikzcd}
     &\ca{P}^{\otimes n}_{h^\prime,\bml^{sab}(g)}\otimes ([n]^*\ca{O}_{-\iota_n^\vee\circ\phi_n(\frac{1}{n}y)})|_{h^\prime}\arrow[r,"\mrm{can}."]\arrow[dd,"\mrm{can}."]&(\ca{O}_{\iota^\vee_n\circ\phi_n(y)})|_{h^\prime}\arrow[d,phantom,"\otimes"]\\
   \ca{P}^{\otimes n}_{h^\prime,(\bml^{sab}\circ\iota_n)(\frac{1}{n}y)}\otimes\ca{P}^{\otimes n}_{h^\prime,(-\iota^\vee_n\circ\phi_n)(\frac{1}{n}y)}\arrow[ur,"\mrm{can}."]&&  (\ca{O}_{-\iota^\vee_n\circ\phi_n(y)})|_{h^\prime}\arrow[d,"\mrm{can. multi}."description]\\
   \ca{P}^{\otimes n}_{h^\prime,{(\bml^{sab}\circ\iota_n-\iota^\vee_n\circ\phi_n)(\frac{1}{n}y)}}\arrow[u,"\mrm{can}."']\arrow[r,"\mrm{can}."]\arrow[rr,bend left=10,"\mrm{can}." description]&\ca{O}_{\bml^{sab}\circ\iota_n(\frac{1}{n}y)}|_0\otimes\ca{O}_{-\iota_n^\vee\circ\phi_n(\frac{1}{n}y)}|_0\arrow[dd,"{\mrm{rig}_{\bml^{sab}\circ\iota_n(\frac{1}{n}y)-\iota^\vee_n\circ \phi_n(\frac{1}{n}y)}}"]& \ca{O}_{\G^{\natural,\prime}[n]}|_{h^\prime}\arrow[dd,"\mrm{rig}_{h^\prime}"]\\\\
    &\ca{O}_{S^\prime}&\ca{O}_{S^\prime}\arrow[l,"{\e(h^\prime,(\bml^{sab}\circ\iota_n-\iota_n^\vee\circ\phi_n)(\frac{1}{n}y))}"'].
    \end{tikzcd}
\end{equation}
In the diagram above, the canonical homomorphism $\ca{O}_{\iota_n^\vee\circ\phi_n(y)}|_{h^\prime}\otimes \ca{O}_{-\iota_n^\vee\circ\phi_n(y)}|_{h^\prime}\xrightarrow{\mrm{can.multi}.}\ca{O}_{\G^{\natural,\prime}[n]}|_{h^\prime}$ is deduced from the multiplication on $A^\vee$; since $\varrho_1^{\ca{R}}$ is a trivialization of biextensions, its multiplicative structure implies that $\mrm{rig}_{h^\prime}\circ (\mrm{can.multi.})=\varrho_1(h^\prime,\phi_n(y))\cdot\varrho_1(h^\prime,-\phi_n(y))$. Then the first (resp. second) column of (\ref{dia-d10-prime}), tensoring with $\varrho_1^{\ca{R}}(0,-\phi_n(\frac{1}{n}y))$ (resp. with $\varrho_1^{\ca{R}}(h^\prime,-\phi_n(y))$), is the second (resp. third) column of the diagram (\ref{dia-e10}).\par
Then the diagram (\ref{dia-e10}) implies that $\e^{\bml}(\varsigma_n(a),\varsigma_n(\overline{y}))=\e_{\G^{\natural,\prime}[n]}(h^\prime,(\bml^{sab}\circ\iota_n-\iota^\vee_n\circ\phi_n)(\frac{1}{n}y))=\e_{A[n]}(a,(\lambda_A\circ c^\vee_n-c_n\circ \phi_n)(\frac{1}{n}y))$, where the last equality is from (\ref{eq-wp-g}).
\end{proofof}
\begin{proofof}[Proposition \ref{prop-deg-lvl}, \ref{2-*}$\Longleftrightarrow$\ref{3-**}: $\mbf{d}_{20,n}=\e_{20,n}$]\par
Let $t\in T[n](S^\prime)$ and $\overline{y}\in Y/nY$. As before, set $\varsigma_n(\overline{y})=(y,g)\in\Q[n](S^\prime)$, for some $y\in Y$ lifting $\overline{y}$ and some $g\in\G^\natural(S^\prime)$ such that $\iota(y)=ng$.\par
Then $\e_{20}(t,(y,g))$ is computed by the following diagram:
\begin{equation}
    \begin{tikzcd}
&\ca{P}^{\otimes n}_{t,\bml^{sab}(g)}\arrow[dl,"\mrm{can}."']\arrow[dr,"\mrm{can}."]&\\
\ca{P}_{0,\bml^{sab}(g)}\arrow[d,"{\varrho_2(0,\bml^{sab}(g))}"']&&\ca{P}_{t,\iota^\vee\circ\phi(y)}\arrow[d,"{\varrho_1(t,\phi(y))}"]\\
\ca{O}_{S^\prime}&&\ca{O}_{S^\prime}\arrow[ll,"{\e^{\bml}(\varsigma_n(t),\varsigma_n(\overline{y}))}"'].
    \end{tikzcd}
\end{equation}
Let $\wdtd{t}\in \G^{\natural,\prime}(S^\prime)$ lifting $t$. Then $\wdtd{t}\in\G^{\natural,\prime}[n^2](S^\prime)$. Note that there is a canonical morphism $\ca{O}^{\otimes n}_{\bml^{sab}(g)}|_{T}\xrightarrow{\mrm{can}.}([n]_T^*\ca{O}_{\bml^{sab}(g)})|_T\iso[n]^*(\ca{O}_{\bml^{sab}(g)}|_T)|_T\xrightarrow{\mrm{can}.}\ca{O}_{\bml^{sab}(g)}|_T$, we denote this canonical morphism also by $[n]$; restricting this morphism to $\wdtd{t}$, there is a canonical morphism $[n]:\ca{P}_{\wdtd{t},\bml^{sab}(g)}^{\otimes n}\iso\ca{O}^{\otimes n}_{\bml^{sab}(g)}|_{\wdtd{t}}\xrightarrow{\mrm{can}.} \ca{O}_{\bml^{sab}(g)}|_t\iso\ca{P}_{t,\bml^{sab}(g)}$.\par
Moreover, since $\bml^{sab}(g)=\iota^\vee_n\circ\phi_n(\frac{1}{n}y)$, we know $\e_{\ca{R}[n^2]}(\wdtd{t},\bml^{sab}(g))=\e_{\ca{R}[n^2]}(\wdtd{t},\phi_n(\frac{1}{n}y))=0$; that is, $\varrho_1^{\ca{R}}(n\wdtd{t},\phi_n(y))=\varrho_1^{\ca{R}}(\wdtd{t},\phi_n(ny))=\varrho_2^{\ca{R}}(0,\bml^{sab}(g))$. Also, since $[n]$ commutes with tensor products, 
\begin{equation*}
    \begin{tikzcd}
    \ca{P}^{\otimes n^2}_{\wdtd{t},\iota_n^\vee\circ\phi_n(\frac{1}{n}y)}\arrow[rr,"\mrm{can}."]\arrow[d,"{[n]}"]&&\ca{P}^{\otimes n}_{\wdtd{t},\iota_n^\vee\circ\phi_n(y)}\arrow[d,"{[n]}"]\\
    \ca{P}^{\otimes n}_{t,\bml^{sab}(g)}\arrow[rr,"\mrm{can}."]&&\ca{P}_{t,\bml^{sab}(y)}
    \end{tikzcd}
\end{equation*}
commutes.
As a result, we have the following commutative diagram, which computes $\e_{20}$:
\begin{equation}
    \begin{tikzcd}
&\ca{P}^{\otimes n^2}_{\wdtd{t},\iota_n^\vee\circ\phi_n(\frac{1}{n}y)}\arrow[ddl,"\mrm{can}."']\arrow[d,"{[n]}"]\arrow[ddr,"\mrm{can}." description]\arrow[drr,"\mrm{can}."]&&\\
&\ca{P}^{\otimes n}_{t,\bml^{sab}(g)}\arrow[dl,"\mrm{can}."']\arrow[drr,"\mrm{can}."]&&\ca{P}^{\otimes n}_{\wdtd{t},\iota^\vee_n\circ\phi_n(y)}\arrow[d,"{[n]}"]\\
\ca{P}_{0,\bml^{sab}(g)}\arrow[d,"{\varrho_2}"]&&\ca{P}_{\wdtd{t},\iota_n^\vee\circ\phi_n(ny)}\arrow[dll,"{\varrho_1}"]&\ca{P}_{t,\iota^\vee\circ\phi(y)}\arrow[d,"{\varrho_1}"]\\
\ca{O}_{S^\prime}&&&\ca{O}_{S^\prime}\arrow[lll,"{\e^{\bml}(\varsigma_n(t),\varsigma_n(\overline{y}))}"'].
    \end{tikzcd}
\end{equation}
Then, from Lemma \ref{lem-wp-ele}, $\e_{20}$ is computed by the commutative diagram
\begin{equation}
    \begin{tikzcd}
\ca{P}^{\otimes n^2}_{\wdtd{t},\iota_n^\vee\circ\phi_n(\frac{1}{n}y)}\arrow[rr,"\mrm{can}."]\arrow[dd,"\mrm{can}."]&&\ca{P}^{\otimes n}_{\wdtd{t},\iota^\vee_n\circ\phi_n(y)}\arrow[dd,"{[n]}"]\arrow[dl,"\mrm{can}."]\\
&\ca{P}_{n\wdtd{t},\iota^\vee_n\circ\phi_n(y)}\arrow[d,"\mrm{can}."]&\\
\ca{P}_{n\wdtd{t},\iota^\vee_n\circ\phi_n(y)}\arrow[r,"\mrm{can}."]\arrow[d,"{\varrho_1}"]&\ca{P}_{0,\iota^\vee\circ\phi(y)}\arrow[dl,"{\varrho_1}"]&\ca{P}_{t,\iota^\vee\circ\phi(y)}\arrow[d,"{\varrho_1}"]\\
\ca{O}_{S^\prime}&&\ca{O}_{S^\prime}\arrow[ll,"{\e^{\bml}(\varsigma_n(t),\varsigma_n(\overline{y}))}"'].
    \end{tikzcd}
\end{equation}
Finally, viewing $T$ as a subgroup of $\G^{\natural,\prime}$, we conclude that $\e^{\bml}(t,(y,g))$ is the difference between the evaluation of $\phi_n(y):T\to \bb{G}_m$ at $\wdtd{t}^n$ and the evaluation of $\phi_n(y)$ at $\wdtd{t}^{n^2}$. By the sign conventions in Lemma \ref{lem-tri}, $\e^{\bml}(t,(y,g))=\phi^*(y)(t)$.
\end{proofof}
Note that the Cartier dual of 
\begin{equation*}
    \begin{tikzcd}
    \ul{Y}\arrow[rr,"\iota"]\arrow[d,hook]&&\G^\natural\\
    \frac{1}{n}\ul{Y}\arrow[urr,"\iota_n"]&&
    \end{tikzcd}
\end{equation*}
is a commutative diagram
\begin{equation*}
    \begin{tikzcd}
\ul{X}\arrow[rr,"\wdtd{\iota}^\vee"]\arrow[drr,"\iota^\vee"]&&\G^{\natural,\wedge}\arrow[d,"n^\vee"]\\
&& \G^{\natural,\vee}.
    \end{tikzcd}
\end{equation*}
Denote by $\Q^\wedge:=[\ul{X}\xrightarrow{\wdtd{\iota}^\vee}\G^{\natural,\wedge}]$ the Cartier dual of $\Q_n:=[\frac{1}{n}\ul{Y}\xrightarrow{\iota_n}\G^\natural]$. The homomorphism $n^\vee:\G^{\natural,\wedge}\to \G^{\natural,\vee}$ is $T^\vee\xrightarrow{x\mapsto x^n}T^\vee$ when restricted to torus parts. \par
More precisely, we have the following commutative diagram which is similar to (\ref{dia-new-isog}):
\begin{equation}\label{dia-new-wedge}
    \begin{tikzcd}
\ul{X}\arrow[rr,"\wdtd{\iota}^\vee"]\arrow[drr,"\iota^\vee"]\arrow[d,"\phi^c"]&&\G^{\natural,\wedge}\arrow[d,"n^\vee"]\\
\ul{Y}\arrow[drr,"\iota"]\arrow[d,hook]&&\G^{\natural,\vee}\arrow[d,"{\bml}^{c,sab}"]\\
\frac{1}{n}\ul{Y}\arrow[rr,"\iota_n"]&&\G^{\natural}.        
    \end{tikzcd}
\end{equation}
To form the diagram above, we can find a polarization $\bml^c:\Q^\vee\to \Q$ such that $\bml\circ \bml^c=[N]_{\Q^\vee}$ and $\bml^c\circ \bml=[N]_{\Q}$ for some integer $N>0$; then we denote $\bml^c:=(\bml^{c,et},\bml^{c,sab})$ and $\bml^{c,et}:=\phi^c$.\par
\begin{proofof}[Proposition \ref{prop-deg-lvl}, \ref{2-*}$\Longleftrightarrow$\ref{3-**}: $\mbf{d}_{00,n}=\e_{00,n}$]\par
Suppose that $\varsigma_n(\overline{y}_i)=(y_i,g_i)\in \Q[n]$ for some $y_i\in Y$ lifting $\overline{y}_i$ and some $g_i\in\G^\natural(S^\prime)$ such that $\iota(y_i)=ng_i$, for $i= 1,2$. 
Then $\e^{\bml}(\varsigma_n(\overline{y}_1),\varsigma_n(\overline{y}_2))$ is computed by the following commutative diagram:
\begin{equation}
    \begin{tikzcd}
\ca{P}^{\otimes n}_{g_1,\bml^{sab}(g_2)}\arrow[r,"\mrm{can}.","\sim"']\arrow[d,"\mrm{can}.","\sim"']&
\ca{P}_{\iota(y_1),\bml^{sab}(g_2)}\arrow[r,"\mrm{can}.","\sim"']\arrow[dr,"{\varrho_2(y_1,\bml^{sab}(g_2))}"']&\ca{P}_{\iota^\vee\circ\phi(y_1),g_2}\arrow[r,"\mrm{can}.","\sim"']&
\ca{P}_{\wdtd{\iota}^\vee\circ\phi(y_1),\iota_n(\frac{1}{n}y_2)}\arrow[dl,"{\tau_{\Q^\wedge}(\phi(y_1),\frac{1}{n}y_2)}" description]\arrow[d,"\mrm{can}.","\sim"']\\
\ca{P}_{g_1,\iota^\vee\circ\phi(y_2)}\arrow[d,"\mrm{can}.","\sim"']\arrow[drr,"{\varrho_1(g_1,\phi(y_2))}"]&&\ca{O}_{S^\prime}&\ca{P}_{\iota_n(\frac{1}{n}y_2),\wdtd{\iota}^\vee\circ\phi(y_1)}\arrow[l,"{\tau_{n}(\frac{1}{n}y_2,\phi(y_1))}"]\\
\ca{P}_{\iota_n(\frac{1}{n}y_1),\wdtd{\iota}^\vee\circ\phi(y_2)}\arrow[rr,"{\tau_{n}(\frac{1}{n}y_1,\phi(y_2))}"]&&\ca{O}_{S^\prime}\arrow[u,"{\e^{\bml}(\varsigma_n(\overline{y}_1),\varsigma_n(\overline{y}_2))}"'].&
    \end{tikzcd}
\end{equation}
More precisely, the commutativity of $\tau_{n}(\frac{1}{n}y_1,\phi(y_2))\circ\mrm{can}.=\varrho_1(g_1,\phi(y_2))$ and $\tau_{\Q^\wedge}(\phi(y_1),\frac{1}{n}y_2)\circ\mrm{can}.\circ\mrm{can}.=\varrho_2(y_1,\bml^{sab}(g_2))$ follows from Lemma \ref{lem-wp-ele}; the commutativity of $\tau_{n}(\frac{1}{n}y_2,\phi(y_1))\circ\mrm{can}.=\tau_{\Q^\wedge}(\phi(y_1),\frac{1}{n}y_2)$ follows from \ref{iota}$\Longleftrightarrow$\ref{eq-44}
Hence, we have shown the desired result.
\end{proofof}
Now we have completed the proof of Proposition \ref{prop-deg-lvl}.\hfill{$\square$}
\subsection{Boundary mixed Shimura varieties of Hodge type, continued}\label{subsec-mxhg-cont}
Fix a cusp label representative $\Phi_0=(Q,X^+_0,g_0)$ of $(G_0,X_0)$. With the conventions in \S\ref{subsec-deg-data-lvl}, suppose that $\Phi_0$ is sent to a cusp label representative $\Phi^\ddag$ of $(G^\ddag,X^\ddag)$ under $\iota$. Let us describe the moduli interpretation of $\sh_{K_{\Phi^\ddag}}(\bb{C})$ as a family of polarized $1$-motives with level structures.
\subsubsection{}
Recall that $D_{\Phi_0}:=D_{Q,X^+_0}$ is defined to be the $P_{\Phi_0}(\bb{R})U_{\Phi_0}(\bb{C})$-orbit of $([x],u_x^{Q}\circ h_\infty)\in \pi_0(X_0)\times \Hom_\bb{C}(\bb{S}_\bb{C},P_{\Phi_0,\bb{C}})$, for any $x\in X_0^+$. 
Moreover, by Proposition \ref{pure shimura}, $D_{\Phi_0}$ is in fact the $P_{\Phi_0}(\bb{R})U_{\Phi_0}(\bb{C})$-orbit of $u_x^Q\circ h_\infty\in\Hom_\bb{C}(\bb{S}_\bb{C},P_{\Phi_0,\bb{C}})$, and the morphism $D_{\Phi_0}\to D_{\Phi^\ddag}\sbst\Hom_\bb{C}(\bb{S}_\bb{C},P_{\Phi^\ddag,\bb{C}})$ induced by $\iota$ is injective. \par
For any $x^\ddag\in D_{\Phi^\ddag}$, $x^\ddag$ determines a rational mixed Hodge structure $\mxh(x^\ddag):=(V,F^\bullet_{x^\ddag},W_\bullet^{x^\ddag})$ on the $\bb{Q}$-vector space $V$, where $F^\bullet_{x^\ddag}$ is a descending filtration on $V_\bb{C}$ and where $W_\bullet^{x^\ddag}$ is an ascending filtration on $V$; by \cite[Prop. 1.16(c)]{Pin89}, the association $x^\ddag\mapsto \mxh(x^\ddag)$ is injective. 
For any $\bb{Z}$-lattice $\mathfrak{V}\sbst V$, $\mxh(x^\ddag)$ determines an integral mixed Hodge structure $\imxh(\mathfrak{V},x^\ddag):=(\mathfrak{V},F^\bullet_{x^\ddag}|_\mathfrak{V},W_\bullet^{x^\ddag}|_\mathfrak{V})$; 
we will abusively write $\imxh(\mathfrak{V})$ or $\imxh(x^\ddag)$ if $x^\ddag$ or $\mathfrak{V}$ is clear in the context; we apply similar shorthand for $\mxh$. Note that if $V$ and $\mathfrak{V}$ are fixed, then $\mxh(x^\ddag)$ determines and is determined by $\imxh(x^\ddag)$.\par
On the other hand, $x^\ddag$ induces on $V$ a mixed Hodge structure of type (see \cite[Ex. 4.25]{Pin89})
$$\{(0,0);(-1,0),(0,-1);(-1,-1)\}.$$ 
The category of integral mixed Hodge structures of this type is closely related to the category of $1$-motives over $\bb{C}$. In fact, Deligne showed that one can geometrize this type of integral mixed Hodge structures.
\begin{thm}[{\cite[10.1.3]{Del74}}]\label{thm-deligne-iii}
For any $1$-motive $\Q=(\ul{Y},\G^\natural,T,A,\iota,c^\vee)$ over $\bb{C}$, there is a canonical finitely generated free $\bb{Z}$-module $T_{\bb{Z}}\Q$ with Hodge filtration $F^\bullet$ and weight filtration $W_\bullet$ on it, such that $(T_\bb{Z}\Q,F^\bullet,W_\bullet)$ is an integral mixed Hodge structure of type $\{(0,0);(-1,0),(0,-1);(-1,-1)\}$, and such that, $\gr^W_{-2}T_\bb{Z}\Q=H_1(T,\bb{Z})=X^\vee\iso\mbf{X}_*(T)$, $\gr^W_{-1}T_\bb{Z}\Q=H_1(A,\bb{Z})$ and $\gr^W_0T_\bb{Z}\Q=Y$.\par
Moreover, the functor $T_\bb{Z}(-)$, from the category of $1$-motives over $\bb{C}$ to the category of integral mixed Hodge structures of type $\{(0,0);(-1,0),(0,-1);(-1,-1)\}$ whose weight $-1$ graded pieces are polarizable, is an equivalence of category.
\end{thm}
The module $T_\bb{Z}\Q$, equipped with its weight filtration and Hodge filtration, is called the \textbf{Hodge realization} of $\Q$. There are canonical comparison isomorphisms between the Hodge realization and the Tate modules, respecting weight filtrations and its graded pieces: $T_\bb{Z}\Q\otimes \wat{\bb{Z}}\iso \wat{T}\Q$, $T_\bb{Z}\Q\otimes \zhp\iso T^\p \Q$ and $T_\bb{Z}\Q\otimes \bb{Z}_p\iso T_p\Q$. For this, see \cite[10.1.5-10.1.10]{Del74}.\par
One can also associate with $\Q$ a canonical \textbf{de Rham realization} $H_\dr \Q$. There is a canonical exact sequence
\begin{equation}\label{eq-hodge-ex}
    0\to F^0H_\dr\Q\iso \ull_{\G^{\natural,\vee}/\bb{C}}^\vee\to H_\dr\Q\to \ull_{\G^\natural/\bb{C}}\to 0.
\end{equation}
There is also a canonical comparison $H_\dr\Q\iso T_\bb{Z}\Q\otimes \bb{C}$.
For this, see \emph{loc. cit.}, and see \cite{Ber05}.
\subsubsection{}
Fix any polarized $1$-motive $(\Q,\bml)$ over $\bb{C}$. Denote $\Q=(\ul{Y},\G^\natural,T,A,\iota,c^\vee)$. Let $\mathfrak{V}=T_\bb{Z}\Q$. The Hodge realization of Cartier dual $\Q^\vee$, $\mathfrak{V}^\vee:=T_\bb{Z}(\Q^\vee)$, can be naturally identified with $T_\bb{Z}(\Q^\vee)=\Hom(\mathfrak{V},\bb{Z}(1))$, induced by a canonical pairing $\Phi^{\can}: \mathfrak{V}\times \mathfrak{V}^\vee\to \bb{Z}(1)$.
Then $\bml:\Q\to \Q^\vee$ induces a pairing $\psi:\mathfrak{V}\times \mathfrak{V}\to \bb{Z}(1)$ by defining $\psi:=\Phi^\can\circ(\id,T_\bb{Z}\bml)$. Its base change to $V:=\mathfrak{V}\otimes\bb{Q}$ is a pairing $\psi: V\times V\to \bb{Q}(1)$. By \cite[Rmk. 10.2.4]{Del74}, $\psi$ is an alternating pairing.\par
The weight $-2$ graded piece $\mathfrak{V}_{-2}:=\gr^W_{-2}\mathfrak{V}\iso H_1(T,\bb{Z})$ is a totally isotropic space of $\mathfrak{V}$ with respect to $\psi$, whose orthogonal complement is $\mathfrak{V}^\perp=W_{-1}\mathfrak{V}$. The alternating pairing $\psi$ induces a pairing $\psi^{et}:\mathfrak{V}_0\times \mathfrak{V}_{-2}\to \bb{Z}(1)$, where $\mathfrak{V}_0:=\gr^W_0\mathfrak{V}\iso Y$; note that $\psi^{et}$ does not depend on the choice of liftings for $\mathfrak{V}_0$. 
Moreover, since $\psi(\mathfrak{V}_{-2},\mathfrak{V}^\perp)=0$, there is a canonical alternating pairing $\psi^{ab}:\mathfrak{V}_{-1}\times\mathfrak{V}_{-1}\to \bb{Z}(1)$, defined by projecting to $\mathfrak{V}_{-1}:=\gr^W_{-1}\mathfrak{V}\iso H_1(A,\bb{Z})$ of the restriction of $\psi$ to $\mathfrak{V}_{-1}$. By Riemann's Theorem, $\psi^{ab}$ is a polarization of $\mathfrak{V}_{-1}$, an integral Hodge structure of type $\{(-1,0),(0,-1)\}$.\par
Note that $\psi^{et}$ defines an embedding $Y\iso \mathfrak{V}_0\hookrightarrow X\iso \mathfrak{V}_{-2}^\vee$ with finite cokernel, which we abusively denote by $\psi^{et}$. Together with $\psi^{ab}$, and by Riemann's Theorem, the pair $(\psi^{et},\psi^{ab})$ determines and is determined by the polarization $\bml=(\bml^{et},\bml^{ab})$. 
We will call such a $\psi$ induced by $\bml$ a \textbf{polarization} of the integral mixed Hodge structure $\mathfrak{V}$.\par
\begin{lem}\label{lem-mx-hdg}
Let us fix $(V,\psi)$ and a cusp label representative $\Phi^\ddag=(Q^\ddag,X^{\ddag,+},g^\ddag)$ of $G^\ddag$. Fix a $\bb{Z}$-lattice $\mathfrak{V}\sbst V$. Let $\imxh(\mathfrak{V})$ be an integral mixed Hodge structure defined by $x^\ddag\in D_{\Phi^\ddag}^+$.
Suppose that there is another pair $(\imxh(\mathfrak{U}),\psi_0)$ of integral mixed Hodge structures of type $\{(0,0);(-1,0),(0,-1);(-1,-1)\}$ with polarization $\psi_0$.
Suppose that there is an isomorphism of lattices $f:\mathfrak{V}\to \mathfrak{U}$, which matches $\psi$ and $\psi_0$, and matches the weight filtrations of $\mathfrak{V}$ and $\mathfrak{U}$.\par
Then there is $g\in P_{\Phi^\ddag}(\bb{R})U_{\Phi^\ddag}(\bb{C})$, such that the mixed Hodge structure on $V_\bb{C}:=\mathfrak{V}\otimes \bb{C}$ induced by the action of $g$ on the mixed Hodge structure $\imxh(\mathfrak{V})\otimes \bb{C}$ coincides with the pullback of the mixed Hodge structure on $U_\bb{C}:=\mathfrak{U}\otimes \bb{C}$ to $V_\bb{C}$ under $f$.
\end{lem}
\begin{proof}
Denote by $\mathfrak{V}_{-i}$ and $\mathfrak{U}_{-i}$ for $i=0,1,2$ the graded pieces of $\mathfrak{V}$ and $\mathfrak{U}$, with respect to the corresponding weights $-i$. Denote by $W^\ddag$ the unipotent radical of $Q^\ddag$, $P^\ddag:=P_{\Phi^\ddag}$ and $U^\ddag:=U_{\Phi^\ddag}$. After pulling back the mixed Hodge structure on $U_\bb{C}$ to $V_\bb{C}$ via $f$, we can assume $f=\id$ without loss of generality.
By \cite[Prop. 1.16 (b)]{Pin89}, it suffices to show that there is a $g\in P^\ddag(\bb{R})W^\ddag(\bb{C})$, such that, the weight filtrations and Hodge filtrations on $V_\bb{C}$ induced by $\imxh(\mathfrak{V})$ and induced by (the pullback of) $\imxh(\mathfrak{U})$ are matched after an action of $g$. The weight filtrations are assumed to be equal; we denote it by $W$.\par
By \cite[Lem. 1.2.11]{Del71b}, there is a canonical bigrading $A^{p,q}_V$ for $\mxh(V)$, such that, $V_{\bb{C}}=\bigoplus\limits_{(p,q)=(-1,-1),(-1,0),(0,-1),(0,0)}A^{p,q}_V$, and such that $A^{p,q}_V$ can recover the weight filtration $W_\bullet$ and the Hodge filtration $F^\bullet_V$ of $V_\bb{C}$. See also \cite[Prop. 1.2]{Pin89}. In our case, we have: $A^{-1,-1}_V=\mathfrak{V}_{-2}\otimes \bb{C}$; the projection of $A_V^{-1,0}\oplus A_V^{0,-1}$ to $\gr^W_{-1}V_\bb{C}$ is the Hodge decomposition of $\gr^W_{-1}V_\bb{C}$; $A^{0,0}_V\iso \gr^W_0 V_\bb{C}$.
Then the Hodge structure is determined by a homomorphism $h_V:\bb{S}_\bb{C}\to \mrm{GL}(V_\bb{C})$, such that $(z_1,z_2)\mapsto z_1^{-p}z_2^{-q}$ on $A^{p,q}_V$. By the same lemma, there are similar complex spaces $A^{p,q}_U$ for $U_\bb{C}$ with corresponding property, replacing ``$V$'' with ``$U$''. Note that the bigrading pieces $A^{p,q}_V$ and $A^{p,q}_U$ are canonically determined by the mixed Hodge structure. To show the statement, it suffices to show that $h_V$ and $h_U$ are conjugate under some element $g\in P^\ddag(\bb{R})W^\ddag(\bb{C})$.\par
Let $\varsigma_V:\bigoplus\limits_{(p,q)}A_V^{p,q}\to V_\bb{C}$ be the canonical splitting, and let $\varsigma_U$ be the corresponding object for $U$. The pullback of $\psi$ to $\bigoplus_{(p,q)}A_U^{p,q}$ (resp. $\bigoplus_{(p,q)}A_V^{(p,q)}$) under $\varsigma_U$ (resp. $\varsigma_U$) is a $3\times 3$ anti-diagonal blocked matrix. In fact, we can see this via the perfect pairing of the de Rham realizations of $1$-motives over $\bb{C}$. See \cite[10.2.7]{Del74} and \cite{Ber05}. Since $A_V^{-1,-1}=A_U^{-1,-1}=\mathfrak{V}_{-2}\otimes \bb{C}$ and since $A_V^{0,0}=A_U^{0,0}=\gr_0^WV_\bb{C}$, $\varsigma_U^{-1}\circ\varsigma_V$ matches the $0$-th and $-2$-th graded pieces of the weight filtration $W$. 
After taking the $W$-graded pieces $\bigoplus\limits_{i=-2}^0\gr_{i}^W V_\bb{C}$, there is a $\overline{g}\in G_{\Phi^\ddag,h}(\bb{R})$ such that the Hodge filtrations on each graded piece induced by $\imxh(\mathfrak{V})$ and $\imxh(\mathfrak{U})$ are matched after the action of $\overline{g}$. Also, we know that $\overline{g}$ lifts to an element $p\in P^\ddag(\bb{R})$.
Then after replacing $\varsigma_V$ with $p\circ\varsigma_V$, we can identify $A_V^{p,q}$ with $A_U^{p,q}$ for each $(p,q)$. We can also choose an identification $\bigoplus_{(p,q)}A^{p,q}_U\iso \bigoplus_{(p,q)}A^{p,q}_V$ such that it matches $\varsigma_U^*\psi$ with $\varsigma_V^*\psi$. 
We then see that $\varsigma_U^{-1}\circ \varsigma_V$ corresponds to a matrix in $W^\ddag(\bb{C})$.
\end{proof}
\subsubsection{}
With the background in previous {subsections}, we can summarize the desired moduli interpretation now. 
Recall that $\p$ is a finite set of prime numbers, which, in our case, is $\{p\}$ or $\emptyset$. Let us fix a cusp label representative $\Phi^\ddag=(Q^\ddag,X^{\ddag,+},g^\ddag)$ of $G^\ddag$. 
Recall that in \S\ref{subsec-siegel} that $(V_{\bb{Z}},\psi_{\bb{Z}})$ is a pair of $\bb{Z}$-lattice $V_\bb{Z}$ with $\psi_\bb{Z}:V_\bb{Z}\times V_\bb{Z}\to \bb{Z}(1)$, the restriction of $\psi$ to $V_\bb{Z}$ which takes its values in $\bb{Z}(1)$.\par
Let $\bm{\xi}_{{(V_{\bb{Z}},\psi_\bb{Z})},\Phi^\ddag,K^{\ddag,\p}}:(Sch/\zbkpp)^{\mrm{op}}\to(Grpds)$ be a moduli problem of isomorphism classes of $1$-motives with additional structures, which can be described as follows:\par
For any locally Noetherian and connected $\zbkpp$-scheme $S$, $\bm{\xi}_{(V_{\bb{Z}},\psi_\bb{Z}),\Phi^\ddag,K^{\ddag,\p}}(S)$ is a groupoid whose objects consist of the tuples $(\Q,\bml,[u]_{(K^\p_{\Phi^\ddag})^{g^\ddag}})$, where 
\begin{itemize}
    \item $(\Q,\bml)$ is a polarized $1$-motive over $S$ such that $\ul{X}_\Q$ and $\ul{Y}_\Q$ split, and
    \item $[u]_{(K^\p_{\Phi^\ddag})^{g^\ddag}}\in \Gamma(S,\isom(V_{\wat{\bb{Z}}^\p},T^\p\Q)/\rcj{(K^\p_{\Phi^\ddag})}{g^\ddag})$ is a level structure of type $(V_\bb{Z},\psi_\bb{Z},\Phi^\ddag,K^{\ddag,\p})$. See Definition \ref{def-deg-lvl}.
    \end{itemize}
Any pair of such tuples $(\Q,\bml,[u]_{(K^\p_{\Phi^\ddag})^{g^\ddag}})$ and $(\Q^\prime,\bml^\prime,[u']_{(K^\p_{\Phi^\ddag})^{g^\ddag}})$ are equivalent in $\bm{\xi}_{(V_\bb{Z},\psi_\bb{Z}),\Phi^\ddag,K^{\ddag,\p}}(S)$, denoted as $(\Q,\bml,[u]_{(K^\p_{\Phi^\ddag})^{g^\ddag}})\sim_{\mrm{iso}}(\Q^\prime,\bml^\prime,[u']_{(K^\p_{\Phi^\ddag})^{g^\ddag}})$, if there is an isomorphism $f:\Q\to\Q^\prime$ such that $f^\vee\circ\lambda^\prime\circ f=\lambda$ and such that $T^\p f\circ u=u^\prime$ modulo the action of $(K^\p_{\Phi^\ddag})^{g^\ddag}$. This also describes the objects in $\bm{\xi}_{(V_\bb{Z},\psi_\bb{Z}),\Phi^\ddag,K^{\ddag,\p}}(S)$ for a general object $S$ in $\ob(Sch/\zbkpp)$, since objects $(\Q,\bml,[u]_{(K^\p_{\Phi^\ddag})^{g^\ddag}})$ are finitely presented.\par 

In fact, for any neat open compact $K^\ddag$, $\xi_{(V_\bb{Z},\psi_\bb{Z}),\Phi^\ddag,K^{\ddag},\bb{C}}$ is representable by $\sh_{K_{\Phi^\ddag},\bb{C}}$. Let us explain this now.
\begin{lem}\label{lem-cons}
Let $(\Q,\bml)$ be a polarized $1$-motive over a $\zbkpp$-scheme $S$. Let $\ca{U}\sbst V^\p\Q$ be an open compact subgroup of $V^\p\Q$ such that $\ca{U}_{\overline{s}}$ is $\pi_1(S,\overline{s})$-invariant for any $\overline{s}\in S$. Then there is a $1$-motive $\Q^\prime$ with a $\zbkppt$-isogeny $f:\Q\to \Q^\prime$ and a $\zbkppt$-polarization $\bml^\prime: \Q'\to \Q^{\prime,\vee}$ such that $f^\vee\circ\bml^\prime\circ f=\bml$.
\end{lem}
\begin{proof}
This is similar to Construction \ref{cons-ab}; see also \cite[Lem. 4.22]{Mum91} and \cite[Cor. 1.3.5.4]{Lan13}. In fact, there is a prime-to-$\p$ positive integer $N$ such that ${N}T^\p\Q\sbst \ca{U}$. Since $\frac{1}{N}\ca{U}/T^{\p}\Q$ is a finite prime-to-$\p$-torsion subgroup of $\Q$, define a homomorphism between $1$-motives $g: \Q\to\Q'$ to be the quotient of $\Q$ by $\frac{1}{N}\ca{U}/T^{\p}\Q$. 
The quotient is a $1$-motive by Lemma \ref{lem-isog-con}. Then $f: \Q\xleftarrow{[N]}\Q\xrightarrow{g}\Q^\prime$ is the desired $\zbkppt$-isogeny.\par
For the second sentence, let $\bml':=f^{\vee,-1}\circ \bml\circ f^{-1}$. See Appendix \ref{comm-lem} for the definition of inverse $\zbkppt$-isogenies. 
\end{proof}
Then we can associate with any $\wp=[(x,p)]\in \sh_{K_{\Phi^\ddag}}(\bb{C})$ a polarized $1$-motive $(\Q_\wp,\bml_\wp)$ with a level structure $[u]_{(K_{\Phi^\ddag}^\p)^{g^\ddag}}$.
\begin{construction}\label{const-1-mot}\upshape
(1) Since $(V_\bb{Z},\psi_\bb{Z})$ is fixed, $x$ determines an integral mixed Hodge structure $\imxh(x)$ of $V_\bb{Z}$. Then $\psi$ or $-\psi$ is a polarization of $\imxh(x)$, which depends on $x$ is in $D_{\Phi^\ddag}^+$ or not. Since $P_{\Phi^\ddag}(\bb{R})$ acts transitively on $\pi_0(D_{\Phi^\ddag})$, we can assume that $x\in D_{\Phi^\ddag}^+$. Then $x$ determines a polarized $1$-motive $(\Q_x,\bml_x)$ such that there is an isomorphism $\alpha_x:T_\bb{Z}\Q_x\xrightarrow{\sim} V_\bb{Z}$; for a different $x^\prime$ in $P_\Phi(\bb{Q})$-orbit of $x$, $\alpha_{x^\prime}$ and $\alpha_{x}$ differ by an action of $P_{\Phi}(\bb{Q})$. That is, if $t\in P_{\Phi^\ddag}(\bb{Q})$ and $tx=x^\prime$, there is a commutative diagram
\begin{equation}\label{eq-dphi-isog}
    \begin{tikzcd}
        T_\bb{Z}\Q_x\otimes \bb{Q}\arrow[rr,"\alpha_x"]\arrow[d,"\sim"]&& V_\bb{Q}\arrow[d,"t"]\\
    T_\bb{Z}\Q_{x^\prime}\otimes\bb{Q}\arrow[rr,"\alpha_{x^\prime}"]&& V_\bb{Q},
    \end{tikzcd}
\end{equation}
such that $V_\bb{Z}$ corresponds to $T_\bb{Z}\Q_x$ (resp. $T_\bb{Z}\Q_{x^\prime}$) under $\alpha_x$ (resp. $\alpha_{x^\prime}$).
\par
Let us fix a choice of $x\in D_{\Phi^\ddag}^+$ and $p\in P_{\Phi^\ddag}(\A)$. Then there is a left-multiplication \begin{equation}\label{eq-lft-mult}c_{p^{\ddag,\p}}:(V_{\zhpp},\psi_{\zhpp})\xrightarrow{p^{\ddag,\p}}(V_{\zhpp}^{(p^{\ddag,\p})},\nu(p^{\ddag,\p})^{-1}\psi_{\App}|_{V_{\zhpp}^{(p^{\ddag,\p})}}),\end{equation}
where $p^{\ddag,\p}$ is the $\App$-component of $p{g^\ddag}$.
Then $c_{p^{\ddag,\p}}V_{\zhpp}$ maps to an open compact subgroup of $V_{\App}\xrightarrow{\alpha^{-1}_x}T_\bb{Z}\Q_x\otimes \App\iso V^\p \Q_x$. By Lemma \ref{lem-cons}, there is a $1$-motive $\Q_\wp$ and a $\zbkppt$-isogeny $f_x: \Q_x\to \Q_\wp$, such that $V^\p f_x\circ c_{p^{\ddag,\p}}(V_{\zhpp})=T^\p\Q_\wp$. Note that $\Q_\wp$ indeed only depends on $\wp$ but not the choice of liftings $(x,p)$. Let $u_\wp:=V^\p f_x\circ c_{p^{\ddag,\p}}|_{V_{\zhpp}}$. Then the 
$(K^\p_{\Phi^\ddag})^{g^\ddag}$-orbit $u^\p:=[u_\wp]_{(K_{\Phi^\ddag}^\p)^{g^\ddag}}$ is the desired level structure.\par
Note that since this process determines a $\zbkppt$-polarization $\bml'$ by \emph{loc. cit.}, there is a unique polarization $\bml_\wp:=r(p^{\ddag,\p})^{-1}\bml'$ for $\Q_\wp$, such that, $u^\p$ sends $\psi_{\zhpp}$ to a $\zhppt$-multiple of $\e^{\bml_\wp}$.\par
(2) Conversely, let us fix any object $(\Q,\bml,[u]_{(K_{\Phi^\ddag}^\p)^{g^\ddag}})$ over $\bb{C}$. 
There is a map $\alpha:T_\bb{Z}\Q\xrightarrow{\sim} V^{(g^\ddag)}_\bb{Z}$, such that $\gr_0 \alpha=\alpha_0^{(g^\ddag),-1}$ and $\gr_{-2}\alpha=\alpha_{-2}^{(g^\ddag),-1}$, and such that $T_\bb{Z}\bml$ corresponds to $\psi_\bb{Z}^{(g^\ddag)}$ under $\alpha$. 
In fact, fix a splitting $\delta_\bb{Z}:\gr_\bullet T_\bb{Z}\Q\to T_\bb{Z}\Q$ and a splitting $\zeta_\bb{Z}:\gr_\bullet^W V_\bb{Z}^{(g^\ddag)}\to V_\bb{Z}^{(g^\ddag)}$, where $W$ is the filtration defined by $Q^\ddag$. Then there is such an $\alpha$, since there is an $\alpha_{-1}: \gr_{-1}T_\bb{Z}\Q\to \gr^W_{-1}V_\bb{Z}$ sending $\bml^{ab}$ to $\psi^{(g^\ddag),ab}_\bb{Z}$ and since $\alpha_0^{(g^\ddag)}$ and $\alpha_{-2}^{(g^\ddag)}$ are determined by the definition of level structures. By Lemma \ref{lem-mx-hdg}, we find a point $x\in D_{\Phi^\ddag}$, which determines an isomorphism of integral mixed Hodge structures $\alpha_x: T_\bb{Z}\Q\xrightarrow{\sim} V_\bb{Z}^{(g^\ddag)}$. \par
Consider the diagram
\begin{equation}\label{diag-level-cons}
    \begin{tikzcd}
 (V_{\App},\psi,Z^{(g^\ddag)}_{\App})\arrow[rr,"u"]\arrow[d,"g^{\ddag,\p}"]&& V^\p\Q\iso T_\bb{Z}\Q\otimes \App\arrow[d,"\alpha_x\otimes 1"]\\
 (V_{\App},\nu(g^{\ddag,\p})^{-1}\psi,W_{\App})\arrow[rr,dashed,"p"]&& (V_{\App},r(g^\ddag)^{-1}\psi,W_{\App}).
    \end{tikzcd}
\end{equation}
The dashed arrow is induced by multiplying some $p\in P_{\Phi^\ddag}(\App)$ since it sends $W$ to $W$, sends $\psi$ to a $\zhppt$-multiple of $\psi$, sends $\gr^W_0 V_{\App}$ to $\gr^W_0 V_{\App}$ identically, and sends $\gr^W_{-2}V_{\App}$ to $\gr^W_{-2}V_{\App}$ by multiplying some element in $\zhppt$. Then we have found a pair $(x,p)$.
\end{construction}
Hence, by Construction \ref{const-1-mot} and Lemma \ref{lem-mx-hdg}, $\sh_{K_{\Phi^\ddag}}(\bb{C})\iso \bm{\xi}_{(V_\bb{Z},\psi_\bb{Z}),\Phi^\ddag,K^{\ddag}}(\bb{C})$. \par 
On the other hand, by \cite[2.2.14-2.2.17]{Mad19}, or by Proposition \ref{prop-deg-lvl} and the proof of \cite[Prop. 6.2.4.7]{Lan13}, assuming that $K_p^\ddag$ is hyperspecial when $\p=\{p\}$, $\bm{\xi}_{(V_\bb{Z},\psi_\bb{Z}),\Phi^\ddag,K^{\ddag,\p}}$ is representable by a separated smooth scheme $\ca{S}_{K_{\Phi^\ddag},\zbkpp}(P_{\Phi^\ddag},D_{\Phi^\ddag})$ over $\zbkpp$. Write \gls{SKPhiddag}$:=\ca{S}_{K_{\Phi^\ddag},\zbkp}(P_{\Phi^\ddag},D_{\Phi^\ddag})$. \par
We then have $\sh_{K_{\Phi^\ddag},\bb{C}}\iso \bm{\xi}_{(V_\bb{Z},\psi_\bb{Z}),\Phi^\ddag,K^{\ddag},\bb{C}}.$ Furthermore, $\sh_{K_{\Phi^\ddag},\bb{Q}}\iso \bm{\xi}_{(V_\bb{Z},\psi_\bb{Z}),\Phi^\ddag,K^{\ddag},\bb{Q}}$ as $\bb{Q}$-schemes, by \cite[2.2.5 and 2.2.6]{Mad19}.\par

\subsubsection{} Recall that $\Phi_0=(Q,X^+_0,g_0)$ maps to $\Phi^\ddag=(Q^\ddag,X^{\ddag,+},g^\ddag)$, that is, 
$Q\sbst Q^\ddag$, $X_0^+\sbst X^{\ddag,+}$, and $g_0=g^\ddag$.\par
The embedding $\iota:(G_0,X_0)\hookrightarrow (G^\ddag,X^\ddag)$ induces a morphism between mixed Shimura data $\iota_{\Phi_0}:(P_{\Phi_0},D_{\Phi_0})\to (P_{\Phi^\ddag},D_{\Phi^\ddag})$, which is an embedding by Corollary \ref{cor-pure-embedding}.\par
Fix a neat open compact $K=K_{0,p}K^p_0$. We can find a neat open compact $K^\ddag=K_p^\ddag K^{p,\ddag}$, such that $K_{0,p}=K^\ddag_p\cap G_0(\bb{Q}_p)$. By Lemma \ref{lem-emb-mix}, $\sh_{K_{\Phi_0,p}}\hookrightarrow \sh_{K_{\Phi^\ddag,p}}$ is an embedding.\par 
For integral models, let $K^\ddag$ and $K_0$ be defined as in \S\ref{subsec-siegel}, and let $\ca{S}_{K_{\Phi_0}}(P_{\Phi_0},D_{\Phi_0})$ be the normalization in $\sh_{K_{\Phi_0}}$ of $\ca{S}_{K_{\Phi^\ddag}}$ under the morphism $\iota_{\Phi_0}:\sh_{K_{\Phi_0}}\to \sh_{K_{\Phi^\ddag}}\hookrightarrow\ca{S}_{K_{\Phi^\ddag}}$. Write \gls{SKPhi0}$:=\ca{S}_{K_{\Phi_0}}(P_{\Phi_0},D_{\Phi_0})$.
Note that the morphism 
$\iota_{\Phi_0}:\ca{S}_{K_{\Phi_0}}\to\ca{S}_{K_{\Phi^\ddag}}$ induced by the normalization is finite, but not a closed embedding in general. When $\Phi_0=(G_0,X_0^+,1)$, a trivial cusp label representative, $\sh_{K_{\Phi_0}}=\sh_{K_0}$. The normalization in $\sh_{K_0}$ of $\ca{S}_{K^\ddag}$ under the morphism $\iota:\sh_{K_0}\to \sh_{K^\ddag}\hookrightarrow \ca{S}_{K^\ddag}$ is denoted by \gls{SK0}, and the generic fiber of $\ca{S}_{K_0}$ is $\sh_{K_0}$.\par 
Let $v|p$ be a place of $E_0$ over $p$, and let $R$ be an $\ca{O}_{(v)}:=\ca{O}_{E_0,(v)}$-algebra. 
For any choice of neat open compact $K^{\ddag,p}$ subgroup containing $K^{p}_0$, 
we can associate an object $(\Q,\bml,[u]_{(K^p_{\Phi^\ddag})^{g^\ddag}})$ in $\bm{\xi}_{(V_\bb{Z},\psi_\bb{Z}),\Phi^\ddag,K^{\ddag,p}}(R)$ for any $R$-valued point in $\ca{S}_{K_{\Phi_0}}(R)$ by pulling back the tautological family via $\iota_{\Phi_0}$.
Define $[u]_{(K_{\Phi_0}^p)^{g^\ddag}}:=\varprojlim_{K_{\Phi_0}^p\sbst K_{\Phi^\ddag}^p}[u]_{(K_{\Phi^\ddag}^p)^{g^\ddag}}$. Then we have associated a tuple $(\Q,\bml,[u]_{(K_{\Phi_0}^p)^{g^\ddag}})$ to any $R$-point. See Construction \ref{const-1-mot} for the case when $R=\bb{C}$.
\subsection{Toroidal compactifications of integral models of Hodge type}\label{subsec-tor-hdg-review}
For the convenience of the readers, let us briefly summarize Chai and Faltings', Lan's and Madapusi's theorems on toroidal compactifications of integral models of Shimura varieties of Siegel/PEL/Hodge type.
\begin{thm}[{Faltings-Chai \cite{FC90}, for Siegel-type case; Lan \cite{Lan13}, \cite{Lan16b} and \cite{Lan17}, for PEL-type case; Madapusi \cite{Mad19}, for Hodge-type case}]\label{thm-tor-hdg}
Let $(G,X)$ be a Shimura datum of Siegel, PEL or Hodge type. Let $\Sigma$ be an admissible polyhedral cone decomposition for $\sh_K(G,X)$. Let $K:=K_pK^p$, where $K_p$ is open compact in $G(\bb{Q}_p)$ and $K^p$ is neat open compact in $G(\Ap)$. Fix a place $v|p$ of the reflex field $E$ of $(G,X)$. For different cases, we have the following assumptions for $\Sigma$:
\begin{enumerate}
    \item If $(G,X)$ is of Siegel or PEL type, we assume that $\Sigma$ is any smooth or projective cone decomposition.
    \item If $(G,X)$ is of Hodge type, we fix an embedding: $\iota:(G,X)\hookrightarrow(G^\ddag,X^\ddag)$. We assume that $\Sigma$ is a refinement of an admissible cone decomposition $\Sigma'$ for $\sh_K$, where $\Sigma'$ is induced by a smooth, or smooth and projective cone decomposition $\Sigma^\ddag$ for $\sh_{K^\ddag}$.
\end{enumerate}
Under the assumptions for $\Sigma$ as above, we have the following statements:
\begin{enumerate}
\item There is a normal algebraic space $\ca{S}_K^\Sigma$ that is proper over $\ca{O}_{E,(v)}$ extending the toroidal compactification $\sh_K^\Sigma$ over $E$ constructed by \cite{Pin89}. The algebraic space $\ca{S}_K^{\Sigma}$ is a projective scheme if $\Sigma$ is induced by a projective $\Sigma^\ddag$; if $\Sigma$ is smooth, $\ca{S}_K^\Sigma$ is smooth if and only if $\ca{S}_K$ is smooth.
\item We have a good stratification $\ca{S}^\Sigma_K=\disju\limits_{\Upsilon}\ca{Z}_{\Upsilon,K}$, where each $\ca{Z}_{\Upsilon,K}$ is a locally closed normal subspace of $\ca{S}^\Sigma_K$ that is flat over $\ca{O}_{E,(v)}$. The strata \gls{ZcaUpsilon} are indexed by the same set of cusp labels with cones with elements of the form $\Upsilon=[(\Phi,\sigma)]\in \cusp_K(G,X,\Sigma)$ defined in the characteristic zero theory. Moreover, $\overline{\ca{Z}}_{\Upsilon,K}=\disju\limits_{\Upsilon^\prime\leq \Upsilon}\ca{Z}_{\Upsilon',K}$.
\item\label{part3-thm-tor-hdg} Fix any $\Upsilon=[(\Phi,\sigma)]$. There is an integral model $\ca{S}_{K_\Phi}$ of the boundary mixed Shimura variety $\sh_{K_\Phi}$, defined over $\ca{O}_{E,(v)}$. There is a tower $\ca{S}_{K_\Phi}\to \overline{\ca{S}}_{K_\Phi}\to \ca{S}_{K_\Phi,h}$ over $\ca{O}_{E,(v)}$ extending the tower $\sh_{K_\Phi}\to \overline{\sh}_{K_\Phi}\to \sh_{K_\Phi,h}$ over $E$, where the first morphism is an $\mbf{E}_{K_\Phi}$-torsor and the second morphism is proper. Denote by $\ca{S}_{K_\Phi}(\sigma)$ the toric scheme defined by the twisted affine torus embedding $\ca{S}_{K_\Phi}(\sigma):=\ca{S}_{K_\Phi}\times^{\mbf{E}_{K_\Phi}}\mbf{E}_{K_\Phi}(\sigma)$. Denote by $\ca{S}_{K_\Phi,\sigma}$ the $\sigma$-stratum of $\ca{S}_{K_\Phi}(\sigma)$.
\item There is an isomorphism $\ca{S}_{K_\Phi,\sigma}\iso \ca{Z}_{[(\cc)],K}$ between the $\sigma$-stratum in $\ca{S}_{K_\Phi}(\sigma)$ and the stratum $\ca{Z}_{[(\Phi,\sigma)],K}$ in $\ca{S}^\Sigma_K$.
Moreover, there is a strata-preserving isomorphism
    $\cpl{\ca{S}_{K_\Phi}(\sigma)}{\ca{S}_{K_\Phi,\sigma}}\iso \cpl{\ca{S}^\Sigma_K}{\ca{Z}_{[(\cc)],K}}$. Here, ``strata-preserving'' has the following meaning: For any affine open formal subscheme $\mathfrak{U}=\spf(A,I)$ of $\cpl{\ca{S}_{K_\Phi}(\sigma)}{\ca{S}_{K_\Phi,\sigma}}\iso \cpl{\ca{S}^\Sigma_K}{\ca{Z}_{[(\cc)],K}}$ which canonically induces morphisms $\mathsf{c}_1:\spec A\to \ca{S}_{K_\Phi}(\sigma)$ and $\mathsf{c}_2:\spec A\to \ca{S}_K^\Sigma$, the two stratifications on $\spec A$ defined by the pullback of the stratification on $\ca{S}_{K_\Phi}(\sigma)$ under $\mathsf{c}_1$ and the pullback of the stratification on $\ca{S}_K^\Sigma$ under $\mathsf{c}_2$ coincide.
\end{enumerate}
\end{thm}
\begin{convention}
To be compatible with the conventions of the previous {subsections}, when $(G,X)$ is of Siegel type, in the theorem above: $(G,X)=(G^\ddag,X^\ddag)$, $K=K^\ddag$, and all the symbols ``$\Phi$'' and ``$\sigma$'' in the theorem above should be added superscripts ``$\ddag$''. Similarly, when $(G,X)$ of Hodge type, then in the theorem above: $(G,X)=(G_0,X_0)$, $K=K_0$, and all ``$\Phi$'' and ``$\sigma$'' above should be added subscripts $0$. That is, $\ca{S}_{K_{\Phi^\ddag}}$ (resp. $\ca{S}_{K_{\Phi_0}}$) introduced in previous {subsections} are the integral models of boundary mixed Shimura varieties $\sh_{K_{\Phi^\ddag}}$ (resp. $\sh_{K_{\Phi_0}}$) described in Part \ref{part3-thm-tor-hdg} of Theorem \ref{thm-tor-hdg} for the compactification $\ca{S}_{K^\ddag}^{\Sigma^\ddag}$ (resp. $\ca{S}_{K_0}^{\Sigma_0}$).
\end{convention}
The assumptions in \cite{FC90} are: (1)$(G,X)=(G^\ddag,X^\ddag)$ is of Siegel type, (2)$K^\ddag=K^\ddag(n)$ is a principal $n$ level subgroup with $(n,p)=1$, and (3)$\Sigma^\ddag$ is smooth with respect to $K^\ddag$ (and is admissible, finite, complete and without self-intersections). In \cite{Lan13}, the $\ca{S}_K^\Sigma$ as above is constructed when $(G,X)$ is of PEL type, $K=K_pK^p$ is defined as Theorem \ref{thm-tor-hdg} with the assumption that $K_p$ is hyperspecial, and $\Sigma$ is smooth with respect to $K$; the hyperspecial assumption for $K_p$ is removed in \cite{Lan16b}. In \cite{Lan17}, the $\ca{S}_K^\Sigma$ as above is constructed for any projective $\Sigma$. In \cite{Mad19}, there is also no assumption on the level $K_{0,p}$ at $p$ when $(G,X)=(G_0,X_0)$ is of Hodge type. In fact, with the same settings as in \S\ref{subsec-siegel}, if $\Sigma=\Sigma_0$ is induced by some smooth cone decomposition $\Sigma^\ddag$ for $(G^\ddag,X^\ddag)$ via the embedding $\iota:(G_0,X_0)\hookrightarrow(G^\ddag,X^\ddag)$, then \gls{SK0Sigma0}, the toroidal compactification of the integral model $\ca{S}_{K_0}$ associated with the cone decomposition $\Sigma_0$, is constructed as the normalization in $\sh_K^\Sigma$ of the toroidal compactification \gls{SKddagSigmaddag} of the Siegel-type integral model $\ca{S}_{K^\ddag}$. When $\Sigma'_0$ is a refinement of $\Sigma_0$, there is also such a $\ca{S}_K^{\Sigma_0'}$ as in the theorem above by \cite[Rmk. 4.1.6]{Mad19}.\par
When $(G,X)=(G^\ddag,X^\ddag)$ is of Siegel type and $K^\ddag$ is chosen as in \S\ref{subsec-siegel}, there is a degenerating family $(\G^\ddag,\lambda^\ddag,[\varepsilon_{\zhp}]_{K^{\ddag,p}})$ of type $\mbf{M}^{\mrm{iso}}_{(V_\bb{Z},\psi_\bb{Z}),K^{\ddag,p}}$ extending the universal family $(\ca{A}^\ddag,\lambda^\ddag|_{\ca{A}},[\varepsilon_{\zhp}]_{K^{\ddag,p}})$ over $\ca{S}_{K^\ddag}$.\par
As we will work over an extension of the reflex field later, let us record the following result:
\begin{prop}\label{prop-normalized-base-change}
Let $F$ be any finite field extension of $E$. Let $\ca{S}_{K,\ca{O}_F}^\Sigma$ be the normalization of the base change $\ca{S}_{K}^\Sigma\otimes_{\ca{O}_{E,(v)}}(\ca{O}_{E,(v)}\otimes_{\ca{O}_E}\ca{O}_F)$. Then the same results in Theorem \ref{thm-tor-hdg} also hold for $\ca{S}^\Sigma_{K,\ca{O}_F}$, with the strata and toric torsors replaced also by the normalizations of their base change from $\ca{O}_{E,(v)}$ to $\ca{O}_{E,(v)}\otimes_{\ca{O}_E}\ca{O}_F$. 
\end{prop}
\begin{proof}
This can be seen from construction, as the arguments in \cite{Mad19} still work if the reflex field is replaced with a finite field extension of it, but let us give a direct proof here. We temporarily denote all normalizations of the base changes of $\ca{O}_{E,(v)}$-algebraic spaces to $\ca{O}_F\otimes_{\ca{O}_E}\ca{O}_{E,(v)}$ by adding subscript $\ca{O}_F$ and denote corresponding usual base changes by adding $[\ca{O}_F]$. The first assertion in Theorem \ref{thm-tor-hdg} is clear; the third one is true because torus parts commute with normalization and base change. There is a quasi-finite morphism $\ca{Z}_{\Upsilon,K,\ca{O}_F}\to \ca{S}_{K,\ca{O}_F}^\Sigma$ with locally closed image denoted by $[\ca{Z}_\Upsilon]$; the induced morphism $\ca{Z}_{\Upsilon,K,\ca{O}_F}\to[\ca{Z}_\Upsilon]$ is finite and birational. 
We just need to show that $[\ca{Z}_\Upsilon]$ is normal and isomorphic to $\ca{Z}_{\Upsilon,K,\ca{O}_F}$, and that $\cpl{\ca{S}_{K_\Phi,\ca{O}_F}(\sigma)}{\ca{S}_{K_\Phi,\ca{O}_F,\sigma}}\iso \cpl{\ca{S}^\Sigma_{K,\ca{O}_F}}{[\ca{Z}_\Upsilon]}$. \par
By \cite[Sec. 3.1]{LS18b}, the assertions except for the normality are true for the usual base change to $\ca{O}_{E,(v)}\otimes_{\ca{O}_E}\ca{O}_F$. Let $\geom{x}$ be any $\overline{\bb{F}}_p$-point of $[\ca{Z}_\Upsilon]$ which maps to $\geom{y}$ in $\ca{Z}_{\Upsilon,K,[\ca{O}_F]}\hookrightarrow\ca{S}_{K,[\ca{O}_F]}^\Sigma$ and lifts to an $\overline{\bb{F}}_p$-point $\geom{x}_0$ in $\ca{Z}_{\Upsilon,K,\ca{O}_F}$. This induces a finite homomorphism from the complete local ring at $\geom{y}$, $\wat{\ca{O}}_{\geom{y}}$, to the complete local ring at $\geom{x}$, $\wat{\ca{O}}_{\geom{x}}$. Moreover, from \emph{loc. cit.} we have $\wat{\ca{O}}_{\geom{y}}\iso \wat{\ca{O}}_{\ca{S}_{K_\Phi,[\ca{O}_F]}(\sigma),\geom{y}}$ (since pullbacks still have a good description of formal completions).
Suppose that $\geom{y}$ maps to a point $\geom{z}$ in $\ca{S}_{K_\Phi,\sigma,[\ca{O}_F]}$ under the canonical projection to $\sigma$-stratum. Let $\wat{\ca{O}}^*_{\geom{z}}$ be the normalization in $\wat{\ca{O}}_{\geom{x}}$ of the complete local ring $\wat{\ca{O}}_{\ca{S}_{K_\Phi,\sigma,[\ca{O}_F]},\geom{z}}$ via the ring homomorphism $\wat{\ca{O}}_{\ca{S}_{K_\Phi,\sigma,[\ca{O}_F]},\geom{z}}\to \wat{\ca{O}}_{\geom{y}}\iso \wat{\ca{O}}_{\ca{S}_{K_\Phi,[\ca{O}_F]}(\sigma),\geom{y}}\to \wat{\ca{O}}_{\geom{x}}$. Then we have $\wat{\ca{O}}_{\ca{Z}_{\Upsilon,K,\ca{O}_F},\geom{x}_0}\iso \wat{\ca{O}}^*_{\geom{z}}$ by Zariski's main theorem. This implies that $\wat{\ca{O}}_{\ca{Z}_{\Upsilon,K,\ca{O}_F},\geom{x}_0}\iso \wat{\ca{O}}_{[\ca{Z}_\Upsilon],\geom{x}}$; therefore, $[\ca{Z}_{\Upsilon}]\iso \ca{Z}_{\Upsilon,K,\ca{O}_F}$. Moreover, we have $\wat{\ca{O}}^*_{\geom{z}}\otimes_{\wat{\ca{O}}_{\ca{S}_{K_\Phi,\sigma,[\ca{O}_F]},\geom{z}}}\wat{\ca{O}}_{\geom{y}}\iso \wat{O}_{\geom{x}}$ by Zariski's main theorem again. Thanks to \cite[Lem. A.3.2]{Mad19}, we now obtain the desired description of the formal completion at the boundary.
\end{proof}
\newpage
\section{Twisting constructions}\label{sec-twt-construction}
We continue with the notation in {Section} \ref{sec-1-motives-and-deg}. Denote by $G^\ad_0$ the adjoint group of $G_0$. Denote by $G^\ad_0(\bb{R})_1$ the subgroup of $G^\ad_0(\bb{R})$ that stabilizes the subspace $X_0$ in $X_0^\ad$. Denote $G^\ad_0(\bb{Q})_1:=G^\ad_0(\bb{R})_1\cap G^\ad_0(\bb{Q})$. \par
By \cite{Del79}, there is a left action of $G^\ad_0(\bb{Q})_1$ on $\sh(G_0,X_0)$. The goal of this {section} is to construct ``moduli descriptions'' (called ``twisting'' in \cite{Kis10} and \cite{KP15}) in some special cases for this action on (integral models of) compactifications. From \S\ref{twt-ab-sch} to \S\ref{subsec-twt-1-mot}, we extend the twisting construction of Kisin-Pappas \cite{KP15} to boundaries. In \S\ref{subsec-ext} and in a different setting, we construct a ``twisting'' related to the quasi-isogeny twisting of Lan \cite{Lan16b} in the Siegel case and conclude that twisting constructions admit extensions to toroidal compactifications (if the cone decompositions are compatible with the actions). 
\subsection{Twisting abelian schemes}\label{twt-ab-sch}
We explain Kisin and Pappas' twisting construction for abelian schemes and for algebraic tori. The main difference is that the original construction due to Kisin and Pappas is established on the level of abelian schemes up to prime-to-$p$ isogenies. 
All materials are from \cite[Sec. 3.1 and 3.2]{Kis10} and \cite[Sec. 4.4 and 4.5]{KP15}. 
\subsubsection{}\label{subsubsec-twt-general-sheaf}
Let $Z$ be a faithfully flat affine
group scheme of finite type over $A$, where $A$ is a commutative ring. For any $A$-module $M$, a $Z$-action on $M$ is defined by a morphism between fppf sheaves $\z:Z\to \aut_AM$. Then $\z$ induces a map $M\xrightarrow{\z^*}M\otimes_A\ca{O}_Z$ that endows $M$ with an $\ca{O}_Z$-comodule structure.
Let $\ca{P}$ be a $Z$-torsor over $A$. Then $M\otimes_A \ca{O}_\ca{P}$ is an $\ca{O}_\ca{P}$-module with semilinear $Z$-action. Then $(M\otimes_A\ca{O}_\ca{P})^Z$, the submodule of $N:=M\otimes_A \ca{O}_\ca{P}$ fixed by $Z$, can be viewed as the equalizer of the natural $\ca{O}_Z$-comodule structure $\mbf{z}^*:N\to N\otimes_A\ca{O}_Z$ of $N$ and the map $p^*:N\to N\otimes_A\ca{O}_Z$ determined by the natural projection $Z\times_{\spec A} M\to M$.
By \cite[Lem. 4.4.3]{KP15},
\begin{equation}\label{twt-comp}
    (M\otimes_A\ca{O_P})^Z\otimes_A \ca{O_P}\iso M\otimes_A\ca{O_P}.
\end{equation}
Following \cite{Kis10} and \cite{KP15}, define the \textbf{twisting of $M$ by $\ca{P}$} by an exact functor $$\tw{(-)}:M\mapsto\tw{M}:=(M\otimes_A\ca{O_P})^Z,$$ where the exactness follows from (\ref{twt-comp}) and the faithful flatness of $Z$ (and $\ca{O}_{\ca{P}}$) over $A$.\par
Now let $A=\zbkp$, let $S$ be a $\spec\zbkp$-scheme and also assume that $Z$ is commutative. Set $\ca{O}_{F,(p)}:=\ca{O}_F\otimes\bb{Z}_{(p)}$. Since $\ca{P}$ is flat of finite type over $\zbkp$, by \cite[I, Thm. 1.6]{MB89}, $\ca{P}(\ca{O}_{F,(p)})$ is not empty for some finite Galois extension $F$ over $\bb{Q}$. For this, see also the proof of \cite[Lem. 4.4.6]{KP15}.\par
Since $Z(\ca{O}_{F,(p)})$ acts on $\ca{P}(\ca{O}_{F,(p)})$ freely and transitively, for any $\tilde{\gamma}\in\ca{P}(\ca{O}_{F,(p)})$ and any $\sigma\in\mrm{Gal}(F/\bb{Q})$, there exists a unique $z_\sigma^{\tilde{\gamma}}\in Z(\ca{O}_{F,(p)})$ such that  $\sigma\tilde{\gamma}=z_\sigma^{\tilde{\gamma}}\cdot\tilde{\gamma}$. Since $z_{\sigma\tau}^{\tilde{\gamma}}\tilde{\gamma}=\sigma\tau\tilde{\gamma}=\sigma(z_\tau^{\tilde{\gamma}}\cdot \tilde{\gamma})=\sigma(z_\tau^{\tilde{\gamma}})\cdot \sigma\tilde{\gamma}=\sigma(z_\tau^{\tilde{\gamma}})\cdot z_\sigma^{\tilde{\gamma}}\cdot \tilde{\gamma}$, we have associated with $\tilde{\gamma}$ a $1$-cocycle $\sigma\mapsto z_\sigma^{\tilde{\gamma}}$. 
For another $\tilde{\gamma}^\prime\in\ca{P}(\ca{O}_{F,(p)})$, there is a unique $z\in Z(\ca{O}_{F,(p)})$ such that $z\tilde{\gamma}=\tilde{\gamma}^\prime$, and $\sigma\tilde{\gamma}^\prime=\sigma(z)z_\sigma^{\tilde{\gamma}}\cdot \tilde{\gamma}=\sigma(z)z^{-1}z_\sigma^{\tilde{\gamma}}\cdot\tilde{\gamma}^\prime$. This implies that $z_\sigma^{\tilde{\gamma}^\prime}=\sigma(z)z^{-1}z_\sigma^{\tilde{\gamma}}$.\par
Specializing (\ref{twt-comp}) at $\tg$, we have an isomorphism 
\begin{equation}
    \iota^{\tg}: M^{\ca{P}}\otimes_{\zbkp}\ca{O}_{F,(p)}\xrightarrow{\sim}M\otimes_{\zbkp}\ca{O}_{F,(p)}.
\end{equation}
The Galois group $\gal(F/\bb{Q})$ acts on $\ca{O}_{F,(p)}$ on the left-hand side and acts on both factors on the right-hand side. The action of $\sigma\in\gal(F/\bb{Q})$ on $M\otimes_\zbkp \ca{O}_{F,(p)}$ is via $\sigma\mapsto z^{\tg}_\sigma\otimes \sigma\mapsto \mbf{z}(z^{\tg}_\sigma)\otimes \sigma$, and we will omit $\mbf{z}$ if it is clear in the context.
In particular, we have $M^{\ca{P}}\iso (M\otimes_{\zbkp}\ca{O}_{F,(p)})^{\gal(F/\bb{Q})}$.\par
\subsubsection{}\label{subsubsec-twt-vec-func}
Let $R$ be a faithfully flat algebra over $\bb{Z}$. Denote $R_{(p)}:=R\otimes \zbkp$. Suppose that $M$ is a free $R_{(p)}$-module of finite rank and that $M_R$ is an $R$-lattice in $M$.
Let $f^{\tg}:M_R\otimes\ca{O}_F\otimes\zbkp\to(M_R\otimes{\ca{O}_F})^{\oplus |\gal(F/\bb{Q})|-1}\otimes\zbkp$ be a homomorphism between $\zbkp$-modules defined by $f^{\tg}:=\bigoplus\limits_{\sigma\in\gal(F/\bb{Q})\text{ and }\sigma\neq\mrm{id}}(z_\sigma^{\tilde{\gamma}}\sigma-\mrm{id})$. Then there is a sufficiently large prime-to-$p$ integer $N>0$, such that $Nz_\sigma^{\tg}\sigma$ is an endomorphism of $M_R\otimes \ca{O}_F$ for all $\sigma\in\gal(F/\bb{Q})$. Define $f^{\tg}_N:=N\circ f^{\tg}|_{M_R\otimes \ca{O}_F}$, and define $M_R^{\tg}:=\ker f^{\tg}_N$. \par
Let $V$ and $W$ be two finite free $R_{(p)}$-modules as $M$ above. Suppose that there is an $R_{(p)}$-isomorphism $g:V\to W$. Suppose that $V$ (resp. $W$) is equipped with a $Z$-action denoted by $\z^V$ (resp. $\z^W$).
We say $\z^V$ and $\z^W$ are compatible with $g$ if, for any $\zbkp$-algebra $R_0$ and any $z\in Z(\ca{O}_{F,(p)})$, $\z^W(z)=g\circ \z^V(z)\circ g^{-1}\in \Endo_{R_{(p)}\otimes_{\zbkp} R_0}(W\otimes_{\zbkp} R_0)$. We summarize the following observation:
\begin{lem}\label{twt-isog}
Let $V_{R}$ (resp. $W_{R}$) be an $R$-lattice of $V$ (resp. $W$). Let $g:V_{R}\to W_{R}$ be a homomorphism between $R$-lattices of $V$ and $W$ such that $f\otimes\zbkp$ is an isomorphism. Let $F$ be a Galois extension of $\bb{Q}$. 
Suppose that $W$ is equipped with a compatible $Z$-action induced by $g$ as above. Fix any $\tg\in \ca{P}(\ca{O}_{F,(p)})$. Let $N$ be any positive prime-to-$p$ integer such that $N\z^V(z_\sigma^{\tg})$ (resp. $N\z^W(z_\sigma^{\tg})$) lies in $\Endo(V_R\otimes{\ca{O}_F})$ (resp. $\Endo(W_R\otimes{\ca{O}_F})$) for any $\sigma\in\gal(F/\bb{Q})$. 
Then there is a unique homomorphism $f_N^{W,\tg}:W_R\otimes{\ca{O}_F}\to W_R\otimes{\ca{O}_F}^{\oplus|\mrm{Gal}(F/\bb{Q})|-1}$ such that the diagram
\begin{equation}\begin{tikzcd}
V_R\otimes{\ca{O}_F}\arrow[rr,"g"]\arrow[d,"f_N^{\tg}"]&&W_R\otimes{\ca{O}_F}\arrow[d,"f^{W,\tg}_N"]\\
(V_R\otimes{\ca{O}_F})^{\oplus|\mrm{Gal}(F/\bb{Q})|-1}\arrow[rr,"g^{\oplus |\mrm{Gal}(F/\bb{Q})|-1}"]&&(W_R\otimes{\ca{O}_F})^{\oplus| \mrm{Gal}(F/\bb{Q})|-1}
\end{tikzcd}\end{equation}commutes. By taking kernels of the two vertical arrows of the diagram above, there is a natural homomorphism  $g^{\tg}:V^{\tg}_{R} \to W^{\tg}_{R}$ which induces an isomorphism after tensoring $\zbkp$. For different choices of $N$, $g^{\tg}$ varies by isomorphisms induced by multiplying prime-to-$p$ integers.\par 
If $g$ is defined by multiplying $r\in R^\times$, then so is $g^{\tg}$.
\end{lem}
\begin{proof}
In fact, we let $f^{W,{\tg}}_N:=\bigoplus\limits_{\sigma\in\mrm{Gal}(F/\bb{Q})\text{ and }\sigma\neq\mrm{id}}((N\z^W(z_\sigma^{\tg})\sigma-N\cdot\mrm{id})$ since $\mbf{z}^V$ and $\mbf{z}^W$ are compatible. The homomorphism $g^{\tg}\otimes \zbkp$ is an isomorphism since $g\otimes \zbkp$ is so.
\end{proof}
\subsubsection{}
Denote by $S$ a $\zbkp$-scheme and denote by $\ca{A}$ a commutative group scheme over $S$. Let $\ca{N}$ be a finite free $\bb{Z}$-algebra. Suppose that there is an $\ca{N}$-action on $\ca{A}$, that is, there is a homomorphism $m: \ca{N}\to \uend_{S}(\ca{A})$. Let $\ca{M}$ be a finite $\ca{N}$-module. Define \emph{Serre's construction} to be $\uhom_\ca{N}(\ca{M},\ca{A})^\circ$, the fiberwise identity component of the {\'e}tale sheaf $\uhom_\ca{N}(\ca{M},\ca{A})$.\par 
If $\ca{A}$ is an abelian scheme, $\uhom_\ca{N}(\ca{M},\ca{A})^\circ$ is representable by an abelian scheme over $S$. When $\ca{M}$ is a finite projective $\ca{N}$-module, $\uhom_\ca{N}(\ca{M},\ca{A})$ is fiberwise connected and therefore is representable by an abelian scheme over $S$. See, e.g., \cite[Prop. 5.2.3.9]{Lan13}.\par 
Assuming that $\ca{M}$ is finite projective, we can write $\ca{A}^{\ca{M}^\vee}:=\uhom_\ca{N}(\ca{M},\ca{A})\iso \ca{M}^\vee\otimes_{\ca{N}}\ca{A}$, where $\ca{M}^\vee:=\uhom_{\ca{N}}(\ca{M},\ca{N})$. For any scheme $T$ over $S$, $\ca{A}^{\ca{M}^\vee}(T)=\Hom_{\ca{N}}(\ca{M},\ca{A}(T))\iso \ca{M}^\vee\otimes_{\ca{N}}\ca{A}(T)$.\par 
For our purpose, we can further assume that $\ca{N}=\bb{Z}$ and that $\ca{M}$ is a finite free $\bb{Z}$-algebra of finite rank. Then $\ca{F}^{\ca{M}^\vee}:=\uhom(\ca{M},\ca{F})=\ca{F}\otimes \ca{M}^\vee$ for any fppf sheaf $\ca{F}$ over $S$; this is just $\ca{F}^{\oplus r}$ equipped with an action of $\ca{M}^\vee$ (or rather, $\ca{M}$). In this case, $\ca{F}^{\ca{M}^\vee}$ is representable by a torus (resp. an {\'e}tale locally constant sheaf of finite rank, resp. an abelian scheme, resp. a semi-abelian scheme) if $\ca{F}$ is representable by a torus (resp. an {\'e}tale locally constant sheaf of finite rank, resp. an abelian scheme, resp. a semi-abelian scheme). If $\Q=[\ul{Y}\to\G^\natural]$ is a $1$-motive, so is $\Q^{\ca{M}^\vee}:=[\ul{Y}^{\ca{M}^\vee}\to \G^{\natural,\ca{M}^\vee}]$.\par
Now suppose that we are given a $1$-motive $\Q$ and a finite Galois extension $F$ of $\bb{Q}$. There is a canonical perfect reduced trace pairing between $\ca{O}_F$ and its inverse different $\mrm{Diff}^{-1}_{\ca{O}_F/\bb{Z}}$:
$$\mrm{Tr}_{F/\bb{Q}}:\ca{O}_F\times \mrm{Diff}^{-1}_{\ca{O}_F/\bb{Z}}\to \bb{Z}.$$
This perfect pairing induces an isomorphism $\ca{O}_F^\vee:=\Hom_\bb{Z}(\ca{O}_F,\bb{Z})\xrightarrow{\sim}\mrm{Diff}^{-1}_{\ca{O}_F/\bb{Z}}\sbst F.$\par
Let $\bml:\Q\to \Q^\vee$ be an isogeny.
Then $$\bml\otimes 1: \Q\otimes \ca{O}_F\lra\Q^\vee\otimes\ca{O}_F$$
is an isogeny. Any embedding of $\bb{Z}$-modules $i_F:\ca{O}_F\hookrightarrow \ca{O}_F^\vee$ 
induces an isogeny
$$(\bml\otimes 1)(i_F): \Q\otimes \ca{O}_F\xrightarrow{\bml\otimes 1}\Q^\vee\otimes\ca{O}_F\xrightarrow{i_F} \Q^\vee\otimes\ca{O}_F^\vee$$
between $\Q\otimes \oo$ and its Cartier dual $\Q^\vee\otimes \ca{O}_F^\vee$.
\subsubsection{}\label{twt-abelian}
Let $F$, $\tg$ and $z^{\tg}_\sigma$ be defined as in \S\ref{subsubsec-twt-general-sheaf}. Let $\ca{A}$ be an abelian scheme or an algebraic torus over $S$. Let $[\ca{A}]:=\ca{A}\otimes \zbkp$: it is the sheaf that represents $\ca{A}$ up to prime-to-$p$ isogenies. More precisely, for any $S$-scheme $T$, define the $T$-value points of $[\ca{A}]$ to be $[\ca{A}](T):=\mrm{Hom}_S(T,\ca{A})\otimes_{\bb{Z}}\zbkp$.
Assume that $Z$ acts on $[\ca{A}]$ via a morphism $\z:Z\to \aut_{\zbkp}([\ca{A}])$: For any $\zbkp$-algebra $R$, there is a homomorphism $\z(R): Z(R)\to \aut_{\zbkp}([\ca{A}])(R):=(\Endo_S(\ca{A})\otimes_\bb{Z} R)^\times$.
From the discussion above, we have a map $\iota_{\gal}:\mrm{Gal}(F/\bb{Q})\to \aut_{\zbkp}([\ca{A}])(\ca{O}_{F,(p)})$ defined by $\sigma\mapsto z^{\tilde{\gamma}}_\sigma\mapsto \z(z^{\tilde{\gamma}}_\sigma)$.
The homomorphism $\z(\ca{O}_{F,(p)}):Z(\op)\to \aut_{\zbkp}([\ca{A}])(\op)$ is equivariant with respect to the actions induced by the Galois action of $\mrm{Gal}(F/\bb{Q})$ on $\ca{O}_F$.\par 
Since $\ca{O}_{F}$ is finite free over $\bb{Z}$, $\aut_{\zbkp}([\ca{A}])(\ca{O}_{F,(p)})=(\Endo_{S}(\ca{A})\otimes_{\bb{Z}} \ca{O}_F\otimes_{\bb{Z}}\zbkp)^\times\iso (\Endo_{S}(\ca{A}^{\ca{O}_F})\otimes_\bb{Z}\zbkp)^\times\iso \aut_{\zbkp}([\ca{A}^{\ca{O}_F}])(\zbkp)$. So the map $\iota_\gal$ above is equivalent to a map $$\iota_{\gal}:\mrm{Gal}(F/\bb{Q})\lra \aut_\zbkp([\ca{A}^{\ca{O}_F}])(\zbkp).$$
The isomorphism (\ref{twt-comp}) can be read as $$[\ca{A}]^{\ca{P}}\otimes_{\bb{Z}_{(p)}}\ca{O}_{\ca{P}}\iso [\ca{A}]\otimes_{\bb{Z}_{(p)}}\ca{O}_{\ca{P}}.$$
Evaluating at $\tilde{\gamma}$, we have an isomorphism between fppf sheaves over $S$:
\begin{equation}\label{twt-comp-a}\iota^{\tg}:[\ca{A}]^{\ca{P}}\otimes_{\bb{Z}_{(p)}}\ca{O}_{F,(p)}\iso [\ca{A}]\otimes_{\bb{Z}_{(p)}}\ca{O}_{F,(p)}.\end{equation}
As discussed in \S\ref{subsubsec-twt-general-sheaf}, we have that $[\ca{A}]^{\ca{P}}\iso(\ca{A}\otimes\ca{O}_{F,(p)})^{\mrm{Gal}(F/\bb{Q})}.$\par
Since $\gal(F/\bb{Q})$ is a finite group, there is a positive prime-to-$p$ integer $N$ such that $N\iota_{\gal}(\sigma)\in \Endo_S(\ca{A}^{\ca{O}_F})$ for all $\sigma\in \gal(F/\bb{Q})$. 
There is a homomorphism between abelian schemes or algebraic tori
$$f^{\tg}_N:\ca{A}^{\ca{O}_F}\lra(\ca{A}^{\ca{O}_F})^{\oplus |\mrm{Gal}(F/\bb{Q})|-1},$$
defined by $f^{\tg}_N:=\bigoplus\limits_{\sigma\in\mrm{Gal}(F/\bb{Q})\text{ and }\sigma\neq \mrm{id}}((Nz_\sigma^{\tilde{\gamma}})\sigma-N\cdot\mrm{id}).$ 
Define $\ca{A}^{\tilde{\gamma}}_N:=(\ker f^{\tg}_N)^\circ$ to be the fiberwise (geometrically) identity component of $\ker f^{\tg}_N$. Hence, $\ker f^{\tg}_N$ and $\ca{A}^{\tilde{\gamma}}_N$ are representable by closed subgroup schemes of $\ca{A}^{\ca{O}_F}$.\par 
Let $f^{\tg}:\ca{A}^{\ca{O}_F}\otimes\zbkp\to(\ca{A}^{\ca{O}_F})^{\oplus |\mrm{Gal}(F/\bb{Q})|-1}\otimes\zbkp$ be a morphism between sheaves defined by $f^{\tg}:=\bigoplus\limits_{\sigma\in\mrm{Gal}(F/\bb{Q})\text{ and }\sigma\neq\mrm{id}}(z_\sigma^{\tilde{\gamma}}\sigma-\mrm{id})$. Since $(p,N)=1$, $[\ca{A}]^{\ca{P}}\iso\ker f^{\tg}\iso\ker Nf^{\tg}\iso\ker f^{\tg}_N\otimes\zbkp$.\par
For different points in $\ca{P}(\ca{O}_{F,(p)})$, we have the following lemma:
\begin{lem}\label{another-gamma} Let $\tg\in \ca{P}(\op)$. Suppose that $\tilde{\gamma}^\prime=z\tilde{\gamma}$ for some $z\in Z(\op)$. Let $N'$ be an integer such that $N^\prime z\in \Endo_S(\ca{A}^{\ca{O}_F})$ and $(N^\prime,p)=1$. 
The homomorphism $N^\prime z\in \Endo_S(\ca{A}^{\ca{O}_F})$ induces a prime-to-$p$ isogeny between $\ker f^{{\tg}^\prime}_N$ and $\ker f^{\tg}_N$.
\end{lem}
\begin{proof}
This is true because $(Nz_\sigma^{\tilde{\gamma}}\sigma-N)\cdot N^\prime z=N^\prime z\cdot(Nz^{\tilde{\gamma}^\prime}_\sigma\sigma-N)$, where $z_\sigma^{\tilde{\gamma}^\prime}=\sigma(z)z^{-1}z_\sigma^{\tilde{\gamma}}$ for any $\sigma\in\mrm{Gal}(F/\bb{Q})$.
\end{proof}
\subsubsection{}
Since $[\ca{A}]^{\ca{P}}\iso \ker f^{\tg}_N\otimes\zbkp$, the isomorphism (\ref{twt-comp-a}) is the base change to $\zbkp$ of the following composition of morphisms between schemes $$\ion:\ker f^{\tg}_N\otimes\ca{O}_{F}\lra \ca{A}^{\ca{O}_F}\otimes\ca{O}_{F}\lra \ca{A}^{\ca{O}_F}\otimes_{\ca{O}_F}\ca{O}_F\iso \ca{A}^{\ca{O}_F},$$
where the first arrow is induced by the embedding $\ker f^{\tg}_N\hookrightarrow \ca{A}^{\ca{O}_F}$ and the second arrow is induced by the multiplication of $\ca{O}_F$ to the second factor of $\ca{A}^{\ca{O}_F}\iso \ca{A}\otimes \ca{O}_F$.\par
Now we explain the proof of the following proposition:
\begin{prop}\label{twt-kp}With the conventions above, the closed subgroup scheme $\ca{A}^{\tilde{\gamma}}_N$ of $\ca{A}^{\ca{O}_F}$ is an abelian scheme (resp. an algebraic torus) over $S$ if $\ca{A}$ is an abelian scheme (resp. an algebraic torus) over $S$.
\end{prop}
\begin{proof}
Since this question is Zariski local on the base, and all the objects and morphisms being studied are finitely presented, we assume that $S$ is affine and of finite type over $\zbkp$. 
If $\ca{A}$ is an abelian scheme, we see that $\ca{A}^{\tilde{\gamma}}_N$ is proper since it is a closed subgroup scheme of $\ca{A}^{\ca{O}_F}$, and that it is fiberwise geometrically connected by the construction. Hence, for either case where $\ca{A}$ is an abelian scheme or an algebraic torus, we only need to show that $\ca{A}_N^{\tilde{\gamma}}$ is smooth. By the fiber-by-fiber criterion of flatness \cite[\href{https://stacks.math.columbia.edu/tag/039E}{Lem. 039E}]{stacks-project}, to show the flatness of $f^{\tg}_N$, it suffices to show the fibers of $f^{\tg}_N$ at all $x\in S$ are flat; by generic flatness \cite[\href{https://stacks.math.columbia.edu/tag/052A}{Prop. 052A}]{stacks-project} and the group scheme structures on $\ca{A}^{\ca{O}_F}_x$, the homomorphisms $f_{\tilde{\gamma},N,x}$ are flat. 
We see that the homomorphism $f^{\tg}_N$ is flat, so is its kernel. 
Then we only need to show that $\ker f^{\tg}_N$, or equivalently, $\ca{A}_{N}^{\tilde{\gamma}}$, is fiberwise smooth, and it suffices to check it at all closed points of $S$. Note that for any point $x^\prime\in S$ whose residue field is of characteristic $0$, $\ker f_{N,x^\prime}^{\tg}$ and $\ca{A}^{\tilde{\gamma}}_{N,x^\prime}$ are smooth by Cartier's theorem; see \cite[\href{https://stacks.math.columbia.edu/tag/047N}{Lem. 047N}]{stacks-project}.\par
For any closed point $x\in S$, the residue field $\kappa_x$ of $x$ is perfect. We write $K:=\ker f^{\tg}_N\otimes \ca{O}_F$ and $A:=\ca{A}^{\ca{O}_F}$ for simplicity. We claim that $\ionx:K_x\to A_{x}$, the fiber of $\ion$ at $x$, is an isogeny, i.e., is surjective with finite kernel. Granting this, we will know that $\ker \ionx$ is a finite group. Since $\ker\ionx\otimes\zbkp$ is trivial, $\ionx$ is a prime-to-$p$ isogeny. Then the identity component of $K_x$ is an abelian variety or algebraic torus and in particular is smooth. So $K$ is smooth. As $K^\circ$ is isomorphic to a product of copies of $\ca{A}^{\tilde{\gamma}}_N$, the latter is also smooth.\par 
Now we show the claim in the last paragraph. We apply Lemma \ref{lem-ab-ker} to $\iota_\red:K_{x,\red}^\circ\hookrightarrow K_x\xrightarrow{\iota_{N,x}^{\tg}} A_x$.
If $(\ker \iota_\red)_\red^\circ$ is non-trivial, it will have nontrivial $p^n$-torsion points since it is an abelian variety or a torus. But this gives us a contradiction since $\ker \iota^{\tg}_{N,x}\otimes \bb{Z}_{(p)}$ is trivial. Applying Lemma \ref{lem-ab-ker} to $f^{\tg}_N$, we see that $\ker \ionx$ is finite. Then $\ionx(K):=K/\ker\ionx$ can be decomposed as $K_1\cdot K'$, where $K_1$ is a finite group and $K'$ is the image of $K_{x,\mrm{red}}^\circ$ in $A_x$ via $\ionx$. If $K'\neq A_x$, there is an abelian subvariety or a subtorus $B$ in $A_x$ such that $K^\prime\cap B$ is finite and $K' \cdot B=A$. 
Then the cokernel sheaf $\coker \ionx$ of $\ionx$ is representable by a quotient of $B$ by a finite group variety, and therefore is representable by an abelian variety or a torus. Since abelian varieties or tori have non-trivial $p^n$-torsion points, $\coker \ionx\otimes\zbkp$ is nontrivial, which leads to a contradiction. So $\iota^{\tg}_{N,x}$ is surjective.  
\end{proof}
\begin{lem}\label{lem-ab-ker}
Let $\kappa$ be a field. Let $A$ and $B$ be two abelian varieties (resp. two algebraic tori) over $\kappa$. Let $f:A\to B$ be a homomorphism with kernel $K$. Let $K_{\mrm{red}}$ be the reduced scheme associated with $K$. Suppose that $K_{\mrm{red}}$ is geometrically reduced (for example, when $\kappa$ is perfect), and therefore is endowed with a group scheme structure induced by $K$. Then: 
\begin{enumerate}
\item $K$ satisfies an extension
$$0\lra K_1\lra K\lra K_2\lra 0,$$
where $K_1$ is a finite flat group variety and $K_2$ is an abelian variety (resp. algebraic torus). 
\item Denote by $K^\circ_{\mrm{red}}$ the identity component of $K_{\mrm{red}}$. Then $K^\circ_{\mrm{red}}$ is an abelian variety (resp. algebraic torus) and the embedding $K_{\mrm{red}}\hookrightarrow K$ induces an isogeny $K_{\mrm{red}}^\circ\to K_2$.
\end{enumerate}
\end{lem}
\begin{proof}Without loss of generality, assume that $f$ is surjective.
If $A$ is an abelian variety (resp. torus), $K_\mrm{red}^\circ$ is an abelian subvariety (resp. a subtorus) of $A$. Then there is an abelian subvariety (resp. a subtorus) $A'$ of $A$ such that $K_{\mrm{fin}}:=K_\mrm{red}^\circ\cap A^\prime$ is finite and such that $K_{\mrm{red}}^\circ\cdot A^\prime=A$; it follows from the Poincar{\'e} reducibility theorem in the abelian variety case. The quotient $K_2:=K_{\mrm{red}}^\circ/K_{\mrm{fin}}$ is an abelian variety (resp. a torus). Then there is a surjective homomorphism $\tilde{f}:A\to K^\circ_\red/K_{\mrm{fin}}\times_\kappa A'/K_{\mrm{fin}}\to K_2\times_\kappa B$ between abelian varieties (resp. tori), which is an isogeny since the dimensions of the source and the target are the same. Let $p_2$ be the projection $p_2:K_2\times_\kappa B$ to the second factor. Then $f=p_2\circ\tilde{f}$. Define $K_1:=\ker \tilde{f}$. Then we have the desired results.
\end{proof}
From the proof of Proposition \ref{twt-kp}, we see that the cokernel sheaf of the inclusion $\ca{A}_N^{\tg}\hookrightarrow\ker\iota^{\tg}_N$ is finite and the orders of its fibers are prime to $p$. Then $\ca{A}_N^{\tg}\otimes\zbkp\iso \ker\iota^{\tg}_N\otimes \zbkp$. Combining this with the discussion before Lemma \ref{another-gamma}, we have
\begin{cor}\label{twt-iso-cor}
The morphism $\iota^{\tg}_N$ is a prime-to-$p$ isogeny and $\ca{A}^{\tilde{\gamma}}_N\otimes\zbkp\iso[\ca{A}]^\ca{P}$ as fppf sheaves over $S$. Moreover, $\iota^{\tg}$ is the base change to $\zbkp$ of the following $\zbkpt$-isogeny $$\ca{A}^{\tg,\ca{O}_F}_N\xleftarrow{[N]}\ca{A}^{\tg,\ca{O}_F}_N\xrightarrow{\ion}\ca{A}^{\ca{O}_F},$$
which we still denote by $\iota^{\tg}$.
\end{cor}
Since, for any positive integer $M$ such that $N|M$ and $(M,p)=1$, $f^{\tg}_M=\frac{M}{N}f^{\tg}_N$. We then have $\ca{A}_N^{\tilde{\gamma}}\iso\ca{A}_M^{\tilde{\gamma}}/\ca{A}_M^{\tilde{\gamma}}[\frac{M}{N}]\iso\ca{A}_M^{\tilde{\gamma}}.$ So it makes sense to define $\ca{A}^{\tilde{\gamma}}:=\ca{A}^{\tilde{\gamma}}_N$, where $N$ is any prime-to-$p$ positive integer such that $Nz^{\tilde{\gamma}}_\sigma\in \Endo_S(\ca{A}^{\ca{O}_F})$. \par
\begin{rk}We may say $\ca{A}^{\tg}$ is the twisting of $\ca{A}$ by $\tg$.
By Lemma \ref{another-gamma}, for another choice $\tilde{\gamma}^\prime\in\ca{P}(\ca{O}_{F,(p)})$, $\ca{A}^{\tilde{\gamma}^\prime}$ and $\ca{A}^{\tilde{\gamma}}$ lie in the same prime-to-$p$ isogeny class, but might not be isomorphic. We will see later that this ambiguity can be eliminated by adding additional structures; see Proposition \ref{twt-sum}.
\end{rk}
\subsection{Twisting degenerating families}\label{subsec-twt-deg}
In this {subsection}, we construct twisting of degenerating families. As a preparation, materials in \S\ref{z-equiv}-\S\ref{twt-pol-revised} are again from \cite[Sec. 3.1 and 3.2]{Kis10}, \cite[Sec. 4.2]{Kis17} and \cite[Sec. 4.4 and 4.5]{KP15}.
\subsubsection{}\label{z-equiv}
Let $(G_0,X_0)$, $(G^\ddag,X^\ddag)$, $K_0$, $K^\ddag$ and $V_\bb{Z}$ be defined as in \S\ref{subsec-siegel}. Fix a place $v|p$ and denote $\ca{O}_{(v)}:=\ca{O}_{E_0,(v)}$. Recall that $\ca{S}_{K_0}:=\ca{S}_{K_0}(G_0,X_0)$ is the normalization of $\ca{S}_{K^\ddag}$ in $\sh_{K_0}:=\sh_{K_0}(G_0,X_0)$ via $\iota:\sh_{K_0}\to \sh_{K^\ddag}\to\ca{S}_{K^\ddag}$, which is a normal integral model of $\sh_{K_0}$; we will abusively denote by $\iota$ the finite morphism $\ca{S}_{K_0}\to\ca{S}_{K^\ddag}$ induced by taking the relative normalization. Recall that $\p=\emptyset$ or $\{p\}$.\par
Denote by $G_{0,\zbkp}$ the schematic closure of $G_0$ in $G_{\zbkp}^\ddag:=\stb_{G^\ddag}(V_{\zbkp})$. Let $Z_0$ be the center of $G_0$. Let $Z_{0,\zbkp}$ be the schematic closure of $Z_0$ in $G_{0,\zbkp}$. Then $Z_{0,\zbkp}$ is a flat affine group scheme of finite type over $\zbkp$. When $\p=\emptyset$, $G_{0,\zbkpp}=G_0$, $Z_{0,\zbkpp}=Z_0$, etc.\par
Let $\ca{A}$ be the pullback to $\ca{S}_{K_0}$ of the universal abelian scheme over $\ca{S}_{K^\ddag}$. By \cite[Lem. 4.5.2]{KP15}, there is a natural embedding 
$$\z:Z_{0,\zbkp}\hookrightarrow \aut_{\zbkp}(\ca{A});$$
on the other hand, there is a natural embedding $Z_{0,\zbkp}\hookrightarrow G_{0,\zbkp}\hookrightarrow \mrm{GL}(V_{\zbkp})$ by the construction. Then there are natural $Z_{0,\zbkpp}$-actions on both $V^\p\ca{A}$ and $V_{\App}$ induced by the two embeddings above. \par
We denote by $(\ca{A},\lambda,[\varepsilon_{\zhpp}]_{K^{\ddag,\p}})$ the pullback to $\ca{S}_{K_0}$ (resp. $\sh_{K_0}$) of the universal family over $\ca{S}_{K^\ddag}$ (resp. $\sh_{K^\ddag}$) obtained from the moduli problem $\mbf{M}^{\mrm{iso}}_{(V_\bb{Z},\psi_\bb{Z}),K^{\ddag,\p}}$ when $\p=\{p\}$ (resp. $\emptyset$). Let $[\varepsilon_{\zhpp}]_{K^\p_0}:=\varprojlim_{K_0^{\p}\sbst K^{\ddag,\p}}[\varepsilon_{\zhpp}]_{K^{\ddag,\p}}$.
\begin{lem}[{\cite[Lem. 4.2.4]{Kis17}}]\label{z-equi-lem}
With the conventions above, fix any representative $\varepsilon_{\zhp}$ in $[\varepsilon_{\zhp}]_{K^p_0}$.
Let $R$ be a $\zbkp$-algebra contained in $\bb{C}$. Then the isomorphism $$\varepsilon_{\zhp}:V_{\Ap}\otimes_{\zbkp} R\lra V^p\ca{A}\otimes_{\zbkp} R$$ 
is $Z_{0,\zbkp}(R)$-equivariant.
\end{lem}
\subsubsection{}\label{twt-level}
We equip $V_{\zbkp}$ with the action of $Z_{0,\zbkp}$ induced by the embedding $Z_{0,\zbkp}\hookrightarrow G_{0,\zbkp}\hookrightarrow \mrm{GL}(V_{\zbkp})$. \par
We denote by $\V_{\zbkp}$ the $\zbkp$-module $V_{\zbkp}$ equipped with the trivial $Z_{0,\zbkp}$-action. Let $\V_\bb{Z}:=\V_\zbkp\cap V_\bb{Z}$.
We denote by $\V:=\V_{\zbkp}\otimes_{\zbkp}\bb{Q}$ (resp. $V:=V_{\zbkp}\otimes_{\zbkp}\bb{Q}$) the vector space $V$ with $Z_0$-action induced by the $Z_{0,\zbkp}$-action on $\V_{\zbkp}$ (resp. $V_{\zbkp}$).\par
Let $G^{\ad}_{0,\zbkp}:=G_{0,\zbkp}/Z_{0,\zbkp}$. Let $G_{0,\zbkp}^\ad(\zbkp)_1$ be the subgroup of $G_{0,\zbkp}^\ad(\bb{Q})=G^\ad_0(\bb{Q})$ that stabilizes $X_0\hookrightarrow X_0^\ad$, that is, $G_{0,\zbkp}^\ad(\zbkp)_1:=G_{0,\zbkp}^\ad(\zbkp)\cap G^\ad(\bb{R})_1$. Let $\gamma\in G_{0,\zbkp}^{\ad}(\zbkp)_1$ and let $\ca{P}$ be the fiber of $G_{0,\zbkp}\to G^{\ad}_{0,\zbkp}$ at $\gamma$. 
As in \S\ref{twt-abelian}, there is a finite Galois extension $F$ over $\bb{Q}$ such that $\ca{P}(\ca{O}_{F,(p)})$ is non-empty.\par 
Let $\tg$ be any element in $\ca{P}(\ca{O}_{F,(p)})$ lifting $\gamma$. Then the $Z_{0,\zbkp}$-action on $\V_{\zbkp}$ (resp. $V_{\zbkp}$) specializes to a $\ca{O}_{F,(p)}$-semilinear $\gal(F/\bb{Q})$-action on $\V_{\zbkp}\otimes_{\zbkp}\ca{O}_{F,(p)}$ (resp. $V_{\zbkp}\otimes_{\zbkp}\ca{O}_{F,(p)}$).\par
With the conventions above, the isomorphism between $\ca{O}_{F,(p)}$-modules $(\tg^{-1}\cdot):\V_{\zbkp}\otimes_{\zbkp} \ca{O}_{F,(p)}\to V_{\zbkp}\otimes\ca{O}_{F,(p)}$ formed by left-multiplying $\tg^{-1}$ on the underlying $\ca{O}_{F,(p)}$-modules is $\gal(F/\bb{Q})$-equivariant: Indeed, $z_\sigma^{\tg}(1\otimes\sigma)\tg^{-1}=(z_\sigma^{\tg} z_\sigma^{\tg^{-1}}\tg^{-1})(1\otimes\sigma)=\tg^{-1}(1\otimes \sigma)$ since $\tg{\tg^{-1}}=1\in G_{0,\zbkp}(\ca{O}_{F,(p)})$. 
Together with Lemma \ref{z-equi-lem}, we obtain that
\begin{prop}\label{prop-comp-gal-equi}The composition
\begin{equation}\label{comp-gal-equi-eq}\begin{tikzcd}
{\V_{\Ap}\otimes_{\zbkp}\ca{O}_{F,(p)}}\arrow[r,"{({\tg}^{-1}\cdot)}","\sim"']&
{V_{\Ap}\otimes_{\zbkp}\ca{O}_{F,(p)}}\arrow[r,"\varepsilon_{\zhp}\otimes 1","\sim"']& V^p\ca{A}^{\ca{O}_F}
\end{tikzcd}\end{equation}
is $\gal(F/\bb{Q})$-equivariant.\end{prop}
Then $\lcj{K^p_0}{\gamma}:=\gamma K^p_0\gamma^{-1}$ (see the convention list after the introduction) acts on $(\varepsilon_{\zhp}\otimes 1)\circ(\tg^{-1}\cdot)$ by right-composition, i.e., $\gamma g\gamma^{-1}\in \lcj{K_0^p}{\gamma}$ sends $$(\varepsilon_{\zhp}\otimes 1)\circ(\tg^{-1}\cdot)\mapsto (\varepsilon_{\zhp}\otimes 1)\circ(\tg^{-1}\cdot)\circ \gamma g\gamma^{-1}.$$
Then
$(\varepsilon_{\zhp}\otimes 1)\circ{\tg}^{-1}(\V_{\zhp}\otimes F)$ is $\gal(F/\bb{Q})$-stable, and the $K_0^p$-orbit of it is $\pi_1(S,\overline{s})$-invariant when specializing to any geometric point $\overline{s}$ of $S$. \par
Assume now that $\lcj{K_0^p}{\gamma}$ stabilizes $\V_{\zhp}$. Consider the image of $\V_{\zhp}\otimes \ca{O}_F\sbst\V_{\Ap}\otimes_{\zbkp}\ca{O}_{F,(p)}$ in $V^p\ca{A}^{\ca{O}_F}$ via $(\varepsilon_{\zhp}\otimes 1)\circ{\tg}^{-1}$.
Then there is a $\zbkpt$-isogeny $$f:\ca{A}^{\ca{O}_F}\lra\ca{A}'$$ over $S$ which induces an isomorphism $V^pf:V^p\ca{A}^{\ca{O}_F}\xrightarrow{\sim}V^p\ca{A}^\prime$ such that $V^pf\circ(\varepsilon_{\zhp}\otimes 1)\circ{\tg}^{-1}(\V_{\zhp}\otimes \ca{O}_F)=T^p\ca{A}'$; the sheaf $[\ca{A}^\prime]$ is equipped with a $\gal(F/\bb{Q})$-action induced by $f$.\par
By Proposition \ref{prop-comp-gal-equi}, we can take the $\gal(F/\bb{Q})$-invariant part so that there is an isomorphism $$\varepsilon_{\zhp}^{\tg}:\V_{\zhp}\lra (T^p\ca{A}')^{\gal(F/\bb{Q})}.$$ 
The open compact subgroup $V^pf\circ (\varepsilon_{\zhp,\overline{s}}\otimes 1)\circ {\tg}^{-1}(\V_{\zhp}\otimes 1)$ of $V^p\ca{A}^{\ca{O}_F}_{\overline{s}}$ is $\gal(F/\bb{Q})$-invariant and $\pi_1(S,\overline{s})$-invariant (by the assumption on $\lcj{K_0^p}{\gamma}$), so it is contained in $V^p\ca{A}_{\overline{s}}^{\tg}$, where 
$\ca{A}^{\tg}$ is the twisted abelian scheme constructed in \S\ref{twt-ab-sch}. Then, from the construction above, there is a $\zbkpt$-isogeny $f^{\tg}:\ca{A}^{\tg}\to\ca{A}^{\tg,\prime}$ such that $T^p\ca{A}^{\tg,\prime}=(T^p\ca{A}')^{\gal(F/\bb{Q})}$.
Hence, we write $\varepsilon^{\tg}_{\zhp}$ as 
$$\varepsilon^{\tg}_{\zhp}: \V_{\zhp}\lra T^p\ca{A}^{\tg,\prime}.$$
\begin{lem}\label{lem-indep}
The isomorphism $\varepsilon^{\tg}_{\zhp}$ is independent of the choice of $\tg\in\ca{P}(\ca{O}_{F,(p)})$ lifting $\gamma$. 
\end{lem}
\begin{proof}
Let $\wdtd{\gamma}^\prime\in\ca{P}(\ca{O}_{F,(p)})$ such that ${\tg}^\prime=z\tg$ for some $z\in Z_{0,\zbkp}(\ca{O}_{F,(p)})$. There is a commutative diagram
\begin{equation}
    \begin{tikzcd}
    \V_{\Ap}\otimes_{\zbkp}\ca{O}_{F,(p)}\arrow[r,"{({\tg}^{-1}\cdot)}"]\arrow[dr,"{({\tg}^{\prime,-1}\cdot)}"]& V_{\Ap}\otimes_{\zbkp}\ca{O}_{F,(p)}\arrow[d,"z^{-1}"]\arrow[rr,"{\varepsilon_{\zhp}\otimes 1}"]&& V^p\ca{A}^{\ca{O}_F}\arrow[d,"z^{-1}"]\arrow[r,"V^pf"]&V^p\ca{A}'\\
    & V_{\Ap}\otimes_{\zbkp}\ca{O}_{F,(p)}\arrow[rr,"{\varepsilon_{\zhp}\otimes 1}"]&&V^p\ca{A}^{\ca{O}_F}\arrow[ur,"{V^pf'}"]& 
    \end{tikzcd}
\end{equation}
The composition $(\varepsilon_{\zhp}\otimes 1)\circ({\tg}^{-1}\cdot)$ (resp. $(\varepsilon_{\zhp}\otimes 1)\circ({\tg}^{\prime,-1}\cdot)$) induces a $\zbkpt$-isogeny $f:\ca{A}^{\ca{O}_F}\to\ca{A}'$ (resp. $f':\ca{A}^{\ca{O}_F}\to\ca{A}^\prime$), such that $f=f^\prime\circ z^{-1}$, since the images of $\V_{\zhp}\otimes\ca{O}_F$ in $V^p\ca{A}^{\ca{O}_F}$ via $z^{-1}\circ(\varepsilon_{\zhp}\otimes 1)\circ({\tg}^{-1}\cdot)$ and $(\varepsilon_{\zhp}\otimes 1)\circ({\tg}^{\prime,-1}\cdot)$ coincide with each other. Then we have a commutative diagram above by \cite[Cor. 1.3.5.4]{Lan13} and Lemma \ref{z-equi-lem}, and we have the desired claim by taking the $\gal(F/\bb{Q})$-invariant part.
\end{proof}
Denote $\varepsilon_{\zhp}^{\tg}$ as $\varepsilon_{\zhp}^{\gamma}$. We then call the $\lcj{K_0^p}{\gamma}$-orbit $[\varepsilon_{\zhp}^\gamma]_{\lcj{K_0^p}{\gamma}}$ a \textbf{twisted (integral) level structure}.
\subsubsection{}\label{twt-pol-revised}
Recall that $\nu:G^\ddag\to\bb{G}_m$ denotes the similitude character. Denote by $\psi^\V_\bb{Z}:\V_\bb{Z}\times \V_\bb{Z}\to \bb{Z}(1)$ the pairing $\psi_\bb{Z}$ viewed as a pairing for $\V_\bb{Z}$. We can also obtain a pairing 
$$\psi_\bb{Z}\otimes 1: (V_\bb{Z}\otimes\ca{O}_{F,(p)})\times (V_\bb{Z}\otimes\ca{O}_{F,(p)})\to \ca{O}_{F,(p)}(\nu)$$
by the base change of $\psi_\bb{Z}$ to $\ca{O}_{F,(p)}$. 
Similarly, define $$\psi^{\V}_\bb{Z}\otimes 1: (\V_\bb{Z}\otimes\ca{O}_{F,(p)})\times (\V_\bb{Z}\otimes\ca{O}_{F,(p)})\to \ca{O}_{F,(p)}(\nu).$$
Note that $\gal(F/\bb{Q})$ acts on $\ca{O}_{F,(p)}(\nu)$ via $\sigma\mapsto \nu(z^{\tg}_\sigma)\sigma$ for $\psi_{\bb{Z}}\otimes 1$ and via $\sigma$ for $\psi^{\V}_\bb{Z}\otimes 1$.  
Hence, there is a $\gal(F/\bb{Q})$-equivariant map 
$$\tg^{-1}: (\V\otimes{\oo}\otimes \zbkp,\psi^\V_\bb{Z}\otimes 1)\to (V\otimes{\oo}\otimes \zbkp,\nu(\tg)\circ\psi_\bb{Z}\otimes 1).$$

On the other hand, we do the same procedures for $\ca{A}^\vee$ as above. First, we obtain an abelian scheme $\ca{A}^{\vee,\tg}$, which is defined by the $Z_{0,\zbkp}$-action given by the canonical isomorphism (\ref{eq-aut-dual-iso}).\par
We then consider the diagram below.
\begin{equation}\label{diag-twt-pol-ab}
    \begin{tikzcd}
    {\V_{\Ap}\otimes_{\zbkp} \ca{O}_{F,(p)}}\arrow[r,"{(\tg^{-1}\cdot)}"]\arrow[d,"\psi^{\V}\otimes 1"]&V_{\Ap}\otimes_{\zbkp}\ca{O}_{F,(p)}\arrow[r,"\varepsilon_{\zhp}\otimes 1"]\arrow[d,"\psi\otimes 1"]& V^p\ca{A}^{\oo}\arrow[r,"V^pf"]\arrow[d,"{V^p\lambda\otimes 1}"]& V^p\ca{A}'\arrow[d,dashed,"{V^p\lambda^{\V,\oo}}"]\\
{\V_{\Ap}^\vee\otimes_{\zbkp} \ca{O}_{F,(p)}(\nu)}\arrow[r]&{V_{\Ap}^\vee\otimes_{\zbkp}\ca{O}_{F,(p)}(\nu)}\arrow[r]& V^p\ca{A}^{\vee}\otimes_{\zbkp}\ca{O}_{F,(p)}(\nu)\arrow[r,dashed,"{V^pf'}"]& V^p\ca{A}^*(\nu)
    \end{tikzcd}
\end{equation}
The first two vertical maps are canonically induced by the perfect pairings marked on the arrows. 
Let us explain how to produce the arrows in the lower-right corner. 
In fact, let $\ca{U}^\vee$ be the image of $\V_{\zhp}^\vee\otimes \oo$ in $V^p\ca{A}^\vee\otimes_{\zbkp}\ca{O}_{F,(p)}(\nu)$ via the composition $(V^p\lambda\otimes 1)\circ(\varepsilon_{\zhp}\otimes 1)\circ (\tg^{-1}\cdot)\circ (\psi^{\V}\otimes 1)^{-1}$. Then, by \cite[Cor. 1.3.5.4]{Lan13} (see also \cite[Lem. 4.22]{Mum91}), $\ca{U}^\vee$ determines a $\zbkpt$-isogeny $f':(\ca{A}^\vee)^{\oo}\to \ca{A}^*$ and a prime-to-$p$ isogeny $\lambda^{\V,\oo}:\ca{A}'\to \ca{A}^*$ so that the right square of the diagram above commutes.\par
By Lemma \ref{z-equi-lem}, $\varepsilon_{\zhp}\otimes 1$ is $Z_{0,\zbkp}(\ca{O}_{F,(p)})$-equivariant. Combining this with the facts that $\varepsilon_{\zhp}\otimes 1$ sends $\psi_{\zhp}$ to a $\wat{\bb{Z}}^{p,\times}$-multiple of $\e^\lambda$, and that $\psi_{\zhp}$ is $Z_{0,\zbkp}(\ca{O}_{F,(p)})$-equivariant, we obtain that $V^p\lambda\otimes 1$ is $Z_{0,\zbkp}(\ca{O}_{F,(p)})$-equivariant and therefore is $\gal(F/\bb{Q})$-equivariant. Hence, this is also true for $V^pf'$ and $V^p\lambda^{\V,\oo}$.\par
Moreover, since $$V^p(\ca{A}^\vee)^{\tg}= (\lambda^{\V,\oo}\circ V^pf\circ(\varepsilon_{\zhp}\otimes 1)\circ({\tg}^{-1}\cdot)\circ(\psi^{\V}\otimes 1)^{-1})(\V^\vee_{\Ap}\otimes_{\zbkp}\ca{O}_{F,(p)}(\nu))^{\gal(F/\bb{Q})},$$
the image of $\V^\vee_{\zhp}$ in $V^p(\ca{A}^\vee)^{\tg}$ induces a $\zbkpt$-isogeny $(\ca{A}^\vee)^{\tg}\to \ca{A}''$ such that $T^p\ca{A}''=(T^p\ca{A}^*)^{\gal(F/\bb{Q})}$. There is also a prime-to-$p$ isogeny $\lambda^{\tg,\prime}:\ca{A}^{\tg,\prime}\to \ca{A}''$ which is the $\gal(F/\bb{Q})$-invariant of $\lambda^{\V,\oo}:\ca{A}'\to \ca{A}^*$.\par
\begin{lem}\label{lem-twt-pol-new}
    With the conventions above, the abelian scheme $\ca{A}''$ is isomorphic to the dual of $\ca{A}^{\tg,\prime}$, and the isogeny $$\lambda^{\tg,\prime}:\ca{A}^{\tg,\prime}\to\ca{A}''$$ is a polarization.
\end{lem}
\begin{proof}
Choose any positive definite perfect pairing $\e_0:\oo\times \oo\to \bb{Z}$ which canonically induces an isomorphism of $\bb{Z}$-modules $i_F:\oo\xrightarrow{\sim}\mrm{Hom}_\bb{Z}(\oo,\bb{Z})$. Since $\lambda(i_F): \ca{A}^{\oo}\xrightarrow{\lambda}\ca{A}^{\vee}\otimes \oo\xrightarrow{i_F}\ca{A}^\vee\otimes\ca{O}_F^\vee$ is a polarization (see \cite{AK18}), so is $\lambda^{\V,\oo}(i_F):\ca{A}'\xrightarrow{\lambda^{\V,\oo}}\ca{A}^*\xrightarrow{i_F}\ca{A}^*\otimes_{\oo,i_F}\ca{O}_F^\vee$. Hence, we have obtained a polarization $$\ca{A}^{\tg,\prime}\otimes \oo\xrightarrow{\lambda^{\tg,\prime}\otimes 1}\ca{A}''\otimes \oo\xrightarrow{\sim}\ca{A}''\otimes\ca{O}_F^\vee.$$

Taking the $\gal(F/\bb{Q})$-invariant part, we obtain a polarization $\ca{A}^{\tg,\prime}\otimes \bb{Z}\xrightarrow{\lambda^{\tg,\prime}\otimes 1}\ca{A}''\otimes \bb{Z}\xrightarrow[\sim]{i_F}\ca{A}''\otimes \mrm{Hom}_\bb{Z}(\bb{Z},\bb{Z})$ (where the second isomorphism is induced by the unique positive definite perfect duality $\bb{Z}\xrightarrow{\sim} \mrm{Hom}_\bb{Z}(\bb{Z},\bb{Z})$ and does not depend on $i_F$), which implies that $\lambda^{\tg,\prime}$ is also a polarization.
\end{proof}
Let $\ca{A}^\gamma:=\ca{A}^{\tg,\prime}$ and $\lambda^\gamma:=\lambda^{\tg,\prime}$.
\begin{prop}\label{twt-sum}
Let $\gamma\in G_{0,\zbkp}^{\ad}(\zbkp)_1$. For the tuple $(\ca{A},\lambda,[\varepsilon_{\zhp}]_{K^p_0})$, we can associate with it a tuple $(\ca{A}^{\gamma},\lambda^\gamma,[\varepsilon_{\zhp}^\gamma]_{\lcj{K^p_0}{\gamma}})$, which is independent of the choice of $\tg$ lifting $\gamma$. For any open compact subgroup $K^{\ddag,p,\prime}$ of $G^\ddag(\Ap)$ containing $\lcj{K^p}{\gamma}$ that stabilizes $\V_{\zhp}$, the tuple $(\ca{A}^\gamma,\lambda^\gamma,[\varepsilon_{\zhp}^\gamma]_{K^{\ddag,p,\prime}})$ is an object in $\mbf{M}_{(\V_{\bb{Z}},\psi^\V_\bb{Z}),K^{\ddag,p,\prime}}^{\mrm{iso}}$. \par
This association defines a morphism $\gamma^{-1}:\ca{S}_K\to\ca{S}_{\lcj{K}{\gamma}}$, whose generic fiber $\gamma^{-1}:\sh_{K_0}\to\sh_{\lcj{K_0}{\gamma}}$ is determined by the map $[(x,g)]\mapsto [(\gamma(x),\gamma g\gamma^{-1})]$ over complex points.\par
\end{prop}
\begin{proof}
The first paragraph is from the main construction in \S\ref{twt-level} and the paragraphs above. The tuple $(\ca{A}^{\gamma},\lambda^\gamma,[\varepsilon^\gamma_{\zhp}]_{K^{\ddag,p,\prime}})$ viewed as an object in $\mbf{M}^{\mrm{isog}}_{K^{\ddag,p,\prime}}$ is exactly the twisting construction in \cite[Sec. 4.5]{KP15}. Combining it with \cite[Lem. 3.2.6]{Kis10}, we obtain the second paragraph as in \cite[Lem. 4.5.7]{KP15} by normalization. 
\end{proof}

\subsubsection{}\label{deg-fam-notion}Let $\Sigma_0$ be any admissible polyhedral cone decomposition for $\ca{S}_{K_0}$ that is compatible with an admissible polyhedral cone decomposition $\Sigma^\ddag$ for $\ca{S}_{K^\ddag}$. Let $\ca{S}^{\Sigma_0}_{K_0}$ (resp. $\ca{S}^{\Sigma^\ddag}_{K^\ddag}$) be the toroidal compactification of $\ca{S}_{K_0}$ (resp. $\ca{S}_{K^\ddag}$) associated with $\Sigma_0$ (resp. $\Sigma^\ddag$). Let $S:=\ca{S}^{\Sigma_0}_{K_0}$ and $U:=\ca{S}_{K_0}$.\par
For any $K^{p}_0$ contained in $K^{\ddag,p}$ as in \S\ref{subsec-siegel}, we denote by $(\G,\lambda_S,[\varepsilon_{\zhp}]_{K^{\ddag,p}})$ the pullback of the degenerating family over $\ca{S}_{K^\ddag}^{\Sigma^\ddag}$ to $S$. Then $(\G,\lambda_S,[\varepsilon_{\zhp}]_{K^{\ddag,p}})$ is a degenerating family in $\mbf{DEG}_{(V_\bb{Z},\psi_\bb{Z}),K^{\ddag,p}}(S,U)$ extending $(\ca{A},\lambda,[\varepsilon_{\zhp}]_{K^{\ddag,p}})$. From \S\ref{z-equiv}, we obtain a tuple $(\G,\lambda_S,[\varepsilon_{\zhp}]_{K^{p}_0})$.\par
Denote $[\G]:=\G\otimes\zbkp$, which represents $\G$ up to prime-to-$p$ isogenies.
\begin{lem}\label{ext-end-deg}There is a natural isomorphism $\aut_{\zbkp}([\G])\iso \aut_{\zbkp}([\ca{A}])$ induced by restricting automorphisms of $[\G]$ from $S$ to $U$. In particular, $\z:Z_{0,\zbkp}\to \aut_\zbkp([\ca{A}])$ extends uniquely to a homomorphism 
    $$\z^{\G}:Z_{0,\zbkp}\lra \aut_{\zbkp}([\G]).$$
\end{lem}
\begin{proof}
The second statement follows from the first. By \cite[I. Prop. 2.7]{FC90}, any homomorphism $f_U\in \Endo_U(\ca{A})$ extends uniquely to a homomorphism $f_S\in \Endo_S(\G)$. Moreover, $f_S$ is a prime-to-$p$ isogeny if and only if $f_U$ is so. 
(Indeed, if $f_U\circ g_U=g_U\circ f_U=[M]$ for some $g_U\in \Endo_U(\ca{A})$ and $(M,p)$=1, then $g_U$ extends uniquely to $g_S\in \Endo_S(\G)$ such that $f_S\circ g_S=g_S\circ f_S=[M]$, so the claim is true since $[M]$ is a prime-to-$p$ isogeny for $\G$.) Hence, $\Endo_S(\G)\iso \Endo_U(\ca{A})$, which induces an isomorphism $\Endo_S([\G])\iso \Endo_U([\ca{A}])$. So $(\Endo_S([\G])\otimes_{\zbkp}R)^\times\iso(\Endo_U([\ca{A}])\otimes_{\zbkp} R)^\times$, as desired.
\end{proof}
For any $\gamma\in G_{0,\zbkp}^\ad(\zbkp)_1$ as in \S\ref{twt-level}, we describe the twisting of $(\G,\lambda,[\varepsilon_{\zhp}])$ by $\gamma$. Let $F$ and $\tg$ be as in \S\ref{twt-ab-sch}. Denote by $f^{\tg}_N$ the homomorphism defining $\ca{A}^{\tg}$ as in \S\ref{twt-ab-sch}. Then Lemma \ref{ext-end-deg} implies that $f^{\tg}_N$ extends uniquely to a homomorphism $$f^{\tg}_{N,S}:\G^{\ca{O}_F}\lra (\G^{\ca{O}_F})^{\oplus |\gal (F/\bb{Q})|-1}$$ over $S$, which can still be written as $f^{\tg}_{N,S}=\bigoplus\limits_{\sigma\in\mrm{Gal}(F/\bb{Q})\text{ and }\sigma\neq\mrm{id}}((Nz_\sigma^{\tg})\sigma-N\cdot\mrm{id})$. Denote by $(\ker f^{\tg}_{N,S})^\circ$ the identity component of the kernel of $f^{\tg}_{N,S}$.
\begin{lem}\label{lem-semiabelian-twist}
$(\ker f^{\tg}_{N,S})^\circ$ is a semi-abelian scheme over $S$.
\end{lem}
\begin{proof}
Arguing as in the first paragraph of the proof of Proposition \ref{twt-kp}, we only need to show that the fiber $(\ker f^{\tg}_{N,S})^\circ_t$ at $t$ is an extension of an abelian variety by a torus for $t\in S\bss U$.\par 
Write $\G_{t}$ as an extension $0\to T_t\to \G_{t}\to \ca{A}_t\to 0.$ Let $f$ be any \emph{surjective} homomorphism between $\G^{\ca{O}_F}_t$. Then $f$ maps $T_t^{\ca{O}_F}$ to itself and also induces an endomorphism of $\ca{A}_t^{\ca{O}_F}$ by taking the quotient of $\G_{t}^{\oo}$ by $T_t^{\oo}$: To see this, note that $\G_{t}^{\ca{O}_F}/f(T_t^{\ca{O}_F})$ only has abelian part since it is the image of the proper scheme $\ca{A}_t^{\ca{O}_F}$. \par 
By applying the last paragraph to $Nz_\sigma^{\tg}\sigma$ and $N\cdot\mrm{id}$, we see that $f^{\tg}_{N,t}$ maps the torus part of $\G_{t}^{\ca{O}_F}$ to the torus part of its target, and induces a homomorphism $f^{\tg,ab}_{N,t}$ from the abelian part of $\G_{t}^{\ca{O}_F}$ to the abelian part of its target. By Proposition \ref{twt-kp}, $(\ker f^{\tg,ab}_{N,t})^\circ$ is an abelian variety, and $(\ker f^{\tg}_{N,t}|_{T_t^{\ca{O}_F}})^\circ$ is a torus. 
Finally, we note that there is an exact sequence $$0\lra \ker f^{\tg}_{N,t}|_{T_t^{\ca{O}_F}}\lra \ker f^{\tg}_{N,t}\lra \ker f^{\tg,ab}_{N,t}\lra \Lambda$$ with $\Lambda$ some finite algebraic group. Indeed, $\ker f^{\tg,ab}_{N,t}$ is proper and $\coker f^{\tg}_{N,t}|_{T_t^{\ca{O}_F}}$ is an affine group, so it follows by snake lemma. So the homomorphism $(\ker f^{\tg}_{N,t})^\circ\to(\ker f^{\tg,ab}_{N,t})^\circ$ must be surjective: $(\ker f^{\tg,ab}_{N,t})^\circ$ is covered by its fiber in $\ker f^{\tg}_{N,t}$, which is a union of connected components of $\ker f^{\tg}_{N,t}$. Then the images of any two connected components of the fiber are equal. \par 
After taking the identity component of $\ker f^{\tg}_{N,t}$, we have an exact sequence
$$1\lra \wdtd{T}\lra (\ker f^{\tg}_{N,t})^\circ\lra \ca{A}_t^{\tg}\lra 1,$$
where the first term $\wdtd{T}$ is an algebraic group of multiplicative type. Then we see that the quotient $(\ker f^{\tg}_{N,t})^\circ/\wdtd{T}^\circ$ is a connected commutative group variety which is an {\'e}tale cover of $\ca{A}^{\tg}_t$, and therefore is also an abelian variety. Then we have proved the desired result.
\end{proof}
We denote $\G^{\tg}:=(\ker f^{\tg}_{N,S})^\circ$.
\begin{prop}\label{prop-ext-semiabelian-twist}
With the conventions and constructions above, $(\ca{A}^\gamma,\lambda^{\gamma})$ extends uniquely to a pair $(\G^\gamma,\lambda^\gamma_S)$ with $\G^\gamma$ a semi-abelian scheme over $S$ and $\lambda^\gamma_S$ an isogeny $\lambda^\gamma_S:\G^\gamma\to \G^{\gamma,\vee}$.
\end{prop}
\begin{proof}
By \cite[I. Prop. 2.7]{FC90}, $\lambda^{\tg}$ extends uniquely to an isogeny $\lambda_S^\gamma$ if $\ca{A}^{\gamma}$ and $\ca{A}^{\gamma,\vee}$ extend to $\G^\gamma$ and $\G^{\gamma,\vee}$; the uniqueness of such $\G^\gamma$ and $\G^{\gamma,\vee}$ also follows from \emph{loc. cit}.\par
There is a natural $\zbkpt$-isogeny $f^{\tg}:(\ca{A}^{\tg},\lambda^{\tg})\to(\ca{A}^{\gamma},\lambda^\gamma)$ by Proposition \ref{twt-sum}. Choose any prime-to-$p$ integer $M>0$ such that $[M]\circ f^{\tg}$ is a prime-to-$p$ isogeny. Denote by $d$ the degree of $[M]\circ f^{\tg}$. Then $\ker [M]\circ f^{\tg}\sbst \ca{A}^{\tg}[d]$.
Denote by $\overline{K}$ the closure of $\ker [M]\circ f^{\tg}$ in $\G^{\tg}[d]$. By the following general fact Lemma \ref{primep-et}, $\overline{K}$ is quasi-finite and {\'e}tale. Then the quotient $\G^\gamma:=\G^{\tg}/\overline{K}$ is a semi-abelian scheme extending $\ca{A}^\gamma$ by \cite[IV, 7.1.2]{MB85} (see also \cite[Lem. 3.4.3.1]{Lan13}).
\end{proof}
\begin{lem}\label{primep-et}
Let $X$ be a quasi-finite separated commutative group scheme over a locally Noetherian normal scheme $S$ over $\zbkp$. Suppose that the order of $X$ over every fiber is prime-to-$p$, and that there is a dense open subscheme $U$ of $S$ over which $X$ is quasi-finite and {\'e}tale. Then the schematic closure $\overline{X}$ of $X_U$ in $X$ is {\'e}tale over $S$.
\end{lem}
\begin{proof}
We can show this by passing to the strict Henselization $\ca{O}_t^{sh}:=\wat{\ca{O}}_{S,t}^{sh}$ of the complete local ring $\wat{\ca{O}}_{S,t}$, for any $t\in S\bss U$. 
Moreover, we can replace $U$ with the {\'e}tale locus of $\overline{X}$ over $S$. If the lemma were not true, the complement $D:=S\bss U$ would be closed and \emph{nonempty}. 
We can further let $t$ be a minimal prime ideal in $D$. Denote the special fiber of $\ca{O}_t^{sh}$ by $\overline{t}$ and the complement of $\overline{t}$ in $\ca{O}_t^{sh}$ by $U_t$. To find the contradiction, we show the claim that  $\overline{X}_{\ca{O}^{sh}_{t}}$ is {\'e}tale over $\ca{O}_t^{sh}$.\par 
As $\overline{X}$ is quasi-finite and separated, it can be decomposed as a disjoint union of schemes $\overline{X}_{\ca{O}^{sh}_t}=\overline{X}^f\disju \overline{X}^{\emptyset}$, where $\overline{X}^f$ is a finite, open, and closed subgroup of $\overline{X}_{\ca{O}^{sh}_{t}}$ and $\overline{X}^{\emptyset}_{\overline{t}}=\emptyset$ (see, e.g., \cite[3.4.1]{Lan13}). Over $U_t$, $\overline{X}_{U_t}=X_{U_t}$ is {\'e}tale (and therefore quasi-finite {\'e}tale over $\ca{O}_t^{sh}$). 
So it suffices to show the claim for $\overline{X}^f$. Since the commutative group scheme $\overline{X}^f_{\overline{t}}$ is contained in $X_{\overline{t}}$, it is of prime-to-$p$ order, so each of its geometrically connected components must be trivial. Hence, $\overline{X}^f_{\overline{t}}\iso \pi_0(\overline{X}^f_{\overline{t}})$ is finite {\'e}tale over $\overline{t}$. By \cite[18.8.3, d)]{EGA4}, for any point $x_0\in \overline{X}^f_{\overline{t}}$, there is an affine open neighborhood $U_{x_0}$ of $x_0$ in $\overline{X}^f$ such that the restriction to $U_{x_0}$ of the structure morphism to $\spec\ca{O}_t^{sh}$ is a closed embedding. Hence, we can write $U_t$ as a disjoint union of open and closed subschemes $U_1$ and $U_2$ such that $U_1$ is the intersection of $U_t$ and the image of $U_{x_0}$ under the structure morphism. However, the closures of $U_1$ and $U_2$ in $\ca{O}^{sh}_{t}$ are two closed subschemes. Since $\ca{O}_{S,t}$ is a normal domain, so is $\ca{O}_t^{sh}$. So either $U_1$ or $U_2$ is empty. But the former case would imply that $x_0$ is not in the closure $\overline{X}_{\ca{O}^{sh}_t}$. So $U_2$ is empty, and therefore the morphism of $U_{x_0}$ to $\spec \ca{O}_t^{sh}$ is a surjective closed immersion. So the claim is true and the lemma is proved.
\end{proof}
Let $[\varepsilon^\gamma_{\zhp}]_{\lcj{K_0}{\gamma}}$ be the twisted level structure induced from $[\varepsilon_{\zhp}]_{K_0}$ by Proposition \ref{twt-sum}. Then the tuple $(\G^\gamma, \lambda^\gamma_S,[\varepsilon^\gamma_{\zhp}]_{\lcj{K_0}{\gamma}})$ is the desired \textbf{twisted degenerating family}.
\subsubsection{}  
Let us recall the construction of Hecke twists on $\ca{S}_{K_0}$ and $\ca{S}_{K_0}^{\Sigma_0}$; see \cite[6.4.3]{Lan13}.
Choose any $g\in G_0(\Ap)$ such that $K_0^{p,g}:=\rcj{(K_0^p)}{g}$ stabilizes some $\zhp$-lattice $V'_{\zhp}\sbst V_{\Ap}$. Recall that $\nu(g)$ can be uniquely decomposed as $\nu(g)=r(g)h(g)$ for some $r(g)\in \bb{Z}_{(p),>0}^\times$ and $h(g)\in \wat{\bb{Z}}^{p,\times}$.\par
Consider the sequence 
$$\begin{tikzcd}V_{\zhp}'\sbst V_{\Ap}\arrow[r,"g\cdot"] & V_{\Ap}\arrow[r,"\varepsilon_{\zhp}\otimes 1"]& V^p\ca{A}.\end{tikzcd}$$
Then the image $(\varepsilon_{\zhp}\otimes 1)\circ g(V_{\zhp}')$ of $V_{\zhp}$ determines a $\zbkpt$-isogeny $f:(\ca{A},\lambda)\lra (\ca{A}^g,\lambda'')$. 
Since the map
$$g:(V_{\Ap},\psi_{\Ap})\lra(V_{\Ap},\nu(g)^{-1}\psi_{\Ap})$$
sends $\psi_{\Ap}$ to $\nu(g)^{-1}\psi_{\Ap}$, we adjust $\lambda''$ to $\lambda^g:=r(g)^{-1}\lambda''$, then $f\circ(\varepsilon_{\zhp}\otimes 1)\circ (g\cdot)$ sends $\psi_{\zhp}$ to a $\wat{\bb{Z}}^{p,\times}$-multiple of $\mrm{e}^{\lambda^g}$.  
By Lemma \ref{primep-et} and \cite[I. Prop. 2.7]{FC90} and as the proof of Proposition \ref{prop-ext-semiabelian-twist}, $(\ca{A}^g,\lambda^g)$ extends uniquely to a pair $(\G^g,\lambda^g_S)$ of a semi-abelian scheme $\G^g$ and an isogeny $\lambda^g_S$ over $S$. The tuple $(\G^g, \lambda^g_S,[\varepsilon^g_{\zhp}]_{K^{p,g}_0})$ is the desired Hecke twist of the degenerating family $(\G,\lambda_S,[\varepsilon_{\zhp}]_{K^p_0})$.\par
Let $g\in G_0(\Ap)$ and let $\gamma\in G_0^\ad(\zbkp)_1$ be chosen as above. 
Following the sign convention in \cite[3.4.4]{Kis10} and \cite[4.5.6]{KP15}, we convert the left action of $G_0^\ad(\zbkp)_1$ to a right action: the right action of $(g,\gamma^{-1})\in G_0(\Ap)\times G_0^\ad(\zbkp)_1$ is given by $$(\G,\lambda_S,[\varepsilon_{\zhp}]_{K^p_0})\mapsto ((\G^g)^\gamma,(\lambda^g)^\gamma,[(\varepsilon^{g}_{\zhp})^\gamma]_{\lcj{(K^{p,g}_0)}{\gamma}}).$$
\subsection{Twisting $1$-motives}\label{subsec-twt-1-mot}
In this {subsection}, we describe the twisting action on $1$-motives with additional structures.
Note that this construction is similar to, but is not a direct consequence of, the last {subsection}.
\subsubsection{}
Let us first record Lemma \ref{lem-kisin}, which serves as the counterpart to \cite[Lem. 3.2.2]{Kis10} and \cite[Lem. 4.5.2]{KP15} for $1$-motives. 
Let $Z_{0,\zbkp}$ and $G_{0,\zbkp}$ be the groups defined in \S\ref{z-equiv}. Recall that $Z_{0}(\bb{C})\hookrightarrow G_0(\bb{C})\hookrightarrow \mrm{GL}(V_\bb{C})$ induces an action of  $Z_{0}(\bb{C})$ on $V_\bb{C}$. Fixing an $x\in D_{\Phi_0}\sbst D_{\Phi^\ddag}$, the integral mixed Hodge structure $\imxh(x)$ corresponds to a $1$-motive $\Q_x$.  The action of $Z_0(\bb{Q})$ on $\mxh(x)$ induces an embedding $$Z_0(\bb{Q})\hookrightarrow \aut_\bb{Q}\mxh(x)\xrightarrow[Thm. \ref{thm-deligne-iii}]{\sim}\aut_\bb{Q}\Q_x.$$ 
The action of $Z_0(\bb{C})$ on $\mxh(x)\otimes_\bb{Q}\bb{C}$ induces an embedding $Z_0(\bb{C})\hookrightarrow \aut_\bb{C}\Q_x$. \par
Let us fix a $\wp=[(x,\mbf{p})]$, where $\mbf{p}\in P_{\Phi_0}(\A)$. Since there is a $\bb{Q}^\times$-isogeny $f:\Q_x\to \Q_\wp$, $f$ induces an isomorphism $\mxh(x)\iso T_\bb{Z}\Q_x\otimes\bb{Q}\iso T_\bb{Z}\Q_\wp\otimes\bb{Q}$. So there is an embedding $Z_0(\bb{Q})\hookrightarrow\aut_\bb{Q}\mxh(x)\iso \aut_\bb{Q} T_\bb{Z}\Q_\wp\iso \aut_\bb{Q} \Q_\wp$; similarly, there is an embedding for $\bb{C}$-points, $Z_0(\bb{C})\hookrightarrow\aut_\bb{C}\Q_\wp$. Those embeddings do not depend on the choice of liftings $(x,\mbf{p})$ of $\wp$. Moreover, we have
\begin{lem}\label{lem-kisin}Let $\Phi_0$ and $\Phi^\ddag$ be the cusp label representatives defined as in \S\ref{subsec-mxhg-cont}.
    Let $\Q$ be the pullback of the tautological $1$-motive $\Q^\ddag$ from $\ca{S}_{K_{\Phi^\ddag}}$ to $\ca{S}_{K_{\Phi_0}}$. Then there is an embedding $\z_{\Phi_0}:Z_{0,\zbkp}\hookrightarrow \aut_\zbkp\Q$, such that the pullback to $\wp$ of the $\bb{C}$-points of $\z_{\Phi_0}$ is the embedding above.
\end{lem}
\begin{proof}The arguments are similar to those in \cite[Lem. 3.2.2]{Kis10} and \cite[Lem. 4.5.2]{KP15}.
For any $1$-motive $\Q$ over $\ca{S}_{K_{\Phi_0}}$, denote by $\Q_\eta$ the generic fiber of $\Q$ over $\sh_{K_{\Phi_0}}$.  \par 
First, let us show that there is an embedding $Z_0\hookrightarrow \aut_\bb{Q}\Q_\eta$. The idea is completely the same as \cite[Lem. 3.2.2]{Kis10}.
Let $\wp=[(x,\mbf{p})]\in \sh_{K_{\Phi_0}}(\bb{C})$. For any $(x,\mbf{p})$ lifting $\wp$, $x$ defines an embedding $e_x:G_0(\bb{C})\hookrightarrow\aut_\bb{C}\Q_x\iso \aut_\bb{C}\mxh(x)$. This embedding depends only on $x$. For any other $x^\prime\in D_{\Phi_0}$, $e_x$ and $e_{x^\prime}$ differ by a conjugation of an element in $P_{\Phi_0}(\bb{R})U_{\Phi_0}(\bb{C})$. Hence, the embedding of $Z_0(\bb{C})\hookrightarrow G_0(\bb{C})$ composed with $e_x$ does not depend on the choice of $x\in D_{\Phi_0}$. Thus, there is a well-defined embedding $Z(\bb{C})\hookrightarrow\aut_\bb{C}\Q_\eta(\bb{C})$, which induces an embedding $Z_0\hookrightarrow\aut_\bb{Q}\Q_\eta(\bb{C})$.\par
Let $E_0:=E(G_0,X_0)=E(P_{\Phi_0},D_{\Phi_0})$ (see \cite[Prop. 12.1]{Pin89}).  
Since $\Q_\eta[n](\bb{C})$ and $\aut\Q_\eta[n](\bb{C})$ are defined over $\sh_{K_{\Phi_0},\overline{E}_0}$, so are $T_p\Q_\eta(\bb{C})$ and $\aut T_p\Q_\eta(\bb{C})$. So the action of $\Aut(\bb{C}/E_0)$ on $T_p\Q_\eta(\bb{C})$ and $\aut T_p\Q_\eta(\bb{C})$ factors through $\Aut(\overline{E}_0/E_0)$. Also, since $\aut(\Q_{\wp}(\bb{C}))\hookrightarrow \aut (T_p\Q_{\wp}(\bb{C}))$ at any $\wp$ (see Lemma \ref{lem-tate-inj}), $\aut\Q_\eta(\bb{C})\iso\aut\Q_\eta(\overline{E}_0)$.\par
Let $(T,X_T)$ be a Shimura datum of a $\bb{Q}$-torus $T$ with an embedding $\mbf{sp}_T:(T,X_T)\to (P_{\Phi_0},D_{\Phi_0})$. Denote by $E_T$ the reflex field of $(T,X_T)$. The projection of $\mbf{sp}_T$ under $P_{\Phi_0}\to P_{\Phi_0,h}$ defines a special point $\mbf{sp}_{T,h}:(T,X_T)\to (P_{\Phi_0,h},D_{\Phi_0,h})$, and we can assume $T$ is maximal without loss of generality. Denote by $\Q_{X_T}$ the $1$-motive determined by $X_T\hookrightarrow D_{\Phi_0}$. For any $\tau\in \gal(\overline{E}_T/E_T)$, by the description in \cite[11.4, p.188]{Pin89}, $Z_0(R)$ commutes with $\tau $ for any finite $\bb{Q}$-algebra $R$ when we view them as elements in $\aut_{\bb{Q}}\Q_{X_T,\overline{E}_T}$. 
Consider the orbit of conjugation of $\mbf{sp}_T$ by $P_{\Phi_0}(\bb{Q})$ in the orbit of conjugation of it by $P_{\Phi_0}(\bb{R})U_{\Phi_0}(\bb{C})$. As $W_{\Phi_0}(\bb{Q})$ is dense in $W_{\Phi_0}(\bb{R})U_{\Phi_0}(\bb{C})$ under the topology induced by Zariski topology on $W_{\Phi_0}$, as $P_{\Phi,h}(\bb{Q})$ is dense  (and in particular, Zariski dense) in $P_{\Phi_0,h}(\bb{R})$, and as $P_{\Phi_0}(\bb{R})U_{\Phi_0}(\bb{C})^+$ acts transitively on $D_{\Phi_0}^+$, $\gal(\overline{E}_T/E_T)$ commutes with $Z_0(R)$ on the connected component of $\sh_{K_{\Phi_0},\overline{E}_0}$ containing $\mbf{sp}_T$; since the action of $\mbf{p}\in P_{\Phi_0}(\A)$ is transitive on connected components of $\sh_{K_{\Phi_0},\overline{E}_0}$ and is defined over $E_0$, this is true for all connected components. 
By \cite[5.1]{Del71} and \cite[Lem. 11.6]{Pin89} (by varying special points), the action of $\gal(\overline{E}_0/E_0)$ commutes with that of $Z_0(R)$ in $\aut_\bb{Q}\Q_{\eta,\overline{E}_0}$.\par
Finally, let us show the statement over the integral model $\ca{S}_{K_{\Phi_0}}$. 
By Lemma \ref{lem-ext-1-mot} below, any endomorphism of $\Q_\eta$ extends uniquely to an endomorphism of $\Q$. Then there is an embedding $Z_0\hookrightarrow \aut_\bb{Q}\Q$. For any point $\wp=[(x,\mbf{p})]\in \sh_{K_{\Phi_0}}(\bb{C})$ as above, there is a map $Z_{0,\zbkp}\to G_{0,\zbkp}\to \aut_{\zbkp} \Q_x\iso \aut_{\zbkp} \imxh(x)\to \mrm{GL}(V_{\zbkp})$; since the composition is injective, the map $Z_{0,\zbkp}\hookrightarrow \aut_{\zbkp}\Q_x$ is injective. Then there is an embedding $Z_{0,\zbkp}\hookrightarrow \aut_\zbkp \Q$. 
\end{proof}
\begin{rk}\label{rk-twt-z-graded-piece}
    Since any homomorphism between $1$-motives respects weight filtrations and graded pieces, $\z_{\Phi_0}$ induces actions $\z_{\Phi_0,-i}:Z_{0,\zbkp}\to \aut_\zbkp W_{-i}\Q$ and $\gr_{-i}\z_{\Phi_0}:Z_{0,\zbkp}\to \aut_\zbkp \gr^W_{-i}\Q$.
\end{rk}
\begin{lem}[{See also \cite[I. Prop. 2.7]{FC90}}]\label{lem-ext-1-mot}
Let $\Q_1$ and $\Q_2$ be two $1$-motives over $S$, where $S$ is a locally Noetherian and normal scheme. Let $U$ be an open dense subscheme of $S$. Suppose that $f_U:\Q_{1,U}\to \Q_{2,U}$ is a homomorphism between $1$-motives over $U$. Then $f_U$ uniquely extends to a homomorphism $$f:\Q_1\lra \Q_2.$$
\end{lem}
\begin{proof}$f_U^{et}:\ul{Y}_{\Q_1,U}\to \ul{Y}_{\Q_2,U}$ extends uniquely to a homomorphism $f^{et}$ over $S$. By \cite[Prop. 2.7]{FC90}, $f^{sab}_U:\G^\natural_{\Q_1,U}\to \G^\natural_{\Q_2,U}$ extends uniquely to $f^{sab}:\G^\natural_{\Q_1}\to \G^\natural_{\Q_2}$. Since $\G^\natural_{\Q_2}$ is separated and $\ul{Y}_1$ is reduced, $f^{sab}\circ\iota_{\Q_1}=\iota_{\Q_2}\circ f^{et}$. So $f=(f^{et},f^{sab})$ is a homomorphism.
\end{proof}
\subsubsection{}\label{twt-1-mot}
Let $G_{0,\zbkp}^\ad$, $\gamma\in G_{0,\zbkp}^\ad(\zbkp)_{1}$, the $Z_{0,\zbkp}$-torsor $\ca{P}$ at $\gamma$, the finite Galois extension $F$ over $\bb{Q}$, the lifting $\tg\in\ca{P}(\ca{O}_{F,(p)})$, the $1$-cocycle $\{z^{\tg}_\sigma\}_\sigma$, and the $\bb{Z}$-module $\V_\bb{Z}$ be defined as in \S\ref{twt-level} and \S\ref{subsubsec-twt-general-sheaf}.\par
Denote by $\Q^{\ca{O}_F}:=\Hom(\ca{O}_F^\vee,\Q)=(\ul{Y}^{\ca{O}_F},\G^{\natural,\oo},T^\oo,A^\oo,\iota^\oo,c^{\vee,\oo})$ the \emph{Serre tensor construction} for $1$-motives. As before, let $\Q$ be the pullback to $\ca{S}_{K_{\Phi_0}}$ of the tautological family $\Q^\ddag$ on $\ca{S}_{K_{\Phi^\ddag}}$. It makes sense to write $\Q^\oo$ as $\Q\otimes \oo$ viewing them as complexes in $D^b_{fppf}(\ca{S}_{K_{\Phi_0}})$.\par 
As in \S\ref{subsubsec-twt-general-sheaf}, there is an action $\iota_{\Phi,\gal}:\gal(F/\bb{Q})\to \aut_{\op}(\Q)\iso \aut_\zbkp(\Q\otimes\oo)\iso \aut_\zbkp(\Q^\oo)$ induced by $\z_{\Phi_0}$ and $\tg$, sending $\sigma\in\gal(F/\bb{Q})$ to $\z_{\Phi_0}(z_\sigma^{\tg})$; as before, we will omit $\z_{\Phi_0}$ if it is clear in the context.
\begin{construction}\label{const-twt-1-mot}\upshape
Fix an $N>0$ such that $(N,p)=1$ and $N z^{\tg}_\sigma\in \aut \Q^{\ca{O}_F}$. We consider a similar homomorphism $$\ff:=\bigoplus\limits_{\sigma\in\mrm{Gal}(F/\bb{Q})\text{ and }\sigma\neq \mrm{id}}((Nz_\sigma^{\tilde{\gamma}})\sigma-N\cdot\mrm{id}),$$
which is a homomorphism from $\Q^{\oo}$ to $(\Q^{\oo})^{|\gal(F/\bb{Q})|-1}$.\par
By Proposition \ref{twt-kp}, $(\ker \ffab)^\circ=: A^{\tg}$ is an abelian scheme, and $(\ker \fftr)^\circ:=T^{\tg}$ is a torus. Since the kernel of any homomorphism between locally constant groups is locally constant, $\ker \ffet$ is a locally constant abelian group of finite rank. We see that $\G^{\natural,\tg}:=(\ker \ffsab)^\circ$ is a semi-abelian scheme, which is an extension of $A^{\tg}$ by $T^{\tg}$ by the proof of Lemma \ref{lem-semiabelian-twist}. Let $\ul{Y}^{\tg}:= \iota^{\oo,-1}\G^{\natural,\tg}\cap \ker \ffet$.\par
Then the homomorphism $\iota^{\oo}$ induces a homomorphism $\iota^{\tg}:\ul{Y}^{\tg}\to \G^{\natural,\tg}$. This defines a $1$-motive $\Q^{\tg}:=(\ul{Y}^{\tg},\G^{\natural,\tg},T^{\tg},A^{\tg},\iota^{\tg})$, which does not depend on the choice of $N$ satisfying the requirements above: To see this, let $M$ be another positive integer, such that $N|M$ and $(M,p)=1$. Then $\mbf{f}^{\tg}_M=\frac{M}{N}\mbf{f}^{\tg}_N$. We have that $\Q^{\tg}/\Q^{\tg}[\frac{M}{N}]\iso\Q^{\tg}.$ See Lemma \ref{lem-isog-con}.
\end{construction}
\begin{lem}\label{lem-ff-isog}
The kernel of $\ff$ induces an isomorphism $\Q^{\tg}\otimes \op\iso \Q^{\oo}\otimes \zbkp$ of complexes in $D^b_{fppf}(\ca{S}_{K_{\Phi_0}})$. Moreover, $$\bm{\iota}^{\tg}_N:\Q^{\tg,\oo}\lra\Q^{\oo}\otimes{\oo}\lra \Q^{\oo}\otimes_{\ca{O}_F}\ca{O}_F\iso \Q^{\ca{O}_F}$$ 
is a prime-to-$p$ isogeny.
\end{lem}
\begin{proof}
By (\ref{twt-comp}) and Corollary \ref{twt-iso-cor}, the quotient $(\ker {\mbf{f}}_N^{\tg,sab})/\G^{\natural,\tg}$ is trivial after tensoring $\zbkp$, so it is prime-to-$p$. Then the index of $\ul{Y}^{\tg}$ in $\ker {\mbf{f}}_N^{\tg,et}$ is prime-to-$p$.
Thus, by decomposing $\ff$ to $\ffab$, $\fftr$ and $\ffet$, and by Corollary \ref{twt-iso-cor} (resp. Proposition \ref{twt-kp}), we have the first (resp. second) statement. 
\end{proof}
The following statement is an immediate corollary of (the proof of) Lemma \ref{lem-kisin}; we have re-interpreted a similar statement in Lemma \ref{z-equi-lem}.
\begin{cor}\label{cor-equiv-zphi}
Let $(\Q,\bml,[u]_{(K^p_{\Phi_0})^{g^\ddag}})$ be the pullback of the tautological family on $\ca{S}_{K_{\Phi^\ddag}}$ to $\ca{S}_{K_{\Phi_0}}$. Let $R$ be a $\zbkp$-algebra contained in $\bb{C}$.\par 
On the one hand, the embedding $Z_{0,\zbkp}\hookrightarrow G_{0,\zbkp}\hookrightarrow\mrm{GL}(V_{\zbkp})$ induces an action of $Z_{0,\zbkp}(R)$ on $V_{\Ap}\otimes_{\zbkp}R$; on the other hand, the embedding $\z_{\Phi_0}:Z_{0,\zbkp}\hookrightarrow \aut_{\zbkp}\Q$ induces an action of $Z_{0,\zbkp}(R)$ on $V^p\Q\otimes_{\zbkp}R$.\par
Then $u\otimes 1: V_{\Ap}\otimes_\zbkp R\xrightarrow{\ \sim\ }V^p\Q\otimes_\zbkp R$ is $Z_{0,\zbkp}(R)$-equivariant.
\end{cor}
\begin{proof}
We can assume that $R=\bb{C}$ and can just consider the generic fiber $\Q_\eta$ over $\sh_{K_{\Phi_0}}$ without loss of generality. \par
In the notation of Construction \ref{const-1-mot}, for any lifting $(x,\mbf{p})$ of $\wp=[(x,\mbf{p})]\in\sh_{K_{\Phi_0}}(\bb{C})$, the map $u_\wp$ fits into a diagram
\begin{equation}
    \begin{tikzcd}
        V_{\Ap}\otimes_{\zbkp} \bb{C}\arrow[rr,"u_\wp\otimes 1"]\arrow[drr,"{p^{\ddag,p}}"]&&V^p\Q_\wp\otimes_\zbkp \bb{C}\iso H_\dr \Q_x\otimes_\zbkp \Ap\arrow[d,"{\alpha_x\otimes 1}"]\\
        && V_{\Ap}\otimes_\zbkp \bb{C}\iso V_\bb{C}\otimes_\zbkp \Ap.
    \end{tikzcd}
\end{equation}
Since $p^{\ddag,p}$ commutes with $Z(\bb{C})$, it suffices to show that $\alpha_x\otimes 1$ is $Z(\bb{C})$-equivariant. By Lemma \ref{lem-kisin} and its proof, the $\bb{C}$-value of $\z_{\Phi_0}$ is given by the variation of mixed Hodge structure on $V_\bb{C}$.
Fix any point $x_0\in D_{\Phi_0}$. We can choose $\alpha_{x_0}$ to be the identity. Then $\alpha_x$ is induced by multiplying some element in $P_{\Phi_0}(\bb{C})$, which commutes with $Z(\bb{C})$. 
\end{proof}
Now, we can construct the twisting of $(\Q,\bml,[u]_{(K_{\Phi_0})^{g^\ddag}})$.
\begin{construction}\label{const-twt-1-mot-lvl}\upshape
Consider 
\begin{equation}\label{dia-twt-1-mot-lvl}
    \begin{tikzcd}
{\V_{\zhp}\otimes \ca{O}_F\sbst\V_{\Ap}\otimes_{\zbkp}\ca{O}_{F,(p)}}\arrow[r,"{({\tg}^{-1}\cdot)}","\sim"']&
{V_{\Ap}\otimes_{\zbkp}\ca{O}_{F,(p)}}\arrow[r,"{u\otimes 1}","\sim"']& V^p\Q^{\oo}\arrow[r,"V^pf"]&V^p\Q^\prime,
\end{tikzcd}
\end{equation}
where $f:\Q^{\oo}\to \Q^\prime$ is a $\zbkpt$-isogeny determined by Lemma \ref{lem-cons}.\par
\begin{enumerate}[leftmargin=*,label=(\arabic*)]
\item We can claim that $(T^p\Q^\prime)^{\gal(F/\bb{Q})}$ comes from the Tate module of a $1$-motive $\Q^{\tg,\prime}$, such that $\Q^{\tg,\prime}$ is $\zbkpt$-isogenous to $\Q^{\tg}$. In fact, by Corollary \ref{cor-equiv-zphi}, we take the $\gal(F/\bb{Q})$-invariant part of both $\V_{\zhp}\otimes\oo$ and $T^p\Q^\prime$, and we have that
$\V_{\zhp}$ is an open compact subgroup of $V^p\Q^{\tg}$ under $(u\otimes 1)\circ(\tg^{-1})$, such that $(u\otimes 1)\circ(\tg^{-1})(\V_{\zhp})_{\overline{s}}$ is $\pi_1(S,\overline{s})$-invariant. Hence, there is a $\zbkpt$-isogeny $$f^{\tg}:\Q^{\tg}\lra \Q^{\tg,\prime}$$ by Lemma \ref{lem-cons}.\par
This construction only depends on $\gamma$, but does not depend on the lifting $\tg$. This is by the same argument as Lemma \ref{lem-indep}. In fact , if $\tg^\prime=z\tg$, there is a commutative diagram
\begin{equation}
    \begin{tikzcd}
    \V_{\Ap}\otimes_{\zbkp}\ca{O}_{F,(p)}\arrow[rr,"{({\tg}^{-1}\cdot)}"]\arrow[drr,"{({\tg}^{\prime,-1}\cdot)}"]&& V_{\Ap}\otimes_{\zbkp}\ca{O}_{F,(p)}\arrow[d,"z^{-1}"]\arrow[rr,"{u\otimes 1}"]&& V^p\Q^{\ca{O}_F}\arrow[d,"z^{-1}"]\\
    && V_{\Ap}\otimes_{\zbkp}\ca{O}_{F,(p)}\arrow[rr,"{u\otimes 1}"]&&V^p\Q^{\ca{O}_F}.
    \end{tikzcd}
\end{equation}
The images of $\V_{\zhp}$ in $V^p\Q^{\oo}$ under $(u\otimes 1)\circ(\tg^{-1}\cdot)$ and $(u\otimes 1)\circ(\tg^{\prime,-1}\cdot)$ differ by $z^{-1}$. Thus, they determine the same $\Q^\prime$ and $\zbkpt$-isogenies $f_1:\Q^{\oo}\to \Q'$ and $f_2:\Q^{\oo}\to \Q'$ such that $f_1=f_2\circ z^{-1}$.\par
Then we can define the twisted level structure as $[u^{\gamma}]_{(K_{\Phi_0})^{g^{\ddag}\gamma^{-1}}}$, where $u^{\gamma}:=(V^pf\circ (u\otimes 1)\circ(\tg^{-1}\cdot))^{\gal(F/\bb{Q})}|_{\V_{\zhp}}$. This is indeed independent of choice of $\tg$ lifting $\gamma$ by the paragraph above. Let $\Q^\gamma:=\Q^{\tg,\prime}$. We will see later that this is indeed a level structure as Definition \ref{def-deg-lvl}.
\item 
We now define the twisted polarization. The treatment is similar to that in \S\ref{twt-pol-revised}, so we will use the notation there.
More precisely, consider the diagram 
\begin{equation}\label{diag-twt-pol-1-mot}
    \begin{tikzcd}
    {\V_{\Ap}\otimes_{\zbkp} \ca{O}_{F,(p)}}\arrow[r,"{(\tg^{-1}\cdot)}"]\arrow[d,"\psi^{\V}\otimes 1"]&V_{\Ap}\otimes_{\zbkp}\ca{O}_{F,(p)}\arrow[r,"u\otimes 1"]\arrow[d,"\psi\otimes 1"]& V^p\Q^{\oo}\arrow[r,"V^pf"]\arrow[d,"{V^p{\bml}\otimes 1}"]& V^p\Q'\arrow[d,dashed,"{V^p{\bml}^{\V,\oo}}"]\\
{\V_{\Ap}^\vee\otimes_{\zbkp} \ca{O}_{F,(p)}(\nu)}\arrow[r]&{V_{\Ap}^\vee\otimes_{\zbkp}\ca{O}_{F,(p)}(\nu)}\arrow[r]& V^p\Q^{\vee}\otimes_{\zbkp}\ca{O}_{F,(p)}(\nu)\arrow[r,dashed,"{V^pf'}"]& V^p\Q^*(\nu)
    \end{tikzcd}
\end{equation}
The image of $\V_{\zhp}\otimes \ca{O}_F$ (resp. $\V_{\zhp}^\vee\otimes \ca{O}_F$) in $V^p\Q^{\oo}$ (resp. $V^p\Q^\vee\otimes_{\zbkp}\ca{O}_{F,(p)}$) via $(u\otimes 1)\circ(\tg^{-1}\cdot)$ (resp. $(V^p\bml\otimes 1)\circ(u\otimes 1)\circ (\tg^{-1}\cdot)\circ (\psi^{\V}\otimes 1)$) determines a 
$\zbkpt$-isogeny $f$ (resp. $f'$) together with a prime-to-$p$ isogeny $\bml^{\V,\oo}:\Q'\to \Q^*$ by Lemma \ref{lem-cons} such that the right square commutes. \par
By Corollary \ref{cor-equiv-zphi}, $u\otimes 1$ is $Z_{0,\zbkp}(\ca{O}_{F,(p)})$-equivariant. Since $u$ sends $\psi_{\zhp}$ to a $\wat{\bb{Z}}^{p,\times}$-multiple of $\e^{\bml}$, $V^p\bml\otimes 1$, $V^pf'$ and $V^p\bml^{\V,\oo}$ are $Z_{0,\zbkp}(\ca{O}_{F,(p)})$-equivariant and therefore are $\gal(F/\bb{Q})$-equivariant.\par
Note that, for $\Q^\vee$, we can also define a $1$-motive $(\Q^{\vee})^{\tg}$ by the $Z_{0,\zbkp}$-action given by the canonical isomorphism (\ref{eq-aut-dual-iso}). Since $$V^p(\Q^\vee)^{\tg}= (V^p{\bml}^{\V,\oo}\circ V^p f\circ(u\otimes 1)\circ({\tg}^{-1}\cdot)\circ(\psi^{\V}\otimes 1)^{-1})(\V^\vee_{\Ap}\otimes_{\zbkp}\ca{O}_{F,(p)}(\nu))^{\gal(F/\bb{Q})},$$
the image of $\V^\vee_{\zhp}$ in $V^p(\Q^\vee)^{\tg}$ induces a $\zbkpt$-isogeny $(\Q^\vee)^{\tg}\to \Q''$ such that $T^p\Q''=(T^p\Q^*)^{\gal(F/\bb{Q})}$. There is also a prime-to-$p$ isogeny $\bml^{\tg,\prime}:\Q^{\gamma}\to \Q''$ which is the $\gal(F/\bb{Q})$-invariant of $\bml^{\V,\oo}:\Q'\to \Q^*$.\par
We now check that $\bml^{\tg,\prime}$ is a polarization. Indeed, choose any positive definite perfect pairing $\e_0:\oo\times\oo\to \bb{Z}$ that canonically induces $i_F:\ca{O}_F\xrightarrow{\sim}\mrm{Hom}_\bb{Z}(\oo,\bb{Z})$. As in Lemma \ref{lem-twt-pol-new}, $\bml(i_F)=i_F\circ\bml$ is a polarization implies that $\bml^{\V,\oo}(i_F):\Q'\to\Q^*\to\Q^*\otimes_{\oo,i_F}\ca{O}_F^\vee$ is a polarization. By construction, this implies that $\Q^{\gamma}\otimes \oo\xrightarrow{\bml^{\tg,\prime}\otimes 1}\Q''\otimes\oo\xrightarrow{\sim}\Q''\otimes\ca{O}_F^\vee$ is a polarization. Taking the $\gal(F/\bb{Q})$-invariant part, we know that $\Q^{\tg,\prime}\otimes \bb{Z}\xrightarrow{\bml^{\tg,\prime}\otimes 1}\Q''\otimes \bb{Z}\xrightarrow{\sim}\Q''\otimes \mrm{Hom}_\bb{Z}(\bb{Z},\bb{Z})$ and $\bml^{\tg,\prime}$ are polarizations.\par
From now on, denote $\bml^\gamma:=\bml^{\tg,\prime}$.
\item Let us check that $[u^{\gamma}]_{(K_{\Phi_0})^{g^\ddag\gamma^{-1}}}$ is a level structure defined as in Definition \ref{def-deg-lvl}.
By the previous part, $u^{\gamma}$ sends $\V_{\zhp}$ to $T^p\Q^\gamma$, and sends $\psi^\V_{\zhp}$ to a $\wat{\bb{Z}}^{p,\times}$-multiple of $\e^{\bml^\gamma}$. \par
Note that $f:\Q^{\ca{O}_F}\to \Q^\prime$ respects weight filtrations $W^{\oo}$ and $W^\prime$ on $\Q^{\oo}$ and $\Q^\prime$, and respects their graded pieces. Then $f$ induces $\zbkpt$-isogenies $\gr_0f:\uhom(\ul{X}^{\oo},\bb{Z}(1))\to\uhom(\ul{X}^\prime,\bb{Z}(1))$ and $\gr_{-2}f:\ul{Y}^{\oo}\to \ul{Y}^\prime$. 
\end{enumerate}
Now we can construct the data we need:
\begin{itemize}[leftmargin=*]
\item $\mbf{Z}^{(\gamma g^\ddag \gamma^{-1}),\ca{O}_F}:= \V_{\zhp}\otimes \ca{O}_F\cap (\tg^{-1})^* 
    Z^{(g^\ddag)}_{\Ap}\otimes \ca{O}_F$; we take $\gal(F/\bb{Q})$-invariant part and obtain a filtration $\mbf{Z}^{(\gamma g^\ddag \gamma^{-1})}$ of $\V_{\zhp}$. It is a filtration determined by the parabolic group $\gamma g^{\ddag,-1}Qg^{\ddag}\gamma^{-1}$, which explains the superscript in our notation. We denote $\gamma g^\ddag \gamma^{-1}$ by $g^\dagger$ for simplicity.\par
\item Consider $\alpha_0^{\tg,\ca{O}_{F,(p)}}:\gr_0^{F^{(g^\dagger)}}\V_\bb{Z}^{(g^\dagger)}\otimes\ca{O}_{F,(p)}\xrightarrow{\gr_0\tg^{-1}}\gr_0^{F^{(g^\ddag)}}V_\bb{Z}^{(g^\ddag)}\otimes \ca{O}_{F,(p)}\xrightarrow{\alpha_0^{(g^\ddag)}}\uhom(\ul{X},\bb{Z}(1))\otimes\ca{O}_{F,(p)}\xrightarrow{\gr_0f}\uhom(\ul{X}^\prime,\bb{Z}(1))\otimes \zbkp$. 
The composition $\alpha^{\tg,\ca{O}_{F,(p)}}_{0,\Ap}\circ (\gr_0 g^\dagger)_{\Ap}$ is $\gr_0 u^{\gamma,\ca{O}_F}_{\Ap}$, so it sends $\gr_0^{\mbf{Z}^{(g^\dagger)}}\V_{\zhp} \otimes \oo$ to $\uhom(\ul{X}',\bb{Z}(1))\otimes \zhp$. 
We then see that the restriction of $\alpha_0^{\tg,\ca{O}_{F,(p)}}$ defines a map $\alpha_0^{\tg,\oo}:\gr_0^{F^{(g^\dagger)}}\otimes\oo\to \uhom(\ul{X}',\bb{Z}(1))$. Taking the $\gal(F/\bb{Q})$-invariant part, there is a map over $\bb{Z}$, $\alpha_0^{\tg}:\gr_0^{F^{(g^\dagger)}}\to \uhom(\ul{X}^\gamma,\bb{Z}(1))$. Denote $\alpha_0^{(g^\dagger)}:=\alpha_0^{\tg}$.
\item Consider $\alpha_{-2}^{\tg,\ca{O}_{F,(p)}}:\gr_{-2}^{F^{(g^\dagger)}}\V_\bb{Z}^{(g^\dagger)}\otimes\ca{O}_{F,(p)}\xrightarrow{\gr_{-2}\tg^{-1}}\gr_{-2}^{F^{(g^\ddag)}}V_\bb{Z}^{(g^\ddag)}\otimes \ca{O}_{F,(p)}\xrightarrow{\alpha_0^{(g^\ddag)}}\ul{Y}\otimes\ca{O}_{F,(p)}\xrightarrow{\gr_{-2}f}\ul{Y}'\otimes \zbkp$. 
The composition $\alpha^{\tg,\ca{O}_{F,(p)}}_{-2,\Ap}\circ (\gr_{-2} g^\dagger)_{\Ap}$ is a $\wat{\bb{Z}}^{p,\times}$-multiple of $\gr_{-2} u^{\gamma,\ca{O}_F}_{\Ap}$, so it sends $\gr_{-2}^{\mbf{Z}^{(g^\dagger)}}\V_{\zhp} \otimes \oo$ to $\ul{Y}'\otimes \zhp$. 
We then see that the restriction of $\alpha_{-2}^{\tg,\ca{O}_{F,(p)}}$ defines a map $\alpha_{-2}^{\tg,\oo}:\gr_{-2}^{F^{(g^\dagger)}}\otimes\oo\to \ul{Y}'$. Taking the $\gal(F/\bb{Q})$-invariant part, there is a map over $\bb{Z}$, $\alpha_{-2}^{\tg}:\gr_{-2}^{F^{(g^\dagger)}}\to \ul{Y}^\gamma$. Denote $\alpha_{-2}^{(g^{\dagger})}:=\alpha_{-2}^{\tg}$.
\end{itemize}
Then $[u^{\gamma}]_{(K_{\Phi_0})^{g^\ddag\gamma^{-1}}}$ satisfies Definition \ref{def-deg-lvl}. See Proposition \ref{prop-twt-clr} below. \par
Let $\Psi_0=\lcj{\Phi_0}{\gamma}:=(\gamma Q\gamma^{-1},\gamma(X_0^+),\gamma g^\ddag\gamma^{-1})$; $\Psi_0$ maps to a cusp label representative $\Psi^\ddag$ of $(G^\ddag,X^\ddag)$.
Let $K^{\odot,p}:=K^{\ddag,p,\prime}$ be any neat open compact subgroup of $G^\ddag(\Ap)$ containing $\lcj{K^p_0}{\gamma}$ and stabilizing $\V_{\zhp}$. Let $K^\odot:=K_p^\ddag K^{\odot,p}$. 
Denote $K_{\Psi^\ddag}^{\odot,p}:= g^\dagger K^{\odot,p} g^{\dagger,-1}\cap P_{\Psi^\ddag}(\Ap)$. Then we see that $K_{\Psi^\ddag}^{\odot,p}$ contains $\lcj{K_{\Phi_0}^p}{\gamma}=K_{\Psi_0}^p$. In fact, we compute as follows:
$$g^\dagger K^{\odot,p}g^{\dagger,-1}\cap P_{\Psi^\ddag}(\Ap)\supset \gamma g^{\ddag}K_0^pg^{\ddag,-1}\gamma^{-1}\cap P_{\Psi^\ddag}(\Ap)=\gamma(g^\ddag K_0^pg^{\ddag,-1}\cap P_{\Phi^\ddag}(\Ap))\gamma^{-1}\supset \gamma K_{\Phi_0}^p\gamma^{-1}.$$
Similarly, $K_{\Psi^\ddag,p}^\odot$ contains $\lcj{K_{\Phi_0,p}}{\gamma}=K_{\Psi_0,p}$. We then have $(K_{\Phi_0}^p)^{g^\ddag \gamma^{-1}}\sbst(K_{\Psi^\ddag}^{\odot,p})^{g^\dagger}$.
\begin{prop}\label{prop-twt-clr}
With the constructions above, $[u^{\gamma}]_{(K^{\odot,p}_{\Psi^\ddag})^{g^\dagger}}:=[u^{\gamma}]_{(K_{\Phi_0})^{g^\ddag\gamma^{-1}}}/(K^{\odot,p}_{\Psi^\ddag})^{g^\dagger}$ is a level structure of type $(\V_\bb{Z},\psi^\V_\bb{Z},\Psi^\ddag,K^{\odot,p})$. 
\end{prop}
\begin{proof}
    This is a consequence of the construction above.
\end{proof}
This ends the Construction \ref{const-twt-1-mot-lvl}.\hfill $\square$
\end{construction}
\subsubsection{}Twisting of $1$-motives sends one cusp label representative of $(G_0,X_0)$ to another one. \par
Denote by $\ca{S}_{K_{\Psi_0}^p}^{\odot}$ the inverse limit $\varprojlim \ca{S}_{K^\odot_{\Psi^\ddag}}\iso \varprojlim\bm{\xi}_{(\V_\bb{Z},\psi^\V_\bb{Z}),\Psi^\ddag,K^{\odot,p}}$, where the inverse system runs over neat open compact subgroups $K^{\odot,p}_{\Psi^\ddag}$ as above containing $K^p_{\Psi_0}$; similarly, let $\sh_{K_{\Psi_0}^p}^{\odot}:=\ca{S}_{K_{\Psi_0}^p,\bb{Q}}^{\odot}$. \par 
\begin{prop}\label{prop-twt-1}
With the constructions above, 
\begin{enumerate}
\item The twisting construction $$(\Q,\bml,[u]_{(K_{\Phi_0})^{g^\ddag}})\mapsto (\Q^\gamma,\bml^\gamma,[u^{\gamma}]_{(K_{\Phi_0})^{g^\ddag\gamma^{-1}}})$$ induces a morphism $\gamma^{-1}:\ca{S}_{K_{\Phi_0}}\to\ca{S}_{K_{\Psi}^p}^{\odot}$ by the universal property.
\item Over $\sh_{K_{\Phi_0}}(\bb{C})\iso P_{\Phi_0}(\bb{Q})\bss D_{\Phi_0}\times P_{\Phi_0}(\A)g^\ddag K/K$, $\gamma^{-1}$
is identical to a map $\gamma^{-1}:\sh_{K_{\Phi_0}}(\bb{C})\to \sh_{K_\Psi^p}^{\odot}(\bb{C})$, sending $\wp=[(x,\mbf{p}g^\ddag)]$ to $\wp^\gamma=[(\gamma(x),\gamma \mbf{p}g^\ddag\gamma^{-1})]$. Hence, it factors as $\gamma^{-1}:\sh_{K_{\Phi_0}}(\bb{C})\xrightarrow{\sim}\sh_{K_{\Psi_0}}(\bb{C})\xrightarrow{\iota_{\Psi_0}} \sh_{K_{\Psi_0}^p}^\odot(\bb{C})$, where the first map is an isomorphism sending $\wp=[(x,\mbf{p}g^\ddag)]$ to $\wdtd{\wp}^\gamma=[(\gamma(x),\gamma \mbf{p}g^\ddag\gamma^{-1})]$.
\item  The map $\gamma^{-1}:\ca{S}_{K_{\Phi_0}}\to\ca{S}_{K_{\Psi_0}^p}^\odot$ factors as $\gamma^{-1}:\ca{S}_{K_{\Phi_0}}\xrightarrow{\sim}\ca{S}_{K_{\Psi_0}}\xrightarrow{\iota_{\Psi_0}} \ca{S}_{K_{\Psi_0}^p}^\odot$.
\end{enumerate}
\end{prop}
\begin{proof}
Given the first and the second part, the third part follows from the construction of $\ca{S}_{K_{\Psi_0}}$ as the normalization of $\ca{S}_{K_{\Psi^\ddag}^\odot}$ in $\sh_{K_{\Psi_0}}$, and the fact that this normalization is independent of the choice of different $K_{\Psi^\ddag}^{\odot,p}$ as above. The first part follows from Construction \ref{const-twt-1-mot-lvl}.\par 
Let us show the second part. This is the same as Kisin's argument in \cite[Lem. 3.2.6]{Kis10}, and let us explain it now.\par 
Recall that we are working over $\bb{C}$. Let $\Q_\wp$ be the $1$-motive at a complex point $\wp=[(x,\mbf{p}g^\ddag)]\in \sh_{K_{\Phi_0}}(\bb{C})$, and let $\Q_{\wp^\gamma}$ be the $1$-motive corresponding to $\wp^\gamma=[(x^\gamma,\gamma \mbf{p}g^\ddag\gamma^{-1})]$. 
First, there is a diagram:
\begin{equation}\label{diag-key}
    \begin{tikzcd}
T_\bb{Z}\Q^{\oo}_\wp\otimes\zbkp\arrow[rr,"f"]\arrow[d,"{\alpha_x\otimes 1}"]&& T_\bb{Z}(\Q^{\gamma,\oo}_{\wp})\otimes\zbkp\arrow[d,dashed]\\
(V_\bb{Z}\otimes \ca{O}_{F,(p)},r(g^\ddag)^{-1}\psi\otimes 1,W)\arrow[rr,"{\tg}"]&& (\V_\bb{Z}\otimes\ca{O}_{F,(p)},r(g^\dagger)^{-1}\nu(\tg^{-1})\psi^\V\otimes 1,\mbf{W}),
    \end{tikzcd}
\end{equation}
where $W$ (resp. $\mbf{W}$) is the filtration determined by $Q^\ddag$ (resp. $\gamma Q^\ddag\gamma^{-1}$) and note that $r(g^\dagger)=r(g^\ddag)$. 
Denote $\mbf{T}:=\tg\circ(\alpha_x\otimes 1)\circ f^{-1}$. By Construction \ref{const-twt-1-mot-lvl} and Construction \ref{const-1-mot}, $\mbf{T}$ sends that weight filtration of $T_\bb{Z}(\Q_\wp^{\gamma,\oo})\otimes \bb{Q}$ to $\mbf{W}$.  
The mixed Hodge structures on $V\otimes\oo\iso \V\otimes\oo$ determined by $x$ and 
$x^\gamma$ viewed as elements in $\mrm{Hom}(\bb{S}_\bb{C},V_\bb{C}\otimes\oo)$ are conjugated by $\tg$. Then the bottom row of the diagram above respects mixed Hodge structures. Hence, $\mbf{T}$ sends the mixed Hodge structure on $T_\bb{Z}(\Q_\wp^{\gamma,\oo})\otimes\bb{Q}$ to the mixed Hodge structure on $\V\otimes\oo$ which corresponds to $x^\gamma$.\par
Then we combine the diagram above with (\ref{dia-twt-1-mot-lvl}) and (\ref{diag-level-cons}) to get the following commutative diagram:
\begin{equation}
    \begin{tikzcd}
{\V_{\zhp}\otimes\oo}\arrow[r,"{\tg^{-1}}"]&V_{\Ap}\otimes\oo\arrow[r,"u"]\arrow[dr,"p^{\ddag,p}"]& T_\bb{Z}\Q^{\oo}_\wp\otimes\Ap\arrow[r,"f"]\arrow[d,"{\alpha_x\otimes 1}"]& T_\bb{Z}(\Q^{\gamma,\oo}_\wp)\otimes\Ap\arrow[d]\\
&&
(V_{\Ap}\otimes \oo,\wat{\bb{Z}}^{p,\times}r(g^\ddag)^{-1}\nu(\tg)\psi,W)\arrow[r,"{\tg}"]& (\V_{\Ap}\otimes\oo,\wat{\bb{Z}}^{p,\times}r(g^\dagger)^{-1}\psi^\V,\mbf{W}).&
    \end{tikzcd}
\end{equation}
The first term is equipped with the pairing $\psi^\V\otimes 1$ and the second one is equipped with the pairing $\nu(\tg)\psi\otimes 1$. 
By adjusting $\mbf{p}$ by an element in $P_{\Phi^\ddag}(\bb{Q})$, 
we see that the lattice $\V_{\zhp}\otimes \oo$ on the top-left is sent to $T^p(\Q^{\gamma,\oo}_\wp)$ via the series of morphisms on the first row, and is sent to $\V_{\zhp}^{(g^\dagger)}\otimes \oo$ via $\tg^{-1} \circ p^{\ddag,p}\circ\tg$. 
Then the second part follows from taking the $\gal(F/\bb{Q})$-invariant of the diagram above and (\ref{diag-level-cons}).
\end{proof}
\subsection{Some actions of the adjoint group on toroidal compactifications}\label{subsec-ext}
In this subsection, we extend the twisting construction in previous subsections to toroidal compactifications of integral models of Hodge-type Shimura varieties, and relate another ``twisting'' construction on the toroidal compactifications of Hodge-type Shimura varieties $\sh_{K_0}^{\Sigma}$ to quasi-isogeny twists in \cite{Lan16b}. The notation and the construction in this subsection will be used in the proof of one of the main technical results \S\ref{subsubsec-main-3}.\par
Let $\gamma\in G^\ad_0(\bb{Q})$ and $g\in G_0(\A)$. Let $\beta:=(g,\gamma^{-1})$ and $(K_0)^\beta:=\lcj{((K_0)^g)}{\gamma}$. 
Then $\beta$ induces a morphism $\beta:\sh_{K_0}(G_0,X_0)\to \sh_{(K_0)^\beta}(G_0,\gamma\cdot X_0)$ described by 
$$[(x,g_0)]\mapsto [(\gamma(x),\gamma g_0g \gamma^{-1})]$$ 
on complex points.\par 
For any admissible cone decomposition $\Sigma_2$ of $\sh_{(K_0)^\beta}:=\sh_{(K_0)^{\beta}}(G_0,\gamma\cdot X_0)$, let $\Sigma_1:=\beta^* \Sigma_2$ be the admissible cone decomposition of $\sh_{K_0}:=\sh_{K_0}(G_0,X_0)$ induced by $\beta$.
By \cite[Prop. 6.25 and Thm. 12.4]{Pin89}, there is a morphism from $\sh_{K_0}^{\Sigma_1}:=\sh_{K_0}^{\Sigma_1}(G_0,X_0)$ to $\sh_{(K_0)^\beta}^{\Sigma_2}:=\sh_{(K_0)^\beta}^{\Sigma_2}(G_0,\gamma \cdot X_0)$ extending the morphism of interiors $\beta: \sh_{K_0}(G_0,X_0)\to \sh_{(K_0)^\beta}(G_0,\gamma\cdot X_0)$. We still denote this extension by $\beta:\sh_{K_0}^{\Sigma_1}\to \sh^{\Sigma_2}_{(K_0)^\beta}$.\par
Fix a cusp label representative $\Phi_0$ mapping to $\Phi^\ddag$ as before. Then $\beta$ induces an action $\beta:\sh_{K_{\Phi_0}}\xrightarrow{g}\sh_{K_{\Phi^g}}\xrightarrow{\gamma^{-1}} \sh_{K_{\Psi'}}$, corresponding to an action $(Q,X_0^+,g^\ddag)\mapsto \Phi^g:=(Q,X^+_0,g^\ddag g)\mapsto \Psi':=(\gamma Q\gamma^{-1},\gamma(X_0^+),\gamma g^\ddag g\gamma^{-1})$ on cusp label representatives. Over complex points, the action of $\beta$ sends $[(x,\mbf{p}g^\ddag)]$ to $[(\gamma(x),\gamma \mbf{p}g^\ddag g \gamma^{-1})]$. \par
\subsubsection{}\label{subsubsec-extend-tor}
With the constructions in previous subsections, the first result we intended to show is now immediate thanks to a lemma of Madapusi \cite[Lem. A.3.4]{Mad19}. Note that we have fixed a Hodge embedding $\iota:(G_0,X_0,K_0)\hookrightarrow (G^\ddag,X^\ddag,K^\ddag)$ as \S\ref{subsec-siegel}.
\begin{lem}\label{lem-action-gadtrivial}
The action $\beta:\sh_{K_0}^{\Sigma_1}\to \sh^{\Sigma_2}_{(K_0)^\beta}$ extends to an action
$$\beta: \ca{S}_{K_0}^{\Sigma_1}\lra \ca{S}^{\Sigma_2}_{(K_0)^\beta},$$
which is compatible with the actions $\beta:\ca{S}_{K_{\Phi_0}}\to \ca{S}_{K_{\Psi'}}$ on the boundaries if $\gamma\in G_{0,\zbkp}^{\ad}(\zbkp)_1$.
\end{lem}
\begin{proof}
Combine Proposition \ref{prop-twt-1}, \cite[Prop. 6.25 and Thm. 12.4(b)]{Pin89}, \cite[Lem. A.3.4]{Mad19} and \cite[4.1.12]{Mad19}. 
\end{proof}
\subsubsection{}\label{subsubsec-ggamma-trivial-case}
We now \textbf{assume} that $\gamma\in G^\ad_0(\bb{Q})_1$ and that $g\gamma^{-1}$ is trivial in $G^\ad_0(\A)$. 
Under this assumption, we have assumed that $(K_0)^\beta=K_0$, and therefore both $\Phi_0$ and $\Psi'$ are representatives of some cusps in $\cusp_{K_0}(G_0,X_0)$. 
We can fix an embedding $\iota:(G_0,X_0,K_0)\hookrightarrow(G^\ddag,X^\ddag,K^\ddag)$ satisfying \cite[Sec. 3.1]{Mad19}; applying the theory of \cite{Mad19}, this choice of the embedding $\iota$ has \emph{fixed} a construction of integral models $\ca{S}_{K_{\Phi_0}}$ for all $\Phi_0\in\ca{CLR}(G_0,X_0)$.\par 
In the following construction, we will mainly work with $\sh_{K_{\Phi_0}}$'s instead of their integral models. 
\begin{construction}\label{cons-general-twisting}\upshape
Let $F$ be a finite Galois extension of $\bb{Q}$ such that $\gamma$ is lifted to $\tg\in G_0(F)$.\par
We view $V^F:=V\otimes_\bb{Q}F$ as a $\bb{Q}$-vector space by fixing a $\bb{Q}$-isomorphism $V^F\xrightarrow{\sim}V^{\oplus[F:\bb{Q}]}$. 
Denote by $\psi^F$ the symplectic pairing on $V^{\oplus [F:\bb{Q}]}$ which is $\psi$ on each factor, and elements from distinct factors are perpendicular.

Define $G^F:=\mrm{GSp}(V^F,\psi^F)$ and let $X^F$ be the union of Siegel upper and lower half spaces defined by $G^F$.\par
Then there is a sequence of embeddings of Shimura data:
\begin{equation*}
\begin{tikzcd}
\iota^+:(G_0,X_0)\arrow[r,hook,"\iota"]& (G^\ddag,X^\ddag)\arrow[r,hook,"\Delta"]&(G^F,X^F).
\end{tikzcd}
\end{equation*}
The morphism $(G^\ddag,X^\ddag)\to (G^F,X^F)$ is defined by the fixed diagonal embedding.\par
On the other hand, let $\nu: R_{F/\bb{Q}}G^\ddag\to R_{F/\bb{Q}}\bb{G}_m$ be the similitude character. The value of $\nu(\tg)$ remains unchanged if we extend $F$. We then choose $F$ to be a finite Galois extension such that it contains the square roots of $\nu(\tg)$. We choose a square root $c(\tg)$ of $\nu(\tg)$. We view $c(\tg)$ as a $\bb{Q}$-point of $R_{F/\bb{Q}}\bb{G}_m\hookrightarrow R_{F/\bb{Q}}G^\ddag$. Denote $\tg_1:=\tg^{-1}c(\tg)$. We view $\tg_1$ as an element in $\mrm{GL}(V^F).$ \par
There is a commutative diagram of algebraic groups:
\begin{equation}\label{diag-group-conj}
    \begin{tikzcd}
    G \arrow[r,"{\mrm{int}(\tg_1^{-1})}"]\arrow[d,"\iota"]& G\arrow[dd,"\mrm{int}(\tg_1^{-1})\iota^+"]\\
G^\ddag\arrow[d,"\Delta"]&\\
G^F\arrow[r,"\mrm{int}(\tg_1^{-1})"]& \mrm{int}(\tg_1^{-1})G^F.    
    \end{tikzcd}
\end{equation}

The operation $(\tg_1\cdot)\circ (g\cdot)$ on $V^F_{\A}$ induces a collection of morphisms 
$$ \sh_{K_{\Phi^F}}\lra \sh_{K_{(\Phi^F)'}},$$
which sends any cusp label representative $\Phi^F:=(Q^F,X^{F,+},g^F)$ of $(G^F,X^F)$ to 
$$(\Phi^F)':=(\tg_1^{-1}Q^F{\tg_1},\tg_1^{-1}(X^{F,+}),\tg_1^{-1}g^{F}g\tg_1)$$ 
of $\mrm{int}(\tg_1^{-1})(G^F, X^F)$.\par
Over complex points, the map is $[(x^F,g^F_0g^F)]\mapsto [(\tg_1^{-1}x^F,\tg_1^{-1}g^F_0g^Fg\tg_1)]$. Note that this action remains unchanged if we change $\tg_1$ to $\tg_1\cdot c$ for any $c\in \bb{G}_m(\bb{Q})\hookrightarrow G^F(\bb{Q})$ in the center of $G^F$. In particular, it remains unchanged if we choose another square root of $\nu(\tg)$.\par
Suppose that $\Phi_0\mapsto\Phi^\ddag\mapsto \Phi^F$ via $\iota^+$ and $\Psi'\mapsto \Psi^{\prime,\ddag}\mapsto (\Phi^F)'$ via $\mrm{int}(\tg_1^{-1})\iota^+$, and if we choose a suitable neat open compact subgroup $K^F\sbst G^F(\A)$ containing $K^\ddag$, we have a commutative diagram
\begin{equation}\label{diag-bd-mix-sh}
    \begin{tikzcd}
{\sh_{K_{\Phi_0}}}\arrow[rr,"{\beta}"]\arrow[d,"\iota_\Phi"]&& {\sh_{K_{\Psi'}}}\arrow[dd,"(\mrm{int}(\tg_1^{-1})\iota^+)_{\Psi'}"]\\
\sh_{K_{\Phi^\ddag}}\arrow[d]&&\\
\sh_{K_{\Phi^F}}\arrow[rr,"{\mrm{int}(\tg_1^{-1})\circ g}"]&& \sh_{K_{(\Phi^F)'}}.
    \end{tikzcd}
\end{equation}
We choose a hyperspecial subgroup $K^F_p$ of $G^F(\bb{Q}_p)$ such that $K^F_p\cap G^\ddag(\bb{Q}_p)=K^\ddag_p$ and choose a neat open compact subgroup $K^{F,p}\sbst G^F(\Ap)$ containing $K^{\ddag,p}$ and stabilizing $V^F_{\zhp}:=V_{\zhp}\otimes\oo$.\par
Since, by construction, $g^\ddag$ maps to $g^F$ and $Q^F$ is the unique proper admissible $\bb{Q}$-parabolic that contains $Q$ and $Q^\ddag$, we have $Q\sbst Q^\ddag\sbst Q^F$. Since the conjugation of $\tg_1^{-1}$ on $G_0$, viewed as a subgroup of $G^F$, is the conjugation of $\tg$, we have that $\tg_1^{-1}Q^F\tg_1$ contains $\tg Q \tg^{-1}$. Hence, $\tg_1^{-1} Q^F \tg_1$ is the unique proper admissible $\bb{Q}$-parabolic that contains $\tg Q \tg^{-1}$ since $\tg_1^{-1}Q^F \tg_1$ depends only on the cocharacter that defines $\tg Q\tg^{-1}$ (see \cite[4.16]{Pin89}).\par 
We have the following commutative diagram induced by (\ref{diag-bd-mix-sh}):
\begin{equation}\label{diag-twisting-check}
\begin{tikzcd}
{P_{\Phi_0}(\A)\cap g^\ddag K_0 g^{\ddag,-1}}\arrow[r]\arrow[d]& {\gamma P_{\Phi_0}\gamma^{-1}(\A)\cap \lcj{K_0}{\gamma g^\ddag g \gamma^{-1}}=\gamma P_{\Phi_0}\gamma^{-1}(\A)\cap \lcj{K_0}{\gamma g^{\ddag}}}\arrow[dd]\\
{P_{\Phi^\ddag}(\A)\cap g^{\ddag} K^\ddag g^{\ddag,-1}}\arrow[d]&  \\
{P_{\Phi^F}(\A)\cap g^\ddag K^F g^{\ddag,-1}}\arrow[r]& P_{(\Phi^F)'}(\A)\cap \lcj{K^F}{\tg_1^{-1}g^\ddag }.
\end{tikzcd}
\end{equation}
It follows from the assumption that $\Psi'=(\gamma Q\gamma,\gamma(X^+_0),\gamma g^\ddag g \gamma^{-1})$ on the top-right corner of (\ref{diag-bd-mix-sh}) is a cusp label representative in $\ca{CLR}(G_0,X_0)$ and that $K_{\Psi'}=\gamma P_{\Phi_0}\gamma^{-1}(\A)\cap \lcj{K_0}{\gamma g^\ddag g \gamma^{-1}}$. So $\sh_{K_{\Psi'}}$ is the boundary mixed Shimura variety associated with $\Psi'$ and $(G_0,X_0,K_0)$.\par
By \cite[Prop. 4.3.5]{Mad19}, the bottom arrow of the diagram (\ref{diag-bd-mix-sh}) above extends to $\mrm{int}(\tg_1^{-1})\circ g: \ca{S}_{K_{\Phi^F,p}}\to \ca{S}_{K_{(\Phi^F)',p}}$. Since there is an equivariant $(P_{\Phi^F}(\Ap)\to P_{(\Phi^F)'}(\Ap))$-action on this map, we have a map $\mrm{int}(\tg_1^{-1})\circ g: \ca{S}_{K_{\Phi^F}}\to \ca{S}_{K_{(\Phi^F)'}}$. This ends Construction \ref{cons-general-twisting}.\hfill $\square$
\end{construction}
Combining the construction above with \cite{Lan12b} and \cite[Sec. 2 and Sec. 3]{Lan16b}, we obtain the following statement:
\begin{lem}\label{lem-lan-qisog-twt}
The bottom arrow in (\ref{diag-bd-mix-sh}) can be viewed as a collection of morphisms between boundary mixed Shimura varieties induced by a quasi-isogeny twist in the sense of \cite[Sec. 2]{Lan16b}.
Let $\Q_{\Phi^F,\eta}$ (resp. $\Q_{(\Phi^F)',\eta}$) be the universal $1$-motive on $\sh_{K_{\Phi^F}}$ (resp. $\sh_{K_{(\Phi^F)'}}$); let $\Q_{\Phi^F}$ (resp. $\Q_{(\Phi^F)'}$) be its extension over $\ca{S}_{K_{\Phi^F}}$ (resp. $\ca{S}_{K_{(\Phi^F)'}}$), which is unique if it exists by Lemma \ref{lem-ext-1-mot}. We then have $\bb{Q}^\times$-isogenies between 
$\Q_{\Phi^F}$ and (the pullback to $\ca{S}_{K_{\Phi^F}}$ of) $\Q_{(\Phi^F)'}$, and between the corresponding components $\ul{X}_?$, $\ul{Y}_?$, $A_?$, $A^\vee_?$, $\G^\natural_?$, $\G^{\natural,\vee}_?$, $T_?$ and $T^\vee_?$ of $?=\Q_{\Phi^F}$ and $\Q_{(\Phi^F)'}$.\par
Furthermore, for any integer $n\geq 1$, this bottom arrow naturally induces \textbf{isomorphisms} over $\ca{S}_{K_{\Phi^F}}$:
$$\uhom(\frac{1}{n}X_{\Q_{\Phi^F}},A_{\Q_{\Phi^F}}^\vee)\iso\uhom(\frac{1}{n}X_{\Q_{(\Phi^F)'}},A_{\Q_{(\Phi^F)'}}^\vee),$$
$$\uhom(\frac{1}{n}Y_{\Q_{\Phi^F}},A_{\Q_{\Phi^F}})\iso\uhom(\frac{1}{n}Y_{\Q_{(\Phi^F)'}},A_{\Q_{(\Phi^F)'}}),$$
$$\uhom(\frac{1}{n}X_{\Q_{\Phi^F}},\G^{\natural,\vee}_{\Q_{\Phi^F}})\iso\uhom(\frac{1}{n}X_{\Q_{(\Phi^F)'}},\G^{\natural,\vee}_{\Q_{(\Phi^F)'}}),$$
$$\uhom(\frac{1}{n}Y_{\Q_{\Phi^F}},\G^{\natural}_{\Q_{\Phi^F}})\iso\uhom(\frac{1}{n}Y_{\Q_{(\Phi^F)'}},\G^{\natural}_{\Q_{(\Phi^F)'}}),$$
$$\uhom^{\mrm{symm}}(\frac{1}{n}X_{\Q_{\Phi^F}},\G^{\natural,\vee}_{\Q_{\Phi^F}})\iso\uhom^{\mrm{symm}}(\frac{1}{n}X_{\Q_{(\Phi^F)'}},\G^{\natural,\vee}_{\Q_{(\Phi^F)'}}),$$
and
$$\uhom^{\mrm{symm}}(\frac{1}{n}Y_{\Q_{\Phi^F}},\G^{\natural}_{\Q_{\Phi^F}})\iso\uhom^{\mrm{symm}}(\frac{1}{n}Y_{\Q_{(\Phi^F)'}},\G^{\natural}_{\Q_{(\Phi^F)'}}).$$
The objects on the right-hand sides are viewed as objects over $\ca{S}_{K_{\Phi^F}}$ by pulling back via $\mrm{int}(\tg_1^{-1})\circ g$.\par
Moreover, the isomorphisms above can be pulled back to $\ca{S}_{K_{\Phi_0}}$.
\end{lem}
\begin{proof}
The situation satisfies the assumptions in \cite[Sec. 2]{Lan16b} to which \cite[Sec. 3]{Lan16b} applies. In fact, the isogeny twists are induced by the difference of two symplectic lattices. In our case, they are $V_\bb{Z}\otimes\ca{O}_F$ and the pullback to $V_{\A}\otimes\ca{O}_F$ of $\V_{\bb{Z}}\otimes \ca{O}_F$. Then we have the first paragraph. For the second paragraph, the best way to see this is to use the description of the tower $\sh_{K_{\Phi^F}}\to \overline{\sh}_{K_{\Phi^F}}\to \sh_{K_{\Phi^F},h}$ at principal level $n$ in the characteristic zero theory (see, e.g., \cite[Lem. 3.5.5, Lem. 3.6.3, and Lem. 3.6.8]{Lan12b}) and the fact that $\mrm{int}(\tg_1^{-1})\circ g$ (or rather $\mrm{int}(\tg_1^{-1})$) induces an isomorphism between $D_{\Phi^F}$ and $D_{(\Phi^F)'}$ and an isomorphism between $P_{\Phi^F}(\bb{R})U_{\Phi^F}(\bb{C})$ and $P_{(\Phi^F)'}(\bb{R})U_{(\Phi^F)'}(\bb{C})$. Then the lemma follows from Lemma \ref{lem-ext-1-mot}.
\end{proof}
\newpage
\section{Toroidal compactifications of integral models of abelian type}\label{sec-tor-abelian}
In this {section}, our main goal is to construct good toroidal compactifications of integral models of $\sh_{K_2}(G_2,X_2)$, given some open compact subgroup $K_{2,p}$ of $G_2(\bb{Q}_p)$ and some neat open compact $K_2^p$. First, we generalize Deligne's induction method in \cite{Del79} to boundary mixed Shimura varieties of abelian type (see \S\ref{sec-del-ind-bd}). Then, in Cases (HS), ($\mrm{STB}_n$) and (DL), we construct toroidal compactifications for $\sh_{K}(G,X_b)$ (see \S\ref{subsec-con-stra-comp-quo}). Finally, we complete the proof of the main theorem in \S\ref{subsec-main-thm}.
\subsection{Deligne's induction for boundary components}\label{sec-del-ind-bd}
\subsubsection{}\label{sss-gp-gen}Let us recall some group-theoretic definitions following \cite[2.0.1]{Del79} and \cite[3.3.1]{Kis10}.
Let $H$ be a group. Suppose that there is an action of $\Delta$ on $H$, that is, a homomorphism $\Delta\to \Aut H$. Suppose that $\Gamma\sbst H$ is a $\Delta$-stable subgroup. Suppose that there is a $\Delta$-equivariant map $\varphi:\Gamma\to \Delta$, where $\Delta$ acts on itself by inner automorphisms. Suppose that, for any $\gamma\in \Gamma$, the action of $\varphi(\gamma)\in \Delta$ on $H$ is the same as the inner automorphism of $\gamma$ on $H$.\par
Then we denote $H*_\Gamma \Delta:=H\rtimes \Delta/(\gamma,\varphi(\gamma)^{-1})_{\gamma\in\Gamma}$; the group structure of $H*_\Gamma \Delta$ is inherited from that of $H\rtimes \Delta$.\par
There is a well-defined natural projection $H*_\Gamma\Delta\to \Gamma\bss H$, sending any representative $(h,d)\in H\times \Delta$ of an element in $H*_\Gamma \Delta$ to $\Gamma\bss H.$\par 
The following lemma is immediate from the definition above:
\begin{lem}\label{lem-agp-imm}
   We have an isomorphism of cosets $\Gamma*_{\Gamma} \Delta\bss H*_{\Gamma} \Delta\iso \Gamma\bss H$ via the well-defined natural projection $H*_\Gamma\Delta\to \Gamma\bss H$. Let $H_1$ be a $\Delta$-stable subgroup of $H$ containing $\Gamma$. Then there is an isomorphism of cosets $H_1*_\Gamma \Delta\bss H*_\Gamma \Delta\iso H_1\bss H$.
\end{lem}
\subsubsection{}\label{sss-del-ind-gen}
Let us recall the general setup of the so-called Deligne's induction following \cite[Sec. 2.7]{Del79} and \cite[3.3.5]{Kis10}. 
Let $O$ be a base scheme. Let $\Gamma$ be a locally profinite group. Let $\{S_\ca{H}\}_{\ca{H}}$ be an inverse system of quasi-projective $O$-schemes $S_\ca{H}$ labeled by a cofinal collection of open compact subgroups $\ca{H}\sbst \Gamma$. Let $S:=\varprojlim_\ca{H}S_\ca{H}$ be the inverse limit.\par
Suppose that there is a right action of $\Gamma$ given by $\gamma_{\ca{H}}:S_\ca{H}\xrightarrow{\sim} S_{\gamma^{-1}\ca{H}\gamma}$. Assume that 
\begin{itemize}
\item $\gamma_{\ca{H}}=\mrm{id}$ if $\gamma\in \ca{H}$. 
\item For any pair of open compact subgroups $\ca{H}_1\sbst \ca{H}\sbst \Gamma$ such that $\ca{H}_1$ is a normal subgroup of $\ca{H}$, the quotient $S_{\ca{H}_1}/(\ca{H}/\ca{H}_1)$ given by the action of $\gamma_{\ca{H}_1}$ for $\gamma\in \ca{H}$ induces an isomorphism $S_{\ca{H}_1}/(\ca{H}/\ca{H}_1)\xrightarrow{\sim}S_\ca{H}$; the transition map $S_{\ca{H}_1}\to S_\ca{H}$ determined by this action is finite and surjective. This implies that the limit $S$ exists as an $O$-scheme.
\end{itemize}
We call such an action of $\Gamma$ a \textbf{continuous (right) $\Gamma$-action on $S$}. One can also similarly define continuous left actions of $\Gamma$ by assigning to $\gamma\in \Gamma$ the inverse $\gamma^{-1}_{\ca{H}}$ for any $\ca{H}$ above.\par
Let $f:\Gamma\to \Gamma'$ be a continuous map of locally profinite groups with a compact kernel $\Delta(\Gamma,\Gamma')$. Define \textbf{Deligne's induction functor} as $\mrm{Ind}_\Gamma^{\Gamma'}S:=S\times \Gamma'/\Gamma$, where the action of $\Gamma$ on $S\times \Gamma'$ is defined by $(s,\gamma')\cdot \gamma=(s\cdot \gamma,\gamma^{-1}\gamma')$. For any compact (but \emph{not necessarily open}) subgroup $\ca{H}'\sbst \Gamma'$, \begin{equation}\label{eq-induction-basic}\begin{split}
&\mrm{Ind}_\Gamma^{\Gamma'}S/\ca{H}'\\
&\iso( S\times \Gamma'/\ca{H}')/\Gamma\\
&\iso \disju_{j\in J}[S\times (f(\Gamma) j \ca{H}'/\ca{H}')]/\Gamma\\
&\iso \disju_{j\in J}S/\ker(\Gamma\to \Gamma'/j\ca{H}' j^{-1}),
\end{split}\end{equation}
where $j$ runs over a complete set of representatives in $\Gamma'$ of $J:=f(\Gamma)\bss \Gamma'/\ca{H}'$.\par
Let $\wat{\pi}$ be a profinite set such that $\Gamma'$ acts on $\wat{\pi}$ continuously and transitively. Let $\Gamma_e'$ be the stabilizer of $e\in\wat{\pi}$ such that $\Gamma_e'\bss \Gamma'\xrightarrow{\sim}\wat{\pi}$ is a homeomorphism. Suppose that the map $\Gamma\to \Gamma'$ factors through a surjective map $\Gamma\twoheadrightarrow \Gamma_e'$. Then $S\times \Gamma'/\Gamma\iso (S/\Delta(\Gamma,\Gamma')\times \Gamma')/\Gamma'_e.$ Note that the quotient $S/\Delta(\Gamma,\Gamma')$ exists: Indeed, for any open compact subgroup $\ca{H}$ of $\Gamma$ such that the conjugation of $\Delta(\Gamma,\Gamma')$ stabilizes $\ca{H}$, the quotient $S_\ca{H}/\Delta(\Gamma,\Gamma')$ factors through a finite quotient since $\Delta(\Gamma',\Gamma)$ is compact.
\subsubsection{}\label{sss-del-ind-shi}
Let $(G,X)$ be any Shimura datum. Denote $Z:=Z_G$ and $\sh(G,X):=\varprojlim_{K\text{ open compact}}\sh_K(G,X)$. Let $$\ag(G):=\frac{G(\A)}{Z(\bb{Q})^{\overline{\ }}}*_{G(\bb{Q})_+/Z(\bb{Q})}G^\ad(\bb{Q})^+\iso
\frac{G(\A)}{Z(\bb{Q})^{\overline{\ }}}*_{G(\bb{Q})/Z(\bb{Q})}G^\ad(\bb{Q})_1.$$
Let 
$$\agsb(G):=\frac{G(\bb{Q})_+^{\overline{\ }}}{Z(\bb{Q})^{\overline{\ }}}*_{G(\bb{Q})_+/Z(\bb{Q})}G^\ad(\bb{Q})^+.$$
The ``${}^{\overline{\ }}$'' denotes the closure of corresponding groups of $\bb{Q}$-points in $G(\A)$.
By \cite[2.1.13-2.1.15]{Del79}, $\ag(G)$ acts transitively on $\wat{\pi}_0(G):=\varprojlim_{K}\pi_0(\sh_K(G,X)_{\overline{\bb{Q}}})\iso G(\A)/G(\bb{Q})_+^{\overline{\ }}$. The group $\wat{\pi}_0(G)$ is an abelian group. The stabilizer of this action is $\agsb(G)$.\par
Let $K$ be a neat open compact subgroup of $G(\A)$. For any $g\in G(\A)$ and fix a connected component $X^+$ of $X$, the connected component of $\sh_{K}(G,X)(\bb{C})$ represented by $g$ is $\Gamma_g\bss X^+$, where $\Gamma_g:=G(\bb{Q})_+\cap g K g^{-1}$ is a congruence subgroup. Moreover, $\Gamma_g$ acts on $X^+$ through its image in $G^\ad(\bb{Q})^+$. When $K$ runs over all open compact subgroup of $G(\A)$, the image of $\Gamma_g$ in $G^\ad(\bb{Q})^+$ forms a fundamental system of neighborhoods of the identity element $1\in G^\ad(\bb{Q})^+$, and defines a topology on $G^\ad(\bb{Q})^+$.\par
The left action of $G^\ad(\bb{Q})^+$ on $\sh^{+}(\bb{C}):=\varprojlim_K \sh_K^{+,(g)}(\bb{C})$, the inverse limit of $\sh_K^{+,(g)}(\bb{C}):=\Gamma_g\bss X^+$ for any $g$, induces an action of the completion $G^\ad(\bb{Q}){}^{+,\wat{\ }}$(rel. $G(\bb{Q})_+$) of $G^\ad(\bb{Q})^+$ with respect to the topology defined by congruence subgroups of $G(\bb{Q})_+$ as above. By \cite[2.0.13 and 2.1.15]{Del79}, $G^\ad(\bb{Q}){}^{+,\wat{\ }}$(rel. $G(\bb{Q})_+$) is canonically isomorphic to $G^\ad(\bb{Q}){}^{+,\wat{\ }}$(rel. $G^\der(\bb{Q})_+$), the completion of $G^\ad(\bb{Q}){}^{+}$ with respect to the topology defined by the congruence subgroups of $G^\der(\bb{Q})_+$. Moreover, we have $\agsb(G)\iso G^\ad(\bb{Q}){}^{+,\wat{\ }}\text{(rel. }G^\der(\bb{Q})_+\text{)}$. \par
Since $\sh_K^{+,(g)}(\bb{C})$ admits a canonical algebraization over $\bb{C}$ denoted by $\sh_{K,\bb{C}}^{+,(g)}$, we denote $\sh^{+}_{\bb{C}}:=\varprojlim_K \sh_{K,\bb{C}}^{+,(g)}$, and the action of $\agsb(G)$ canonically algebraizes to an action on $\sh^{+}_{\bb{C}}$. Moreover, the variety $\sh^{+}_{K,\bb{C}}$ descends to a variety $\sh^{+,(g)}_{K,\qbar}$ over $\qbar$, and the action of $\agsb(G)$ is also defined over $\qbar$. Denote $\sh^{+}_{\qbar}:=\varprojlim_K \sh_{\qbar}^{+,(g)}$.\par
Hence, we have that (see \cite[2.1.15]{Del79})
$$\sh(G,X)_{\qbar}\iso (\sh_{\qbar}^+\times \ag(G))/\agsb(G).$$
Let $(G_2,X_2)$ be a Shimura datum such that there is a central isogeny $G^\der\to G^\der_2$ that induces an isomorphism $(G^\ad,X^\ad)\iso(G_2^\ad,X_2^\ad)$. As any pullback of a congruence subgroup of $G_2^\der(\bb{Q})$ contains a congruence subgroup of $G^\der(\bb{Q})$, we see that there is a natural continuous surjective map 
$$\agsb(G)\lra \agsb(G_2).$$
Let $\Delta^\circ(G,G_2)$ be the kernel of the map above. Following almost the same proof as \cite[Lem. E. 6]{Kis17}, the group $\Delta^\circ(G,G_2)$ is compact: To see this, it suffices to replace the Galois cohomology input there with the fact that 
$\ker (H^1(\bb{Q},Z_G)\to \prod H^1(\bb{Q}_l,Z_G))$ is finite. Hence, the kernel of $\ag(G)\to \ag(G^\ad)$ contains $Z_G(\A)/Z_G(\bb{Q})^{\overline{\ }}$ as a finite index subgroup.\par
From \cite[2.7.11]{Del79} (see also \cite[Prop. 3.3.10]{Kis10} and \ref{sss-del-ind-gen}), we have that 
$$\sh(G_2,X_2)_{\qbar}\iso (\sh_{\qbar}^+\times \ag(G_2))/\agsb(G).$$
\subsubsection{}\label{sss-agp-mixed}
Let $(G,X)$ be any Shimura datum. Let $\Phi=(Q,X^+,g)$ be any cusp label representative. Let $K$ be a neat open compact subgroup of $G(\A)$.\par 
Let $Y_\Phi$ be any connected normal $\bb{Q}$-subgroup of $Q$ such that $P_\Phi\sbst Y_\Phi\sbst ZP_\Phi$. Let $\sh_{K_\Phi^Y}(\bb{C}):=\sh_{K_\Phi^Y}(Y_\Phi(\A),Y_\Phi(\bb{Q})D_\Phi)(\bb{C})$, where $K_\Phi^Y:=Y_\Phi(\A)\cap g K g^{-1}$.\par
Fix a connected component $D_\Phi^+\sbst D_\Phi$, then this choice determines an isomorphism between the group $\pi_0(Y_\Phi)_{K}$ of connected components of the mixed Shimura variety $\sh_{K_\Phi^Y}(\bb{C})$ and $Y_\Phi(\bb{Q})_+\bss Y_\Phi(\A)/K_\Phi$, where $Y_\Phi(\bb{Q})_+:=\stb_{Y_\Phi(\bb{Q})}(D_\Phi^+)$. Then $\sh_{K_\Phi^Y}(\bb{C})$ is isomorphic to a disjoint union 
$$\disju_{p_f\in Y_\Phi(\bb{Q})_+\bss Y_\Phi(\A)/K_\Phi} \Gamma_{p_f}\bss D_\Phi^+.$$
In the equation above, $\Gamma_{p_f}:=Y_\Phi(\bb{Q})_+\cap p_f K_\Phi^Y p^{-1}_f.$\par
Fix a $p_f\in Y_\Phi(\A)$. Denote $\sh_{K_\Phi^Y}^{+,(p_f)}(\bb{C}):=\Gamma_{p_f}\bss D_\Phi^+$. By \cite[Ch.9]{Pin89}, $\sh_{K_\Phi^Y}^{+,(p_f)}(\bb{C})$ uniquely algebraizes to an algebraic variety $\sh_{K_\Phi^Y,\bb{C}}^{+,(p_f)}$; by \cite[Ch.11]{Pin89}, this algebraic variety descends to an algebraic variety $\sh_{K_\Phi^Y,\qbar}^{+,(p_f)}$ over $\qbar$. Let $\sh_{Y_\Phi,\qbar}^{+}:=\varprojlim_K \sh_{K_\Phi^Y,\qbar}^{+,(p_f)}$, which does not depend on $p_f$. Denote $\sh_{Y_\Phi}:=\varprojlim_{K\text{ open compact}}\sh_{K_\Phi^Y}\iso \varprojlim_{K_\Phi^Y\text{ open compact in }Y_\Phi(\A)}\sh_{K_\Phi^Y}$.\par
We set 
$$\ag(Y_\Phi):=\frac{Y_\Phi(\A)}{Z^Y_{\Phi}(\bb{Q})^{\overline{\ }}}*_{Y_{\Phi}(\bb{Q})_+/Z^Y_{\Phi}(\bb{Q})}Y_{\Phi}^{\ad}(\bb{Q})^+,$$
and set 
$$\agsb(Y_\Phi):=\frac{Y_{\Phi}(\bb{Q})_+^{\overline{\ }}}{Z^Y_{\Phi}(\bb{Q})^{\overline{\ }}}*_{Y_{\Phi}(\bb{Q})_+/Z^Y_{\Phi}(\bb{Q})}Y_{\Phi}^\ad(\bb{Q})^+.$$
In the two equations above, $Z^Y_\Phi:=Z_G\cap Y_\Phi$ and $Y_\Phi^\ad:=Y_\Phi/Z^Y_\Phi$. So $P^\ad_\Phi=Y_\Phi^\ad=ZP_\Phi^\ad$. The ``$\overline{\ }$'' denotes the closure in $Y_\Phi(\A)$. The group $\ag(Y_\Phi)$ acts on $\wat{\pi}_0(Y_\Phi):=\varprojlim \pi_0(Y_\Phi)_K$ transitively with stabilizer $\agsb(Y_\Phi)$.\par 
For any closed connected unipotent normal subgroup $W\sbst Y_\Phi$, we can also define
$$\ag(Y_\Phi/W):=\frac{Y_\Phi(\A)/W(\A)}{Z^Y_{\Phi}(\bb{Q})^{\overline{\ }}}*_{(Y_{\Phi}(\bb{Q})_+/W(\bb{Q}))/Z^Y_{\Phi}(\bb{Q})}Y_{\Phi}^{\ad}(\bb{Q})^+/W(\bb{Q}),$$
and define 
$$\agsb(Y_\Phi/W):=\frac{Y_{\Phi}(\bb{Q})_+^{\overline{\ }}/W(\bb{Q})^{\overline{\ }}}{Z^Y_{\Phi}(\bb{Q})^{\overline{\ }}}*_{(Y_{\Phi}(\bb{Q})_+/W(\bb{Q}))/Z^Y_{\Phi}(\bb{Q})}Y_{\Phi}^\ad(\bb{Q})^+/W(\bb{Q}).$$
\begin{lem}\label{lem-agy-basic}
With the notation above, 
\begin{enumerate}
\item $\ag(Y_\Phi)$ and $\agsb(Y_\Phi)$ fit into the following exact sequences induced by Levi quotients:
\begin{equation}
\begin{split}
0\lra W_\Phi(\A)\lra \ag(Y_\Phi)\lra \ag(Y_\Phi/W_\Phi)\lra 0 ,   \\
0\lra W_\Phi(\A)\lra \agsb(Y_\Phi)\lra \agsb(Y_\Phi/W_\Phi)\lra 0.
\end{split}
\end{equation}
\item $\agsb(Y_\Phi)$ is the completion of $Y_\Phi^\ad(\bb{Q})^+$ with respect to the topology generated by the images of congruence subgroups of $Y_\Phi(\bb{Q})_+$ in $Y_\Phi^\ad(\bb{Q})^+.$
\item Let $Y'_\Phi$ be the smallest normal $\bb{Q}$-subgroup of $Y_\Phi$ such that there is a product $Y'_\Phi\cdot Z_\Phi^Y\iso Y_\Phi$ and such that $Y_\Phi'$ and $Z^Y_\phi$ have finite intersection. Then the definition of $\agsb(Y_\Phi)$ depends only on $Y'_\Phi$.
\end{enumerate}
\end{lem}
\begin{proof}
By \cite[Lem. 2.7 and Prop. 6.6]{PR94}, $W_\Phi(\bb{A})\iso W_\Phi(\bb{Q})^{\overline{\ }}$. Since $W_\Phi(\A)*_{W_\Phi(\bb{Q})}W_\Phi(\bb{Q})\iso W_\Phi(\A)$, the first part follows. For any open compact subgroup $K\sbst Y_\Phi(\A)$, by \cite[Prop. 6.5, p.296]{PR94}, the image of $K$ in the quotient $Y_\Phi/W_\Phi(\A)$ is open and compact. Combining this fact with strong approximation of unipotent groups, and choosing a Levi decomposition representing $Y_\Phi$ as a semidirect product of $W_\Phi$ and $Y_\Phi/W_\Phi$, we see that the family of congruence subgroups of $Y_\Phi(\bb{Q})$ contains a cofinal system consisting of the products of congruence subgroups of $W_\Phi(\bb{Q})$ and $Y_\Phi(\bb{Q})/W_\Phi(\bb{Q})$. Then the second part follows from \cite[2.0.9]{Del79} and the first part. To show the third part, we can assume that $Y_\Phi$ is reductive by the first part. Then the third part is the consequence of the following Lemma \ref{lem-del-cong-gen} that can be directly derived from \cite[2.0.13]{Del79}.
\end{proof}
\begin{lem}\label{lem-del-cong-gen}
    Let $Y$ be a reductive $\bb{Q}$-group. Let $Y_1$ be a $\bb{Q}$-subgroup containing $Y^\der$. Let $Z_1$ be a subgroup of $Z_Y$ such that $Z_1$ and $Y_1$ have finite intersection and such that $Z_1\cdot Y_1\iso Y$. Then the product of a congruence subgroup in $Y_1$ and a congruence subgroup in $Z_1$ is a congruence subgroup of $Y$. 
\end{lem}
\begin{proof}
By \cite[2.0.13]{Del79}, any congruence subgroup $U_1\sbst Y_1(\bb{Q})$ contains a finite index subgroup $U'$, where $U'$ is a product of a congruence subgroup $U^\der$ of $Y^\der$ and a congruence subgroup $U_Z^1$ of $Z^\circ_{Y_1}(\bb{Q})$, which is a finite index subgroup of group of units of $Z^\circ_{Y_1}(\bb{Q})$ by \cite[2.0.10]{Del79}. By \emph{loc. cit.} again, the product of the congruence subgroup $U_Z^1$ and a congruence subgroup $U_Z\sbst Z_1$ is a congruence subgroup in $Z_1\cdot Z_{Y_1}(\bb{Q})\iso Z_Y(\bb{Q})$. Hence, $U'\cdot U_Z$ is a congruence subgroup of $Y(\bb{Q})$ by \cite[2.0.13]{Del79} again, and so is $U_1\cdot U_Z$.
\end{proof}
\subsubsection{}
By part 2 of Lemma \ref{lem-agy-basic} and the paragraphs in \S\ref{sss-agp-mixed} before that, we have an $\ag(Y_\Phi)$-equivariant isomorphism
$$\sh_{Y_\Phi,\qbar}\iso \sh^+_{Y_\Phi,\qbar}\times \ag(Y_\Phi)/\agsb(Y_\Phi).$$
Let $E'$ be a finite extension of $E:=E(G,X)=E(P_\Phi,D_\Phi)$ (see \cite[Prop. 12.1]{Pin89}).\par
The action of $\gal(E^{\prime,ab}/E^\prime)$ on $\wat{\pi}_0(Y_{\Phi})\iso Y_{\Phi}(\A)/Y_{\Phi}(\bb{Q})_+^{\overline{\ }}\iso (Y_{\Phi}/W_{\Phi})(\A)/(Y_{\Phi}/W_{\Phi})(\bb{Q})_+^{\overline{\ }}$ is given by
\begin{equation}\label{eq-gal-action-yphih}
\begin{tikzcd}
r_{E'}(Y_{{\Phi}}):\gal(E^{\prime,ab}/E^\prime)\arrow[r,"\mrm{rec_{E^\prime}}"]&\pi_0(\bb{G}_m(\bb{A}_{E^\prime})/\bb{G}_m(E^\prime))\arrow[r,"\sim"]&{\frac{\bb{A}_{E'}^\times}{E^{\prime,\times}(E'\otimes_\bb{Q}\bb{R})^{\times,+}}}\\
\iso \bb{A}_{E',f}^\times/E^{\prime,\times,+}
\arrow[r,"{[\mu_{{\Phi},x}]}"]& Y_{\Phi}(\bb{A}_{E',f})/Y_{\Phi}(E')_+^{\overline{\ }}\arrow[r,"N_{E'/\bb{Q}}"]&\frac{Y_{\Phi}(\A)}{Y_{\Phi}(\bb{Q})_+^{\overline{\ }}}.
\end{tikzcd}
\end{equation}
In the diagram (\ref{eq-gal-action-yphih}) above, $[\mu_{\Phi,x}]$ is the morphism determined by the conjugacy class of Hodge cocharacter $\mu_{\Phi,x}$ for any $x\in Y_{\Phi}(\bb{Q})D_{\Phi}$ over $E'$; $N_{E'/\bb{Q}}$ is the map determined by norm map $N_{E'/\bb{Q}}:R_{E'/\bb{Q}}Y_{\Phi,E'}\to Y_{\Phi}$. 
Let $\wdtd{r}_{E'}(Y_\Phi)$ be the composition of $r_{E'}(Y_\Phi)$ with the natural map $\gal(\overline{E}'/E')\to \gal(E^{\prime,ab}/E')$.\par 
Let $\eg(Y_\Phi)$ be the pullback of $$0\lra \agsb(Y_\Phi)\lra \ag(Y_\Phi)\lra \agsb(Y_\Phi)\bss \ag(Y_\Phi)\lra 0$$ via $-\wdtd{r}_{E'}(Y_\Phi)$. Then $\eg(Y_\Phi)$ is the stabilizer of $\sh_{Y_\Phi,\qbar}^+$ in $\ag(Y_\Phi)\times \gal(\overline{E}'/E')$.\par
So there is an action of $\ag(Y_\Phi)\times \gal(\overline{E}'/E')$ on $\sh_{Y_\Phi,\qbar}$ such that the natural projection $\wat{\pi}_0:\sh_{Y_\Phi,\qbar}\to \wat{\pi}_0(Y_\Phi)$ from $\sh_{Y_\Phi,\qbar}$ to its group of connected components is $\ag(Y_\Phi)\times \gal(\overline{E}'/E')$-equivariant.
\subsubsection{}\label{sss-del-ind-mix}
Let $\Phi=(Q,X^+,g)$ be a cusp label representative of $(G,X)$. The cusp label representative $\Phi$ maps to a cusp label representative $\Phi^\ad=(Q^\ad,X^{\ad,+},g)$ of $(G^\ad,X^\ad)$, where $Q^\ad=Q/Z_G$, $X^{\ad,+}$ is the image of $X^+$ and $g$ here abusively denotes the image of $g$.
\begin{lem}\label{lem-dq-quotient}The map $D_{Q,X^+}\to D_{Q^\ad,X^{\ad,+}}$ induced by $\Phi\to \Phi^\ad$ is injective and is a homeomorphism on each connected component of $D_{Q,X^+}$.
\end{lem}
\begin{proof}
We have a $P_{Q}(\bb{R})\cdot U_{Q}(\bb{C})$-equivariant diagram
    \begin{equation*}
        \begin{tikzcd}
            X\arrow[r]\arrow[d]&\pi_0(X)\times\Hom(\bb{S}_\bb{C}, P_{Q,\bb{C}})\arrow[d]\\
            X^\ad\arrow[r]& \pi_0(X^\ad)\times\Hom(\bb{S}_\bb{C},P_{Q^\ad,\bb{C}}).
        \end{tikzcd}
    \end{equation*}
By \cite[Prop. 2.9]{Pin89}, we can write $D_{Q,X^+}$ as a quotient $P_Q(\bb{R})U_Q(\bb{C})/\stb(x)$ where $\stb(x)$ is the stabilizer of a point in $D_{Q,X^+}$. The map $D_{Q,X^+}\to D_{Q^\ad,X^{\ad,+}}$ is induced by the quotient of $P_Q(\bb{R})U_Q(\bb{C})$ by $(P_Q\cap Z_G)(\bb{R})$, and the latter group is contained in $\stb(x)$.
\end{proof}
Let $(G_2,X_2)$ be the Shimura datum as in \S\ref{sss-del-ind-shi}. 
Let $\Phi$ (resp. $\Phi_2$) be a cusp label representative of $(G,X)$ (resp. $(G_2,X_2)$). Recall that $Y'_\Phi$ is the smallest normal $\bb{Q}$-subgroup of $Y_\Phi$ such that $Y'_\Phi$ and $Z^Y_\Phi$ have finite intersection and such that $Y'_\Phi\cdot Z_\Phi^Y\iso Y_\Phi$.\par 
We assume that there is an isogeny $Y'_\Phi\to Y'_{\Phi_2}$ whose kernel is in $Z_G\cap Y'_\Phi$. \par
From Part 2 and Part 3 of Lemma \ref{lem-agy-basic}, we see that there is a natural continuous and surjective map
$$\agsb(P_\Phi)\to \agsb(Y_{\Phi_2}).$$
\begin{lem}\label{lem-ker-cpt}
The kernel $\Delta^\circ(P_\Phi,Y_{\Phi_2}):=\ker(\agsb(P_\Phi)\to\agsb(Y_{\Phi_2}))$ is compact. The intersection $(Z_G\cap P_\Phi)(\A)/(Z_G\cap P_\Phi)(\bb{Q})^{\overline{\ }}\cap \Delta^\circ (P_\Phi,Y_{\Phi_2})$ is a finite index subgroup in $\Delta^\circ(P_\Phi,Y_{\Phi_2})$. 
\end{lem}
\begin{proof}
    Combine Part 1 of Lemma \ref{lem-agy-basic} and \cite[Lem. E. 6]{Kis17}, with the Galois cohomology input there changed to the fact that
$\ker (H^1(\bb{Q},Z_G\cap P_\Phi)\lra \prod H^1(\bb{Q}_l,Z_G\cap P_\Phi))$ is finite.
\end{proof}
Note that $\sh_{Y_{\Phi_2},\qbar}$ is equipped with a Hecke action of $Y_{\Phi_2}(\A)$ on the right (see \cite[11.5]{Pin89}). From the arguments above, we see that
\begin{prop}\label{prop-del-ind-mix}There is an $\ag(Y_{\Phi_2})$-equivariant isomorphism  
$$\sh_{Y_{\Phi_2},\qbar}\iso \sh^+_{P_\Phi,\qbar}\times \ag(Y_{\Phi_2})/\agsb(P_\Phi).$$
The $\ag(Y_{\Phi_2})$-action of the right-hand side is induced by the right multiplication on its second factor.
\end{prop}
From the Proposition above, we see that $\eg(Y_{\Phi_2})$ is the pushout of $0\to \agsb(P_{\Phi})\to \eg(P_{\Phi})\to \gal(\overline{E}'/E')\to 0$ under the natural map $\agsb(P_{\Phi})\to \agsb(Y_{\Phi_2})$ given by Part 3 of Lemma \ref{lem-agy-basic}.
So $\ag(Y_{\Phi_2})*_{\agsb(P_{\Phi})}\eg(P_{\Phi})\iso \ag(Y_{\Phi_2})\times \gal(\overline{E}'/E')$. \par
So there is an action of $\ag(Y_{\Phi_2})\times \gal(\overline{E}'/E')$ on $\sh_{Y_{\Phi_2},\overline{\bb{Q}}}$ such that the natural projection $\wat{\pi}_0:\sh_{Y_{\Phi_2},\overline{\bb{Q}}}\to \wat{\pi}_0(Y_{\Phi_2})$ from $\sh_{Y_{\Phi_2},\overline{\bb{Q}}}$ to its group of connected components is $\ag(Y_{\Phi_2})\times \gal(\overline{E}'/E')$-equivariant.\par
Recall that any $g_2\in G_2(\bb{Q})$ induces an isomorphism $\sh_{Y_{\Phi_2},\qbar}\xrightarrow{\sim}\sh_{\lcj{Y_{\Phi_2}}{g_2},\qbar}$, where $\sh_{\lcj{Y_{\Phi_2}}{g_2},\qbar}$ is the inverse limit of a tower of mixed Shimura varieties $\sh_{\lcj{K^Y_{\Phi_2}}{g_2}}(\lcj{Y_{\Phi_2}}{g_2}(\A),\lcj{Y_{\Phi_2}(\bb{Q})}{g_2}g_2\cdot D_{\Phi_2})_{\qbar}$. Also, there is a right action of $h_2\in G_2(\A)$, mapping $\sh_{K_{\Phi_2}^Y}(Y_{\Phi_2},Y_{\Phi_2}(\bb{Q})D_{\Phi_2})_{\qbar}$ to $\sh_{\rcj{(K_{\Phi_2}^Y)}{h_2}}(Y_{\Phi_2},Y_{\Phi_2}(\bb{Q})D_{\Phi_2})_{\qbar}$. Therefore, there is a right action $ h_2: \sh_{Y_{\Phi_2},\qbar}\xrightarrow{\sim} \sh_{Y_{\Phi_2},\qbar}$ on the inverse limit $\sh_{Y_{\Phi_2},\qbar}$. 
\subsubsection{}
As in \cite{Kis10} and \cite{KP15}, we want to decompose $\ag(Y_{\Phi_2})$ to smaller groups.\par
Let us define more ``$\ag$-type'' and ``$\agsb$-type'' groups.\par
Fix any general $(G,X)$ and $Y_\Phi$. Let $K:=K_{p}K^p$ be a neat open compact subgroup of $G(\A)$ such that $K^p\sbst G(\Ap)$ is neat open compact. We define 
$$\atd_{K_p}(Y_\Phi):=\frac{Y_{\Phi}(\Ap)\times K_{\Phi,p}^Y}{Z_{\Phi,K_p}^Y(\zbkp)^{\overline{\ }}}*_{Y_{\Phi,K_p}(\zbkp)_+/Z_{\Phi,K_p}^Y(\zbkp)}P_{\Phi,K_p}^{\ad}(\zbkp)^+.$$
In the equation above, $Y_{\Phi,K_p}(\zbkp)_+:=Y_\Phi(\bb{Q})_+\cap \K_{\Phi,p}^Y$, $Z_{\Phi,K_p}^Y(\zbkp):=Z_G(\bb{Q})\cap K_{\Phi,p}^Y$, and $P^{\ad}_{\Phi,K_p}(\zbkp)^+:=(Y_\Phi^\ad)(\bb{Q})^+\cap\im K_{\Phi,p}^Y=(ZP_\Phi/Z_G)(\bb{Q})^+\cap \im K_{\Phi,p}^Y$, where $\im K_{\Phi,p}^Y$ is the image of $K_{\Phi,p}^Y$ in $G^\ad(\bb{Q}_p)$.\par
We define
$$\atdsb_{K_p}(Y_\Phi):=\frac{Y_{\Phi,K_p}(\zbkp)_+^{\overline{\ }}}{Z^Y_{\Phi,K_p}(\zbkp)^{\overline{\ }}}*_{Y_{\Phi,K_p}(\zbkp)_+/Z^Y_{\Phi,K_p}(\zbkp)}P_{\Phi,K_p}^{\ad}(\zbkp)^+.$$
Therefore, $\atdsb(Y_\Phi)$ is the closure of $P_{\Phi,K_p}^\ad(\zbkp)^+$ in $\agsb(Y_\Phi)$ (with respect to the topology on $\agsb(Y_\Phi)$ defined in \S\ref{sss-agp-mixed}).
Moreover, let
$$\agsb_{K_p}(Y_\Phi):=\frac{Y_{\Phi,K_p}(\zbkp)_+^{\overline{p}}}{Z^Y_{\Phi,K_p}(\zbkp)^{\overline{p }}}*_{Y_{\Phi,K_p}(\zbkp)_+/Z^Y_{\Phi,K_p}(\zbkp)}P_{\Phi,K_p}^{\ad}(\zbkp)^+.$$
In the definition above, ``$\overline{p}$'' denotes the closure in $Y_\Phi(\Ap)$.
Then $\agsb_{K_p}(Y_\Phi)$ is the completion of $P_{K_p}^\ad(\zbkp)^+$ with respect to the topology generated by the image of the sets of the form $Y'_{K_p}(\zbkp)_+\cap K^p$ in $P_\Phi^\ad(\bb{Q})^+$ by Lemma \ref{lem-agy-basic} and \cite[2.0.12]{Del79}.\par
Define 
$$\ag_{K_p}(Y_\Phi):=\frac{Y_\Phi(\Ap)}{Z^Y_{\Phi,K_p}(\zbkp)^{\overline{p }}}*_{Y_{\Phi,K_p}(\zbkp)_+/Z^Y_{\Phi,K_p}(\zbkp)}P_{\Phi,K_p}^{\ad}(\zbkp)^+.$$

Now, let $(G,X)$, $(G_2,X_2)$, $\Phi$, $\Phi_2$, $P_\Phi$ and $Y_{\Phi_2}$ be the symbols defined as in \S\ref{sss-del-ind-mix}. Let $K_2=K_{2,p}K^p_2$ be a neat open compact subgroup of $G_2(\A)$ with $K_2^p$ neat open compact. So we have defined $\atdsb_{K_{2,p}}(Y_{\Phi_2})$ and $\agsb_{K_{2,p}}(Y_{\Phi_2})$ as above.\par
Define $\atdsb_{K_{2,p}}(P_\Phi)$ as the closure of $P^\ad_{\Phi_2,K_{2,p}}(\zbkp)^+$ in $\agsb(P_\Phi)$. 
Since there is a central isogeny $P'_\Phi\to Y'_{\Phi_2}$, we can choose a neat open compact subgroup $K_p\sbst G(\bb{Q}_p)$ such that $g_\Phi K_pg_\Phi^{-1}\cap P'_\Phi(\bb{Q}_p)$ is a normal subgroup contained in the pullback of $g_{\Phi_2}K_{2,p}g_{\Phi_2}^{-1}\cap Y'_{\Phi_2}(\bb{Q}_p)$ under this central isogeny.\par
Define $\agsb_{K_{2,p}}(P_\Phi)$ to be the completion of $P^\ad_{\Phi_2,K_{2,p}}(\zbkp)^+$ with respect to the topology 
generated by the image of the sets of the form $P'_{K_p}(\zbkp)_+\cap K^p$ in $P_{\Phi_2}^\ad(\bb{Q})^+$.\par
Explicitly, we can write 
$$\atdsb_{K_{2,p}}(P_\Phi):=\frac{P_{\Phi,K_p}(\zbkp)_+^{\overline{\ }}}{Z_{\Phi,K_p}(\zbkp)^{\overline{\ }}}*_{P_{\Phi,K_p}(\zbkp)_+/Z_{\Phi,K_p}(\zbkp)}P_{\Phi_2,K_{2,p}}^{\ad}(\zbkp)^+,$$
and write
$$\agsb_{K_{2,p}}(P_\Phi):=\frac{P_{\Phi,K_p}(\zbkp)_+^{\overline{p}}}{Z_{\Phi,K_p}(\zbkp)^{\overline{p }}}*_{P_{\Phi,K_p}(\zbkp)_+/Z_{\Phi,K_p}(\zbkp)}P_{\Phi_2,K_{2,p}}^{\ad}(\zbkp)^+.$$
The maps $\atdsb_{K_{2,p}}(P_\Phi)\to \atdsb_{K_{2,p}}(Y_{\Phi_2})$ and $\agsb_{K_{2,p}}(P_\Phi)\to \agsb_{K_{2,p}}(Y_{\Phi_2})$ are surjective; their kernels are compact by Lemma \ref{lem-ker-cpt}.
\subsubsection{}
For a general $(G,X)$,
define $$\jg_{K^Y_{\Phi,p}}(Y_\Phi):=Y_{\Phi}(\bb{Q})_+\bss Y_\Phi(\bb{Q}_p)/K_{\Phi,p}^Y.$$
We simplify this notation to $\jg_{K_{p}}(Y_\Phi)$ if we know how to form the group $K^Y_{\Phi,p}$ from some $K_p\sbst G(\bb{Q}_p)$.
Then $\jg_{K_p}(Y_\Phi)\iso (Y_\Phi/W_\Phi(\bb{Q})_+)\bss (Y_\Phi/W_\Phi(\bb{Q}_p))/K_{\Phi,h,p}^Y$ is an abelian group by \cite[2.1.3]{Del79}, where $K_{\Phi,h,p}^Y$ is the image of $K^Y_{\Phi,p}$ in $Y_\Phi/W_\Phi(\bb{Q}_p)$. If we assume that $K_p$ is hyperspecial, $K^Y_{\Phi,h,p}$ is also hyperspecial. By \cite[Lem. 2.2.6]{Kis10}, $\jg_{K_p}(Y_\Phi)$ is trivial. 
\begin{lem}\label{lem-jgp}
With the previous definitions, there is a natural injective map $$\atdsb_{K_p}(Y_\Phi)\bss \atd_{K_p}(Y_\Phi)/K^Y_{\Phi,p}\hookrightarrow \agsb(Y_\Phi)\bss \ag(Y_\Phi)/K^Y_{\Phi,p},$$
where $\atdsb_{K_p}(Y_\Phi)\bss \atd_{K_p}(Y_\Phi)/K^Y_{\Phi,p}\iso \agsb_{K_p}(Y_\Phi)\bss \ag_{K_p}(Y_\Phi)$. The index of the left-hand side in the right-hand side is isomorphic to $\jg_{K_p}(Y_\Phi)$.\par
Therefore, we can write $\ag(Y_\Phi)$ as a fiber bundle over $\jg_{K_p}(Y_\Phi)$; at each $j\in \jg_{K_p}(Y_\Phi)$, the fiber is ${\agsb_{\lcj{K_p}{j}}(Y_\Phi)}\bss {\ag_{\lcj{K_p}{j}}(Y_\Phi)}\iso \agsb_{K_p}(Y_\Phi)\bss\ag_{K_p}(Y_\Phi).$ When $K_p$ is hyperspecial, we have that $\agsb_{K_p}(Y_\Phi)\bss\ag_{K_p}(Y_\Phi)\iso \agsb(Y_\Phi)\bss \ag(Y_\Phi)/K^Y_{\Phi,p}$.
\end{lem}
\begin{proof}
By Lemma \ref{lem-agp-imm}, $\agsb(Y_\Phi)\bss \ag(Y_\Phi)\iso Y_\Phi(\bb{Q})_+^{\overline{\ }}\bss Y_\Phi(\A)$, $\atdsb_{K_p}(Y_\Phi)\bss \atd_{K_p}(Y_\Phi)\iso Y_{\Phi,K_p}(\zbkp)_+^{\overline{\ }}\bss Y_\Phi(\Ap)\times K_{\Phi,p}^Y$ and $\agsb_{K_p}(Y_\Phi)\bss \ag_{K_p}(Y_\Phi)\iso Y_{\Phi,K_p}(\zbkp)_+^{\overline{p}}\bss Y_\Phi(\Ap)$. 
Moreover, 
\begin{equation*}
\begin{split}
&Y_{\Phi,K_p}(\zbkp)_+^{\overline{p}}\bss Y_\Phi(\Ap)\\
&\iso\varprojlim_{K^p}Y_{\Phi,K_p}(\zbkp)_+\bss Y_\Phi(\Ap)/K^{Y,p}_\Phi\\
&\iso\varprojlim_{K^p}Y_{\Phi,K_p}(\zbkp)_+\bss Y_\Phi(\Ap)\times K^Y_{\Phi,p}/K^Y_{\Phi,p}K^{Y,p}_\Phi\\
&\iso \varprojlim_{K_p'K^p,\text{ }K_p'\sbst K_p} \{[Y_{\Phi,K_p}(\zbkp)_+\bss Y_\Phi(\Ap)\times K^{Y}_{\Phi,p}/K^{\prime,Y}_{\Phi,p}K^{Y,p}_\Phi]/K^Y_{\Phi,p}\}\\
&\iso \{\varprojlim_{K_p'K^p,\text{ }K_p'\sbst K_p} [Y_{\Phi,K_p}(\zbkp)_+\bss Y_\Phi(\Ap)\times K^{Y}_{\Phi,p}/K^{\prime,Y}_{\Phi,p}K^{Y,p}_\Phi]\}/K^Y_{\Phi,p}\\
&\iso Y_{\Phi,K_p}(\zbkp)_+^{\overline{\ }}\bss Y_\Phi(\Ap)\times K^Y_{\Phi,p}/K^Y_{\Phi,p}. 
\end{split}
\end{equation*}
The fourth line to the fifth line is by the Mittag-Leffler criterion.
So the lemma follows.
\end{proof}
The ``$\eg$-type'' groups also decompose over $\jg_{K_p}(Y_\Phi)$.\par
By weak approximation for $\bb{G}_m$, the map defined by the quotient of $r_{E'}(Y_\Phi)$ by $K_{\Phi,p}^Y$, namely,
\begin{equation}\label{eq-gal-kyp}
    \begin{tikzcd}
r_{E',K_p}(Y_\Phi):\gal(E^{\prime,ab}/E^\prime)\arrow[r,"\mrm{rec_{E^\prime}}"]&\pi_0(\bb{G}_m(\bb{A}_{E^\prime})/\bb{G}_m(E^\prime))\arrow[r,"\sim"]&\frac{\bb{A}_{E'}^\times}{E^{\prime,\times}(E'\otimes_\bb{Q}\bb{R})^{\times,+}}\\
\iso \bb{A}_{E',f}^\times/E^{\prime,\times,+}
\arrow[r,"{[\mu_{\Phi,x}]}"]& Y_{\Phi}(\bb{A}_{E',f})/Y_{\Phi}(E')_+^{\overline{\ }}\arrow[r,"N_{E'/\bb{Q}}/K^Y_{\Phi,p}"]&\frac{Y_{\Phi}(\A)}{Y_{\Phi}(\bb{Q})_+^{\overline{\ }}K^Y_{\Phi,p}},
    \end{tikzcd}
\end{equation}
factors through $$Y_{\Phi}(\Ap)/Y_{\Phi,K_p}(\zbkp)_+^{\overline{p}}\iso Y_{\Phi}(\Ap)\times K^Y_{\Phi,p}/Y_{\Phi,K_p}(\zbkp)_+^{\overline{\ }}K^Y_{\Phi,p}\sbst \frac{Y_{\Phi}(\A)}{Y_{\Phi}(\bb{Q})_+^{\overline{\ }}K^Y_{\Phi,p}}.$$\par
We abusively denote the resulting map $\gal(E^{\prime,ab}/E')\to Y_\Phi(\Ap)/Y_{\Phi,K_p}(\zbkp)_+^{\overline{p}}$ by $r_{E',K_p}(Y_\Phi)$. Let $\wdtd{r}_{E',K_p}(Y_\Phi)$ be the composition of $r_{E',K_p}(Y_\Phi)$ with the natural map $\gal(\overline{E}'/E')\to \gal(E^{\prime,ab}/E')$.
Let $\eg_{K_p}(Y_\Phi)$ be the pullback of the extension $$0\lra \agsb_{K_p}(Y_\Phi)\lra \ag_{K_p}(Y_\Phi)\lra Y_{\Phi,K_p}(\zbkp)_+^{\overline{p}}\bss Y_{\Phi}(\Ap)\lra 0$$
along $-\wdtd{r}_{E',K_p}(Y_\Phi)$.
\begin{lem}\label{lem-egp-j}
With the previous definitions, for any $j\in \jg_{K_p}(Y_\Phi)$, there is a natural map $\eg_{\lcj{K_p}{j}}(Y_\Phi)\to \eg(Y_\Phi)$ that fits the following diagram:
\begin{equation}
    \begin{tikzcd}
        0\arrow[r]&\agsb_{\lcj{K_p}{j}}(Y_\Phi)\arrow[r]\arrow[d]&\eg_{\lcj{K_p}{j}}(Y_\Phi)\arrow[r]\arrow[d]&\gal(\overline{E}^{\prime}/E')\arrow[r]\arrow[d,equal]&0\\
        0\arrow[r]&\agsb(Y_\Phi)\arrow[r]&\eg(Y_\Phi)\arrow[r]&\gal(\overline{E}'/E')\arrow[r]&0.
    \end{tikzcd}
\end{equation}
Moreover, the group $\eg_{\lcj{K_p}{j}}(Y_\Phi)$ is the pullback of the extension 
$$0\lra \agsb_{\lcj{K_p}{j}}(Y_\Phi)\lra {\ag_{\lcj{K_p}{j}}(Y_\Phi)}\lra Y_{\Phi,\lcj{K_p}{j}}(\zbkp)_+^{\overline{p}}\bss Y_\Phi(\Ap) \lra 0$$
along $-\wdtd{r}_{E',\lcj{K_p}{j}}(Y_\Phi)$.
\end{lem}
\begin{proof}
    This is self-explanatory from the constructions above.
\end{proof}
\subsubsection{}
Let us return to the settings for $(G_0,X_0)$ and $(G_2,X_2)$.  
Let $\Phi_0$ be a cusp label representative of $(G_0,X_0)$ and let $\Phi_2$ be a cusp label representative of $(G_2,X_2)$. 
Assume that there is an isogeny $P'_{\Phi_0}\to Y'_{\Phi_2}$ whose kernel is in $Z_{G_0}\cap P'_{\Phi_0}$.\par
We assume that $K_{\Phi_0,p}\cap \agsb(P_{\Phi_0})$ maps to $K_{\Phi_2,p}^Y\cap \agsb(Y_{\Phi_2})$ and that $g_\Phi K_p g_\Phi^{-1}\cap P'_\Phi(\bb{Q}_p)$ is the pullback of $g_{\Phi_2}K_{2,p}g_{\Phi_2}^{-1}\cap Y'_{\Phi_2}(\bb{Q}_p)$ under the isogeny $P'_{\Phi_0}\to Y'_{\Phi_2}$.\par
Combining Lemma \ref{lem-jgp}, Lemma \ref{lem-egp-j} and Proposition \ref{prop-del-ind-mix}, we obtain a $Y_{\Phi_2}(\Ap)\times\gal(\overline{E}'/E')$-equivariant isomorphism (cf. (\ref{eq-induction-basic}))
\begin{equation}\label{eq-j-decompose-gen}
\begin{split}
&\sh_{K_{\Phi_2,p}^Y,\overline{\bb{Q}}}\iso (\sh^+_{P_{\Phi_0},\overline{\bb{Q}}}\times (\ag(Y_{\Phi_2})/K_{\Phi_2,p}^Y))/\agsb(P_{\Phi_0})\\
&\iso\disju_{j\in\jg_{K_{2,p}}(Y_{\Phi_2})}(\sh^+_{\lcj{K_{\Phi_0,p}}{j},\overline{\bb{Q}}}\times (\atd_{\lcj{K_{2,p}}{j}}(Y_{\Phi_2})/jK_{\Phi_2,p}^Yj^{-1}))/{\atdsb_{\lcj{K_{2,p}}{j}}(P_{\Phi_0})}\\
&\iso\disju_{j\in\jg_{K_{2,p}}(Y_{\Phi_2})}(\sh^+_{\lcj{K_{\Phi_0,p}}{j},\overline{\bb{Q}}}\times \ag_{\lcj{K_{2,p}}{j}}(Y_{\Phi_2}))/\agsb_{\lcj{K_{2,p}}{j}}(P_{\Phi_0}).  
\end{split}
\end{equation}
Since every component in the disjoint union above descends to $E'$, we can write $\sh_{K_{\Phi_0,p}}$ as a disjoint union over $E'$, 
$$\sh_{K_{\Phi_0,p}}\iso \disju_{j\in \jg_{K_{0,p}}(P_{\Phi_0})}\lcj{\sh}{j}_{K_{\Phi_0,p}}.$$
Let $\lcj{\ca{S}}{j}_{K_{\Phi_0,p}K_{\Phi_0}^p}$ be the normalization in $\lcj{\sh}{j}_{K_{\Phi_0,p}K_{\Phi_0}^p}$ of $\ca{S}_{K_{\Phi_0,p}K_{\Phi_0}^p}$.\par 
Define $\lcj{\ca{S}}{j}_{K_{\Phi_0,p}}:=\varprojlim_{K_{\Phi_0}^p\text{\ neat open compact}}\lcj{\ca{S}}{j}_{K_{\Phi_0,p}K_{\Phi_0}^p}$.\par
We also record a lemma.
\begin{lem}\label{lem-field-conn}
Let $K_p$ be an open compact subgroup of $G(\bb{Q}_p)$. Let $\Phi$ be a cusp label representative of $(G,X)$. Then there is a finite Galois extension $F$ of the reflex field $E:=E(G,X)$ depending on $K_p$ and $\Phi$, such that all connected components of $\sh_{K_{\Phi,p}}(G,X)_{\bb{C}}$ are defined over $F\cdot E^p$, where $E^p$ is the union of all finite extensions of $E$ that are unramified at $p$.
\end{lem}
\begin{proof}
By \cite[Prop. 11.2 (c)]{Pin89}, it suffices to show the lemma after replacing $(G,X)$ with $(P_{\Phi,h},\hbar (D_{\Phi,h}))$. So we assume that $P_\Phi=G$ is reductive without loss of generality.\par
Recall that by \cite[Thm. 2.6.3]{Del79}, the Galois group $\gal(E^{ab}/E)$ acts on the connected components by 
\begin{equation}\label{eq-conn-components-action}\begin{split}\gal(E^{ab}/E)\xrightarrow{rec_E}\pi_0(\bb{G}_m(\bb{A}_E)/\bb{G}_m(E))\xrightarrow{\sim}\bb{A}^\times_E/E^\times(E\otimes_\bb{Q}\bb{R})^{\times,+}\xrightarrow{[\mu_x]}\\
G(\bb{A}_{E,f})/G(E)^{\overline{\ }}_+\xrightarrow{N_{E/\bb{Q}}}G(\A)/G(\bb{Q})^{\overline{\ }}_+K_p.\end{split}\end{equation}
This map sends open compact groups of $\bb{G}_m(E\otimes_\bb{Q} \bb{Q}_p)$ to compact subgroups of $G(\bb{Q}_p)$. So there is a possibly ramified finite field extension $F$ of $E$, such that the image of $\gal(E^{ab}/F)$ under the map above is trivial at $p$. So we have the desired result.
\end{proof}
\subsection{Construction of strata and compactifications by quotients}\label{subsec-con-stra-comp-quo}
We continue with the settings in \S\ref{sss-conclusion}. Let $(G,X_a)$, $(G,X_b)$ and $K=K_pK^{p}\sbst G(\A)$ be the Shimura data and the neat open compact subgroup chosen there. Since in most occasions only $(G,X_b)$ will appear in this {subsection}, we will omit the superscript $b$ from many symbols in the first {section}.\par 
Let $\Phi:=\Phi_b=(Q,X_b^+,g^b)$ be any cusp label representative of $(G,X_b)$. Set $P_\Phi:=P_\Phi^b$ and $ZP_\Phi:=ZP_\Phi^b$. Let $K_\Phi:=K_\Phi^b=P_{\Phi}^b(\A)\cap g^b K g^{b,-1}$ and let $\wdtd{K}_\Phi:=\wdtd{K}^b_\Phi=ZP_\Phi^b(\A)\cap g^bK g^{b,-1}$. Recall that $\sh_{\K_\Phi}:=\sh_{\wdtd{K}^b_\Phi}(ZP^b_\Phi,ZP^b_\Phi(\bb{Q})D_{\Phi_b})$ is the mixed Shimura variety defined by $(ZP^b_\Phi,ZP^b_\Phi(\bb{Q})D_{\Phi_b})$. In our conventions $\sigma$ are cones in $ \Sigma(\Phi)$, $\sigma_0^\alpha$ are cones in $\Sigma_0^\alpha(\Phi_0^\alpha)$, etc. Recall that we chose $\Sigma_0^\alpha$ to be \emph{induced} by $\Sigma$ (see Proposition \ref{zp-cones}). Let us explain how to use quotients to construct integral models $\ca{S}_{\wdtd{K}_\Phi}$,  $\ca{S}_{\wdtd{K}_\Phi,\sigma}$ and $\ca{S}_{\wdtd{K}_\Phi}(\sigma)$ for $\sh_{\wdtd{K}_\Phi}$, $\sh_{\wdtd{K}_\Phi,\sigma}$ and $\sh_{\wdtd{K}_\Phi}(\sigma)$, respectively.\par
We will study the following cases, where each of them is a special case of the one following it:
\begin{itemize}
\item (Hyperspecial levels; \textbf{HS}) See \cite{Kis10} and \cite{KM15}. Suppose that $G_2$ is quasi-split and unramified at $p$. By \cite[Lem. 3.4.13]{Kis10}, we can assume that $G_0$ is chosen such that $G_{0,\bb{Q}_p}$ is unramified. Then we can choose $G$ such that the center of $G$ is quasi-split and unramified at $p$ (and $G^\der=G_2^\der$). Hence, we can and we shall choose $G_0$ and $G$ such that $G_0$ and $G$ are both quasi-split and unramified at $p$. Furthermore, we choose all of $g_\alpha$ in $G(\Ap)$ by \cite[Lem. 2.2.6]{Kis10}.\par 
Let $G_{2,\zbkp}$ be a smooth reductive model over $\zbkp$ such that $G_{2,\zbkp}(\bb{Z}_p)=K_{2,p}$. As in the proof of \cite[Cor. 3.4.14]{Kis10}, the center $Z'$ of $G$ extends to a smooth reductive model $Z'_{\zbkp}$ over $\zbkp$, and the kernel of the isogeny $G_2^\der\times Z'\to G$ also extends to a finite flat group over $\zbkp$. Hence, we choose a smooth reductive model $G_{\zbkp}$ of $G$ over $\zbkp$ such that $G_{2,\zbkp}$ embeds into $G_\zbkp$ and $K_p=G_{\zbkp}(\bb{Z}_p)$. 
By \emph{loc. cit}., we can choose a hyperspecial level $K_{0,p}^\alpha= G_{0,\bb{Z}_p}(\bb{Z}_p)$, where $G_{0,\bb{Z}_p}$ is a reductive group scheme over $\bb{Z}_p$, such that the kernel of $G_{0,\bb{Z}_p}\to G_{\bb{Z}_p}$ is finite and is contained in the center of $G^{\der}_{0,\bb{Z}_p}$. Here the subgroups $K_{0,p}^\alpha$ are chosen to be the same. Write this group as $K_{0,p}:=K_{0,p}^\alpha$. Then we can choose $K^{\ddag,\alpha}_p=:K^\ddag_p$ to be the same. We also require the embeddings to be \textbf{$p$-integral} as in \cite[Sec. 4.3]{Mad19}.
\item (Bruhat-Tits stabilizer levels; \textbf{STB}) Fix a Bruhat-Tits stabilizer group scheme $\mathscr{G}_{x_2}$ corresponding to a point $x_2$ in the extended Bruhat-Tits building $\ca{B}(G_{2,\bb{Q}_p},\bb{Q}_p)$. Choose a Hodge-type Shimura datum $(G_0,X_0)$ such that there is a central isogeny $G_0^\der\to G_2^\der$ which induces an isomorphism $(G_2^\ad,X_2^\ad)\iso (G_0^\ad,X_0^\ad)$. We further require that $(G_0,X_0)$ satisfies the condition (3) in \cite[Lem. 4.6.22]{KP15}. By making this choice, we have that, for any place $v_2|p$ of $E_2$ and any place $v'|v_2$ of $E'$, $\ca{O}_{E_{2,v_2}}\iso \ca{O}_{E'_{v'}}$. Note that the argument in the proof of \emph{loc. cit.} for only part (3) still works for general cases.  
Let $x$ be the image of $x_2$ in $\ca{B}(G_{\bb{Q}_p},\bb{Q}_p)$ and $x$ determines a Bruhat-Tits stabilizer group scheme $\mathscr{G}_x$. By our construction and by the definition of Bruhat-Tits stabilizer group schemes, we see that $\mathscr{G}_x(\bb{Z}_p)\cap G_2(\bb{Q}_p)=\mathscr{G}_{x_2}(\bb{Z}_p)$. Let $K_{2,p}:=\mathscr{G}_{x_2}(\bb{Z}_p)$ and let $K_p:=\mathscr{G}_{x}(\bb{Z}_p)$. \par
For any $g_\alpha$, let $K_{0,p}^\alpha$ be be the preimage of $g_\alpha K_pg_\alpha^{-1}$ under $\pi^b$; since the kernel of $G_0\to G$ is finite and is contained in the center of $G_0^\der$, $K_{0,p}^\alpha$ is a Bruhat-Tits stabilizer subgroup of $G_0(\bb{Q}_p)$. Write $K_{0,p}^{\alpha}:=\mathscr{G}_0^\alpha(\bb{Z}_p)$, where $\mathscr{G}_0^\alpha$ is a Bruhat-Tits stabilizer group scheme corresponding to a point $x_0^\alpha$ in the extended building $\ca{B}(G_{0,\bb{Q}_p},\bb{Q}_p)$. By \cite{Lan00} and as explained in \cite[4.5.2, p. 86]{PR24}, for any Hodge embedding $(G_0,X_0)\hookrightarrow (\mrm{GSp}(V,\psi),X^\ddag)$, up to replacing $V$ with $V^{\oplus n}$ for some $n$, we can choose a $\bb{Z}_p$-lattice $V_{\bb{Z}_p}^\alpha$ such that $K_{0,p}^\alpha$ is the stabilizer in $G_0(\bb{Q}_p)$ of $V_{\bb{Z}_p}^\alpha$. 
By Zarhin's trick, we can further replace $V$ with $V^{\oplus 8}$, so that there is a self-dual $\bb{Z}$-lattice $V_{\bb{Z}}^\alpha$, such that $\stb_{G_0(\bb{Q}_p)}(V_{\bb{Z}_p}^\alpha)=K_{0,p}^\alpha$ and $K_p^{\ddag,\alpha}:=\stb_{G^\ddag(\bb{Q}_p)}(V_{\bb{Z}_p}^\alpha)$. 
Now we have chosen for any $g_\alpha$ a Hodge embedding $\iota^\alpha: (G_0,X_0)\hookrightarrow (\mrm{GSp}((V,\psi)^{\perp r(\alpha)}),X^{\ddag,r(\alpha)})$ for some positive integer $r(\alpha)$ and a $\bb{Z}$-lattice $V_\bb{Z}^\alpha$ satisfying the setup in \cite[Sec. 3.1]{Mad19} and \S\ref{subsec-siegel}. 
We can further assume that $r(\alpha)=r(\alpha')$ for different $\alpha\neq \alpha'\in I_{G/G_0}$ by further replacing all $r(\alpha)$ with their lcm. Hence, we can and we will write $(\mrm{GSp}((V,\psi)^{\perp r(\alpha)}),X^{\ddag,\alpha})$ as $(G^\ddag,X^\ddag)$.
\item (Finite Intersections of Bruhat-Tits stabilizers; $\mbf{STB}_n$) Slightly generalizing the case above, we still assume that $K_p\cap G_2(\bb{Q}_p)=K_{2,p}$ but only require that $K_p:=\cap_{i=1}^n \mathscr{G}_{x_i}(\bb{Z}_p)$, where $\{x_i\}_{i=1}^n$ is a set of $n$ points in $\ca{B}(G_{\bb{Q}_p},\bb{Q}_p)$ and $\mathscr{G}_{x_i}$ are corresponding Bruhat-Tits stabilizer group schemes of $x_i$.
Pulling back to $G_0$, we have an intersection $K_{0,p}^\alpha:=\pi^{b,-1}(g_\alpha\cap_{i=1}^n\mathscr{G}_{x_i}(\bb{Z}_p)g_\alpha^{-1})=\cap_{i=1}^n\pi^{b,-1}(g_\alpha\mathscr{G}_{x_i}(\bb{Z}_p)g_\alpha^{-1})$. Each $\pi^{b,-1}(g_\alpha\mathscr{G}_{x_i}(\bb{Z}_p)g_\alpha^{-1})=\mathscr{G}_{x_{i,0}^\alpha}(\bb{Z}_p)$, where $\mathscr{G}_{x_{i,0}^\alpha}$ is a Bruhat-Tits stabilizer group scheme correspoinding to $x_{i,0}^\alpha\in \ca{B}(G_{0,\bb{Q}_p},\bb{Q}_p)$ and $x_{i,0}^\alpha$ corresponds to the same point as $g_\alpha\cdot x_i$ in the (reduced) Bruhat-Tits Building of $G^\ad$ over $\bb{Q}_p$. 
Write $K_{0,p}^{(\alpha,i)}:=\mathscr{G}_{x^\alpha_{i,0}}(\bb{Z}_p)$.
Exactly as above, by \cite[4.5.2]{PR24} and by Zarhin's trick, we find a Hodge embedding $\iota^{(\alpha,i)}:(G_0,X_0)\hookrightarrow (\mrm{GSp}((V,\psi)^{\perp r(\alpha,i)}),X^{\ddag,r(\alpha,i)})$ for some positive integer $r(\alpha,i)$ for each $(\alpha,i)$ and a $\bb{Z}$-lattice $V_\bb{Z}^{(\alpha,i)}\sbst V^{\oplus r(\alpha,i)}$ satisfying the setup in \cite[Sec. 3.1]{Mad19} and \S\ref{subsec-siegel}, such that $K_{0,p}^{(\alpha,i)}$ is exactly the stabilizer of $V_{\bb{Z}_p}^{(\alpha,i)}$ in $G_0(\bb{Q}_p)$.
Let $r_{\mrm{sum}}(\alpha):=\sum_{i=1}^nr(\alpha,i)$ and $V_{\mrm{sum},\bb{Z}}^\alpha:=\oplus_{i=1}^n V_\bb{Z}^{(\alpha,i)}\sbst V^{\oplus r_{\mrm{sum}}(\alpha)}$. Let $r$ be the lcm of all $r_{\mrm{sum}}(\alpha)$ for all $\alpha$. Hence, up to replacing $(V,\psi)$ with its $r$-th orthogonal direct sum, and replacing $V^\alpha_{\mrm{sum},\bb{Z}}$ with its $r/r_{\mrm{sum}}(\alpha)$-th copies for each $\alpha$, we have found a Hodge embedding $(G_0,X_0)\hookrightarrow(G^\ddag,X^\ddag)$ and a lattice $V_\bb{Z}^\alpha$ satisfying \cite[Sec. 3.1]{Mad19} for each $\alpha$. Let $K_p^{\ddag,\alpha}:=\stb_{G^\ddag(\bb{Q}_p)}(V_{\bb{Z}_p}^\alpha)$ for each $\alpha$. 
By definition, (STB)$=$($\mrm{STB}_1$).
\item (Deeper levels in (HS) and ($\mrm{STB}_n$); \textbf{DL}) Under the assumptions of either (HS) or ($\mrm{STB}_n$), we choose an open compact subgroup $K_p\sbst K_p^*$ where $K_p^*$ is a hyperspecial subgroup in Case (HS), or an intersection of Bruhat-Tits stabilizer subgroups in Case ($\mrm{STB}_n$). 
Let $I_{G/G_0}$ (resp. $I_{G/G_0}^*$) be the double coset $I_{G/G_0,K}$ (resp. $I_{G/G_0,K^*}$) defined in \S\ref{sss-conclusion} with $K=K_pK^p$ (resp. $K^*:=K_p^*K^p$). Fix a choice of objects $(G_0,X_0)$ and $(G^\ddag,X^\ddag)$; for any $\alpha^*\in I_{G/G_0}^*$, fix a choice of representatives $g_{\alpha^*}$, and open compact subgroups $K_{0,p}^{\alpha^*}$ and $K_p^{\ddag,\alpha^*}$ as in (HS) or ($\mrm{STB}_n$). 
For any $\alpha\in I_{G/G_0}$ mapping to $\alpha^*\in I_{G/G_0}^*$, we choose the representative $g_{\alpha}$ to be $g_{\alpha}=g_{\alpha^*}\cdot k_{\alpha^*}$ for some $k_{\alpha^*}\in K_p^*$.
Let $K_{0,p}^{\alpha}$ be the preimage of $g_{\alpha} K_p g_{\alpha}^{-1}$ under $\pi^b$. Then $K_{0,p}^{\alpha}\sbst K_{0,p}^{\alpha^*}\sbst K_p^{\ddag,\alpha^*}$ by construction. Hence, the conditions in \cite[Sec. 3.1]{Mad19} and \S\ref{subsec-siegel} are still satisfied. Note that for different choices of intersections of Bruhat-Tits stabilizer subgroups $K_p^*$ of $G(\bb{Q}_p)$ containing $K_p$, the integral models constructed later might be different.
\end{itemize}
In all four cases, let $K_0^{\alpha,p}$ be a neat open compact subgroup of $G_0(\Ap)$ contained in $\pi^{b,-1}(g_\alpha K^pg_\alpha^{-1})$. Let $K_0^\alpha=K_{0,p}^\alpha K_0^{\alpha,p}$. Since any open compact subgroup of $G_2(\bb{Q}_p)$ has a fixed point in the extended building $\ca{B}(G_{2,\bb{Q}_p},\bb{Q}_p)$ (see, e.g., \cite[Thm. 1.4.4]{Lan00}), for any open compact subgroup $K_{2,p}$, we can choose $n$ Bruhat-Tits stabilizer subgroups in $G_2(\bb{Q}_p)$ for some $n$ and implement one of the constructions above.\par
In summary, for all four cases, we have chosen $(G_0,X_0,K_{0,p}^\alpha,K_{0}^{\alpha,p})\to (G^\ddag,X^\ddag,K_p^{\ddag,\alpha},K^{\ddag,\alpha,p})$ from $(G_2,X_2,K_{2,p},K^p_2)$ for each $\alpha\in I_{G/G_0}$ with the settings in \cite{Mad19}. 
In Case (HS) (resp. (STB)), we shall additionally require the embeddings to satisfy the settings in \cite[4.3.2]{Mad19} (resp. \cite{DvHKZ24}). Let $\ca{S}_{K_0^\alpha}$ be the relative normalization of $\ca{S}_{K^{\ddag,\alpha}}$ in $\sh_{K_0^\alpha}$, and by the main theorem of \cite{Mad19} (cf. Theorem \ref{thm-tor-hdg}), $\ca{S}_{K_0^\alpha}$ admits good toroidal compactifications $\ca{S}_{K_0^\alpha}^{\Sigma_0^\alpha}$ with integral models of boundary mixed Shimura varieties $\ca{S}_{K_{\Phi_0^\alpha}}$. The integral model $\ca{S}_{K_{\Phi_0^\alpha}}$ (resp. $\overline{\ca{S}}_{K_{\Phi_0^\alpha}}$ and resp. $\ca{S}_{K_{\Phi_0^\alpha},h}$) is constructed as the relative normalization in $\sh_{K_{\Phi_0^\alpha}}$ (resp. $\overline{\sh}_{K_{\Phi_0^\alpha}}$ and resp. $\sh_{K_{\Phi_0^\alpha},h}$) of $\ca{S}_{K_{\Phi^{\ddag,\alpha}}}$ (resp. $\overline{\ca{S}}_{K_{\Phi^{\ddag,\alpha}}}$ and resp. $\ca{S}_{K_{\Phi^{\ddag,\alpha}},h}$) where $\Phi_0^\alpha\in \ca{CLR}(G_0,X_0)$ maps to $\Phi^\ddag\in \ca{CLR}(G^\ddag,X^\ddag)$.\par
Fix any place $v|p$ of $E_0$. Let $\ca{O}_{(v)}':= \ca{O}_{E'}\otimes_{\ca{O}_{E_0}}\ca{O}_{E_0,(v)}$ and let $\ca{O}_{v}':= \ca{O}_{E'}\otimes_{\ca{O}_{E_0}}\ca{O}_{E_0,v}$. 
We will mainly work over $\ca{O}_{(v)}'$; but we will finally state the main results over $\ca{O}_{E_2,(v_2)}$ in Case (HS) and its deep-level case, and over $\ca{O}_{E_2,v_2}$ in ($\mrm{STB}_n$) and its deep-level case.
We will clarify the base ring if it is needed. Let $E^{\prime,p}$ be the maximal field extension of $E'$ that is unramified at all primes dividing $p$. Write $\ca{O}_{(v)}^{\prime,ur}:=\ca{O}_{E^{\prime,p}}\otimes_{\ca{O}_{E'}}\ca{O}_{(v)}'$.
\subsubsection{}\label{subsubsec-hodge-gad1q}We will omit and the readers can ignore the superscript $\alpha$ before introducing Lemma \ref{lem-gadq1-ext}. Recall that $\pi:=\pi^b:(G_0,X_0)\to (G,X_b)$ is a map such that the kernel of $G_0\to G$ is finite and is in the center of $G_0^\der$. Let $\Phi_0:=(Q_0,X_0^+,g_0)$ be a cusp label representative in $\ca{CLR}(G_0,X_0)$ mapping to $\Phi^\ddag:=(Q^\ddag, X^{\ddag,+},g^\ddag)\in \ca{CLR}(G^\ddag,X^\ddag)$.
Recall that for $\gamma\in G^\ad_0(\bb{Q})$, there is a natural morphism $\mrm{int}(\gamma):(G_0,X_0)\to (G_0,\gamma\cdot X_0)$, which induces isomorphisms (see \cite[4.16]{Pin89})
$$\sh_{K_{\Phi_0}}\to \sh_{\lcj{K_{\Phi_0}}{\gamma}},$$
which maps $[(x,p)]_{P_{\Phi_0}(\A)\cap g_{0}K g_{0}^{-1}}$ to $[(\gamma(x),\gamma p \gamma^{-1})]_{\lcj{(P_{\Phi_0}(\A)\cap g_{0}K g_{0}^{-1})}{\gamma}}$ over complex points. Similarly, we have morphisms induced by $\gamma$ for $\overline{\sh}_{K_{\Phi_0}}$ and $\sh_{K_{\Phi_0},h}$.\par
By \cite[Thm. 3.58(3)]{Mao25}, if $K_{0,p}$ is a Bruhat-Tits stabilizer subgroup, so is $K_{\Phi_0,h,p}$, and we have $K_{\Phi_0,h,p}=G_0(\bb{Q}_p)\cap K_{\Phi^\ddag,h,p}$. Note that if $K_{0,p}$ is parahoric, we can only assume that $K_{\Phi_0,h,p}$ is quasi-parahoric (see Notation and Conventions). \par
Note that $\gamma$ maps $\Phi_0$ to $\gamma\cdot \Phi_0:=(\gamma Q_0 \gamma^{-1},\gamma(X_0^+),\gamma g_0 \gamma^{-1})\in \ca{CLR}(G_0,\gamma\cdot X_0)$ (see \S\ref{cusp-label}). If $K_{0,p}$ is a Bruhat-Tits stabilizer subgroup, so is $\lcj{K_{0,p}}{\gamma}:=\gamma K_{0,p} \gamma^{-1}$. \par
As in (STB), we can choose a Hodge embedding $(G_0,\gamma\cdot X_0)\hookrightarrow (G^{\ddag,\prime},X^{\ddag,\prime})$ and a Bruhat-Tits stabilizer subgroup $K_p^{\ddag,\prime}$ of $G^{\ddag,\prime}(\bb{Q}_p)$, such that $K_p^{\ddag,\prime}$ is the stabilizer in $G^{\ddag,\prime}(\bb{Q}_p)$ of a self-dual lattice satisfying the settings in \cite[Sec. 3.1]{Mad19} and such that $\lcj{K_{0,p}}{\gamma}=K_p^{\ddag,\prime}\cap G_0(\bb{Q}_p)$. Suppose that $\gamma\cdot \Phi_0$ maps to a cusp label representative $\Phi^{\ddag,\prime}\in \ca{CLR}(G^{\ddag,\prime},X^{\ddag,\prime})$. Then one can define an integral model $\ca{S}_{\lcj{K_{\Phi_0}}{\gamma}}$ (resp. $\overline{\ca{S}}_{\lcj{K_{\Phi_0}}{\gamma}}$ and resp. $\ca{S}_{\lcj{K_{\Phi_0},h}{\gamma}}$) of $\sh_{\lcj{K_{\Phi_0}}{\gamma}}$ (resp. $\overline{\sh}_{\lcj{K_{\Phi_0}}{\gamma}}$ and resp. $\sh_{\lcj{K_{\Phi_0},h}{\gamma}}$) as the relative normalization of $\ca{S}_{K_{\Phi^{\ddag,\prime}}}$ (resp. $\overline{\ca{S}}_{K_{\Phi^{\ddag,\prime}}}$ and resp. $\ca{S}_{K_{\Phi^{\ddag,\prime}},h}$) in $\sh_{\lcj{K_{\Phi_0}}{\gamma}}$ (resp. $\overline{\sh}_{\lcj{K_{\Phi_0}}{\gamma}}$ and resp. $\sh_{\lcj{K_{\Phi_0},h}{\gamma}}$) as \cite{Mad19}.
\begin{lem}\label{lem-gadq1-ext} Let $v$ be any place of $E_0$ over $p$. 
Under the assumptions in (HS), (STB) or (DL) of (HS)/(STB), the action of $\gamma\in G_0^\ad(\bb{Q})$, $\gamma: \sh_{K_{\Phi_0^\alpha}}\to \sh_{\lcj{K_{\Phi_0^\alpha}}{\gamma}}$, extends to an isomorphism $\ca{S}_{K_{\Phi_0^\alpha}}\to \ca{S}_{\lcj{K_{\Phi_0^\alpha}}{\gamma}}$ over $\ca{O}_{E_0,(v)}$ for any cusp label representative $\Phi_0^\alpha$ and neat open compact subgroup $K_0^\alpha=K_{0,p}^\alpha K_0^{\alpha,p}$. 
Moreover, the morphism $\gamma: \ca{S}_{K_{\Phi_0^\alpha}}\to \ca{S}_{\lcj{K_{\Phi_0^\alpha}}{\gamma}}$ as a morphism between torus torsors is equivariant under the action of $\mbf{E}_{K_{\Phi_0^\alpha}}\to \mbf{E}_{\lcj{K_{\Phi_0^\alpha}}{\gamma}}$. Similarly, there are morphisms $\gamma:\overline{\ca{S}}_{K_{\Phi_0^\alpha}}\to\overline{\ca{S}}_{\lcj{K_{\Phi_0^\alpha}}{\gamma}}$ and $\gamma: \ca{S}_{K_{\Phi_0^\alpha},h}\to \ca{S}_{\lcj{K_{\Phi_0^\alpha},h}{\gamma}}$ extending the conjugation of $\gamma$ over the generic fiber.\par
The integral models we constructed depend only on the choice of a Bruhat-Tits stabilizer subgroup in $G_0(\bb{Q}_p)$ and are independent of the choice of the Hodge embedding in the construction made above.
\end{lem}
\begin{proof} Write $\Phi_0:=\Phi_0^\alpha$ and $K_0:=K_0^\alpha$ to simplify the notation. 
We only have to show the first two cases since the last case follows from taking relative normalizations. In Case (HS), it follows from \cite[Prop. 4.3.5]{Mad19} that the action of $\gamma$ extends to $\ca{S}_{K_{\Phi_0},p}\to \ca{S}_{\lcj{K_{\Phi_0},p}{\gamma}}$. With the same reference, we see that this extended map is $P_{\Phi_0}(\Ap)$-equivariant, and this proves the desired statement.\par 
Since the schemes involved are normal and separated, and by gluing two morphisms over $\ca{O}_{E_0,v}$ and $E_0$ that are identified over $E_{0,v}$, it suffices to show the statement over $\ca{O}_{E_0,v}$. In Case (STB), by \cite[Thm. 4.1.5(1)]{Mad19}, the tower $\ca{S}_{K_{\Phi_0}}\to \overline{\ca{S}}_{K_{\Phi_0}}\to \ca{S}_{K_{\Phi_0},h}$ is a torus torsor $\ca{S}_{K_{\Phi_0}}\to \overline{\ca{S}}_{K_{\Phi_0}}$ under $\mbf{E}_{K_{\Phi_0}}$ composing with an abelian scheme torsor $\overline{\ca{S}}_{K_{\Phi_0}}$ over $\ca{S}_{K_{\Phi_0},h}$ under $\ca{A}_{K_{\Phi_0}}$.\par 
From the theory of canonical integral models of Shimura varieties of Hodge type (see \cite[Cor. 4.3.2]{PR24}, \cite[Thm. 4.14]{DY25} and \cite[Cor. 4.1.10]{DvHKZ24}) and Corollary \ref{cor-pure-shimura}, there is a morphism $\gamma:\ca{S}_{K_{\Phi_0},h}\to \ca{S}_{\lcj{K_{\Phi_0},h}{\gamma}}$ extending the one over the generic fiber. 
By the valuative criterion of properness and Weil's extension theorem (see \cite[Sec. 4.4, Thm. 1]{BLR90}), there is an extension $\overline{\ca{S}}_{K_{\Phi_0}}\to \overline{\ca{S}}_{\lcj{K_{\Phi_0}}{\gamma}}$. For the extension of the whole $\ca{S}_{K_{\Phi_0}}$, we adopt the following method (with which we can also obtain the extension $\gamma:\overline{\ca{S}}_{K_{\Phi_0}}\to \overline{\ca{S}}_{\lcj{K_{\Phi_0}}{\gamma}}$).\par 
Suppose that $\ca{S}_{K_0}$ and $\ca{S}_{K_{\Phi_0}}$ are constructed from a Hodge embedding $\iota:(G_0,X_0,K_{0,p})\hookrightarrow (\mrm{GSp}(V,\psi),X^\ddag,K^\ddag_p)$, and $\ca{S}_{\lcj{K_0}{\gamma}}$ and $\ca{S}_{\lcj{K_{\Phi_0}}{\gamma}}$ are constructed from another Hodge embedding $\iota':(G_0,\gamma\cdot X_0,\lcj{K_{0,p}}{\gamma})\hookrightarrow (\mrm{GSp}(V',\psi'),X^{\ddag,\prime},K_p^{\ddag,\prime})$. Let $V'':=V\oplus V'$ with symplectic pairing $\psi''$ on it defined by the direct sum of the two factors. Let $\iota'':=\iota\times (\iota'\circ\mrm{int}(\gamma))$. Suppose that $K_p^{\ddag}$ (resp. $K_p^{\ddag,\prime}$) is the stabilizer in $\mrm{GSp}(V,\psi)(\bb{Q}_p)$ (resp. $\mrm{GSp}(V',\psi')(\bb{Q}_p)$) of $\Lambda$ (resp. $\Lambda'$), where $\Lambda$ (resp. $\Lambda'$) is the base change to $\bb{Z}_p$ of some self-dual $\bb{Z}$-lattice $V_\bb{Z}$ (resp. $V_\bb{Z}'$). Set $\Lambda'':=\Lambda\oplus \Lambda'$, $V''_\bb{Z}:=V_\bb{Z}\oplus V_\bb{Z}'$ and $K_p^{\ddag,\prime\prime}:=\stb_{\mrm{GSp}(V'',\psi'')(\bb{Q}_p)}(\Lambda'')$. Let $G^\ddag_{\zbkp}$ (resp. $G^{\ddag,\prime}_{\zbkp}$ and resp. $G^{\ddag,\prime\prime}_\zbkp$) be the smooth reductive group scheme over $\zbkp$ whose group of $\bb{Z}_p$-points is $K_p^\ddag$ (resp. $K^{\ddag,\prime}_p$ and resp. $K^{\ddag,\prime\prime}_p$).
Then there are integral models $\ca{S}_{K_0}''$ and $\ca{S}_{K_{\Phi_0}}''$ of $\sh_{K_0}$ and $\sh_{K_{\Phi_0}}$, respectively, constructed from $\iota''$ and $\Lambda''$ as in \cite{PR24}, \cite{DvHKZ24} and \cite{Mad19} by taking relative normalizations.\par 
We claim that we have a diagram
\begin{equation}\label{eq-change-of-lattice-mixsh}\ca{S}_{K_{\Phi_0}}\longleftarrow\ca{S}''_{K_{\Phi_0}}\lra \ca{S}_{\lcj{K_{\Phi_0}}{\gamma}}.\end{equation}
In fact, let $H'':=(\mrm{GSp}(V,\psi)\times \mrm{GSp}(V',\psi'))\times_{\bb{G}_m\times\bb{G}_m}\bb{G}_m$, where the left factor maps to $\bb{G}_m\times\bb{G}_m$ by a product of similitude characters, and $\bb{G}_m$ maps to $\bb{G}_m\times\bb{G}_m$ by diagonal embedding. Then the embedding $G_0\hookrightarrow \mrm{GSp}(V'',\psi'')$ factors through $H''$, and this induces a Shimura datum $(H'',X'')$ associated with $H''$. The morphism $(H'',X'')\hookrightarrow (\mrm{GSp}(V'',\psi''),X^{\ddag,\prime\prime})$ is an embedding, and there are natural projections (see Remark \ref{rk-kisin-zhou-deleted-part}) 
$$(\mrm{GSp}(V,\psi),X^\ddag)\longleftarrow (H'',X'')\longrightarrow (\mrm{GSp}(V',\psi'),X^{\ddag,\prime}).$$
Let $h$ be any element in $X''$. Its weight cocharacter $\omega_h$ sends $\bb{G}_m$ to $\omega_h(\bb{G}_m)$, and the intersection $C''$ of $\omega_h(\bb{G}_m)$ and $\mrm{Sp}(V,\psi)\times\mrm{Sp}(V',\psi')$ is isomorphic to $\mu_2$. 
The split torus $\omega_{h}(\bb{G}_m)$ canonically extends to a split torus $\omega_h(\bb{G}_m)_{\zbkp}$ over $\zbkp$, and the finite group $C'$ extends to a finite flat group scheme $C'_{\zbkp}$ over $\zbkp$. Let $H''_{\zbkp}:=(G^{\ddag,\der}_{\zbkp} \times G^{\ddag,\prime,\der}_{\zbkp}\times \omega_h(\bb{G}_m)_{\zbkp})/C'_{\zbkp}$ be a reductive model over $\zbkp$ extending $H''$. By \cite[Cor. 1.3]{PY06} and the proof of \cite[Lem. 4.7]{KM15}, there is an embedding
$H''_{\zbkp}\hookrightarrow G^{\ddag,\prime\prime}_\zbkp$
between reductive group schemes over $\zbkp$ extending the embedding $H''\hookrightarrow \mrm{GSp}(V'',\psi'')$. Let $K''_p=H''_\zbkp(\bb{Z}_p)$. Then $\iota_H'':(H'',X'',K_p'')\hookrightarrow(\mrm{GSp}(V'',\psi''),X^{\ddag,\prime\prime},K_p^{\ddag,\prime\prime})$ is $p$-integral in the sense of \cite[4.3.2]{Mad19}. 
Suppose that $\Phi_0$ maps to $\Phi''\in\ca{CLR}(H'',X'')$, and then to $\Phi^{\ddag,\prime\prime}\in\ca{CLR}(G^{\ddag,\prime\prime},X^{\ddag,\prime\prime})$.
Then by \cite[Prop. 4.3.5]{Mad19}, the integral model $\ca{S}_{K_{\Phi'',p}}:=\ca{S}_{K_{\Phi'',p}}(P_{\Phi''},D_{\Phi''})$ associated with any cusp label representative $\Phi''\in\ca{CLR}(H'',X'')$ (constructed from taking relative normalization of $\ca{S}_{K_{\Phi^{\ddag,\prime\prime}}}$ in $\sh_{K_{\Phi''}}$ as in \cite[Sec. 4.1]{Mad19}) has the extension property. 
Then $\ca{S}''_{K_{\Phi_0}}$ constructed from $\Phi^{\ddag,\prime\prime}$ can also be constructed from taking relative normalization of $\ca{S}_{K_{\Phi''}}$ in $\sh_{K_{\Phi_0}}$ for suitable away-from-$p$ levels. By the extension property and by taking quotient of away-from-$p$ levels as in the first paragraph, we have a diagram
\begin{equation*}
   \begin{tikzcd}
       \sh_{K_{\Phi_0}}\arrow[d]&\sh_{K_{\Phi_0}}\arrow[d]\arrow[l,"\mrm{id}"']\arrow[r,"\gamma"]&\sh_{\lcj{K_{\Phi_0}}{\gamma}}\arrow[d]\\
       \ca{S}_{K_{\Phi^\ddag}}&\ca{S}_{K_{\Phi''}} \arrow[l]\arrow[r]&\ca{S}_{K_{\Phi^{\ddag,\prime}}}.
   \end{tikzcd} 
\end{equation*}
Now we have obtained the desired diagram by functoriality of relative normalizations.
\par 
By \cite[Thm. 4.1.5 (4)]{Mad19}, the two morphisms in (\ref{eq-change-of-lattice-mixsh}) are morphisms between $\mbf{E}_{K_{\Phi_0}}\iso \mbf{E}_{\lcj{K_{\Phi_0}}{\gamma}}$-torsors over $\overline{\ca{S}}_{K_{\Phi_0}}\leftarrow\overline{\ca{S}}''_{K_{\Phi_0}} \to\overline{\ca{S}}_{\lcj{K_{\Phi_0}}{\gamma}}$ since this is so over generic fibers. Similarly, $\overline{\ca{S}}_{K_{\Phi_0}}''$ (resp. $\overline{\ca{S}}_{\lcj{K_{\Phi_0}}{\gamma}}$) is also an abelian scheme torsor under some $\ca{A}''$ (resp. $\lcj{\ca{A}}{\gamma}$) over $\ca{S}_{K_{\Phi_0},h}''$ (resp. $\ca{S}_{\lcj{K_{\Phi_0},h}{\gamma}}$) by \cite[Thm. 4.1.5(1)]{Mad19}. More precisely, we have a commutative diagram of schemes
\begin{equation}\label{diag-mixsh-functoriality}
\begin{tikzcd}
\ca{S}_{K_{\Phi_0}}\arrow[d]&\ca{S}_{K_{\Phi_0}}''\arrow[d]\arrow[l]\arrow[r,"\gamma"]&\ca{S}_{\lcj{K_{\Phi_0}}{\gamma}}\arrow[d]\\
\overline{\ca{S}}_{K_{\Phi_0}}\arrow[d]&\overline{\ca{S}}_{K_{\Phi_0}}''\arrow[d]\arrow[l]\arrow[r,"\gamma"]&\overline{\ca{S}}_{\lcj{K_{\Phi_0}}{\gamma}}\arrow[d]\\
\ca{S}_{K_{\Phi_0},h}&\ca{S}_{K_{\Phi_0},h}''\arrow[l,"\iso"]\arrow[r,"\gamma","\iso"']&{\ca{S}_{\lcj{K_{\Phi_0}}{\gamma},h}}
\end{tikzcd}
\end{equation}
By \cite[I. Prop. 2.7]{FC90} and characteristic zero theory, there is a diagram of isomorphisms between abelian schemes $\ca{A}_{K_{\Phi_0}}\xleftarrow{\sim}\ca{A}''\xrightarrow[\sim]{\gamma}\lcj{\ca{A}}{\gamma}$; checking over an {\'e}tale cover of the third row such that the abelian scheme torsors in the second row all split, the second row of (\ref{diag-mixsh-functoriality}) are isomorphisms between abelian scheme torsors. Similarly, the arrows in the first row of (\ref{diag-mixsh-functoriality}) are isomorphisms between torus torsors. Now we have obtained the desired extended isomorphism $\ca{S}_{K_{\Phi_0}}\iso \ca{S}''_{K_{\Phi_0}}\to \ca{S}_{\lcj{K_{\Phi_0}}{\gamma}}$.\par
The last sentence follows from taking $\gamma=1$.
\end{proof}
The proof above essentially implies the following result:
\begin{cor}\label{cor-gadq1-ext}
Assume that $(G_2,X_2)$ is of Hodge type and choose $(G_0,X_0)=(G_2,X_2)$. In the second case (STB) or its (DL) case, choose any open compact subgroup $K_{0,p}$ contained in a \textbf{fixed} Bruhat-Tits stabilizer subgroup $K_{0,p}^*$ of $G_0(\bb{Q}_p)$ and any $K_0^p$ neat open compact in $G_0(\Ap)$ as before. \textbf{Fix} an admissible cone decomposition $\Sigma_0$ of $(G_0,X_0,K_0)$.\par 
The good toroidal compactification $\ca{S}_{K_0}^{\Sigma_0}$ of $\ca{S}_{K_0}$ over $\ca{O}_{E_0,(v)}$ for the fixed cone decomposition $\Sigma_0$ with associated integral models of boundary mixed Shimura varieties $\ca{S}_{K_{\Phi_0}}$ satisfying \cite[Thm. 4.1.5]{Mad19} and constructed from a Hodge embedding chosen as in \cite{DvHKZ24} (and (STB)), if it exists, is independent of the choice of such embedding.   
\end{cor}
\begin{proof}
    The argument above works for any isomorphism between Shimura data $(G_0,X_0)\to (G_0',X_0')$. Let $\gamma=1$. Since the integral models with deeper level structures are constructed by taking relative normalizations from the integral models defined with the fixed Bruhat-Tits level structures, we assume $K_{0,p}=K_{0,p}^*$. We see that $\ca{S}_{K_{\Phi_0}}\to \overline{\ca{S}}_{K_{\Phi_0}}\to \ca{S}_{K_{\Phi_0},h}$ is independent of the choice of the embedding as above. Combining this with \cite[Lem. A.3.4]{Mad19}, we see the uniqueness of such toroidal compactifications.
\end{proof}
\begin{rk}\label{rk-kisin-zhou-deleted-part}
In \cite[Prop. 4.3]{Mao25}, this result was also proved to study the change-of-parahoric morphisms of Hodge-type compactifications. We learned the construction of $(H'',X'')$ from the second version of \cite[Prop. 4.3.2]{KZ21}, but this part is now replaced by a new proof in the published version. The same argument above can extend the \emph{full} functoriality result in \emph{loc. cit.} to toroidal compactifications. We do not need to pursue this generality in this paper, but hope to record it in a future project with Shengkai Mao (see Remark \ref{rk-future-work}).
\end{rk}
\subsubsection{}\label{subsubsec-nstb}
Now we consider Case ($\mrm{STB}_n$). Suppose $K_{0,p}=\cap_{i=1}^n K^i_{0,p}$, where $K^i_{0,p}$ are Bruhat-Tits stabilizer subgroups of $G_0(\bb{Q}_p)$. As in the discussion under the third item ($\mrm{STB}_n$), 
we can choose a Hodge embedding $(G_0,X_0)\hookrightarrow (G^{\ddag,i,\prime},X^{\ddag,i,\prime})$ where $G^{\ddag,i,\prime}=\mrm{GSp}(V^i,\psi^i)$ and a self-dual lattice $V_\bb{Z}^i$ for each $i$, such that $K_{0,p}^i$ is exactly the stabilizer of $V_{\bb{Z}_p}^i$ in $G_0(\bb{Q}_p)$. Let $K_p^{\ddag,i,\prime}$ be the stabilizer in $G^{\ddag,i,\prime}(\bb{Q}_p)$ of $V_{\bb{Z}_p}^i$. Let $V_{\mrm{sum}}:=\oplus_{i=1}^nV^i$ and $V_{\mrm{sum},\bb{Z}}:=\oplus_{i=1}^n V_\bb{Z}^i$. Let $K_p^{\ddag,\prime}$ be the stabilizer of $V_{\mrm{sum},\bb{Z}_p}$ in $G^{\ddag,\prime}(\bb{Q}_p)$, where $G^{\ddag,\prime}:=\mrm{GSp}(V_{\mrm{sum}},\perp\psi^i)$. Let $X^{\ddag,\prime}$ be the union of Siegel upper and lower half-spaces induced by a $G^{\ddag,\prime}(\bb{R})$-conjugacy class of an element in $X_0$. 
Define $H':=\prod_{i=1}^n G^{\ddag,i,\prime}\times_{\prod_{i=1}^n\bb{G}_m}\bb{G}_m$, where the homomorphism from $\prod_{i=1}^nG^{\ddag,i,\prime}$ to $\prod_{i=1}^n\bb{G}_m$ is given by the product of similitude characters of $G^{\ddag,i,\prime}$, and the homomorphism from $\bb{G}_m$ to $\prod_{i=1}^n\bb{G}_m$ is given by the diagonal embedding. Note that the embedding $G_0\hookrightarrow \prod_{i=1}^n G^{\ddag,i,\prime}$ factors through $H'$: Indeed, $G_0$ is generated by $G_0^\der$ and $\omega(\bb{G}_m)$ (see Lemma \ref{injective}) and the restriction of the similitude character to $\omega(\bb{G}_m)$ is independent of the choice of Hodge embedding. 
Let $X_{H'}$ be the $H'(\bb{R})$-conjugacy class generated by an element in $X_0$ and the embedding $G_0\to H'$. The Shimura datum $(H',X_{H'})$ is of Hodge type. Moreover, there is a diagram
\begin{equation}\label{eq-emb-shimura-data}
    \begin{tikzcd}
    (G_0,X_0)\arrow[r,hook]& (H',X_{H'})\arrow[r,hook]\arrow[d,hook]& (\prod_{i=1}^nG^{\ddag,i,\prime},\prod_{i=1}^nX^{\ddag,i,\prime})\\
    &(G^{\ddag,\prime},X^{\ddag,\prime}).&
    \end{tikzcd}
\end{equation}
Let $G^{\ddag,i,\prime}_{\bb{Z}_p}$ (resp. $G^{\ddag,\prime}_{\bb{Z}_p}$) be a smooth reductive model of $G^{\ddag,i,\prime}_{\bb{Q}_p}$ (resp. $G^{\ddag,\prime}_{\bb{Q}_p}$) over $\bb{Z}_p$ such that $G^{\ddag,i,\prime}_{\bb{Z}_p}(\breve{\bb{Z}}_p)$ (resp. $G^{\ddag,\prime}_{\bb{Z}_p}(\breve{\bb{Z}}_p)$) is exactly the stabilizer of $V^i_{\breve{\bb{Z}}_p}$ (resp. $V_{\mrm{sum},\breve{\bb{Z}}_p}$).
By \cite[Prop. 1.7.6]{BT84} and \cite[Cor. 1.3]{PY06} (and the proof of \cite[Lem. 4.7]{KM15}), there are closed embeddings between reductive group schemes over $\bb{Z}_p$, $\prod_{i=1}^nG^{\ddag,i,\prime}_{\bb{Z}_p}\hookrightarrow \mrm{GL}(V_{\mrm{sum},\bb{Z}_p})$ and $G^{\ddag,\prime}_{\bb{Z}_p}\hookrightarrow \mrm{GL}(V_{\mrm{sum},\bb{Z}_p})$ extending $\prod_{i=1}^nG_{\bb{Q}_p}^{\ddag,i,\prime}\hookrightarrow \mrm{GL}(V_{\mrm{sum},\bb{Q}_p})$ and $G_{\bb{Q}_p}^{\ddag,\prime}\hookrightarrow \mrm{GL}(V_{\mrm{sum},\bb{Q}_p})$, respectively; moreover, the flat base change from $\zbkp$ to $\bb{Z}_p$ of the schematic closure $G^{\ddag,i,\prime}_{\zbkp}$ (resp. $G^{\ddag,\prime}_{\zbkp}$) of $G^{\ddag,i,\prime}$ (resp. $G^{\ddag,\prime}$) in $\mrm{GL}(V^i_{\zbkp})$ (resp. $\mrm{GL}(V_{\mrm{sum},\zbkp})$) is isomorphic to $G_{\bb{Z}_p}^{\ddag,i,\prime}$ (resp. $G_{\bb{Z}_p}^{\ddag,\prime}$) (cf. \cite[2.3.1 and 2.3.2]{Kis10}). 
Let $h$ be any element in $X_{H'}$. Its weight cocharacter $\omega_h$ sends $\bb{G}_m$ to $\omega_h(\bb{G}_m)$, and the intersection $C'$ of $\omega_h(\bb{G}_m)$ and $\prod_{i=1}^n G^{\ddag,i,\prime,\der}$ is isomorphic to $\mu_2$. 
The split torus $\omega_{h}(\bb{G}_m)$ canonically extends to a split torus $\omega_h(\bb{G}_m)_{\zbkp}$ over $\zbkp$. The finite group $C'$ extends to a finite flat group scheme $C'_{\zbkp}$ over $\zbkp$. So $H'_{\zbkp}:=\prod_{i=1}^n G^{\ddag,i,\prime,\der}_{\zbkp}\times \omega_h(\bb{G}_m)_{\zbkp}/C'_{\zbkp}$ is a reductive model over $\zbkp$ extending $H'$. By \cite[Cor. 1.3]{PY06} and the proof of \cite[Lem. 4.7]{KM15} again, there is a diagram of embeddings between reductive group schemes over $\zbkp$ induced by (\ref{eq-emb-shimura-data})
\begin{equation*}
  \begin{tikzcd}
H'_{\zbkp}\arrow[r,hook]\arrow[d,hook]&\prod_{i=1}^nG^{\ddag,i,\prime}_{\zbkp}\arrow[d,hook]\\
G^{\ddag,\prime}_{\zbkp}\arrow[r,hook]&\mrm{GL}(V_{\mrm{sum},\zbkp}).
  \end{tikzcd}  
\end{equation*}
Let $\ca{H}_{p}:=H'_{\zbkp}(\bb{Z}_p)$. Then the embedding $(H',X_{H'},\ca{H}_p)\hookrightarrow (G^{\ddag,\prime},X^{\ddag,\prime},K_p^{\ddag,\prime})$ is $p$-integral (see \cite[4.3.2]{Mad19}). 
Choose any cusp label representative $\Phi\in \ca{CLR}(H',X_{H'})$ mapping to $\Phi^{\ddag,\prime}\in\ca{CLR}(G^{\ddag,\prime},X^{\ddag,\prime})$ and to $\Phi^i\in\ca{CLR}(G^{\ddag,i,\prime},X^{\ddag,i,\prime})$. Choose a neat open compact subgroup $K^{\ddag,\prime,p}\sbst G^{\ddag,\prime}(\Ap)$ (resp. $K^{\ddag,i,\prime}$) stabilizing $V_{\mrm{sum},\zhp}$ (resp. $V_{\zhp}^i$) and choose a neat open compact subgroup $\ca{H}^p\sbst H'(\Ap)$ contained in both $K^{\ddag,\prime,p}$ and $\prod_{i=1}^nK^{\ddag,i,\prime}$. Define $\ca{H}:=\ca{H}_p\ca{H}^p$, $K^{\ddag,\prime}:=K^{\ddag,\prime}_pK^{\ddag,\prime,p}$ and $K^{\ddag,i,\prime}:=K^{\ddag,i,\prime,p}K_p^{\ddag,i,\prime}$. 
Let $\sh_{\ca{H}_{\Phi}}:=\sh_{P_\Phi(\A)\cap g_\Phi \ca{H}g_\Phi^{-1}}(P_\Phi,D_\Phi)$, and let $\ca{S}_{K_{\Phi^{\ddag,\prime}}}$ (resp. $\ca{S}_{K_{\Phi^i}}$) be the integral model corresponding to $\bm{\xi}_{(V_{\mrm{sum},\bb{Z}},\perp \psi^i_\bb{Z}),\Phi^\ddag,K^{\ddag,\prime,p}}$ (resp. $\bm{\xi}_{(V_\bb{Z}^i,\psi^i_\bb{Z}),\Phi^i,K^{\ddag,i,\prime,p}}$). 
Let $\ca{S}_{\ca{H}_{\Phi}}$ be the normalization in $\sh_{\ca{H}_\Phi}$ of $\ca{S}_{K_{\Phi^\ddag}}$. By \cite[Prop. 4.3.5]{Mad19}, the inverse limit $\ca{S}_{\ca{H}_{\Phi,p}}:=\varprojlim_{\ca{H}^p}\ca{S}_{\ca{H}_\Phi}$ has the extension property. 
By Proposition \ref{ext-imp} and \cite[Lem. 2.1.4]{Lov17}, $\ca{S}_{\ca{H}_{\Phi,p}}$ is canonically isomorphic to the inverse limit $\varprojlim_{\prod_{i=1}^nK^{\ddag,i,\prime,p}}\prod_{i=1}^n\ca{S}_{K_{\Phi^i}}$. Moreover, let $\ca{S}_{\ca{H}}$ be the normalization of $\ca{S}_{K^{\ddag,\prime}}$ in $\sh_{\ca{H}}$. By \cite{Kis10} and \cite{KM15}, $\ca{S}_{\ca{H}_p}:=\varprojlim_{\ca{H}^p}\ca{S}_{\ca{H}}$ has the extension property. 
Similarly, by Proposition \ref{ext-imp} and \cite[Lem. 2.1.4]{Lov17}, $\ca{S}_{\ca{H}_p}$ is an open and closed subscheme of $\varprojlim_{\prod^n_{i=1}K^{\ddag,i,\prime,p}}\prod_{i=1}^n\ca{S}_{K^{\ddag,i,\prime}}$. Hence, $\ca{S}_{\ca{H}_\Phi}$ is also the normalization in $\sh_{\ca{H}_\Phi}$ of $\prod_{i=1}^n\ca{S}_{K_{\Phi^i}}$, and $\ca{S}_{\ca{H}}$ is also the normalization in $\sh_{\ca{H}}$ of $\prod_{i=1}^n\ca{S}_{K^{\ddag,i,\prime}}$.\par
Hence, the discussion above implies the following lemma:
\begin{lem}\label{lem-n-stb-normalization-iso}
With the conventions above, let $\Phi_0'\in\ca{CLR}(G_0,X_0)$ be a cusp label representative mapping to $\Phi$. Choose any open compact subgroup $K_{0,p}'\sbst G_0(\bb{Q}_p)$ contained in $K^{\ddag,\prime}_p$ (and therefore in $\ca{H}_p$) and any neat open compact subgroup $K_0^{\prime,p}\sbst G_0(\Ap)$. Choose neat open compact subgroups $K^{\ddag,\prime,p}$ and $\prod_{i=1}^nK^{\ddag,i,\prime,p}$ containing $K_0^{\prime,p}$. Let $K_0':=K_{0,p}'K_{0}^{\prime,p}$. 
Then the normalizations in $\sh_{K_0'}(G_0,X_0)$ of $\ca{S}_{K^{\ddag,\prime}}$ and of $\prod_{i=1}^n\ca{S}_{K^{\ddag,i,\prime}}$ are isomorphic, and the normalizations in $\sh_{K_{\Phi_0'}}$ (resp. $\overline{\sh}_{K_{\Phi_0'}}$ and resp. $\sh_{K_{\Phi_0'},h}$) of $\ca{S}_{K_{\Phi^{\ddag,\prime}}}$ (resp. $\overline{\ca{S}}_{K_{\Phi^{\ddag,\prime}}}$ and resp. $\ca{S}_{K_{\Phi^{\ddag,\prime}},h}$) and of $\prod_{i=1}^n\ca{S}_{K_{\Phi^i}}$ (resp. $\prod_{i=1}^n\overline{\ca{S}}_{K_{\Phi^i}}$ and resp. $\prod_{i=1}^n\ca{S}_{K_{\Phi^i},h}$) are isomorphic.
\end{lem}
\begin{proof}
Choose $\ca{H}^p$ as above and containing $K_0^{\prime,p}$. By the discussion above, the normalizations in the statement are isomorphic to the normalization in $\sh_{K_0'}$ of $\ca{S}_{\ca{H}}$ and the normalization in $\sh_{K_{\Phi_0'}}$ (resp. $\overline{\sh}_{K_{\Phi_0'}}$ and resp. $\sh_{K_{\Phi_0'},h}$) of $\ca{S}_{\ca{H}_\Phi}$ (resp. $\overline{\ca{S}}_{\ca{H}_\Phi}$ and resp. $\ca{S}_{\ca{H}_\Phi,h}$), respectively.
\end{proof}
\subsubsection{}\label{subsubsec-stbndl}
Now we come back to the situation that $\gamma$ maps $\Phi_0$ to $\gamma\cdot \Phi_0:=(\gamma Q_0 \gamma^{-1},\gamma(X_0^+),\gamma g_0 \gamma^{-1})$ in Case ($\text{STB}_n$) and its (DL).\par
Suppose that $K_{0,p}$ is contained in $K_{0,p}'$, an intersection of $n$ Bruhat-Tits stabilizer subgroups of $G_0(\bb{Q}_p)$. Then $\lcj{K_{0,p}}{\gamma}\sbst \lcj{K_{0,p}'}{\gamma}$. 
We can now follow ($\text{STB}_n$) to choose a Hodge embedding $(G_0,\gamma\cdot X_0)\hookrightarrow(G^{\ddag,\gamma},X^{\ddag,\gamma})$ as follows:
Suppose $K_{0,p}'=\cap_{i=1}^n K^i_{0,p}$, where $K^i_{0,p}$ are Bruhat-Tits stabilizer subgroups of $G_0(\bb{Q}_p)$. As in the discussion under the third item ($\mrm{STB}_n$), 
we can choose a Hodge embedding $(G_0,\gamma\cdot X_0)\hookrightarrow (G^{\ddag,i,\gamma},X^{\ddag,i,\gamma})$ where $G^{\ddag,i,\gamma}=\mrm{GSp}(V^i,\psi^i)$ and a self-dual lattice $V_\bb{Z}^i$ for each $i$, such that $\lcj{K_{0,p}^i}{\gamma}$ is exactly the stabilizer of $V_{\bb{Z}_p}^i$ in $G_0(\bb{Q}_p)$. 
Let $V_{\mrm{sum}}:=\oplus_{i=1}^nV^i$ and $V_{\mrm{sum},\bb{Z}}:=\oplus_{i=1}^n V_\bb{Z}^i$. Let $K_p^{\ddag,\gamma}$ be the stabilizer of $V_{\mrm{sum},\bb{Z}_p}$ in $G^{\ddag,\gamma}(\bb{Q}_p)$, where $G^{\ddag,\gamma}:=\mrm{GSp}(V_{\mrm{sum}},\perp\psi^i)$.\par

Suppose $\gamma\cdot \Phi_0$ maps to a cusp label representative $\Phi^{\ddag,\gamma}\in \ca{CLR}(G^{\ddag,\gamma},X^{\ddag,\gamma})$. Let $K^{\ddag,\gamma}:=K_p^{\ddag,\gamma}K^{\ddag,\gamma,p}$ for suitable neat open compact $K^{\ddag,\gamma,p}\sbst G^{\ddag,\gamma}(\A)$ containing $\lcj{K_{0}^p}{\gamma}$. Denote $K_{\Phi^{\ddag,\gamma}}=g_{\Phi^{\ddag,\gamma}}K^{\ddag,\gamma}g_{\Phi^{\ddag,\gamma}}^{-1}\cap P_{\Phi^{\ddag,\gamma}}(\A)$.
We can still define an integral model $\ca{S}_{\lcj{K_{\Phi_0}}{\gamma}}$ (resp. $\overline{\ca{S}}_{\lcj{K_{\Phi_0}}{\gamma}}$ and resp. $\ca{S}_{\lcj{K_{\Phi_0},h}{\gamma}}$) of $\sh_{\lcj{K_{\Phi_0}}{\gamma}}$ (resp. $\overline{\sh}_{\lcj{K_{\Phi_0}}{\gamma}}$ and resp. $\sh_{\lcj{K_{\Phi_0},h}{\gamma}}$) as the normalization of $\ca{S}_{K_{\Phi^{\ddag,\gamma}}}$ (resp. $\overline{\ca{S}}_{K_{\Phi^{\ddag,\gamma}}}$ and resp. $\ca{S}_{K_{\Phi^{\ddag,\gamma}},h}$) in $\sh_{\lcj{K_{\Phi_0}}{\gamma}}$ (resp. $\overline{\sh}_{\lcj{K_{\Phi_0}}{\gamma}}$ and resp. $\sh_{\lcj{K_{\Phi_0},h}{\gamma}}$).
\begin{lem}\label{lem-gadq1-ext-n}
    Under the assumptions in ($\text{STB}_n$) or its (DL), we have the same result as Lemma \ref{lem-gadq1-ext}. The integral models we constructed depend only on the choice of the collection of Bruhat-Tits stabilizers $\{K^i_{0,p}\}_{i=1}^n$ and are independent of the choice of the collection of Hodge embeddings made above. 
\end{lem}
\begin{proof}
    Combine Lemma \ref{lem-gadq1-ext} with Lemma \ref{lem-n-stb-normalization-iso}.
\end{proof}
\subsubsection{}
Note that $\sh_K(G,X_b)(\bb{C})$ is covered by $\disju_{\alpha\in I_{G/G_0}} \sh_{K_0^\alpha}(G_0,X_0)(\bb{C})$.
We can write $\sh_K(G,X_b)_{\overline{\bb{Q}}}$ as
\begin{equation}\label{split-big-G}
\begin{split}
&\sh^+(G_0,X_0)_{\overline{\bb{Q}}}\times (\ag(G)/K)/\agsb(G_0)\\
&\iso \disju_{\alpha\in I_{G/G_0}} (\sh^+(G_0,X_0)_{\overline{\bb{Q}}}\times (\agsb(G)\pi(G_0(\A))g_\alpha K/K))/\agsb(G_0)\\
&\iso (\disju_{\alpha\in I_{G/G_0}} \sh_{K_0^\alpha}(G_0,X_0)_{\overline{\bb{Q}}}/\Delta(G_0,G)\cdot g_\alpha K)/K\\
&\iso  \disju_{\alpha\in I_{G/G_0}} (\sh_{K_0^\alpha}(G_0,X_0)_{\overline{\bb{Q}}}/\Delta(G_0,G))/\lcj{K}{g_\alpha}.
\end{split}
\end{equation}
In the equation above, $\Delta(G_0,G):=\ker (\ag(G_0)\to \ag(G))$. 
Moreover, (\ref{split-big-G}) is defined over $E_0$.\par
Now we construct the toroidal compactification $\ca{S}_{K}^\Sigma:=\ca{S}_K^\Sigma(G,X_b)$ as a disjoint union of quotients. \par
Denote $\lcj{K}{g_\alpha}=g_\alpha K g_\alpha^{-1}$ as $\lcj{K}{\alpha}$. Similarly, write $\lcj{K^p}{\alpha}=\lcj{K^p}{g_\alpha}$ and $\lcj{K_p}{\alpha}=\lcj{K_p}{g_\alpha}$. 
Denote $\Delta_{\lcj{K}{\alpha}}(G_0,G):=\ker(\ag(G_0)\to\ag(G)/\lcj{K}{\alpha})$.
Choose a neat open compact subgroup $K_0^{\alpha,p}$ for each $\alpha$ such that it is stabilized by $\Delta_{\lcj{K}{\alpha}}(G_0,G)$. This can be achieved because 
\begin{lem}\label{lem-deltak-cpt}
$\Delta_{\lcj{K}{\alpha}}(G_0,G)$ is compact.
\end{lem}
\begin{proof}
Without loss of generality, assume $\alpha$ is trivial. For any open compact $K_0\sbst G_0(\A)$, write $\ag(G_0)=\disju_{[g_0]\in\pi_0(\sh_{K_0,\bb{C}}(G_0,X_0))}\agsb(G_0)g_0 K_0$. Also, write $\ag(G)=\disju_{[g]\in \pi_0(\sh_{K,\bb{C}}(G,X_b))}\agsb(G)g K$. We now assume that $K_0$ maps to $K$ under $\pi$. We have decompositions
$$\ag(G_0)/K_0=\disju_{[g_0]\in\pi_0(\sh_{K_0,\bb{C}}(G_0,X_0))}\agsb(G_0)/\agsb(G_0)\cap\lcj{K_0}{g_0}$$
and
$$\ag(G)/K=\disju_{[g]\in\pi_0(\sh_{K,\bb{C}}(G,X_b))}\agsb(G)/\agsb(G)\cap\lcj{K}{g}.$$
Both $\agsb(G_0)\cap\lcj{K_0}{g_0}$ and $\agsb(G)\cap \lcj{K}{g}$ are arithmetic subgroups of $G^\ad(\bb{Q})^+$, and $\agsb(G_0)/\agsb(G_0)\cap\lcj{K_0}{g_0}\iso G^\ad(\bb{Q})^+/\agsb(G_0)\cap\lcj{K_0}{g_0}$ and $\agsb(G)/\agsb(G)\cap\lcj{K}{g}\iso G^\ad(\bb{Q})^+/\agsb(G)\cap\lcj{K}{g}$.
Then we see that the preimage of $[K]$ under $\ag(G_0)/K_0\to \ag(G)/K$ is finite.
\end{proof}
\begin{construction}\label{const-tor-Gb}\upshape
Recall that we choose $\Sigma_0^\alpha$ induced by $\Sigma$ and we can choose them to be projective as in Proposition \ref{zp-cones}. 
Let
\begin{equation}\label{eq-def-const-tor-Gb}
\ca{S}^\Sigma_{K}:=
\disju_{\alpha\in I_{G/G_0}} \ca{S}^{\Sigma_0^\alpha}_{K_0^\alpha}/\Delta_{\lcj{K}{\alpha}}(G_0,G).\end{equation}
We shall check that:
\begin{lem}\label{lem-const-tor-Gb}
In (\ref{eq-def-const-tor-Gb}) above, the action of $\Delta_{\lcj{K}{\alpha}}(G_0,G)$ on $\ca{S}^{\Sigma_0^\alpha}_{K_0^\alpha}$ and the quotient of it are well defined. Hence, (\ref{eq-def-const-tor-Gb}) defines a normal algebraic space that is proper over $\ca{O}'_{(v)}$, and is representable by a normal scheme which is projective over $\ca{O}'_{(v)}$ if $\Sigma_0^\alpha$ are induced by projective cone decompositions $\Sigma^{\ddag,\alpha}$ for all $\alpha\in I_{G/G_0}$. 
\end{lem}
\begin{proof}
By Lemma \ref{lem-gadq1-ext}, Lemma \ref{lem-gadq1-ext-n} and \cite[4.1.3 and 4.1.12]{Mad19}, there is an action of $\Delta_{\lcj{K}{\alpha}}:=\Delta_{\lcj{K}{\alpha}}(G_0,G)$ on the disjoint union of integral models $\ca{S}_{K_{\Phi_0^\alpha}}$ for all $\Phi_0^\alpha\in \ca{CLR}(G_0,X_0)$. Since the cone decompositions $\Sigma_0^\alpha$ are \emph{induced} by $\Sigma$ and since any element $d\in \Delta_{\lcj{K}{\alpha}}$ maps to $g_\alpha d'g_\alpha^{-1}\in\lcj{K}{\alpha}$, the action of $d'$ on $\Sigma$ is trivial by admissibility of $\Sigma$, and therefore the action of $\Delta_{\lcj{K}{\alpha}}$ on the integral models $\ca{S}_{K_{\Phi_0^\alpha}}$ induces an action of it on the disjoint union of strata $\ca{S}_{K_{\Phi_0^\alpha},\sigma_0^\alpha}\iso \ca{Z}_{\Upsilon_0^\alpha,K_0^\alpha}$ for all $\Upsilon_0^\alpha:=[(\Phi_0^\alpha,\sigma_0^\alpha)]\in \mrm{Cusp}_{K_0^\alpha}(G_0,X_0,\Sigma_0^\alpha)$. 
Finally, by \cite[Lem. A.3.4]{Mad19}, we have an action of $\Delta_{\lcj{K}{\alpha}}$ on $\ca{S}_{K_0^\alpha}^{\Sigma_0^\alpha}$ that extends the actions on the individual strata.\par
Next, we shall explain why this quotient makes sense. By Lemma \ref{lem-deltak-cpt} above, the action of $\Delta_{\lcj{K}{\alpha}}$ on $\sh_{K_0^\alpha}$ factors through a finite quotient, and therefore the same is true for $\ca{S}_{K_0^\alpha}^{\Sigma_0^\alpha}$ since $\sh_{K_0^\alpha}$ is dense in $\ca{S}_{K_0^\alpha}^{\Sigma_0^\alpha}$.\par 
Consequently, the quotients in (\ref{eq-def-const-tor-Gb}) exist as algebraic spaces by the Keel-Mori theorem \cite{KM97}. (See also \cite[Thm. 3.1.13]{CLO12} and \cite[Intro. p.17]{Knu71}.) 
The algebraic space $\ca{S}_K^\Sigma$ is separated over the base by \cite[Thm. 3.1.12(i)]{CLO12} (cf. \cite[\href{https://stacks.math.columbia.edu/tag/05Z2}{Lem. 05Z2}]{stacks-project}), and the morphism $\disju_{\alpha\in I_{G/G_0}}\ca{S}^{\Sigma_0^\alpha}_{K_0^\alpha}\to\ca{S}_K^\Sigma$ is finite by \cite[Thm. 3.1.13(i)(ii)]{CLO12} and \cite[\href{https://stacks.math.columbia.edu/tag/04NX}{Lem. 04NX}]{stacks-project}. 
Combining it with \cite[\href{https://stacks.math.columbia.edu/tag/08AJ}{Lem. 08AJ}]{stacks-project} and \cite[Thm. 3.1.13(iii)]{CLO12}, we have that the algebraic space $\ca{S}_K^\Sigma$ is proper over $\ca{O}'_{(v)}$. Finally, $\ca{S}_K^\Sigma$ is normal as we can check normality over {\'e}tale covers and combine \emph{loc. cit.} again.\par
For the projectivity, note that $\ca{S}_{K_0^\alpha}^{\Sigma_0^\alpha}$ will be projective if all $\Sigma_0^\alpha$'s for varying $\alpha$ satisfy the condition in the second sentence of the statement (see \cite[Ch. V., Thm. 5.8]{FC90} and \cite[Thm. 7.3.3.4]{Lan13}). Then $\ca{S}_K^\Sigma$ is projective over the base since the quotient of a projective scheme by a finite group is projective (see \cite[Ch. 4, Prop. 1.5]{Knu71}).
\end{proof}
The action of $\Delta_{\lcj{K}{\alpha}}(G_0,G)$ is trivial on $\sh_K^\Sigma$. So the finite morphism $\disju_{\alpha\in I_{G/G_0}}\sh_{K_0^\alpha}^{\Sigma_0^\alpha}\to \sh_{K}^\Sigma$ factors through the generic fiber of (\ref{eq-def-const-tor-Gb}). By Zariski's main theorem, the generic fiber of (\ref{eq-def-const-tor-Gb}) is isomorphic to $\sh_K^\Sigma$. Hence, (\ref{eq-def-const-tor-Gb}) is indeed an integral model
of $\sh_K^\Sigma$ over $\ca{O}_{(v)}'$.\hfill$\square$
\end{construction}
\subsubsection{}
Next, let us construct the strata for $\ca{S}_K^\Sigma$ over $\ca{O}_{(v)}'$.
Let $Q_0:=\pi^{-1}(Q)$ and $X^+_0:=\pi^{-1}(X^+_b)$. Up to an action of $Q_{\Phi}(\bb{Q})$, we can assume that $X^+_0$ is nonempty. 

Write $P_0:=P_{Q_0}$. There is a surjective map $P_0\to P_{\Phi}$. 
We can write $G(\bb{Q})_+ZP_\Phi(\A)g^b K$ as a disjoint union
$$\disju_{\alpha\in  I_{G/G_0}}G(\bb{Q})_+\pi(G_0(\A))\alpha K\cap G(\bb{Q})_+ZP_\Phi(\A)g^bK,$$
where $ I_{G/G_0}\iso G(\bb{Q})_+\pi(G_0(\A))\bss G(\A)/K$.\par
\begin{lem}\label{lem-db-cst-g}
    $G(\bb{Q})_+ZP_\Phi(\A)g^bK=G(\bb{Q})_+^{\overline{\ }}ZP_\Phi(\A)g^bK=$ 
    $$\disju_{\alpha\in I_{G/G_0}}\disju_{\pi(g_0^\alpha)\alpha\sim g^b} G(\bb{Q})_+\pi(P_0(\A)g_0^\alpha)g_\alpha K.$$
    In the equation above, the second index runs over a complete set of representatives $\{\pi(g_0^\alpha)\}$ in $\pi(G_0(\A))$ of $J_{G_0/P_0}^\alpha:=G(\bb{Q})_+^{\overline{\ }}\pi(P_0(\A))\bss G(\bb{Q})_+^{\overline{\ }} \pi(G_0(\A))g_\alpha K/K$ such that $\pi(g_0^\alpha)g_\alpha$ and $g^b$ are equivalent in $G(\bb{Q})_+^{\overline{\ }}ZP_\Phi(\A)\bss G(\A)/K$.
\end{lem}
\begin{proof}The first equation follows from the fact that $G(\bb{Q})_+^{\overline{\ }}\ca{H}=G(\bb{Q})_+\ca{H}$ for any open compact subgroup $\ca{H}$ in $G(\A)$. 
    Write $G_0(\A)=\disju_{g_0^\alpha}P_0(\A)g_0^\alpha K_0^\alpha$. Then, for any $g_0^\alpha$, $$G(\bb{Q})_+^{\overline{\ }}\pi(P_0(\A)g_0^\alpha K_0^\alpha)g_\alpha K=G(\bb{Q})_+^{\overline{\ }}\pi(P_0(\A)g_0^\alpha)g_\alpha K$$ 
is contained in $G(\bb{Q})_+^{\overline{\ }}ZP_\Phi(\A)g^bK$ if and only if $\pi(g_0^\alpha)g_\alpha\sim g^b$.
\end{proof}
Here we add a ``$\overline{\ }$'' to $G(\bb{Q})_+ZP_\Phi(\A)$ since $G(\bb{Q})_+^{\overline{\ }}ZP_\Phi(\A)$ is a group but $G(\bb{Q})_+ZP_\Phi(\A)$ is not.\par  
By Lemma \ref{lem-dq-quotient}, Lemma \ref{lem-db-cst-g} and \cite[Thm. 12.4]{Pin89}, we have that
\begin{prop}\label{prop-shimura-sur}
    With the notation as above, there are surjective maps 
    \begin{equation}\label{eq-sur-str-min}
 \disju_{\alpha\in{I}_{G/G_0}}\disju_{[\Phi_0^\alpha]=[(P_0,X_0^+,g_0^\alpha)]\xrightarrow{\pi(\alpha)}[ZP^b(\Phi)]}  \mrm{Z}_{[\Phi_0^\alpha],K_0^\alpha,\bb{C}}\lra 
 \mrm{Z}_{[ZP^b(\Phi)],K,\bb{C}},     
    \end{equation}
    \begin{equation}\label{eq-sur-sh-min}
 \disju_{\alpha\in{I}_{G/G_0}}\disju_{[\Phi_0^\alpha]=[(P_0,X_0^+,g_0^\alpha)]\xrightarrow{\pi(\alpha)}[ZP^b(\Phi)]}  \Delta_{\Phi_0^\alpha,K^\alpha_0}\bss\sh_{K_{\Phi_0^\alpha},h,\bb{C}}\lra \Delta^{ZP}_{\Phi,K}\bss\sh_{\K_{\Phi},h,\bb{C}},         
    \end{equation}
     \begin{equation}\label{eq-sur-str-tor}
\disju_{\alpha\in{I}_{G/G_0}}\disju_{[(\Phi_0^\alpha,\sigma_0^\alpha)]\xrightarrow{\pi(\alpha)}[ZP^b(\Phi,\sigma)]}  \mrm{Z}_{[(\Phi_0^\alpha,\sigma_0^\alpha)],K_0^\alpha,\bb{C}}\lra
\mrm{Z}_{[ZP^b(\Phi,\sigma)],K,\bb{C}},     
    \end{equation}
    \begin{equation}\label{eq-sur-sh-tor}
\disju_{\alpha\in{I}_{G/G_0}}\disju_{[(\Phi_0^\alpha,\sigma_0^\alpha)]\xrightarrow{\pi(\alpha)}[ZP^b(\Phi,\sigma)]}  \Delta_{\Phi_0^\alpha,K_0^\alpha}^\circ\bss\sh_{K_{\Phi_0^\alpha},\bb{C}}(\sigma_0^\alpha)\lra \Delta_{\Phi,K}^{ZP,\circ}\bss\sh_{\K_\Phi,\bb{C}}(\sigma),    
    \end{equation}
    and
    \begin{equation}\label{eq-sur-clstr-tor}
\disju_{\alpha\in I_{G/G_0}}\disju_{[(\Phi_0^\alpha,\sigma_0^\alpha)]\xrightarrow{\pi(\alpha)}[ZP^b(\Phi,\sigma)]}  \Delta_{\Phi_0^\alpha,K_0^\alpha}^\circ\bss\sh_{K_{\Phi_0^\alpha},\sigma_0^\alpha,\bb{C}}\lra \Delta_{\Phi,K}^{ZP,\circ}\bss\sh_{\K_\Phi,\sigma,\bb{C}},    
    \end{equation}
such that the images to right-hand sides of any two components in the disjoint unions of the left-hand sides of (\ref{eq-sur-str-min}) (\ref{eq-sur-sh-min}) (\ref{eq-sur-str-tor}), (\ref{eq-sur-sh-tor}) and (\ref{eq-sur-clstr-tor}) are disjoint or the same. All of the maps above algebraize and descend to the reflex field of $(G_0,X_0)$.
\end{prop}
\begin{construction}\label{const-strata-zpb}\upshape
By Lemma \ref{lem-db-cst-g} and Proposition \ref{prop-shimura-sur}, $\mrm{Z}_{[ZP^b({\Phi},\sigma)],K,\bb{C}}\iso \Delta^{ZP,\circ}_{\Phi,K}\bss \sh_{\wdtd{K}_\Phi,\sigma,\bb{C}}$ is isomorphic to the quotient 
$$ \disju_{\alpha\in I_{G/G_0}}\disju_{\pi(g_0^\alpha)\alpha\sim g^b} (\Delta_{\lcj{K}{\alpha}}(G_0,G)(G_0(\bb{Q})_+\bss G_0(\bb{Q})_+ \sh_{K_{\Phi_0^\alpha},\sigma_0^\alpha,\bb{C}}))/\Delta_{\lcj{K}{\alpha}}(G_0,G),$$
where $\sigma_0^\alpha$ is the unique cone in $\mbf{P}_{\Phi_0^\alpha}$ induced by the cone $\sigma$. 
Hence, we define an $\ca{O}'_{(v)}$-scheme associated with $[ZP^b(\Phi,\sigma)]$:
$$\ca{Z}_{[ZP^b(\Phi,\sigma)],K}:=  \disju_{\alpha\in I_{G/G_0}}\disju_{\pi(g_0^\alpha)\alpha\sim g^b}(\Delta_{\lcj{K}{\alpha}}(G_0,G) (G_0(\bb{Q})_+\bss G_0(\bb{Q})_+ \ca{S}_{K_{\Phi_0^\alpha},\sigma_0^\alpha}))/\Delta_{\lcj{K}{\alpha}}(G_0,G).$$ 
It is isomorphic to $$\disju_{\alpha\in I_{G/G_0}}\disju_{\pi(g_0^\alpha)\alpha\sim g^b}(\Delta_{\lcj{K}{\alpha}}(G_0,G) \Delta_{\Phi_0^\alpha,K_0^\alpha}^\circ\bss \ca{S}_{K_{\Phi_0^\alpha},\sigma_0^\alpha})/\Delta_{\lcj{K}{\alpha}}(G_0,G)$$ and to $$ 
\disju_{\alpha\in I_{G/G_0}}\disju_{\pi(g_0^\alpha)\alpha\sim g^b}(\Delta_{\lcj{K}{\alpha}}(G_0,G)\ca{Z}_{[(\Phi_0^\alpha,\sigma_0^\alpha)],K_0^\alpha})/\Delta_{\lcj{K}{\alpha}}(G_0,G)$$ by \cite[Thm. 4.1.5(5)]{Mad19}. When $K_0^\alpha$ are neat, $\Delta^\circ_{\Phi_0^\alpha,K_0^\alpha}$ are trivial. 
Since all the strata in the second disjoint union are locally closed in $\ca{S}_{K_0^\alpha}^{\Sigma_0^\alpha}$ and since the actions of $\Delta_{\lcj{K}{\alpha}}(G_0,G)$ on those strata and toroidal compactifications are equivariant, we see that the action of $\Delta_{\lcj{K}{\alpha}}(G_0,G)$ on those strata factors through a finite quotient, and that we have a finite morphism $\ca{Z}_{[ZP^b(\Phi,\sigma)],K}\hookrightarrow \ca{S}^{\Sigma}_K$. 
We shall check the following statement:
\begin{lem}\label{lem-strata-zpb}
    The induced map $i:\ca{Z}_{[ZP^b(\Phi,\sigma)],K}\to \ca{S}_K^\Sigma$ is a locally closed embedding. Thus we obtain a stratification on $\ca{S}_K^\Sigma$ by all $ZP$-cusp labels with cones.
\end{lem}
\begin{proof}
Denote by $[\ca{Z}]$ the image of $i$. (The image here is inductively defined, which means the schematic image defined in $\ca{S}_{K}^\Sigma$ removing the images defined by other $\ca{Z}_{[ZP^b(\Phi,\sigma)],K}$'s with strictly smaller dimensions. In particular, it is automatically locally closed.) The map $i:\ca{Z}_{[ZP^b(\Phi,\sigma)],K}\to[\ca{Z}]$ is finite, birational and a bijection over geometric points (see, e.g., \cite[Cor. A.7.2.2]{KM85}).
To check this map is an embedding, we can check over the complete local ring at every geometric point.
At any geometric point $\geom{y}$ of $\ca{Z}_{[(\Phi_0^\alpha,\sigma_0^\alpha)],K_0^\alpha}$, note that there is a structural morphism
$$\cpl{\ca{S}_{K_0^\alpha}^{\Sigma_0^\alpha}}{\geom{y}}\iso \cpl{\ca{S}_{K_{\Phi_0^\alpha}}(\sigma_0^\alpha)}{\geom{y}}\twoheadrightarrow\cpl{\ca{Z}_{[(\Phi_0^\alpha,\sigma_0^\alpha)],K_0^\alpha}}{\geom{y}}$$
whose pre-composition with the natural inclusion $\cpl{\ca{Z}_{[(\Phi_0^\alpha,\sigma_0^\alpha)],K_0^\alpha}}{\geom{y}}\hookrightarrow \cpl{\ca{S}_{K_0^\alpha}^{\Sigma_0^\alpha}}{\geom{y}}$ is the identity. The structural morphism is defined by the canonical projection of a twisted affine torus embedding with respect to a cone $\sigma_0^\alpha$ to its $\sigma_0^\alpha$-stratum.\par
Let $\geom{x}$ be any geometric point on $\ca{Z}_{[ZP^b(\Phi,\sigma)],K}$. Denote the image of $\geom{x}$ on $[\ca{Z}]$ by $[\geom{x}]$, and suppose $[\geom{x}]$ can be lifted to a geometric point on $\ca{Z}_{[(\Phi_0^\alpha,\sigma_0^\alpha)],K_0^\alpha}$. 
We then obtain a sequence
$$\cpl{\ca{Z}_{[ZP^b(\Phi,\sigma)],K}}{\geom{x}}\to\cpl{\ca{S}_K^\Sigma}{[\geom{x}]}\iso (\cpl{\ca{S}^{\Sigma_0^\alpha}_{K_0^\alpha}}{\pi^{b,-1}([\geom{x}])})^{\Delta_{\lcj{K}{\alpha}}(G_0,G)}\to \cpl{\ca{Z}_{[ZP^b(\Phi,\sigma)],K}}{\geom{x}}$$
whose composition is the identity. This induces an isomorphism $\cpl{\ca{Z}_{[ZP^b(\Phi,\sigma)],K}}{\geom{x}}\iso \cpl{[\ca{Z}]}{[\geom{x}]}$. The desired assertion that $\ca{Z}_{[ZP^b(\Phi,\sigma)],K}\iso [\ca{Z}]$ is now proved.
\end{proof}
Since the quotient map from the disjoint union of $\ca{Z}_{[(\Phi_0^\alpha,\sigma_0^\alpha)],K_0^\alpha}$ to $\ca{Z}_{[ZP^b(\Phi,\sigma)],K}$ is finite, by induction on the dimension of $\ca{Z}_{[ZP^b(\Phi,\sigma)],K}$ and by Lemma \ref{lem-strata-zpb}, we have a stratification of $\ca{S}_K^\Sigma$ by all $\ca{Z}_{[ZP^b(\Phi,\sigma)],K}$ extending that defined by $ZP$-cusp labels with cones on the generic fiber, and we have that those strata are normal, are locally closed in $\ca{S}_K^\Sigma$, and are flat over $\ca{O}'_{(v)}$. \par
We then define finer strata $\ca{Z}_{[(\Phi,\sigma)],K}$ as relative normalizations of $\ca{Z}_{[ZP^b(\Phi,\sigma)],K}$ in $\mrm{Z}_{[(\Phi,\sigma)],K}$. Since the embedding of  $\mrm{Z}_{[(\Phi,\sigma)],K}$ in $\mrm{Z}_{[ZP^b(\Phi,\sigma)],K}$ is open and closed, $\ca{Z}_{[(\Phi,\sigma)],K}$ is open and closed in $\ca{Z}_{[ZP^b(\Phi,\sigma)],K}$. Hence, $\ca{Z}_{[(\Phi,\sigma)],K}$ are locally closed in $\ca{S}_K^{\Sigma}$. Moreover, by construction, we have that
$$\ca{S}_K^{\Sigma}\iso\disju_{[ZP^b(\Phi,\sigma)]\in \mrm{Cusp}_K^{ZP}(G,X_b,\Sigma)}\ca{Z}_{[ZP^b(\Phi,\sigma)],K}\iso \disju_{\Upsilon=[(\Phi,\sigma)]\in\mrm{Cusp}_K(G,X_b,\Sigma)} \ca{Z}_{[(\Phi,\sigma)],K},$$
which gives us stratifications of $\ca{S}_K^{\Sigma}$ given by $ZP$-cusp labels and cusp labels with cones, extending the stratifications of both types of cusp labels with cones over the generic fiber.\hfill$\square$
\end{construction}
\subsubsection{}
Now we are ready to find an integral model of $\sh_{\wdtd{K}_\Phi}:=\sh_{\K_\Phi^b}(ZP^b_\Phi,ZP^b_\Phi(\bb{Q})D_\Phi)$.
Set $$I_{ZP_{\Phi}/P_{0},K}:=\stb_{ZP_{\Phi}(\bb{Q})}(D_{Q_0,X^+_0})\pi(P_0(\A))\bss ZP_{\Phi}(\A)g^bK/K.$$ 
Since $P_{0}(\bb{Q})$ acts on connected components of $D_{Q_0,X^+_0}$ transitively, we have $$I_{ZP_{\Phi}/P_{0},K}\iso \stb_{ZP_{\Phi}(\bb{Q})}(D_{Q_0,X_0^+}^+)\pi(P_{0}(\A))\bss ZP_{\Phi}(\A)g^bK/K.$$
For any $\beta\in I_{\Phi}:=I_{ZP_{\Phi}/P_{0},K}$, we consider the image of $\beta$ in $$J_{\Phi}:=G(\bb{Q})_+^{\overline{\ }}\pi(P_{0}(\A))\bss G(\bb{Q})_+^{\overline{\ }}ZP_{\Phi}(\A)g^b K/K.$$ 
By Proposition \ref{prop-shimura-sur}, for any $\beta$, we choose any representative $e_\beta\in ZP_\Phi(\A)g^bK$ for $\beta$, and there is a $\gamma_\beta\in G(\bb{Q})_+$ and a cusp label representative $\Phi_\beta=(Q_0^{\gamma_\beta},X^+_0,g_\beta)$ in $\ca{CLR}(G_0,X_0)$ such that $\gamma_\beta g_\beta g_\alpha=e_\beta k$ for some $k\in K$,
$$\sh_{\wdtd{K}_\Phi}\iso \disju_{\beta\in I_{\Phi}} \gamma_\beta\cdot \sh_{K_{\Phi_\beta}}/\Delta(P_0,ZP_{\Phi}/ZP_{\Phi}(\A)\cap \lcj{K}{\gamma_\beta g_\beta g_\alpha}),$$
where $\Delta(P_0,ZP_{\Phi}/ZP_{\Phi}(\A)\cap \lcj{K}{\gamma_\beta g_\beta g_\alpha}):=\ker (\ag(P_0)\to \ag(ZP_{\Phi})/ZP_{\Phi}(\A)\cap \lcj{K}{\gamma_\beta g_\beta g_\alpha})$ and $g_\alpha$ is determined by the image of $\beta$ in $I_{G/G_0}$. Set $K^\beta:=ZP_{\Phi}(\A)\cap \lcj{K}{\gamma_\beta g_\beta g_\alpha}$ and $\Delta_{K^\beta}(P_0,ZP_\Phi):=\Delta(P_0,ZP_{\Phi}/ZP_{\Phi}(\A)\cap \lcj{K}{\gamma_\beta g_\beta g_\alpha})$.\par
\begin{lem}\label{lem-deltakp-cpt}
$\Delta_{\rcj{K}{\beta}}(P_0,ZP_\Phi)$ is compact.
\end{lem}
\begin{proof}
For any open compact subgroup $K_{\Phi_0}\sbst P_0(\A)$ mapping to $K^\beta$ under $\pi$, we have decompositions
$$\ag(P_0)/K_{\Phi_0}=\disju_{[p_0]\in P_0(\bb{Q})_+\bss P_0(\A)/K_{\Phi_0}}\agsb(P_0)/\agsb(P_0)\cap\lcj{K_{\Phi_0}}{p_0}$$
and
$$\ag(ZP_\Phi)/K^\beta=\disju_{[p]\in ZP_\Phi(\bb{Q})_+\bss ZP_\Phi(\A)/K^\beta}\agsb(ZP_\Phi)/\agsb(ZP_\Phi)\cap\lcj{(K^\beta)}{p}.$$
Both $\agsb(P_0)\cap\lcj{K_{\Phi_0}}{p_0}$ and $\agsb(ZP_\Phi)\cap \lcj{(K^\beta)}{p}$ are arithmetic subgroups of $P_0^\ad(\bb{Q})^+$, and $\agsb(P_0)/\agsb(P_0)\cap\lcj{K_{\Phi_0}}{p_0}\iso P_0^\ad(\bb{Q})^+/\agsb(P_0)\cap\lcj{K_{\Phi_0}}{p_0}$ and $\agsb(ZP_\Phi)/\agsb(ZP_\Phi)\cap\lcj{(K^\beta)}{p}\iso P_0^\ad(\bb{Q})^+/\agsb(ZP_\Phi)\cap\lcj{(K^\beta)}{p}$.
Then we see that the preimage of $[K^\beta]$ under $\ag(P_0)/K_{\Phi_0}\to \ag(G)/K^\beta$ is finite.
\end{proof}
Hence, we can write $\sh_{\wdtd{K}_\Phi}$ as 
$$\disju_{\beta\in I_{\Phi}} \gamma_\beta\cdot \sh_{K_{\Phi_\beta}}/\Delta_{K^\beta}(P_0,ZP_{\Phi}).$$

\begin{construction}\label{const-zp-mix-sh}\upshape
    Set $\ca{S}_{\K_\Phi}:=\ca{S}_{\wdtd{K}^b_\Phi}(ZP_\Phi^b,ZP^b_\Phi(\bb{Q})D_{\Phi}):=$
\begin{equation*}
\disju_{\beta\in I_{\Phi}}\gamma_\beta\cdot \ca{S}_{K_{\Phi_\beta}}/\Delta_{K^\beta}(P_0,ZP_{\Phi})
\end{equation*}
The quotient makes sense since the action of $\Delta_{K^\beta}(P_0,ZP_\Phi)$ makes sense by Lemma \ref{lem-gadq1-ext}, Lemma \ref{lem-gadq1-ext-n}, and \cite[4.1.3]{Mad19}, and factors through a finite group by Lemma \ref{lem-deltakp-cpt}.\par
Similarly, we construct $$\overline{\ca{S}}_{\wdtd{K}_\Phi}:=\overline{\ca{S}}_{\K_\Phi^b}(ZP^b_\Phi,ZP^b_\Phi(\bb{Q})D_{\Phi_b}):=\disju_{\beta\in I_{\Phi}}\gamma_\beta\cdot \overline{\ca{S}}_{K_{\Phi_\beta}}/\Delta_{K^\beta}(P_0,ZP_{\Phi}),$$ 
and $$\ca{S}_{\wdtd{K}_{\Phi},h}:=\ca{S}_{\K_\Phi^b,h}(ZP^b_\Phi,ZP^b_\Phi(\bb{Q})D_{\Phi_b}):=\disju_{\beta\in I_{\Phi}}\gamma_\beta\cdot \ca{S}_{K_{\Phi_\beta},h}/\Delta_{K^\beta}(P_0,ZP_{\Phi}).$$ Consequently, there is a tower of $\ca{O}'_{(v)}$-schemes 
$$\ca{S}_{\wdtd{K}_\Phi}\lra \overline{\ca{S}}_{\wdtd{K}_\Phi} \lra \ca{S}_{\wdtd{K}_{\Phi},h}.$$
By the construction above and by \cite[4.1.5(4)]{Mad19} (cf. Theorem \ref{thm-tor-hdg} (3)), the first morphism above is affine, while the second morphism is proper.\hfill$\square$
\end{construction}
Denote $\Phi_{\beta'}:=\gamma_\beta\cdot \Phi_\beta$.\par
Now we shall check the following statement:
\begin{prop}\label{prop-mixsh-independent}
In Cases (HS), ($\text{STB}_n$) and (DL), $\ca{S}_{\wdtd{K}_\Phi}$, $\overline{\ca{S}}_{\wdtd{K}_\Phi}$ and $\ca{S}_{\wdtd{K}_{\Phi},h}$ are independent of the choice of representatives $e_\beta$ and the pairs $(\gamma_\beta,g_\beta)$ for all $\beta$ in the construction.\par 
Moreover, for any $\gamma\in G(\bb{Q})$ and $q'\in ZP_{\Phi_2}(\A)$ such that $\Phi_1\sim_{ZP}\Phi_2$ under the equivalence $\Phi_1\xrightarrow[\sim]{(\gamma,q')_K}\Phi_2$, the isomorphism of mixed Shimura varieties induced by this equivalence extends to an isomorphism between integral models $(\gamma,q'):\ca{S}_{\wdtd{K}_{\Phi_1}}\xrightarrow{\sim}\ca{S}_{\wdtd{K}_{\Phi_2}}$.\par 
The statement in the last paragraph also holds if we replace $\ca{S}_?$ with $\overline{\ca{S}}_?$ and $\ca{S}_{?,h}$. The isomorphisms induced by $(\gamma,q')$ are compatible with the towers $\ca{S}_{\wdtd{K}_\Phi}\to \overline{\ca{S}}_{\wdtd{K}_\Phi} \to \ca{S}_{\wdtd{K}_{\Phi},h}$.
\end{prop}
\begin{proof} By Lemma \ref{lem-gadq1-ext}, Lemma \ref{lem-gadq1-ext-n} and \cite[4.1.3]{Mad19}, different choices of $(e_\beta;\gamma_\beta,g_\beta)$ for the same $\beta$ will induce isomorphisms between the models $\ca{S}_{K_{\Phi_{\beta'}}}$ constructed as above. Now we show the second paragraph. \par
By Lemma \ref{lem-gadq1-ext-n} and the constructions above, the statements hold for $(\gamma,1)$ since the morphism $\gamma:\ca{S}_{K_{\beta'}}\to \ca{S}_{\lcj{K_{\beta'}}{\gamma}}$ is well defined. We only need to show the case where $\gamma=1$ and $q'$ is non-trivial.
In fact, we can directly define a morphism of $q'$; if we can do so, then the statement follows from separatedness and normality of the schemes involved. 
Suppose $\Phi\xrightarrow[\sim]{(1,q')_K}\Phi_1$. 
Fix a complete set of representatives $\{e_\beta\}$ for $I_{ZP_\Phi/P_{0}}$. For any $q'\in ZP_\Phi(\A)$, $\{e_{\beta}^{q'}:=e_\beta q'\}$ is a complete set of representatives for $ \stb_{ZP_\Phi(\bb{Q})}(D_{Q_0,X_0^+})\pi(P_{0}(\A))\bss ZP_\Phi(\A)(q')^{-1}g^b K/K$. We still assign $(\gamma_\beta,g_\beta)$ to $e_\beta^{q'}$ and construct $\ca{S}_{\K_{\Phi_1}}$ as above. By the last paragraph, the construction of $\ca{S}_{\K_{\Phi_1}}$ is independent of the choice of such assignment. 
Then there is an isomorphism $q': \ca{S}_{\K_\Phi}\xrightarrow{\sim} \ca{S}_{\K_{\Phi_1}}$ induced by the map $e_\beta\mapsto e_\beta^{q'}$ extending the isomorphism over the generic fiber sending the scheme $\ca{S}_{K_{\Phi_{\beta'}}}$ constructed for the index $e_\beta$ identically to the scheme constructed for the index $e_\beta^{q'}$. 
\end{proof}
\subsubsection{}\label{subsubsec-intmod-study} Let us study the integral models we just constructed. \par
For any admissible cone decomposition $\Sigma_0$ of $(G_0,X_0,K_0^\alpha)$ and $\Phi_{\beta'}$, define $\Delta_{K^\beta}^\circ(P_0,G)$ to be 
the group consisting of the elements $(h,\gamma^{-1})\in\Delta_{K^\beta}(P_0,G):=\ker(\ag(P_0)\times^{P_0^\ad(\bb{Q})^+}Q_0^\ad(\bb{Q})^+\to \ag(G)/\lcj{K}{\gamma_\beta g_\beta g_\alpha})$ such that $\gamma$ stabilizes $\Sigma(\Phi_{\beta'})$, where we denote $Q_0^\ad(\bb{Q})^+:=Q^\ad_0(\bb{Q})\cap G^\ad_0(\bb{Q})^+$.\par 
In fact, $\Delta_\beta:=\Delta_{K^\beta}(P_0,G)$ is a generalization of the group $\Delta_{\Phi,K}$ defined by Pink (see \cite[6.3 and 7.3]{Pin89}). To explain this similarity, we convert the right action of $(h,\gamma^{-1})\in \Delta_{\beta}$ to $[\gamma h\gamma^{-1}]\circ [\gamma]$, where $[\gamma]$ is the left action of $\gamma$ and $[\gamma h\gamma^{-1}]$ is the right action of $\gamma h\gamma^{-1}$. 
Let \begin{equation*}
\begin{split}
&\wdtd{\Delta}_\beta:=\{(h_1,\gamma_1)\in P_0(\A)\times Q^\ad_\Phi(\bb{Q})^+|\gamma_1\text{ is lifted to an element }\gamma^*_1\in Q_\Phi(\bb{Q})_+,\\
&\text{ and }\gamma_1^*=h_1\text{ modulo }\lcj{K}{\gamma_\beta g_\beta g_\alpha}\text{ on the right in }G(\A)\},
\end{split}
\end{equation*}
where $Q_\Phi^\ad(\bb{Q})^+:=Q_\Phi^\ad(\bb{Q})\cap G^\ad(\bb{Q})^+=Q_0^\ad(\bb{Q})^+$.
We have that $(h,\gamma^{-1})\in\Delta_\beta$ if and only if there is a lifting $\gamma_1^*$ of $\gamma$ to $Q_\Phi(\bb{Q})_+$ such that $(h,\gamma^{-1})(\gamma_1^{*,-1},\gamma)=(h\gamma^{-1}\gamma_1^{*,-1}\gamma,1)=(h\gamma_1^{*,-1},1)=(k,1)$ for some $k\in \lcj{K}{\gamma_\beta g_\beta g_\alpha}$. We then have $(h,\gamma^{-1})\in \Delta_\beta$ if and only if $\gamma_1^*=\gamma h\gamma^{-1}$ modulo $\lcj{K}{\gamma_\beta g_\beta g_\alpha}$. As a result, there is a surjective homomorphism 
$$\wdtd{\Delta}_\beta\lra \Delta_\beta$$
with the assignment $(\gamma h\gamma^{-1},\gamma)\mapsto (h,\gamma^{-1})$, which converts the right action of $\Delta_\beta$ to $[\gamma h\gamma^{-1}]\circ[\gamma]$. Note that the action of $\wdtd{\Delta}_\beta$ modulo $P_0(\bb{Q})$ coincides with the action of Pink's group $\Delta_{\Phi,K}$ when $(G_0,X_0)=(G,X_b)$, $P_0=P_\Phi$ and $\gamma_\beta=g_\alpha=1$. \par
By \cite[Thm. 6.19(a)]{Pin89} and \cite[2.1.19]{Mad19}, the definition of $\Delta_{K^\beta}^\circ(P_0,G)$ is independent of the choice of the cone decomposition $\Sigma_0$.
\begin{lem}\label{lem-twt-free-action}
Recall that we assume $K^p$ is \textbf{neat}. With the conventions and assumptions above, and in all cases, 
\begin{enumerate}
\item\label{lem-twt-free-1} 
We have $\Delta^\circ_{K^\beta}(P_0,G)=\Delta_{K^\beta}(P_0,ZP_\Phi)$. 
\item\label{lem-twt-free-2} The group $\Delta_{K^\beta}(P_0,ZP_{\Phi})$ acts equivariantly on the tower of schemes $\ca{S}_{K_{\Phi_{\beta'}}}\to \overline{\ca{S}}_{K_{\Phi_{\beta'}}}\to \ca{S}_{K_{\Phi_{\beta'}},h}$; the quotient induced by this action on each scheme in the tower is finite {\'e}tale. More precisely, the action of $\Delta_{K^\beta}(P_0,ZP_\Phi)$ on $\ca{S}_{K_{\Phi_{\beta'}}}$ (resp. $\overline{\ca{S}}_{K_{\Phi_{\beta'}}}$ and resp. $\ca{S}_{K_{\Phi_{\beta'}},h}$) factors through some finite group $H_{\beta}$ (resp. $\overline{H}_\beta$ and resp. $H_{\beta,h}$), and the action of this finite group is free. Consequently, for $K^{\prime,p}$ any normal subgroup of $K^p$, we have a tower $\ca{S}_{\K_{\Phi}'}\to \overline{\ca{S}}_{\K_\Phi'}\to \ca{S}_{\K_{\Phi}',h}$ associated with $\Phi$ and $K'=K_pK^{\prime,p}$ by Construction \ref{const-zp-mix-sh}. The quotient of each scheme in the tower by $\K_\Phi^p$ is finite {\'e}tale.
\item\label{lem-twt-free-3} 
The action of $\Delta_{K^\beta}(P_0,G)$ on $\ca{S}_{K_{\Phi_{\beta'}},h}$ also factors through a finite group, and the action of it on $\ca{S}_{K_{\Phi_{\beta'}},h}$ is also free.
\end{enumerate}

\end{lem}
\begin{proof}
First, we check the actions factor through finite groups. 
To check this action factors through a finite group, we argue as in \cite{Pin89}. Indeed, it suffices to check over the generic fiber. The action of $(\gamma h\gamma^{-1},\gamma)\in \wdtd{\Delta}_\beta$ is determined by $\gamma$. There is an almost semi-direct product $Q_0\iso G^*_h\cdot P_0$ with $G^*_h$ a reductive group over $\bb{Q}$ which is a finite cover of $Q_0/P_0$. The group $G^*_h(\bb{Q})$ acts on $\sh_{K_{\Phi_{\beta'}},h}$ trivially. Hence, the action of $\Delta_{\beta}$ on $\ca{S}_{K_{\Phi_{\beta'}},h}$ factors through a finite quotient group by projecting $\wdtd{\Delta}_\beta$ to $Q_0/P_0(\bb{Q})$. From Construction \ref{const-zp-mix-sh}, we see that the action of $\Delta^\circ_{K^\beta}(P_0,G)$ also factors through a finite group if Part \ref{lem-twt-free-1} is shown. \par
By Lemma \ref{lem-pin214-str}, the group $\Delta_{K^\beta}(P_0,ZP_\Phi)$ is included in $\Delta^\circ_{K^\beta}(P_0,G)$. The projection of $\Delta^\circ_{K^\beta}(P_0,G)$ to $Q_0^\ad/ZP_\Phi^\ad(\bb{Q})$ is an arithmetic subgroup contained in the image of $\lcj{K}{\gamma_\beta g_\beta g_\alpha}$. Moreover, this projection is contained in $P_\Phi''/ZP_\Phi(\bb{Q})$ (see \S\ref{subsubsec-Harris-group} for the notation) by the neatness of $K$ and by \cite[Thm. 6.19(a)]{Pin89}. Then the $Q_\Phi(\bb{Q})_+$-component of $\Delta^\circ_{K^\beta}(P_0,G)$ centralizes $U_\Phi$ and is contained in $ZP_\Phi(\bb{Q})$ by Lemma \ref{lem-Harris-group-compare}. This proves Part \ref{lem-twt-free-1}.\par
Next, we show the freeness of actions. 
Let $\overline{x}\in \ca{S}_{K_{\Phi_{\beta'}}}(\kappa)$ be a geometric point of $\ca{S}_{K_{\Phi_{\beta'},p}}$, where 
$\kappa$ is an algebraically closed field. Suppose that the action of $\tau=(h,\gamma^{-1})\in\Delta_{K^\beta}(P_0,ZP_\Phi)$ on $\overline{x}$ is trivial. The aim is to show that $\tau$ acts on every point of $\ca{S}_{K_{\Phi_{\beta'}}}$ trivially.
Then $\tau$ stabilizes the component $\lcj{\ca{S}}{j}_{K_{\Phi_{\beta'},p}K_{\Phi_{\beta'}}^p}$ containing $\overline{x}$. Up to replacing $K_{\Phi_{\beta'},p}$ and $K_p^\beta$ with the right conjugate of a representative in $P_0(\bb{Q}_p)$ of $j$, it suffices to assume that $j=1$ and to consider the subgroup of $\wdtd{\Delta}_\beta$ that stabilizes $\lcj{\ca{S}}{1}_{K_{\Phi_{\beta'},p}K_{\Phi_{\beta'}}^p}$. 
This stabilizer is 
\begin{equation*}
\begin{split}
&\wdtd{\Delta}_\beta^1:=\{(h_1,\gamma_1)\in P_0(\bb{Q})_+P_0(\Ap)K_{\Phi_{\beta'},p}\times Q^\ad(\bb{Q})^+|\gamma_1\text{ is lifted to an element }\gamma^*_1\in Q_\Phi(\bb{Q})_+,\\
&\text{ and }\gamma_1^*=h_1\text{ modulo }\lcj{K}{\gamma_\beta g_\beta g_\alpha}\text{ on the right in }G(\A)\}.
\end{split}
\end{equation*}
Modulo the trivial action of $P_0(\bb{Q})_+$, we can assume that $\gamma_1$'s in the second component of elements in $\wdtd{\Delta}^1_\beta$ are all lifted to $Q_\Phi(\bb{Q})_+\cap \lcj{K_p}{\gamma_\beta g_\beta g_\alpha}$.
Since we can assume $(h,\gamma^{-1})\in \Delta_{K^\beta}(P_0,ZP_\Phi)$, we assume that $\gamma$ can be lifted to $ZP_\Phi(\bb{Q})_+\cap \lcj{K_p}{\gamma_\beta g_\beta g_\alpha}$. By Lemma \ref{lem-twist-admissible} below, the action of the last group can be described by twisting constructions for some $G_{0,\zbkp}$ as in Section \ref{sec-twt-construction}.\par 
Write $G$ as an almost-direct product $G\iso (G_0/\ker(\pi^b)\times S_G)/Z_G'$, where $S_G\sbst Z_G$ and $Z_G'$ is a finite group. Then Lemma \ref{lem-twist-admissible} shows that $\gamma$ lifts to $P_0^\dagger(F)\cap G_{0,\zbkp}(\ca{O}_{F,(p)})$, where $P_0^\dagger:=\pi^{b,-1}( Z_G')\cdot P_0$. \par
In all cases, we can pull back the $1$-motive with additional structures to $\ca{S}_{K_{\Phi_{\beta'}}}$ from the integral model of some Siegel-type boundary mixed Shimura variety. 
Suppose that the image of $\beta$ is some $\alpha\in I_{G/G_0}$. Suppose that $\Phi_{\beta'}$ maps to $\Phi^{\ddag,\prime}\in \ca{CLR}(G^\ddag,X^\ddag)$ as in the notation of \S\ref{subsubsec-hodge-gad1q} and \S\ref{subsubsec-nstb}. With the notation in the above-mentioned subsections, in Cases (HS), (STB) and their deep-level cases, we pull back the tautological $1$-motive with additional structures over $\ca{S}_{K_{\Phi^{\ddag,\prime}}}$. In Case ($\text{STB}_n$) and its deep-level case, we choose an $i$ and pull back the tautological family over $\ca{S}_{K_{\Phi^i}}$.\par
By Proposition \ref{prop-twt-1}, there is a commutative diagram 
\begin{equation}
    \begin{tikzcd}
    V_{\zhp}\otimes\ca{O}_F\arrow[r,"\wdtd{\gamma}^{-1}\circ h"]\arrow[dr,"k\cdot"]&V_{\Ap}\otimes_{\zbkp} \ca{O}_{F,(p)}\arrow[r,"u_{\overline{x}}"]& V^p\Q_{\overline{x}}^{\ca{O}_F}\arrow[r,"V^pf"]&V^p\Q_{\overline{x}}^{\tau,\ca{O}_F}\\
   & V_{\Ap}\otimes_{\zbkp} \ca{O}_{F,(p)}\arrow[r,"u_{\overline{x}}"] & V^p\Q_{\overline{x}}^{\ca{O}_F}\arrow[ur,"V^ps"]&
    \end{tikzcd}
\end{equation}
where $s$ is induced by an isomorphism between $1$-motives $s:\Q_{\overline{x}}\to \Q^\tau_{\overline{x}}$ preserving polarizations, $k\in K_{\Phi_{\beta'}}^p$, $\tg\in P_0^\dagger(F)\cap G_{0,\zbkp}(\ca{O}_{F,(p)})$ is a lifting of $\gamma$ and $h\in P_0(\Ap)$.\par
Since $h\tg^{-1}$ is in the preimage in $G_0(\A)$ of the image of $\lcj{K}{\gamma_\beta g_\beta g_\alpha}$ in $G/S_{G}(\A)$, $(\tg)^{-1}hk^{-1}\in G_0(\Ap\otimes_\bb{Q} F)$ lies in an open compact subgroup of $G_0(\Ap\otimes_\bb{Q} F)$. 
We replace $V_{\zhp}\otimes \ca{O}_F$ with another $\zhp$-lattice $V_{\zhp}'\sbst V_{\Ap}\otimes_{\zbkp} \ca{O}_{F,(p)}$ and replace $\Q^{\oo}_{\overline{x}}$ with $\Q'$ with a $\zbkpt$-isogeny $\Q'\to \Q_{\overline{x}}^{\oo}$ such that the commutative diagram above induces a commutative diagram preserving polarizations
\begin{equation}\label{diag-auto}
    \begin{tikzcd}
    V_{\zhp}'\arrow[rr,"u_{\overline{x}}"]\arrow[d,"kh^{-1}\tg"]&& T^p\Q'\arrow[d]\\
    V_{\zhp}'\arrow[rr,"u_{\overline{x}}"]&& T^p\Q'.
    \end{tikzcd}
\end{equation}
The automorphism on the right is induced by an automorphism in $\mrm{Aut} \Q'$. Note that we have 
$$\mrm{Aut}(\Q')\sbst \mrm{Aut}(Y_{\Q'})\times \mrm{Aut}(T_{\Q'})\times \mrm{Aut}(A_{\Q'}).$$
Indeed, we can check that $\mrm{Aut}(\Q')\sbst \mrm{Aut}(Y_{\Q'})\times \mrm{Aut}(\G^{\natural}_{\Q'})$ by definition; and $\mrm{Aut}(\G^{\natural}_{\Q'})\iso \mrm{Aut}(((\Q')^\vee)^\circ)\sbst \mrm{Aut}(X_{\Q'})\times \mrm{Aut}(A_{\Q'}^\vee)\iso \mrm{Aut}(T_{\Q'})\times \mrm{Aut}(A_{\Q'})$. Moreover, since the projection to $\mrm{Aut}(A_{\Q'})\to \mrm{Aut}(T^p A_{\Q'})$ of the right vertical arrow of (\ref{diag-auto}) is induced by an isomorphism of polarized abelian schemes with level structures, this projection is some multiplication of roots of unity (see, e.g., \cite[Cor. 2.3.3.2]{Lan13}). We will use these facts later.\par
Set $ZP_0:=Z_{G_0}P_0$ and $W_0$ the unipotent radical of $P_0$. Set $ZW_0:=Z_{G_0}W_0$. We claim that the image of $q:=k^{-1}h\tg^{-1}$ in $ZP_0/ZW_0(\Ap\otimes_\bb{Q} F)$ is trivial. Then it will follow that $q\in ZW_0(\Ap\otimes_\bb{Q}F)$. Let us show the claim.\par
Let $P_\Phi^\dagger$ be the quotient of $P_0^\dagger$ by $\bb{K}:=\ker(G_0\to G)$. The projection of $q$ to $P_\Phi^\dagger/P_\Phi(\Ap\otimes_\bb{Q}F)$ via $P_0^\dagger\to P_\Phi^\dagger\to P_\Phi^\dagger/P_\Phi$ lies in $P_\Phi^\dagger/P_\Phi(F)$, which is a finite group. 
Hence, we can replace $q$ with $q^m$ for some positive integer $m$ and assume that $q^m$ is in $P_0(\Ap\otimes_\bb{Q} F)$. The quotient $P_0/W_0(\Ap\otimes F)$ acts \emph{faithfully} on $\gr_{-1} V_{\Ap}\otimes_\bb{Q} F$. So, by the diagram (\ref{diag-auto}) above, the image of $q^m$ in $P_0/W_0(\Ap\otimes_\bb{Q} F)$ satisfies the equation $X^M=1$ for some positive integer $M$, and therefore it is also true for the image of $q$ in $ZP_0/ZW_0(\Ap\otimes_\bb{Q} F)$. But it maps into the image in $ZP_\Phi/(Z_GW_0)(\Ap)=ZP_0/ZW_0(\Ap)$ of a \emph{neat} open compact subgroup $ZP(\Ap)\cap\lcj{K}{\gamma_\beta g_\beta g_\alpha}$. Consequently, we have the claim for $q$.\par
Since $\mrm{Aut}(\Q')\sbst \mrm{Aut}(Y_{\Q'})\times \mrm{Aut}(T_{\Q'})\times \mrm{Aut}(A_{\Q'})$, $q$ is in $Z_{G_0}(\Ap\otimes_\bb{Q} F)$. Now we can follow the argument in \cite[Prop. 3.6.4]{Kis10} and \cite[Cor. 4.6.15]{KP15}. Since $q\in \mrm{Aut}(\Q')\cap Z_0(\Ap\otimes_\bb{Q} F)$, we have $q\in Z_{0,\zbkp}(\ca{O}_{F,(p)})$. Replacing $\tg$ with another lifting, we can assume that the prime-to-$p$ factor of $q$ is trivial.
Then we know that $\tg\in ZP_0(\bb{Q})$ and $\tg$ maps into $\lcj{K_p}{\gamma_\beta g_\beta g_\alpha}$ in $G(\bb{Q}_p)$. This implies that $\tg=\gamma \in ZP_0(\bb{Q})$ maps to $ K^\beta_{p}$. Combining this fact with the paragraph above, we see that $q$ acts on every point of $\ca{S}_{K_{\Phi_{\beta'}}}$ trivially, since over the generic fiber, $q$ acts as the right action of some element in $K_{\Phi_{\beta'},p}$. Note that the last sentence is true because of the condition that $K_{0,p}^\alpha=\pi^{-1}(\lcj{K_p}{\alpha})$.\par 
The argument for $\overline{\ca{S}}_{K_{\Phi_{\beta'}}}$ is similar to the argument above. The differences are: for $\overline{\ca{S}}_{K_{\Phi_{\beta'}}}$, one replaces $V_?$ with $V_?/W_{-2}V_?$, and replaces objects related to $1$-motives $?\Q$ with $?\Q/W_{-2}?\Q$.\par
Let us show Part \ref{lem-twt-free-3}, and the statement in Part \ref{lem-twt-free-2} for $\ca{S}_{K_{\Phi_{\beta'}},h}$ will follow immediately. 
The idea is still similar to the paragraphs above with some changes. 
We can only assume that $\gamma$ lifts to $Q_\Phi(\bb{Q})_+\cap \lcj{K_p}{\gamma_\beta g_\beta g_\alpha}$ and to $\tg\in Q_0(F)\cap G_{0,\zbkp}(\ca{O}_{F,(p)})$. See the definition of $\wdtd{\Delta}^1_\beta$. \par
Let $G_l$ be the smallest connected normal subgroup in $Q_0/W_0$ such that there is an almost-direct product decomposition $Q_0/W_0=(P_0/W_0)\cdot G_l$.
Let $G_h':=(Q_0/W_0)/(G_l\cdot \bb{K})$ and $G_{0,h}:=P_0/W_0$. Then there is a natural central isogeny $G_{0,h}\to G_h'$, and this isogeny factors through $P_{\Phi,h}:=P_\Phi/W_\Phi$ since $G_{0,h}/\bb{K}=P_{\Phi,h}$. 
Recall that we have fixed an $i$ and $\lcj{K_p}{\gamma_\beta g_\beta g_\alpha}$ is contained in a Bruhat-Tits stabilizer subgroup $K_p^{i,\prime}:=\lcj{\mathscr{G}_{x_i}(\bb{Z}_p)}{\gamma_\beta g_\beta g_\alpha}$ of $G(\bb{Q}_p)$. The projection to $G_{h}'(\bb{Q}_p)$ of the intersection $K_p^{i,\prime}\cap Q_\Phi(\A)$ is compact and therefore contained in a Bruhat-Tits stabilizer subgroup $K_{\Phi,h,p}^{i,\prime}$. Denote the image of $\gamma$ in $P^\ad_{\Phi,h}(\bb{Q})^+=G_{0,h}^\ad(\bb{Q})^+$ by $\overline{\gamma}$. In particular, by Lemma \ref{lem-twist-admissible}, the action of $\overline{\gamma}$ on $\ca{S}_{K_{\Phi_{\beta'}},h}$ can be described by twisting abelian schemes. Moreover, we can check over generic fiber and by reducedness and separatedness that the action of $\gamma$ on $\ca{S}_{K_{\Phi_{\beta'}},h}$ is determined by the image $\overline{\gamma}$.
We now choose a finite Galois extension $F$ of $\bb{Q}$ and lift $\overline{\gamma}$ to $\tg_{11}\in G_{0,h,\zbkp}(\ca{O}_{F,(p)})$ via $G_{0,h}\to G'_h$ (in fact by lifting the projection to $G'_h(\bb{Q})$).\par 
If $(h,\gamma^{-1})$ acts on some geometric point of $\ca{S}_{K_{\Phi_{\beta'}},h}$ trivially, it follows that the $(-1)$-graded piece of the diagram (\ref{diag-auto}) is commutative. Denote by $\overline{h}$ the projection of $h$ to $G_{0,h}(\Ap)$. With similar conventions as before, we then see that the projection to $G_{0,h}'(\Ap)$ of $q:= k^{-1}\overline{h}\tg_{11}^{-1}\in G_{0,h}(\Ap\otimes_{\bb{Q}} F)$ is a neat element, where in this case $k\in K_{\Phi_{\beta'},h}$. 
On the other hand, $q$ satisfies the equation $X^m=1$ for some positive integer $m$, as it is a multiplication of a root of unity by \cite[Cor. 2.3.3.2]{Lan13} and the $(-1)$-graded piece of (\ref{diag-auto}).
This implies that $q\in Z_{G_{0,h}}(\Ap\otimes_\bb{Q} F)$ again. Finally, by the argument of \cite[Prop. 3.6.4]{Kis10} and \cite[Cor. 4.6.15]{KP15} as above, we see that $q\in Z_{G_{0,h}}(\Ap\otimes_\bb{Q}F)\cap\Aut(A_{\Q'})$, and this implies that $q\in Z_{G_{0,h},\zbkp}(\ca{O}_{F,(p)})$. Up to replacing $\tg_{11}$ with another lifting again, we can assume that the prime-to-$p$ factor of $q$ is trivial. Again, this implies that $q$ acts on every point of $\ca{S}_{K_{\Phi_{\beta'}},h}$ trivially.
\end{proof}
Note that the action of the whole group $G^\ad(\bb{Q})_1$ cannot be described by the twisting construction over integral models. But we have the following statement:
\begin{lem}\label{lem-twist-admissible}
With the conventions in the proof above, if $\gamma\in Q_\Phi(\bb{Q})^+\cap \lcj{K_p}{\gamma_\beta g_\beta g_\alpha}$, then we can describe the action of $\gamma$ as a (left) twisting action. 
\end{lem}
\begin{proof}
It suffices to deal with Case (STB). Without loss of generality, we assume $\gamma_\beta=1$, $g_\beta=1$ and $g_\alpha=1$. We show that the elements in $G(\bb{Q})^+\cap K_p$ can be described by twisting. 
Write $G\iso (G_0/\ker(\pi^b)\times S_G)/Z'_G$ as in the last proof, where $S_G\sbst Z_G$, and $Z'_G$ is a finite group.
By assumption, $\gamma$ can be projected to $\overline{\gamma}\in G':=G/S_G(\bb{Q})$, and $\overline{\gamma}$ can be lifted to $\tg\in G_0(F)$ for some finite Galois extension $F$ over $\bb{Q}$.\par
Suppose that $K_p=\mathscr{G}_x(\bb{Z}_p)$ is a Bruhat-Tits stabilizer subgroup corresponding to some point $x$ in the extended building $\ca{B}(G_{\bb{Q}_p},\bb{Q}_p)$. 
The point $x$ corresponds to a point $y$ in the Bruhat-Tits building of $G^\der_{\bb{Q}_p}$, and $y$ corresponds to a point $y'$ in the Bruhat-Tits building of $G'_{\bb{Q}_p}$ and a point $y_0$ in the Bruhat-Tits building of $G_{0,\bb{Q}_p}$.\par 
Then $\gamma$ lies in $\mathscr{G}_{y'}(\bb{Z}_p)$ under the projection, where $\mathscr{G}_{y'}$ is the Bruhat-Tits stabilizer group scheme corresponding to $y'$. Let $F\otimes_\bb{Q}\bb{Q}_p=\prod_i F_{v_i}$ for the places $v_i$ of $F$ over $p$. Then we can view $\gamma$ in $\prod_i R_{\ca{O}_{F_{v_i}}/\bb{Z}_p}\mathscr{G}_{y',\ca{O}_{F_{v_i}}}(\bb{Z}_p)$. Moreover, $R_{\ca{O}_{F_{v_i}}/\bb{Z}_p}\mathscr{G}_{y',\ca{O}_{F_{v_i}}}$ is the Bruhat-Tits stabilizer group scheme corresponding to the image of $y'$ in $\ca{B}(G_{F_{v_i}}',F_{v_i})\iso \ca{B}(R_{F_{v_i}/\bb{Q}_p}G_{F_{v_i}}',\bb{Q}_p)$. We can also view $\tg\in G_0(F)$ as an element in $R_{F\otimes_\bb{Q}\bb{Q}_p/\bb{Q}_p}G_0$. Since the map $R_{F\otimes_\bb{Q}\bb{Q}_p/\bb{Q}_p}G_0 \to R_{F\otimes_\bb{Q} \bb{Q}_p/\bb{Q}_p}G'$ is a central isogeny between reductive groups, $\tg\in \prod_i R_{\ca{O}_{F_{v_i}}/\bb{Z}_p}\mathscr{G}_{y_0,\ca{O}_{F_{v_i}}}(\bb{Z}_p)\cap G_0(F)\iso \mathscr{G}_{y_0}(\bb{Z}_p\otimes \ca{O}_F)\cap G_0(F)$.
Hence, $\tg$ lies in $G_{0,\zbkp}(\bb{Z}_p\otimes \ca{O}_F)\cap G_0(F)\sbst G_{0,\zbkp}(\zbkp\otimes \ca{O}_F)$, as desired.
\end{proof}

When $(G_0,X_0)=(G,X_b)$ and when $\gamma_\beta=g_\alpha=1$, we have reproved (and generalized) a result due to Madapusi.
\begin{cor}[{cf. the proof of \cite[Thm. 5.3.1]{Mad19}}]\label{cor-twt-free-action}
The quotient map $$\ca{S}_{K_{\Phi},h}\to\Delta_{\Phi,K}\bss\ca{S}_{K_{\Phi},h}$$ is finite {\'e}tale.
\end{cor}
\begin{rk} Note that we show this without assumptions on $K_p=K_{0,p}$. 
In fact, Madapusi has already remarked in the proof of \emph{loc. cit.} that one can use Kisin's theory of twisting abelian varieties to study the action of $\Delta_{\Phi,K}$, and has also remarked that one can drop the hyperspecial assumption added in \emph{loc. cit.} by doing this.
\end{rk}
The lemma above also generalizes the argument in the proof of \cite[Cor. 4.6.15]{KP15}; in particular, we have shown the following result:
\begin{prop}\label{prop-kisin-pappas-free-action}
When $P_0$ is a normal subgroup of $G_0$ and when $ZP_\Phi=G$, Part \ref{lem-twt-free-2} of Lemma \ref{lem-twt-free-action} implies that the morphism $\ca{S}_{K_0^\alpha}\to \ca{S}_K$ induced by taking the quotient of ${\Delta}_{\lcj{K}{\alpha}}(G_0,G)$ is finite {\'e}tale and that $\{\ca{S}_{K_pK^p}\}_{K^p}$ has finite {\'e}tale transition maps. Taking inverse limits with finite {\'e}tale transition maps, the quotient map $\ca{S}_{K_{0,p}^\alpha}\to \ca{S}_{K_p}$ induced by taking the quotient of $\Delta_{\lcj{K_p}{\alpha}}(G_0,G)$ is pro-finite {\'e}tale.
\end{prop}
Now, we see that the tower associated with $\ca{S}_{\K_\Phi}$ has the desired structure.
\begin{prop}\label{prop-mixsh-canonicity}Under the conventions and assumptions above, we have the following statements for different cases:
\begin{enumerate}
\item\label{prop-mixsh-canonicity-1}  In all cases, $\ca{S}_{\wdtd{K}_\Phi}$ is a torus torsor over $\overline{\ca{S}}_{\wdtd{K}_\Phi}$ under $\mbf{E}_{\K_\Phi}$ (see the paragraph below Definition \ref{def-zp-mix-shimura}).
Moreover, in Cases (HS) and ($\mrm{STB}_n$), 
$\overline{\ca{S}}_{\wdtd{K}_\Phi}$ is an abelian scheme torsor over $\ca{S}_{\wdtd{K}_{\Phi},h}$.
\item\label{prop-mixsh-canonicity-2} In all cases, the action of $\Delta^{ZP,\circ}_{\Phi,K}$ on $\ca{S}_{\K_\Phi}$ is trivial.
\item\label{prop-mixsh-canonicity-3} In Case (HS), the inverse limit formed by $\ca{S}_{\wdtd{K}_{\Phi,p}}:=\varprojlim_{\wdtd{K}^{p}_{\Phi}}\ca{S}_{\wdtd{K}_\Phi}$ satisfies the extension property.
\end{enumerate}
\end{prop}
\begin{proof}
By \cite[3.13]{Pin89} and Lemma \ref{lem-large-lattices}, we know that $\sh_{\K_\Phi}\to \overline{\sh}_{\K_\Phi}$ is a torsor under $\mbf{E}_{\K_\Phi}$. By \cite[Thm. 4.1.5 (4)]{Mad19}, $\ca{S}_{K_{\Phi_{\beta'}}}\to\overline{\ca{S}}_{K_{\Phi_{\beta'}}}$ is an $\mbf{E}_{K_{\Phi_{\beta'}}}$-torsor. From the construction of $\ca{S}_{\K_\Phi}\to \overline{\ca{S}}_{\K_\Phi}$ as a quotient, we have a finite morphism $\disju_{\beta\in I_\Phi}\mbf{E}_{K_{\Phi_{\beta'}}}\times \ca{S}_{K_{\Phi_{\beta'}}}\to \ca{S}_{\K_\Phi}\times_{\overline{\ca{S}}_{\K_\Phi}}\ca{S}_{\K_\Phi}$. From Lemma \ref{lem-twt-free-action}(2), we see that the target of this morphism is normal. This morphism factors as $\disju_{\beta\in I_\Phi}\mbf{E}_{K_{\Phi_{\beta'}}}\times \ca{S}_{K_{\Phi_{\beta'}}}\to \mbf{E}_{\K_\Phi}\times \ca{S}_{\K_\Phi}\to \ca{S}_{\K_\Phi}\times_{\overline{\ca{S}}_{\K_\Phi}}\ca{S}_{\K_\Phi}$, and the last morphism $\mbf{E}_{\K_\Phi}\times \ca{S}_{\K_\Phi}\to \ca{S}_{\K_\Phi}\times_{\overline{\ca{S}}_{\K_\Phi}}\ca{S}_{\K_\Phi}$ is an isomorphism by Zariski's main theorem. Finally, since, over the fppf cover $\disju_{\beta\in I_\Phi}\ca{S}_{K_{\Phi_{\beta'}}}\to \overline{\ca{S}}_{\K_\Phi}$, $\ca{S}_{\K_\Phi}\to \overline{\ca{S}}_{\K_\Phi}$ admits a section, we conclude that it is an $\mbf{E}_{\K_\Phi}$-torsor. 
The claim for abelian scheme torsors is proved in a similar manner. By \cite[Thm. 4.1.5 (1)]{Mad19}, $\overline{\ca{S}}_{K_{\Phi_{\beta'}}}\to \ca{S}_{K_{\Phi_{\beta'}},h}$ is a torsor under some abelian scheme $\ca{A}_{K_{\Phi_{\beta'}}}$ over $\ca{S}_{K_{\Phi_{\beta'}},h}$. By \cite[3.13]{Pin89} again, $\overline{\sh}_{\K_\Phi}\to\sh_{\K_\Phi,h}$ is a torsor under some abelian scheme $A_{\K_{\Phi}}$ over $\sh_{\K_\Phi,h}$, and $\overline{\sh}_{K_{\Phi_{\beta'}}}\to \sh_{K_{\Phi_{\beta'}},h}$ is a torsor under $A_{K_{\Phi_{\beta'}}}$, where $A_{K_{\Phi_{\beta'}}}:= \ca{A}_{K_{\Phi_{\beta'}}}\times_{\ca{S}_{K_{\Phi_{\beta'}},h}}\sh_{K_{\Phi_{\beta'}},h}$ is an abelian scheme over $\sh_{K_{\Phi_{\beta'}},h}$. 
By Lemma \ref{lem-twt-free-action}(2), $\disju_{\beta\in I_\Phi}\overline{\ca{S}}_{K_{\Phi_{\beta'}}}\to\overline{\ca{S}}_{\K_\Phi}$ and $\disju_{\beta\in I_\Phi}\ca{S}_{K_{\Phi_{\beta'}},h}\to\ca{S}_{\K_\Phi,h}$ are finite {\'e}tale surjective morphisms between normal schemes. Then the induced morphism $\disju_{\beta\in I_\Phi}\overline{\ca{S}}_{K_{\Phi_{\beta'}}}\to\disju_{\beta\in I_\Phi}\ca{S}_{K_{\Phi_{\beta'}},h}\times_{\ca{S}_{\K_\Phi,h}} \overline{\ca{S}}_{\K_\Phi}$ is finite {\'e}tale and surjective. Let $\ca{K}$ be the schematic closure in $\ca{A}:=\disju_{\beta\in I_\Phi} \ca{A}_{K_{\Phi_{\beta'}}}$ of the kernel of $\disju_{\beta\in I_\Phi}A_{K_{\Phi_{\beta'}}}\to A_{\K_{\Phi}}$. Then $\ca{K}$ is finite over $\disju_{\beta\in I_\Phi}\ca{S}_{K_{\Phi_{\beta'}},h}$ as it is a closed subscheme in $\ca{A}[n]$ for some positive integer $n$. 
By Zariski's main theorem, the quotient of $\disju_{\beta\in I_\Phi}\overline{S}_{K_{\Phi_{\beta'}}}$ by $\ca{K}$ is isomorphic to $\disju_{\beta\in I_\Phi}\ca{S}_{K_{\Phi_{\beta'}},h}\times_{\ca{S}_{\K_\Phi,h}} \overline{\ca{S}}_{\K_\Phi}$. Hence, $\ca{K}$ is a finite {\'e}tale group scheme over $\disju_{\beta\in I_\Phi}\ca{S}_{K_{\Phi_{\beta'}},h}$. We then have that $\ca{A}/\ca{K}$ is an abelian scheme over $\disju_{\beta\in I_\Phi}\ca{S}_{K_{\Phi_{\beta'}},h}$ extending $A_{\K_\Phi}$. By Lemma \ref{lem-twt-free-action}(2) again, $\overline{\ca{S}}_{\K_\Phi}$ is an abelian scheme torsor over $\ca{S}_{\K_\Phi,h}$ since we have checked it over an {\'e}tale cover $\disju_{\beta\in I_\Phi}\ca{S}_{K_{\Phi_{\beta'}},h}$ of $\ca{S}_{\K_{\Phi},h}$.\par
Now let us show Part \ref{prop-mixsh-canonicity-2}. We argue as in \cite[7.3]{Pin89} and the proof of Lemma \ref{lem-twt-free-action}(1). 
The projection of $\Delta_{\Phi,K}^{ZP,\circ}$ to $Q_\Phi/ZP_\Phi(\bb{Q})$ is a neat arithmetic subgroup. By \cite[Thm. 6.19(a)]{Pin89} and Lemma \ref{lem-Harris-group-compare}, this projection is in a neat arithmetic subgroup in $P''_\Phi(\bb{Q})$, so it has to be in $ZP_\Phi(\bb{Q})$, as desired.\par
Now assume that we are in Case (HS). By \cite[Lem. 2.1.3]{Lov17} and Proposition \ref{prop-par-geo-conn}, we can and we will check the extension property over $\ca{O}^p:=\ca{O}_{E^p}\otimes_{\ca{O}_{E_0}}\ca{O}_{E_0,(v)}$, where $E^p$ is the maximal field extension of $E'$ that is unramified at all primes dividing $p$.
It follows from essentially the same argument as in \cite[E. 6]{Kis17}. 
By \cite[A.3.5 and A.3.6]{Mad19}, we have that, 
for any rational prime $l\neq p$ and any neat open compact subgroup $K_{\Phi_0}^{p,l}\sbst P_{\Phi_0}(\A^{p,l})$, 
$\ca{S}_{K_{\Phi_0,p}K^{p,l}_{\Phi_0}}$ has the extension property.
Part \ref{prop-mixsh-canonicity-3} will follow from \cite[Lem 2.1.5 and Lem. 2.1.4]{Lov17} and Lemma \ref{lem-twt-free-action}, if we show that the quotient of $\ca{S}_{K_{\Phi_0,p}K_{\Phi_0}^{p,l},\ca{O}^p}$ by $\Delta_{K^{\beta,l}}$ factors through a finite group, where
$$\Delta_{K^{\beta,l}}:=\ker(\ag(P_0)\to \ag(ZP_\Phi)/K^{\beta,l}),$$
$K^{\beta,l}:=ZP_\Phi(\A)\cap \lcj{K^l}{\gamma_\beta g_\beta g_\alpha}$ and $K^l=K_pK^{l,p}$ for some neat compact open $K^{l,p}\sbst G(\A^{p,l})$.\par
Write $\ag(P_0)$ as $P_0(\A)\cdot \agsb(P_0)$ and $\ag(ZP_\Phi)$ as $ZP_\Phi(\A)\cdot \agsb(ZP_\Phi)$. Since 
$$\ker(P_0(\bb{Q})_+^-\bss P_0(\A)/K_{\Phi_0,p}K_{\Phi_0}^{p,l}\to ZP_\Phi(\bb{Q})_+^-\bss ZP_\Phi(\A)/K^{\beta,l})$$ 
is finite (since the kernel of $P_0\to ZP_\Phi$ is finite by construction, one can see this by decomposing to $l$- and away-from-$l$ places), and by Lemma \ref{lem-ker-cpt},  
the quotient of $\ca{S}_{K_{\Phi_0,p}K_{\Phi_0}^{p,l},\ca{O}^p}$ by $\Delta_{K^{\beta,l}}$ indeed factors through a finite group.
\end{proof}
\begin{construction}\label{const-cone-mix-int-model}\upshape
By Proposition \ref{prop-mixsh-canonicity}, $\ca{S}_{\wdtd{K}_\Phi}\to \overline{\ca{S}}_{\wdtd{K}_\Phi}$ is an $\mbf{E}_{\K_\Phi}$-torsor. Let $\Sigma$ be an admissible cone decomposition and $\sigma\in \Sigma(\Phi)$ such that $\sigma^\circ\sbst \mbf{P}^+_\Phi$. It makes sense to define
$$\ca{S}_{\wdtd{K}_\Phi}(\Sigma):=\ca{S}_{\wdtd{K}_\Phi}\times^{\mbf{E}_{\K_\Phi}}\mbf{E}_{\K_\Phi}(\Sigma),$$
where $\mbf{E}_{\K_\Phi}(\Sigma)$ is the affine torus embedding of $\mbf{E}_{\K_\Phi}$ defined by $\Sigma^+(\Phi)$,
$$\ca{S}_{\wdtd{K}_\Phi}(\sigma):=\ca{S}_{\wdtd{K}_\Phi}\times^{\mbf{E}_{\K_\Phi}}\mbf{E}_{\K_\Phi}(\sigma),$$
and the $\sigma$-stratum of $\ca{S}_{\wdtd{K}_\Phi}(\sigma)$:
$$\ca{S}_{\wdtd{K}_\Phi,\sigma}:=\ca{S}_{\wdtd{K}_\Phi}\times^{\mbf{E}_{\K_\Phi}}\mbf{E}_{\K_\Phi,\sigma}.$$
Note that there is a canonical projection from $\mbf{E}_{\K_\Phi}$ to $\mbf{E}_{\K_\Phi,\sigma}$.\hfill$\square$
\end{construction}
Finally, let us explicitly check the following statement:
\begin{prop}\label{prop-zpb-completion-iso}
With the assumptions and constructions above, the isomorphism $\mrm{Z}_{[ZP^b({\Phi},\sigma)],K}\iso \Delta^{ZP,\circ}_{\Phi,K}\bss \sh_{\wdtd{K}_\Phi,\sigma}$ extends to 
\begin{equation}\label{eq-zp-strata-compatible-int}
    \ca{Z}_{[ZP^b({\Phi},\sigma)],K}\iso \Delta^{ZP,\circ}_{\Phi,K}\bss \ca{S}_{\wdtd{K}_\Phi,\sigma}.
\end{equation}
Moreover, this isomorphism extends to an isomorphism
\begin{equation}\label{eq-zp-strata-completion-int}
   \cpl{\ca{S}^{\Sigma}_K}{ \ca{Z}_{[ZP^b({\Phi},\sigma)],K}} \iso \Delta^{ZP,\circ}_{\Phi,K}\bss \cpl{\ca{S}_{\K_\Phi}(\sigma)}{\ca{S}_{\K_\Phi,\sigma}}.
\end{equation}
Since we assume that $K$ is neat, the quotients of $\Delta^{ZP,\circ}_{\Phi,K}$ in the right-hand sides of (\ref{eq-zp-strata-compatible-int}) and (\ref{eq-zp-strata-completion-int}) are isomorphic to $\ca{S}_{\K_\Phi,\sigma}$ and $\cpl{\ca{S}_{\K_\Phi}(\sigma)}{\ca{S}_{\K_\Phi,\sigma}}$, respectively.
\end{prop}
\begin{proof}The last sentence follows from Proposition \ref{prop-mixsh-canonicity}(2). 
Write the RHS of (\ref{eq-zp-strata-compatible-int}) as $\Delta_{\Phi,K}^{ZP}\bss \Delta_{\Phi,K}^{ZP}\ca{S}_{\K_\Phi,\sigma}\iso G(\bb{Q})_+\bss G(\bb{Q})_+\times^{\wdtd{\Delta}_{\Phi,K}^{ZP}\cap G(\bb{Q})_+}\ca{S}_{\K_\Phi,\sigma}$. We write the last term as $G(\bb{Q})_+\bss G(\bb{Q})_+ \ca{S}_{\K_\Phi,\sigma}.$ 
From the construction of $\ca{S}_{\K_\Phi}$, we have
\begin{equation*}
\begin{split}
&G(\bb{Q})_+\bss G(\bb{Q})_+ \ca{S}_{\K_\Phi,\sigma}\\ 
&\iso G(\bb{Q})_+\bss G(\bb{Q})_+\disju_{\beta \in I_\Phi} \gamma_\beta\cdot \ca{S}_{K_{\Phi_\beta},\sigma_\beta}/\Delta_{K^\beta}(P_0,ZP_\Phi)\\
&\iso \disju_{\overline{\beta} \in J_\Phi} G(\bb{Q})_+\bss G(\bb{Q})_+ \ca{S}_{K_{\Phi_{\beta'}},\sigma_{\beta'}}\\
&\iso \disju_{\alpha\in I_{G/G_0}}\disju_{\pi(g_0^\alpha)\alpha\sim g^b} G(\bb{Q})_+\bss G(\bb{Q})_+ \ca{S}_{K_{\Phi_0^\alpha},\sigma_0^\alpha}\\
&\iso \disju_{\alpha\in I_{G/G_0}}\disju_{\pi(g_0^\alpha)\alpha\sim g^b} \Delta_{\lcj{K}{g_0^\alpha g_\alpha}}(P_{0},G)\ca{S}_{K_{\Phi_0^\alpha},\sigma_0^\alpha}/\Delta_{\lcj{K}{g_0^\alpha g_\alpha}}(P_{0},G)\\
&\iso \disju_{\alpha\in I_{G/G_0}}\disju_{\pi(g_0^\alpha)\alpha\sim g^b} (\Delta_{\lcj{K}{\alpha}}(G_0,G)\ca{Z}_{[(\Phi_0^\alpha,\sigma_0^\alpha)],K_0^\alpha})/\Delta_{\lcj{K}{\alpha}}(G_0,G)\\
& \iso \mrm{LHS}.
\end{split}
\end{equation*}
In the equations above, we write $\Phi_{\beta'}=\gamma_\beta\cdot \Phi_\beta$ and $\sigma_{\beta'}:=\gamma_\beta\cdot\sigma_\beta$. Write $G(\bb{Q})_+\ca{S}_{K_{\Phi_{\beta'}},\sigma_{\beta'}}:=G(\bb{Q})_+\times^{\Delta_{K^\beta}(P_0,G)}\Delta_{K^\beta}(P_0,G)\ca{S}_{K_{\Phi_{\beta'}},\sigma_{\beta'}}$, where $\Delta_{K^\beta}(P_0,G):=\ker (\ag(P_0)\times^{P_0^\ad(\bb{Q})^+} Q_0^\ad(\bb{Q})^+\to \ag(G)/\lcj{K}{\gamma_\beta g_\beta g_\alpha})$ (see the beginning of \S\ref{subsubsec-intmod-study} for the definition of $Q_0^\ad(\bb{Q})^+$). Similarly, we write $G(\bb{Q})_+\ca{S}_{K_{\Phi_0^\alpha},\sigma_0^\alpha}:=G(\bb{Q})_+\times^{\Delta_{\lcj{K}{g_0^\alpha g_\alpha}}(P_0,G)}\Delta_{\lcj{K}{g_0^\alpha g_\alpha}}(P_0,G)\ca{S}_{K_{\Phi_0^\alpha},\sigma_0^\alpha}$, where $\Delta_{\lcj{K}{g_0^\alpha g_\alpha}}(P_0,G):=\ker (\ag(P_0)\times^{P_0^\ad(\bb{Q})^+} Q_0^\ad(\bb{Q})^+\to \ag(G)/\lcj{K}{g_0^\alpha g_\alpha})$.\par
Let us explain the equations above: The second line comes from the definition. Any $(p,\gamma^{-1})\in \Delta_{K^\beta}(P_0,ZP_\Phi)$ is determined by the element $\gamma$ that is liftable to $ZP_\Phi(\bb{Q})_+/Z(\bb{Q})$ and to $\gamma^*\in ZP_\Phi(\bb{Q})_+$. Then $\Delta_{K^\beta}(P_0,ZP_\Phi)\sbst \Delta_{K^\beta}(P_0,G)$. 
We see the third line is isomorphic to $$ \disju_{\overline{\beta} \in J_\Phi} G(\bb{Q})_+\bss G(\bb{Q})_+ (\ca{S}_{K_{\Phi_{\beta'}},\sigma_{\beta'}}/\Delta_{K^\beta}(P_0,ZP_\Phi)).$$ Moreover, the stabilizer in $\wdtd{\Delta}^{ZP}_{\Phi,K}\cap G(\bb{Q})_+$ of an element in $J_\Phi$ is isomorphic to $\stb_{G(\bb{Q})_+}(D_{\Phi_\beta})\cap \pi(P_0(\A))\lcj{K}{\gamma_\beta g_\beta g_\alpha}$; the image of $\Delta_{K^\beta}(P_0,G)$, viewed as elements in $(p,\gamma)\in P_0(\A)\times G(\bb{Q})_+$, is surjective to this stabilizer. Hence, we have the isomorphism from the second line to the third line. 
The fourth line comes from the discussion in Lemma \ref{lem-db-cst-g}, and the fifth line is the definition. The isomorphism from the fifth line to the sixth line follows from the fact that $\Delta_{\Phi_0^\alpha,K_0^\alpha}$ maps to $\Delta_{\lcj{K}{g_0^\alpha g_\alpha}}(P_{0},G)$ and the fact that $\Delta_{\lcj{K}{g_0^\alpha g_\alpha}}(P_{0},G)$ is contained in $\Delta_{\lcj{K}{\alpha}}(G_0,G)$ via the map sending $(h,\gamma)$ to $((g_0^\alpha)^{-1}hg_0^\alpha,\gamma)$; in addition, an element $((g_0^\alpha)^{-1}hg_0^\alpha,\gamma)\in \Delta_{\lcj{K}{\alpha}}(G_0,G)$ stabilizes $\ca{Z}_{[(\Phi_0^\alpha,\sigma_0^\alpha)]}$ if and only if it is in the image of $\Delta_{\lcj{K}{g_0^\alpha g_\alpha}}(P_0,G)$. \par
Let us show the second isomorphism. Note that we have similar isomorphisms as above if we replace closed strata with the corresponding twisted affine torus embeddings. Hence, we have 
\begin{equation*}
\begin{split}
&\mrm{RHS}\\ 
&\iso\disju_{\alpha\in I_{G/G_0}}\disju_{\pi(g_0^\alpha)\alpha\sim g^b} (\Delta_{\lcj{K}{\alpha}}(G_0,G)\cpl{\ca{S}_{K_0^\alpha}(\sigma_0^\alpha)}{\ca{S}_{K_{\Phi_0^\alpha},\sigma_0^\alpha}})/\Delta_{\lcj{K}{\alpha}}(G_0,G)\\
&\iso \disju_{\alpha\in I_{G/G_0}}\disju_{\pi(g_0^\alpha)\alpha\sim g^b} (\Delta_{\lcj{K}{\alpha}}(G_0,G) \cpl{\ca{S}_{K_0^\alpha}^{\Sigma_0^\alpha}}{\ca{Z}_{[(\Phi_0^\alpha,\sigma_0^\alpha)],K_0^\alpha}})/\Delta_{\lcj{K}{\alpha}}(G_0,G)\\
& \iso \mrm{LHS}.
\end{split}
\end{equation*}
The third line to the last line requires that the actions of $\Delta_{\lcj{K}{\alpha}}(G_0,G)$ on boundaries and on the compactification exist and are compatible, which follows from Lemma \ref{lem-gadq1-ext} and \cite[Lem. A.3.4]{Mad19}.
\end{proof}
\subsubsection{}\label{subsubsec-const-GXa}
Assume that we are in Case ($\text{STB}_n$) or (HS). From $\ca{S}_{\K_\Phi^b}(G,X_b)$, we can switch it to $\ca{S}_{\K_\Phi^a}:=\ca{S}_{\K_\Phi^a}(G,X_a)$ by multiplying the Galois action induced by $c$. 

\begin{construction}\label{const-mixed-a}\upshape
Let $\ca{S}_{\K_{\Phi,p},\ca{O}^{\prime,ur}_{(v)}}:=\varprojlim_{\K_\Phi^p\text{\ neat open compact}}\ca{S}_{\K_\Phi,\ca{O}^{\prime,ur}_{(v)}}$. Similarly, define $\overline{\ca{S}}_{\K_{\Phi,p},\ca{O}^{\prime,ur}_{(v)}}$ and $\ca{S}_{\K_{\Phi,h,p},\ca{O}^{\prime,ur}_{(v)}}$.\par
Recall that $Z:=Z_G$. Let $c$ be the homomorphism $c:\bb{S}\to Z_\bb{R}$ defined as in \S\ref{a,b}. Then $(Z,\{c\})$ is a $0$-dimensional Shimura datum. 
Let $r_{E',{K_{Z,p}K_Z^p}}(Z,\{c\})$ be the homomorphism 
$\gal(\overline{E}'/E')\to \gal(E^{\prime,ab}/E')\to Z(\bb{Q})\bss Z(\A)/K_{Z,p}K_Z^p$ defined as in \S\ref{subsubsec-can-model} describing the action of $\gal(\overline{E}'/E')$ on $\sh_{K_{Z,p}K_Z^p,\overline{\bb{Q}}}\iso \pi_0(\sh_{K_{Z,p}K_Z^p,\overline{\bb{Q}}})\iso Z(\bb{Q})\bss Z(\A)/K_{Z,p}K_Z^p$. Here we can take $K_{Z,p}$ as $K_p\cap Z(\bb{Q}_p)$ and take $K_Z^p\sbst K^p\cap Z(\Ap)$. 
By Proposition \ref{prop-mixsh-independent}, the left action of $Z(\A)$ on $\sh_{\K_\Phi}$ extends to an action on $\ca{S}_{\K_\Phi}$, and $Z(\bb{Q})K_{Z,p}K_Z^p$ is in the kernel of this action since it is so over the generic fiber.\par
Then we can define a new Galois action as follows:
Let $r_b: \gal(E^{\prime,ab}/E')\times \ca{S}_{\K_{\Phi,p},\ca{O}^{\prime,ur}_{(v)}}\to \ca{S}_{\K_{\Phi,p},\ca{O}^{\prime,ur}_{(v)}}$ be the action of $\gal(E^{\prime,ab}/E')$ on $\ca{S}_{\K_{\Phi,p},\ca{O}^{\prime,ur}_{(v)}}$ which determines the descent datum of $\ca{S}_{\K_{\Phi,p},\ca{O}^{\prime,ur}_{(v)}}$ from $\ca{O}^{\prime,ur}_{(v)}$ to $\ca{O}'_{(v)}$. Define $r_a:\gal(E^{\prime,ab}/E')\times \ca{S}_{\K_{\Phi,p},\ca{O}^{\prime,ur}_{(v)}}\to \ca{S}_{\K_{\Phi,p},\ca{O}^{\prime,ur}_{(v)}}$ by sending each $\sigma\in \gal(E^{\prime,ab}/E')$ to $r_a(\sigma):=r_{E',{K_{Z,p}K_Z^p}}(Z,\{c\})(\sigma)\circ r_b(\sigma)=r_b(\sigma)\circ r_{E',{K_{Z,p}K_Z^p}}(Z,\{c\})(\sigma)$. We can exchange the order since the left action of $Z(\A)$ is defined over $\ca{O}_{(v)}'$. Note that the descent datum factors through $\gal(E^{p}/E')$ since the intersection of $n$ quasi-parahoric subgroups in $G(\bb{Q}_p)$ still contains the parahoric subgroup of $Z(\bb{Q}_p)$.\par
The tower induced by the new descent datum determined by $r_a$ is a model of $\sh_{\K_{\Phi,p}^a}:=\sh_{\K_{\Phi,p}^a}(G,X_a)$ over $\ca{O}_{(v)}'$, we denote it by $\ca{S}_{\K^a_{\Phi,p}}$ or $\ca{S}_{\K^a_{\Phi,p}}(G,X_a)$. Similarly, we define the towers of $\ca{O}_{(v)}'$-schemes $\overline{\ca{S}}_{\K^a_{\Phi,p}}$ and $\ca{S}_{\K^a_{\Phi,p},h}$ which are models of $\overline{\sh}_{\K^a_{\Phi,p}}$ and $\sh_{\K^a_{\Phi,p},h}$, respectively.\par
By {\'e}tale descent of morphisms and since the schemes involved are normal and separated, we see that for any neat open compact $K^p\sbst G(\Ap)$, there is also a tower of $\ca{O}'_{(v)}$-schemes $\ca{S}_{\K_{\Phi}^a}\to \overline{\ca{S}}_{\K_\Phi^a}\to \ca{S}_{\K_{\Phi}^a,h}$. Also from {\'e}tale descent, we know that this tower also satisfies Proposition \ref{prop-mixsh-canonicity} and Proposition \ref{prop-mixsh-independent}.
\hfill$\square$
\end{construction}
\begin{rk}
It is possible to write the tower $\ca{S}_{\K^i_{\Phi,p},\ca{O}^{\prime,ur}_{(v)}}$ in a similar form as \cite[Cor. 4.6.15]{KP15}, for $i=a,b$ and for every $\Phi$. However, the complexity will increase when one wants to study the relation between the compactifications and the boundaries with this method. This partially explains why we chose to use a different strategy here to study compactifications.
\end{rk}
\subsection{Main theorem on toroidal compactifications}\label{subsec-main-thm}
\subsubsection{}\label{sss-main-setup-again} For the rest of this {section}, let us complete the main constructions.\par
Let us fix an open compact subgroup $K_{2,p}\sbst G_2(\bb{Q}_p)$ and a neat open compact subgroup $K_2^p$ of $G_2(\Ap)$. Let $\Sigma_2$ be an admissible rational polyhedral cone decomposition of $(G_2,X_2)$. By Proposition \ref{ext-imp}, Corollary \ref{ext-cpt-imp} and Proposition \ref{zp-cones}, there is a neat open compact subgroup $K=K_pK^p$ of $G(\A)$ such that, $K_p\cap G_2(\bb{Q}_p)=K_{2,p}$, and such that $\sh_{K_2}\to \sh_{K}$ is an open and closed embedding; moreover, there is a refinement $\Sigma_2'$ of $\Sigma_2$ and a $ZP$-invariant admissible cone decomposition $\Sigma$ associated with $(G,X_a,K)$ and $(G,X_b,K)$, such that $\Sigma_2'$ is induced by $\Sigma$, and such that the morphisms between toroidal compactifications, boundary mixed Shimura varieties and boundary strata from the ones associated with $(G_2,X_2,K_2)$ to those associated with $(G,X_a,K)$ are all open and closed embeddings. For different cases (HS), ($\text{STB}_n$) and (DL), we choose $(G_0,X_0)$, $(G,X)$, $g_\alpha$, $K=K_pK^p$, $(G^{\ddag},X^{\ddag})$, $K_0^\alpha=K_{0,p}^\alpha K_0^{\alpha,p}$, $K^{\ddag,\alpha}=K_p^{\ddag,\alpha}K^{\ddag,\alpha,p}$, etc., according to the discussion at the beginning of \S\ref{subsec-con-stra-comp-quo}. By Proposition \ref{zp-cones}, $\Sigma$ can be and will be chosen such that $\Sigma$ also induces $\Sigma_0^\alpha$ by pulling back under $\pi^b(g_\alpha)$, and such that good toroidal compactifications $\ca{S}_{K_0^\alpha}^{\Sigma_0^\alpha}$ associated with $(G_0,X_0,K_0^\alpha,\Sigma_0^\alpha)$ satisfying \cite[Thm. 4.1.5]{Mad19} can be constructed.\par
Let us construct integral models of toroidal compactifications of $(G,X_a)$ and $(G_2,X_2)$.\par
\begin{construction}\label{const-strata-compactification}\upshape
Let $\Sigma$ and $\Sigma_0^\alpha$ be cone decompositions chosen as mentioned in Construction \ref{const-tor-Gb} and at the beginning of \S\ref{sss-main-setup-again}. Let $\Phi_2$ be a cusp label representative of $(G_2,X_2)$ that maps to $\Phi\in \ca{CLR}(G,X_a)$ under $\pi^a$. By \S\ref{subsec-zp-cusp}, we can view $\Phi$ as a cusp label representative of $(G,X_b)$ in $\ca{CLR}(G,X_b)$.\par 
Let us first treat cases (HS) and ($\text{STB}_n$).\par
We construct $\ca{S}_{\K_{\Phi}^a}\to \overline{\ca{S}}_{\K_\Phi^a}\to \ca{S}_{\K_{\Phi}^a,h}$ as in Construction \ref{const-mixed-a}. 
Since $\ca{S}_{\K^a_\Phi}\to \overline{\ca{S}}_{\K^a_\Phi}$ is an $\mbf{E}_{K_\Phi}$-torsor, we define $\ca{S}_{\K_\Phi^a}(\Sigma)$, $\ca{S}_{\K_\Phi^a}(\sigma)$ and $\ca{S}_{\K_\Phi^a,\sigma}$ as before. More precisely, 
let $\Sigma$ be an admissible cone decomposition and $\sigma\in \Sigma(\Phi)$ such that $\sigma^\circ\sbst \mbf{P}^+_\Phi$. Define
$$\ca{S}_{\wdtd{K}^a_\Phi}(\Sigma):=\ca{S}_{\wdtd{K}^a_\Phi}\times^{\mbf{E}_{\K_\Phi}}\mbf{E}_{\K_\Phi}(\Sigma),$$
where $\mbf{E}_{\K_\Phi}(\Sigma)$ is the affine torus embedding of $\mbf{E}_{\K_\Phi}$ defined by $\Sigma^+(\Phi)$,
$$\ca{S}_{\wdtd{K}^a_\Phi}(\sigma):=\ca{S}_{\wdtd{K}^a_\Phi}\times^{\mbf{E}_{\K_\Phi}}\mbf{E}_{\K_\Phi}(\sigma),$$
and the $\sigma$-stratum of $\ca{S}_{\wdtd{K}^a_\Phi}(\sigma)$:
$$\ca{S}_{\wdtd{K}^a_\Phi,\sigma}:=\ca{S}_{\wdtd{K}^a_\Phi}\times^{\mbf{E}_{\K_\Phi}}\mbf{E}_{\K_\Phi,\sigma}.$$
Let $\sigma_2\in \Sigma_2'(\Phi_2)$ be the cone such that $\sigma_2^\circ$ maps to $\sigma^\circ$. Let $\ca{S}_{K_{\Phi_2}}(\sigma_2)$ be the normalization of $\ca{S}_{\K_{\Phi}^a}(\sigma)$ in $\sh_{K_{\Phi_2}}(\sigma_2)$ via the composition of morphisms $\sh_{K_{\Phi_2}}(\sigma_2)\to \sh_{\K_{\Phi}^a}(\sigma)\hookrightarrow \ca{S}_{\K_{\Phi}^a}(\sigma)$. We define $\ca{S}_{K_{\Phi_2},\sigma_2}$, $\overline{\ca{S}}_{K_{\Phi_2}}$ and $\ca{S}_{K_{\Phi_2},h}$ similarly. Hence, this construction induces towers of schemes $\ca{S}_{K_{\Phi_2}}(\sigma_2)\to \overline{\ca{S}}_{K_{\Phi_2}}\to \ca{S}_{K_{\Phi_2},h}$, and $\ca{S}_{K_{\Phi_2},\sigma_2}\to \overline{\ca{S}}_{K_{\Phi_2}}\to \ca{S}_{K_{\Phi_2},h}$.\hfill$\square$
\end{construction}
For all the cases, we can continue with the following construction:
\begin{construction}\label{const-strata-compactification-2}\upshape 
Let $K_{2,p}$ be an open compact subgroup contained in $K_{2,p}'\sbst G_2(\bb{Q}_p)$, where $K_{2,p}'$ is an intersection of $n$ Bruhat-Tits stabilizer subgroups. Then $K_{2,p}'$ is contained in $K'_p\sbst G(\bb{Q}_p)$, where $K'_p$ is an intersection of $n$ Bruhat-Tits stabilizer subgroups of $ G(\bb{Q}_p)$. Set $K_2:=K_{2,p}K^p_2$ and $K':=K_p'K^{\prime,p}$ where $K^p_2\sbst K^{\prime,p}$ are neat open compact.
\begin{enumerate}
\item Let $\ca{S}_{K_{\Phi_2}}$ (resp. $\overline{\ca{S}}_{K_{\Phi_2}}$, and resp. $\ca{S}_{K_{\Phi_2,h}}$) be the normalization of  $\ca{S}_{\K^{\prime,a}_\Phi}$ (resp. $\overline{\ca{S}}_{\K^{\prime,a}_\Phi}$ and resp. $\ca{S}_{\K^{\prime,a}_{\Phi},h}$) in $\sh_{K_{\Phi_2}}$ (resp. $\overline{\sh}_{K_{\Phi_2}}$, and resp. $\sh_{K_{\Phi_2,h}}$).
\item Let $\ca{S}_{\K^{a}_\Phi}$ (resp. $\overline{\ca{S}}_{\K^{a}_\Phi}$ and resp. $\ca{S}_{\K^{a}_{\Phi},h}$) for (DL) be the normalization of $\ca{S}_{\K^{\prime,a}_\Phi}$ (resp. $\overline{\ca{S}}_{\K^{\prime,a}_\Phi}$ and resp. $\ca{S}_{\K^{\prime,a}_{\Phi},h}$) in (HS) or ($\text{STB}_n$) case in $\sh_{\K^{a}_\Phi}$ (resp. $\overline{\sh}_{\K^{a}_\Phi}$ and resp. $\sh_{\K^{a}_{\Phi},h}$).
\item Recall that $F_{K_Z}$ is a number field such that every geometrically connected component of $\sh_{K_Z}$ is defined over $F_{K_Z}$, where $K_Z$ is an open compact subgroup of $Z(\A)$ contained in $K$ and $\sh_{K_Z}:=\sh_{K_Z}(Z,\{c\})$ (see \S\ref{a,b}). We denote $E_{K_Z}':=E'\cdot F_{K_Z}$ and choose $F_{K_Z}$ such that $E_{K_Z}'/E'$ is finite and Galois. Note that $E'_{K_Z}$ can not be chosen to be unramified over $E'$ at $p$ in general. We can still find a $F_{K_Z}$ so that $E'_{K_Z}$ is unramified over $E'$ at $p$ if $K_p$ contains the parahoric subgroup of $Z(\bb{Q}_p)$. 
Let $\ca{O}_{K_Z}:=\ca{O}_{E'_{K_Z}}\otimes_{\ca{O}_{E'}}\ca{O}'_{(v)}$. 
Let $\ca{S}^{\Sigma_2'}_{K_2,\ca{O}_{K_Z}}$ be the normalization of $\ca{S}_{K,\ca{O}_{K_Z}}^\Sigma$ in $\sh_{K_2,E'_{K_Z}}^{\Sigma_2'}$ via the composition of morphisms $\sh_{K_2,E'_{K_Z}}^{\Sigma_2'}\to \sh_{K}^\Sigma(G,X_a)_{E'_{K_Z}}\xrightarrow{\sim}\sh^\Sigma_K(G,X_b)_{E'_{K_Z}}\hookrightarrow\ca{S}_{K,\ca{O}_{K_Z}}^{\Sigma}$, where the last subscript $\ca{O}_{K_Z}$ denotes the normalized base change (i.e., the normalization of the base change) of $\ca{S}^\Sigma_{K}$ from $\ca{O}_{(v)}'$ to $\ca{O}_{K_Z}$. The middle isomorphism is from Corollary \ref{zp-ab}. Then $\ca{S}^{\Sigma_2'}_{K_2,\ca{O}_{K_Z}}$ is a normal algebraic space that is proper and flat over $\ca{O}_{K_Z}$ and is representable by a projective scheme if $\Sigma$ is projective. 
Define $\ca{Z}_{[(\Phi_2,\sigma_2)],K_2,\ca{O}_{K_Z}}$ to be the normalization of the normalized base change $\ca{Z}_{[ZP^b(\Phi,\sigma)],K,\ca{O}_{K_Z}}$ in $\mrm{Z}_{[(\Phi_2,\sigma_2)],K_2,E'_{K_Z}}$ via the composition of morphisms $\mrm{Z}_{[(\Phi_2,\sigma_2)],K_2,E'_{K_Z}}\to \mrm{Z}_{[ZP^a(\Phi,\sigma)],K,E'_{K_Z}}\xrightarrow{\sim}\mrm{Z}_{[ZP^b(\Phi,\sigma)],K,E'_{K_Z}}\hookrightarrow \ca{Z}_{[ZP^b(\Phi,\sigma)],K,\ca{O}_{K_Z}}$. The middle isomorphism is also obtained from Corollary \ref{zp-ab}. Let $\ca{Z}_{[ZP^a(\Phi,\sigma)],K,\ca{O}_{K_Z}}$ be the normalization of $\ca{Z}_{[ZP^b(\Phi,\sigma)],K,\ca{O}_{K_Z}}$ in $\mrm{Z}_{[ZP^a(\Phi,\sigma)],K,\ca{O}_{K_Z}}$, which is isomorphic to $\ca{Z}_{[ZP^b(\Phi,\sigma)],K,\ca{O}_{K_Z}}$.
\end{enumerate}
\end{construction}
\begin{convention}
Let $X$ be a normal algebraic space over $\ca{O}_{(v)}'$. The symbol $X_{\ca{O}_{K_Z}}$ always denotes the normalized base change of $X$ from $\ca{O}_{(v)}'$ to $\ca{O}_{K_Z}$; these are the usual base changes if we can choose $F_{K_Z}$ such that $E'_{K_Z}$ is unramified over $E'$ at $p$.
\end{convention}
\begin{lem}\label{lem-strata-compactification-2}
With the construction above, if $\Phi_2\xrightarrow[\sim]{(\gamma,q)_{K_2}}\Phi_2'$ for $\Phi_2,\Phi_2'\in\ca{CLR}(G_2,X_2)$, $\gamma\in G_2(\bb{Q})$ and $q\in G_2(\A)$, we obtain an isomorphism $\ca{S}_{K_{\Phi_2}}(P_{\Phi_2},D_{\Phi_2})\to \ca{S}_{K_{\Phi_2}'}(P_{\Phi_2'},D_{\Phi_2'})$ extending the one over the generic fiber (see \cite[2.1.15]{Mad19}). In particular, $\Delta_{\Phi_2,K_2}$ acts equivariantly on the tower $\ca{S}_{K_{\Phi_2}}\to\overline{\ca{S}}_{K_{\Phi_2}}\to\ca{S}_{K_{\Phi_2},h}$, and the action of $\Delta_{\Phi_2,K_2}^\circ$ factors through a finite quotient.
\end{lem}
\begin{proof}
    Combine Proposition \ref{prop-mixsh-independent}, the construction above, and the following commutative diagram 
    \begin{equation*}
        \begin{tikzcd}
        \Phi_2\arrow[rr,"{(\gamma,q)_{K_2}}"]\arrow[d,"\pi^a"]&&\Phi_2'\arrow[d,"\pi^a"]\\
        \pi^a_*(\Phi_2)\arrow[rr,"{(\gamma,q)_K}"]&&\pi^a_*(\Phi_2').
        \end{tikzcd}
    \end{equation*}
The last claim is true since one can check this over the generic fiber and $\Delta^\circ_{\Phi_2,K_2}$ acts on the tower $\sh_{K_{\Phi_2}}\to\overline{\sh}_{K_{\Phi_2}}\to\sh_{K_{\Phi_2},h}$ through a finite group.
\end{proof}
Let us state the main theorem of this paper. 
\begin{thm}\label{maintheorem}
The triple $(G_2,X_2,K_{2,p})$ can be regarded as one of the Cases (HS), ($\text{STB}_n$) and (DL) in \S\ref{subsec-con-stra-comp-quo} by fixing some $K_{2,p}'$ containing $K_{2,p}$, where $K_{2,p}'$ is a hyperspecial subgroup or an intersection of $n$ Bruhat-Tits stabilizer subgroups containing $K_{2,p}$ as in Construction \ref{const-strata-compactification-2}. Let $K_2^p$ be a neat open compact subgroup of $G_2(\Ap)$ and write $K_2:=K_{2,p}K_2^p$. Let $\Sigma_2$ be an admissible cone decomposition of $(G_2,X_2,K_2)$. There is an admissible cone decomposition $\Sigma'_2$ refining $\Sigma_2$, which can be chosen to be smooth, projective, or both smooth and projective. Let $\ca{O}_2$ be $\ca{O}_{E_2,(v_2)}$ in Case (HS) and its deep-level case, and let $\ca{O}_2$ be $\ca{O}_{E_2,v_2}$ in Case ($\text{STB}_n$) and its deep-level case. With the constructions above, we have the following statements:
\begin{enumerate}[label=(\textrm{\ref{maintheorem}}.\arabic*)]
\item\label{maintheorem-1} There is an integral model $\ca{S}_{K_2}^{\Sigma_2^\prime}$ over $\ca{O}_2$ extending 
$\sh_{K_2}^{\Sigma_2^\prime}$ 
which is a normal algebraic space that is proper and flat over $\ca{O}_2$. For a suitable refinement $\Sigma_2'$, $\ca{S}_{K_2}^{\Sigma_2'}$ is representable by a normal scheme that is projective over $\ca{O}_2$.
\item\label{maintheorem-2} The compactification $\ca{S}_{K_2}^{\Sigma_2'}$ has a good stratification $$\ca{S}^{\Sigma_2'}_{K_2}=\disju\limits_{\Upsilon_2\in \cusp_{K_2}(G_2,X_2,\Sigma_2')}\ca{Z}_{\Upsilon_2,K_2}$$ 
such that each stratum $\ca{Z}_{\Upsilon_2,K_2}$ is a locally closed normal subscheme labeled by an $\Upsilon_2=[(\Phi_2,\sigma_2)]$, which is flat over $\ca{O}_2$. Recall that there is a partial order ``$\preceq$'' on $\cusp_{K_2}(G_2,X_2,\Sigma'_2)$ (see \cite[2.1.18]{Mad19} and \S\ref{tor-comp}). The closure $\overline{\ca{Z}}_{\Upsilon_2,K_2}$ of 
$\ca{Z}_{\Upsilon_2,K_2}$ in $\ca{S}^{\Sigma_2'}_{K_2}$ satisfies $\overline{\ca{Z}}_{\Upsilon_2,K_2}=\disju\limits_{\Upsilon_2^\prime\preceq \Upsilon_2}\ca{Z}_{\Upsilon_2',K_2}$. In particular, there is a natural open embedding $J^{\Sigma_2'}:\ca{S}_{K_2}\hookrightarrow\ca{S}_{K_2}^{\Sigma_2'}$ making $\ca{S}_{K_2}$ an open dense subscheme of $\ca{S}_{K_2}^{\Sigma_2'}$.
\item\label{maintheorem-3} For each cusp label representative $\Phi_2$ of $(G_2,X_2)$, the scheme $\ca{S}_{K_{\Phi_2}}$ is a normal scheme that is flat over $\ca{O}_2$ whose generic fiber is $\sh_{K_{\Phi_2},\ca{O}_2\otimes \bb{Q}}$. There is a tower of normal $\ca{O}_2$-schemes
\begin{equation}\label{eq-2-step}\begin{tikzcd}\ca{S}_{K_{\Phi_2}}\arrow[r,"\mbf{p}_1"]& \overline{\ca{S}}_{K_{\Phi_2}}\arrow[r,"\mbf{p}_2"]& \ca{S}_{K_{\Phi_2},h}.\end{tikzcd}\end{equation} 
The first morphism $\ca{S}_{K_{\Phi_2}}\to\overline{\ca{S}}_{K_{\Phi_2}}$ represents an $\mbf{E}_{K_{\Phi_2}}$-torsor and the second morphism is proper. 
In particular, since $\mbf{p}_1$ is an $\mbf{E}_{K_{\Phi_2}}$-torsor, we can construct $\ca{S}_{K_{\Phi_2}}(\sigma_2)$ as the twisted affine torus embedding of $\ca{S}_{K_{\Phi_2}}$ over $\overline{\ca{S}}_{K_{\Phi_2}}$ with respect to the affine torus embedding $\mbf{E}_{K_{\Phi_2}}\hookrightarrow \mbf{E}_{K_{\Phi_2}}(\sigma_2)$, and we can construct $\ca{S}_{K_{\Phi_2},\sigma_2}$ as the $\sigma_2$-stratum of $\ca{S}_{K_{\Phi_2}}(\sigma_2)$. Also, it makes sense to define $\ca{S}_{K_{\Phi_2}}(\Sigma'_2)$.\par
Moreover, let $\ca{S}^*_{K_{\Phi_2}}$ (resp. $\overline{\ca{S}}^*_{K_{\Phi_2}}$ and resp. $\ca{S}^*_{K_{\Phi_2},h}$) be the quotient of $\ca{S}_{K_{\Phi_2}}$ (resp. $\overline{\ca{S}}_{K_{\Phi_2}}$ and resp. $\ca{S}^*_{K_{\Phi_2},h}$) by $\Delta^\circ_{\Phi_2,K_2}$. Then $\ca{S}_{K_{\Phi_2}}^*\to\overline{\ca{S}}^*_{K_{\Phi_2}}$ is also a torsor under a split torus.
\item\label{maintheorem-4} For any cusp label representative $\Phi_2$ of $(G_2,X_2)$, the group $\Delta^\circ_{\Phi_2,K_2}$ acts on $\ca{S}_{K_{\Phi_2},\sigma_2}$ and $\ca{S}_{K_{\Phi}}(\sigma_2)$ through a finite quotient. We have an isomorphism $\Delta_{\Phi_2,K_2}^\circ\bss \ca{S}_{K_{\Phi_2},\sigma_2}\iso \ca{Z}_{[(\Phi_2,\sigma_2)],K_2}$, where $\ca{S}_{K_{\Phi_2},\sigma_2}\hookrightarrow \ca{S}_{K_{\Phi_2}}(\sigma_2)$ is the closed stratum of $\sigma_2$. The scheme $\ca{S}_{K_{\Phi_2}}(\sigma_2)$ has a natural stratification defined by the cone decomposition $\Sigma_2'(\Phi_2)$.\par
We have a strata-preserving isomorphism 
    $$\Delta^\circ_{{\Phi_2},K_2}\bss\cpl{\ca{S}_{K_{\Phi_2}}(\sigma_2)}{\ca{S}_{K_{\Phi_2},\sigma_2}}\iso \cpl{\ca{S}^{\Sigma_2'}_{K_2}}{\ca{Z}_{[(\Phi_2,\sigma_2)],K_2}}.$$
Moreover, let $\ca{S}^*_{K_{\Phi_2}}(\sigma_2)$ (resp. $\ca{S}^*_{K_{\Phi_2},\sigma_2}$) be the quotient of $\ca{S}_{K_{\Phi_2}}(\sigma_2)$ (resp. $\ca{S}_{K_{\Phi_2,\sigma_2}}$) by $\Delta^\circ_{\Phi_2,K_2}$. We have a strata-preserving isomorphism 
  $$\Delta^\circ_{{\Phi_2},K_2}\bss\cpl{\ca{S}_{K_{\Phi_2}}(\sigma_2)}{\ca{S}_{K_{\Phi_2},\sigma_2}}\iso \cpl{\ca{S}^{\Sigma_2'}_{K_2}}{\ca{Z}_{[(\Phi_2,\sigma_2)],K_2}}\iso \cpl{\ca{S}^*_{K_{\Phi_2}}(\sigma_2)}{\ca{S}^*_{K_{\Phi_2},\sigma_2}}.$$
Here, ``strata-preserving'' has the following meaning: For any affine open formal subscheme $\spf(A,I)$ of $\cpl{\ca{S}^{\Sigma_2'}_{K_2}}{\ca{Z}_{[(\Phi_2,\sigma_2)],K_2}}$ which canonically induces morphisms $\mathsf{c}_1:\spec A\to \ca{S}_{K_{\Phi_2}}^*(\sigma_2)$ and $\mathsf{c}_2:\spec A\to \ca{S}^{\Sigma_2'}_{K_2}$, the two stratifications on $\spec A$ defined by the pullback of the stratification on $\ca{S}_{K_{\Phi_2}}(\sigma_2)$ under $\mathsf{c}_1$ and the pullback of the stratification on $\ca{S}_{K_2}^{\Sigma_2'}$ under $\mathsf{c}_2$ coincide. 
\item\label{maintheorem-5} In Case (HS), each $$\ca{S}_{K_{\Phi_2,p}}:=\varprojlim_{K_{\Phi_2}^p\text{ neat open compact}}\ca{S}_{K_{\Phi_2,p}K_{\Phi_2}^p}$$ is an inverse limit of schemes that are smooth over $\ca{O}_{2}$ and has the extension property. In particular, the tower $\ca{S}_{K_{\Phi_2,p}}$ does not depend on the choices made at the beginning of \S\ref{sss-main-setup-again}. 
\end{enumerate}
\end{thm}
\begin{rk}\label{rk-field-of-def}
\begin{itemize}
\item Our construction works for any choice of $(G_0,X_0)$ that satisfies Condition (3) in \cite[Lem. 4.6.22]{KP15}; therefore, if, for some $(G_2,X_2)$, one can choose a $(G_0,X_0)$ such that $E_0\sbst E_2$, then $E'=E_2$, and $\ca{S}_{K_2}^{\Sigma_2'}$ and $\ca{S}_{K_{\Phi_2}}$ in the theorem are automatically defined over $\spec \ca{O}_{E_2,(v_2)}$. 
\item The group $\Delta^\circ_{\Phi_2,K_2}$ and the action of it on $\ca{S}_{K_{\Phi_2}}$ are not trivial for a general abelian-type Shimura datum $(G_2,X_2)$ and a general $K_2$. Our result \ref{maintheorem-3} shows that both $\ca{S}_{K_{\Phi_2}}$ and $\ca{S}_{K_{\Phi_2}}^*$ are torsors under split tori. Note that instead of $\{\ca{S}_{K_{\Phi_2}}(\sigma_2)\}$, $\{\ca{S}_{K_{\Phi_2}}^*(\sigma_2)\}$ is the collection of affine toric schemes attached to the toroidal embedding $\ca{S}_{K_2}\hookrightarrow \ca{S}_{K_2}^{\Sigma_2'}$ in the sense of \cite[p.54, Def. 1]{KKMS73}.
\item One can also formulate and prove a refinement of \ref{maintheorem-4} similar to \cite[Prop. 2.1.3]{LS18b}. Note that this refined assertion was used in \cite{LS18b} and \cite{Mao25b} to study the well-positionedness of some subschemes on compactifications of integral models of PEL- and Hodge-type Shimura varieties. 
\end{itemize}
\end{rk}
\begin{proofof}[Theorem \ref{maintheorem}: first reductions] It suffices to show all of the statements above over $\ca{O}'_{(v)}$. 
The statement \ref{maintheorem-2} will be shown along the way in the proof of 
\ref{maintheorem-1} and \ref{maintheorem-4} (see below). 
For \ref{maintheorem-3}, it suffices to show that $\ca{S}_{K_{\Phi_2}}\to \overline{\ca{S}}_{K_{\Phi_2}}$ and $\ca{S}^*_{K_{\Phi_2}}\to\overline{\ca{S}}^*_{K_{\Phi_2}}$ are torsors under $\mbf{E}_{K_{\Phi_2}}$ and $\mbf{E}_{\K_{\Phi_2}}$, respectively; other claims in it follow from Construction \ref{const-strata-compactification}, \ref{const-strata-compactification-2} and Lemma \ref{lem-strata-compactification-2}. 
\end{proofof}
\subsubsection{Proof of \ref{maintheorem-3}.}\label{subsubsec-main-3}
\S\ref{subsubsec-main-3} and \S\ref{subsubsec-s*true} are devoted to the proof of \ref{maintheorem-3}.\par 
First, we need the following preparation. 
\begin{lem}[{cf. \cite[Lem. 3.7.2]{Mad19} and \cite[Prop. 8.7]{Lan16b}}]\label{lem-torus-hodge}
Let $K_{0,p}$ be an open compact subgroup contained in some intersection $K_{0,p}'\sbst G_0(\bb{Q}_p)$ of
$n$ Bruhat-Tits stabilizer subgroups. Let $K^p_0$ be a neat open compact subgroup in $G_0(\Ap)$. Let $K_0:=K_{0,p}K_0^p$ and $K_0':=K_{0,p}'K_0^p$. 
Let $\Phi_0$ be a cusp label representative of $(G_0,X_0)$ mapping to $\Phi^\ddag\in \ca{CLR}(G^\ddag,X^\ddag)$. Let $K_{\Phi_0}$ (resp. $K_{\Phi_0}'$) be the intersection $\lcj{K_0}{g_{\Phi_0}}\cap P_{\Phi_0}(\A)$ (resp. $\lcj{K'_0}{g_{\Phi_0}}\cap P_{\Phi_0}(\A)$). As before, suppose that $\ca{S}_{K_{\Phi_0}}$ is constructed from the relative normalization of $\ca{S}_{K_{\Phi_0}'}$ in $\sh_{K_{\Phi_0}}$, and that $\ca{S}_{K_{\Phi_0}'}$ is constructed by choosing $K_p^\ddag\sbst G^\ddag(\bb{Q}_p)$ and a Hodge embedding as in ($\text{STB}_n$), and by taking the relative normalization of $\ca{S}_{K_{\Phi^\ddag}}$ in $\sh_{K_{\Phi_0}'}$. Then there is a finite flat cover $U$ of $\overline{\ca{S}}_{K_{\Phi_0}'}$ and an $\mbf{E}_{K_{\Phi_0}}$-torsor $\ca{S}$ over $U$, such that the pullback of the torsor $\ca{S}\to U$ to $U':=U\times_{\overline{\ca{S}}_{K'_{\Phi_0}}}\overline{\ca{S}}_{K_{\Phi_0}}$ is isomorphic to the pullback of the torsor $\ca{S}_{K_{\Phi_0}}\to\overline{\ca{S}}_{K_{\Phi_0}}$ to the finite flat cover $U'$ of the target.
\end{lem}
\begin{proof}
We combine \cite[Lem. 3.7.2]{Mad19} with \cite[Thm. 4.1.5]{Mad19}. 
Choose some neat open compact $K^{\ddag,\prime}:=K_p^{\ddag,\prime}K^{\ddag,\prime,p}\sbst K^\ddag$ containing $K_0$ such that $\sh_{K_{\Phi_0}}\to\sh_{K_{\Phi^\ddag}'}$ is a closed embedding, where $K_{\Phi^\ddag}':=P_{\Phi^\ddag}(\A)\cap \lcj{K^{\ddag,\prime}}{g_{\Phi^\ddag}}$. This can be achieved by Corollary \ref{cor-pure-embedding} and Lemma \ref{lem-emb-mix}. 
Then $\mbf{E}_{K_{\Phi_0}}\to \mbf{E}_{K_{\Phi^\ddag}'}$ is a closed embedding. By \cite[Thm. 4.1.5(4)]{Mad19}, $\ca{S}_{K_{\Phi_0}}\to \ca{S}_{K_{\Phi^\ddag}'}$ is equivariant under $\mbf{E}_{K_{\Phi_0}}\to \mbf{E}_{K_{\Phi^\ddag}'}$.
By the proof of \cite[Lem. 3.7.2]{Mad19}, there is a finite flat cover $U^\ddag$ of $\overline{\ca{S}}_{K_{\Phi^\ddag}}$, a split torus $\mbf{E}$ and an $\mbf{E}$-torsor $\ca{S}^\ddag$ over $U^\ddag$ with an isogeny between split tori $\mbf{E}\to \mbf{E}_{K_{\Phi^\ddag}'}$, such that the pullback of $\ca{S}^\ddag\times^{\mbf{E}}\mbf{E}_{K'_{\Phi^\ddag}}\to U^\ddag$ to $U^{\ddag,\prime}:=\overline{\ca{S}}_{K'_{\Phi^\ddag}}\times_{\overline{\ca{S}}_{K_{\Phi^\ddag}}} U^\ddag$ is isomorphic to the pullback of $\ca{S}_{K'_{\Phi^\ddag}}\to\overline{\ca{S}}_{K'_{\Phi^\ddag}}$ to $U^{\ddag,\prime}$.
There is a projection $\mbf{S}_{K'_{\Phi^\ddag}}\to\mbf{S}_{K_{\Phi_0}}$ between character groups associated with $\mbf{E}_{K'_{\Phi^\ddag}}$ and $\mbf{E}_{K_{\Phi_0}}$.\par 
Let $$U:=U^\ddag\times_{\overline{\ca{S}}_{K_{\Phi^\ddag}}}\overline{\ca{S}}_{K'_{\Phi_0}},$$ and let $$\ca{S}^\ddag_1:=\ca{S}_{K_{\Phi^\ddag}}\times_{\overline{\ca{S}}_{K_{\Phi^\ddag}}} \overline{\ca{S}}_{K'_{\Phi_0}}.$$
Moreover, let $\ca{S}_0^\ddag$ be the pullback of the torsor $\ca{S}^\ddag\times^{\mbf{E}}\mbf{E}_{K'_{\Phi^\ddag}}$ over $U^\ddag$ along $U\to U^\ddag$.\par
Then there is a diagram of torsors
\begin{equation*}
    \begin{tikzcd}
    &\ca{S}_0^\ddag\arrow[d]\\
    \ca{S}_{K'_{\Phi_0}}\arrow[r]&\ca{S}^\ddag_1,
    \end{tikzcd}
\end{equation*}
where the bottom one is equivariant under $\mbf{E}_{K_{\Phi_0}'}\to \mbf{E}_{K_{\Phi^\ddag}}$ over $\overline{\ca{S}}_{K'_{\Phi_0}}$ (and the last map between tori can be chosen as an embedding by Corollary \ref{cor-pure-embedding} and Lemma \ref{lem-emb-mix} again), and the vertical arrow is equivariant under $\mbf{E}_{K'_{\Phi^\ddag}}\to \mbf{E}_{K_{\Phi^\ddag}}$.\par
We now take the inverse image of the $\mbf{E}_{K'_{\Phi_0}}$-torsor $\ca{S}_{K'_{\Phi_0}}$ in the $\mbf{E}_{K'_{\Phi^\ddag}}$-torsor $\ca{S}_0^\ddag$ by pulling back under the vertical map in the displayed diagram. This is an $\mbf{E}_{K_{\Phi_0}}$-torsor $\ca{S}$ over $U$.\par
Then the torsor $\ca{S}\to U$ is the desired one, as $\ca{S}_{K_{\Phi_0}}$ is an $\mbf{E}_{K_{\Phi_0}}$-torsor by \cite[Thm. 4.1.5]{Mad19} and has an $\mbf{E}_{K_{\Phi_0}}$-equivariant morphism to $\ca{S}$ by checking using the universal property of fiber products in the construction.
\end{proof}
\begin{rk}\label{rk-torus-hodge}
If there is another open compact subgroup $K_{0,p}''\sbst G_0(\bb{Q}_p)$ such that $K_{0,p}\sbst K_{0,p}''\sbst K_{0,p}'$. Let $K'':=K_{0,p}''K^p_0$ and $K''_{\Phi_0}:=\lcj{K_0''}{g_{\Phi_0}}\cap P_{\Phi_0}(\A)$. Then we can pull $U$ back to $U'':=\overline{\ca{S}}_{K_{\Phi_0}''}(P_{\Phi_0},D_{\Phi_0})\times_{\overline{\ca{S}}_{K'_{\Phi_0}}}U$ and also pull the torsor $\ca{S}$ over $U$ back to $U''$. We denote this pullback of $\ca{S}$ by $\ca{S}''$. Then Lemma \ref{lem-torus-hodge} above implies that the pullback of $\ca{S}''$ to $U'$ is isomorphic to the pullback of $\ca{S}_{K_{\Phi_0}}\to \overline{\ca{S}}_{K_{\Phi_0}}$ to the finite flat cover $U'$.
\end{rk}
\subsubsection{}\label{subsubsec-s*true}
Let us show the following lemma. The proof of it occupies all paragraphs below until Lemma \ref{lem-torus}.
\begin{lem}\label{lem-maintorus-s*}
The morphism $\ca{S}_{\K_{\Phi}^a}\to \overline{\ca{S}}_{\K_\Phi^a}$ is an $\mbf{E}_{\K_\Phi}$-torsor. The assertion in \ref{maintheorem-3} that $\ca{S}^*_{K_{\Phi_2}}\to\overline{\ca{S}}^*_{K_{\Phi_2}}$ is a torsor under a split torus is true.  
\end{lem}
By Proposition \ref{ext-imp}, we see that it suffices to show the statement for $\ca{S}_{\K^a_\Phi}\to \overline{\ca{S}}_{\K^a_\Phi}$ since $\ca{S}_{K_{\Phi_2}}^*\iso \overline{\ca{S}}_{K_{\Phi_2}}^*\times_{\overline{\ca{S}}_{\K^a_\Phi}}\ca{S}_{\K^a_\Phi}$ by our choice of $K_pK^p$ satisfying the proposition.\par
Recall that we have already shown that  $\ca{S}_{\K^b_\Phi}\to\overline{\ca{S}}_{\K^b_\Phi}$ is an $\mbf{E}_{\K_\Phi}$-torsor for all cases in Part \ref{prop-mixsh-canonicity-1} of Proposition \ref{prop-mixsh-canonicity}. As in Construction \ref{const-strata-compactification}, we have proved that $\ca{S}_{\K^a_\Phi}\to \overline{\ca{S}}_{\K^a_\Phi}$ is an $\mbf{E}_{\K_\Phi}$-torsor, when we can choose $F_Z$ such that $\ca{O}_{K_Z}$ is unramified over $\ca{O}'_{(v)}$. This is the case when $K_p$ contains the parahoric subgroup of $Z(\bb{Q}_p)$. For general $K_p$, we know from Construction \ref{const-strata-compactification-2} that $\ca{S}_{\K^a_\Phi,\ca{O}_{K_Z}}\to\overline{\ca{S}}_{\K^a_\Phi,\ca{O}_{K_Z}}$ is an $\mbf{E}_{\K_\Phi}$-torsor.\par
Recall that we fixed some intersection of Bruhat-Tits stabilizer subgroups $K_p'$ containing $K_p$. \par
Denote $K':=K_p'K^p$. Let $K_{Z,p}^+\sbst Z_G(\bb{Q}_p)$ be the unique parahoric subgroup of $Z_G(\bb{Q}_p)$. Denote $K_p^+:=K_pK_{Z,p}^+$ and $K^+:=K_p^+K^p$. Let $\K^+_\Phi:= ZP_\Phi(\A)\cap g_\Phi K^+g_\Phi^{-1}$. Note that $\K^+_\Phi$ contains the parahoric subgroup of $Z_G(\bb{Q}_p)$.
Then the scheme $\ca{S}_{\K^{+,i}_\Phi}:=\ca{S}_{\K^+_\Phi}(ZP^i_\Phi,ZP^i_\Phi(\bb{Q})D_{\Phi})$ obtained as in Construction \ref{const-strata-compactification-2} is a torsor over $\overline{\ca{S}}_{\K^{+,i}_\Phi}:=\overline{\ca{S}}_{\K^+_\Phi}(ZP^i_\Phi,ZP^i_\Phi(\bb{Q})D_\Phi)$ under $\mbf{E}_{\K^+_\Phi}$. 
Note that $\mbf{E}_{\K^+_\Phi}$ is determined by the lattice $\mbf{\Lambda}_{\K_\Phi^+}$, the projection to the second factor of $(Z_G(\bb{Q})\times U_\Phi(\bb{Q}))\cap \K^+_\Phi$, and there is an isogeny $\mbf{E}_{\K_\Phi}\to \mbf{E}_{\K^+_\Phi}$.
We then have a finite map $\ca{S}_{\K_\Phi^a}\to  \overline{\ca{S}}_{\K_\Phi^a}\times_{\overline{\ca{S}}_{\K^{+}_\Phi}} \ca{S}_{\K^{+}_\Phi}$ that is equivariant under $\mbf{E}_{\K_\Phi}\to\mbf{E}_{\K^+_\Phi}$.\par 
Recall that by construction we have a finite index $I_{ZP_\Phi/P_0,K}$ mapping to another finite index $I_{ZP_\Phi/P_0,K^+}$. Without loss of generality, we pick some $e_\beta$ for $\beta\in I_{ZP_\Phi/P_0,K}$ mapping to $\beta_1\in I_{ZP_\Phi/P_0,K^+}$, and let $e_{\beta}=e_{\beta_1}$.

Denote $\Delta^\beta:=\Delta_{\rcj{K}{\beta}}(P_0,ZP_\Phi)$ and $\Delta^{\beta_1}:=\Delta_{{K^{\prime,\beta_1}}}(P_0,ZP_\Phi)$.
The quotient $\ca{S}_{\K^b_\Phi,\beta}:=\ca{S}_{K_{\Phi_{\beta'}}}/\Delta^\beta$ (resp. $\overline{\ca{S}}_{\K^b_\Phi,\beta}:=\overline{\ca{S}}_{K_{\Phi_{\beta'}}}/\Delta^{\beta}$) is an open and closed subscheme of $\ca{S}_{\K^b_\Phi}$ (resp. $\overline{\ca{S}}_{\K^b_\Phi}$). 
Similarly, the quotient $\ca{S}_{\K^{+,b}_\Phi,\beta_1}:=\ca{S}_{K_{\Phi_{\beta_1'}}}/\Delta^{\beta_1}$ (resp. $\overline{\ca{S}}_{\K^{+,b}_\Phi,\beta_1}:=\overline{\ca{S}}_{K_{\Phi_{\beta_1'}}}/\Delta^{\beta_1}$) is an open and closed subscheme of $\ca{S}_{\K_\Phi^{+,b}}$.
Then $\ca{S}_{\K_\Phi^b,\beta}$ (resp. $\overline{\ca{S}}_{\K_\Phi^b,\beta}$) maps surjectively to $\ca{S}_{\K_\Phi^{+,b},\beta_1}$ (resp. $\overline{\ca{S}}_{\K_\Phi^{+,b},\beta_1}$).\par
The action of $\Delta^{\beta_1}$ on $\overline{\ca{S}}_{K_{\Phi_{\beta_1'}}}$ factors through a finite quotient $\overline{H}_{\beta_1}$, and the action of the latter group is free by Lemma \ref{lem-twt-free-action}; denote by $\ca{K}_{\beta_1}$ the kernel of this action. Let $\Delta_1:=\Delta^{\beta}\cap \ca{K}_{\beta_1}$. Let $\Delta^1:=\Delta^{\beta}/\Delta_1$.
Define $\overline{\ca{S}}_1:=\overline{\ca{S}}_{K_{\Phi_{\beta'}}}/\Delta_1$. Define $\overline{\ca{S}}_2:=\overline{\ca{S}}_{K_{\Phi_{\beta_1'}}}/\Delta^1.$
\begin{lem}\label{lem-torus-cartesian-diag}
There is a Cartesian commutative diagram
\begin{equation}\label{diag-cartesian}
   \begin{tikzcd}
    \overline{\ca{S}}_1\arrow[r,"f_1"]\arrow[d,"g_1"]&\overline{\ca{S}}_{K_{\Phi_{\beta_1'}}}\arrow[d,"g_2"]\\
    \overline{S}_{\K^b_\Phi,\beta}\arrow[r,"f_2"]&\overline{\ca{S}}_2,
   \end{tikzcd} 
\end{equation}
which fits into a commutative diagram
\begin{equation}\label{diag-quotient-decomposition}
   \begin{tikzcd}
\overline{\ca{S}}_{K_{\Phi_{\beta'}}}\arrow[rr,"h_1"]\arrow[rrrr,bend left=20,"\mrm{can.}"]\arrow[d,"\pi^b"]&&\overline{\ca{S}}_1\arrow[rr,"f_1"]\arrow[dll,"g_1"]&&\overline{\ca{S}}_{K_{\Phi_{\beta_1'}}}\arrow[dll,"g_2"]\arrow[d,"\pi^b"]\\
\overline{\ca{S}}_{\K^b_\Phi,\beta}\arrow[rr,"f_2"]\arrow[rrrr,bend right=20,"\mrm{can.}"]&&\overline{
\ca{S}}_2\arrow[rr,"h_2"]&&\overline{\ca{S}}_{\K_\Phi^{+,b},\beta_1}.
   \end{tikzcd} 
\end{equation}
The compositions of horizontal arrows of the diagram above are canonical transition maps induced by relative normalizations. 
\end{lem}
\begin{proof}
The compositions of finite surjective morphisms between normal schemes $\overline{\ca{S}}_{K_{\Phi_{\beta'}}}\to\overline{\ca{S}}_1\xrightarrow{g_1}\overline{\ca{S}}_{\K^b_\Phi,\beta}$ and $\overline{\ca{S}}_{K_{\Phi_{\beta_1'}}}\xrightarrow{g_2}\overline{\ca{S}}_2\to \overline{\ca{S}}_{K^{+,b}_\Phi,\beta_1}$ are finite {\'e}tale. So $g_1$ and $g_2$ are finite {\'e}tale. Moreover, (\ref{diag-cartesian}) is Cartesian since $g_1$ and $g_2$ are locally finite free of the same degree by Lemma \ref{lem-twt-free-action}.
\end{proof}
By construction, $\Delta^{\beta_1}$ is generated by $\Delta^\beta$ and $\Delta^+:=\ker(\ag(P_0)\to\ag(ZP_\Phi)/K_{Z,p}^+)$, and both $\Delta^\beta$ and $\Delta^+$ are normal subgroups of $\Delta^{\beta_1}$; $\Delta^+$ is finite modulo $K_{\Phi_{\beta_1'}}\cap \Delta^+$ and satisfies the assumption in \S\ref{subsubsec-ggamma-trivial-case} (that $h\gamma^{-1}$ is trivial in $G^\ad_0(\A)$).\par
By Lemma \ref{lem-torus-hodge} and Remark \ref{rk-torus-hodge} above and by taking contracted products of $\mbf{E}_{K_{\Phi_{\beta'}}}$-torsors along $\mbf{E}_{K_{\Phi_{\beta'}}}\to \mbf{E}_{\K_{\Phi}}$, there is a finite flat cover $V(\beta_1)$ of $\overline{\ca{S}}_{\K_{\Phi_{\beta_1'}}}$ and an $\mbf{E}_{\K_\Phi}$-torsor $\ca{T}(\beta_1)$ over $V(\beta_1)$, such that the pullback of the torsor $\ca{T}(\beta_1)\to V(\beta_1)$ to $V'(\beta):=\overline{\ca{S}}_{K_{\Phi_{\beta'}}}\times_{\overline{\ca{S}}_{K_{\Phi_{\beta_1'}}}}V(\beta_1)$ is isomorphic to the pullback of $\ca{S}_{K_{\Phi_{\beta'}}}\times^{\mbf{E}_{K_{\Phi_{\beta'}}}}\mbf{E}_{\K_\Phi}$ over $\overline{\ca{S}}_{K_{\Phi_{\beta'}}}$ to $V'(\beta)$.\par
Let $V(\beta):=V(\beta_1)\times_{\overline{\ca{S}}_{K_{\Phi_{\beta_1'}}},f_1}\overline{\ca{S}}_1$. 
Taking the quotient of $\Delta_1$ and by Lemma \ref{lem-torus-cartesian-diag} above, we have that 
the pullback of $\ca{S}_{\K^b_\Phi,\beta}\to \overline{\ca{S}}_{\K^b_\Phi,\beta}$ along the finite flat cover $V(\beta)\to\overline{\ca{S}}_1\xrightarrow{g_1}\overline{\ca{S}}_{\K^b_\Phi,\beta}$ is isomorphic to the pullback of $\ca{T}(\beta_1)\to V(\beta_1)$ to $V(\beta)$. 
Moreover, we have $V(\beta)\times_{V(\beta_1)}\ca{T}(\beta_1)\iso \overline{\ca{S}}_{\K^b_\Phi,\beta}\times_{f_2,\overline{\ca{S}}_2}\ca{T}(\beta_1)$ by Lemma \ref{lem-torus-cartesian-diag}. (The morphism from $\ca{T}(\beta_1)$ to $\overline{\ca{S}}_2$ is given by $\ca{T}(\beta_1)\to V(\beta_1)\to \overline{\ca{S}}_{K_{\Phi_{\beta_1'}}}\xrightarrow{g_2}\overline{\ca{S}}_2$.)\par
In conclusion, we have a canonical isomorphism 
$$\ca{S}_{\K^b_\Phi,\beta}\times_{\overline{\ca{S}}_{\K^b_\Phi,\beta}}V(\beta)\iso \overline{\ca{S}}_{\K^b_\Phi,\beta}\times_{f_2,\overline{\ca{S}}_2}\ca{T}(\beta_1).$$
This isomorphism is defined over $\ca{O}'_{(v)}$. 
\begin{lem}\label{lem-cover-invariant}
With the conventions above, there is a finite flat cover $V(\beta_1)$ over $\overline{\ca{S}}_{\K_{\Phi_{\beta_1'}}}$ and an $\mbf{E}_{\K_\Phi}$-torsor $\ca{T}(\beta_1)$ over $V(\beta_1)$ such that the sequence $\ca{T}(\beta_1)\to V(\beta_1)\to \overline{\ca{S}}_{\K_{\Phi_{\beta_1'}}}$ is equivariant under the action of $\Delta^+$.
\end{lem}
\begin{proof}
Since $\Delta^+\sbst \Delta^{\beta_1}$ satisfies the assumption in \S\ref{subsubsec-ggamma-trivial-case} and is finite modulo $K_{\Phi_{\beta_1'}}\cap \Delta^+$, we can fix a sufficiently large finite Galois extension $F$ over $\bb{Q}$, and replace $(G^\ddag,X^\ddag)$ in Lemma \ref{lem-torus-hodge} with $(G^F,X^F)$. In the notation of \S\ref{subsubsec-ggamma-trivial-case}, $\Delta^+$ acts on the integral models $\ca{S}_{K_{\Phi^F}}\to\overline{\ca{S}}_{K_{\Phi^F}}\to\ca{S}_{K_{\Phi^F},h}$ associated with the cusp label representative $\Phi^F$ of $(G^F,X^F)$ and the universal $1$-motive on $\ca{S}_{K_{\Phi^F}}$. 
We now adopt the notation in Lemma \ref{lem-torus-hodge}. Suppose $\Phi_{\beta_1'}$ maps to $\Phi^F$. By \cite[Lem. 3.7.2]{Mad19}, the finite flat cover $U^\ddag$ over $\ca{S}_{K_{\Phi^F}}$ is $\Hom(\frac{1}{n}Y_{\Q_{\Phi^F}}, A_{\Q_{\Phi^F}})$, where $n$ is some positive integer and $Y_{\Q_{\Phi^F}}$ (resp. $A_{Q_{\Phi^F}}$) is the $(-2)$-graded piece (resp. $(-1)$-graded piece) of the universal $1$-motive $\Q_{\Phi^F}$ over $\ca{S}_{K_{\Phi^F}}$, and the torsor $\ca{S}_0^\ddag$ is a quotient of $\uhom^{\mrm{symm}}(\frac{1}{n}Y_{\Q_{\Phi^F}},\G^\natural_{\Q_{\Phi^F}})$ (see \cite[3.7.2.2]{Mad19}).
Then, by Lemma \ref{lem-lan-qisog-twt}, the sequence $\ca{S}_0^\ddag\times^{\mbf{E}_{K_{\Phi_0}}}\mbf{E}_{\K_\Phi}\to U^\ddag\to \overline{\ca{S}}_{K_{\Phi^F}}$ is equivariant under the action of $\Delta^+$ after pulling back to $\ca{S}_{\Phi_{\beta_1'}}$. We now have the desired assertion.
\end{proof}
Now we can and we will assume that $\ca{T}(\beta_1)\to V(\beta_1)$ satisfies both Lemma \ref{lem-cover-invariant} and Remark \ref{rk-torus-hodge}.
\begin{lem}\label{lem-torsor-iso}
With the conventions above, the pullback of $V(\beta_1)\to \overline{\ca{S}}_{K_{\Phi_{\beta_1'}}}\xrightarrow{\pi^b}\overline{\ca{S}}_{\K^{+,b}_\Phi,\beta_1}$ along the transition map $\overline{\ca{S}}_{\K^b_\Phi,\beta}\to \overline{\ca{S}}_{\K_\Phi^{+,b},\beta_1}$ is a finite flat cover $U(\beta)$ of $\overline{\ca{S}}_{\K^b_\Phi,\beta}$. The pullback of $\ca{T}(\beta_1)$ to $V_\beta$ is canonically isomorphic to the pullback of $\ca{S}_{\K^b_\Phi,\beta}$ to $U(\beta)$.
\end{lem}
\begin{proof}
The first claim is immediate from construction. 
The diagram (\ref{diag-quotient-decomposition}) is equivariant under any $\delta\in \Delta^+$. Combining this with Lemma \ref{lem-cover-invariant}, we have the following commutative diagram:
\begin{equation}\label{diag-quotient-decomposition-delta}
 \begin{tikzcd}
 &&&&&&\ca{T}(\beta_1)\arrow[d,"\pi_1"]\arrow[dr,"\delta"]&\\
&&&&&&V(\beta_1) \arrow[d,"\pi_2"]\arrow[dr,"\delta"]&\ca{T}(\beta_1)\arrow[d,"\pi_1"]\\
\overline{\ca{S}}_{K_{\Phi_{\beta'}}}\arrow[rrr,"h_1"]\arrow[dd,"\pi^b"]\arrow[dr,"\delta"]&&&\overline{\ca{S}}_1\arrow[rrr,"f_1"]\arrow[ddlll]\arrow[dr,"\delta"]&&&\overline{\ca{S}}_{K_{\Phi_{\beta_1'}}}\arrow[ddlll,"g_2"]\arrow[dd,"\pi^b"]\arrow[dr,"\delta"]&V(\beta_1)\arrow[d,"\pi_2"]\\
&\overline{\ca{S}}_{K_{\Phi_{\beta'}}}\arrow[rrr,"h_1"]\arrow[dd,"\pi^b"']&&&\overline{\ca{S}}_1\arrow[rrr,"f_1"]\arrow[ddlll,"g_1"']&&&\overline{\ca{S}}_{K_{\Phi_{\beta_1'}}}\arrow[ddlll]\arrow[dd,"\pi^b"]\\
\overline{\ca{S}}_{\K^b_\Phi,\beta}\arrow[rrr,"f_2"]\arrow[dr,"\delta"]&&&\overline{
\ca{S}}_2\arrow[rrr,"h_2"]\arrow[dr,"\delta"]&&&\overline{\ca{S}}_{\K_\Phi^{+,b},\beta_1}\arrow[dr,"\mrm{Id}"]&\\
&\overline{\ca{S}}_{\K^b_\Phi,\beta}\arrow[rrr,"f_2"]&&&\overline{
\ca{S}}_2\arrow[rrr,"h_2"]&&&\overline{\ca{S}}_{\K_\Phi^{+,b},\beta_1}.
   \end{tikzcd} 
\end{equation}
From now on, we will freely use the conventions appearing in (\ref{diag-quotient-decomposition-delta}). 
The quotient $\Delta_0:=\Delta^{\beta_1}/\Delta^\beta\cdot\ca{K}_{\beta_1}$ acts on $\overline{\ca{S}}_2$ freely. Since $\Delta^+$ and $\Delta^\beta$ generate $\Delta^{\beta_1}$ and are normal subgroups, $\Delta_0$ is a quotient of $\Delta^+$. Hence, the fiber product $\overline{\ca{S}}_{\K_\Phi^b,\beta}\times_{\overline{\ca{S}}_{\K_\Phi^{+,b},\beta_1}}\ca{T}(\beta_1)$ is isomorphic to a disjoint union
$$\disju_{[\delta]\in \Delta_0}\overline{\ca{S}}_{\K^b_\Phi,\beta}\times_{f_2,\overline{\ca{S}}_2,[\delta]\circ g_2\circ\pi_2\circ\pi_1}\ca{T}(\beta_1).$$
Similarly, the fiber product $U(\beta)=\overline{\ca{S}}_{\K_\Phi^b,\beta}\times_{\overline{\ca{S}}_{\K_\Phi^{+,b},\beta_1}}V(\beta_1)$ is isomorphic to a disjoint union
$$\disju_{[\delta]\in \Delta_0}\overline{\ca{S}}_{\K^b_\Phi,\beta}\times_{f_2,\overline{\ca{S}}_2,[\delta]\circ g_2\circ\pi_2}V(\beta_1).$$
From (\ref{diag-quotient-decomposition-delta}) and Lemma \ref{lem-torus-cartesian-diag}, we see that $\overline{\ca{S}}_{\K_\Phi^b,\beta}\times_{\overline{\ca{S}}_{\K_\Phi^{+,b},\beta_1}}\ca{T}(\beta_1)$ is isomorphic to a disjoint union
$$\disju_{\{\delta\}}\delta^*(\overline{\ca{S}}_{\K^b_\Phi,\beta}\times_{f_2,\overline{\ca{S}}_2}\ca{T}(\beta_1))\iso \disju_{\{\delta\}}\delta^*(\ca{S}_{\K^b_\Phi,\beta}\times_{\overline{\ca{S}}_{\K^b_\Phi,\beta}}V(\beta) )\iso\disju_{\{\delta\}}  \ca{S}_{\K^b_\Phi,\beta}\times_{\overline{\ca{S}}_{\K^b_\Phi,\beta}}\delta^*V(\beta) ,$$
for a complete set of representatives $\{\delta\}$ in $\Delta^+$ of the quotient group $\Delta^{\beta_1}/\Delta^\beta\cdot \ca{K}_{\beta_1}$. For a different choice of representatives, the left-hand side of the displayed equation above remains unchanged since the composition of the map from $\ca{T}(\beta_1)$ to $\overline{\ca{S}}_2$ is always equal to $[\delta]\circ g_2\circ\pi_2\circ\pi_1$. 
We continue to compute $$\disju_{\{\delta\}}  \ca{S}_{\K^b_\Phi,\beta}\times_{\overline{\ca{S}}_{\K^b_\Phi,\beta}}\delta^*V(\beta)\iso \disju_{\{\delta\}}  \ca{S}_{\K^b_\Phi,\beta}\times_{\overline{\ca{S}}_{\K^b_\Phi,\beta},[\delta]\circ g_1\circ f_1^*\pi_2}V(\beta).$$
By Lemma \ref{lem-torus-cartesian-diag}, the last expression is isomorphic to 
$$\disju_{[\delta]\in\Delta_0}  \ca{S}_{\K^b_\Phi,\beta}\times_{f_2\circ\mbf{p}_1,\overline{\ca{S}}_2,[\delta]\circ g_2\circ\pi_2}V(\beta_1)\iso \ca{S}_{\K^b_\Phi,\beta}\times_{\overline{\ca{S}}_{\K^b_\Phi,\beta}}U(\beta).$$
The second claim has been shown.
\end{proof}
Write $\ca{S}_{\K^b_\Phi}$ (resp. $\overline{\ca{S}}_{\K^b_\Phi}$) as a disjoint union 
$\disju_{\beta\in I_{ZP_\Phi/P_0,K}}\ca{S}_{\K^b_\Phi,\beta}$ (resp. $\disju_{\beta\in I_{ZP_\Phi/P_0,K}}\overline{\ca{S}}_{\K^b_\Phi,\beta}$). Let $\ca{U}:=\disju_{\beta\in I_{ZP_\Phi/P_0,K}}U(\beta)$, $\ca{V}:=\disju_{\beta_1\in I_{ZP_\Phi/P_0,K^+}}V(\beta_1)$ and $\ca{T}:=\disju_{\beta_1\in I_{ZP_\Phi/P_0,K^+}}\ca{T}(\beta_1)$.
\begin{lem}\label{lem-torsor-iso-union}
From the constructions above, we have a finite flat cover $\ca{V}$ over $\overline{\ca{S}}_{\K^{+,b}_\Phi}$ and an $\mbf{E}_{\K_\Phi}$-torsor $\ca{T}$ over $\ca{V}$, such that the pullback of $\ca{V}$ along the transition map $\overline{\ca{S}}_{\K_\Phi^b}\to\overline{\ca{S}}_{\K^{+,b}_\Phi}$ is $\ca{U}$ and such that there is an isomorphism of $\mbf{E}_{\K_\Phi}$-torsors
\begin{equation*}
\ca{S}_{\K^b_\Phi}\times_{\overline{\ca{S}}_{\K^b_\Phi}}\ca{U}\iso \overline{\ca{S}}_{\K_\Phi^b}\times_{\overline{\ca{S}}_{\K^{+,b}_\Phi}}\ca{T}.\end{equation*}
\end{lem}
\begin{proof}
Combine Lemma \ref{lem-torsor-iso} with the definitions above. Note that for distinct elements $\beta$ and $\beta^*$ in $I_{ZP_\Phi/P_0,K}$ mapping to the same $\beta_1\in I_{ZP_\Phi/P_0,K^+}$, both $U(\beta)$ and $U(\beta^*)$ are the pullbacks of the same $V(\beta_1)$ by (\ref{diag-quotient-decomposition-delta}) again since $\beta^*\equiv \delta\beta$ in $I_{ZP_\Phi/P_0,K}$ for some $\delta\in \Delta^+$ lifting an element in $K_{Z,p}^+$.
\end{proof}
\begin{proofof}[Lemma \ref{lem-maintorus-s*}]
By Lemma \ref{lem-torsor-iso-union}, we have \begin{equation}\label{eq-separate-torus-part}\ca{S}_{\K^b_\Phi}\times_{\overline{\ca{S}}_{\K^b_\Phi}}\ca{U}\iso \overline{\ca{S}}_{\K_\Phi^b}\times_{\overline{\ca{S}}_{\K^{+,b}_\Phi}}\ca{T},\end{equation} and this isomorphism is equivariant under $K_{Z,p}^+$-action by (\ref{diag-quotient-decomposition-delta}), where $K_{Z,p}^+$ acts equivariantly on $\ca{S}_{\K_\Phi^b}\to \overline{\ca{S}}_{\K_\Phi^b}$ and trivially on $\overline{\ca{S}}_{\K^{+,b}_\Phi}$, and this induces $K_{Z,p}^+$-action on both sides of (\ref{eq-separate-torus-part}).\par
Let $\ca{U}^*:=\overline{\ca{S}}_{\K^b_\Phi,\ca{O}_{K_Z}}\times_{\overline{\ca{S}}_{\K^b_\Phi}}\ca{U}$. Pulling back (\ref{eq-separate-torus-part}) along $\overline{\ca{S}}_{\K^b_\Phi,\ca{O}_{K_Z}}\to\overline{\ca{S}}_{\K^b_\Phi}$, we have an isomorphism 
$$\ca{S}_{\K^b_\Phi,\ca{O}_{K_Z}}\times_{\overline{\ca{S}}_{\K^b_\Phi,\ca{O}_{K_Z}}}\ca{U}^*\iso \overline{\ca{S}}_{\K^b_\Phi,\ca{O}_{K_Z}}\times_{\overline{\ca{S}}_{\K^{+,b}_\Phi}}\ca{T}.$$
Let $E^{\prime,ur}_{K_Z}:=E'_{K_Z}\cap E^p$ and $\ca{O}_{K_Z}^{ur}:=\ca{O}_{E^{\prime,ur}_{K_Z}}\otimes_{\ca{O}_{E'}}\ca{O}_{E',(v)}$. The $\gal(E'_{K_Z}/E')$-action on the left-hand side is given by the descent datum of $\sh_{\K^b_\Phi,E'_{K_Z}}$ from $E'_{K_Z}$ to $E'$, and it differs from the descent datum of $\sh_{\K^a_\Phi,E'_{K_Z}}$ from $E'_{K_Z}$ to $E'$ by $r_{E',K_{Z,p}}$, where $r_{E',K_{Z,p}}$ assigns to any $\sigma\in \gal(E'_{K_Z}/E')$ the multiplication of an element in $Z_G(\bb{Q})\bss Z_G(\A)/K_{Z,p}$ on $\ca{S}_{\K^b_\Phi,\ca{O}_{K_Z}}$. 
Since $\ca{S}_{\K^b_\Phi,\ca{O}_{K_Z}}\iso \ca{S}_{\K^a_\Phi,\ca{O}_{K_Z}}$, $\overline{\ca{S}}_{\K^b_\Phi,\ca{O}_{K_Z}}\iso\overline{\ca{S}}_{\K^a_\Phi,\ca{O}_{K_Z}}$ and $\overline{\ca{S}}_{\K^{+,b}_\Phi,\ca{O}_{K_Z}^{ur}}\iso\overline{\ca{S}}_{\K^{+,a}_\Phi,\ca{O}_{K_Z}^{ur}}$,
by taking the $\gal(E'_{K_Z}/E^{\prime,ur}_{K_Z})$-quotient of the right-hand side of  
\begin{equation}\label{eq-separate-torus-part-a}
\ca{S}_{\K^a_\Phi,\ca{O}_{K_Z}}\times_{\overline{\ca{S}}_{\K^a_\Phi,\ca{O}_{K_Z}}}\ca{U}^*\iso \overline{\ca{S}}_{\K^a_\Phi,\ca{O}_{K_Z}}\times_{\overline{\ca{S}}_{\K^{+,a}_\Phi,\ca{O}_{K_Z}^{ur}}}\ca{T}_{\ca{O}_{K_Z}^{ur}}
\end{equation}
we see $\ca{S}_{\K^a_\Phi,\ca{O}^{ur}_{K_Z}}\to\overline{\ca{S}}_{\K^a_\Phi,\ca{O}^{ur}_{K_Z}}$ is an fppf $\mbf{E}_{\K_\Phi}$-torsor. The desired result now follows since it has been shown that $\ca{S}_{\K^a_\Phi}$ is an fppf $\mbf{E}_{\K_\Phi}$-torsor and $\mbf{E}_{\K_\Phi}$ is smooth.
\end{proofof}

\begin{lem}\label{lem-torus}
The statement \ref{maintheorem-3} is true.
\end{lem}
\begin{proof}
It remains to show that $\ca{S}_{K_{\Phi_2}}\to \overline{\ca{S}}_{K_{\Phi_2}}$ is an $\mbf{E}_{K_{\Phi_2}}$-torsor. 
For a fixed $K_{\Phi_2}=K_{\Phi_2,p}K^p_{\Phi_2}$, there is a neat open compact subgroup $K^{p,\uparrow}_{\Phi}$ of $ZP_{\Phi}(\Ap)$ such that $\sh_{K_{\Phi_2}}\to \sh_{\K_{\Phi,p}K_{\Phi}^{p,\uparrow}}(ZP^a_{\Phi},ZP^a_{\Phi}(\bb{Q})D_{\Phi})$ is an open and closed embedding by Lemma \ref{lem-emb-mix} and Proposition \ref{ext-imp}. Moreover, we can choose a sufficiently small $K^p\sbst G(\Ap)$ containing $g_{\Phi_2}^{-1}K_{\Phi}^{p,\uparrow}g_{\Phi_2}$ such that $g_{\Phi_2}K^p g_{\Phi_2}^{-1}\cap ZP_{\Phi}(\Ap)=K_{\Phi}^{p,\uparrow}$. We can then deduce the desired result from the result for $\ca{S}_{\K_\Phi^a}$: 
In fact, under this choice of $K^p$, $\sh_{\K^a_{\Phi}}\to \overline{\sh}_{\K^a_\Phi}$ is an $\mbf{E}_{K_{\Phi_2}}$-torsor since $\sh_{K_{\Phi_2}}\to \sh_{\K_{\Phi,p}K_{\Phi}^{p,\uparrow}}(ZP^a_{\Phi},ZP^a_{\Phi}(\bb{Q})D_{\Phi})$ is a closed embedding. 
Then $\ca{S}_{K_{\Phi_2}}\iso \overline{\ca{S}}_{K_{\Phi_2}}\times_{\overline{\ca{S}}_{\K^a_\Phi}}\ca{S}_{\K_\Phi^a}$ by Zariski's main theorem.
\end{proof}
\subsubsection{End of the proof.}
\begin{lem}\label{lem-completion-comparison}
The statements \ref{maintheorem-1}, \ref{maintheorem-2} and \ref{maintheorem-4} are true.
\end{lem}
\begin{proof}
From Proposition \ref{prop-normalized-base-change}, exactly the same proof as Proposition \ref{prop-zpb-completion-iso} and Construction \ref{const-strata-compactification-2}, we have that 
$$\ca{Z}_{[ZP^a({\Phi},\sigma)],K,\ca{O}_{K_Z}}\iso\ca{S}_{\wdtd{K}^a_\Phi,\sigma,\ca{O}_{K_Z}}$$ and that
\begin{equation}\label{eq-kaphi-completion-iso}
\cpl{\ca{S}_{\K_{\Phi}^a}(\sigma)_{\ca{O}_{K_Z}}}{\ca{S}_{\K_{\Phi}^a,\sigma,\ca{O}_{K_Z}}}\iso \cpl{\ca{S}^{\Sigma}_{K,\ca{O}_{K_Z}}}{\ca{Z}_{[ZP^a(\Phi,\sigma)],K,\ca{O}_{K_Z}}}.\end{equation}
By Zariski's main theorem, Lemma \ref{zp-iso} and Corollary \ref{zp-ab}, we have $\Delta^\circ_{\Phi,K}\bss \ca{S}_{K_{\Phi}^a,\sigma,\ca{O}_{K_Z}}\iso \ca{Z}_{[(\Phi,\sigma)],K,\ca{O}_{K_Z}}$ and $\disju_{[(\Phi',\sigma')]\in [ZP^a(\Phi',\sigma')]}\ca{Z}_{[(\Phi',\sigma')],K,\ca{O}_{K_Z}}\iso \ca{Z}_{[ZP^a(\Phi,\sigma)],K,\ca{O}_{K_Z}}$.
Taking the completion with respect to $\ca{Z}_{[(\Phi,\sigma)],K,\ca{O}_{K_Z}}$ on both sides of (\ref{eq-kaphi-completion-iso}), as we can check that $$\cpl{\ca{S}_{\K_{\Phi}^a}(\sigma)_{\ca{O}_{K_Z}}}{\ca{S}_{\K_{\Phi}^a,\sigma,\ca{O}_{K_Z}}}\iso \disju_{[(\Phi',\sigma')]\in [ZP^a(\Phi',\sigma')]} \Delta^\circ_{\Phi',K}\bss \cpl{\ca{S}_{K_{\Phi'}}(\sigma')_{\ca{O}_{K_Z}}}{\ca{S}_{K_{\Phi'},\sigma',\ca{O}_{K_Z}}},$$
we have $$\cpl{\ca{S}^\Sigma_{K,\ca{O}_{K_Z}}(G,X_a)}{\ca{Z}_{[(\Phi,\sigma)],K,\ca{O}_{K_Z}}}\iso \Delta^\circ_{\Phi,K}\bss \cpl{\ca{S}_{K_{\Phi}^a}(\sigma)_{\ca{O}_{K_Z}}}{\ca{S}_{K_{\Phi},\sigma,\ca{O}_{K_Z}}}.$$

Since $\sh_{K_2}^{\Sigma_2'}(G_2,X_2)\lra \sh_{K}^\Sigma(G,X_a)$ is an open and closed embedding by Corollary \ref{ext-cpt-imp}, we see that  
$$\Delta^\circ_{{\Phi_2},K_2}\bss\cpl{\ca{S}_{K_{\Phi_2}}(\sigma_2)_{\ca{O}_{K_Z}}}{\ca{S}_{K_{\Phi_2},\sigma_2,\ca{O}_{K_Z}}}\iso \cpl{\ca{S}^{\Sigma_2'}_{K_2,\ca{O}_{K_Z}}}{\ca{Z}_{[(\Phi_2,\sigma_2)],K_2,\ca{O}_{K_Z}}}.$$
Applying \ref{maintheorem-3}, the action of $\gal(E'_{K_Z}/E')$ on $\ca{S}_{K_{\Phi_2},\ca{O}_{K_Z}}$ naturally induces actions of $\gal(E'_{K_Z}/E')$ on $\ca{S}_{K_{\Phi_2},\sigma_2,\ca{O}_{K_Z}}$ and $\ca{S}_{K_{\Phi_2}}(\sigma_2)_{\ca{O}_{K_Z}}$.\par
The actions of $\Delta^\circ_{\Phi_2,K_2}$ and $\gal(E'_{K_Z}/E')$ commute with each other since it is so over the generic fiber. Finally, by \cite[Lem. A.3.4]{Mad19}, the action of $\gal(E'_{K_Z}/E')$ on $\ca{S}_{K_2,\ca{O}_{K_Z}}$ extends to $\ca{S}^{\Sigma_2'}_{K_2,\ca{O}_{K_Z}}$.\par 
Define $\ca{S}^{\Sigma_2'}_{K_2}$ as the quotient of $\ca{S}^{\Sigma_2'}_{K_2,\ca{O}_{K_Z}}$ by $\gal(E'_{K_Z}/E')$ (note that the quotient exists and is proper (resp. projective) over $\ca{O}'_{(v)}$ by \cite[Thm. 3.1.13]{CLO12} (resp. \cite[Ch. 4, Prop. 1.5]{Knu71}) again). Then $\ca{S}^{\Sigma'_2}_{K_2}$ satisfies \ref{maintheorem-1}.\par 
Finally, taking the quotient by $\gal(E'_{K_Z}/E')$ of the isomorphism above, let us show \ref{maintheorem-4} holds and show \ref{maintheorem-2} along the way. Taking Galois quotient induces a quasi-finite morphism $f_{\Upsilon_2}:\ca{Z}_{[(\Phi_2,\sigma_2)],K_2}\to \ca{S}_{K_2}^{\Sigma_2'}$ which is finite onto its image (which is inductively defined as remarked in Lemma \ref{lem-strata-zpb}). Denote by $[\ca{Z}_{\Upsilon_2}]$ the image of $f_{\Upsilon_2}$. There is a structural morphism induced by canonical projection to the closed stratum of the torus embedding:
$$\cpl{\ca{S}^{\Sigma_2}_{K_2,\ca{O}_{K_Z}}}{\ca{Z}_{\Upsilon_2,K_2,\ca{O}_{K_Z}}}\iso\cpl{\ca{S}_{K_{\Phi_2},\ca{O}_{K_Z}}^*(\sigma_2)}{\ca{S}^*_{K_{\Phi_2},\sigma_2,\ca{O}_{K_Z}}}\twoheadrightarrow\ca{Z}_{\Upsilon_2,K_2,\ca{O}_{K_Z}},$$
and the pre-composition of it with $\ca{Z}_{\Upsilon_2,K_2,\ca{O}_{K_Z}}\to \cpl{\ca{S}^{\Sigma_2'}_{K_2,\ca{O}_{F_Z}}}{\ca{Z}_{\Upsilon_2,K_2,\ca{O}_{F_Z}}}$ is identity.
Then taking Galois quotient induces an isomorphism $f_{\Upsilon_2}$. In conclusion, we obtained a stratification of normal and flat subschemes for $\ca{S}_{K_2}^{\Sigma_2'}$ and also the assertion \ref{maintheorem-4} over $\ca{O}'_{(v)}$. 
\end{proof}

\begin{lem}\label{cor-hs-def-ring}
The statement \ref{maintheorem-5} is true.    
\end{lem}
\begin{proof}
It follows from Proposition \ref{prop-mixsh-canonicity} (3) and Construction \ref{const-mixed-a} that the towers $\ca{S}_{\K^a_{\Phi,p}}$, $\overline{\ca{S}}_{\K^a_{\Phi,p}}$ and $\ca{S}_{\K^a_{\Phi,h,p}}$ over $\ca{O}'_{(v)}$ have the extension property. By Lemma \ref{lem-emb-mix}, we see that $\ca{S}_{K_{\Phi_2,p}}$ (resp. $\overline{\ca{S}}_{K_{\Phi_2,p}}$ and resp. $\ca{S}_{K_{\Phi_2,h,p}}$) is open and closed in $\ca{S}_{\K^a_{\Phi,p}}$ (resp. $\overline{\ca{S}}_{\K^a_{\Phi,p}}$ and resp. $\ca{S}_{\K^a_{\Phi,h,p}}$) when $\K^{a,p}_\Phi$ is sufficiently small. Once we have the statement over $\ca{O}'_{(v)}$, the descent data from $\ca{O}'_{(v)}$ to $\ca{O}_2$ will extend to the integral models. Thus, we have the statement over $\ca{O}_2$.
\end{proof}
In Case (HS) and its deeper levels, by the extension property \ref{maintheorem-5} and \cite[Lem. A.3.4]{Mad19} and since we have shown all statements over $\ca{O}'_{(v)}$, we then all results descend to $\ca{O}_2=\ca{O}_{E_2,(v_2)}$.
In Case ($\text{STB}_n$) and its deeper levels, by construction at the beginning of \S \ref{subsec-con-stra-comp-quo}, for any place $v_2|p$ of $E_2$, we can find a place $v'$ of $E'$ over $v_2$ and find a place $v$ of $E_0$ such that $v'|v|p$, and we can establish the results over $\ca{O}_2=\ca{O}_{E_2,v_2}$. \par
Now we have finished the proof of Theorem \ref{maintheorem}.\hfill$\square$
\subsubsection{}
Let us collect some consequences of Theorem \ref{maintheorem}.\par
\begin{cor}\label{cor-maintheorem-immediate}
With Theorem \ref{maintheorem} and all constructions in this {section}, the following statements are true: 
\begin{enumerate}
    \item (See \cite[Cor. 4.1.7]{Mad19}.) When $G_2$ is $\bb{Q}$-anisotropic (i.e. $G_2$ does not have proper parabolic subgroups over $\bb{Q}$, or equivalently, $G_2$ does not contain noncentral split torus), $\ca{S}_{K_2}$ is projective over $\ca{O}_2$. The converse is also true.
    \item (See \cite[Prop. 14.1 and Prop. 14.2]{Lan16b}.) Suppose that $\Sigma_2'$ is smooth. Let $\ca{P}$ be one of the following properties: (geometrically) reduced, (geometrically) normal, Cohen-Macaulay, (geometrically) $(R_i)$ each one for $i\geq 0$, $(S_i)$ each one for $i\geq 0$, (geometrically) regular. Then a fiber of $\ca{S}^{\Sigma_2'}_{K_2}\to\spec\ca{O}_2$ satisfies the property $\ca{P}$ if and only if the corresponding fiber of $\ca{S}_{K_2}\to\spec\ca{O}_2$ has the property $\ca{P}$. Moreover, $\ca{S}_{K_2}^{\Sigma_2'}$ is regular (resp. smooth over $\spec \ca{O}_2$) if and only if $\ca{S}_{K_2}$ is regular (resp. smooth over $\spec\ca{O}_2$).
    \item (See \cite[Cor. 10.18 and Cor. 14.4]{Lan16b}.) The scheme $\ca{S}_{K_2}$ is fiberwise dense in $\ca{S}_{K_2}^{\Sigma_2'}$. If $\ca{S}_{K_2}\to\spec\ca{O}_2$ has geometrically reduced fibers, so does $\ca{S}_{K_2}^{\Sigma_2'}\to\spec\ca{O}_2$; in this case, the geometric fibers of the morphisms $\ca{S}_{K_2}\to \spec\ca{O}_2$ and $\ca{S}_{K_2}^{\Sigma_2'}\to\spec\ca{O}_2$ have the same number of connected components.
 
\end{enumerate}
\end{cor}
\begin{proof}
The proof is almost identical to that of the references cited. If $G_2$ is $\bb{Q}$-anisotropic, there are no proper $\bb{Q}$-parabolic subgroups, and therefore there are no extra strata in toroidal compactifications by \ref{maintheorem-2}. Conversely, if $\ca{S}_{K_2}$ is proper, it is identical to any of its toroidal compactifications by Zariski's main theorem. Hence, it does not have extra strata in its compactifications and $G_2^\ad$ does not contain admissible $\bb{Q}$-parabolic subgroups. Since a proper admissible $\bb{Q}$-parabolic subgroup is a product of maximal proper $\bb{Q}$-parabolic subgroups in some $\bb{Q}$-simple factors of $G^\ad_2$ with the whole groups of other $\bb{Q}$-simple factors of $G^\ad_2$, there are no proper $\bb{Q}$-parabolic subgroups contained in $G^\ad_2$ if there are no maximal ones. So part one is true.\par
Part two follows from exactly the same arguments as in \cite[Prop. 14.1 and 14.2]{Lan16b}. Note that, with \ref{maintheorem-3} and \ref{maintheorem-4} in hand, one can replace $\Vec{\Xi}_{\Phi_{\ca{H}},\delta_{\ca{H}}}\to \Vec{C}_{\Phi_{\ca{H}},\delta_{\ca{H}}}$ and $\Vec{\Xi}_{\Phi_{\ca{H}},\delta_{\ca{H}}}(\sigma)\to \Vec{C}_{\Phi_{\ca{H}},\delta_{\ca{H}}}$ in the proof of \emph{loc. cit.} with $\ca{S}_{K_{\Phi_2}}^*\to \overline{\ca{S}}_{K_{\Phi_2}}^*$ and $\ca{S}_{K_{\Phi_2}}^*(\sigma_2)\to \overline{\ca{S}}_{K_{\Phi_2}}^*$.\par
It follows from \ref{maintheorem-3} and \ref{maintheorem-4} that $\ca{S}_{K_2}$ is fiberwise dense in $\ca{S}_{K_2}^{\Sigma_2'}$. Then the rest of the third part follows from 
 \cite[\href{https://stacks.math.columbia.edu/tag/0E0N}{Lem. 0E0N}]{stacks-project} and \ref{maintheorem-1}.\par 
 Note that for both the second and the third parts, we needed to use the fact that the isomorphisms in \ref{maintheorem-4} preserve the open dense strata corresponding to $\ca{S}_{K_2}\sbst \ca{S}_{K_2}^{\Sigma_2'}$ and $\ca{S}_{K_{\Phi_2}}^*\sbst \ca{S}_{K_{\Phi_2}}^*(\sigma_2)$, and this fact follows from our quotient construction and from the corresponding fact in the Hodge-type case.
\end{proof}
\begin{rk}
For the definition of being $\bb{Q}$-anisotropic, we follow \cite[Prop. 16.2.2]{Spr09}. This definition is equivalent to ``being $\bb{Q}$-anisotropic modulo center'' in many other references where a reductive group is called $\bb{Q}$-anisotropic if it does not contain any split $\bb{Q}$-torus.
\end{rk}
\subsubsection{}The following statements follow from a combination of the main results and the Artin approximation technique (and its generalizations).\par
The next proposition is a formal consequence of having a \emph{strata-preserving} toroidal embedding.
\begin{prop}\label{prop-artin-approximation-strata-preserving}
(See \cite[Prop. 2.2(9)]{LS18i} and \cite[Prop. 2.1.2(9)]{LS18b}) Let $x$ be a point of $\ca{Z}_{[(\Phi_2,\sigma_2)],K_2}\iso \ca{S}^*_{K_{\Phi_2},\sigma_2}$. There is an {\'e}tale neighborhood $\mathsf{e}_1:\overline{U}\to \ca{S}_{K_2}^{\Sigma_2'}$ of $x$, and an {\'e}tale morphism $\mathsf{e}_2:\overline{U}\to \ca{S}_{K_{\Phi_2}}^*(\sigma_2)$, such that the stratifications on $\overline{U}$ defined by pulling back those on $\ca{S}_{K_2}^{\Sigma_2'}$ and $\ca{S}^*_{K_{\Phi_2}}(\sigma_2)$ via $\mathsf{e}_1$ and $\mathsf{e}_2$ coincide. In particular, we have $\mathsf{e}_1^{-1}(\ca{S}_{K_2})=\mathsf{e}_2^{-1}(\ca{S}_{K_{\Phi_2}}^*)$.     
\end{prop}
\begin{proof}
Let $A$ be the local ring $\ca{O}_{\ca{S}^*_{K_{\Phi_2}}(\sigma_2),x}$ at $x$ with Henselization denoted by $A^h$. Denote by $\wat{A}$ the complete local ring of $A^h$ at $x$. Let $\wdtd{A}$ be the finitely generated $A^h$-algebra generated by the local ring $\ca{O}_{\ca{S}^{\Sigma_2'}_{K_2},x}$ via the last isomorphism in \ref{maintheorem-4}. 
Then there is a sequence of inclusions $A\xhookrightarrow{i} A^h\xhookrightarrow{\alpha} \wdtd{A}\xhookrightarrow{\beta} \wat{A}$. By Artin approximation \cite[Thm. 1.10]{Art69} (cf. \cite[Prop. 6.3.2.2]{Lan13}), for any integer $k>0$, there is a section $\alpha_k: \wdtd{A}\to A^h$, such that $\alpha_k\circ\alpha=\mrm{Id}$ and $\beta\circ\alpha\circ\alpha_k-\beta$ is trivial modulo $\mathfrak{m}^k_x$, where $\mathfrak{m}_x$ is the maximal ideal in $\wat{A}$ corresponding to $x$. For $k\geq 2$, $\alpha_k$ defines a scheme $\overline{U}$ with {\'e}tale morphisms $\mathsf{e}_1$ and $\mathsf{e}_2$ to $\ca{S}_{K_2}^{\Sigma_2'}$ and $\ca{S}^*_{K_{\Phi_2}}(\sigma_2)$ by \cite[Cor. 2.6]{Art69}.\par
Now we modify the construction above to get strata-preserving $(\overline{U},\mathsf{e}_1,\mathsf{e}_2)$. Note that using induction it is enough to show it for closures of strata. Let $\overline{\ca{M}}$ be the closure of a stratum $\ca{M}$ in $\ca{S}^{\Sigma_2'}_{K_2}$ which under the last isomorphism in \ref{maintheorem-4}, corresponds to $\overline{\ca{N}}$, the closure in $\ca{S}_{K_{\Phi_2}}(\sigma_2)$ of a stratum $\ca{N}$. Suppose that the pullback of $\overline{\ca{N}}$ to $\spec A$ determines an ideal $I\sbst A$ and the pullback of $\overline{\ca{M}}$ to $\spec \wdtd{A}$ determines an ideal $J\sbst \wdtd{A}$. We have $\beta\circ\alpha\circ i(I)\wat{A}=\beta(J)\wat{A}$. 
Since $\wat{A}$ is Noetherian, choose a finite set of elements $\{x_t^{\ca{M}}\}_{t\in \Lambda_\ca{M}}\sbst \wat{A}$ such that $\beta(J)\sbst \sum_{t\in \Lambda_\ca{M}} I\wdtd{A} x_t^{\ca{M}}$, and a finite set of $\{y_s^{\ca{M}}\}_{s\in \ca{I}_\ca{M}}\sbst \wat{A}$ such that $\beta\circ\alpha\circ i(I)\sbst \sum_{s\in\ca{I}_\ca{M}}J \wdtd{A}y_s^{\ca{M}}$. Let $\wdtd{B}$ be the $A^h$-algebra generated by $\wdtd{A}$ and all $\{x_t^{\ca{M}}\}_{t\in \Lambda_\ca{M}}$ and $\{y_s^{\ca{M}}\}_{s\in \ca{I}_\ca{M}}$ running over all strata $\ca{M}$ of $\ca{S}^{\Sigma_2'}_{K_2}$. The algebra $\wdtd{B}$ is still finitely generated over $A^h$, and there is a sequence of inclusions $A\xhookrightarrow{i} A^h\xhookrightarrow{\wdtd{\alpha}} \wdtd{B}\xhookrightarrow{\wdtd{\beta}} \wat{A}$. By Artin approximation as above, for any integer $k\geq 2$, there is a section $\wdtd{\alpha}_k: \wdtd{B}\to A^h$, such that $\wdtd{\alpha}_k\circ\wdtd{\alpha}=\mrm{Id}$, $\wdtd{\beta}\circ\wdtd{\alpha}\circ\wdtd{\alpha}_k-\wdtd{\beta}$ is trivial modulo $\mathfrak{m}^k_x$, and $\alpha_k:= \wdtd{\alpha}_k|_{\wdtd{A}}$ defines a scheme $\overline{U}$ with {\'e}tale morphisms $\mathsf{e}_1$ and $\mathsf{e}_2$ to $\ca{S}_{K_2}^{\Sigma_2'}$ and $\ca{S}^*_{K_{\Phi_2}}(\sigma_2)$. 
Since for any $j\in J$, $ (\beta\circ\alpha\circ\alpha_k-\beta)(j)\sbst (\wdtd{\beta}\circ\wdtd{\alpha}\circ\wdtd{\alpha}_k-\wdtd{\beta})(\sum I \wdtd{A}x_t^\ca{M})=\sum I (\wdtd{\beta}\circ\wdtd{\alpha}\circ\wdtd{\alpha}_k-\wdtd{\beta})(\wdtd{A}x_t^{\ca{M}})\sbst I\wat{A}$ and $\alpha\circ\alpha_k\circ\alpha\circ i(I)=\alpha\circ i(I)$, we have $I\wat{A}=\beta(J)\wat{A}=\beta\circ\alpha\circ\alpha_k(J)\wat{A}$, as desired. 
\end{proof}
From Proposition \ref{prop-artin-approximation-strata-preserving} above, we see that
\begin{cor}\label{cor-artin-approximation-strata-preserving}
 The complement of $\ca{S}_{K_2}$ in $\ca{S}_{K_2}^{\Sigma_2'}$ is a relative effective Cartier divisor.   
\end{cor}
\begin{cor}\label{cor-maintheorem}
With the constructions and conventions above,
after taking the normalized base change to $\ca{O}_{K_Z}$, $\ca{S}_{K_2,\ca{O}_{K_Z}}^{\Sigma_2'}$ is an open and closed subspace of $\ca{S}_{K,\ca{O}_{K_Z}}^{\Sigma}(G,X_b)$, and $\ca{S}_{K,\ca{O}_{K_Z}}^{\Sigma}(G,X_b)$, equipped with a log structure defined by the complement of $\ca{S}_{K,\ca{O}_{K_Z}}(G,X_b)$ in $\ca{S}_{K,\ca{O}_{K_Z}}^\Sigma(G,X_b)$, admits a finite Kummer {\'e}tale cover of disjoint union of Hodge-type toroidal compactifications $\ca{S}^{\mrm{TOR}}_{K_0,\ca{O}_{K_Z}}:=\disju_{\alpha\in I_{G/G_0}}\ca{S}^{\Sigma_0^\alpha}_{K_0^\alpha,\ca{O}_{K_Z}}$. The induced map between interiors $\disju_{\alpha\in I_{G/G_0}}\ca{S}_{K_0^\alpha,\ca{O}_{K_{Z}}}\to\ca{S}_{K,\ca{O}_{K_Z}}(G,X_b)$ is finite {\'e}tale.
\end{cor}
\begin{proof}
It remains to show that the morphism from $\ca{S}_{K_0,\ca{O}_{K_Z}}^{\mrm{TOR}}$ to $\ca{S}^\Sigma_{K,\ca{O}_{K_Z}}$ is Kummer {\'e}tale.
First, we apply the nested approximation theorem due to Teissier and Spivakovsky, see \cite[Cor. 2.10]{Tei95} and \cite[Thm. 11.5]{Spi99}, to find compatible approximations of a morphism of completions. Let us explain the exact meanings of the symbols in \emph{loc. cit.} specialized to our situation. 
The $m$ in \emph{loc. cit.} is $2$ in our situation. Pick any point $x\in \ca{Z}:=\ca{Z}_{[ZP^b(\Phi,\sigma)],K,\ca{O}_{K_Z}}$ and let $A_1$ in \emph{loc. cit.} be the Henselization of the local ring of $\ca{S}^{\Sigma}_{K,\ca{O}_{K_Z}}$ at $x$. The pullback of $\ca{S}_{K,\ca{O}_{K_Z}}$ to $\spec A_1$ determines an ideal $I\sbst A_1$. Let $B_1$ be the completion of $A_1$ with respect to the maximal ideal $\mathfrak{m}_x\sbst A_1$ corresponding to $x$.
The point $x$ corresponds to a point $x'$ in $\ca{S}_{\K_\Phi^b,\sigma,\ca{O}_{K_Z}}$, and let $\ca{O}$ be the local ring of $\ca{S}_{\K^b_\Phi,\ca{O}_{K_Z}}(\sigma)$ at $x'$. 
The pullback of $\ca{S}_{\K^b_\Phi,\ca{O}_{K_Z}}$ to $\spec \ca{O}$ determines an ideal $J\sbst \ca{O}$. We have $J B_1= I B_1$. Choose a finite set of elements $\{x_t\}_{t\in \Lambda}$ in $B_1$ such that $J\sbst \sum_{t\in \Lambda}I\ca{O} \cdot x_t$ and a finite set of elements $\{y_s\}_{s\in \ca{I}}$ in $B_1$ such that $I\sbst \sum_{s\in \ca{I}}J\ca{O} \cdot y_s$.
Let $C_1$ be the finitely generated $A_1$-algebra generated by $\{x_t\}$ and $\{y_t\}$ and $\ca{O}$ via (\ref{eq-zp-strata-completion-int}). Let $A_2$ (resp. $B_2$) be the algebra obtained by tensoring $\otimes_{\ca{O}_{\ca{S}^\Sigma_{K,\ca{O}_{K_Z}}}}\ca{O}_{\ca{S}^{\mrm{TOR}}_{K_0,\ca{O}_{K_Z}}}$ to $A_1$ (resp. $B_1$). 
Let $C_2$ be the finitely generated $A_2$-algebra generated by $\{x_t\}$, $\{y_s\}$ and a semi-local ring of the disjoint union of $\ca{S}_{K_{\Phi_0^\alpha},\ca{O}_{K_Z}}(\sigma_0^\alpha)$, corresponding to all cusp labels with cones mapping to $[ZP^b(\Phi,\sigma)]$ (cf. (\ref{eq-sur-str-tor})), at the inverse image of $x$ in $\ca{S}^{\mrm{TOR}}_{K_0,\ca{O}_{K_Z}}$. Since completion and Henselization commute with finite ring maps, we can apply \emph{loc. cit.} to find compatible approximations. 
More precisely (cf. \cite[Cor. 3.8]{Lan19}), there is an {\'e}tale neighborhood $\overline{U}$ of $x$ and an {\'e}tale morphism $\overline{U}\to \ca{S}_{\K^b_\Phi,\ca{O}_{K_Z}}(\sigma)$ such that the preimage $U$ in $\overline{U}$ of $\ca{S}_{K,\ca{O}_{K_Z}}$ coincides with the preimage of $\ca{S}_{\K^b_\Phi,\ca{O}_{K_Z}}$. Moreover, the {\'e}tale morphisms $\overline{U}\to \ca{S}^{\Sigma}_{K,\ca{O}_{K_Z}}$ and $\overline{U}\to \ca{S}_{\K^b_\Phi,\ca{O}_{K_Z}}(\sigma)$ induce {\'e}tale morphisms $\overline{U}_0\to \ca{S}^{\mrm{TOR}}_{K_0,\ca{O}_{K_Z}}$ and $\overline{U}_0\to \disju \ca{S}_{K_{\Phi_0^\alpha},\ca{O}_{K_Z}}(\sigma_0^\alpha)$ by pulling back via $\ca{S}^{\mrm{TOR}}_{K_0,\ca{O}_{K_Z}}\to \ca{S}^\Sigma_{K,\ca{O}_{K_Z}}$ and $\disju \ca{S}_{K_{\Phi_0^\alpha},\ca{O}_{K_Z}}(\sigma_0^\alpha)\to\ca{S}_{\K^b_\Phi,\ca{O}_{K_Z}}(\sigma)$, respectively. The two {\'e}tale morphisms obtained by pullback as before have the property that the preimages of $\disju_{I_{G/G_0}}\ca{S}_{K_0^\alpha,\ca{O}_{K_Z}}$ and $\disju \ca{S}_{K_{\Phi_0^\alpha},\ca{O}_{K_Z}}$ coincide with the preimage $U_0$ of $U\sbst \overline{U}$ via $\overline{U}_0\to \overline{U}$.\par
Now, thanks to Proposition \ref{prop-kisin-pappas-free-action}, it suffices to show that each $\ca{S}_{K_{\Phi_0^\alpha}}(\sigma_0^\alpha)\to \ca{S}_{\K^b_\Phi}(\sigma)$ is Kummer {\'e}tale, where the log structures on the source and the target are defined by the torus embeddings. By Lemma \ref{lem-twt-free-action} and Proposition \ref{prop-mixsh-canonicity}, the morphism $\mbf{E}_{K_{\Phi_0^\alpha}}\to\mbf{E}_{\K_{\Phi}}$ between tori induced by $\ca{S}_{K_{\Phi_0^\alpha}}\to \ca{S}_{\K^b_\Phi}$ is a prime-to-$p$ isogeny. Since $\sigma_0^\alpha$ is induced by $\sigma$, the map $\ca{S}_{K_{\Phi_0^\alpha}}(\sigma_0^\alpha)\to \ca{S}_{\K^b_\Phi}(\sigma)$ is log {\'e}tale and Kummer, as desired.
\end{proof}
\subsubsection{}\label{subsubsec-abelian-scheme}The main result in \S\ref{subsubsec-abelian-scheme} concerns when the proper morphism $\overline{\ca{S}}_{K_{\Phi_2}}^*\to\ca{S}_{K_{\Phi_2},h}^*$ is an abelian scheme torsor.\par
Consider (DL) of ($\mrm{STB}_n$). Suppose $K_{2,p}\sbst K_{2,p}'$, where $K_{2,p}'$ is an intersection of $n$ Bruhat-Tits stabilizer subgroups $K_{2,p}^i\sbst G_2(\bb{Q}_p)$ for $1\leq i\leq n$.
Denote by $K_{2,p}^{i,\circ}$ the parahoric subgroups in $G_2(\bb{Q}_p)$ associated with $K_{2,p}^i$. \textbf{Assume} in \S\ref{subsubsec-abelian-scheme} that \begin{equation}\label{eq-assumption-n-quasi-para} K_{2,p}^\circ:=\cap_{i=1}^nK_{2,p}^{i,\circ}\sbst K_{2,p}\sbst K_{2,p}'.\end{equation}
\begin{lem}[{See \cite[Lem. 3.24 and Cor. 3.56]{Mao25}}]\label{lem-abelian-torsor-group}
Assume that $(G_2,X_2)=(G_0,X_0)$. Let $K^p_0$ be any neat open compact subgroup of $G_0(\Ap)$, $g$ be any element in $G_0(\A)$ and $\Phi_0\in\ca{CLR}(G_0,X_0)$. Under the assumption above, we have 
$$\overline{g K_{0,p}^\circ K_{0}^p g^{-1}\cap W_{\Phi_0}(\A)}=\overline{g K_{0,p}'K_0^pg^{-1}\cap W_{\Phi_0}(\A)},$$ where the ``$\overline{\ }$'' denotes the projection to $V_{\Phi_0}(\A):=W_{\Phi_0}(\A)/U_{\Phi_0}(\A)$.
\end{lem}
\begin{proof}
In \cite[Lem. 3.24]{Mao25}, it is shown when $n=1$ that $g K_{0,p}^\circ K_{0}^p g^{-1}\cap W_{\Phi_0}(\A)=g K_{0,p}'K_0^pg^{-1}\cap W_{\Phi_0}(\A)$. The rest of the lemma follows from taking the intersection and projection.
\end{proof}
We need the following technical lemma:
\begin{lem}\label{lem-parahoric-auxiliary}Fixing $K_{2,p}$, 
the construction of $\ca{S}_{K_{\Phi_2}}\to \overline{\ca{S}}_{K_{\Phi_2}}\to\ca{S}_{K_{\Phi_2},h}$ in Theorem \ref{maintheorem} does not depend on the choice of $K_p\sbst G(\bb{Q}_p)$ made at the beginning of \S\ref{sss-main-setup-again}. \par
Under Assumption (\ref{eq-assumption-n-quasi-para}), 
\begin{itemize}
\item we can choose an open compact subgroup $K_p$ such that, for any $\alpha$, $\pi^{b,-1}(\lcj{K_p}{\alpha})$ is an intersection of quasi-parahoric subgroups corresponding to the same points as $\lcj{K^i_{2,p}}{\alpha}$'s in the Bruhat-Tits building of $G^\ad_{\bb{Q}_p}$; 
\item with this $K_p$, we can choose $\ca{O}_{K_Z}$ to be unramified over $\ca{O}'_{(v)}$.
\end{itemize}
\end{lem}
\begin{proof}
The first sentence follows from the construction that the integral models such as $\ca{S}_{K_{\Phi_2}}$ are normalizations from the corresponding integral models in Case ($\mrm{STB}_n$). \par
We now construct a $K_p$ to satisfy the first item. In Case ($\mrm{STB}_n$), $K_p$ is chosen to be an intersection of Bruhat-Tits stabilizers that correspond to the same points as $K^i_{2,p}$'s in the Bruhat-Tits building of $G^\ad_{\bb{Q}_p}$. Under only Assumption (\ref{eq-assumption-n-quasi-para}), it is reduced to find such a $K_p$ when $K_{2,p}$ is quasi-parahoric as one can take an intersection after this. When $K_{2,p}$ is quasi-parahoric, denote by $K_{0,p}$ the parahoric subgroup of $G_0(\bb{Q}_p)$ such that $K_{2,p}$ and $K_{0,p}$ correspond to the same point in the Bruhat-Tits building of $G^\ad_{\bb{Q}_p}$. Note that $K_p^\bullet:=\pi^b(K_{0,p})\cdot \pi^a(K_{2,p})$ is a compact group in $G(\bb{Q}_p)$ and $K_p^\bullet\cap G_2(\bb{Q}_p)=K_{2,p}$. We then find an open compact subgroup $K_p$ in the Bruhat-Tits stabilizer subgroup of $G(\bb{Q}_p)$ corresponding to $K_{2,p}$ containing $K_p^\bullet$ such that $K_p\cap G_2(\bb{Q}_p)=K_{2,p}$. The preimage $\pi^{b,-1}(K_p)$ contains $K_{0,p}$ and is contained in the Bruhat-Tits stabilizer corresponding to $K_{0,p}$. Taking conjugates of $g_\alpha$, we get the desired result.\par
We now explain the second item. Again, we assume that $K_{2,p}$ is quasi-parahoric first. Pick any point $x\in X_a\iso X_b\iso X^\ad$ and denote by $\mu_a$ (resp. $\mu_b$) its associated Hodge cocharacter regarded as a point on $X_a$ (resp. $X_b$). By the proof of \cite[Lem. 4.3.2]{KP15}, the image of $\ca{O}_{E'_{v'}}^\times$ under (\ref{eq-conn-components-action}) induced by $\mu_a$ (resp. $\mu_b$) is contained in $\pi^a(K_{2,p})$ (resp. $\pi^b(K_{0,p})$) up to a conjugation. As $\mu_a\cdot\mu_b^{-1}$ is central, the images of $\ca{O}_{E'_{v'}}^\times$ under (\ref{eq-conn-components-action}) induced by $\mu_a$ and $\mu_b$ are contained in $\pi^a(K_{2,p})$ and $\pi^b(K_{0,p})$ respectively after the same conjugation, and the image of it under (\ref{eq-conn-components-action}) induced by $\mu_a\cdot \mu_b^{-1}$ is contained in $K_p^\bullet\cap Z_G(\bb{Q}_p)\sbst K_p\cap Z_G(\bb{Q}_p)$. By taking intersection, we have shown the second item.
\end{proof}
\begin{rk}
When $p>2$, we expect that one can also show the statement above using the machinery developed in Daniels-Youcis \cite{DY25} and Kisin-Pappas-Zhou \cite{KPZ24}. See, e.g., \cite[Prop. 2.15 and Prop. 5.4]{DY25}. 
\end{rk}
\begin{prop}\label{prop-abelian-torsor}
Under Assumption (\ref{eq-assumption-n-quasi-para}), $\overline{\ca{S}}_{K_{\Phi_2}}^*\to\ca{S}_{K_{\Phi_2},h}^*$ is an abelian scheme torsor.
\end{prop}
\begin{proof}
Under the assumption, we can choose $E'_{K_Z}$ unramified over $E'$ at $p$ by Lemma \ref{lem-parahoric-auxiliary}. Hence, by Part \ref{prop-mixsh-canonicity-2} of Proposition \ref{prop-mixsh-canonicity} it is enough to show the statement for $\overline{\ca{S}}_{\K^b_\Phi}\to \ca{S}_{\K^b_\Phi,h}$. 
By \cite[2.1.10, Thm. 4.1.5(1)]{Mad19}, Lemma \ref{lem-abelian-torsor-group} and Zariski's main theorem, all $\overline{\ca{S}}_{K_{\Phi_{\beta'}}}\to\ca{S}_{K_{\Phi_{\beta'}},h}$ used to construct $\overline{\ca{S}}_{\K^b_\Phi}\to\ca{S}_{\K^b_\Phi,h}$ are abelian scheme torsors. 
We then employ the proof of Part \ref{prop-mixsh-canonicity-1} of Proposition \ref{prop-mixsh-canonicity} to obtain the desired result.
\end{proof}

\newpage
\section{Minimal compactifications of integral models of abelian type}\label{sec-min-cpt}
In this {section}, we construct the minimal compactification $\ca{S}^{\mrm{min}}_{K_2}$ of $\ca{S}_{K_2}$ along with the stratification on it. We still assume that $K_2=K_{2,p}K^{p}_2$, where $K^p_2$ is a neat open compact subgroup of $G_2(\Ap)$ and $K_{2,p}$ is an open compact subgroup of $G_2(\bb{Q}_p)$ and regard $(G_2,X_2,K_{2,p})$ as one of the Cases (HS), ($\text{STB}_n$) and (DL). To avoid overloaded notation, we assume $\Sigma_2=\Sigma_2'$.
\subsection{Minimal compactifications of integral models of Hodge type}\label{subsection-min-hodge}
Fix a Hodge-type Shimura datum $(G_0,X_0)$ with an embedding $\iota:(G_0,X_0)\hookrightarrow (G^\ddag,X^\ddag)$. Recall that $E_0=E(G_0,X_0)$, $v$ is a place of $E_0$ over $p$, and all integral models related to $(G_0,X_0)$ are defined over $\ca{O}_{E_0,(v)}$. Let us first briefly recall the construction of the minimal compactification $\ca{S}_{K_0}^\mmin$ for $\ca{S}_{K_0}$ (see \cite[Sec. 5]{Mad19}). We then deduce a functoriality result analogous to Corollary \ref{cor-gadq1-ext} and Lemma \ref{lem-gadq1-ext-n} for minimal compactifications.
\subsubsection{}\label{subsubsec-review-min-hodge}
As explained in \S\ref{subsec-siegel}, the integral model $\ca{S}_{K^{\ddag}_pK^{\ddag,p}}$ of Siegel type is associated with a moduli problem $\mbf{M}^{\mrm{iso}}_{(V_\bb{Z},\psi_\bb{Z}),K^{\ddag}_pK^{\ddag,p}}$. Denote by $(\ca{A}^{\ddag},\lambda,[\varepsilon]_{K^{\ddag,p}})$ the universal polarized abelian scheme with level structure on $\ca{S}_{K^{\ddag}_pK^{\ddag,p}}$. Denote the \emph{Hodge invertible sheaf} associated with $\ca{A}^{\ddag}$ by  $\omega_{\ca{A}^{\ddag}/\ca{S}_{K^{\ddag}}}:=(\wedge^{\mrm{top}}\ull_{\ca{A}^{\ddag}/\ca{S}_{K^{\ddag}}})^{\otimes -1}$; the pullback $\omega_{K_0}:=\iota^*\omega_{\ca{A}^\ddag/\ca{S}_{K^\ddag}}$ of it to $\ca{S}_{K_0}$ is called a Hodge invertible sheaf on $\ca{S}_{K_0}$. Then $\omega_{\ca{A}^\ddag/\ca{S}_{K^\ddag}}$ extends to an invertible sheaf $\omega_{K^\ddag}^{\mrm{tor}}:=(\wedge^{\mrm{top}}\ull_{\ca{G}^\ddag/\ca{S}_{K^\ddag}})^{\otimes -1}$, where $\ca{G}^\ddag$ is the degenerating semi-abelian scheme over $\ca{S}_{K^\ddag}^{\Sigma^\ddag}$ extending $\ca{A}^\ddag$; the pullback of it to $\ca{S}_{K_0}^{\Sigma_0}$ via $\iota:\ca{S}_{K_0}^{\Sigma_0}\to \ca{S}_{K^\ddag}^{\Sigma^\ddag}$ is an invertible sheaf extending $\omega_{K_0}$, denoted by $\omega_{K_0}^{\mrm{tor}}$. 
\begin{construction}[{\cite[Sec. 5.2]{Mad19}}]\label{con-min-hodge}\upshape
For a sufficiently large integer $k_0$, the tensor power $\omega_{K_0}^{\mrm{tor},\otimes k_0}$ is generated by its global sections. Then there is a morphism $\int_{K_0}:\ca{S}_{K_0}^{\Sigma_0}\to \bb{P}^M$ defined by the global sections of $\omega_{K_0}^{\mrm{tor},\otimes k_0}$. 
The minimal compactification of $\ca{S}_{K_0}$ is defined to be 
 $$\ca{S}_{K_0}^\mmin:= \mrm{Proj}(\bigoplus\limits_{k=0}^\infty\Gamma(\ca{S}_{K_0}^{\Sigma_0},\omega^{\mrm{tor},\otimes k}_{K_0}))\iso \mrm{Proj}(\bigoplus\limits_{k=0}^\infty\Gamma(\ca{S}_{K_0}^{\Sigma_0},\omega^{\mrm{tor},\otimes k_0k}_{K_0}))\iso \underline{\mrm{Spec}}_{\bb{P}^M}\int_{K_0,*}\ca{O}_{\ca{S}_{K_0}^{\Sigma_0}}.$$
Then there is a proper morphism $$\oint_{K_0}^{\Sigma_0}:\ca{S}_{K_0}^{\Sigma_0}\lra \ca{S}_{K_0}^\mmin$$
and a finite morphism 
$$\overline{\oint}_{K_0}:\ca{S}_{K_0}^\mmin\lra \bb{P}^M.$$
By general properties of Stein factorizations (see, e.g., \cite[\href{https://stacks.math.columbia.edu/tag/0A1B}{Thm. 0A1B}]{stacks-project}), we have that $\oint_{K_0}^{\Sigma_0}$ is proper and surjective, the fibers of $\oint_{K_0}^{\Sigma_0}$ are geometrically connected and $\oint_{K_0}^{\Sigma_0,\sharp}:\ca{O}_{\ca{S}^\mmin_{K_0}}\to \oint_{K_0,*}^{\Sigma_0}\ca{O}_{\ca{S}_{K_0}^{\Sigma_0}}$ is an isomorphism. Moreover, this construction is independent of the choice of $\Sigma_0$.\hfill$\square$
\end{construction}
\begin{thm}[{\cite[Thm. 5.2.11]{Mad19}; cf. \cite[Thm. V.2.5]{FC90} and \cite[Thm. 7.2.4.1]{Lan13}}]\label{thm-cpt-min-hodge}
With the conventions and constructions above, there is a natural stratification on $\ca{S}_{K_0}^\mmin$ consisting of locally closed normal subschemes that are flat over $\ca{O}_{E_0,(v)}$, that is, $\ca{S}_{K_0}^\mmin=\disju_{[\Phi_0]\in\cusp_{K_0}(G_0,X_0)}\ca{Z}_{[\Phi_0],K_0}$. The strata are labeled by equivalence classes in $\cusp_{K_0}(G_0,X_0)$; each $\ca{Z}_{[\Phi_0],K_0}$ is naturally isomorphic to $\Delta_{\Phi_0,K_0}\bss\ca{S}_{K_{\Phi_0},h}$. A stratum $\ca{Z}_{[\Phi_0'],K_0}$ lies in the closure of $\ca{Z}_{[\Phi_0],K_0}$ if and only if $[\Phi_0]'\preceq [\Phi_0]$. \par
Moreover, the morphism $\oint_{K_0}^{\Sigma_0}$ is compatible with stratifications in the following sense: The restriction of $\oint_{K_0}^{\Sigma_0}$ to each $\ca{Z}_{[(\Phi_0,\sigma_0)],K_0}$ is the natural projection $\ca{Z}_{[(\Phi_0,\sigma_0)],K_0}\iso \Delta^\circ_{\Phi_0,K_0}\bss\ca{S}_{K_{\Phi_0},\sigma_0}\to \ca{S}_{K_{\Phi_0},h}\to\Delta_{\Phi_0,K_0}\bss \ca{S}_{K_{\Phi_0},h}\iso \ca{Z}_{[\Phi_0],K_0}$; the inverse image of each $\ca{Z}_{[\Phi_0],K_0}$ under $\oint_{K_0}^{\Sigma_0}$ is $\disju_{\Upsilon_0\in\cusp_{K_0}(G_0,X_0)\mapsto [\Phi_0]}\ca{Z}_{\Upsilon_0,K_0}$.\par
The Hodge invertible sheaf $\omega_{K_0}$ extends to an ample invertible sheaf $\omega_{K_0}^\mmin$ on $\ca{S}_{K_0}^\mmin$ such that $\oint_{K_0}^{\Sigma_0,*}\omega_{K_0}^\mmin\iso \omega_{K_0}^\mrm{tor}$.\end{thm}
For normalized base change, we have
\begin{prop}\label{prop-normalizedbasechange-min}
Let $F$ be any finite field extension of $E_0$. Let $\ca{S}^\mmin_{K_0,\ca{O}_F}$ be the normalization of the base change $\ca{S}_{K_0}^\mmin\otimes_{\ca{O}_{E,(v)}}\ca{O}_{E_0,(v)}\otimes_{\ca{O}_{E_0}}\ca{O}_F$. Then $\ca{S}^\mmin_{K_0,\ca{O}_F}$ is constructed as 
$$\mrm{Proj}(\bigoplus\limits_{k=0}^\infty\Gamma(\ca{S}_{K_0,\ca{O}_F}^{\Sigma_0},\omega^{\mrm{tor},\otimes k}_{K_0,\ca{O}_F}))$$
where $\ca{S}_{K_0,\ca{O}_F}^{\Sigma_0}$ denotes the normalized base change of $\ca{S}_{K_0}^{\Sigma_0}$ to $\ca{O}_{E_0,(v)}\otimes_{\ca{O}_{E_0}}\ca{O}_F$ and $\omega^{\mrm{tor}}_{K_0,\ca{O}_F}$ denotes the pullback of $\omega^{\mrm{tor}}_{K_0}$ to $\ca{S}^{\Sigma_0}_{K_0,\ca{O}_F}$. Moreover, the properties in Theorem \ref{thm-cpt-min-hodge} are also true for $\ca{S}^\mmin_{K_0,\ca{O}_F}$, with all strata replaced with normalized base change to $\ca{O}_{E_0,(v)}\otimes_{\ca{O}_{E_0}}\ca{O}_F$.
\end{prop}
\begin{proof}
Repeat the construction in \cite[Sec. 5]{Mad19}. Note that one can also use the argument in Proposition \ref{prop-normalized-base-change} with \cite[Cor. 5.2.8]{Mad19} instead of \cite[Thm. 4.1.5 (4)]{Mad19} as an input.
\end{proof}
\subsubsection{}\label{subsubsec-review-min-gen}Let $(G_0,X_0)$ be any Shimura datum. The canonical sheaf $\varpi_0$ of $\sh_{K_0}(G_0,X_0)$ extends to an ample invertible sheaf $\varpi_0^\mmin$ on $\sh_{K_0}^\mmin$ and an invertible sheaf $\varpi_0^\mrm{tor}$ on $\sh_{K_0}^{\Sigma_0}$, which is locally defined as in \cite[5.26]{Pin89}. In fact, $\varpi_0^{\mrm{tor}}=\wedge^{\mrm{top}}\Omega^1_{\sh_{K_0}^{\Sigma_0}/\spec E_0}(\log D)$, the sheaf of meromorphic differential forms with top degree and logarithmic poles along the boundary divisor, and $\varpi_0^{\mrm{tor}}$ is isomorphic to the pullback of $\varpi_0^\mmin$ along the natural projection $\sh_{K_0}^{\Sigma_0}\to \sh_{K_0}^\mmin$ (see \cite[Prop. 3.4(b)]{Mum77}). 
Moreover, by a similar construction as in \ref{con-min-hodge}, there is a natural isomorphism $\sh_{K_0}^\mmin\iso \mrm{Proj}\bigoplus_{k=0}^\infty\Gamma(\sh_{K_0}^{\Sigma_0},\varpi_0^{\mrm{tor},\otimes k})$ (see \cite[8.2]{Pin89} and \cite[Thm. 5.2.11 (5)]{Mad19}; see also the proof of \cite[Thm. 5.1.1]{Lan12b}). In particular, the natural projection $\sh_{K_0}^{\Sigma_0}\to\sh_{K_0}^\mmin$ is proper, surjective and of connected geometric fibers by general properties of Stein factorization again.\par
The adjoint Shimura datum $(G^\ad,X^\ad)$ admits a decomposition $(G^\ad,X^\ad)=\prod_{i=1}^k(G^i,X^i)$, where $G_i$ are $\bb{Q}$-simple factors of $G^\ad$. Choose suitable neat open compact subgroups $K^i\sbst G^i(\A)$ such that there is a (finite {\'e}tale) morphism $\pi_\ad:\sh_{K_0}(G_0,X_0)\to \prod_{i=1}^k\sh_{K^i}(G^i,X^i)$. 
Denote by $\varpi_{i}$ the canonical sheaf of $\sh_{K^i}(G^i,X^i)$.\par
Now resume the assumptions in \S\ref{subsubsec-review-min-hodge} that $(G_0,X_0)$ is of Hodge type. \par
It is explained in \cite[Thm. 5.1.1]{Lan12b} that, for all PEL-type $(G^\ddag,X^\ddag)$, the construction of $\sh_{K^\ddag}^\mmin$ as above using the extension of the canonical sheaf with log poles on $\sh^{\Sigma^\ddag}_{K^\ddag}\bss\sh_{K^\ddag}$ coincides with the generic fiber of $\ca{S}_{K^\ddag}^\mmin$ constructed as in \cite{Lan13} using the extension of the Hodge invertible sheaf. Hence, there is a natural finite morphism $\sh_{K_0}^\mmin\to \ca{S}_{K^\ddag}^\mmin$ defined as post-composing the natural morphism $\sh_{K_0}^\mmin\to \sh_{K^\ddag}^\mmin$ with the natural embedding $\sh_{K^\ddag}^\mmin\hookrightarrow \ca{S}_{K^\ddag}^\mmin$.\par
The following property follows directly from $\iota$ and Construction \ref{con-min-hodge} itself; in other words, it can be shown without knowing the main theorem \cite[Thm. 5.2.11]{Mad19} or its proof:
\begin{prop}\label{prop-gen-classical}
With the conventions above, there are positive integers $N$ and $N_i$ for $1\leq i \leq k$ such that $\pi_{\ad}^*(\boxtimes_{i=1}^k(\varpi_{i}^{\otimes N_i}))\iso \omega_{K_0,\bb{Q}}^{\otimes N}$.\par 
Moreover, we can also construct $\ca{S}_{K_0}^{\mmin}$ as the normalization in 
$\sh_{K_0}$, $\sh_{K_0}^{\Sigma_0}$ or $\sh_{K_0}^\mmin$ of the minimal compactification $\ca{S}_{K^\ddag}^\mmin$ of $\ca{S}_{K^\ddag}$ constructed in \cite[V.2]{FC90} and \cite[7.2.3]{Lan13}, and construct $\oint_{K_0}^{\Sigma_0}$ as the morphism induced by the normalization in $\sh_{K_0}$ of the morphism $\oint_{K^\ddag}^{\Sigma^\ddag}:\ca{S}_{K^\ddag}^{\Sigma^\ddag}\to \ca{S}_{K^\ddag}^\mmin$. In particular, the generic fiber of $\ca{S}_{K_0}^\mmin$ is $\sh_{K_0}^\mmin$.
\end{prop}
\begin{proof}
The first assertion follows from the fifth paragraph of the proof of \cite[Prop. 2.4]{LS18}.\par 
Let us show the second paragraph. Indeed, the morphism $\oint_{K_0}^{\Sigma_0}$ in Construction \ref{con-min-hodge} is induced by the normalization in $\ca{S}^{\Sigma_0}_{K_0}$ of $\ca{S}_{K^\ddag}^{\mmin}$; in other words, $\ca{S}_{K_0}^\mmin$ is the normalization in $\ca{S}_{K_0}^{\Sigma_0}$ of $\ca{S}_{K^\ddag}^\mmin$: To see this, denote this normalization by $\ca{S}^\mmin_1$. The normalization induces a Stein factorization of $\oint_{K^\ddag}^{\Sigma^\ddag}\circ\iota$, $$\ca{S}_{K_0}^{\Sigma_0}\xrightarrow{\pi_2}\ca{S}_1^\mmin\xrightarrow{\pi_1} \ca{S}_{K^\ddag}^\mmin.$$ Since $\pi_2^*\pi_1^*\omega_{K^\ddag}^\mmin\iso \iota^*\oint_{K^\ddag}^{\Sigma^\ddag,*}\omega_{K^\ddag}^\mmin$, we have 
$$\mrm{Proj}\bigoplus_{k=0}^\infty\Gamma(\ca{S}_{K_0}^{\Sigma_0},\omega_{K_0}^{\mrm{tor},\otimes k})\iso \mrm{Proj}\bigoplus_{k=0}^\infty\Gamma(\ca{S}_1^\mmin,(\pi_1^*\omega_{K^\ddag}^\mmin)^{\otimes k})\iso \ca{S}_1^\mmin.$$
The first isomorphism is due to the projection formula and the second one holds because $\pi_1$ is finite and therefore $\pi_1^*\omega_{K^\ddag}^{\mmin}$ is ample.\par 
Then the normalization in $\sh_{K_0}$ of $\ca{S}_{K^\ddag}^\mmin$ has a finite birational morphism to $\ca{S}_{K_0}^\mmin$ and is isomorphic to $\ca{S}_{K_0}^\mmin$ by Zariski's main theorem. Similarly, by functoriality of relative normalizations and Zariski's main theorem, $\ca{S}_{K_0}^\mmin$ is also isomorphic to the normalization in $\sh_{K_0}^{\Sigma_0}$ or $\sh^\mmin_{K_0}$ of $\ca{S}_{K^\ddag}^\mmin$.
\end{proof}
\subsubsection{}\label{subsubsec-func-min}Now we resume with the notation and conventions in \S\ref{sss-conclusion} and \S\ref{subsec-con-stra-comp-quo}. 
Fix any $\gamma\in G^\ad_0(\bb{Q})_1$. Recall that by Lemma \ref{lem-gadq1-ext-n} (cf. Lemma \ref{lem-gadq1-ext}), there is a morphism $\gamma:\ca{S}_{K_{\Phi_0^\alpha},h}\to \ca{S}_{\lcj{K_{\Phi_0^\alpha},h}{\gamma}}=\ca{S}_{K_{\lcj{\Phi_0^\alpha}{\gamma}},h}$ extending the conjugation $\gamma:\sh_{K_{\Phi_0^\alpha},h}\to \sh_{\lcj{K_{\Phi_0^\alpha},h}{\gamma}}$ defined over the generic fiber.\par
As in \S\ref{subsubsec-stbndl}, we choose a Hodge embedding $\iota^\gamma:(G_0,X_0)\hookrightarrow(G^{\ddag,\gamma},X^{\ddag,\gamma})$. Define $\ca{S}_{\lcj{K_0}{\gamma}}^\mmin:=\ca{S}_{\lcj{K_0}{\gamma}}^\mmin(G_0,X_0)$ to be the normalization in $\sh_{\lcj{K_0}{\gamma}}(G_0,X_0)$ of $\ca{S}_{K^{\ddag,\gamma}}(G^{\ddag,\gamma},X^{\ddag,\gamma})$ (see the notation there).
\begin{lem}\label{lem-gadq1-min}
The morphisms $\gamma:\ca{S}_{K_{\Phi_0^\alpha},h}\to\ca{S}_{\lcj{K_{\Phi_0^\alpha},h}{\gamma}}$ for integral models of Hodge-type Shimura varieties induce morphisms $\gamma:\ca{Z}_{[\Phi_0^\alpha],K_0^\alpha}\to\ca{Z}_{[\lcj{\Phi_0^\alpha}{\gamma}],\lcj{K_0^\alpha}{\gamma}}$ between strata of $\ca{S}_{K_0^\alpha}^\mmin$ and $\ca{S}_{\lcj{K_0^\alpha}{\gamma}}^\mmin$. These morphisms between strata and the morphism $\gamma:\sh_{K_0}^\mmin\to \sh_{\lcj{K_0}{\gamma}}^\mmin$ induced by conjugation of $\gamma$ (see \cite[12.3]{Pin89}) together (uniquely) induce a morphism between integral models of minimal compactifications $\gamma:\ca{S}_{K_0}^\mmin\to\ca{S}_{\lcj{K_0}{\gamma}}^\mmin$. In particular, the construction of $\ca{S}_{K_0}^\mmin$ and $\ca{S}_{\lcj{K_0}{\gamma}}^\mmin$ is independent of the choice of Hodge embedding as in ($\mrm{STB}_n$).
\end{lem}
\begin{proof}
For the first sentence, we can use \cite[Thm. 5.2.11(4)]{Mad19} and check that $\sh_{K_{\Phi_0^\alpha},h}\to\Delta_{\lcj{\Phi_0^\alpha}{\gamma},\lcj{K_0}{\gamma}}\bss \sh_{\lcj{K_{\Phi_0^\alpha}}{\gamma},h}$ factors through $\Delta_{\Phi_0^\alpha,K_0}\bss \sh_{K_{\Phi_0^\alpha},h}$ and that $\gamma \Delta_{\Phi_0^\alpha,K}\gamma^{-1}=\Delta_{\lcj{\Phi_0^\alpha}{\gamma},\lcj{K}{\gamma}}$.
The second sentence follows from \cite[Lem. A.3.4]{Mad19} (note that the proof of \emph{loc. cit.} only used the condition that there are morphisms between strata that are compatible with a morphism over the generic fiber for the whole space).
The third sentence also follows from \cite[Lem. A.3.4]{Mad19} and Lemma \ref{lem-gadq1-ext-n}.
\end{proof}
\subsection{Construction of minimal compactifications}\label{sec-con-min}
Our construction for minimal compactifications is similar to the method in \S\ref{subsec-con-stra-comp-quo} and \S\ref{subsec-main-thm} for toroidal compactifications. We will freely use the notation there.
\subsubsection{}We first begin with the construction of $\ca{S}_K^\mmin(G,X_b)$.
\begin{construction}\label{const-min-quo-b}\upshape
Since $\sh_{K_0^\alpha}(G_0,X_0)$ is open dense in $\ca{S}^\mmin_{K_0^\alpha}(G_0,X_0)$, as Lemma \ref{lem-const-tor-Gb}, $\Delta_{\lcj{K}{\alpha}}(G_0,G)$ acts on $\ca{S}^\mmin_{K_0^\alpha}(G_0,X_0)$ through a finite group. 
Let 
\begin{equation}\ca{S}_K^\mmin(G,X_b):=\disju_{\alpha\in I_{G/G_0}}\ca{S}_{K_0^\alpha}^\mmin(G_0,X_0)/\Delta_{\lcj{K}{\alpha}}(G_0,G).\end{equation} 
The existence of $\Delta_{\lcj{K}{\alpha}}(G_0,G)$-action on $\ca{S}_{K_0^\alpha}^\mmin(G_0,X_0)$ follows from Lemma \ref{lem-gadq1-min} above and \cite[Prop. 5.2.13]{Mad19}.
Since $\sh_{K}(G,X_b)=\disju_{\alpha\in I_{G/G_0}}\sh_{K_0^\alpha}(G_0,X_0)/\Delta_{\lcj{K}{\alpha}}(G_0,G)$ and since the finite surjective morphism $\disju_{\alpha\in I_{G/G_0}}\sh^\mmin_{K_0^\alpha}(G_0,X_0)\to \sh^\mmin_K(G,X_b)$ factors through a finite surjective morphism 
$$\sh^\mmin_{K_0^\alpha}(G_0,X_0)/\Delta_{\lcj{K}{\alpha}}(G_0,G)\to \sh^\mmin_{K}(G,X_b),$$ 
the generic fiber of $\ca{S}_K^\mmin(G,X_b)$ is $\sh^\mmin_K(G,X_b)$ by Zariski's main theorem.
By construction and \cite[Ch. 4, Prop. 1.5]{Knu71}, $\ca{S}_K^\mmin(G,X_b)$ is a normal scheme which is flat and projective over $\ca{O}'_{(v)}$.\hfill$\square$
\end{construction}
\subsubsection{}\label{subsubsec-x0toxb}
As in \S\ref{subsec-main-thm}, the subscript $\ca{O}_{K_Z}$ denotes normalized base change to $\ca{O}_{K_Z}$. 
The diagram formed by taking quotient of $\Delta_{\lcj{K}{\alpha}}(G_0,G)$ 
\begin{equation}\label{diag-quo-tor-min}
    \begin{tikzcd}
    \disju_{I_{G/G_0}}\ca{S}^{\Sigma_0^\alpha}_{K_0^\alpha,\ca{O}_{K_Z}}\arrow[d,"\pi^b"]\arrow[rr,"{\disju\oint_{K_0^\alpha,\ca{O}_{K_Z}}^{\Sigma_0^\alpha}}"]&&\disju_{I_{G/G_0}}\ca{S}^\mmin_{K_0^\alpha,\ca{O}_{K_Z}}\arrow[d,"\pi^b"]\\
    \ca{S}^\Sigma_{K,\ca{O}_{K_Z}}\arrow[rr,dashed,"\oint_{K,\ca{O}_{K_Z}}^\Sigma"]&&\ca{S}^\mmin_{K,\ca{O}_{K_Z}}
    \end{tikzcd}
\end{equation}
induces a morphism $\oint^\Sigma_{K,\ca{O}_{K_Z}}:\ca{S}^\Sigma_{K,\ca{O}_{K_Z}}\to \ca{S}^\mmin_{K,\ca{O}_{K_Z}}$. (Again we can check the actions of $\Delta_{\lcj{K}{\alpha}}(G_0,G)$ on both source and target in the first row are compatible as they all contain a dense subscheme $\sh_{K_0^\alpha}$.) This morphism is proper and surjective since both source and target are proper over $\ca{O}_{K_Z}$ and since all solid arrows in (\ref{diag-quo-tor-min}) are surjective.
\begin{lem}\label{lem-xb-geom-conn-fib}
The normalization of $\ca{S}^\mmin_{K,\ca{O}_{K_Z}}$ in $\ca{S}^\Sigma_{K,\ca{O}_{K_Z}}$ is $\ca{S}^\mmin_{K,\ca{O}_{K_Z}}$ itself. The morphism $\oint_{K,\ca{O}_{K_Z}}^\Sigma$ has geometrically connected fibers.    
\end{lem}
\begin{proof}
Let $\ca{S}'$ be the normalization of $\ca{S}^\mmin_{K,\ca{O}_{K_Z}}$ in $\ca{S}^\Sigma_{K,\ca{O}_{K_Z}}$. Taking Stein factorization commutes with flat base change; combing with the discussion in the first prargraph of \S\ref{subsubsec-review-min-gen}, we check over generic fiber that the induced finite map $\ca{S}'\to \ca{S}^\mmin_{K,\ca{O}_{K_Z}}$ is birational. We then have the conclusion by Zariski's main theorem and \cite[\href{https://stacks.math.columbia.edu/tag/0A1B}{Thm. 0A1B}]{stacks-project}. 
\end{proof}
\subsubsection{} Fix any index $\alpha\in I_{G/G_0}$. Suppose that $(\Phi_0^\alpha,\sigma_0^\alpha)$ maps to $(\Phi,\sigma)$ as before. \par
Denote by $\oint_{[(\Phi_0^\alpha,\sigma_0^\alpha)],K_0^\alpha}$ the restriction $\oint_{K_0^\alpha}^{\Sigma_0^\alpha}|_{\ca{Z}_{[(\Phi_0^\alpha,\sigma_0^\alpha)],K_0^\alpha}}$. Similar conventions will be applied to restrictions to other strata or their (normalized) base changes.\par
By Theorem \ref{thm-cpt-min-hodge} and (\ref{diag-quo-tor-min}), the image of $\pi^{b,-1}(\ca{Z}_{[ZP^b(\Phi,\sigma)],K,\ca{O}_{K_Z}})=$
$$\disju_{\alpha\in I_{G/G_0}}\disju_{\pi^b(g_0^\alpha)\alpha\sim g^b}(\Delta_{\lcj{K}{\alpha}}(G_0,G)\ca{Z}_{[(\Phi_0^\alpha,\sigma_0^\alpha)],K_0^\alpha,\ca{O}_{K_Z}})$$
under $\disju \oint^{\Sigma_0^\alpha}_{K_0^\alpha,\ca{O}_{K_Z}}$ is 
$$\disju_{\alpha\in I_{G/G_0}}\disju_{\pi^b(g_0^\alpha)\alpha\sim g^b}(\Delta_{\lcj{K}{\alpha}}(G_0,G)\ca{Z}_{[\Phi_0^\alpha],K_0^\alpha,\ca{O}_{K_Z}}).$$
\begin{construction}\label{const-zpstrata-min}\upshape
From the paragraph above, it makes sense to define 
$$ \ca{Z}_{[ZP^b(\Phi)],K,\ca{O}_{K_Z}}\iso \ca{Z}_{[ZP^a(\Phi)],K,\ca{O}_{K_Z}}:=
\disju_{\alpha\in I_{G/G_0}}\disju_{\pi^b(g_0^\alpha)\alpha\sim g^b}(\Delta_{\lcj{K}{\alpha}}(G_0,G)\ca{Z}_{[\Phi_0^\alpha],K_0^\alpha,\ca{O}_{K_Z}})/\Delta_{\lcj{K}{\alpha}}(G_0,G).$$
\end{construction}
Define $\oint_{[ZP^b(\Phi,\sigma)],K,\ca{O}_{K_Z}}$ to be the projection induced by the restriction $\oint_{[(\Phi_0^\alpha,\sigma_0^\alpha)],K_0^\alpha,\ca{O}_{K_Z}}$.
\begin{lem}\label{lem-zpstrata-min} The following statements about $\ca{Z}_{[ZP^b(\Phi)],K,\ca{O}_{K_Z}}$ are true:
\begin{enumerate}
\item The generic fiber of the scheme above is $\mrm{Z}_{[ZP^a(\Phi)],K,E'_{K_Z}}\iso \mrm{Z}_{[ZP^b(\Phi)],K,E'_{K_Z}}$.
\item It is a locally closed normal flat subscheme of $\ca{S}_{K,\ca{O}_{K_Z}}^\mmin$ fitting into the commutative diagram
\begin{equation}\label{diag-rest-strata-zpb}
\begin{tikzcd}
\ca{Z}_{[(\Phi_0^\alpha,\sigma_0^\alpha)],K_0^\alpha,\ca{O}_{K_Z}}\arrow[rrr,"{\oint_{[(\Phi_0^\alpha,\sigma_0^\alpha)],K_0^\alpha,\ca{O}_{K_Z}}}"]\arrow[dd]\arrow[dr,hook]&&&\ca{Z}_{[\Phi_0^\alpha],K_0^\alpha,\ca{O}_{K_Z}}\arrow[rd,hook]\arrow[dd]&\\
&\ca{S}^{\Sigma_0^\alpha}_{K_0^\alpha,\ca{O}_{K_Z}}\arrow[rrr,"{\oint^{\Sigma_0^\alpha}_{K_0^\alpha,\ca{O}_{K_Z}}}"]\arrow[dd]&&&\ca{S}_{K_0^\alpha,\ca{O}_{K_Z}}^\mmin\arrow[dd]\\
\ca{Z}_{[ZP^b(\Phi,\sigma)],K,\ca{O}_{K_Z}}\arrow[rrr,"{\oint_{[ZP^b(\Phi,\sigma)],K,\ca{O}_{K_Z}}}"]\arrow[dr,hook,"{\text{Lem. \ref{lem-strata-zpb}}}"']&&&\ca{Z}_{[ZP^b(\Phi)],K,\ca{O}_{K_Z}}\arrow[rd,hook,"j"]&\\
&\ca{S}_{K,\ca{O}_{K_Z}}^\Sigma\arrow[rrr,"{\oint_{K,\ca{O}_{K_Z}}^\Sigma}"]&&&\ca{S}_{K,\ca{O}_{K_Z}}^\mmin.
\end{tikzcd}    
\end{equation}
\end{enumerate}
\end{lem}
\begin{proof}
The first statement follows from Proposition \ref{prop-shimura-sur} and Corollary \ref{zp-ab}. The scheme $\ca{Z}_{[ZP^b(\Phi)],K,\ca{O}_{K_Z}}$ is normal and is flat over $\ca{O}'_{(v)}$ by construction. The bottom square is induced by quotient construction, so (\ref{diag-rest-strata-zpb}) is commutative. We now show $j$ is an embedding. Denote the image of $j$ by $[\ca{Z}_\Phi]$; it is a locally closed subscheme by induction on the dimension of strata. Moreover, $j:\ca{Z}_{[ZP^b(\Phi)],K,\ca{O}_{K_Z}}\to [\ca{Z}_\Phi]$ is finite, surjective and a bijection on geometric points by construction.
From \cite[Cor. 5.2.8]{Mad19}, we see that for any geometric point $\geom{y}$ on $\ca{Z}_{[\Phi_0^\alpha],K_0^\alpha,\ca{O}_{K_Z}}$, there is a structural morphism
$\cpl{\ca{S}^\mmin_{K_0^\alpha,\ca{O}_{K_Z}}}{\geom{y}}\to \cpl{\ca{Z}_{[\Phi_0^\alpha],K_0^\alpha,\ca{O}_{K_Z}}}{\geom{y}}$
whose pre-composition with the morphism $\cpl{\ca{Z}_{[\Phi_0^\alpha],K_0^\alpha,\ca{O}_{K_Z}}}{\geom{y}}\to \cpl{\ca{S}^\mmin_{K_0^\alpha,\ca{O}_{K_Z}}}{\geom{y}}$ induced by the embedding $\ca{Z}_{[\Phi_0^\alpha],K_0^\alpha,\ca{O}_{K_Z}}\to\ca{S}^\mmin_{K_0^\alpha,\ca{O}_{K_Z}}$ is an identity. \par
For any geometric point $\geom{x}$ on $\ca{Z}_{[ZP^b(\Phi)],K,\ca{O}_{K_Z}}$ mapping to $[\geom{x}]$ on $[\ca{Z}_\Phi]$, we have a sequence (cf. Lemma \ref{lem-strata-zpb})
$$\cpl{\ca{Z}_{[ZP^b(\Phi)],K,\ca{O}_{K_Z}}}{\geom{x}}\to\cpl{\ca{S}^\mmin_{K,\ca{O}_{K_Z}}}{[\geom{x}]}\to\cpl{\ca{Z}_{[ZP^b(\Phi)],K,\ca{O}_{K_Z}}}{\geom{x}}$$
who factors through $\cpl{[\ca{Z}_\Phi]}{[\geom{x}]}$ with composition an identity. This forces $j$ to be an embedding.
\end{proof}
From the proof above, we have 
\begin{lem}\label{lem-str-morphism-gxb}
For any $\geom{x}\in \ca{Z}_{[ZP^b(\Phi,\sigma)],K,\ca{O}_{K_Z}}(\overline{\bb{F}}_p)$, the homomorphism $$\wat{\ca{O}}_{\ca{S}^\mmin_{K,\ca{O}_{K_Z}},\geom{x}}\to \wat{\ca{O}}_{\ca{Z}_{[ZP^b(\Phi)],K,\ca{O}_{K_Z}},\geom{x}}$$ between complete local rings induced by the embedding $\ca{Z}_{[ZP^b(\Phi,\sigma)],K,\ca{O}_{K_Z}}\hookrightarrow \ca{S}^\mmin_{K,\ca{O}_{K_Z}}$ admits a section
$$\wat{\ca{O}}_{\ca{Z}_{[ZP^b(\Phi)],K,\ca{O}_{K_Z}},\geom{x}}\to \wat{\ca{O}}_{\ca{S}^\mmin_{K,\ca{O}_{K_Z}},\geom{x}}$$
such that the pre-composition of the homomorphism with the section is an identity.
\end{lem}
Recall that we constructed $$\ca{S}_{\wdtd{K}_{\Phi}^b,h}:=\disju_{\beta\in I_{\Phi}}\gamma_\beta\cdot \ca{S}_{K_{\Phi_\beta},h}/\Delta_{K^\beta}(P_0,ZP_{\Phi})$$
in Construction \ref{const-zp-mix-sh}.\par 
\begin{lem}\label{lem-zpb-completion-iso-min}
There is an isomorphism $\ca{Z}_{[ZP^b(\Phi)],K,\ca{O}_{K_Z}}\iso \Delta_{\Phi,K}^{ZP}\bss \ca{S}_{\K^b_\Phi,h,\ca{O}_{K_Z}}$.    
\end{lem}
\begin{proof}
    This follows from the same proof as Proposition \ref{prop-zpb-completion-iso}.
\end{proof}
Then $\oint_{[ZP^b(\Phi,\sigma)],K,\ca{O}_{K_Z}}$ is identified with the canonical projection $\ca{S}_{\K^b_\Phi,\sigma,\ca{O}_{K_Z}}\to \Delta^{ZP}_{\Phi,K}\bss \ca{S}_{\K^b_\Phi,h,\ca{O}_{K_Z}}$ by density and the corresponding result in characteristic zero theory.
\subsubsection{}\label{subsubsec-xbtoxa}
Since $\ca{S}_{\K_\Phi^a,h,\ca{O}_{K_Z}}\iso \ca{S}_{\K_\Phi^b,h,\ca{O}_{K_Z}}$, $\ca{Z}_{[ZP^a(\Phi)],K,\ca{O}_{K_Z}}\iso \ca{Z}_{[ZP^b(\Phi)],K,\ca{O}_{K_Z}}$ and $\ca{S}^{\Sigma}_{K}(G,X_a)_{\ca{O}_{K_Z}}\iso\ca{S}^\Sigma_K(G,X_b)_{\ca{O}_{K_Z}}$,
let $\ca{S}^\mmin_{K}(G,X_a)_{\ca{O}_{K_Z}}:=\ca{S}^\mmin_K(G,X_b)_{\ca{O}_{K_Z}}$.\par
From the Galois action of $\gal(E'_{K_Z}/E')$ on strata with subscript $b$, we can obtain the action of it for corresponding strata with subscript $a$ as follows:
Let $r_b: \gal(E'_{K_Z}/E')\times \ca{Z}_{[ZP^b(\Phi)],K,\ca{O}_{K_Z}}\to \ca{Z}_{[ZP^b(\Phi)],K,\ca{O}_{K_Z}}$ be the action of $\gal(E'_{K_Z}/E')$ on $\ca{Z}_{[ZP^b(\Phi)],K,\ca{O}_{K_Z}}$ which determines the descent datum of $\mrm{Z}_{[ZP^b(\Phi)],K,E'_{K_Z}}$ from $E'_{K_Z}$ to $E'$. Define $r_a:\gal(E'_{K_Z}/E')\times \ca{Z}_{[ZP^b(\Phi)],K,\ca{O}_{K_Z}}\to \ca{Z}_{[ZP^a(\Phi)],K,\ca{O}_{K_Z}}$ by sending each $\sigma\in \gal(E'_{K_Z}/E')$ to $r_a(\sigma):=r_{E',{K_{Z,p}}}(Z,\{c\})(\sigma)\circ r_b(\sigma)=r_b(\sigma)\circ r_{E',{K_{Z,p}}}(Z,\{c\})(\sigma)$. 
Note that we can do this since the action of $Z(\A)$ is well defined by Proposition \ref{prop-mixsh-independent} and Lemma \ref{lem-zpb-completion-iso-min}. By \cite[Lem. A.3.4]{Mad19}, there is an action of $Z(\A)$ on $\ca{S}^\mmin_{K,\ca{O}_{K_Z}}$ extending the Hecke action on $\sh_K^\mmin$ and the strata $\ca{Z}_{[ZP^b(\Phi)],K,\ca{O}_{K_Z}}$. Note that, to run the proof in 
\cite[Lem. A.3.4]{Mad19}, one only needs the condition that there are morphisms between strata that are compatible with a morphism over the generic fiber for the whole $\sh_{K}^\mmin$.\par
We now do the quotient of $\ca{S}^\mmin_{K,\ca{O}_{K_Z}}$ and $\ca{Z}_{[ZP^b(\Phi)],K,\ca{O}_{K_Z}}$ by $\gal(E'_{K_Z}/E')$ via $r_a$. 
\begin{construction}\label{const-min-quo}\upshape
Let $\ca{S}^\mmin_K(G,X_a)$ (resp. $\ca{Z}_{[ZP^a(\Phi)],K}$) be the quotient of $\ca{S}^\mmin_{K,\ca{O}_{K_Z}}$ (resp. $\ca{Z}_{[ZP^b(\Phi)],K,\ca{O}_{K_Z}}$) by $\gal(E'_{K_Z}/E')$ via $r_a$. Then $\ca{S}_K^\mmin(G,X_a)$ is normal and is projective (see \cite[Ch. 4, Prop. 1.5]{Knu71}) and flat over $\ca{O}'_{(v)}$, and $\ca{Z}_{[ZP^a(\Phi)],K}$ is normal and is flat over $\ca{O}'_{(v)}$.
\end{construction}
\begin{prop}\label{prop-sum-min-xa}
Taking the quotient by $\gal(E'_{K_Z}/E')$ via $r_a$ induces a proper morphism $\oint_{K}^{\Sigma}:\ca{S}_K^{\Sigma}(G,X_a)\to \ca{S}^\mmin_K(G,X_a)$ with geometrically connected fibers. 
The normal scheme $\ca{Z}_{[ZP^a(\Phi)],K}$ is a locally closed subscheme of $\ca{S}^\mmin_K(G,X_a)$.\par
The restriction $\oint_{[ZP^a(\Phi)],K}$ of $\oint_K^\Sigma$ to $\ca{Z}_{[ZP^a(\Phi,\sigma)],K}$ is isomorphic to the natural projection $\ca{Z}_{[ZP^a(\Phi,\sigma)],K}\iso \ca{S}_{\K_\Phi^a,\sigma}\to \Delta_{\Phi,K}^{ZP}\bss \ca{S}_{K_{\Phi},h}\iso \ca{Z}_{[ZP^a(\Phi)],K}$. 
\end{prop}
\begin{proof}
Since this is still doing a quotient, all arguments will be similar to what we have seen. To show $\oint^\Sigma_{K}$ has geometrically connected fibers, we let $\ca{S}''$ be the normalization in $\ca{S}^{\Sigma}_{K}(G,X_a)$ of $\ca{S}_K^\mmin(G,X_a)$. But the induced morphism $\ca{S}''\to \ca{S}_K^\mmin(G,X_a)$ between normal schemes is finite and birational. We have the conclusion by Zariski's main theorem and \cite[\href{https://stacks.math.columbia.edu/tag/0A1B}{Thm. 0A1B}]{stacks-project}.\par
The quotient process induces a quasi-finite birational morphism $j: \ca{Z}_{[ZP^a(\Phi)],K}\to \ca{S}^\mmin_K(G,X_a)$ which is a bijection over geometric points onto its image $[\ca{Z}'']$.
By induction on the dimension of strata, this image $[\ca{Z}'']$ is locally closed in $\ca{S}^\mmin_K(G,X_a)$. To show that $j$ is an embedding, we again reduce the question to studying complete local rings at geometric points. By Lemma \ref{lem-str-morphism-gxb}, 
for any geometric point $\geom{z}\in \ca{Z}_{[ZP^a(\Phi)],K}(\overline{\bb{F}}_p)$ mapping to $[\geom{z}]\in \ca{S}^\mmin_K(G,X_a)(\overline{\bb{F}}_p)$, the homomorphism between complete local rings induced by $j$ 
$$\wat{\ca{O}}_{\ca{S}^\mmin_K(G,X_a),[\geom{z}]}\iso (\wat{\ca{O}}_{\ca{S}^\mmin_{K,\ca{O}_{K_Z}},[\geom{z}]})^{\gal(E'_{K_Z}/E')}\to \wat{\ca{O}}_{\ca{Z}_{[ZP^a(\Phi)],K},\geom{z}}$$
admits a section $\wat{\ca{O}}_{\ca{Z}_{[ZP^a(\Phi)],K,\geom{z}}}\to\wat{\ca{O}}_{\ca{S}^\mmin_K(G,X_a),[\geom{z}]}$.
This forces $\wat{\ca{O}}_{\ca{Z}_{[ZP^a(\Phi)],K},\geom{z}}\to\wat{\ca{O}}_{[\ca{Z}''],[\geom{z}]}$ to be an isomorphism. Here $\wat{\ca{O}}_{\ca{S}^\mmin_{K,\ca{O}_{K_Z}},[\geom{z}]}$ denotes the complete semi-local ring at the inverse image of $[\geom{z}]$ under the quotient map.\par
The second paragraph follows from the following diagram, which is commutative since it is so over the generic fiber:
\begin{equation}\label{diag-rest-strata-zpb-gal}
\begin{tikzcd}
\ca{Z}_{[ZP^a(\Phi,\sigma)],K,\ca{O}_{K_Z}}\arrow[rrr,"{\oint_{[ZP^a(\Phi,\sigma)],K,\ca{O}_{K_Z}}}"]\arrow[dd]\arrow[dr,hook]&&&\ca{Z}_{[ZP^a(\Phi)],K,\ca{O}_{K_Z}}\arrow[rd,hook]\arrow[dd]&\\
&\ca{S}^{\Sigma}_{K,\ca{O}_{K_Z}}\arrow[rrr,"{\oint^{\Sigma}_{K,\ca{O}_{K_Z}}}"]\arrow[dd]&&&\ca{S}_{K,\ca{O}_{K_Z}}^\mmin\arrow[dd]\\
\ca{Z}_{[ZP^a(\Phi,\sigma)],K}\arrow[rrr,"{\oint_{[ZP^a(\Phi,\sigma)],K}}"]\arrow[dr,hook,"{\text{Lem. \ref{lem-completion-comparison}}}"']&&&\ca{Z}_{[ZP^a(\Phi)],K}\arrow[rd,hook,"j"]&\\
&\ca{S}_K^\Sigma(G,X_a)\arrow[rrr,"{\oint_K^\Sigma}"]&&&\ca{S}_{K}^\mmin(G,X_a).
\end{tikzcd}    
\end{equation}
\end{proof}
\subsubsection{}\label{subsubsec-xatox2} Recall that by Corollary \ref{ext-cpt-imp} and our construction, $K$ is chosen such that $\sh_{K_2}^\mmin(G_2,X_2)\to \sh_K^\mmin(G,X_a)$ and $\mrm{Z}_{[\Phi_2],K_2}\to \mrm{Z}_{[\Phi],K}$ are open and closed embeddings for each $\Phi_2\mapsto \Phi$.
\begin{construction}\label{const-xatox2}\upshape
As the final step, construct $\ca{S}_{K_2}^\mmin(G_2,X_2)$ as the normalization of $\ca{S}_K^\mmin(G,X_a)$ in $\sh_{K_2}^\mmin$, and construct $\ca{Z}_{[\Phi_2],K_2}$ as the normalization of $\ca{Z}_{[ZP^a(\Phi)],K}$ in $\mrm{Z}_{[\Phi_2],K_2}$.
\end{construction}
\subsection{Main theorem on minimal compactifications}\label{subsec-main-theorem-min}Let us now summarize the main theorem:
\begin{thm}\label{thm-main-theorem-min}
With the conventions in Theorem \ref{maintheorem}, there is a normal scheme $\ca{S}_{K_2}^\mmin$ that is projective and flat over $\ca{O}_2$, which satisfies the following properties:
\begin{enumerate}[label=(\textrm{\ref{thm-main-theorem-min}}.\arabic*)]
\item\label{min-1}The scheme $\ca{S}_{K_2}$ is an open dense subscheme of $\ca{S}^\mmin_{K_2}$ with a natural open embedding $J^\mmin:\ca{S}_{K_2}\hookrightarrow\ca{S}_{K_2}^\mmin$.
\item\label{min-2}The compactification $\ca{S}_{K_2}^\mmin$ has a stratification of locally closed normal subschemes 
$$\ca{S}_{K_2}^\mmin=\disju_{[\Phi_2]\in\cusp_{K_2}(G_2,X_2)}\ca{Z}_{[\Phi_2],K_2},$$
such that, for each $[\Phi_2]$, $\ca{Z}_{[\Phi_2],K_2}$ is flat over $\ca{O}_2$ extending the stratum $\mrm{Z}_{[\Phi_2],K_2}$ defined in Pink's theory. The closure $\overline{\ca{Z}}_{[\Phi_2],K_2}$ of $\ca{Z}_{[\Phi_2],K_2}$ in $\ca{S}_{K_2}^\mmin$ satisfies $\overline{\ca{Z}}_{[\Phi_2],K_2}=\disju_{[\Phi_2']\preceq [\Phi_2]}\ca{Z}_{[\Phi_2'],K_2}$.
\item\label{min-3}There is a proper surjective morphism 
$$\oint_{K_2}^{\Sigma_2}:\ca{S}_{K_2}^{\Sigma_2}\lra \ca{S}^\mmin_{K_2}$$
with geometrically connected fibers such that $\oint_{K_2}\circ J^{\Sigma_2}=J^\mmin$. The morphism $\oint_{K_2}^{\Sigma_2}$ is compatible with stratifications on the source and the target in the sense that 
$$\oint^{\Sigma_2,-1}_{K_2}(\ca{Z}_{[\Phi_2],K_2})=\disju_{\Upsilon_2\in \cusp_{K_2}(G_2,X_2,\Sigma_2)\mapsto [\Phi_2]}\ca{Z}_{\Upsilon_2,K_2}.$$
\item\label{min-4} For each $[\Phi_2]$, $\ca{Z}_{[\Phi_2],K_2}$ is isomorphic to $\Delta_{\Phi_2,K_2}\bss \ca{S}_{K_{\Phi_2},h}$. The restriction of $\oint_{K_2}^{\Sigma_2}$ to any stratum $\ca{Z}_{\Upsilon_2,K_2}$ of $\ca{S}_{K_2}^{\Sigma_2}$ labeled by $\Upsilon_2=[(\Phi_2,\sigma_2)]$ is the natural projection (see (\ref{eq-2-step}) for conventions) $$\ca{Z}_{\Upsilon,K_2}\iso \Delta^\circ_{\Phi_2,K_2}\bss \ca{S}_{K_{\Phi_2},\sigma_2}\xrightarrow{\mbf{p}_2\circ\mbf{p}_1}\Delta_{\Phi_2,K_2}^\circ\bss\ca{S}_{K_{\Phi_2},h}\to \Delta_{\Phi_2,K_2}\bss \ca{S}_{K_{\Phi_2},h}\iso \ca{Z}_{[\Phi_2],K_2}.$$
\end{enumerate}
\end{thm}
\begin{proof}It suffices to show all statements above over $\ca{O}'_{(v)}$.
Let $\ca{Z}_{[\Phi],K}$ be the normalization of $\ca{Z}_{[ZP^a(\Phi)],K}$ in $\mrm{Z}_{[\Phi],K}$. By Corollary \ref{zp-ab}, $\ca{Z}_{[\Phi],K}$ is an open and closed subscheme of $\ca{Z}_{[ZP^a(\Phi)],K}$. By construction, $\ca{Z}_{[\Phi_2],K_2}\xrightarrow{\sim} \ca{Z}_{[\Phi],K}$, and $\ca{S}_{K_2}^\mmin\hookrightarrow\ca{S}_K^\mmin$ is an open and closed embedding. So $\ca{Z}_{[\Phi_2],K_2}$ is a locally closed normal subscheme of $\ca{S}_{K_2}^\mmin$ and is flat over $\ca{O}_{(v)}'$ for each $[\Phi_2]$. So \ref{min-1} is proved.\par
It follows from Proposition \ref{zp-ab-min}, Proposition \ref{prop-sum-min-xa} and Zariski's main theorem that $\Delta_{\Phi_2,K_2}\bss\ca{S}_{K_{\Phi_2},h}\iso \ca{Z}_{[\Phi_2],K_2}$.\par
The morphism $\oint_{K_2}^{\Sigma_2}$ and its restriction $\oint_{[(\Phi_2,\sigma)],K_2}$ to $\ca{Z}_{[\Phi_2,\sigma_2],K_2}$ is given by the following diagram
\begin{equation}\label{diag-rest-strata-final}
\begin{tikzcd}
\ca{Z}_{[(\Phi_2,\sigma_2)],K_2}\arrow[rrr,"{\oint_{[(\Phi_2,\sigma_2)],K_2}}"]\arrow[dd,"\simeq"]\arrow[dr,hook]&&&\ca{Z}_{[\Phi_2],K_2}\arrow[rd,hook]\arrow[dd,"\simeq"]&\\
&\ca{S}^{\Sigma_2}_{K_2}\arrow[rrr,"{\oint^{\Sigma_2}_{K_2}}"]\arrow[dd,hook]&&&\ca{S}_{K_2}^\mmin\arrow[dd,hook]\\
\ca{Z}_{[(\Phi,\sigma)],K}\arrow[rrr,"{\oint_{[(\Phi,\sigma)],K}}"]\arrow[dr,hook]&&&\ca{Z}_{[\Phi],K}\arrow[rd,hook,"j"]&\\
&\ca{S}_K^\Sigma(G,X_a)\arrow[rrr,"{\oint_K^\Sigma}"]&&&\ca{S}_{K}^\mmin(G,X_a)
\end{tikzcd}    
\end{equation}
by pulling back the bottom square along $\ca{S}_{K_2}^\mmin\hookrightarrow\ca{S}_K^\mmin(G,X_a)$.
Then \ref{min-3} and \ref{min-4} follow.
For \ref{min-2}, it remains to show the statement about the closure of a stratum. But this follows from \ref{min-3} and \ref{maintheorem-2}.
\end{proof}
\subsection{Pink's formula}\label{subsec-pink-formula}
We show Pink's formula when $K_{2,p}$ is an intersection of $n$ quasi-parahoric subgroups.
\subsubsection{}\label{subsubsec-pink-formula-finite-etale}It is not hard to verify the following statements.
\begin{lem}\label{lem-finite-etale-quotient}
Under Assumption (\ref{eq-assumption-n-quasi-para}), the natural quotient $\ca{S}_{\K_{\Phi}^a,h}\to \Delta^{ZP}_{\Phi,K}\bss \ca{S}_{\K_\Phi^a,h}$ is finite {\'e}tale.
\end{lem}
\begin{proof}
By assumption and Lemma \ref{lem-parahoric-auxiliary}, it suffices to prove the statement for $\ca{S}_{\K_{\Phi}^b,h}\to \Delta^{ZP}_{\Phi,K}\bss \ca{S}_{\K_\Phi^b,h}$. This follows from Part \ref{lem-twt-free-3} of Lemma \ref{lem-twt-free-action}, together with the fact that
$$\Delta^{ZP}_{\Phi,K}\bss\ca{S}_{\K_\Phi^b,h}\iso \disju_{\alpha\in I_{G/G_0}}\disju_{\pi(g_0^\alpha)\alpha\sim g^b} \Delta_{\lcj{K}{g_0^\alpha g_\alpha}}(P_{\Phi_0^\alpha},G)\bss\ca{S}_{K_{\Phi_0^\alpha},h}$$
from the proof of Proposition \ref{prop-zpb-completion-iso} (or Lemma \ref{lem-zpb-completion-iso-min}).
\end{proof}
We immediately obtain that
\begin{lem}[{cf. \cite[Assumption 4.3.1]{LS18b}}]\label{lem-lanstroh-assumption}
Under Assumption (\ref{eq-assumption-n-quasi-para}), $\overline{\ca{S}}^*_{K_{\Phi_2}}\to \ca{S}^*_{K_{\Phi_2},h}$ is an abelian scheme torsor and $\ca{S}^*_{K_{\Phi_2},h}\to \Delta_{\Phi_2,K_2}\bss \ca{S}_{K_{\Phi_2},h}$ is finite {\'e}tale.    
\end{lem}
\begin{proof}
This follows from Lemma \ref{lem-finite-etale-quotient} and Proposition \ref{prop-abelian-torsor}.
\end{proof}
Now we convert the left action of $\Delta_{\Phi_2,K_2}$ on $\ca{S}_{K_{\Phi_2},h}$ and $\ca{S}_{K_{\Phi_2},h}^*$ to a right action to match the conventions in \cite{Pin92}. 
Following \cite[3.7.4]{Pin92}, denote
$$H_{\Phi_2,K_2}:=g_{\Phi_2}K_2g_{\Phi_2}^{-1}\cap \stb_{Q_{\Phi_2}(\bb{Q})}(D_{\Phi_2,h})P_{\Phi_2}(\A),$$
$$H_{\Phi_2,K_2}^C:=g_{\Phi_2}K_2g_{\Phi_2}^{-1}\cap \mrm{Cent}_{Q_{\Phi_2}(\bb{Q})}(D_{\Phi_2,h})W_{\Phi_2}(\A),$$
and 
$$K_{W,\Phi_2}:=g_{\Phi_2}K_2g_{\Phi_2}^{-1}\cap W_{\Phi_2}(\A).$$
Recall that any representative $(h,\gamma)\in P_{\Phi_2}(\A)\times \stb_{Q_{\Phi_2}(\bb{Q})}(D_{\Phi_2})$ of $\Delta_{\Phi_2,K_2}$ acts on $\sh_{K_{\Phi_2}}\to\overline{\sh}_{K_{\Phi_2}}\to\sh_{K_{\Phi_2},h}$ by the left multiplication of $\gamma$, which is equivalent to the left conjugation of $\gamma$ followed by the right action of $h$; converting it to a right action means a right conjugation of $\gamma^{-1}$ followed by the right action of $h$. For any $(\gamma^{-1},h)\in \stb_{Q_{\Phi_2}(\bb{Q})}(D_{\Phi_2})\times P_{\Phi_2}(\A)$ such that $\gamma^{-1}h\in H_{\Phi_2,K_2}$, we define a right action of this pair given by the converted right action we just mentioned. It is easy to verify that for any $p\in P_{\Phi_2}(\bb{Q})=\stb_{Q_{\Phi_2}(\bb{Q})}(D_{\Phi_2})\cap P_{\Phi_2}(\A)$, the action of $(\gamma^{-1}p^{-1},ph)$ is equal to that of $(\gamma^{-1},h)$. Hence, the left action of $\Delta_{\Phi_2,K_2}$ is equivalent to the right action of $H_{\Phi_2,K_2}/K_{\Phi_2}$.
\subsubsection{}\label{subsubsec-pink-formula}
Now we state Pink's formula for general abelian-type Shimura varieties. Note that, since the assumption \cite[3.1.5]{Pin92} that the connected center of $G_2$ is an extension of a split $\bb{Q}$-torus and a compact-type torus is not true for some abelian type Shimura data that are not of Hodge type, the formulation is slightly different from those in \cite{Pin92} and \cite{LS18b}.\par
For any connected algebraic group $H$ over $\bb{Q}$, define \gls{Hc} to be the $\bb{Q}$-group formed by the quotient of $H$ by $Z_{ac}(H)$, where $Z_{ac}(H)$ is the minimal $\bb{Q}$-subgroup in the connected multiplicative center $Z(H)^\circ$ of $H$ such that $Z(H)^\circ/Z_{ac}(H)$ has equal $\bb{Q}$-rank and $\bb{R}$-rank. (Here, the subscript ``ac'' means ``anti-cuspidal''.) As remarked in \cite[footnote p.34]{KSZ21}, if we assume that $H$ is reductive and splits over a CM field, the definition coincides with that in \cite[Ch. III]{Mil90}.\par
By \cite{Del79} and \cite[1.5.8.2]{KSZ21}, we know that 
$$\gal(\sh(G_2,X_2)/\sh_{K_2}(G_2,X_2))\iso \frac{K_2}{Z_{G_2}(\bb{Q})^{\overline{\ }}\cap K_2},$$ 
where $Z_{G_2}(\bb{Q}){}^{\overline{\ }}\cap K_2$ is trivial when $K_2$ is neat and $Z_{G_2}$ is isogenous to a product of $\bb{Q}$-split tori and compact-type tori.\par 
We check if the Hecke actions away from $p$ are finite {\'e}tale:
\begin{lem}\label{lem-etale-hecke-action}
For all cases and for $?=*$ or $\emptyset$, consider prime-to-$p$ Hecke actions on $\ca{S}_{K_{2}}$, $\ca{S}_{K_{\Phi_2}}^?$, $\overline{\ca{S}}_{K_{\Phi_2}}^?$ and $\ca{S}_{K_{\Phi_2},h}^?$. More precisely, for any normal open compact subgroup $K_2^{p,\prime}\sbst K_2^p$, denote by $K_{\Phi_2}^{p,\prime}:=P_{\Phi_2}(\Ap)\cap g_{\Phi_2}K_2^{p,\prime}g_{\Phi_2}^{-1}$ a normal subgroup of $K_{\Phi_2}$. We have the following statements:
\begin{enumerate}
    \item\label{lem-etale-hecke-action-1} If $K_{2,p}$ satisfies Assumption (\ref{eq-assumption-n-quasi-para}), the natural morphisms $\ca{S}_{K_{2,p}K_2^{p,\prime}}\to\ca{S}_{K_{2,p}K_2^p}$, $\ca{S}_{K_{\Phi_2,p}K_{\Phi_2}^{p,\prime}}^?\to \ca{S}_{K_{\Phi_2,p}K_{\Phi_2}^{p}}^?$, $\overline{\ca{S}}_{K_{\Phi_2,p}K_{\Phi_2}^{p,\prime}}^?\to \overline{\ca{S}}_{K_{\Phi_2,p}K_{\Phi_2}^{p}}^?$ and    $\ca{S}_{K_{\Phi_2,p}K_{\Phi_2}^{p,\prime},h}^?\to \ca{S}_{K_{\Phi_2,p}K_{\Phi_2}^{p},h}^?$ 
    are all finite {\'e}tale for any neat open compact subgroup $K_2^p$ and $K_{\Phi_2}^p$.
    \item\label{lem-etale-hecke-action-2} In general, the statement above holds for sufficiently small $K_2^p$ and $K_{\Phi_2}^p$.
\end{enumerate}
\end{lem}
\begin{proof}
As in the proof of Lemma \ref{lem-torus}, by adjusting prime-to-$p$ level group, it is enough to show the statement for $?=*$, and this is further reduced to the statement for $\ca{S}_{\K^a_\Phi}$, $\overline{\ca{S}}_{\K^a_\Phi}$ and $\ca{S}_{\K^a_\Phi,h}$.\par 
Write $\K^a_\Phi$ as $\K_{\Phi,p}^a\K_{\Phi}^{a,p}$. 
If $K_{2,p}$ satisfies (\ref{eq-assumption-n-quasi-para}), by Lemma \ref{lem-parahoric-auxiliary}, the first part follows directly from Lemma \ref{lem-twt-free-action} by reducing to the related statement about $\ca{S}_{\K^b_\Phi}$.\par 
Now we show the second part. 
Denote $\K^+_{\Phi,p}:=ZP_\Phi(\bb{Q}_p)\cap g_\Phi K_p^+ g_\Phi^{-1}$ where $K_p^+$ is defined in \S\ref{subsubsec-s*true} above. 
Let $\K_\Phi^{a,p,\prime}$ be a normal open compact subgroup of $\K_\Phi^{a,p}$. Pick any $g^p\in \K_\Phi^{a,p}$ such that the action of it on $\ca{S}_{\K^{a,\prime}_\Phi}:=\ca{S}_{\K_{\Phi,p}^a\K_\Phi^{a,p,\prime}}$ is not free on every geometric point. 
By Lemma \ref{lem-twt-free-action}, the Hecke action of $g^p\in \K_\Phi^{a,p}$ on $\ca{S}_{\K_{\Phi,p}^{+,a}\K_\Phi^{a,p,\prime}}$ is not free if and only if its action on $\sh_{\K_{\Phi,p}^{+,a}\K_\Phi^{a,p,\prime}}$ is trivial. Writing $\K_{\Phi,p}^{+,a}\K_\Phi^{a,p,\prime}=:\K_\Phi^{+,a,\prime}$, the last condition is equivalent to $g^p\in \mrm{Cent}_{ZP_\Phi(\bb{Q})}(D_{\Phi})\K_\Phi^{+,a,\prime}$. 
Replacing $g^p$ with a multiplication of it by an element in $\K_\Phi^{a,p,\prime}$, assume that $g^p\in  \mrm{Cent}_{ZP_\Phi(\bb{Q})}(D_{\Phi})\K_{\Phi,p}^{+,a}$. Then $g^p$ comes from an element $c\in \mrm{Cent}_{ZP_\Phi(\bb{Q})}(D_\Phi)$ and $c\in  \mrm{Cent}_{ZP_\Phi(\bb{Q})}(D_{\Phi})\cap\K_{\Phi}^{+,a}= \mrm{Cent}_{ZP_\Phi(\bb{Q})}(D_{\Phi})\cap\K_{\Phi}^{a}K^+_{Z,p}$. That is, the prime-to-$p$ part of $c$, viewed as an $\Ap$-point, is equal to $g^p$.\par
Recall that $\mrm{Cent}_{ZP_\Phi(\bb{Q})}(D_{\Phi})$ is multiplicative and recall Chevalley's theorem \cite[Thm. 1]{Che51}. We have that $$\mrm{Cent}_{ZP_\Phi(\bb{Q})}(D_{\Phi})\cap\K_{\Phi}^{a}K^+_{Z,p}=\mrm{Cent}_{ZP_\Phi(\bb{Q})}(D_{\Phi})\cap\K_{\Phi}^{a}$$ when $\K_\Phi^{a,p}$ is sufficiently small. Choose such a sufficiently small $\K_\Phi^{a,p}$, 
we have that the action of $g^p$ is trivial on all geometric points of $\ca{S}_{\K^{a,\prime}_\Phi}$.\par
The proof for $\overline{\ca{S}}_{K_{\Phi}^a}$ and $\ca{S}_{\K^a_\Phi,h}$ is exactly the same as above.
\end{proof}
Fix a prime number $\ell\neq p$. Denote by $\ca{G}(K_2)_l$ the quotient 
$$\frac{K_2}{(K_2\cap G_2(\A^l))(Z_{G_2}(\bb{Q}){}^{\overline{\ }}\cap K_2)};$$
this is a profinite group with a natural homomorphism to $G_2^c(\bb{Q}_l)$. 
Let $\xi$ be a continuous representation of $\ca{G}(K_2)_l$ on a finite-dimensional $\overline{\bb{Q}}_l$-vector space $V_\xi$. By \cite[1.5.2]{KSZ21} (see also \cite[Lem. 3.6]{LS18} and \cite[III. 2]{HT00}), there is a lisse $\overline{\bb{Q}}_l$-sheaf $\ca{V}_\xi$ on $\ca{S}_{K_2}$ associated with $V_\xi$.\par 
Now assume that $\xi$ is an algebraic $\overline{\bb{Q}}_l$-representation of $G_2^c$ on a finite-dimensional $\overline{\bb{Q}}_l$-vector space $V_\xi$, then $\xi$ induces a continuous representation $\ca{G}(K_2)_l\to G^c_2(\bb{Q}_l)\xrightarrow{\xi} \mrm{GL}(V_\xi)$, which we abusively denote by $\xi$. We attach a lisse $\overline{\bb{Q}}_l$-sheaf $\ca{V}_\xi$ to $V_\xi$ as in the last paragraph.\par
On the other hand, 
the group $H^C_{\Phi_2,K_2}$ is by definition the kernel of the right action of $H_{\Phi_2,K_2}$ on $\sh_{\Phi_2,h}:=\varprojlim_{K_2}\sh_{K_{\Phi_2},h}$. 
There is a continuous representation of $H_{\Phi_2,K_2}$ on $V_\xi$ induced by
$$\xi_{\Phi_2}:H_{\Phi_2,K_2}\hookrightarrow g_{\Phi_2}K_2g_{\Phi_2}^{-1}\xrightarrow[\sim]{g_{\Phi_2}^{-1}(\cdot)g_{\Phi_2}}K_2\hookrightarrow G_2(\A)\to G_2^c(\bb{Q}_l)\xrightarrow{\xi}\mrm{GL}(V_\xi).$$
Viewing $V_\xi$ as a $H_{\Phi_2,K_2}$-representation, it makes sense to define the following object 
$$V^\natural_{\xi,\Phi_2,K_2}:=R\mrm{Inv}(H^C_{\Phi_2,K_2},R\mrm{Inv}(\lie W_{\Phi_2}(\bb{Q}_l),V_\xi))$$
in the derived category of bounded complexes of finite-dimensional $\overline{\bb{Q}}_l$-vector spaces with $H_{\Phi_2,K_2}/H^C_{\Phi_2,K_2}$-actions. By Part \ref{lem-etale-hecke-action-1} of Lemma \ref{lem-etale-hecke-action}, we can and will denote by $\ca{V}_{\xi,\Phi_2,K_2}$ the lisse $\overline{\bb{Q}}_l$-sheaf on $\ca{Z}_{\Phi_2,K_2}$ associated with the $H_{\Phi_2,K_2}/H^C_{\Phi_2,K_2}$-representation $V^\natural_{\xi,\Phi_2,K_2}$. 

From Theorem \ref{thm-main-theorem-min}, we have a diagram
$$\ca{S}_{K_2}\xhookrightarrow{J^\mmin}\ca{S}_{K_2}^\mmin\xhookleftarrow{i_{[\Phi_2],K_2}}\ca{Z}_{[\Phi_2],K_2}.$$
\begin{prop}[{Pink; cf. \cite[Thm. 4.2.1]{Pin92} and \cite[Thm. 4.3.10]{LS18b}}]\label{prop-pink-formula}
Under Assumption (\ref{eq-assumption-n-quasi-para}), for any algebraic representation $\xi$ of $G_2^c$ on a finite-dimensional $\overline{\bb{Q}}_l$-vector space $V_\xi$, we have
\begin{equation}\label{eq-pink-formula}
    i^*_{[\Phi_2],K_2}RJ_*^\mmin\ca{V}_\xi\iso \ca{V}^\natural_{\xi,\Phi_2,K_2}.
\end{equation}
\end{prop}
\begin{proof}
Pink's formula (\ref{eq-pink-formula}) can be shown by the same argument in \cite[4.3-4.7]{Pin92} with the geometric input developed in this article, and we will explain it for the convenience of the readers.\par
There is a finite extension $F_l$ of $\bb{Q}_l$ and an $\ca{O}_{F_l}$-lattice $W_\xi\sbst V_\xi$ that is stable under $\ca{G}(K_2)_l$ such that we can write $V_\xi$ as $\overline{\bb{Q}}_l\otimes_{\ca{O}_{F_l}}\varprojlim_{n}W_\xi/l^n W_\xi$ and $V_{\xi,n}:=W_\xi/l^n W_\xi$ are continuous $\bb{Z}/l^n\bb{Z}$-representations of $\ca{G}(K_2)_l$. Let $\ca{V}_{\xi,n}$ be the $\bb{Z}/l^n\bb{Z}$-sheaf on $\ca{S}_{K_2}$ attached to $V_{\xi,n}$. 
Let $\ca{Z}_{[\Phi_2],K_2}^\mrm{tor}$ be the preimage $\oint_{K_2}^{\Sigma_2,-1}(\ca{Z}_{[\Phi_2],K_2})$. There is a commutative diagram
\begin{equation*}
  \begin{tikzcd}  \ca{S}_{K_2}\arrow[rr,hook,"J^{\Sigma_2}"]\arrow[rrd,hook,"J^\mmin"]&&\ca{S}^{\Sigma_2}_{K_2}\arrow[d,"\oint_{K_2}^{\Sigma_2}"]&&\ca{Z}^{\mrm{tor}}_{[\Phi_2],K_2}\arrow[ll,hook',"i^{\mrm{tor}}_{[\Phi_2],K_2}"']\arrow[d,"\oint_{K_2}^{\Sigma_2}"]\\
    &&\ca{S}_{K_2}^\mmin&&\ca{Z}_{[\Phi_2],K_2}.\arrow[ll,hook',"i_{[\Phi_2],K_2}"']
    \end{tikzcd}
\end{equation*}
By Theorem \ref{min-3}, $\oint_{K_2}^{\Sigma_2}$ is proper. By proper base change, there is a canonical isomorphism $$i^*_{[\Phi_2],K_2}RJ_*^\mmin\ca{V}_{\xi,n}\iso R\oint^{\Sigma_2}_{K_2,*}i^{\mrm{tor},*}_{[\Phi_2],K_2}RJ_*^{\Sigma_2}\ca{V}_{\xi,n}.$$
The argument in \cite[Prop. 4.4.3]{Pin92} still works thanks to \cite[2.7]{Pin92} and Proposition \ref{prop-artin-approximation-strata-preserving}. 
We deduce a canonical isomorphism
$$R\oint^{\Sigma_2}_{K_2,*}i^{\mrm{tor},*}_{[\Phi_2],K_2}RJ_*^{\Sigma_2}\ca{V}_{\xi,n}\iso R\oint^{\Sigma_2}_{K_2,*}R\mrm{Inv}(H_{\Phi_2,K_2}, V_{\xi,n}\otimes \varinjlim_{U_l}\phi_{U_l,*}\bb{Z}),$$
where $U_l\sbst G_2(\bb{Q}_l)$ runs over a cofinal system of sufficiently small open compact subgroups of $G_2(\bb{Q}_l)$ such that $K_{2}^p\cap G_2(\A^{p,l})U_l\sbst K_2^p$, and $\phi_{U_l}:\ca{Z}_{[\Phi_2],K_{2,p}K_2^{p,l}U_l}^\mrm{tor}\to \ca{Z}^{\mrm{tor}}_{[\Phi_2],K_2}$ is the natural transition map.\par
The $\ca{S}_{K_{\Phi_2}}^*\to\overline{\ca{S}}_{K_{\Phi_2}}^*\to\ca{S}^*_{K_{\Phi_2},h}\to \Delta_{\Phi_2,K_2}\bss \ca{S}_{K_{\Phi_2},h}$ in our convention corresponds to $\Xi\to C\to C^{\mrm{st}}\to \mathsf{Z}$ in the conventions of \cite[Sec. 4.3]{LS18b}; all terms and morphisms in this sequence satisfy \cite[Prop. 2.1.2 and Assumption 4.3.1(4)]{LS18b} by Theorem \ref{maintheorem}, Theorem \ref{thm-main-theorem-min} and Lemma \ref{lem-lanstroh-assumption}.
Note that \cite[Lem. 4.3.2]{LS18b} follows from characteristic zero theory, Lemma \ref{lem-lanstroh-assumption} and \cite[I. Prop. 2.7]{FC90}. Together with Part \ref{lem-etale-hecke-action-1} of Lemma \ref{lem-etale-hecke-action}, the argument in \cite[4.6 and 4.7]{Pin92} works. We deduce that 
$R\oint^{\Sigma_2}_{K_2,*}R\mrm{Inv}(H_{\Phi_2,K_2}, V_{\xi,n}\otimes \varinjlim_{U_l}\phi_{U_l,*}\bb{Z})$ is isomorphic to the complex of $\bb{Z}/l^n\bb{Z}$-sheaves on $\ca{Z}_{[\Phi_2],K_2}$ in the derived category associated with $R\mrm{Inv}(H^C_{\Phi_2,K_2},V_{\xi,n})$. We now have that $i^*_{[\Phi_2],K_2}RJ_*^\mmin\ca{V}_\xi$ is isomorphic to the complex of $\overline{\bb{Q}}_l$-sheaves on $\ca{Z}_{[\Phi_2],K_2}$ in the derived category associated with $R\mrm{Inv}(H^C_{\Phi_2,K_2},V_{\xi})$. Finally, we obtain (\ref{eq-pink-formula}) by the proof in \cite[5.2 and 5.3]{Pin92}.
\end{proof}
\subsection{Nearby cycles}\label{subsec-nearby}
This subsection is to generalize \cite[Thm. 4.1]{LS18} to the abelian-type case. As the readers will see, the proof requires essentially nothing but (a generalization of) the qualitative descriptions of good compactifications stated there.\par 
\begin{rk}
Unlike the Case (Ab) in \cite{LS18} which only deals with a ``crude'' construction of minimal compactifications for integral models in \cite{KP15}, the proof presented here works for all integral models constructed in this paper without technical assumptions on $G_2$ and the level at $p$.
\end{rk}
\begin{prop}[{Lan-Stroh; cf. \cite[Thm. 4.1 and Cor. 4.6]{LS18}}]\label{prop-nearby}
Fix a prime $l\neq p$. 
Let $\ca{V}$ be either a lisse $\overline{\bb{Q}}_l$-sheaf associated with an algebraic representation $\xi$ of $G^c_2$ on a finite-dimensional $\overline{\bb{Q}}_l$-vector space $V_\xi$, or a finite $\bb{Z}_l$-module equipped with an action of an open compact subgroup $G_2^c(\bb{Z}_l)$ of $G_2^c(\bb{Q}_l)$. Denote by $\eta$ (resp. $s$) the generic (resp. special) point of the base ring of the integral model $\ca{S}_{K_2}$ and denote the corresponding geometric point by adding ``$\overline{\ }$''.\par
In all cases ($\mrm{HS}$), ($\mrm{STB}_n$) and ($\mrm{DL}$), assume that the projection of $K_2$ into $G^c_2(\bb{Q}_l)$ factors through $G^c_2(\bb{Z}_l)$. Then there are natural isomorphisms induced by adjunctions that are equivariant under the actions of the absolute Galois group $\gal(\overline{\eta}/\eta)$:
\begin{equation}\label{eq-nearby-tor}
R\Psi_{\ca{S}_{K_2}^{\Sigma_2}}RJ_{\eta,*}^{\Sigma_2}\ca{V}\xrightarrow{\sim} RJ_{\overline{s},*}^{\Sigma_2}R\Psi_{\ca{S}_{K_2}}\ca{V},
\end{equation}
\begin{equation}\label{eq-nearby-tor-shriek}
J_{\overline{s},!}^{\Sigma_2}R\Psi_{\ca{S}_{K_2}}\ca{V}\xrightarrow{\sim}R\Psi_{\ca{S}_{K_2}^{\Sigma_2}}J_{\eta,!}^{\Sigma_2}\ca{V},   
\end{equation}
\begin{equation}\label{eq-nearby}
  R\Gamma_{et}(\sh_{K_2,\overline{\eta}},\ca{V})\xrightarrow{\sim}R\Gamma_{et}(\ca{S}_{K_2,\overline{s}},R\Psi_{\ca{S}_{K_2}}\ca{V}),  
\end{equation}
and
\begin{equation}\label{eq-nearby-shriek}
  R\Gamma_{et,c}(\ca{S}_{K_2,\overline{s}},R\Psi_{\ca{S}_{K_2}}\ca{V})\xrightarrow{\sim} R\Gamma_{et,c}(\sh_{K_2,\overline{\eta}},\ca{V}).  
\end{equation}
In the equations above, ``$R\Psi$'' denotes nearby cycles (see \cite[Exp. XIII]{SGA7-II}).
\end{prop}
\begin{proof}
The last two isomorphisms (\ref{eq-nearby}) and (\ref{eq-nearby-shriek}) follow from the first two (\ref{eq-nearby-tor}) and (\ref{eq-nearby-tor-shriek}) and the proper base change. (Alternatively, when $K_{2,p}$ is parahoric, we can choose $\ca{O}_{K_Z}$ so that the normalized base changes in Corollary \ref{cor-maintheorem} are base changes by Lemma \ref{lem-parahoric-auxiliary}. Then one can also use the Hochschild-Serre spectral sequence and Corollary \ref{cor-maintheorem} to show (\ref{eq-nearby}) and (\ref{eq-nearby-shriek}) as in the proof of \cite[Cor. 4.6]{LS18}.)\par
We show the first two isomorphisms. The proof is almost identical to that in \cite{LS18} given the established geometric inputs in this paper. We will repeat it and highlight the difference. \par
It suffices to assume that $\ca{V}$ is associated with a finite $\bb{Z}/l^n\bb{Z}$-module $\Lambda_n$. We first choose a prime-to-$p$ level $K^{p,\prime}_2$ such that $K^{p,\prime}_2$ acts trivially on $\Lambda_n$. By Lemma \ref{lem-finite-etale-reduction-step}, we can choose suitable $K_2^{p,\prime}$ to further make the transition maps in Lemma \ref{lem-etale-hecke-action} (1) finite {\'e}tale. 
Then the argument in \cite[Prop. 5.1 and Prop. 5.3]{LS18} still works as summarized below.\par 
Denote the base ring by $\ca{O}$ until the end of this proof. We now fix a point $x$ in a stratum $\ca{Z}_\sigma:=\ca{Z}_{[(\Phi_2,\sigma_2)],K_2}$. For simplicity, write $\sigma:=\sigma_2$. We know that $\ca{S}_{K_2}^{\Sigma_2}$ is {\'e}tale locally isomorphic to $C'\times_{\ca{O}} \mbf{E}$ by Lemma \ref{lem-finite-etale-reduction-step} below, and that we can reduce the problem to $C'\times_{\ca{O}} \mbf{E}$ by assuming that $\ca{V}:=\ca{V}_E\boxtimes \ul{\bb{Z}/l^n\bb{Z}}$ for an {\'e}tale $\bb{Z}/l^n\bb{Z}$-module $\ca{V}_E$ over $\mbf{E}$ that is trivialized over $\mbf{E}'$ and a constant sheaf $\ul{\bb{Z}/l^n\bb{Z}}$ over $C'$. (We first have a constant sheaf $\ul{\Lambda_n}$ on $\ca{S}^*_{K_{\Phi_2,p}K_{\Phi_2}^{p,\prime}}$ and pull it back to $C'\times_{\ca{O}}\mbf{E}'$, and we then descend it to $C'\times_{\ca{O}}\mbf{E}$ with the given $(K_2^{p}\cap (Z_{G_2}\times U_{\Phi_2})(\Ap))$-action.)\par 
Note that if the coefficient we are considering is already trivial itself, we do not need the reduction step in the last two paragraphs and can directly reduce the problem to the problem over $C\times_{\ca{O}}\mbf{E}$ for trivial coefficients.\par 
Now the argument in \cite[pp.18-19]{LS18} gets through verbatim, and it shows that 
$$R\Psi_{C'\times_{\ca{O}} \mbf{E}(\sigma)}Rj_{\eta,*}\ca{V}\xrightarrow{\sim}Rj_{\overline{s},*}R\Psi_{C'\times_{\ca{O}}\mbf{E}}\ca{V}$$ and that
$$j_{\overline{s},!}R\Psi_{C'\times_{\ca{O}}\mbf{E}}\ca{V}\xrightarrow{\sim}R\Psi_{C'\times_{\ca{O}} \mbf{E}(\sigma)}j_{\eta,!}\ca{V}$$
for $j: C'\times_{\ca{O}}\mbf{E}\hookrightarrow C'\times_{\ca{O}}\mbf{E}(\sigma)$.
We now have the desired isomorphisms.
\end{proof}
\begin{lem}\label{lem-finite-etale-reduction-step}
Assume that $\ca{V}$ is associated with a finite $\bb{Z}/l^n\bb{Z}$-module $\Lambda_n$. Suppose that a neat open compact subgroup $K^{p,\prime\prime}_2\sbst K^{p}_2$ acts trivially on $\Lambda_n$ via $K^p_2\to \ca{G}(K_2)_l\to G_2^c(\bb{Q}_l)$. 
\begin{enumerate}
    \item We can choose a prime-to-$p$ level $K^{p,\prime}_2$ such that $K^{p,\prime}_2$ acts trivially on $\Lambda_n$ and the transition maps in Lemma \ref{lem-etale-hecke-action} (1) are finite {\'e}tale. 
    \item Fix any cusp label $[(\Phi_2,\sigma_2)]$. Under the choice in the first item, we have the following commutative diagram:
    \begin{equation}
\begin{tikzcd}
\ca{S}^*_{K_{\Phi_2,p}K_{\Phi_2}^{p,\prime}}\arrow[r]\arrow[d,"{\mbf{E}'}"]&\ca{S}^*_{K_{\Phi_2,p}K_{\Phi_2}^{p}}\arrow[d,"\mbf{E}"]\\
C':=\overline{\ca{S}}^*_{K_{\Phi_2,p}K_{\Phi_2}^{p,\prime}}\arrow[r]&C:=\overline{\ca{S}}^*_{K_{\Phi_2,p}K_{\Phi_2}^{p}}.
\end{tikzcd}  
    \end{equation}
\end{enumerate}
The horizontal maps are finite {\'e}tale, the vertical maps are torsors under split tori $\mbf{E}'$ and $\mbf{E}$, respectively, and the diagram is equivariant under a finite {\'e}tale isogeny between the tori $\mbf{E}'\to \mbf{E}.$
\end{lem}
\begin{proof}
We adopt the notation in the proof of Lemma \ref{lem-etale-hecke-action}. We reduce to the case where $(G_2,X_2)=(G,X_a)$. We see from the proof there that, if we replace $K^{p,\prime}$ with $K^{p,\prime}\cdot (Z_G(\bb{Q})\cap K_{Z,p}^+K)^p$ (where $(-)^p$ denotes viewing $\bb{Q}$-points as $\Ap$-points), then $\ca{S}_{K'}\to \ca{S}_K$ is finite {\'e}tale. By the proof of \cite[Cor. 4.10]{Pin89} and since $\mrm{Cent}_{ZP_\Phi(\bb{Q})}(D_\Phi)\cap \K_\Phi^{+,a}$ is neat, we have that $\mrm{Cent}_{ZP_\Phi(\bb{Q})}(D_\Phi)\cap \K_\Phi^{+,a}\sbst Z_G(\bb{Q})\cap\K_\Phi^{+,a}\sbst Z_G(\bb{Q})\cap K_{Z,p}^+K$. From the proof of Lemma \ref{lem-etale-hecke-action}, the transition maps $\ca{S}_{\K_\Phi^{a,\prime}}\to \ca{S}_{\K_\Phi^a}$ and $\overline{\ca{S}}_{\K_\Phi^{a,\prime}}\to \overline{\ca{S}}_{\K_\Phi^a}$ are finite {\'e}tale if we replace $K^{p,\prime}$ with $K^{p,\prime}\cdot (Z_G(\bb{Q})\cap K_{Z,p}^+K)^p$. Since the action of $K^p$ on $\Lambda_n$ factors through $G^c(\bb{Q}_l)$ and the projection of $Z_G(\bb{Q})\cap K_{Z,p}^+K$ to $G^c(\bb{Q}_l)$ is trivial, we have shown the first item. The second item follows from the first one and Theorem \ref{maintheorem}.
\end{proof}
The corollaries for intersection cohomology now also hold for all cases ($\mrm{HS}$), ($\mrm{STB}_n$) and ($\mrm{DL}$) by proper base change and by the proof of \cite[Thm. 5.26]{LS18i}.
\begin{cor}\label{cor-nearby-intersection}
Let $\ca{V}$ be a lisse $\overline{\bb{Q}}_l$-sheaf, or an $\bb{F}_l$- or $\overline{\bb{F}}_l$-sheaf. We have, in all cases, the following natural isomorphisms that are equivariant under the actions of $\gal(\overline{\eta}/\eta)$:
\begin{equation}\label{eq-nearby-min}
R\Psi_{\ca{S}_{K_2}^{\mmin}}RJ_{\eta,*}^{\mmin}\ca{V}\xrightarrow{\sim} RJ_{\overline{s},*}^{\mmin}R\Psi_{\ca{S}_{K_2}}\ca{V},
\end{equation}
\begin{equation}\label{eq-nearby-min-shriek}
J_{\overline{s},!}^{\mmin}R\Psi_{\ca{S}_{K_2}}\ca{V}\xrightarrow{\sim}R\Psi_{\ca{S}_{K_2}^{\mmin}}J_{\eta,!}^{\mmin}\ca{V},   
\end{equation}
and
\begin{equation}\label{eq-nearby-intersection}
R\Psi_{\ca{S}_{K_2}^\mmin}J_{,\eta,!*}^{\mmin}\ca{V}[d]\xrightarrow{\sim}J_{\overline{s},!*}^\mmin R\Psi_{\ca{S}_{K_2}}\ca{V}[d].
\end{equation}
\end{cor}
\subsubsection{Remarks on future projects}\label{subsubsec-future-projects}
Now we have established a good compactification theory in which $\oint_{K_2}^{\Sigma_2}:\ca{S}_{K_2}^{\Sigma_2}\to \ca{S}_{K_2}^\mmin$ has the desired properties as in \cite[Prop. 2.2]{LS18i} and \cite[Prop. 2.1.2]{LS18b}. 
In some future projects, we plan to study some aspects that have not been addressed in this paper.
\begin{rk}\label{rk-future-work}
In a forthcoming joint work with Shengkai Mao, we will study the extension of $p$-adic shtukas to toroidal compactifications, where a framework generalizing \cite{PR24} to show the well-positionedness of many types of strata on the special fibers of Shimura varieties, and to describe the uniqueness and functoriality of integral models of toroidal compactifications is pursued.
\end{rk}
\begin{rk}\label{rk-fj-expansion}
Unlike Hodge-type Shimura varieties, a general abelian-type Shimura variety does not automatically have a ``Hodge invertible sheaf'' determined by the pullback of some universal abelian variety under a Hodge embedding. Hence, this requires a more systematic treatment which is different from that in \cite[Sec. 5]{Mad19}.\par 
We plan to study the extension of canonical sheaves and ``Hodge invertible sheaves'' to integral models of minimal compactifications, as well as the construction of Fourier-Jacobi expansions. In addition, we also want to extend the main theorems of this paper to more general setups such as allowing all projective cone decompositions for $(G_2,X_2,K_2)$ and constructing toroidal compactifications as normalized blow-ups.
\end{rk}

\newpage

\appendix
\section{Comparison with PEL-type cusp labels}\label{cpr-pel-cl}
Let us review Lan's definition of (PEL-type) cusp labels and its relation, in the case of PEL-type Shimura varieties, with the definition of cusp labels introduced in \cite{Pin89} and \cite{Mad19}, and summarized in \S\ref{cusp-label}. 
See \cite{Lan13}, \cite{Lan12b} and the appendix of \cite{LW15}.
\subsection{PEL-type Shimura varieties} 
Let us fix a PEL-type $\ca{O}$-lattice $\ca{D}:=(V_\bb{Z},\psi_\bb{Z},h_0)$ as defined in \cite[Def. 1.2.1.3]{Lan13}, where $\ca{O}$ is an order of a finite semisimple $\bb{Q}$-algebra with a positive involution $*$, $V_\bb{Z}$ is a $\bb{Z}$-module with $\ca{O}$-action, $\psi_\bb{Z}$ is an alternating pairing of $V_\bb{Z}$ such that $\psi_\bb{Z}(x,by)=\psi_\bb{Z}(b^*x,y)$ and $h_0$ is a Hodge cocharacter $h_0:\bb{C}\lra\Endo_{\ca{O}\otimes\bb{R}}V_\bb{R}$. Fix a \emph{good prime} $p$ as defined in \cite[Def. 1.4.1.1]{Lan13}.\par
Let $G_\bb{Z}$ be the similitude group over $\bb{Z}$ defined by assigning any $\bb{Z}$-algebra $R$ a group $G_\bb{Z}(R):=\{(g,r)\in\mrm{GL}_{\ca{O}\otimes R}(V_R)\times \bb{G}_m(R)|\psi_{R}(gx,gy)=r\psi_{R}(x,y),x,y\in V_R\}$. Let $r:=\nu(g)$, the map $\nu$ is the similitude character of $G_\bb{Z}$. We have $G_\bb{Q}:=G_\bb{Z}\otimes \bb{Q}$ is quasi-split and unramified at $p$ since $\ca{O}\otimes\bb{Q}$ is unramified at $p$ by assumption. \par
Let $X$ be a $G(\bb{R})$-orbit of $h_0$. Then we form a pair $(G_\bb{Q},X)$. Note that $G_\bb{Q}$ is not connected in general when $G^\ad_\bb{Q}$ has simple factors of type $D$, and $G^\ad_\bb{Q}$ might have compact $\bb{Q}$-simple factors; however, we can abusively call such a pair $(G_\bb{Q},X)$ a \textbf{PEL-type Shimura datum}. Define the \textbf{PEL-type Shimura variety with level $K$} by the double coset
$$ \sh_K(G_\bb{Q},X)(\bb{C}):=G(\bb{Q})\bss X\times G(\A)/K.$$
The complex orbifold $\sh_K(G_\bb{Q},X)(\bb{C})$ is canonically the analytification of a smooth quasi-projective algebraic variety $\sh_{K,\bb{C}}$ over $\bb{C}$ if $K$ is neat.\par 
Assume that $K=K_pK^p$ where $K_p=G_\bb{Z}(\bb{Z}_p)$ is the stabilizer of $V_{\bb{Z}_p}$ in $G(\bb{Q}_p)$ and $K^p\sbst G_\bb{Z}(\zhp)$ is a neat open compact subgroup. Those assumptions imply that $K$ stabilizes $V_{\wat{\bb{Z}}}$. Let $\p$ be $\emptyset$ or $\{p\}$.\par
Let us denote by $\mbf{M}_{\ca{D},K^\p}$ the \textbf{PEL-type moduli problem} defined as in \cite[Def. 1.4.1.4]{Lan13} by the PEL-type $\ca{O}$-lattice $\ca{D}$ and the level $K^\p$; if 
$K$ is neat, $\mbf{M}_{\ca{D},K^\p}$ is representable by a smooth and quasi-projective scheme over $\spec \ca{O}_{F_0,(\p)}$, where $F_0$ is the reflex field of $\ca{D}$. For any locally Noetherian and connected scheme $S$ over $\spec \ca{O}_{F_0,(\p)}$, the objects in $\mbf{M}_{\ca{D},K^\p}(S)$ consists of the tuples $(\ca{A}_S,\lambda_S,i_S,[\varepsilon_S]_{K^\p})$, where $\ca{A}_S$ is an abelian scheme over $S$ with a prime-to-$\p$ polarization $\lambda_S$ and an endomorphism structure 
$i_S:\ca{O}\lra \Endo_S(\ca{A}_S)$, and where $[\varepsilon_S]_{K^\p}$ is an (integral) $K^\p$-level structure; the Lie algebra $\ull_{\ca{A}_S/S}$ satisfies the determinantal condition given by $\ca{D}$. 
Denote by $(\ca{A},\lambda,i,[\varepsilon]_{K^p})$ the universal family over $\mbf{M}_{\ca{D},K^\p}$.\par
On the other hand, the variation of Hodge structures defines a holomorphic family $(\ca{A}_{\hol},\lambda_{\hol},i_{\hol},[\varepsilon_{\hol}]_{K})$ over $\sh_K(\bb{C})$. This holomorphic family descends to an algebraic family over $\sh_{K,\bb{C}}$:
\begin{prop}[{\cite[p.16]{Lan12b}}]
There is a canonical open and closed embedding $i:\sh_{K,\bb{C}}\hookrightarrow \mbf{M}_{\ca{D},K,\bb{C}}$, such that the holomorphic family $(\ca{A}_{\hol},\lambda_{\hol},i_{\hol},[\varepsilon_{\hol}]_{K})$ over $\sh_{K}(\bb{C})$ is canonically isomorphic to the analytification of the pullback of the universal family $(\ca{A},\lambda,i,[\varepsilon]_{K})$ defined over $\mbf{M}_{\ca{D},K,\bb{C}}$ by the representability of moduli problem.
\end{prop}
\subsection{PEL-type cusp labels}
\begin{definition}[{\cite[Def. 5.4.1.9 and Def. 5.4.2.4]{Lan13}; see also \cite[A.4]{LW15}}]A \textbf{(PEL-type) cusp label} of $\mbf{M}_{\ca{D},K^\p}$ is denoted by $[(Z_{K^\p},\Phi_{K^\p},\delta_{K^\p})]$, which is an equivalence class of $K^\p$-orbits of tuples of the form $(Z,\Phi,\delta)$, where:
\begin{itemize}
\item $Z=\{Z_i\}_{i={-2}}^0$ is an admissible and fully symplectic (ascending) filtration on $V_{\zhpp}$. We denote by $\gr^Z_i$ the graded pieces defined by $Z$, for $i=0,-1,-2$ (see \cite[Def. 1.2.6.6 and Def. 5.2.7.1]{Lan13});
\item A torus argument $\Phi=(X,Y,\phi,\varphi_{-2},\varphi_{0})$: $X$ and $Y$ are $\bb{Z}$-lattices with $\ca{O}$-action, $\phi:Y\hookrightarrow X$ is an $\ca{O}$-equivariant embedding with finite cokernel, $\varphi_{-2}:\gr^Z_{-2}\xrightarrow{\sim}\Hom (X\otimes {\zhpp},{\zhpp}(1))$ and $\varphi_0:\gr^Z_0\xrightarrow{\sim}Y\otimes{\zhpp}$. The pairing $\langle-,-\rangle_{20}$ induced by $Z$ and $\psi_{\zhpp}$ on $\gr^Z_{-2}\times\gr^Z_0$ coincides with the pullback of the pairing $\langle-,-\rangle_\phi$ induced by $\phi$ under $\varphi_{-2}\times \varphi_0$;
\item $\delta$ is an $\ca{O}$-equivariant splitting $\delta: \bigoplus\limits_{i=-2}^0\gr^Z_i\xrightarrow{\sim}V_{\zhpp}$. 
\end{itemize}
Two orbits $(Z_{K^\p},\Phi_{K^\p},\delta_{K^\p})$ and $(Z_{K^\p}^\prime,\Phi_{K^\p}^\prime,\delta_{K^\p}')$ are equivalent if $Z_{K^\p}=Z'_{K^\p}$ and there is a pair of $\ca{O}$-equivariant isomorphisms $(\gamma_X:X'\xrightarrow{\sim}X,\gamma_Y:Y\xrightarrow{\sim}Y')$ sending $\Phi_{K^\p}$ to $\Phi^\prime_{K^\p}$ in the sense of \cite[Def. 5.4.2.2]{Lan13}.\par
The orbit $(Z_{K^\p},\Phi_{K^\p},\delta_{K^\p})$ itself is called a \textbf{(PEL-type) cusp label representative} of the cusp label $[(Z_{K^\p},\Phi_{K^\p},\delta_{K^\p})]$.
\end{definition}
Denote $V:=V_\bb{Z}\otimes \bb{Q}$. By \cite[Lem. 3.1.1]{Lan12b}, there is a bijection between admissible $\bb{Q}$-parabolic subgroups of $G_\bb{Q}$ and non-trivial minimal admissible symplectic filtrations $W=\{W_i\}_{i=-2}^0$ on $V$; such filtrations are completely determined by their maximal isotropic graded pieces $\gr^W_{-2}$.\par
Let $\mathscr{L}:=(Q,g)$ be a pair of an admissible $\bb{Q}$-parabolic subgroup $Q$ and an element $g\in G(\A)$. Denote by $W$ the filtration on $V$ determined by $Q$. Following \cite[p.21]{Lan12b}, one can associate a PEL-type cusp label with $\mathscr{L}=(Q,g)$:\par
\begin{itemize}
\item $Z^{(g)}_i:=g^{-1}(gV_{\zhpp}\cap W_i\otimes_\bb{Q}\App)=V_{\zhpp}\cap g^{-1}(W_i\otimes_{\bb{Q}}\App)$. Then $Z^{(g)}_\bullet$ is actually a filtration defined by $g^{-1}Qg$.
Let $W^{(g)}_{i}:=g^{-1}(gV_\bb{Z}\cap (W_i\otimes_\bb{Q}\App)).$
\item Let $V^{(g)}_\bb{Z}:=gV_{\wat{\bb{Z}}}\cap V_\bb{Z}$. Let $F^{(g)}_{i}:=V^{(g)}_\bb{Z}\cap W_i.$ 
Then $X^{(g)}:=\Hom(\gr^{F^{(g)}}_{-2},\bb{Z}(1))$ and $Y^{(g)}:=\gr^{F^{(g)}}_0$. Note that $V^{(g)}_\bb{Z}$ is equipped with a symplectic pairing: $\psi_\bb{Z}^{(g)}: V^{(g)}_\bb{Z}\times V^{(g)}_\bb{Z}\to \bb{Z}(1)$, which is defined by $r^{-1}\psi$, where $r$ is the unique $r\in \bb{Q}^\times_{>0}$ such that $r\cdot u=\nu(g)$ for some $u\in \wat{\bb{Z}}^\times$.
\item $\phi^{(g)}:Y^{(g)}\hookrightarrow X^{(g)}$ is induced by the symplectic pairing $\psi^{(g)}_\bb{Z}$;
$\varphi_{-2}^{(g)}$ is defined by $\gr^{Z^{(g)}}_{-2}\xrightarrow{u^{-1}\circ\gr(g)}\gr^{F^{(g)}}_{-2}\otimes \zhpp\xrightarrow{\sim}\Hom(X^{(g)}\otimes \zhpp,\zhpp(1))$ and $\varphi_0^{(g)}$ is defined by $\gr^{Z^{(g)}}_0\xrightarrow{\gr(g)}\gr_0^{F^{(g)}}\otimes {\zhpp}\xrightarrow{\sim}Y^{(g)}\otimes \zhpp$.
\item $\delta^{(g)}$ is any $\ca{O}$-equivariant splitting $\delta^{(g)}:\bigoplus \gr^{Z^{(g)}}_i\xrightarrow{\sim}V_{\zhpp}$.
\end{itemize}\par
Denote by $(Z^{(g)}_{K^\p},\Phi^{(g)}_{K^\p},\delta^{(g)}_{K^\p})$ the $K^\p$-orbit of $(Z^{(g)},\Phi^{(g)},\delta^{(g)})$, and this orbit is a cusp label representative. Then we have constructed a map $\mrm{CL}^\p$, which maps any $\mathscr{L}=(Q,g)$ to $\Psi_\mathscr{L}:=[(Z^{(g)}_{K^\p},\Phi^{(g)}_{K^\p},\delta^{(g)}_{K^\p})]$. \par
Now we fix an admissible $\bb{Q}$-parabolic subgroup $Q$.
Define a normal $\bb{Q}$-subgroup $\wdtd{P}_Q$ of $Q$ as, for any $\bb{Q}$-algebra $R$, 
$$\wdtd{P}_Q(R):=\{g\in Q(R)|\ g|_{\gr_0^W}=\mrm{Id}\text{ and }g|_{\gr_{-2}^W}=\nu(g)\mrm{Id}\}.$$
When $(G_\bb{Q},X)$ is of Siegel type, we have $\wdtd{P}_Q=P_Q$ by \cite[4.25]{Pin89}; in general, $P_Q\neq \wdtd{P}_Q$ since there might be compact factors in $\wdtd{P}_Q$ even in the type A case.\par
For $\p=\emptyset$ or $\{p\}$, define 
$$I_{\PEL}^\p(Q):=Q(\bb{Q})\wdtd{P}_{Q}(\App)\bss G(\App)/K^\p.$$
Let $G_l:= Q/\wdtd{P}_Q$. Let $\wdtd{P}_{Q,h}$ (resp. $Q_h$) be the Levi quotient of $\wdtd{P}_Q$ (resp. $Q$). Then in the PEL-type case, there is a natural exact sequence 
\begin{equation}\label{eq-exact-pel-gp}1\to \wdtd{P}_{Q,h} \to Q_h\to G_l\to 1\end{equation}
and $G_{l}\iso \mrm{GL}(\gr_0^W)\times \mrm{GL}(\gr_{-2}^W)$.\par
From now on, we \textbf{assume} that $G_\bb{Q}$ is connected, that is, $G^\ad_\bb{Q}$ (or equivalently, $\ca{O}\otimes \bb{Q}$) has no simple factors of type $D$.\par
\begin{lem}\label{lem-cusp-lan}
Under the assumption above, every PEL-type cusp label representative comes from a $\mrm{CL}^p(Q,g)=\mrm{CL}(Q,g)$ of some admissible $\bb{Q}$-parabolic subgroup $Q$ and some $g\in G(\A)$.
\end{lem}
\begin{proof}
    This is \cite[Lem. A.4.6 and Prop. A.5.9]{LW15}.
\end{proof}
\begin{lem}\label{lem-comparison-pink-pel}
Under the assumption above, the map 
$$G(\A)\to \mrm{CL}(Q,-)$$ 
factors through a bijection $I_{\PEL}(Q)\xrightarrow{\sim} \mrm{CL}(Q,-)$.\par
If $\wdtd{P}_Q=P_Q$, we have $I_{\PEL}(Q)\iso I(Q)$; this is the case when $(G_\bb{Q},X)$ is of Siegel type by \cite[4.25]{Pin89}. Moreover, in this case, there is a bijection between the set of PEL-type cusp labels and the set of cusp labels defined by \cite{Pin89} and \cite{Mad19}.
\end{lem}
\begin{proof}
It remains to show the first statement (see \cite[B.3]{HLTT16} for a specific case).
Firstly, the map $G(\A)\to\mrm{CL}(Q,-)$ sending $g$ to $\mrm{CL}(Q,g)$ factors through $I_\mrm{PEL}(Q)$: for any $p\in \wdtd{P}_Q(\A)$, $Z^{(pg)}=Z^{(g)}$ by the definition above and $\Phi^{(pg)}=\Phi^{(g)}$ by (\ref{eq-exact-pel-gp}); $q\in Q(\bb{Q})$ gives a pair of $\ca{O}$-equivariant isomorphisms $(\gamma_X,\gamma_Y)$. Secondly, assume that, for some $h\in G(\A)$, $\mrm{CL}(Q,g)=\mrm{CL}(Q,hg)$. Then up to replacing $h$ with $hgkg^{-1}$ for some $k\in K$, we have $Z^{(hg)}=Z^{(g)}$ so $h\in Q(\A)$. Moreover, there is a pair of $\ca{O}$-equivariant isomorphisms $(\gamma_X,\gamma_Y)$ sending $\Phi^{(g)}_K$ to $\Phi^{(hg)}_K$. Then $h\in Q(\A)$ projects to $G_l(\bb{Q})$ via (\ref{eq-exact-pel-gp}) by the construction of $X^{(g)}$ and $Y^{(g)}$. The remaining part of the statement follows from the fact that $Q_h\iso \wdtd{P}_{Q,h}\times G_l$ in the PEL-type case and Lemma \ref{lem-cusp-lan} above.
\end{proof}
\section{Isogenies}\label{comm-lem}
Let $\square$ be a set of prime numbers, which is still $\{p\}$ or $\emptyset$ in this appendix.\par
\subsection{Quasi-isogenies}
Let $\ca{C}$ be the category of abelian schemes over $S$, or the category of $1$-motives over $S$.
Recall that we have defined isogenies for $1$-motives; see Definition \ref{def-isog-1-mot}. Then it is also possible to define $\zbkppt$-isogenies for $1$-motives.
\begin{definition}\label{def-pp-isog}
Let $\Q_1$ and $\Q_2$ be two objects in $\ca{C}$. 
Let $f:\Q_1\to \Q_2$ be an isogeny. Since $\ker f$ is finite flat over $S$, the rank of $\ker f$ is locally constant in both cases. We say $f$ is a \textbf{prime-to-$\p$ isogeny} if the rank of $\ker f$ over any connected component of $S$ is prime to $\p$.\par
A \emph{$\zbkppt$-isogeny} $g$ from $\Q_1$ to $\Q_2$ is represented by a tuple $(f_1,f_2,\Q_3)$ consisting of a $\Q_3\in\ob\ca{C}$ and two prime-to-$\p$ isogenies $f_1:\Q_3\to \Q_1$ and $f_2:\Q_3\to \Q_2$, such that, $g$ is the diagram $g: \Q_1\xleftarrow{f_1}\Q_3\xrightarrow{f_2}\Q_2$. If $f_1=[n]$, we denote $g:=\frac{1}{n}f_2$. Any two tuples $(f_1,f_2,\Q_3)$ and $(h_1,h_2,\ca{R}_3)$ are equivalent if there is another tuple $(c_1,c_2,\T_3)$ such that $c_1$ and $c_2$ are prime-to-$\p$ isogenies, and there are prime-to-$\p$ isogenies $e_1:\T_3\to \ca{R}_3$ and $e_2:\T_3\to \Q_3$, which fit into the following commutative diagram:
\begin{equation}
    \begin{tikzcd}
    & \Q_3\arrow[dl,"f_1"']\arrow[dr,"f_2"]&\\
\Q_1&\T_3\arrow[l,"c_1"']\arrow[r,"c_2"]\arrow[d,"e_1"]\arrow[u,"e_2"']&\Q_2\\
& \ca{R}_3\arrow[ul,"h_1"']\arrow[ur,"h_2"].&
    \end{tikzcd}
\end{equation}
By Lemma \ref{lem-isog-1-mot} in the case of $1$-motives, for any tuple $(f_1,f_2,\Q_3)$ representing $g$, there is a tuple $([N],g_2,\Q_1)$, for some $(N,\p)=1$ and some prime-to-$\p$ $g_2$, that is equivalent to $(f_1,f_2,\Q_3)$.\par
The \textbf{composition} of two $\zbkppt$-isogenies $f=([N_1],f_1,\Q_1):\Q_1\to \Q_2$ and $g=([N_2],g_1,\Q_2):\Q_2\to \Q_3$ is represented by $g\circ f=([N_1N_2],g_1\circ f_1): \Q_1\to \Q_3$.
\end{definition}
By Lemma \ref{lem-isog-1-mot}, for any $\zbkppt$-isogeny between $1$-motives $f:\Q_1\to \Q_2$, there is a prime-to-$\p$ isogeny $g: \Q_2\to \Q_1$ and a prime-to-$\p$ positive integer $M$ such that $f\circ g=[M]$ and $g\circ f=[M]$. Hence, we define the $\zbkppt$-isogeny $f^{-1}:= \frac{1}{M}g$ to be the \textbf{inverse} of $f$. For the case of abelian schemes, see, e.g., \cite[Sec. 1.3.1]{Lan13}.\par
In the categorical language, let $T$ be the set of arrows consisting of arrows consisting of prime-to-$\p$ isogenies in $\mrm{Mor}(\ca{C})$; alternatively, we can let $T$ be the set of arrows consisting of all $[N]$, for 
$N>0$ and $(N,\p)=1$. Then the \emph{isogeny category} is defined to be the localization $\ca{C}^\p:=T^{-1}\ca{C}$.

\subsection{Polarizations}
\subsubsection{}\label{subsubsec-zppol-1-mot}
Let $S$ be a $\zbkp$-scheme. Let $\Q$ be a $1$-motive over $S$. Denote by $[\Q]$ the fppf sheaf representing $\Q$ up to $\zbkppt$-isogenies.\par
For any $f\in \aut_{\zbkpp}([\Q])$, $f$ naturally induces an automorphism $f^\vee\in \aut_{\zbkpp}([\Q^\vee])$ by the definition of Cartier dual and dual abelian schemes.\par 
Then there is a canonical isomorphism:
\begin{equation}\label{eq-aut-dual-iso}
    \begin{split}
    \aut_{\zbkpp}([\Q])&\xrightarrow{\sim}\aut_{\zbkpp}([\Q^\vee]);\\
    f &\mapsto f^\bigstar:=(f^\vee)^{-1}.
    \end{split}
\end{equation}
We say a polarization $\bml:\Q\to\Q^\vee$ is \textbf{prime-to-$\p$} if it is a prime-to-$\p$ isogeny. We say that a $\zbkppt$-isogeny $\bml:\Q\to\Q^\vee$ is a \textbf{$\zbkppt$-polarization} if $N\bml$ is a prime-to-$\p$ polarization for some positive prime-to-$\p$ integer $N$.
\subsubsection{}
Let $\ca{B}$ be another $1$-motive over a $\zbkpp$-scheme $S$. 
Suppose that there is a prime-to-$\p$ isogeny $f:\Q\to \ca{B}$. For any $\zbkpp$-algebra $R$ and $z\in \mrm{End}_{\zbkpp}([\Q])(R)$, $f$ induces an endomorphism $z^\prime:=f\circ z\circ f^{-1}\in \mrm{End}_{\zbkpp}([\ca{B}])(R)$. We will use the notation in \S\ref{subsec-twt-1-mot}.
\begin{lem}\label{twt-isog-gen}
There is a unique homomorphism $f_N^{\tilde{\gamma},\prime}:\ca{B}^{\ca{O}_F}\lra(\ca{B}^{\ca{O}_F})^{\oplus|\mrm{Gal}(F/\bb{Q})|-1}$ such that the diagram
\begin{equation}\begin{tikzcd}
\Q^{\ca{O}_F}\arrow[rr,"f"]\arrow[d,"f^{\tilde{\gamma}}_N"]&&\ca{B}^{\ca{O}_F}\arrow[d,"f^{\tilde{\gamma},\prime}_N"]\\
(\Q^{\ca{O}_F})^{\oplus|\mrm{Gal}(F/\bb{Q})|-1}\arrow[rr,"f^{\oplus |\mrm{Gal}(F/\bb{Q})|-1}"]&&(\ca{B}^{\ca{O}_F})^{\oplus| \mrm{Gal}(F/\bb{Q})|-1}
\end{tikzcd}\end{equation}commutes. The diagram above induces a prime-to-$p$ isogeny $f^{\tilde{\gamma}}:\Q^{\tg} \lra\ca{B}^{\tg}$. The isogeny $f^{\tg}$ does not depend on $N$. \end{lem}
\begin{prop}[{cf.\cite[Lem. 4.4.8]{KP15}}]\label{twt-isog-iso-gen}
With the setting in Lemma \ref{twt-isog-gen}, there is a commutative diagram of isomorphisms between fppf sheaves 
\begin{equation}
    \begin{tikzcd}
{\Q^{\tg}\otimes \ca{O}_{F,(p)}}\arrow[r,"f^{\tg}","\sim"']\arrow[d,sloped,"\sim"]&
    {\ca{B}^{\tg}\otimes \ca{O}_{F,(p)}}\arrow[d,sloped,"\sim"]\\
    {\Q\otimes\ca{O}_{F,(p)}}\arrow[r,"f","\sim"']&{\ca{B}\otimes \ca{O}_{F,(p)}}.
    \end{tikzcd}
\end{equation}
\end{prop}
\begin{proof}
This follows directly from Lemma \ref{twt-isog-gen}.
\end{proof}
\newpage
\printnoidxglossary[type=symbols,title={List of symbols in the first two {sections}}]
\newpage
\bibliography{References}
\end{document}